\documentclass[11pt,a4]{article}
\usepackage[T2A]{fontenc}
\usepackage[cp866]{inputenc}
\usepackage[mathscr]{eucal}
\usepackage{amssymb}
\usepackage{amsmath}
\usepackage{graphicx}
\usepackage[russian]{babel}  % затем убрать
\usepackage{array}  % добавили новый пакет для таблицы

\renewcommand{\leq}{\leqslant}
\renewcommand{\geq}{\geqslant}

\renewcommand{\C}{{\mathbb C}}
\newcommand{\N}{{\mathbb N}}
\newcommand{\R}{{\mathbb R}}

\renewcommand{\k}{\boxtimes}
\usepackage{fancyhdr}
\begin{document}
\sloppy

\lhead
    [\scriptsize V.N. Gorbuzov]
    {\scriptsize V.N. Gorbuzov}
\rhead
    [\it \scriptsize Partial integrals of ordinary differential systems]
    {\it \scriptsize Partial integrals of ordinary differential systems}

{\normalsize

\thispagestyle{empty}

\mbox{}
\\[-0.15ex]
\centerline{
{\large
\bf
PARTIAL INTEGRALS OF ORDINARY DIFFERENTIAL SYSTEMS 
}
}
\\[2.5ex]
\centerline{
\bf 
V.N. Gorbuzov
}
\\[2ex]
\centerline{
\it
Department of Mathematics and Computer Science,
}
\\[1ex]
\centerline{
\it
Yanka Kupala State University of Grodno,
}
\\[1ex]
\centerline{
\it
Ozeshko {\rm 22}, Grodno, {\rm 230023}, Belarus
}
\\[1.5ex]
\centerline{
E-mail: gorbuzov@grsu.by
}
\\[5.5ex]
\centerline{{\large\bf Abstract}}
\\[1ex]
\indent
Properties of partial integrals such as real and complex-valued polynomial, multiple polynomial, exponential, and conditional for ordinary differential systems are studied.
The possibilities of constructing first integrals and last multipliers by known partial integrals are considered.
Applications of partial integrals to solve the Darboux problem and to the exten\-ded Darboux problem are given.
Integral basis of the Jacobi system is built. And 
the in\-ver\-se problem of constructing differential systems on the base of their partial integrals is solved.
\\[1.5ex]
\\[1.25ex]
\\[7.5ex]
\centerline{{\large\bf Contents}}
\\[1.5ex]
{\bf  Introduction}                   \dotfill\ 2
\\[1ex]
{\bf \S 1. Partial integrals}
                                                 \dotfill \ 6
\\[0.5ex]
\mbox{}\hspace{1.35em}
1. Partial integral. Definition. Properties
                                                 \dotfill \ 6
\\[0.5ex]
\mbox{}\hspace{1.35em}
2. Polynomial partial integrals
                                                 \dotfill \ 10
\\[0.5ex]
\mbox{}\hspace{1.35em}
3. Exponential partial integrals
                                                 \dotfill \ 12
\\[0.5ex]
\mbox{}\hspace{1.35em}
4. Conditional partial integrals
                                                 \dotfill \ 16
\\[0.5ex]
\mbox{}\hspace{1.35em}
5. Multiple polynomial partial integrals
                                                 \dotfill \ 17
\\[0.5ex]
\mbox{}\hspace{1.35em}
6. Complex-valued polynomial partial integrals
                                                 \dotfill \ 23
\\[0.5ex]
\mbox{}\hspace{1.35em}
7. Multiple complex-valued polynomial partial integrals
                                                 \dotfill \ 28
\\[1ex]
\noindent
{\bf \S 2. 
Last multipliers
}
                                                 \dotfill \ 31
\\[0.75ex]
\mbox{}\hspace{1.35em}
8. Last multiplier as partial integral
                                                 \dotfill \ 31
\\[0.5ex]
\mbox{}\hspace{1.35em}
9. Building of last multipliers on base of polynomial partial integrals 
                                                 \dotfill \ 35
\\[0.5ex]
\mbox{}\hspace{1.35em}
10. Exponential last multipliers
                                                 \dotfill \ 40
\\[1ex]
\noindent
{\bf \S 3. First integrals
}
                                                 \dotfill \ 46
\\[0.75ex]
\mbox{}\hspace{1.35em}
11. First integrals defined by partial integrals and last multipliers
                                                 \dotfill \ 46
\\[0.5ex]
\mbox{}\hspace{1.35em}
12. Applications
                                                 \dotfill \ 54
\\[1ex]
{\bf References}
                                              \dotfill \ 83

\newpage

\mbox{}
\\[-2.25ex]
\centerline{\large\bf  Introduction}
\\[1.5ex]
\indent
{\bf Object of research and statement of problem}.
Consider an $n\!$-th order normal ordinary differential system 
\\[1.75ex]
\mbox{}\hfill                             %(0.1)
$
\dfrac{dx_i^{}}{dt}=X_i^{}(t,x),
\quad 
i = 1,\ldots,n,
$
\hfill(0.1)
\\[2.25ex]
where $X_{i}^{}\colon \Xi\to\R,\ i=1,\ldots,n,\ \Xi = T\times\R^n,\ T\subset\R,$ 
\vspace{0.35ex}
are polynomials in the dependent variables $x_1^{},\ldots, x_n^{}$ 
\vspace{0.25ex}
with continuously differentiable on the domain $T$ coefficients-functions
in the independent variable $t.$
\vspace{0.25ex}
Moreover, we assume that the degrees of polynomials $X_{i}^{}$ such that
$
d=\max\bigl\{\deg_{\;\!x}^{} X_{i}^{}\colon i=1,\ldots, n\bigr\}\geq 1.
$ 
\vspace{0.5ex}

For the system (0.1) we introduce the notion of {\it partial integral}, 
with respect to which we solve the series of problems, such as: 
1) existence;
2) analytical structure;
3) geometric interpretation;
4) properties;
5) relationship with last multipliers;
6) construction of first integrals by known partial integrals.
\vspace{1ex}

{\bf General provisions}.
For the purpose of an unambiguous interpretation, we will define the notions used, 
formulate the basic statements, specify the accepted terminology, and introduce conventions.
At the same time, we will mainly follow the approaches,
described in the monograph [1] and in the article [2].
\vspace{0.15ex}

We recall that by {\it domain} we mean open arcwise connected set.
\vspace{0.15ex}

Without special provisions, we assume that domains
\vspace{0.25ex}
$\Xi^{\;\!\prime}$ and $\Omega$ are subsets of the domain 
$\Xi = T\times\R^n,\ T\subset\R,$ such that  
\\[1.5ex]
\mbox{}\hfill 
$
\Xi^{\;\!\prime} = T^{\;\!\prime}\times\R^n,\ \ \,
T^{\;\!\prime}\subset T,
$ 
\ and \
$
\Omega = T^{\;\!\prime}\times X^{\;\!\prime},\ \ \,
T^{\;\!\prime}\subset T,\ \ \,
X^{\;\!\prime}\subset\R^n.
\hfill
$ 
\\[1.5ex]
\indent
The set of functions that are polynomials in the variables $x_1^{},\ldots,x_n^{}$
\vspace{0.15ex}
with continuously differentiable coefficients-functions in the variable $t$ on the domain $T$ 
is denoted by $\text{P}_{_{\!\Xi}}.$  
\vspace{0.35ex}

Functions $p$ and $q$ from the set $\text{P}_{_{\!\Xi}}$ are called {\it relatively prime} 
\vspace{0.35ex}
if there is no function $u\in \text{P}_{_{\!\Xi}},\ u\not\equiv \text{const},$ such that
\\[2ex]
\mbox{}\hfill                            
$
\displaystyle
p(t,x)=u(t,x)\;\! p_1^{}(t,x),
\quad 
q(t,x)=u(t,x)\;\! q_1^{}(t,x)
$  
\ for all $(t,x)\in \Xi,
\quad
p_1^{},\, q_1^{}\in \text{P}_{_{\!\Xi}}.
\hfill
$
\\[2ex]
\indent
The system (0.1) is induced the linear differential ope\-ra\-tor of first order
\\[2ex]
\mbox{}\hfill                            
$
\displaystyle
{\frak d}(t,x) =\partial_{{}_{\scriptstyle t}} +
\sum\limits_{i=1}^n\, X_{i}^{}(t,x)\, \partial_{{}_{\scriptstyle x_i^{}}}
$ 
\ for all $(t,x)\in \Xi.
\hfill
$
\\[2ex]
We say that [1, p. 20] this operator is  the 
\vspace{0.5ex}
{\it operator of differentiation by virtue of sys\-tem} (0.1).

The set of continuously differentiable functions on the domain $\Omega$ we denote by $C^{1}\Omega.$
 \vspace{0.35ex}
 
Let $f\in C^{1}\Omega.$ Then the function 
\\[2ex]
\mbox{}\hfill                            
$
\displaystyle
{\frak d}f\colon (t,x)\to\ 
\partial_{{}_{\scriptstyle t}}f(t,x) +
\sum\limits_{i=1}^n\, X_{i}^{}(t,x)\, \partial_{{}_{\scriptstyle x_i^{}}}f(t,x)
$ 
\ for all 
$(t,x)\in \Omega
\hfill
$
\\[2ex]
is said to be the {\it derivative by virtue of system} (0.1) for the function $f;$ 
and if $f(t,x)> 0$ for all $(t,x)\in \Omega,$ then the function  
\\[1.75ex]
\mbox{}\hfill                            
$
\displaystyle
{\frak d}\ln f\colon (t,x)\to\ 
{\frak d}\ln f(t,x)
$ 
\ for all $(t,x)\in \Omega
\hfill
$
\\[1.75ex]
is called the {\it logarithmic derivative by virtue of system} (0.1) for the function $f.$
 \vspace{1ex}
 
{\bf Definition 0.1} [3, p. 256; 4, pp. 129 -- 132].
\vspace{0.35ex}
{\it
A function $F\in C^{1}\Omega$ is said to be a 
\textit{\textbf{first integral on the domain}} $\Omega$ of system {\rm (0.1)}, 
 \vspace{0.35ex}
if this function is constant along any solution  
$x\colon t\to x(t)$ for all $t\in T_{_{0}}^{}\subset T^{\;\!\prime}$ 
 \vspace{0.5ex}
of system {\rm (0.1)} such that the points $(t,x(t))\in \Omega$ for all $t\in T_{_{0}}^{},$  
i.e.,
$F(t,x(t))=C$ for all $t\in T_{_{0}}^{},\ C\in\R.$
}
 \vspace{0.75ex}

By $\text{I}_{_{\Omega}}$ denote the set of all first integrals on the domain $\Omega$ of system (0.1).
\vspace{0.5ex}

Let us give the existence criterion of first integral, 
which is sometimes taken for the definition of first integral.
 \vspace{0.5ex}
 
{\bf Theorem 0.1} [1, p. 26; 5, p. 337].
\vspace{0.35ex}
{\it
A function $F\in C^{1}\Omega$ is a first integral on the domain $\Omega$ of system {\rm (0.1)} 
if and only if the differential by virtue of system {\rm (0.1)} 
\vspace{0.35ex}
for this function is identically zero on the domain $\Omega\colon$}
\\[1.75ex]
\mbox{}\hfill
$
\displaystyle
dF(t,x)_{\displaystyle |_{(0.1)}} = \Bigl(\partial_{{}_{\scriptstyle t}}F(t,x)\;\! dt +
\sum\limits_{i=1}^n\, \partial_{{}_{\scriptstyle x_i^{}}}F(t,x)\;\! dx_i^{}\Bigr)_{\displaystyle |_{(0.1)}} =
\hfill
$
\\
\mbox{}\hfill (0.2)
\\
\mbox{}\hfill
$
\displaystyle
=\Bigl(\partial_{{}_{\scriptstyle t}}F(t,x) +
\sum\limits_{i=1}^n\, X_{i}^{}(t,x)\, \partial_{{}_{\scriptstyle x_i^{}}}F(t,x)\Bigr)dt =0
$ 
\ for all 
$(t,x)\in\Omega.
\hfill
$
\\[2ex]
\indent
Taking into account the relation between 
the differential by virtue of system (0.1) and the derivative by the virtue of (0.1), 
we replace the differential identity (0.2) by the operator identity.
\vspace{0.5ex}

{\bf Theorem 0.2} [1, p. 26; 2].
\vspace{0.35ex}
{\it
A function $F\in C^{1}\Omega$ is a first integral on the domain $\Omega$ of system {\rm (0.1)} 
if and only if the derivative by virtue of system {\rm (0.1)} 
\vspace{0.35ex}
for this function is identically zero on the domain $\Omega\colon$
\\[2ex]
\mbox{}\hfill                             %(0.3)
$
{\frak d}F(t,x)=0
$ 
\ for all} $(t,x)\in \Omega.
$
\hfill(0.3)
\\[2.5ex]
\indent 
{\bf Property 0.1} [1, pp. 28 -- 29; 3, p. 262].
{\it
If $F_\nu^{}\in \text{\rm I}_{_{\Omega}}\;\!,\ \nu=1,\ldots,k,$ then the function 
\\[2ex]
\mbox{}\hfill
$
\Psi\colon (t,x)\to\Phi(F_1^{}(t,x),\ldots,F_k^{}(t,x))
$
\ for all $(t,x)\in\Omega,
\hfill
$ 
\\[2ex]
where $\Phi$ is arbitrary continuously differentiable function,
is also a first integral on the domain $\Omega$ of system}~(0.1).
 \vspace{0.75ex}
 
This property under $k=1$ expresses the functional ambiguity of first integral.
\vspace{1ex}
 
{\bf Definition 0.2} [1, p. 29].
\vspace{0.35ex}
 {\it
A set of the functionally independent on the domain $\Omega$ first integrals 
 \vspace{0.35ex}
 $F_\nu^{}\colon \Omega\to\R,\ \nu=1,\ldots,k,$ of system {\rm (0.1)} is called a 
\textit{\textbf{basis of first integrals}} {\rm(}or 
\textit{\textbf{integral basis}}{\rm)} 
\textit{\textbf{on the domain}} $\Omega$ of system {\rm (0.1)} if
\vspace{0.35ex}
 any first integral $F\colon \Omega\to\R$ of system {\rm (0.1)} 
can be represented on the domain $\Omega$  in the form
\\[1.5ex]
\mbox{}\hfill
$
F(t,x)=\Phi(F_1^{}(t,x),\ldots,F_k^{}(t,x))
$
\ for all $(t,x)\in\Omega,
\hfill
$ 
\\[1.75ex]
where $\Phi$ is some continuously differentiable function. 
\vspace{0.35ex}
The number $k$ is said to be the \textit{\textbf{dimension}} of 
ba\-sis of first integrals on the domain $\Omega$ of system} {\rm (0.1)}.
\vspace{1.25ex} 
 
{\bf Theorem 0.3} [1, p. 54; 3, p. 367].
\vspace{0.5ex} 
{\it
The system {\rm (0.1)} on an neighbourhood of any point from the domain $\Omega$ 
has a basis of first integrals of dimension $n.$
}
\vspace{1.5ex} 
 
{\bf Definition 0.3} [2].
\vspace{0.35ex} 
{\it
A smooth manifold ${\rm g}(t,x)=0$ is said to be an
\textit{\textbf{integral manifold}} of system {\rm (0.1)} if
\vspace{0.35ex} 
the  differential by virtue of system {\rm(0.1)} for the function ${\rm g}\in C^{1}\Omega$ 
is identically zero on the manifold ${\rm g}(t,x)=0\colon$
\\[2ex]
\mbox{}\hfill                           
$
d\;\!{\rm g}(t,x)_{\displaystyle |_{(0.1)}} = \Phi(t,x)\;\! dt
$ 
\ for all $(t,x)\in \Omega,
$
\hfill {\rm (0.4)}
\\[2ex]
where the function $\Phi\colon \Omega \to \R$ such that
\\[2ex]
\mbox{}\hfill                           
$
\Phi(t,x)_{\displaystyle |_{{\rm g}(t,x)=0}} = 0
$ 
\ for all} $(t,x)\in \Omega.
\hfill
$
\\[2.75ex]
\indent
Along with Definition 0.3 we will use the existence criterion of integral manifold for 
\linebreak
system (0.1).

\newpage

{\bf Theorem 0.4.}
\vspace{0.35ex}
{\it
A smooth manifold ${\rm g}(t,x) = 0$ is an integral manifold of system {\rm(0.1)}
if and only if 
the  derivative by virtue of system {\rm(0.1)} for the function ${\rm g}\colon \Omega\to \R$ 
\vspace{0.25ex}
is identically zero on this manifold{\rm:}
\\[2ex]
\mbox{}\hfill                                    
$
{\frak d}\;\! {\rm g}(t,x)  =  \Phi(t,x), 
\quad 
\Phi(t,x)_{\displaystyle |_{{\rm g}(t,x)=0}} = 0
$
\ 
for all} $(t,x)\in \Omega.
$
\hfill  (0.5)
\\[2.75ex]
\indent
{\bf Definition 0.4.}
\vspace{0.35ex}
{\it
A function $\mu\in C^{1}\Omega$ is called  a
\textit{\textbf{last multiplier on the domain}} $\Omega$  of system {\rm (0.1)} if
the  differential by virtue of system {\rm(0.1)} 
\\[2ex]
\mbox{}\hfill                     
$
d\;\!\mu(t,x)_{\displaystyle |_{(0.1)}}={}-\mu(t,x)\;{\rm div}\, {\frak d}(t,x)\, dt
$ 
\ for all $(t,x)\in \Omega,
\hfill
$
\\[1.25ex]
where
\\[0.75ex]
\mbox{}\hfill
$
\displaystyle
{\rm div}\, {\frak d}(t,x) = 
\sum\limits_{i=1}^n\, \partial_{{}_{\scriptstyle x_i^{}}}X_i^{}(t,x)
$ 
\ for all} $(t,x)\in \Xi.
\hfill
$
\\[2ex]
\indent
Let us give the existence criterion of last multiplier, 
which is sometimes taken for the definition of last multiplier.
\vspace{0.5ex}

{\bf Theorem 0.5} [5, pp. 341 -- 346; 6, p. 117].
\vspace{0.35ex}
{\it
A function $\mu\in C^{1}\Omega$ is a 
last multiplier on the domain $\Omega$  of system {\rm (0.1)}
if and only if 
the derivative by virtue of system {\rm(0.1)} for this function
\\[1ex]
\mbox{}\hfill                         %(0.6)
$
{\frak d}\;\!\mu(t,x)={}-\mu(t,x)\;{\rm div}\, {\frak d}(t,x)
$ 
\ for all} $(t,x)\in \Omega.
$
\hfill(0.6)
\\[2.25ex]
\indent
By $\text{M}_{_{\Omega}}$ denote the set of all last multipliers on the domain $\Omega$ of system (0.1).
\vspace{1.25ex}

{\bf Property 0.2} (Jacobi's property of last multipliers).
 \vspace{0.5ex}
 {\it 
Suppose $\mu_1^{},\mu_2^{}\in \text{\rm M}_{_{\Omega}},$ 
the set 
\linebreak
$\Omega_{_0}\subset \Omega$ such that  
$\mu_2^{}(t,x)\ne 0$ for all $(t,x)\in \Omega_{_0}$ and
 \vspace{0.35ex}
 $\mu_2^{}(t,x)= 0$ for all $(t,x)\in {\sf C}_{_\Omega}\Omega_{_0}.$
Then the function 
\\[2ex]
\mbox{}\hfill
$
F \colon (t,x) \to\ 
\dfrac{\mu_1^{}(t,x)}{\mu_2^{}(t,x)}
$  
\ for all $(t,x)\in \Omega_{_0}
\hfill
$
\\[2.25ex]
is a first integral on any domain from the set $\Omega_{_0}$ for system} (0.1).
\vspace{1.25ex}
 
{\bf Remark 0.1.}
Definitions 0.1 -- 0.4, Theorems 0.1 -- 0.5, and also Properties 0.1 and 0.2 
are true in the more general case, when the right parts of system (0.1) such that
$X_i^{}\in C^{1}\Omega,$ $i=1,\ldots, n.$
\vspace{1.5ex}

To simplify the reading, we give the list of symbols:
\vspace{1ex}

\noindent
$T$ is a domain from $\R;$ 
$T^{\;\!\prime}$ is  a subdomain from the domain $T;$
\vspace{1ex}

\noindent
the domain $\Xi = T\times\R^n;$
the domain $\Xi^{\;\!\prime} = T^{\;\!\prime}\times\R^n;$
\vspace{1ex}

\noindent
$X^{\;\!\prime}$ is a domain from $\R^n;$
the domain $\Omega = T^{\;\!\prime}\times X^{\;\!\prime};$
\vspace{1ex}

\noindent
${\sf C}_{_{\Omega^{\;\!\prime}}}\Omega_{_0}$ 
the complement of set $\Omega_{_0}\subset \Omega^{\;\!\prime}$ 
to the set $\Omega^{\;\!\prime};$
\vspace{1.25ex}

\noindent
${\sf D}$ 
\qquad
the Cauchy symbol of derivative for one-variable function;
\vspace{1ex}

\noindent
$\partial_{t}^{}$ 
\qquad
the partial derivative with respect to the variable $t$ 
(the operator of differentiation 
\\[0.5ex]
\mbox{}\hspace{3.3em}
with respect to the variable $t);$
\vspace{1ex}

\noindent
${\frak d}$
\qquad\,
the differential operator by virtue of system (0.1);
\vspace{1ex}

\noindent
$d\;\!f_{\displaystyle |_{(0.1)}}\!$
the differential of function $f$ by virtue of system (0.1);
\vspace{1.25ex}

\noindent
$
d=\max\bigl\{\deg_{\;\!x}^{} X_{i}^{}\colon i=1,\ldots, n\bigr\}\;\!;
$ 
\vspace{1.25ex}

\noindent
$C^1\Omega$ 
\quad
the set of continuously differentiable functions on the domain $\Omega\;\!;$
\vspace{1ex}

\newpage

\noindent
$\text{P}_{_{\!\Xi}}$
\qquad
the set of functions that are polynomials in the variables $x_1^{},\ldots,x_n^{}$ with continuously
\\[0.5ex]
\mbox{}\hspace{3.5em}
differentiable real coefficients-functions in the variable $t$ on the domain $T;$
\vspace{1ex}

\noindent
$\text{Z}_{_{\Xi}}$ 
\qquad
the set of functions that are polynomials in the variables $x_1^{},\ldots,x_n^{}$ with continuously 
\\[0.5ex]
\mbox{}\hspace{3.5em}
differentiable complex coefficients-functions in the variable $t$ on the domain $T;$
\vspace{1ex}

\noindent
$\text{I}_{_{\Omega}}$ 
\qquad
the set of all first integrals on the domain $\Omega$ of system (0.1); 
\vspace{1ex}

\noindent
$\text{J}_{_{\Omega}}$ 
\qquad
the set of all partial integrals on the domain $\Omega$ of system (0.1); 
\vspace{1ex}

\noindent
${\rm(g}, M)\in \text{J}_{_{\Omega}}$ 
\quad
the function ${\rm g}$ is a partial integral with cofactor $M$ on the domain $\Omega$
\\[0.5ex]
\mbox{}\hspace{6.7em}
of system (0.1); 
\vspace{1ex}

\noindent
$\text{A}_{_{\Xi}}$
\qquad
the set of polynomial (real) partial integrals on the domain $\Xi$ of system (0.1);
\vspace{1ex}

\noindent
$(p, M)\in \text{A}_{_{\Xi}}$
\quad
the function $p$ is a polynomial (real) partial integral with cofactor $M$ 
\\[0.5ex]
\mbox{}\hspace{6.8em}
on the domain $\Xi$ of system (0.1);
\vspace{1ex}

\noindent
$\text{B}_{_{\Xi}}$
\qquad
the set of multiple polynomial (real) partial integrals on the domain $\Xi$
\\[0.5ex]
\mbox{}\hspace{3.7em}
of system (0.1);
\vspace{1ex}

\noindent
$\bigl((p, M), (h,q, N)\bigr)\in \text{B}_{_{\Xi}}$
\ \;\!
the polynomial (real) partial integral $p$ with cofactor $M$ on the 
\\[0.35ex]
\mbox{}\hspace{11.4em}
domain $\Xi$ of system (0.1) 
\vspace{1ex}
is multiple and the identity (5.1) holds;

\noindent
$\text{E}_{_{\Omega}}$
\qquad
the set of exponential partial integrals on the domain $\Omega$ of system (0.1);
\vspace{1ex}

\noindent
$(\exp\omega , M)\in \text{E}_{_{\Omega}}$
\quad
the function $\exp\omega$ is an exponential partial integral with cofactor $M$ 
\\[0.5ex]
\mbox{}\hspace{8.7em}
on the domain $\Omega$ of system (0.1);
\vspace{1ex}

\noindent
$\text{F}_{_{\Xi}}$
\qquad
the set of conditional partial integrals on the domain $\Xi$ of system (0.1);
\vspace{1.25ex}

\noindent
$(\exp p, M)\in \text{F}_{_{\Xi}}$
\quad
the function $\exp p$ is a conditional partial integral with cofactor $M$ 
\\[0.5ex]
\mbox{}\hspace{8.4em}
on the domain $\Xi$ of system (0.1);
\vspace{1ex}

\noindent
$\text{H}_{_{\Xi}}$
\qquad
the set of complex-valued polynomial partial integrals on the domain $\Xi$
\\[0.5ex]
\mbox{}\hspace{3.7em}
of system (0.1); 
\vspace{1ex}

\noindent
$(w, W)\in \text{H}_{_{\Xi}}$ 
\quad
the function $w$ is a complex-valued polynomial partial integral with 
\\[0.35ex]
\mbox{}\hspace{7.05em}
cofactor $W$ on the domain $\Xi$ of system (0.1);
\vspace{1ex}

\noindent
$\text{G}_{_{\Xi}}$
\qquad
the set of multiple complex-valued polynomial partial integrals on the domain $\Xi$
\\[0.5ex]
\mbox{}\hspace{3.7em}
 of system (0.1);
\vspace{1ex}

\noindent
$\bigl((w, W), (h,z, Q)\bigr)\in \text{G}_{_{\Xi}}$
\quad
the complex-valued polynomial partial integral $w$ with 
\\[0.5ex]
\mbox{}\hspace{12em}
cofactor $W$ on the domain $\Xi$ of system (0.1) is multiple  
\\[0.5ex]
\mbox{}\hspace{12em}
and the identity (7.1) holds;
\vspace{1ex}

\noindent
$\text{M}_{_{\Omega}}$
\quad \ \
the set of last multipliers on the domain $\Omega$  of system (0.1);
\vspace{1ex}

\noindent
$\text{MA}_{_{\Xi}}$
\quad
the set of polynomial (real) last multipliers on the domain $\Xi$  of system (0.1);
\vspace{1ex}

\noindent
$\text{MB}_{_{\Xi}}$ 
\quad
the set of multiple polynomial (real) last multipliers on the domain $\Xi$ 
\\[0.5ex]
\mbox{}\hspace{3.7em}
of system (0.1);
\vspace{1ex}

\noindent
$\text{ME}_{_{\Omega}}$
\quad
the set of exponential last multipliers on the domain $\Omega$ of system (0.1);
\vspace{1ex}

\noindent
$\text{MF}_{_{\Xi}}$
\quad
the set of conditional last multipliers on the domain $\Xi$ of system (0.1);
\vspace{1ex}

\noindent
$\text{MH}_{_{\Xi}}$
\quad
the set of complex-valued polynomial last multipliers on the domain $\Xi$ 
\\[0.5ex]
\mbox{}\hspace{3.7em}
of system (0.1);
\vspace{1ex}

\noindent
$\text{MG}_{_{\Xi}}$
\quad
the set of multiple complex-valued polynomial last multipliers on the domain $\Xi$ 
\\[0.5ex]
\mbox{}\hspace{3.7em}
of system (0.1).

\newpage

\mbox{}
\\[-0.5ex]
\centerline{
{\bf\large \S\;\!1. Partial integrals}}
\\[2ex]
\centerline{
{\bf  1. 
Partial integral. Definition. Properties
}
}
\\[1.5ex]
\indent
{\bf Definition 1.1.}
\vspace{0.5ex}
{\it
We shall say that a continuously differentiable function
$
{\rm g}\colon \Omega\to\R
$
is a \textit{\textbf{partial integral on the domain}} 
$\Omega$ of system {\rm (0.1)} if  the differential of this function 
by virtue of system {\rm (0.1)} 
\\[2ex]
\mbox{}\hfill                                           % (1.1)
$
\displaystyle
d\!\;{\rm g}(t,x)_{\displaystyle |_{(0.1)}} =
{\rm g}(t,x)\;\!M(t,x)\;\!dt
$ 
\ for all 
$(t,x)\in \Omega,
$
\hfill {\rm(1.1)}
\\[2.25ex]
where the function 
\vspace{0.25ex}
$M\in \text{\rm P}_{_{\!\Xi^{\;\!\prime}}}$ and has the degree 
$\deg_{\;\!x}^{} M\leq d-1.$ 
Moreover, the function $M$ is said to be 
\textit{\textbf{cofactor of partial integral}} ${\rm g}.$
}
\vspace{0.5ex}

By $\text{J}_{_{\Omega}}$ denote the set of all partial integrals on the domain $\Omega$ 
of system  (0.1).
The phrase 
"the function ${\rm g}$ is a partial integral with cofactor $M$ on the domain 
$\Omega$ of system {\rm (0.1)}" is denoted by $({\rm g}, M)\in \text{J}_{_{\Omega}}.$
\vspace{0.25ex}
If a set $\Omega_{_0}\subset \R^{n+1}$ is not a domain, then  
by $({\rm g}, M)\in \text{J}_{_{\Omega_{{}_{\tiny 0}}}}$ denote 
the phrase "the function ${\rm g}$ is a partial integral with cofactor $M$
on any domain from the set $\Omega_{_0}$ of system {\rm (0.1)}". 

Using the operator of differentiation by virtue of system (0.1),
we can write the differential identity (1.1) as one of the operator identities
\\[2ex]
\mbox{}\hfill                             %(1.2)
$
\displaystyle
{\frak d}\;\! {\rm g}(t,x) =
{\rm g}(t,x)\;\!M(t,x)
$ 
\ for all 
$
(t,x)\in \Omega
$
\hfill(1.2)
\\[1ex]
or
\\[1ex]
\mbox{}\hfill                             %(1.3)
$
\displaystyle
{\frak d} \ln \bigl| {\rm g}(t,x) \bigr| = M(t,x)
$ 
\ for all 
$
(t,x)\in \Omega_{_0}\;\!,
$
\hfill(1.3)
\\[2.5ex]
where the set $\Omega_{_0}\subset \Omega$ such that  
\vspace{0.25ex}
${\rm g}(t,x)\ne 0$ for all $(t,x)\in \Omega_{_0}$ and
${\rm g}(t,x)=0$ for all $(t,x)\in {\sf C}_{_\Omega}\Omega_{_0}.$
\vspace{0.5ex}
Thus we have two existence criteria of partial integral.

{\bf Theorem 1.1.}
\vspace{0.35ex}
${\rm (g}, M)\in \text{J}_{_{\Omega}}$
{\it 
if and only if the identity {\rm (1.2)} holds
{\rm(}as well as if and only if the identity {\rm (1.3)} holds{\rm)}.
Moreover, both in the identity {\rm (1.2)} and in the identity {\rm (1.3)}, the function
$M\in \text{\rm P}_{_{\!\Xi^{\;\!\prime}}}$ and has the degree 
$\deg_{\;\!x}^{} M\leq d-1.$ 
}
\vspace{0.5ex}

Using the definition of integral manifold (Definition 0.3)
and the definition of partial integral (Definition 1.1), we get 
the geometric sense of partial integral.
\vspace{0.5ex}

{\bf Theorem 1.2.}
{\it 
If a partial integral ${\rm g}$ of system {\rm (0.1)} defines the manifold
${\rm g}(t,x)=0,$
then this manifold is an integral manifold of system {\rm (0.1)}.
}
\vspace{0.5ex}

{\bf Theorem 1.3.}
{\it
A continuously differentiable function
\\[1.5ex]
\mbox{}\hfill
$
\displaystyle
{\rm g}\colon (t,x)\to\ 
\prod\limits_{j=1}^{m} {\rm g}_j^{}(t,x)
$
\ for all $(t,x)\in\Omega
$
\hfill {\rm (1.4)}
\\[1.5ex]
is a partial integral with cofactor 
\vspace{0.25ex}
$M$ on the domain $\Omega$ of system {\rm (0.1)}
if and only if 
there exist functions 
$M_j^{}\in C^1\Omega,\ j=1,\ldots,m,$ 
such that 
\\[1.5ex]
\mbox{}\hfill
$
\displaystyle
{\frak d}\;\!{\rm g}_j^{}(t,x)  = 
{\rm g}_j^{}(t,x)\;\! M_j^{}(t,x)
$
\ for all $(t,x)\in\Omega_{_0}\;\!,
\quad
j=1,\ldots,m,
$
\hfill {\rm (1.5)}
\\[1ex]
and 
\\[1ex]
\mbox{}\hfill
$
\displaystyle
\sum\limits_{j=1}^{m} M_j^{}(t,x)= M(t,x)
$
\ for all 
$(t,x)\in\Omega_{_0}\;\!,
$
\hfill {\rm (1.6)}
\\[1.5ex]
where the function
$M\in \text{\rm P}_{_{\!\Xi^{\;\!\prime}}}$ and has the degree  
$\deg_{\;\!x}^{} M\leq d-1.$ 
}
\vspace{0.5ex}

{\sl Proof. Necessity}.
Case $m=2,\  ({\rm g}_1^{}{\rm g}_2^{}, M)\in  \text{\rm J}_{_{\Omega}}.$
The identity (1.2) from Theorem (1.2) has the form 
\\[1.5ex]
\mbox{}\hfill                          
$
\displaystyle
{\rm g}_2^{}(t,x)\;\! {\frak d}\;\! {\rm g}_1^{}(t,x) +
{\rm g}_1^{}(t,x)\;\! {\frak d}\;\! {\rm g}_2^{}(t,x) =
{\rm g}_1^{}(t,x)\;\! {\rm g}_2^{}(t,x)\;\!M(t,x)
$
\ for all 
$(t,x)\in\Omega,
\hfill
$
\\[2ex]
where the function
$M\in \text{\rm P}_{_{\!\Xi^{\;\!\prime}}}$ and  
$\deg_{\;\!x}^{} M\leq d-1.$
Hence the identity holds
\\[1.5ex]
\mbox{}\hfill                          
$
\displaystyle
{\frak d}\;\! {\rm g}_1^{}(t,x) =
{\rm g}_1^{}(t,x)\;\! 
\biggl(
M(t,x)-\dfrac{{\frak d}\;\! {\rm g}_2^{}(t,x)}{{\rm g}_2^{}(t,x)}
\biggr)
$
\ for all 
$(t,x)\in\Omega_{_0},
\hfill
$
\\[2ex]
where the set $\Omega_{_0}$ such that   
${\rm g}_2^{}(t,x)\ne 0$ for all $(t,x)\in \Omega_{_0}$ and 
${\rm g}_2^{}(t,x)=0$ for all $(t,x)\in {\sf C}_{_\Omega}\Omega_{_0}.$
\vspace{0.5ex}

Suppose the derivative by virtue of system (0.1)
\\[1.5ex]
\mbox{}\hfill                            
$
M(t,x)-\dfrac{{\frak d}\;\! {\rm g}_2^{}(t,x)}{{\rm g}_2^{}(t,x)}=
M_1^{}(t,x)
$ 
\ for all 
$(t,x)\in \Omega_{_0}.
\hfill 
$
\\[1.5ex]
\indent
Then, 
\\[1.5ex]
\mbox{}\hfill                          
$
\displaystyle
{\frak d}\;\! {\rm g}_1^{}(t,x)=
{\rm g}_1^{}(t,x)\;\! M_1^{}(t,x)
$ 
\ for all $(t,x)\in\Omega_{_0}\;\!,
\hfill
$
\\[2.5ex]
\mbox{}\hfill                          
$
\displaystyle
{\frak d}\;\! {\rm g}_2^{}(t,x)=
{\rm g}_2^{}(t,x)\;\! \bigl(M(t,x)-M_1^{}(t,x)\bigr)
$ 
\ for all 
$(t,x)\in\Omega_{_0}\;\!.
\hfill
$
\\[2ex]
\indent
So the identities (1.5) and (1.6) for $m=2$ are proved.
\vspace{0.35ex}

The identities (1.5) and (1.6) for $m>2$ are proved by induction.
\vspace{0.35ex}

{\sl Sufficiency.}
\vspace{0.25ex}
Suppose the identities (1.5) and (1.6) are true, and the function 
$M\in \text{\rm P}_{_{\!\Xi^{\;\!\prime}}}$ such that the degree  
$\deg_{\;\!x}^{} M\leq d-1.$ Then,
\\[2ex]
\mbox{}\hfill                        
$
\displaystyle
{\frak d}\;\!{\rm g}(t,x)  =
{\frak d}\;\!\prod\limits_{j=1}^{m} {\rm g}_j^{}(t,x) \, = \,
\sum\limits_{\nu=1}^{m}\, 
\prod\limits_{{\,}_{\scriptstyle j\ne \nu}^{\scriptstyle j=1,}}^{m}
{\rm g}_j^{}(t,x)\, {\frak d}\;\! {\rm g}_{\nu}^{}(t,x)  
=
\hfill
$
\\[2ex]
\mbox{}\hfill                        
$
\displaystyle
=\,
\prod\limits_{j=1}^{m}\, 
{\rm g}_j^{}(t,x)\, 
\sum\limits_{\nu=1}^{m} M_{\nu}^{}(t,x)=
{\rm g}(t,x)\;\! M(t,x)
$
\ for all $(t,x)\in \Omega.
\hfill
$
\\[1.5ex]
\indent
By Theorem 1.1, the function (1.4) is a partial integral with cofactor 
$M$ on the domain $\Omega$ of system (0.1). $\k$
\vspace{0.75ex}

{\bf  Property 1.1.}
\vspace{0.5ex}
{\it
Suppose a function $\varphi\in C^1T^{\;\!\prime},$ 
the set $T_{_0}\subset T^{\;\!\prime}$ such that
$\varphi(t)\ne 0$ for all $t\in T_{_0}$ and 
$\varphi(t)=0$ for all $t\in {\sf C}_{{}_{T^{\;\!\prime}}}T_{_0},$
the set $\Omega_{_0}=T_{_0}\times X^{\;\!\prime}.$ Then}
\\[2ex]
\mbox{}\hfill
$
({\rm g}, M)\in \text{J}_{_{\Omega}}
\iff 
\bigl(\varphi\;\! {\rm g}\;\!,\;\! M+{\sf D} \ln |\varphi|\bigr)\in 
\text{J}_{_{\Omega_{{}_{\tiny\;\! 0}}}}.
\hfill
$
\\[2ex]
\indent
{\sl Proof} is based on Theorem 1.1 and the identity  
\\[1.5ex]
\mbox{}\hfill                        
$
{\frak d}\bigl(\varphi(t)\;\!{\rm g}(t,x)\bigr)=
{\rm g}(t,x)\, {\sf D} \varphi(t) +\varphi(t)\, {\frak d}\;\!{\rm g}(t,x) 
=
\hfill
$
\\[2ex]
\mbox{}\hfill
$
=
\varphi(t)\;\! \bigl({\rm g}(t,x)\, {\sf D} \ln |\varphi(t) | + 
 {\frak d}\;\!{\rm g}(t,x)\bigr)
$
\ for all 
$(t,x)\in  \Omega_{_0}.\ \k
\hfill 
$
\\[2ex]
\indent
As a consequence of Property 1.1, we obtain the following
\vspace{1ex}

{\bf Property 1.2.} 
{\it
If $\lambda\in\R\backslash\{0\},$ then}
\\[1.5ex]
\mbox{}\hfill
$
({\rm g}, M)\in \text{J}_{_{\Omega}}
\iff 
(\lambda\;\! {\rm g}\;\!,\;\! M)\in \text{J}_{_{\Omega}}.
\hfill
$
\\[2ex]
\indent
By Property 1.2, 
if we have two or more partial integrals of system (0.1), then 
we assume that they are pairwise linearly independent.
\vspace{1ex}

{\bf Property 1.3.} 
\vspace{0.5ex}
{\it
Suppose the set $\Omega_{_0}\subset\Omega$ such that   
${\rm g}(t,x)\ne 0$ for all $(t,x)\in \Omega_{_0}$ and 
${\rm g}(t,x)=0$ for all $(t,x)\in {\sf C}_{_\Omega}\Omega_{_0}.$
Then}
\\[1.5ex]
\mbox{}\hfill  % (1.7)
$
({\rm g}, M)\in \text{J}_{_{\Omega}}
\iff 
(|{\rm g}|\;\!,\;\! M)\in \text{J}_{_{\Omega_{{}_{\tiny\;\! 0}}}}.
$
\hfill (1.7)
\\[2ex]
\indent
{\sl Proof.}
Since
\\[1ex]
\mbox{}\hfill 
$
({\rm g}, M)\in \text{J}_{_{\Omega}}
\iff 
({\rm g}, M)\in \text{J}_{_{\Omega_{{}_{\tiny\;\! 0}}}},
\hfill
$
\\[2ex]
\mbox{}\hfill 
$
|{\rm g}(t,x)|=\text{sgn}\;\!{\rm g}(t,x)\, {\rm g}(t,x)
$
\ for all 
$(t,x)\in\Omega,
\hfill
$
\\[1.5ex]
by Property 1.2, we  have
\\[1.5ex]
\mbox{}\hfill 
$
({\rm g}, M)\in \text{J}_{_{\Omega_{{}_{\tiny\;\! 0}}}}
\iff 
({\rm sgn\;\!g\,g}, M)\in \text{J}_{_{\Omega_{{}_{\tiny\;\! 0}}}}.
\hfill
$
\\[1.5ex]
Now, using the transitivity of equivalence, we obtain the statement  (1.7). $\k$
\vspace{0.75ex}

By Property 1.3, the connection between the identities (1.2) and (1.3) is established. 
\vspace{1ex}

{\bf Property 1.4.} 
{\it
If the derivative by virtue of system {\rm (0.1)}
\\[1.5ex]
\mbox{}\hfill
$
{\frak d}\;\!{\rm g}(t,x)=
\bigl({\rm g}(t,x)+c\bigr)\;\!M(t,x)
$
\ for all 
$(t,x)\in\Omega,
\hfill
$
\\[2ex]
where $c\in\R,\ M\in \text{\rm P}_{_{\!\Xi^{\;\!\prime}}},\  
\deg_{\;\!x}^{} M\leq d-1,$ then 
$({\rm g}+c,M)\in \text{\rm J}_{_{\Omega}}.$}
\vspace{0.75ex}

{\sl Indeed,} the derivative by virtue of system (0.1)
\\[1.5ex]
\mbox{}\hfill
$
{\frak d}\;\!\bigl({\rm g}(t,x)+c\bigr)=
{\frak d}\;\!{\rm g}(t,x)=
\bigl({\rm g}(t,x)+c\bigr)\;\!M(t,x)
$
\ for all 
$(t,x)\in\Omega.
\hfill
$
\\[1.5ex]
\indent
Therefore, by Theorem 1.1, we have $({\rm g}+c,M)\in \text{\rm J}_{_{\Omega}}.\ \k$
\vspace{1ex}

{\bf Property 1.5.} 
\vspace{0.5ex}
{\it
If $({\rm g}_j^{}, M)\in \text{\rm J}_{_{\Omega}},\ 
\lambda_j^{}\in\R\backslash\{0\},\ j=1,\ldots,m,$ then 
$\biggl(\,\sum\limits_{j=1}^{m}\! \lambda_j^{}\;\!{\rm g}_j^{}, M\biggl) 
\in \text{\rm J}_{_{\Omega}}.$}

{\sl Proof.} If 
$({\rm g}_j^{}, M)\in \text{\rm J}_{_{\Omega}}\;\!,\ j=1,\ldots,m,$
then, by Theorem 1.1,
\\[1.75ex]
\mbox{}\hfill
$
{\frak d}\;\!{\rm g}_j^{}(t,x)=
{\rm g}_j^{}(t,x)\;\!M(t,x)
$
\ for all 
$(t,x)\in\Omega,
\quad
j=1,\ldots,m,
\hfill
$
\\[1.75ex]
where the cofactor $M\in \text{\rm P}_{_{\!\Xi^{\;\!\prime}}}$ and has the degree 
$\deg_{\;\!x}^{} M\leq d-1.$
\vspace{0.75ex}

Consequently the derivative by virtue of system (0.1)
\\[1.5ex]
\mbox{}\hfill                         
$
\displaystyle
{\frak d}\;\!\sum\limits_{j=1}^{m} \lambda_j^{}\;\!{\rm g}_j^{}(t,x)  = 
\sum\limits_{j=1}^{m}
\lambda_j^{}\;\!{\rm g}_j^{}(t,x)\;\!M(t,x)
$
\ for all 
$(t,x)\in\Omega.
\hfill
$
\\[1.5ex]
\indent
Thus, by Theorem 1.1, 
$\biggl(\,\sum\limits_{j=1}^{m}\! \lambda_j^{}\;\!{\rm g}_j^{}, M\biggl) 
\in \text{\rm J}_{_{\Omega}}.\ \k$
\vspace{1.25ex}

{\bf Property 1.6.} 
{\it
Suppose $\gamma\in\R\backslash\{0\},\ {\rm g}^{\gamma}\in C^1\Omega.$
Then}
\\[2ex]
\mbox{}\hfill    %(1.8)
$
({\rm g}, M)\in \text{J}_{_{\Omega}}
\iff 
\bigl({\rm g}^{\gamma},\;\! \gamma\;\!M\bigr)\in \text{J}_{_{\Omega}}.
$
\hfill (1.8)
\\[2ex]
\indent
{\sl Proof}
is based on Theorem 1.1, which applied to the functions ${\rm g}$ and ${\rm g}^{\gamma},$ 
and on the identity 
\\[0.15ex]
\mbox{}\hfill
$
{\frak d}\;\!{\rm g}^{\gamma}(t,x)=
\gamma\;\!{\rm g}^{\gamma-1}(t,x)\,
{\frak d}\;\!{\rm g}(t,x)
$
\ for all 
$(t,x)\in\Omega. \ \k
\hfill
$
\\[2.25ex]
\indent
{\bf Property 1.7.} 
\vspace{0.5ex}
{\it
Suppose $\gamma\in\R\backslash\{0\},$ 
the set $\Omega_{_0}\subset\Omega$ such that   
${\rm g}(t,x)\ne 0$ for all $(t,x)\in \Omega_{_0}$ 
and 
${\rm g}(t,x)=0$ for all $(t,x)\in {\sf C}_{_\Omega}\Omega_{_0}.$
Then}
\\[2ex]
\mbox{}\hfill  % (1.9)
$
({\rm g}, M)\in \text{J}_{_{\Omega}}
\iff 
\bigl(|{\rm g}|^{\gamma},\;\! \gamma\;\!M\bigr)\in 
\text{J}_{_{\Omega_{{}_{\tiny\;\! 0}}}}.
$
\hfill (1.9)
\\[2ex]
\indent
{\sl Proof.}
By Property 1.3, equivalence (1.7) is true.
By Property 1.6,
\\[1.5ex]
\mbox{}\hfill 
$
(|{\rm g}|, M)\in \text{J}_{_{\Omega_{{}_{\tiny\;\! 0}}}}
\iff 
\bigl(|{\rm g}|^{\gamma},\;\! \gamma\;\!M\bigr)\in 
\text{J}_{_{\Omega_{{}_{\tiny\;\! 0}}}}.
\hfill
$
\\[1.5ex]
\indent
Using the transitivity of equivalence, we get the statement  (1.9). $\k$
\vspace{1ex}

{\bf Property 1.8.} 
{\it
If $\rho_j^{}, \lambda_j^{}\in\R\backslash\{0\},
\ ({\rm g}_j^{}, \rho_j^{}M_{0}^{})\in \text{\rm J}_{_{\Omega}},\ 
{\rm g}_j^{{}^{\scriptsize 1/\rho_{\!j}^{}}}\in C^1\Omega, ,\ j=1,\ldots,m,$ then}
\\[1ex]
\mbox{}\hfill
$
\displaystyle
\biggl(\,\sum\limits_{j=1}^{m}\! 
\lambda_j^{}\;\!{\rm g}_j^{{}^{\scriptsize 1/\rho_{\!j}^{}}}, M_0^{}\biggl) 
\in \text{\rm J}_{_{\Omega}}.
\hfill
$
\\[0.75ex]
\indent
{\sl Proof.} 
Since  
$\!\!({\rm g}_j^{}, \rho_j^{}M_0^{})\!\in\! \text{\rm J}_{_{\Omega}},\!$
by Property 1.6, we see that  
$\!\!\Bigl({\rm g}_j^{{}^{\scriptsize 1/\rho_{\!j}^{}}}, M_0^{}\Bigr)\!\in\! 
\text{\rm J}_{_{\Omega}}, 
\, 
j\!=\!1,\ldots,m.$
Then, from Property 1.5 it follows that
$\biggl(\,\sum\limits_{j=1}^{m}\! 
\lambda_j^{}\;\!{\rm g}_j^{{}^{\scriptsize 1/\rho_{\!j}^{}}}, M_0^{}\biggl) 
\in \text{\rm J}_{_{\Omega}}.\ \k$
\vspace{0.75ex}

{\bf Property 1.9.} 
\vspace{0.5ex}
{\it
Let 
$({\rm g}_j^{}, M_{j}^{})\in \text{\rm J}_{_{\Omega}},\ 
\gamma_j^{}\in\R\backslash\{0\},\
{\rm g}_j^{{}^{\scriptsize \gamma_{j}^{}}}\in C^1\Omega, \ j=1,\ldots,m.$ 
Then   
$\biggl(\,\prod\limits_{j=1}^{m} 
{\rm g}_j^{{}^{\scriptsize \gamma_{j}^{}}}, M\biggl) 
\in \text{\rm J}_{_{\Omega}}$
if and only if the cofactors $M,\;\!M_1^{},\ldots,M_m^{}$ such that 
the identity holds}
\\[1ex]
\mbox{}\hfill                      % (1.10)
$
\displaystyle
M(t,x)=\sum\limits_{j=1}^{m}
\gamma_j^{}\;\!M_{j}^{}(t,x)
$
\ for all 
$(t,x)\in\Xi^{\;\!\prime}.
$
\hfill (1.10)
\\[1.5ex]
\indent
{\sl Proof.} 
\vspace{0.35ex}
Since
$({\rm g}_j^{}, M_{j}^{})\in \text{\rm J}_{_{\Omega}},\ j=1,\ldots,m,$
we see that from Property 1.6 and Theorem 1.1 it follows that  
\\[1.75ex]
\mbox{}\hfill   % (1.11)
$
{\frak d}\;\!{\rm g}_j^{{}^{\scriptsize \gamma_{j}^{}}}\!(t,x)=
\gamma_j^{}\;\!{\rm g}_j^{{}^{\scriptsize \gamma_{j}^{}}}\!(t,x)\;\!M_j^{}(t,x)
$
\ for all 
$
(t,x)\in\Omega,
\quad
j=1,\ldots,m,
$
\hfill (1.11)
\\[2.5ex]
where the functions $M_j^{}\in \text{\rm P}_{_{\!\Xi^{\;\!\prime}}}$ and have the degrees
$\deg_{\;\!x}^{} M_j^{}\leq d-1,\ j=1,\ldots, m.$
\vspace{1ex}

By Theorem 1.3, 
\vspace{0.75ex}
$\biggl(\,\prod\limits_{j=1}^{m} 
{\rm g}_j^{{}^{\scriptsize \gamma_{j}^{}}}, M\biggl) 
\in \text{\rm J}_{_{\Omega}}$
if and only if 
there exist functions
$M_j^{\ast}\in C^1\Omega,
\linebreak 
j=1,\ldots, m,$
such that the identities hold
\\[2.25ex]
\mbox{}\hfill   % (1.12)
$
{\frak d}\;\!{\rm g}_j^{{}^{\scriptsize \gamma_{j}^{}}}\!(t,x)=
{\rm g}_j^{{}^{\scriptsize \gamma_{j}^{}}}\!(t,x)\;\!M_j^{\ast}(t,x)
$
\ for all 
$(t,x)\in\Omega,
\quad
j=1,\ldots,m,
$
\hfill (1.12)
\\[1.25ex]
and
\\[1ex]
\mbox{}\hfill                      
$
\displaystyle
M(t,x)=\sum\limits_{j=1}^{m}
M_{j}^{\ast}(t,x)
$
\ for all 
$(t,x)\in\Omega,
\hfill
$
\\[1.5ex]
where the function $M\in \text{\rm P}_{_{\!\Xi^{\;\!\prime}}}$ and has the degree  
$\deg_{\;\!x}^{} M\leq d-1.$
\vspace{1ex}

Using the identities (1.12) and (1.11), we obtain
\\[1.5ex]
\mbox{}\hfill                      
$
\displaystyle
M_{j}^{\ast}(t,x)=
\gamma_j^{}\;\!M_j^{}(t,x)
$
\ for all 
$
(t,x)\in\Xi^{\;\!\prime},
\quad
j=1,\ldots,m,
\hfill
$
\\[2.25ex]
where the functions $M_j^{}\in \text{\rm P}_{_{\!\Xi^{\;\!\prime}}}$ and have the degrees 
$\deg_{\;\!x}^{} M_j^{}\leq d-1,\ j=1,\ldots, m.$
\vspace{1ex}

Thus,
\vspace{1ex}
$\biggl(\,\prod\limits_{j=1}^{m} 
{\rm g}_j^{{}^{\scriptsize \gamma_{j}^{}}}, M\biggl) 
\in \text{\rm J}_{_{\Omega}}$
if and only if the identity (1.10) is true. $\k$

{\bf Corollary 1.1.} 
\vspace{0.5ex}
{\it
Let 
$\rho_j^{},\gamma_j^{}\in\R\backslash\{0\},\
({\rm g}_j^{}, \rho_j^{}M_{0}^{})\in \text{\rm J}_{_{\Omega}},\ 
{\rm g}_j^{{}^{\scriptsize \gamma_{j}^{}}}\in C^1\Omega, \ j=1,\ldots,m.$ 
Then  
$\biggl(\,\prod\limits_{j=1}^{m} 
{\rm g}_j^{{}^{\scriptsize \gamma_{j}^{}}}, M\biggl) 
\in \text{\rm J}_{_{\Omega}}$
if and only if the cofactor}
\\[1ex]
\mbox{}\hfill                      % (1.13)
$
\displaystyle
M(t,x)=\sum\limits_{j=1}^{m}
\rho_j^{}\;\!\gamma_j^{}\;\!M_{0}^{}(t,x)
$
\ for all 
$
(t,x)\in\Xi^{\;\!\prime}.
$
\hfill (1.13)
\\[1.5ex]
\indent
Using Theorem 1.3, from Properties 1.6 and 1.9, we can state the following
\vspace{1.25ex}

{\bf Property\! 1.10.}\! 
\vspace{0.5ex}
{\it
Let 
$\!({\rm g}_\tau^{}, M_{\tau}^{})\!\in\! \text{\rm J}_{_{\Omega}},\, 
\tau\!=\!1,\ldots, m-1,\,
\gamma_j^{}\!\in\!\R\backslash\{0\},\,
{\rm g}_j^{{}^{\scriptsize \gamma_{j}^{}}}\!\in\! C^1\Omega, \,
%\linebreak 
j\!=\!1,\ldots,m.$ 
Then we claim that}
\\[1.25ex]
\mbox{}\hfill                      
$
\displaystyle
\biggl(\,\prod\limits_{j=1}^{m} 
{\rm g}_j^{{}^{\scriptsize \gamma_{j}^{}}}, M\biggl) 
\in \text{\rm J}_{_{\Omega}}
\iff
\biggl({\rm g}_m^{}, \,
\dfrac{1}{\gamma_m^{}}\;\!
\biggl(M-\sum\limits_{\tau=1}^{m-1}
\gamma_{\tau}\;\!M_{\tau}^{}\biggl)\biggl) 
\in \text{\rm J}_{_{\Omega}}\;\!.
\hfill
$
\\[1.5ex]
\indent
From Theorem 1.3 and Property 1.10, we get the following
\vspace{1.25ex}

{\bf Property 1.11.} 
\vspace{0.5ex}
{\it
Let 
$({\rm g}_\nu^{}, M_{\nu}^{})\in \text{\rm J}_{_{\Omega}},\ 
\nu=1,\ldots, s,\ s\leq m-2, \ \
\gamma_j^{}\in\R\backslash\{0\},\
{\rm g}_j^{{}^{\scriptsize \gamma_{j}^{}}}\!\in C^1\Omega, 
\linebreak 
j=1,\ldots,m.$ 
Then we have
\\[1ex]
\mbox{}\hfill                      
$
\displaystyle
\biggl(\,\prod\limits_{j=1}^{m} 
{\rm g}_j^{{}^{\scriptsize \gamma_{j}^{}}}, M\biggl) 
\in \text{\rm J}_{_{\Omega}}
\iff
\biggl(
\,\prod\limits_{k=s+1}^{m} 
{\rm g}_k^{{}^{\scriptsize \gamma_{k}^{}}}, \,
M-\sum\limits_{\nu=1}^{s}
\gamma_{\nu}\;\!M_{\nu}^{}\biggl) 
\in \text{\rm J}_{_{\Omega}}\;\!.
\hfill
$
\\[1.5ex]
Moreover, there exist functions
$M_k^{}\in C^1\Omega,\ k=s+1,\ldots, m,$
such that the identities hold
\\[1.5ex]
\mbox{}\hfill                        % (1.14)
$
{\frak d}\;\!{\rm g}_k^{}(t,x)=
{\rm g}_k^{}(t,x)\;\!M_k^{}(t,x)
$
\ for all 
$(t,x)\in\Omega,
\quad
k=s+1,\ldots,m,
$
\hfill {\rm (1.14)}
\\[1ex]
and}
\\[1ex]
\mbox{}\hfill                      % (1.15)
$
\displaystyle
\sum\limits_{k=s+1}^{m}
\gamma_{k}^{}\;\!M_k^{}(t,x)=
M(t,x)-\sum\limits_{\nu=1}^{s}
\gamma_{\nu}^{}\;\!M_\nu^{}(t,x)
$
\ for all 
$(t,x)\in\Omega.
$
\hfill (1.15)
\\[2ex]
\indent
{\bf Remark 1.1.}\!
\vspace{0.35ex}
Properties\! 1.8 -- 1.11 and Corollary 1.1 were proved on the basis of equiva\-lence (1.8) and Property 1.6.
In case of need we can replace the degree 
$\!{\rm g}_j^{{}^{\scriptsize \gamma_{j}^{}}},\ j\!\in\!\{1,\ldots,m\}$
by the degree 
$\!|{\rm g}_j^{}|^{{}^{\scriptsize \gamma_{j}^{}}}\!\!$
in Properties\! 1.8--1.11 and Corollary\! 1.1\! (with correction on Property\! 1.7).
\\[5.25ex]
\centerline{
{\bf  2. Polynomial partial integrals}
}
\\[1.75ex]
\indent
{\bf Definition 2.1.}
\vspace{0.15ex}
{\it
A partial integral on the domain $\Xi^{\;\!\prime}$ of system {\rm (0.1)} is called 
a \textit{\textbf{polynomial partial integral on the domain}}
$\Xi^{\;\!\prime}$ of system {\rm (0.1)} if
this partial integral is a polynomial in the variables $x_1^{},\ldots, x_n^{}$ with 
\vspace{0.25ex}
continuously differentiable on the domain $T^{\;\!\prime}$
coefficients-functions of the variable $t.$
}
\vspace{0.75ex}

By $\text{A}_{_{\Xi^{\;\!\prime}}}$ denote the set of polynomial partial integrals on the domain 
\vspace{0.75ex}
$\Xi^{\;\!\prime}$ of system  (0.1).

From Definition 2.1 it follows that the set 
$\text{A}_{_{\Xi^{\;\!\prime}}}\subset \text{J}_{_{\Xi^{\;\!\prime}}}.$
\vspace{0.75ex}

The phrase 
\vspace{0.5ex}
"the function $p\in \text{P}_{_{\!\Xi^{\;\!\prime}}}$ is 
a polynomial partial integral with cofactor $M$ on the domain 
$\Xi^{\;\!\prime}$ of system {\rm (0.1)}" is denoted by $(p, M)\in \text{A}_{_{\Xi^{\;\!\prime}}}.$
\vspace{0.75ex}

Using Definition 2.1, we obtain
\\[1.5ex]
\mbox{}\hfill               % (2.1)
$
(p,M)\in \text{A}_{_{\Xi^{\;\!\prime}}}
\iff
(p,M)\in \text{J}_{_{\Xi^{\;\!\prime}}}
\ \&\ 
p\in \text{P}_{_{\!\Xi^{\;\!\prime}}}.
$
\hfill (2.1)
\\[2ex]
\indent
{\bf Theorem 2.1}
\vspace{0.35ex}
(existence criterion of polynomial partial integral).
{\it
$(p, M)\in \text{\rm A}_{_{\Xi^{\;\!\prime}}}$ 
if and only if the derivative by virtue of system} (0.1)
\\[1.5ex]
\mbox{}\hfill               % (2.2)
$
{\frak d}\;\!p(t,x)=
p(t,x)\;\!M(t,x)
$
\ {\it for all} 
$(t,x)\in \Xi^{\;\!\prime},
\quad
p\in \text{P}_{_{\!\Xi^{\;\!\prime}}}.
$
\hfill (2.2)
\\[1.5ex]
\indent
{\sl Proof.} 
\vspace{0.25ex}
Taking into account the equivalence (2.1), from 
the existence criterion of polynomial partial integral (Theorem 1.1), we get the 
statement of Theorem 2.1. $\k$
\vspace{1ex}

{\bf Remark 2.1.}
\vspace{0.75ex}
If the function $p\in \text{P}_{_{\!\Xi^{\;\!\prime}}},$
then from the identity (2.2) it follows that the function 
$M\in \text{P}_{_{\!\Xi^{\;\!\prime}}}$ and the degree 
$\deg_{\;\!x}^{} M\leq d-1.$
\vspace{1ex}

{\bf Theorem 2.2.}\! 
\vspace{0.5ex}
{\it
Suppose 
$p_j^{}\!\in\! \text{\rm P}_{_{\!\Xi^{\;\!\prime}}},\, 
\gamma_j^{}\!\in\!\R\backslash\{0\},\!$ 
a set $\Omega_{_0}\!\subset\!\Xi^{\;\!\prime}\!$ such that  
$p_j^{{}^{\scriptsize \gamma_j^{}}}\!\in\! C^1\Omega_{_0},$ $j=1,\ldots, m.$
Then
\\[1.5ex]
\mbox{}\hfill  
$
\displaystyle
\biggl(\,\prod\limits_{j=1}^{m} 
p_j^{{}^{\scriptsize \gamma_{j}^{}}}, M\biggl) 
\in \text{\rm J}_{_{\Omega_{{}_{\tiny\;\! 0}}}}
\iff
\bigl(p_j^{}, M_j^{}\bigr)\in \text{\rm A}_{_{\Xi^{\;\!\prime}}},
\ \ 
j=1,\ldots, m,
\hfill
$
\\[1.5ex]
where the cofactors $M$ and $M_j^{},\ j=1,\ldots,m,$
\vspace{0.75ex}
such that the identity {\rm (1.10)} holds.}

{\sl Proof.}
\vspace{0.5ex}
By Theorem 1.3,
$
\biggl(\,\prod\limits_{j=1}^{m} 
p_j^{{}^{\scriptsize \gamma_{j}^{}}}, M\biggl) 
\in \text{J}_{_{\Omega_{{}_{\tiny\;\! 0}}}}
$
if and only if  there exist functions 
\linebreak
 $M_j^{\ast}\in C^1\Omega_{_0},\ j=1,\ldots,m,$
such that the identities hold
\\[2ex]
\mbox{}\hfill   % (2.3)
$
{\frak d}\;\!p_j^{{}^{\scriptsize \gamma_{j}^{}}}\!(t,x)=
p_j^{{}^{\scriptsize \gamma_{j}^{}}}\!(t,x)\;\!M_j^{\ast}(t,x)
$
\ for all 
$(t,x)\in\Omega_{_{0}},
\quad
j=1,\ldots,m,
$
\hfill (2.3)
\\[1.5ex]
and
\\[1ex]
\mbox{}\hfill     % (2.4)                 
$
\displaystyle
\sum\limits_{j=1}^{m}
M_{j}^{\ast}(t,x)= M(t,x)
$
\ for all 
$(t,x)\in \Omega_{_{0}},
$
\hfill (2.4)
\\[2ex]
where the function $M\in \text{\rm P}_{_{\!\Xi^{\;\!\prime}}}$ and has the degree 
$\deg_{\;\!x}^{} M\leq d-1.$
%\vspace{0.5ex}

\newpage

Since the derivative by virtue of system (0.1)
\\[1.75ex]
\mbox{}\hfill 
$
{\frak d}\;\!p_j^{{}^{\scriptsize \gamma_{j}^{}}}\!(t,x)=
\gamma_j^{}\;\!p_j^{{}^{\scriptsize \gamma_{j}^{}-1}}\!(t,x)\;\!
{\frak d}\;\!p_j^{{}}\!(t,x)
$
\ for all 
$(t,x)\in\Omega_{_{0}},
\quad
j=1,\ldots,m,
\hfill
$
\\[2.5ex]
we see that the identities (2.3) are true if and only if 
\\[2ex]
\mbox{}\hfill 
$
{\frak d}\;\!p_j^{}(t,x)=
p_j^{}(t,x)\;\!M_j^{}(t,x)
$
\ for all 
$
(t,x)\in \Xi^{\;\!\prime},
\quad
j=1,\ldots,m,
\hfill
$
\\[1.5ex]
where
\\[1ex]
\mbox{}\hfill 
$
M_j^{}(t,x)=
\dfrac{1}{\gamma_j^{}}\, M_j^{\ast}(t,x)
$
\ for all 
$(t,x)\in \Xi^{\;\!\prime},
\quad
j=1,\ldots,m.
\hfill
$
\\[1.5ex]
\indent
Using Theorem 2.1, we have the identities (2.3) are true if and only if 
\vspace{0.75ex}
$
\bigl(p_j^{}, M_j^{}\bigr)\in \text{\rm A}_{_{\Xi^{\;\!\prime}}},
\linebreak 
j=1,\ldots, m.
$
From the identity (2.4) it follows that the identity (1.10) is true.$\k$
\vspace{1ex}

In particular, from Theorem 2.2, we get the following
\vspace{1ex}

{\bf Theorem 2.3}
\vspace{0.35ex}
(existence criterion of rational partial integral). {\it 
Suppose functions
\linebreak
$p_1^{},p_2^{}\in \text{\rm P}_{_{\!\Xi^{\;\!\prime}}}$
are relatively prime,
\vspace{0.75ex}
the set $\Omega_{_0}\subset\Xi^{\;\!\prime}$ such that   
$p_2^{}(t,x)\ne 0$ for all $(t,x)\in \Omega_{_0}$ and
$p_2^{}(t,x)= 0$ for all $(t,x)\in {\sf C}_{_{\Xi^{\;\!\prime}}}\Omega_{_0}.$
Then
\\[2ex]
\mbox{}\hfill  
$
\Bigl(\,\dfrac{p_1^{}}{p_2^{}}\,,\, M\Bigl) 
\in \text{\rm J}_{_{\Omega_{{}_{\tiny\;\! 0}}}}
\iff
\bigl(p_1^{}, M_1^{}\bigr)\in \text{\rm A}_{_{\Xi^{\;\!\prime}}}
\ \& \ 
\bigl(p_2^{}, M_2^{}\bigr)\in \text{\rm A}_{_{\Xi^{\;\!\prime}}}.
\hfill
$
\\[2ex]
Moreover, the cofactors $M,\ M_1^{},\ M_2^{}$ such that}
\\[1.5ex]
\mbox{}\hfill  
$
M(t,x)=M_1^{}(t,x)-M_2^{}(t,x)
$
\ {\it for all}
$(t,x)\in \Xi^{\;\!\prime}.
\hfill
$
\\[2.5ex]
\indent
{\bf  Property 2.1.}
\vspace{0.75ex}
{\it
Suppose $\varphi\in C^1T^{\;\!\prime},$ 
the set $T_{_0}\subset T^{\;\!\prime}$ 
such that  
$\varphi(t)\ne 0$ for all $t\in T_{_0}$ and 
$\varphi(t)=0$ for all $t\in {\sf C}_{{}_{T^{\;\!\prime}}}T_{_0},$
the set $\Xi_{_0}=T_{_0}\times \R^n.$ Then}
$
\bigl(\varphi, {\sf D} \ln |\varphi|\bigr)\in 
\text{\rm A}_{_{\Xi_{{}_{\tiny\;\! 0}}}}.
$
\vspace{0.75ex}

{\sl Indeed}, the derivative by virtue of system (0.1)  
\\[1.75ex]
\mbox{}\hfill                        
$
{\frak d}\;\!\varphi(t) =
{\sf D} \varphi(t) =
\varphi(t)\,\dfrac{{\sf D} \varphi(t)}{\varphi(t)}=
\varphi(t)\;\! {\sf D} \ln |\varphi(t) | 
$
\ for all 
$(t,x)\in  \Xi_{_0},
\hfill 
$
\\[1.75ex]
where $\deg_{\;\!x}^{} {\sf D} \ln |\varphi |=0\leq d-1.$
\vspace{1ex}

By Theorem 2.1, 
$
\bigl(\varphi, {\sf D} \ln |\varphi|\bigr)\in 
\text{\rm A}_{_{\Xi_{{}_{\tiny\;\! 0}}}}\!.\,  \k
$
\vspace{0.5ex}

{\bf  Property 2.2.}
%\vspace{0.5ex}
{\it
Suppose $\lambda_j^{},\;\! c_j^{}\in\R, \ j=1,\ldots,m,\ m\leq n,$ 
and $\sum\limits_{j=1}^{m}|\lambda_j^{}|\ne 0.$
Then  
\\
\mbox{}\hfill
$
\displaystyle
\biggl(\,\sum\limits_{j=1}^{m} 
\lambda_j^{}(x_j^{}+c_j^{}), M\biggl) 
\in\text{\rm A}_{_{\Xi}}
\hfill
$
\\[1ex]
if and only if the identity holds}
\\[1.5ex]
\mbox{}\hfill
$
\displaystyle
\sum\limits_{j=1}^{m}
\lambda_j^{}\;\!X_j^{}(t,x)=
M(t,x)\;\!\sum\limits_{j=1}^{m}\lambda_j^{}(x_j^{}+c_j^{})
$
\ {\it for all}
$(t,x)\in \Xi.
\hfill
$
\\[1.5ex]
\indent
{\sl Proof}\,
is based on Theorem 2.1 and used the derivative by virtue of system (0.1)
\\[1.5ex]
\mbox{}\hfill                        
$
\displaystyle
{\frak d}\;\!\sum\limits_{j=1}^{m}\lambda_j^{}(x_j^{}+c_j^{}) =
\sum\limits_{j=1}^{m}\lambda_j^{}\;\!X_j^{}(t,x)
$
\ for all 
$(t,x)\in  \Xi.
\ \k
\hfill 
$
\\[1.75ex]
\indent
From Property 2.2 under the condition $m=1,$ we have
\vspace{1ex}

{\bf  Property 2.3.}
\vspace{0.5ex}
{\it
If $c\in\R,$ then
$
\bigl(x_i^{}+c, M\bigr) \in \text{\rm A}_{_{\Xi}},\ 
i\in \{1,\ldots, n\},
$
if and only if the function $X_i^{}$ from right side of system {\rm (0.1)} has the form}
\\[1.5ex]
\mbox{}\hfill
$
X_i^{}(t,x)=
(x_i^{}+c)\;\! M(t,x)
$
\ for all 
$(t,x)\in \Xi.
\hfill
$
\\[-2.75ex]

\newpage

\mbox{}
\\[-1.75ex]
\centerline{
{\bf  3. Exponential partial integrals}
}
\\[1.5ex]
\indent
{\bf Definition 3.1.}
\vspace{0.35ex}
{\it
A function $\exp\omega$ is said to be an
\textit{\textbf{exponential partial integral}}
with cofactor $M$ on the domain $\Omega$ of system {\rm (0.1)}
if $\bigl(\exp\omega,M\bigr)\in \text{\rm J}_{_{\Omega}}.$
}
\vspace{0.75ex}

By $\text{E}_{_{\Omega}}$ denote the set of exponential partial integrals on the domain $\Omega$ 
\vspace{0.75ex}
of system (0.1).

From Definition 3.1 it follows that the set
$\text{E}_{_{\Omega}}\subset \text{J}_{_{\Omega}}.$
\vspace{0.5ex}

The phrase "the function $\exp\omega$ 
\vspace{0.35ex}
is an exponential partial integral with cofactor $M$ on
the domain $\Omega$ of system (0.1)" is denoted by 
$\bigl(\exp\omega,M\bigr)\in \text{\rm E}_{_{\Omega}}.$
\vspace{0.75ex}

{\bf Theorem 3.1}
\vspace{0.35ex}
(existence criterion of exponential partial integral).
{\it
$\bigl(\exp\omega,M\bigr)\in \text{\rm E}_{_{\Omega}}$ 
if and only if the identity holds
\\[1.5ex]
\mbox{}\hfill               % (3.1)
$
{\frak d}\;\!\omega(t,x)=M(t,x)
$
\ for all 
$(t,x)\in \Omega,
$
\hfill {\rm(3.1)}
\\[1.5ex]
where the function $M\in \text{\rm P}_{_{\!\Xi^{\;\!\prime}}}$ and has the degree
$\deg_{\;\!x}^{} M\leq d-1.$
}
\vspace{0.5ex}

{\sl Proof.}
Using 
the existence criterion of partial integral (Theorem~1.1) and 
the definition of exponential partial integral (Definition 3.1), we see that
\vspace{0.25ex}
the identity (1.3) (or the identity (1.2)) 
for the function ${\rm g}\colon (t,x)\to\exp\omega$ for all $(t,x)\in\Omega$
has the form (3.1). $\k$
\vspace{0.75ex}

{\bf Theorem 3.2}
\vspace{0.35ex}
(existence criterion of exponential partial integral).
{\it
A function $\exp\omega\in \text{\rm E}_{_{\Omega}}$ 
if and only if the function 
${\frak d}\;\!\omega\in \text{\rm P}_{_{\!\Xi^{\;\!\prime}}}$ 
and has the degree $\deg_{\;\!x}^{} {\frak d}\;\!\omega\leq d-1.$
\vspace{0.35ex}
Moreover, the function ${\frak d}\;\!\omega$ is the cofactor of 
the exponential partial integral $\exp\omega.$
}
\vspace{1ex}

{\bf  Theorem 3.3.}
{\it
Let $\gamma_j^{}\in\R\backslash\{0\},\ c_j^{}\in\R, \ j=1,\ldots,m.$
Then we claim that
\\[1.5ex]
\mbox{}\hfill
$
\displaystyle
\exp\sum\limits_{j=1}^{m}
\gamma_j^{}\;\!(\omega_j^{}+c_j^{})\in \text{\rm E}_{_{\Omega}}
\ \iff \
\sum\limits_{j=1}^{m}\gamma_j^{}\;\!{\frak d}\;\!\omega_j^{}\in 
\text{\rm P}_{_{\!\Xi^{\;\!\prime}}}
\, \ \&\ \,
\deg_{\;\!x}^{} 
\sum\limits_{j=1}^{m}
\gamma_j^{}\;\!{\frak d}\;\!\omega_j^{}\leq d-1.
\hfill
$
\\[1.5ex]
Furthermore, the function
$\sum\limits_{j=1}^{m}\gamma_j^{}\;\!{\frak d}\;\!\omega_j^{}$
is the cofactor of the exponential partial integral 
\\[1ex]
\mbox{}\hfill
$
\displaystyle
\exp\sum\limits_{j=1}^{m}\gamma_j^{}\;\!(\omega_j^{}+c_j^{}).
\hfill
$
}
\\[1.5ex]
\indent
{\sl Proof} follows from Theorem 3.2 under 
$\omega=\exp\sum\limits_{j=1}^{m}\gamma_j^{}\;\!(\omega_j^{}+c_j^{}).\ \k$
\vspace{0.5ex}

Theorem 3.3 is the analog of Theorem 1.3 in the case 
\\[1.5ex]
\mbox{}\hfill
$
{\rm g}_j^{}=\exp\bigl(\gamma_j^{}\;\!(\omega_j^{}+c_j^{})\bigr),
\ \
j=1,\ldots,m.
\hfill
$
\\[2ex]
\indent
{\bf  Property 3.1.}
{\it
Let $\varphi\in C^1T^{\;\!\prime}.$ Then} 
\\[1.5ex]
\mbox{}\hfill
$
\displaystyle
({\rm g}, M)\in \text{\rm J}_{_{\Omega}}
\iff
\bigl({\rm g}\exp\varphi, M+{\sf D} \varphi\bigr)\in \text{\rm J}_{_{\Omega}}.
\hfill
$
\\[1.75ex]
\indent
{\sl Proof}\;\! follows from Theorem 1.1 and 
\\[1.5ex]
\mbox{}\hfill                        
$
{\frak d}\bigl({\rm g}(t,x) e^{\varphi(t)}\!\bigr)\! =\!
e^{\varphi(t)}\;\!{\frak d}\;\!{\rm g}(t,x)\! +
{\rm g}(t,x)\;\!{\sf D}e^{\varphi(t)}\! =\!
e^{\varphi(t)}{\rm g}(t,x)
\bigl(M(t,x)\!+\!{\sf D} \varphi(t)\bigr)
$
 for all $\!(t,x)\!\in\!\Omega. \k
\hfill 
$
\\[2ex]
\indent
{\bf  Property 3.2.}
{\it
Let $\varphi\in C^1T^{\;\!\prime}.$ Then} 
\\[1.5ex]
\mbox{}\hfill
$
\displaystyle
(\exp\omega, M)\in \text{\rm E}_{_{\Omega}}
\iff
\bigl(\exp(\omega+\varphi), M+{\sf D} \varphi\bigr)\in \text{\rm E}_{_{\Omega}}.
\hfill
$
\\[1.75ex]
\indent
{\sl Proof}\;\! follows from Property 3.1 under ${\rm g}=\exp\omega.\ \k$  
\vspace{0.75ex}

In particular, from Property 3.2, we obtain
\vspace{0.75ex}

{\bf  Property 3.3.}
{\it
If $c\in\R,$ then we have} 
\\[1.5ex]
\mbox{}\hfill
$
\displaystyle
(\exp\omega, M)\in \text{\rm E}_{_{\Omega}}
\iff
\bigl(\exp(\omega+c), M\bigr)\in \text{\rm E}_{_{\Omega}}.
\hfill
$
\\[-1.75ex]

\newpage

{\bf  Property 3.4.}
{\it
If  
$\bigl(\exp\omega_j^{}, M\bigr)\in \text{\rm E}_{_{\Omega}},\ 
\lambda_j^{}\in\R\backslash\{0\}, \ j=1,\ldots,m,$ then}
\\[1.25ex]
\mbox{}\hfill
$
\displaystyle
\biggl(\,\sum\limits_{j=1}^{m}\lambda_j^{}\exp\omega_j^{}\;\!,\;\! M\biggr)\in \text{\rm J}_{_{\Omega}}.
\hfill
$
\\[1.5ex]
\indent
{\sl Proof}\;\! follows from Property 1.5. $\k$  
\vspace{0.75ex}

{\bf  Property 3.5.}
{\it
If 
$\rho_j^{},\lambda_j^{}\in\R\backslash\{0\}, \ 
\bigl(\exp\omega_j^{}, \rho_j^{}\;\!M_0^{}\bigr)\in \text{\rm E}_{_{\Omega}},\ 
j=1,\ldots,m,$ then}
\\[1.5ex]
\mbox{}\hfill
$
\displaystyle
\biggl(\,\sum\limits_{j=1}^{m}\lambda_j^{}
\exp\dfrac{\omega_j^{}}{\rho_j^{}}\;\!,\;\! M_0^{}\biggr)\in \text{\rm J}_{_{\Omega}}.
\hfill
$
\\[1.5ex]
\indent
{\sl Proof}\;\! follows from Property 1.8. $\k$  
\vspace{1ex}

{\bf  Property 3.6.}
\vspace{0.35ex}
{\it
Suppose  
$\bigl(\exp\omega_j^{}, M_j^{}\bigr)\in \text{\rm E}_{_{\Omega}}$ and
$\gamma_j^{}\in\R\backslash\{0\}, \ j=1,\ldots,m.$ Then
$
\biggl(\,\exp\sum\limits_{j=1}^{m}\gamma_j^{}\;\!\omega_j^{}\;\!,\;\! M\biggr)\in \text{\rm E}_{_{\Omega}}
$
if and only if the identity {\rm (1.10)} holds}.
\vspace{0.75ex}

{\sl Proof}\;\! follows from Property 1.9. $\k$  
\vspace{0.75ex}

{\bf  Corollary 3.1.}
{\it
Suppose
$\gamma_j^{}, \rho_j^{}\in\R\backslash\{0\}, \ 
\bigl(\exp\omega_j^{}, \rho_j^{}\;\!M_0^{}\bigr)\in \text{\rm E}_{_{\Omega}},\ 
j=1,\ldots,m.$ Then
$
\biggl(\,\exp\sum\limits_{j=1}^{m}
\gamma_j^{}\;\!\omega_j^{}\;\!,\;\! M\biggr)\in \text{\rm E}_{_{\Omega}}
$
if and only if the identity {\rm (1.13)} holds}.
\vspace{1ex}

{\bf  Corollary 3.2.}
{\it
If $\gamma\in\R\backslash\{0\},$ then we claim that} 
\\[1.5ex]
\mbox{}\hfill
$
\displaystyle
(\exp\omega, M)\in \text{\rm E}_{_{\Omega}}
\iff
\bigl(\exp(\gamma\;\!\omega), \gamma\;\!M\bigr)\in \text{\rm E}_{_{\Omega}}.
\hfill
$
\\[1.75ex]
\indent
Corollary 3.2 is the analog of Properties 1.6 and 1.7 in the case of 
exponential partial integral.
\vspace{0.35ex}

{\bf  Property 3.7.}
{\it
If
$
\bigl(\exp\omega_j^{}, M_j^{}\bigr)\in \text{\rm E}_{_{\Omega}},\ 
\gamma,\gamma_j^{}\in\R\backslash\{0\}, \ j=1,\ldots,m, \
{\rm g}^{\gamma}\in C^1\Omega,$ then}
\\[1.5ex]
\mbox{}\hfill
$
\displaystyle
\biggl(\,
{\rm g}^{\gamma}\exp\sum\limits_{j=1}^{m}\gamma_j^{}\;\!\omega_j^{}\;\!,\;\! M\biggr)\in 
\text{\rm J}_{_{\Omega}}
\iff
\biggl(\,
{\rm g},\, \dfrac{1}{\gamma}\,\biggl(M-
\sum\limits_{j=1}^{m}\gamma_j^{}\;\!M_j^{}\biggr)\biggr)\in 
\text{\rm J}_{_{\Omega}}.
\hfill
$
\\[1.25ex]
\indent
{\sl Proof}\;\! follows from Property 1.10. $\k$  
\vspace{0.5ex}

{\bf  Corollary 3.3.}
{\it
If
$
\bigl(\exp\omega_\tau^{}, M_\tau^{}\bigr)\!\in\! \text{\rm E}_{_{\Omega}},\,
\tau\!=\!1,\ldots, m\!-\!1,\ 
\gamma_j^{}\in\R\backslash\{0\}, \, j\!=\!1,\ldots,m,$ then}
\\[1.25ex]
\mbox{}\hfill
$
\displaystyle
\biggl(\,
\exp\sum\limits_{j=1}^{m}\gamma_j^{}\;\!\omega_j^{}\;\!,\;\! M\biggr)\in 
\text{\rm E}_{_{\Omega}}
\iff
\biggl(\,
\exp\omega_{m}^{}\;\!,\,
 \dfrac{1}{\gamma_{m}^{}}\,\biggl(M-
\sum\limits_{\tau=1}^{m-1}\gamma_\tau^{}\;\!M_\tau^{}\biggr)\biggr)\in 
\text{\rm E}_{_{\Omega}}.
\hfill
$
\\[1.75ex]
\indent
{\bf  Property 3.8.}
\vspace{0.25ex}
{\it
Suppose 
$
\bigl(\exp\omega_\nu^{}, M_\nu^{}\bigr)\in \text{\rm E}_{_{\Omega}},\ 
\nu=1,\ldots, s,\ s\leqslant m-2,\  
\gamma_j^{}\in\R\backslash\{0\}, 
\linebreak 
j=1,\ldots,m,\ 
{\rm g}_k^{{}^{\scriptsize \gamma_k^{}}}\!\in C^1\Omega,\ 
k=s+1,\ldots,m.$ Then we have
\\[1.5ex]
\mbox{}\hfill
$
\displaystyle
\biggl(\,
\prod\limits_{k=s+1}^{m}
{\rm g}_k^{{}^{\scriptsize \gamma_k^{}}}
\exp\sum\limits_{\nu=1}^{s}
\gamma_\nu^{}\;\!\omega_\nu^{}\;\!,\;\! M\biggr)\in 
\text{\rm J}_{_{\Omega}}
\iff
\biggl(\,
\prod\limits_{k=s+1}^{m}
{\rm g}_k^{{}^{\scriptsize \gamma_k^{}}},\,
M-\sum\limits_{\nu=1}^{s}\gamma_\nu^{}\;\!M_\nu^{}\biggr)\in 
\text{\rm J}_{_{\Omega}}.
\hfill
$
\\[1.25ex]
Moreover, there exist functions $M_k^{}\in C^1\Omega,\ k=s+1,\ldots, m,$ such that
the identities {\rm (1.14)} and {\rm (1.15)} are true.
}
\vspace{0.35ex}

{\sl Proof}\;\! follows from Property 1.11. $\k$
\vspace{0.75ex}

{\bf  Property 3.9.}
\vspace{0.25ex}
{\it
Suppose 
$
\bigl(\exp\omega_\nu^{}, M_\nu^{}\bigr)\in \text{\rm E}_{_{\Omega}},\ 
\nu=1,\ldots, s,\ s\leqslant m-1,\  
({\rm g}_\xi^{}, M_\xi^{})\in \text{\rm J}_{_{\Omega}},
\linebreak 
{\rm g}_\xi^{{}^{\scriptsize \gamma_\xi^{}}}\!\in C^1\Omega,
\ 
\xi=s+1,\ldots,m,\ 
\gamma_j^{}\in\R\backslash\{0\}, 
\ j=1,\ldots,m. 
$ Then we have
\\[1.5ex]
\mbox{}\hfill
$
\displaystyle
\biggl(\,
\prod\limits_{\xi=s+1}^{m}
{\rm g}_\xi^{{}^{\scriptsize \gamma_\xi^{}}}
\exp\sum\limits_{\nu=1}^{s}
\gamma_\nu^{}\;\!\omega_\nu^{}\;\!,\;\! M\biggr)\in 
\text{\rm J}_{_{\Omega}}
\hfill
$
\\[1.5ex]
if and only if the identity {\rm (1.10)} holds.
}
\vspace{0.5ex}

{\sl Proof}\;\! follows from Property 1.9. $\k$
\vspace{1ex}

{\bf Property 3.10.} 
\vspace{0.5ex}
{\it
Suppose functions 
$p,q\in \text{\rm P}_{_{\!\Xi^{\;\!\prime}}}$ are relatively prime,
the set $\Omega_{_0}\subset\Xi^{\;\!\prime}$ such that   
$
p(t,x)\ne 0$ for all $(t,x)\in\Omega_{_0}$ and
$p(t,x)= 0$ for all $(t,x)\in {\sf C}_{_{\!\Xi^{\;\!\prime}}}\Omega_{_0}\;\!.
$
Then
$
\Bigl(\exp\dfrac{q}{p}\,, N\Bigr)\in 
\text{\rm E}_{_{\Omega_{{}_{\tiny\;\! 0}}}}
$
if and only if
$(p,M)\in \text{\rm A}_{_{\!\Xi^{\;\!\prime}}}$
and the identity holds
\\[1.75ex]
\mbox{}\hfill                            % (3.2)
$
{\frak d}\;\!q(t,x)=
q(t,x)\;\!M(t,x)+p(t,x)\;\!N(t,x)
$
\ for all 
$(t,x)\in \Xi^{\;\!\prime},
$
\hfill {\rm (3.2)}
\\[1.75ex]
where the function $N\in \text{\rm P}_{_{\!\Xi^{\;\!\prime}}}$ and has the degree 
$\deg_{x}^{}N\leq d-1.$
}
\vspace{0.75ex}

{\sl Proof.}
By Theorem 3.1,
$
\Bigl(\exp\dfrac{q}{p}\,, N\Bigr)\in 
\text{\rm E}_{_{\Omega_{{}_{\tiny\;\! 0}}}} 
$
if and only if the derivative by virtue of system (0.1)
\\[0.75ex]
\mbox{}\hfill  
$
{\frak d}\,\dfrac{q(t,x)}{p(t,x)}=N(t,x)
$
\ for all 
$(t,x)\in\Omega_{_{0}}
\hfill 
$
\\[2ex]
under the conditions $N\in \text{\rm P}_{_{\!\Xi^{\;\!\prime}}}$  and $\deg_{x}^{}N\leq d-1.$
\vspace{0.75ex}

This implies that
$
\Bigl(\exp\dfrac{q}{p}\,, N\Bigr)\in 
\text{\rm E}_{_{\Omega}} 
$
if and only if the identity holds
\\[1.5ex]
\mbox{}\hfill   % (3.3)
$
p(t,x)\;\!{\frak d}\;\!q(t,x)-
q(t,x)\;\!{\frak d}\;\!p(t,x)=
p^2(t,x)\;\!N(t,x)
$
\ for all $(t,x)\in \Xi^{\;\!\prime},
$
\hfill (3.3)
\\[1.75ex]
where the function $N\in \text{\rm P}_{_{\!\Xi^{\;\!\prime}}}$ and has the degree 
$\deg_{x}^{}N\leq d-1.$
\vspace{0.75ex}

{\sl Necessity.}
From the identity (3.3) it follows that
\\[1.5ex]
\mbox{}\hfill                       % (3.4)
$
{\frak d}\;\!q(t,x)=
p(t,x)\;\!N(t,x)+
\dfrac{{\frak d}\;\!p(t,x)}{p(t,x)}\ q(t,x)
$
\ for all 
$(t,x)\in\Omega_{_{0}},
$
\hfill  (3.4)
\\[1.5ex]
where the function $N\in \text{\rm P}_{_{\!\Xi^{\;\!\prime}}}$ and has the degree
$\deg_{x}^{}N\leq d-1.$
\vspace{0.75ex}

Since
\vspace{0.5ex}
$p,q, {\frak d}\;\!p, N\in \text{\rm P}_{_{\!\Xi^{\;\!\prime}}}$
and the functions $p,\,q$ are relatively prime, we see that from the identity (3.4) it follows that  
\vspace{0.5ex}
$\dfrac{{\frak d}\;\!p}{p}\in \text{\rm P}_{_{\!\Omega_{{}_{\tiny\;\! 0}}}}.$
Therefore there exists a function
$M\in \text{\rm P}_{_{\!\Xi^{\;\!\prime}}}$ such that
the identity (2.2) is true. 
Thus, by Theorem 2.1, we have 
$(p,M)\in \text{\rm A}_{_{\!\Xi^{\;\!\prime}}}.$
\vspace{0.5ex}
Now, using the identity (3.4) under the condition (2.2), we obtain the identity (3.2), 
\vspace{0.5ex}
where the function $N\in \text{\rm P}_{_{\!\Xi^{\;\!\prime}}}$ and the degree $\deg_{x}^{}N\leq d-1.$
\vspace{0.5ex}

{\sl Sufficiency.}
Using the identity (3.2) under the condition (2.2), we get the identity (3.4).
Further, multiplying both sides of the identity (3.4) by $p,$ we obtain the identity (3.3). 
\vspace{0.5ex}
This proves that 
$
\Bigl(\exp\dfrac{q}{p}\,, N\Bigr)\in 
\text{\rm E}_{_{\Omega_{{}_{\tiny\;\! 0}}}}. \ \k
$
\vspace{1.25ex}

{\bf Property 3.11.} 
\vspace{0.5ex}
{\it
Suppose functions
$u,v\in \text{\rm P}_{_{\!\Xi^{\;\!\prime}}}$ 
are relatively prime, the set
$\Omega_{_0}\subset\Xi^{\;\!\prime}$ such that  
$u(t,x)\ne 0$ for all $(t,x)\in\Omega_{_0}$ and 
$u(t,x)= 0$ for all $(t,x)\in {\sf C}_{_{\!\Xi^{\;\!\prime}}}\Omega_{_0}\;\!.$
Then
\\[1.5ex]
\mbox{}\hfill
$
\Bigl(\exp\arctan\dfrac{v}{u}\,, V\Bigr)\in
\text{\rm E}_{_{\Omega_{{}_{\tiny\;\! 0}}}}
\hfill
$
\\[1.5ex]
if and only if there exists a function 
$U\in \text{\rm P}_{_{\!\Xi^{\;\!\prime}}}$
such that the system of identities holds
\\[1.5ex]
\mbox{}\hfill                            % (3.5)
$
{\frak d}\;\!u(t,x)=
u(t,x)\;\!U(t,x)-v(t,x)\;\!V(t,x)
$
for all 
$(t,x)\in \Xi^{\;\!\prime},
\hfill
$
\\
\mbox{}\hfill {\rm (3.5)}
\\
\mbox{}\hfill                           
$
{\frak d}\;\!v(t,x)=
u(t,x)\;\!V(t,x)+v(t,x)\;\!U(t,x)
$
 for all $(t,x)\in \Xi^{\;\!\prime},
\hfill
$
\\[1.5ex]
where the function $V\in \text{\rm P}_{_{\!\Xi^{\;\!\prime}}}$ and has the degree 
$\deg_{x}^{}V\leq d-1.$
}
\vspace{0.75ex}

{\sl Proof.}
By Theorem 3.1, 
$
\Bigl(\exp\arctan\dfrac{v}{u}\,, V\Bigr)\in 
\text{\rm E}_{_{\Omega_{{}_{\tiny\;\! 0}}}}
$
if and only if 
\\[1.75ex]
\mbox{}\hfill                               % (3.6)
$
{\frak d}\,\arctan\dfrac{v(t,x)}{u(t,x)}=V(t,x)
$
\ for all $(t,x)\in\Omega_{_{0}},
$
\hfill  (3.6)
\\[2ex]
where the function $V\in \text{\rm P}_{_{\!\Xi^{\;\!\prime}}}$ and has the degree
$\deg_{x}^{}V\leq d-1.$
\vspace{0.75ex}

Using the derivative by virtue of system (0.1)
\\[1.5ex]
\mbox{}\hfill                           
$
{\frak d}\,\arctan\dfrac{v(t,x)}{u(t,x)}=
\dfrac{u(t,x)\,{\frak d}\;\!v(t,x)-v(t,x)\,{\frak d}\;\!u(t,x)}{u^2(t,x)+v^2(t,x)}
$
\ for all $(t,x)\in\Omega_{_{0}},
\hfill 
$
\\[2ex]
from the identity (3.6), we get
\\[1.75ex]
\mbox{}\hfill                           
$
{\frak d}\;\!v(t,x)=
u(t,x)\;\!V(t,x)+v(t,x)\,
\dfrac{{\frak d}\;\!u(t,x)+v(t,x)\;\!V(t,x)}{u(t,x)}
$
\ for all $(t,x)\in\Omega_{_{0}}.
\hfill 
$
\\[2ex]
\indent
Since the functions
\vspace{0.5ex}
$u,v, {\frak d}\;\!u,{\frak d}\;\!v, V\in \text{\rm P}_{_{\!\Xi^{\;\!\prime}}}$ and
the functions $u,\,v$ are relatively prime, we see that there exists a function 
$U\in \text{\rm P}_{_{\!\Xi^{\;\!\prime}}}$ such that
\\[2ex]
\mbox{}\hfill 
$
{\frak d}\;\!u(t,x)+
v(t,x)\;\!V(t,x)=
u(t,x)\;\!U(t,x)
$
\ for all $(t,x)\in \Xi^{\;\!\prime}.
\hfill
$
\\[1.75ex]
\indent
Thus the identity (3.6) is true 
\vspace{0.5ex}
if and only if the system of identities (3.5) is true,
where $V\in \text{\rm P}_{_{\!\Xi^{\;\!\prime}}}$ and $\deg_{x}^{}V\leq d-1.\ \k$
\vspace{1.25ex}

{\bf Property 3.12.} 
\vspace{0.5ex}
{\it
Let  
$p\in \text{\rm P}_{_{\!\Xi^{\;\!\prime}}}.$
Then we claim that
\\[1.35ex]
\mbox{}\hfill
$
\bigl(\exp\arctan p, M\bigr)\in
\text{\rm E}_{_{\Xi^{\;\!\prime}}}
\hfill
$
\\[1.5ex]
if and only if the derivative by virtue of system {\rm(0.1)}
\\[1.75ex]
\mbox{}\hfill                            % (3.7)
$
{\frak d}\;\!p(t,x)=
\bigl(1+p^2(t,x)\bigr)\;\!M(t,x)
$
\ for all $(t,x)\in \Xi^{\;\!\prime}.
$
\hfill {\rm (3.7)}
\\[2ex]
Moreover, the cofactor $M$ of exponential partial integral
\vspace{0.35ex}
$\exp\arctan p$ has the degree $\deg_{x}^{}M\leq d-1-\deg_{x}^{}p.$
}
\vspace{0.75ex}

{\sl Proof}\;\! follows from
\vspace{0.25ex}
Property 3.11 under the condition $v=p,\ u=1.$
Besides, using the identity (3.7), we have
$\deg_{x}^{}M\leq d-1-\deg_{x}^{}p.\ \k$
\vspace{1.25ex}

{\bf Property 3.13.} 
{\it
Let $f\in C^1\Omega.$ Then we have}
\\[1.75ex]
\mbox{}\hfill                       
$
\bigl(\exp\arctan f, M\bigr)\in
\text{\rm E}_{_{\Omega}}
\iff
\bigl(\exp{\rm arccotan} f, {}-M\bigr)\in
\text{\rm E}_{_{\Omega}}.
\hfill
$
\\[1.75ex]
\indent
{\sl Proof}\;\!
is based on Theorem 3.1 and on the identity
\\[1.75ex]
\mbox{}\hfill                       
$
{\rm arccotan} f(t,x)=\dfrac{\pi}{2}-\arctan f(t,x)
$
\ for all $(t,x)\in\Omega. \ \k
\hfill
$
\\[2ex]
\indent
{\bf Property 3.14.} 
\vspace{0.5ex}
{\it
Suppose
$k\in\N,\  F\in \text{\rm I}_{_{\Omega}},$ a function
$f\in C^1\Omega$ such that $F$ and $f$ are relatively prime, 
the set $\Omega_{_0}\subset\Omega$ such that
\vspace{0.5ex}
$
F(t,x)\ne 0$ for all $(t,x)\in\Omega_{_0}$ and
$F(t,x)= 0$ for all $(t,x)\in {\sf C}_{_{\Omega}}\Omega_{_0}\;\!.$
Then we claim that
\\[1.5ex]
\mbox{}\hfill
$
\Bigl(\exp\dfrac{f}{F^k}\,, M\Bigr)\in 
\text{\rm E}_{_{\Omega_{{}_{\tiny\;\! 0}}}}
\hfill
$
\\[1.5ex]
if and only if the identity holds
\\[2ex]
\mbox{}\hfill                      
$
{\frak d}\;\!f(t,x)=
F^k(t,x)\;\!M(t,x)
$
\ for all $(t,x)\in \Omega,
\hfill
$
\\[2ex]
where the function $M\in \text{\rm P}_{_{\!\Xi^{\;\!\prime}}}$  and has the degree
$\deg_{x}^{}M\leq d-1.$
}
\vspace{0.75ex}

{\sl Proof.}
\vspace{0.5ex}
By Theorem 3.1, using the identity (0.3) from Theorem 0.2, we have
$
\Bigl(\exp\dfrac{f}{F^k}\,, M\Bigr)\in 
\text{\rm E}_{_{\Omega_{{}_{\tiny\;\! 0}}}}
$
if and only if the identity on the domain $\Omega_{_{0}}$ holds 
\\[2ex]
\mbox{}\hfill                               
$
{\frak d}\;\!\dfrac{f(t,x)}{F^k(t,x)}=
\dfrac{F^k(t,x)\,{\frak d}f(t,x)-kf(t,x)\;\!F^{k-1}(t,x)\,{\frak d}F(t,x)}{F^{\;\!2k}(t,x)}=
\dfrac{{\frak d}f(t,x)}{F^k(t,x)}=M(t,x),
\hfill
$
\\[2.25ex]
where the function $M\in \text{\rm P}_{_{\!\Xi^{\;\!\prime}}}$ has the degree
$\deg_{x}^{}M\leq d-1.\ \k$
\vspace{1.25ex}

The application of Property 3.14 is given in [56, pp. 52 -- 55 ].

\newpage

\mbox{}
\\[-1.75ex]
\centerline{
{\bf  4. Conditional partial integrals}
}
\\[1.5ex]
\indent
{\bf Definition 4.1.}
\vspace{0.75ex}
{\it
A function $\exp p$ is called 
a \textit{\textbf{conditional partial integral}}
with cofactor $M$ on the domain $\Xi^{\;\!\prime}$ of system {\rm (0.1)}
if $p\in \text{\rm P}_{_{\!\Xi^{\;\!\prime}}},\  
\bigl(\exp p, M\bigr)\in \text{\rm E}_{_{\Xi^{\;\!\prime}}}.$
}
\vspace{1ex}

By $\text{\rm F}_{_{\!\Xi^{\;\!\prime}}}$ 
denote the set of conditional partial integrals on the domain $\Xi^{\;\!\prime}$
\vspace{0.5ex}
of system (0.1).
 
From the definitions of  conditional and exponential partial integrals 
(Definitions 4.1 and 3.1) it follows that the set
$\text{\rm F}_{_{\!\Xi^{\;\!\prime}}}\subset
\text{\rm E}_{_{\Xi^{\;\!\prime}}}\subset 
\text{\rm J}_{_{\!\Xi^{\;\!\prime}}}.$
\vspace{1ex}

The phrase "the function $\exp p$ 
\vspace{0.35ex}
is a conditional partial integral with cofactor $M$ on
the domain $\Xi^{\;\!\prime}$ of system (0.1)" is denoted by 
$\bigl(\exp p, M\bigr)\in \text{\rm F}_{_{\!\Xi^{\;\!\prime}}}.$
\vspace{0.75ex}

Definition 4.1 can be expressed by the equivalence
\\[1.5ex]
\mbox{}\hfill               % (4.1)
$
\bigl(\exp p, M\bigr)\in \text{\rm F}_{_{\!\Xi^{\;\!\prime}}}
\iff
\bigl(\exp p, M\bigr)\in \text{\rm E}_{_{\Xi^{\;\!\prime}}}
\ \&\ 
p\in \text{\rm P}_{_{\!\Xi^{\;\!\prime}}}.
$
\hfill {\rm(4.1)}
\\[2ex]
\indent
{\bf Theorem 4.1}\!
\vspace{0.35ex}
(existence criterion of conditional partial integral).\!
{\it
Suppose $p\in \text{\rm P}_{_{\!\!\Xi^{\;\!\prime}}}.$
Then  
$\bigl(\exp p, M\bigr)\in \text{\rm F}_{_{\!\Xi^{\;\!\prime}}}$
if and only if the derivative by virtue of system {\rm(0.1)}
\\[1.5ex]
\mbox{}\hfill               % (4.2)
$
{\frak d}\;\!p(t,x)=M(t,x)
$
\ for all 
$(t,x)\in \Xi^{\;\!\prime}
$
\hfill {\rm(4.2)}
\\[1.5ex]
under the condition $\deg_{\;\!x}^{} M\leq d-1.$
}
\vspace{0.5ex}

{\sl Proof}\;\! uses the equivalence (4.1) and 
follows from Theorem 3.1 under 
\vspace{0.75ex}
$\omega=p\in \text{\rm P}_{_{\!\Xi^{\;\!\prime}}}.\ \k$

Using Theorem 4.1, we have
\vspace{0.5ex}

{\bf Theorem 4.2}\!
\vspace{0.35ex}
(existence criterion of conditional partial integral).\!
{\it
Let $p\in \text{\rm P}_{_{\!\!\Xi^{\;\!\prime}}}.$
Then  
\\[0.5ex]
\mbox{}\hfill
$
\exp p\in \text{\rm F}_{_{\!\Xi^{\;\!\prime}}}
\iff
\deg_{\;\!x}^{}{\frak d}\;\!p\leq d-1.
\hfill
$
\\[1.5ex]
Moreover, the function ${\frak d}\;\!p$ is the cofactor of  
\vspace{0.75ex}
conditional partial integral $\exp p.$}

Theorem 4.2 is the analog of Theorem 3.2 for the case of 
\vspace{0.25ex}
conditional partial integral 
\linebreak
(when $\omega=p\in \text{\rm P}_{_{\!\Xi^{\;\!\prime}}}).$
\vspace{1ex}

{\bf Remark 4.1.}
\vspace{0.75ex}
If $p\in \text{\rm P}_{_{\!\Xi^{\;\!\prime}}},$ then 
${\frak d}\;\!p\in \text{\rm P}_{_{\!\Xi^{\;\!\prime}}}$ and
$\deg_{\;\!x}^{}{\frak d}\;\!p\leq d-1+\deg_{\;\!x}^{}p.$
Thus the requirement $\deg_{\;\!x}^{} M\leq d-1$ 
\vspace{0.75ex}
in the identity (4.2) of Theorem 4.1 
(as the requirement $\deg_{\;\!x}^{}{\frak d}\;\!p\leq d-1$ in Theorem 4.2)
is important.
\vspace{1.25ex}

{\bf  Theorem 4.3.}
\vspace{0.5ex}
{\it
Suppose 
$p_j^{}\in \text{\rm P}_{_{\!\Xi^{\;\!\prime}}},\
\gamma_j^{}\in\R\backslash\{0\},\ c_j^{}\in\R, \ j=1,\ldots,m.$
Then
\\[1.5ex]
\mbox{}\hfill
$
\displaystyle
\exp\sum\limits_{j=1}^{m}
\gamma_j^{}\;\!(p_j^{}+c_j^{})\in \text{\rm F}_{_{\!\Xi^{\;\!\prime}}}
\ \iff \
\deg_{\;\!x}^{} 
\sum\limits_{j=1}^{m}
\gamma_j^{}\,{\frak d}\;\!p_j^{}\leq d-1.
\hfill
$
\\[1.5ex]
Also, the function 
$\sum\limits_{j=1}^{m}\gamma_j^{}\;\!{\frak d}\;\!p_j^{}$
is the cofactor of conditional partial integral  
$\exp\sum\limits_{j=1}^{m}\gamma_j^{}\;\!(p_j^{}+c_j^{}).$}

{\sl Proof}\;\! follows from Theorem 4.2 under the condition 
$p=\sum\limits_{j=1}^{m}\gamma_j^{}\;\!(p_j^{}+c_j^{}).\ \k$
\vspace{0.75ex}

Theorem 4.3 is the analog of Theorem 1.3 
\vspace{0.5ex}
for the case of conditional partial integral 
\linebreak
(when the functions 
${\rm g}_j^{}=\exp\bigl(\gamma_j^{}\;\!(p_j^{}+c_j^{})\bigr),
\ j=1,\ldots,m).$
\vspace{1.5ex}

{\bf  Property 4.1.}
\vspace{1.25ex}
{\it
If $\varphi\in C^1T^{\;\!\prime},$ then} 
$
\bigl(\exp\varphi, {\sf D} \varphi\bigr)\in \text{\rm F}_{_{\!\Xi^{\;\!\prime}}}.
$

{\sl Indeed}, 
\vspace{1ex}
$\varphi\in \text{\rm P}_{_{\!\Xi^{\;\!\prime}}},\
{\frak d}\;\!\varphi(t)={\sf D}\varphi(t)$ for all $(t,x)\in \Xi^{\;\!\prime},\
\deg_{\;\!x}^{}{\sf D}\varphi=0\leq d-1.$
By Theorem~4.1,  
$
\bigl(\exp\varphi, {\sf D} \varphi\bigr)\in \text{\rm F}_{_{\!\Xi^{\;\!\prime}}}.\ \k
$
\vspace{1.25ex}

{\bf  Property 4.2.}
{\it
Let $\varphi\in C^1T^{\;\!\prime}, \
\gamma_j^{}\in\R\backslash\{0\},\ c_j^{}\in\R, \ j=1,\ldots,m,\ m\leq n.$
Then
\\[1.5ex]
\mbox{}\hfill
$
\displaystyle
\exp\biggl(\varphi(t)+\sum\limits_{j=1}^{m}
\gamma_j^{}\;\!(x_j^{}+c_j^{})\biggr)\in \text{\rm F}_{_{\!\Xi^{\;\!\prime}}}
\ \iff \
\deg_{\;\!x}^{} 
\sum\limits_{j=1}^{m}
\gamma_j^{}\,X_j^{}\leq d-1.
\hfill
$
\\[1.5ex]
Moreover, the function  
${\sf D} \varphi+\sum\limits_{j=1}^{m}\gamma_j^{}\;\!X_j^{}$
is the cofactor of conditional partial integral}
\\[1.5ex]
\mbox{}\hfill
$
\displaystyle
\exp\Bigl(\varphi(t)+\sum\limits_{j=1}^{m}
\gamma_j^{}\;\!(x_j^{}+c_j^{})\Bigr).
\hfill
$
\\[1.5ex]
\indent
{\sl Proof}\;\! 
is based on Theorem 4.1 and 
\\[1.5ex]
\mbox{}\hfill
$
\displaystyle
{\frak d}\;\!\biggl(\varphi(t)+\sum\limits_{j=1}^{m}
\gamma_j^{}\;\!(x_j^{}+c_j^{})\biggr)=
{\sf D} \varphi(t)+\sum\limits_{j=1}^{m}\gamma_j^{}\;\!X_j^{}(t,x)
$
\ for all $(t,x)\in \Xi^{\;\!\prime}.\ \k
\hfill
$
\\[1.75ex]
\indent
In particular, if $m=1,$ then we obtain
\vspace{1ex}

{\bf  Property 4.3.}
{\it
Suppose $\varphi\in C^1T^{\;\!\prime}, \
\gamma\in\R\backslash\{0\},\ c\in\R.$
Then we claim that
\\[2ex]
\mbox{}\hfill
$
\displaystyle
\exp\bigl(\varphi(t)+\gamma\;\!(x_j^{}+c)\bigr)\in \text{\rm F}_{_{\!\Xi^{\;\!\prime}}}
\ \iff \
\deg_{\;\!x}^{} X_j^{}\leq d-1,
\quad 
j\in \{1,\ldots, n\}. 
\hfill
$
\\[2ex]
Further, the function
${\sf D} \varphi+\gamma\;\!X_j^{}$
is the cofactor of conditional partial integral}
\\[1.5ex]
\mbox{}\hfill
$
\exp\bigl(\varphi(t)+\gamma\;\!(x_j^{}+c)\bigr).
\hfill
$
\\[5.25ex]
\centerline{
{\bf  5. Multiple polynomial partial integrals}
}
\\[1.5ex]
\indent
{\bf Definition 5.1.}
\vspace{0.35ex}
{\it
A polynomial partial integral $p$ with cofactor $M$ 
on the domain $\Xi^{\;\!\prime}$ of system {\rm (0.1)} is said to be 
\textit{\textbf{multiple}}
\vspace{0.25ex}
if there exist a natural number $h$ and  a function
$q\in \text{\rm P}_{_{\!\Xi^{\;\!\prime}}}$
\vspace{0.5ex}
relatively prime to the function $p$
such that the derivative by virtue of system {\rm (0.1)} 
\\[1.5ex]
\mbox{}\hfill                             % (5.1)
$
\displaystyle
{\frak d}\,\dfrac{q(t,x)}{p^{\;\!h}(t,x)}=N(t,x)
$
\ for all $(t,x)\in \Omega_{_0},
$
\hfill {\rm (5.1)}
\\[2.25ex]
where the function $N\in \text{\rm P}_{_{\!\Xi^{\;\!\prime}}}$ 
has the degree $\deg_{\;\!x}^{} N\leq d-1,$
\vspace{1ex}
the set 
$\Omega_{_0}\subset\Xi^{\;\!\prime}$ such that   
$
p(t,x)\ne 0$ for all $(t,x)\in\Omega_{_0}$ and 
$p(t,x)= 0$ for all $(t,x)\in {\sf C}_{_{\Xi^{\;\!\prime}}}\Omega_{_0}\;\!.$
}
\vspace{1.25ex}

By $\text{\rm B}_{_{\Xi^{\;\!\prime}}}$ 
\vspace{0.5ex}
denote the set of multiple polynomial partial integrals on the domain 
$\Xi^{\;\!\prime}$ of system (0.1).
From Definition 5.1 it follows that the set 
$\text{\rm B}_{_{\Xi^{\;\!\prime}}}\subset
\text{\rm A}_{_{\Xi^{\;\!\prime}}}\subset 
\text{\rm J}_{_{\!\Xi^{\;\!\prime}}}.$
\vspace{1ex}

The phrase "the polynomial partial integral $p$ with cofactor $M$ 
\vspace{0.5ex}
on the domain $\Xi^{\;\!\prime}$ of system {\rm (0.1)} is multiple 
and the identity (5.1) holds" is denoted by 
\vspace{1.25ex}
$\bigl((p, M), (h,q,N)\bigr)\in \text{\rm B}_{_{\Xi^{\;\!\prime}}}.$

{\bf Theorem 5.1}
\vspace{0.35ex}
(existence criterion of multiple polynomial partial integral).
{\it
$\bigl((p, M), (h,q,N)\bigr)\in \text{\rm B}_{_{\Xi^{\;\!\prime}}}$
\vspace{0.5ex}
if and only if the identities {\rm(2.2)} and {\rm(5.1)} are true, 
where the number $h\in\N,$ the functions
\vspace{0.75ex}
$q, N\in \text{\rm P}_{_{\!\Xi^{\;\!\prime}}},$
the functions $p$ and $q$ are relatively prime, $\deg_{\;\!x}^{} N\leq d-1,$
\vspace{0.75ex}
the set
$\Omega_{_0}\subset\Xi^{\;\!\prime}$ such that   
$p(t,x)\ne 0$ for all $(t,x)\in\Omega_{_0}$ and
$p(t,x)= 0$ for all $(t,x)\in {\sf C}_{_{\Xi^{\;\!\prime}}}\Omega_{_0}\;\!.$
}
\vspace{0.75ex}

{\sl Proof}\;\! follows from Definition 5.1 and Theorem 2.1 
(taking into account Remark 2.1). $\k$
\vspace{1ex}

{\bf Theorem 5.2}
\vspace{0.35ex}
(existence criterion of multiple polynomial partial integral).
{\it
$\bigl((p, M), (h,q,N)\bigr)\in \text{\rm B}_{_{\Xi^{\;\!\prime}}}$
if and only if the identities {\rm(2.2)} and 
\\[2ex]
\mbox{}\hfill                             % (5.2)
$
\displaystyle
{\frak d}\;\!q(t,x)=h\;\!q(t,x)\;\!M(t,x)+p^{\;\!h}(t,x)\;\!N(t,x)
$
\ for all 
$(t,x)\in \Xi^{\;\!\prime}
$
\hfill {\rm (5.2)}
\\[2.25ex]
are true, where the number $h\in\N,$ the functions
\vspace{0.5ex}
$q, N\in \text{\rm P}_{_{\!\Xi^{\;\!\prime}}},$
the functions $p$ and $q$ are relatively prime, $\deg_{\;\!x}^{} N\leq d-1.$
}
\vspace{1ex}

{\sl Proof}\;\! is based on Theorem 5.1.
The identity (5.1) under the condition (2.2) is true
if and only if the identity (5.2) is true, because

\newpage

\mbox{}
\\[-2.5ex]
\mbox{}\hfill                           
$
\displaystyle
{\frak d}\,\dfrac{q(t,x)}{p^{\;\!h}(t,x)}=
\dfrac{p^{\;\!h}(t,x)\, {\frak d}\;\!q(t,x)-
h\;\!p^{\;\!h-1}(t,x)\;\!q(t,x)\, {\frak d}\;\!p(t,x)}{p^{\;\!2h}(t,x)}=
\hfill
$
\\[2.25ex]
\mbox{}\hfill                           
$
\displaystyle
=\dfrac{{\frak d}\;\!q(t,x)-
h\;\!q(t,x)\;\!M(t,x)}{p^{\;\!h}(t,x)}
$
\ for all 
$(t,x)\in \Omega_{_0},
\hfill
$
\\[1.75ex]
where the set 
\vspace{0.35ex}
$\Omega_{_0}\subset\Xi^{\;\!\prime}$ such that   
$p(t,x)\ne 0$ for all $(t,x)\in\Omega_{_0},$ and
$p(t,x)= 0$ for all $(t,x)\in {\sf C}_{_{\Xi^{\;\!\prime}}}\Omega_{_0}\;\!.\ \k$
\vspace{1ex}

{\bf Lemma 5.1.}
\vspace{0.5ex}
{\it
Suppose $h\in\N,$ functions
$p, q\in \text{\rm P}_{_{\!\Xi^{\;\!\prime}}}\!$ are relatively prime, 
the set
$\Omega_{_0}\subset\Xi^{\;\!\prime}$ such that   
\vspace{0.5ex}
$p(t,x)\ne 0$ for all $(t,x)\in\Omega_{_0}$ and
$p(t,x)= 0$ for all $(t,x)\in {\sf C}_{_{\Xi^{\;\!\prime}}}\Omega_{_0}.$
Then
$
\Bigl(\exp\dfrac{q}{p^{\;\!h}}\,, N\Bigr)\!\in\! 
\text{\rm E}_{_{\Omega_{{}_{\tiny\;\! 0}}}}
$
if and only if
\vspace{0.5ex}
$(p, M)\in \text{\rm A}_{_{\Xi^{\;\!\prime}}}$
and the identity {\rm (5.2)} holds, 
where the function $N\in \text{\rm P}_{_{\!\Xi^{\;\!\prime}}}$  and has the degree
$\deg_{x}^{}N\leq d-1.$
}
\vspace{0.75ex}

{\sl Proof.}
By Theorem 2.2, we have
\\[1.5ex]
\mbox{}\hfill                           
$
\displaystyle
(p, M)\in \text{\rm A}_{_{\Xi^{\;\!\prime}}}
\iff 
\bigl(p^{\;\!h}, h\;\!M\bigr)\in \text{\rm A}_{_{\Xi^{\;\!\prime}}}.
\hfill
$
\\[1.5ex]
\indent
Now, from Property 3.10 it follows that 
the statement of Lemma 5.1 is true. $\k$
\vspace{0.75ex}

{\bf Theorem 5.3}\!
\vspace{0.15ex}
(existence criterion of multiple polynomial partial integral).\!\!
{\it
Suppose $\!h\!\in\!\N,\!$ functions
$p, q\in \text{\rm P}_{_{\!\Xi^{\;\!\prime}}}\!$ are relatively prime, 
the set
$\Omega_{_0}\subset\Xi^{\;\!\prime}$ such that   
\vspace{0.5ex}
$p(t,x)\ne 0$ for all $(t,x)\in\Omega_{_0}$ and
$p(t,x)= 0$ for all $(t,x)\in {\sf C}_{_{\Xi^{\;\!\prime}}}\Omega_{_0}.$
Then}
\\[1.5ex]
\mbox{}\hfill                           
$
\displaystyle
\bigl((p, M), (h,q,N)\bigr)\in \text{\rm B}_{_{\Xi^{\;\!\prime}}}
\iff
\Bigl(\exp\dfrac{q}{p^{\;\!h}}\,, N\Bigr)\in 
\text{\rm E}_{_{\Omega_{{}_{\tiny\;\! 0}}}}.
\hfill
$
\\[2ex]
\indent
{\sl Proof}.
Taking into account Theorem 2.1, from Theorem 5.2 and Lemma 5.1, we obtain the statement of 
Theorem 5.3 $\k$
\vspace{0.5ex}

By Theorem 5.3 (using Lemma 5.1), it follows that 
\vspace{0.35ex}
$\bigl((p, M), (h,q,N)\bigr)\in \text{\rm B}_{_{\Xi^{\;\!\prime}}}$
has two cofactors. Namely
$M$ is the cofactor of polynomial partial integral $p$ and 
\vspace{0.35ex}
$N$ is the cofactor of exponential partial integral $\exp\dfrac{q}{p^{\;\!h}}\,.$
\vspace{1ex}

{\bf Property 5.1.}
{\it
Suppose $h\in\N,\ \varphi\in C^1T^{\;\!\prime},$ functions
\vspace{0.5ex}
$p, q\in \text{\rm P}_{_{\!\Xi^{\;\!\prime}}}\!$ are relatively prime, 
the set
$\Omega_{_0}\subset\Xi^{\;\!\prime}$ such that  
$p(t,x)\ne 0$ for all $(t,x)\in\Omega_{_0}$ and 
\vspace{0.25ex}
$p(t,x)= 0$ for all $(t,x)\in {\sf C}_{_{\Xi^{\;\!\prime}}}\Omega_{_0}.$
Then we claim that}
\\[1ex]
\mbox{}\hfill                           
$
\displaystyle
\bigl((p, M), (h,q,N)\bigr)\in \text{\rm B}_{_{\Xi^{\;\!\prime}}}
\iff
\biggl(\exp\Bigl(\;\!\dfrac{q}{p^{\;\!h}}+\varphi\Bigr), N+{\sf D}\;\!\varphi\biggr)\in 
\text{\rm E}_{_{\Omega_{{}_{\tiny\;\! 0}}}}.
\hfill
$
\\[2ex]
\indent
{\sl Proof}.
\vspace{1ex}
Using Property 3.1, from Theorem 5.3, we get the statement of this property. $\!\k$

{\bf Corollary 5.1.}
{\it
Suppose numbers $h\in\N$ and $c\in\R,$ functions
\vspace{0.5ex}
$p, q\in \text{\rm P}_{_{\!\Xi^{\;\!\prime}}}\!$ are relatively prime, 
the set
$\Omega_{_0}\subset\Xi^{\;\!\prime}$ such that  
\vspace{0.5ex}
$p(t,x)\ne 0$ for all $(t,x)\in\Omega_{_0}$ and
$p(t,x)= 0$ for all $(t,x)\in {\sf C}_{_{\Xi^{\;\!\prime}}}\Omega_{_0}.$
Then we have}
\\[1.25ex]
\mbox{}\hfill                           
$
\displaystyle
\bigl((p, M), (h,q,N)\bigr)\in \text{\rm B}_{_{\Xi^{\;\!\prime}}}
\iff
\biggl(\exp\Bigl(\;\!\dfrac{q}{p^{\;\!h}}+c\Bigr), N\biggr)\in 
\text{\rm E}_{_{\Omega_{{}_{\tiny\;\! 0}}}}.
\hfill
$
\\[2ex]
\indent
{\bf Property 5.2.}
\vspace{0.75ex}
{\it
Suppose 
$\bigl((p, M), (h,q,N)\bigr)\in \text{\rm B}_{_{\Xi^{\;\!\prime}}},\ 
\gamma_1^{},\gamma_2^{}\in\R,\ \varphi\in C^1T^{\;\!\prime},$ 
the set 
$\Omega_{_0}\subset\Xi^{\;\!\prime}$ such that   
$p(t,x)\ne 0$ for all $(t,x)\in\Omega_{_0}$ and
$p(t,x)= 0$ for all $(t,x)\in {\sf C}_{_{\Xi^{\;\!\prime}}}\Omega_{_0}.$
Then
\\[1.5ex]
\mbox{}\hfill                           
$
\displaystyle
\biggl(p^{{}^{\scriptsize \gamma_1^{}}}
\exp\Bigl(\gamma_2^{}\Bigl(\;\!\dfrac{q}{p^{\;\!h}}+\varphi\Bigr)\Bigr), L\biggr)\in 
\text{\rm J}_{_{\Omega_{{}_{\tiny\;\! 0}}}}
\hfill
$
\\[1.75ex]
if and only if  the cofactors $M,\, N,\, L$ such that
\\[1.5ex]
\mbox{}\hfill                           
$
\displaystyle
L(t,x)=\gamma_1^{}\;\!M(t,x)+\gamma_2^{}\;\!
\bigl(N(t,x)+{\sf D}\;\!\varphi(t)\bigr)
$
\ for all}
$(t,x)\in \Xi^{\;\!\prime}.
\hfill
$
\\[1.5ex]
\indent
{\sl Proof.}
Using Property 5.1, by Property 3.9, we get the statement of Property 5.2. $\k$
\vspace{1ex}

{\bf Corollary 5.2.}
\vspace{0.75ex}
{\it
Suppose 
$\bigl((p, M), (h,q,N)\bigr)\in \text{\rm B}_{_{\Xi^{\;\!\prime}}},\ 
c, \gamma_1^{},\gamma_2^{}\in\R,$ 
the set 
$\Omega_{_0}\subset\Xi^{\;\!\prime}$ such that  
$p(t,x)\ne 0$ for all $(t,x)\in\Omega_{_0}$ and 
$p(t,x)= 0$ for all $(t,x)\in {\sf C}_{_{\Xi^{\;\!\prime}}}\Omega_{_0}.$
Then 
\\[1.5ex]
\mbox{}\hfill                           
$
\displaystyle
\biggl(p^{\gamma_1^{}}
\exp\Bigl(\gamma_2^{}\Bigl(\;\!\dfrac{q}{p^{\;\!h}}+c\Bigr)\Bigr), L\biggr)\in 
\text{\rm J}_{_{\Omega_{{}_{\tiny\;\! 0}}}}
\hfill
$
\\[2ex]
if and only if the cofactors $M,\, N,\, L$ such that
\\[1.5ex]
\mbox{}\hfill                           
$
\displaystyle
L(t,x)=\gamma_1^{}\;\!M(t,x)+\gamma_2^{}\;\!N(t,x)
$
\ for all} 
$(t,x)\in \Xi^{\;\!\prime}.
\hfill
$
\\[2.25ex]
\indent
{\bf Property 5.3.}
\vspace{0.75ex}
{\it
Suppose 
$\!\bigl((p_j^{}, M_j^{}), (h_j^{},q_j^{},N)\bigr)\!\in\! \text{\rm B}_{_{\Xi^{\;\!\prime}}},\, 
\lambda_j^{},\gamma_j^{}\!\in\!\R,\, j\!=\!1,\ldots, m,\, \sum\limits_{j=1}^{m}|\lambda_j^{}|\!\ne\! 0,$ 
a function $\varphi\in C^1T^{\;\!\prime},$ 
the set 
\vspace{0.75ex}
$\Omega_{_0}\subset\Xi^{\;\!\prime}$ such that   
$\prod\limits_{j=1}^{m}p_j^{}(t,x)\ne 0$ for all $(t,x)\in\Omega_{_0}$ and
$\prod\limits_{j=1}^{m}p_j^{}(t,x)= 0$
for all $(t,x)\in {\sf C}_{_{\Xi^{\;\!\prime}}}\Omega_{_0},\ 
p_j^{\gamma_j^{}}\in C^1\Omega_{_0}, \ j=1,\ldots, m.$
Then we have
\\[1ex]
\mbox{}\hfill                           
$
\displaystyle
\biggl(\ \prod\limits_{j=1}^{m}p_j^{\gamma_j^{}}
\sum\limits_{j=1}^{m}\lambda_j^{}
\exp\biggl(\;\!\dfrac{q_j^{}}{p_j^{\;\!h_j^{}}}+\varphi\biggr), L\biggr)\in 
\text{\rm J}_{_{\Omega_{{}_{\tiny\;\! 0}}}}
\hfill
$
\\[1.5ex]
if and only if the cofactors $M_1^{},\ldots, M_{m}^{},\, N,\, L$ such that
\\[1.5ex]
\mbox{}\hfill                           
$
\displaystyle
L(t,x)={\sf D}\;\!\varphi(t)+N(t,x)+
\sum\limits_{j=1}^{m}\gamma_j^{}\;\!M_j^{}(t,x)
$
\ for all} 
$(t,x)\in \Xi^{\;\!\prime}.
\hfill
$
\\[1.25ex]
\indent
{\sl Proof}.
Using Properties 5.1 and 3.4 in consecutive order, we obtain
\\[1.5ex]
\mbox{}\hfill                           
$
\displaystyle
\biggl(
\exp\biggl(\;\!\dfrac{q_j^{}}{p_j^{\;\!h_j^{}}}+\varphi\biggr),\, 
N+{\sf D}\;\!\varphi\biggr)\in 
\text{\rm E}_{_{\Omega_{{}_{\tiny\;\! 0}}}},\ j=1,\ldots, m,
\ \Longrightarrow 
\hfill                           
$
\\[1.5ex]
\mbox{}\hfill                           
$
\displaystyle
\Longrightarrow \
\biggl(\ \sum\limits_{j=1}^{m}\lambda_j^{}
\exp\biggl(\;\!\dfrac{q_j^{}}{p_j^{\;\!h_j^{}}}+\varphi\biggr),\, 
N+{\sf D}\;\!\varphi\biggr)\in 
\text{\rm J}_{_{\Omega_{{}_{\tiny\;\! 0}}}}.
\hfill
$
\\[1.5ex]
\indent
Now, taking into account that
\vspace{0.75ex} 
$\bigl(p_j^{}, M_j^{}\bigr)\in \text{\rm A}_{_{\Xi^{\;\!\prime}}},\ j=1,\ldots, m,$ 
by Property 1.9, it follows that the statement of Property 5.3 is true. $\k$
\vspace{1ex}

{\bf Property 5.4.}
\vspace{0.75ex}
{\it
Suppose 
$\bigl((p_j^{}, M_j^{}), (h_j^{},q_j^{},N)\bigr)\in \text{\rm B}_{_{\Xi^{\;\!\prime}}},\ 
c_j^{}, \lambda_j^{}, \gamma_j^{}\in\R,\ j=1,\ldots, m,$ 
$\sum\limits_{j=1}^{m}|\lambda_j^{}|\ne 0,$ 
the set 
\vspace{0.75ex}
$\Omega_{_0}\subset\Xi^{\;\!\prime}$ such that   
$\prod\limits_{j=1}^{m}p_j^{}(t,x)\!\ne\! 0$ for all $(t,x)\in\Omega_{_0}$ and
$\prod\limits_{j=1}^{m}p_j^{}(t,x)\!=\! 0$
for all $(t,x)\in {\sf C}_{_{\Xi^{\;\!\prime}}}\Omega_{_0},\ 
p_j^{\gamma_j^{}}\in C^1\Omega_{_0}, \ j=1,\ldots, m.$
Then we claim that
\\[1.25ex]
\mbox{}\hfill                           
$
\displaystyle
\biggl(\ \prod\limits_{j=1}^{m}p_j^{\gamma_j^{}}
\sum\limits_{j=1}^{m}\lambda_j^{}
\exp\biggl(\;\!\dfrac{q_j^{}}{p_j^{\;\! h_j^{}}}+c_j^{}\biggr), L\biggr)\in 
\text{\rm J}_{_{\Omega_{{}_{\tiny\;\! 0}}}}
\hfill
$
\\[1.5ex]
if and only if the cofactors $M_1^{},\ldots, M_{m}^{},\, N,\, L$ such that 
\\[1.5ex]
\mbox{}\hfill                           
$
\displaystyle
L(t,x)=N(t,x)+
\sum\limits_{j=1}^{m}\gamma_j^{}\;\!M_j^{}(t,x)
$
\ for all} 
$(t,x)\in \Xi^{\;\!\prime}.
\hfill
$
\\[1.25ex]
\indent
{\sl Proof}.
Using Corollary 5.1 and Property 3.4 in consecutive order, we obtain
\\[1.25ex]
\mbox{}\hfill                           
$
\displaystyle
\biggl(
\exp\biggl(\;\!\dfrac{q_j^{}}{p_j^{\;\!h_j^{}}}+c_j^{}\biggr),\, N\biggr)\in 
\text{\rm E}_{_{\Omega_{{}_{\tiny\;\! 0}}}},\ j=1,\ldots, m,
\ \Longrightarrow \
\biggl(\ \sum\limits_{j=1}^{m}\lambda_j^{}
\exp\biggl(\;\!\dfrac{q_j^{}}{p_j^{\;\!h_j^{}}}+c_j^{}\biggr),\, N\biggr)\in 
\text{\rm J}_{_{\Omega_{{}_{\tiny\;\! 0}}}}.
\hfill
$
\\[1.25ex]
\indent
Finally, since
\vspace{0.75ex} 
$\bigl(p_j^{}, M_j^{}\bigr)\in \text{\rm A}_{_{\Xi^{\;\!\prime}}},\ j=1,\ldots, m,$ 
by Property 1.9, we see that the statement of Property 5.4 is true. $\k$
\vspace{0.75ex} 

{\bf Property 5.5.}
\vspace{0.75ex}
{\it
Suppose 
$\bigl((p_j^{}, M_j^{}), (h_j^{},q_j^{},\rho_{j}^{}N_0^{})\bigr)\in 
\text{\rm B}_{_{\Xi^{\;\!\prime}}},\ 
\rho_j^{}\in\R\backslash\{0\},\ 
\lambda_j^{},\gamma_j^{}\in\R,
\linebreak 
j=1,\ldots, m,\ \sum\limits_{j=1}^{m}|\lambda_j^{}|\ne 0,\ 
\varphi\in C^1T^{\;\!\prime},$ 
the set 
\vspace{0.75ex}
$\Omega_{_0}\subset\Xi^{\;\!\prime}$ such that   
$\prod\limits_{j=1}^{m}p_j^{}(t,x)\ne 0$ for all $(t,x)\in\Omega_{_0}$ and
$\prod\limits_{j=1}^{m}p_j^{}(t,x)= 0$
for all $(t,x)\in {\sf C}_{_{\Xi^{\;\!\prime}}}\Omega_{_0},\ 
p_j^{\gamma_j^{}}\in C^1\Omega_{_0}, \ j=1,\ldots, m.$
Then
\\[1.25ex]
\mbox{}\hfill                           
$
\displaystyle
\biggl(\ \prod\limits_{j=1}^{m}p_j^{\gamma_j^{}}
\sum\limits_{j=1}^{m}\lambda_j^{}
\exp\biggl(\;\!\dfrac{q_j^{}}{\rho_j^{}\;\!p_j^{\;\!h_j^{}}}+\varphi\biggr), L\biggr)\in 
\text{\rm J}_{_{\Omega_{{}_{\tiny\;\! 0}}}}
\hfill
$
\\[1.25ex]
if and only if the cofactor
\\[1.5ex]
\mbox{}\hfill                           
$
\displaystyle
L(t,x)={\sf D}\;\!\varphi(t)+N_0^{}(t,x)+
\sum\limits_{j=1}^{m}\gamma_j^{}\;\!M_j^{}(t,x)
$
\ for all}
$(t,x)\in \Xi^{\;\!\prime}.
\hfill
$
\\[1.5ex]
\indent
{\sl Proof}.
Using Theorem 5.3, Properties 3.5 (under $m=1,$ $\lambda_1^{}=1),$
3.2, and 3.4 in consecutive order, we obtain
\\[1.5ex]
\mbox{}\hfill                           
$
\displaystyle
\biggl(
\exp\;\!\dfrac{q_j^{}}{p_j^{\;\!h_j^{}}}\,,\, \rho_{j}^{}\;\!N_0^{}\biggr)\in 
\text{\rm E}_{_{\Omega_{{}_{\tiny\;\! 0}}}}
\ \Longrightarrow \
\biggl(
\exp\;\!\dfrac{q_j^{}}{ \rho_{j}^{}\;\! p_j^{\;\!h_j^{}}}\,,\, N_0^{}\biggr)\in 
\text{\rm E}_{_{\Omega_{{}_{\tiny\;\! 0}}}}
\ \Longrightarrow \
\hfill                           
$
\\[2ex]
\mbox{}\hfill                           
$
\displaystyle
\Longrightarrow \
\biggl(
\exp\biggl(\;\!\dfrac{q_j^{}}{ \rho_{j}^{}\;\! p_j^{\;\!h_j^{}}}+\varphi\biggr),\, 
N_0^{}+{\sf D}\;\!\varphi\biggr)\in 
\text{\rm E}_{_{\Omega_{{}_{\tiny\;\! 0}}}},\ j=1,\ldots, m,
\ \Longrightarrow 
\hfill                           
$
\\[2ex]
\mbox{}\hfill                           
$
\displaystyle
\Longrightarrow \
\biggl(\ \sum\limits_{j=1}^{m}\lambda_j^{}
\exp\biggl(\;\!\dfrac{q_j^{}}{ \rho_{j}^{}\;\! p_j^{\;\!h_j^{}}}+\varphi\biggr),\, 
N_0^{}+{\sf D}\;\!\varphi\biggr)\in 
\text{\rm J}_{_{\Omega_{{}_{\tiny\;\! 0}}}}.
\hfill
$
\\[1.5ex]
\indent
Further, taking into account that
\vspace{0.75ex} 
$\bigl(p_j^{}, M_j^{}\bigr)\in \text{\rm A}_{_{\Xi^{\;\!\prime}}},\ j=1,\ldots, m,$ 
from Property 1.9 it follows that the statement of Property 5.5 is true. $\k$
\vspace{1.25ex} 

{\bf Property 5.6.}
\vspace{0.75ex}
{\it
Suppose 
$\bigl((p_j^{}, M_j^{}), (h_j^{},q_j^{},\rho_{j}^{}N_0^{})\bigr)\in 
\text{\rm B}_{_{\Xi^{\;\!\prime}}},\ 
\rho_j^{}\in\R\backslash\{0\},\ 
c_j^{}, \lambda_j^{},\gamma_j^{}\in\R,$  
$j=1,\ldots, m,\
\sum\limits_{j=1}^{m}|\lambda_j^{}|\ne 0,$ 
the set 
\vspace{0.75ex}
$\Omega_{_0}\subset\Xi^{\;\!\prime}$ such that   
$\prod\limits_{j=1}^{m}p_j^{}(t,x)\ne 0$ for all $(t,x)\in\Omega_{_0}$ and
$\prod\limits_{j=1}^{m}p_j^{}(t,x)= 0$
for all $(t,x)\in {\sf C}_{_{\Xi^{\;\!\prime}}}\Omega_{_0},\ 
p_j^{\gamma_j^{}}\in C^1\Omega_{_0}, \ j=1,\ldots, m.$
Then we have
\\[1.25ex]
\mbox{}\hfill                           
$
\displaystyle
\biggl(\ \prod\limits_{j=1}^{m}p_j^{\gamma_j^{}}
\sum\limits_{j=1}^{m}\lambda_j^{}
\exp\biggl(\;\!\dfrac{q_j^{}}{\rho_j^{}\;\!p_j^{\;\!h_j^{}}}+c_j^{}\biggr), L\biggr)\in 
\text{\rm J}_{_{\Omega_{{}_{\tiny\;\! 0}}}}
\hfill
$
\\[1.5ex]
if and only if the cofactor
\\[1.5ex]
\mbox{}\hfill                           
$
\displaystyle
L(t,x)=N_0^{}(t,x)+
\sum\limits_{j=1}^{m}\gamma_j^{}\;\!M_j^{}(t,x)
$
\ for all} $(t,x)\in \Xi^{\;\!\prime}.
\hfill
$
\\[1.5ex]
\indent
{\sl Proof}\;\!
is similar to the proof of Property 5.5 and has only one difference.
Namely, Property 3.3 is used instead of Property 3.2. $\k$
\vspace{1.25ex} 

{\bf Property 5.7.}\!
\vspace{0.75ex}
{\it
Suppose 
$\!\bigl((p_j^{}, M_j^{}), (h_j^{},q_j^{},N_j^{})\bigr)\!\in\! \text{\rm B}_{_{\Xi^{\;\!\prime}}},\, 
\gamma_j^{},\xi_j^{}\!\in\!\R,\, \varphi_j^{}\!\in\! C^1T^{\;\!\prime},\, 
j\!=\!1,\ldots, m,\!$ 
the set
\vspace{0.75ex}
$\Omega_{_0}\subset\Xi^{\;\!\prime}$ such that  
$\prod\limits_{j=1}^{m}p_j^{}(t,x)\ne 0$ for all $(t,x)\in\Omega_{_0}$ and
$\prod\limits_{j=1}^{m}p_j^{}(t,x)= 0$
for all $(t,x)\in {\sf C}_{_{\Xi^{\;\!\prime}}}\Omega_{_0},$ 
$p_j^{\gamma_j^{}}\in C^1\Omega_{_0}, \ j=1,\ldots, m.$
Then we claim that
\\[1.5ex]
\mbox{}\hfill                           
$
\displaystyle
\biggl(\ \prod\limits_{j=1}^{m}p_j^{\gamma_j^{}}
\exp\sum\limits_{j=1}^{m}\xi_j^{}
\biggl(\;\!\dfrac{q_j^{}}{p_j^{\;\! h_j^{}}}+\varphi_j^{}\biggr), L\biggr)\in 
\text{\rm J}_{_{\Omega_{{}_{\tiny\;\! 0}}}}
\hfill
$
\\[2ex]
if and only if the cofactor
\\[1.75ex]
\mbox{}\hfill                           
$
\displaystyle
L(t,x)=
\sum\limits_{j=1}^{m}\bigl(
\gamma_j^{}\;\!M_j^{}(t,x)+
\xi_j^{}\;\!\bigl(N_j^{}(t,x)+{\sf D}\;\!\varphi_j^{}(t)\bigr)\bigr)
$
\ for all} $(t,x)\in \Xi^{\;\!\prime}.
\hfill
$
\\[1.5ex]
\indent
{\sl Proof}.
By Property 5.1, we have
\\[1.5ex]
\mbox{}\hfill                           
$
\displaystyle
\biggl(
\exp\biggl(\;\!\dfrac{q_j^{}}{p_j^{\;\!h_j^{}}}+\varphi_j^{}\biggr),\, 
N_j^{}+{\sf D}\;\!\varphi_j^{}\biggr)\in 
\text{\rm E}_{_{\Omega_{{}_{\tiny\;\! 0}}}},
\ \ 
j=1,\ldots, m.
\hfill                           
$
\\[1.5ex]
\indent
Taking into account that
\vspace{0.75ex} 
$\bigl(p_j^{}, M_j^{}\bigr)\in \text{\rm A}_{_{\Xi^{\;\!\prime}}},\ j=1,\ldots, m,$ 
by Property 3.9, it follows that the statement of Property 5.7 is true. $\k$
\vspace{0.75ex} 

By definition 5.1, 
the multiplicity of polynomial partial integral $p$ depends on
the quantity of natural numbers $h$ and on the 
functions $q\in \text{\rm P}_{_{\!\Xi^{\;\!\prime}}}$ 
\vspace{0.25ex} 
(corresponding to these natural numbers) such that 
the identity (5.1) holds.
Moreover, taking into account Theorem 5.3, we obtain 
the multiplicity of polynomial partial integral $p$
\vspace{0.25ex} 
is determined by the quantity of exponential partial integrals
$\exp\;\!\dfrac{q}{p^{\;\!h}}\;\!.$
\vspace{0.75ex} 

{\bf Definition 5.2.}
{\it
A polynomial partial integral  $p$ on the domain $\Xi^{\;\!\prime}$ of system {\rm (0.1)} 
is said to be  \textit{\textbf{multiple with multiplicity}}
$
\varkappa =1 + \sum\limits_{\xi=1}^{\varepsilon} \delta_{\xi}^{}
$
if there exist natural numbers
$h_{\xi}^{},$ 
\linebreak
$\xi =1,\ldots, \varepsilon,$ and functions
\vspace{0.75ex} 
$
q_{_{\scriptstyle h_\xi^{}f_\xi^{}}}\in \text{\rm P}_{_{\!\Xi^{\;\!\prime}}},\ 
f_{\xi}^{}=1,\ldots,\delta_{\xi}^{},\ \xi=1,\ldots,\varepsilon,$
that correspond to these numbers and relatively prime to the function $p,$
such that the identities hold 
\\[2.25ex]
\mbox{}\hfill                    % (5.3)
$
\displaystyle
{\frak d}\, \dfrac{q_{_{\scriptstyle h_\xi^{}  f_\xi^{}} }(t,x)}{\displaystyle  p^{\;\!h_\xi^{}} (t,x)}=
N_{h_\xi^{}  f_\xi^{}}^{} (t,x)
$
\ for all 
$(t,x)\in \Omega_{_0},
\quad 
f_{\xi}^{}=1,\ldots, \delta_{\xi}^{}, \ \   
\xi=1,\ldots,\varepsilon,
$
\hfill {\rm (5.3)}
\\[2.5ex]
where 
the functions  
\vspace{0.75ex}
$N_{h_\xi^{}  f_\xi^{}}^{}\in\text{\rm P}_{_{\!\Xi^{\;\!\prime}}}$ and
have the degrees 
$
\deg_{\;\!x}^{}  N_{h_\xi^{}  f_\xi^{}}^{}\leq d-1, \ 
f_{\xi}^{}=1,\ldots, \delta_{\xi}^{}, 
\linebreak  
\xi=1,\ldots,\varepsilon,
$
the set 
\vspace{0.5ex}
$\Omega_{_0}\subset\Xi^{\;\!\prime}$ such that   
$p(t,x)\ne 0$ for all $(t,x)\in\Omega_{_0}$ and
$p(t,x)= 0$ for all $(t,x)\in {\sf C}_{_{\Xi^{\;\!\prime}}}\Omega_{_0}.$
}
\vspace{1.5ex}

Using Definition 5.2, we get
\vspace{1.25ex}

{\bf Proposition 5.1.}\!
\vspace{0.75ex}
{\it
If  
$\!\Bigl((p, M), 
\Bigl(h_\xi^{},q_{_{\scriptstyle h_\xi^{}  f_\xi^{}}}, N_{h_\xi^{}  f_\xi^{}}^{}\Bigr)\Bigr)
\!\in\! \text{\rm B}_{_{\Xi^{\;\!\prime}}},\,
f_{\xi}^{}\!=\!1,\ldots, \delta_{\xi}^{},\, \xi\!=\!1,\ldots,\varepsilon,\!$
then the polynomial partial integral $p$ with cofactor $M$
\vspace{0.35ex}
on the domain $\Xi^{\;\!\prime}$ of system {\rm (0.1)} is 
multiple with multiplicity
$
\varkappa =1 + \sum\limits_{\xi=1}^{\varepsilon} \delta_{\xi}^{}.
$
}
\vspace{0.75ex}

We stress that 
\vspace{0.5ex}
in Definition 5.2 and in Proposition 5.1 there is a possibility, 
when to one number $\!h_\xi^{}\!\in\!\N\!$ corresponds $\!\delta_{\xi}^{}\!$ functions
$\!q_{_{\scriptstyle h_{\xi}^{}  f_\xi^{}}}\!\!\in\! \text{\rm P}_{_{\!\Xi^{\;\!\prime}}},
\, f_{\xi}^{}\!=\!1,\ldots, \delta_{\xi}^{}.\!$
Also it is not excluded that
\\[2.25ex]
\mbox{}\hfill
$
q_{_{\scriptstyle h_\xi^{}  f_\xi^{}}}(t,x)=
q_{_{\scriptstyle h_\zeta^{}  f_\zeta^{}}}(t,x)
$
for all $(t,x)\in \Xi^{\;\!\prime}
$
\ under \ $\zeta\ne \xi,\ \ \xi,\zeta\in\{1,\ldots,\varepsilon\}.
\hfill
$
\\[2.5ex]
\indent
{\bf Property 5.8.}
\vspace{1ex}
{\it
Suppose 
$\Bigl((p, M), 
\Bigl(h_\xi^{},q_{_{\scriptstyle h_\xi^{}  f_\xi^{}}}, N_{h_\xi^{}  f_\xi^{}}^{}\Bigr)\Bigr)
\in \text{\rm B}_{_{\Xi^{\;\!\prime}}},\
\gamma,\gamma_{_{\scriptstyle h_\xi^{}  f_\xi^{}}}\in\R,\ 
\varphi_{_{\scriptstyle h_\xi^{}  f_\xi^{}}}\in C^1T^{\;\!\prime},$ 
$f_{\xi}^{}=1,\ldots, \delta_{\xi}^{},\ \xi=1,\ldots,\varepsilon,$
\vspace{1ex}
the set
$\Omega_{_0}\subset\Xi^{\;\!\prime}$ such that   
$p(t,x)\ne 0$ for all $(t,x)\in\Omega_{_0}$ and
$p(t,x)= 0$ for all $(t,x)\in {\sf C}_{_{\Xi^{\;\!\prime}}}\Omega_{_0},\ 
p^{\gamma}\in C^1\Omega_{_0}.$
Then
\\[1.5ex]
\mbox{}\hfill                           
$
\displaystyle
\biggl(\ p^{\gamma}
\exp\sum\limits_{\xi=1}^{\varepsilon}\sum\limits_{f_{\xi}^{}=1}^{\delta_{\xi}^{}}
\biggl(\;\!\gamma_{_{\scriptstyle h_\xi^{}  f_\xi^{}}}\biggl(\,
\dfrac{q_{_{\scriptstyle h_\xi^{}  f_\xi^{}} }}{\displaystyle  p^{\;\!h_\xi^{}}}+
\varphi_{_{\scriptstyle h_\xi^{}  f_\xi^{}}}\biggr)\biggr),\, L\biggr)\in 
\text{\rm J}_{_{\Omega_{{}_{\tiny\;\! 0}}}}
\hfill
$
\\[2ex]
if and only if the cofactor
\\[1.75ex]
\mbox{}\hfill                           
$
\displaystyle
L(t,x)=\gamma\;\!M(t,x)+
\sum\limits_{\xi=1}^{\varepsilon}\sum\limits_{f_{\xi}^{}=1}^{\delta_{\xi}^{}}
\biggl(\gamma_{_{\scriptstyle h_\xi^{}  f_\xi^{}}}\Bigl(
N_{h_\xi^{}  f_\xi^{}}^{} (t,x)+{\sf D}\;\!\varphi_{_{\scriptstyle h_\xi^{}  f_\xi^{}}}(t)\Bigr)\biggr)
$
\ for all} $(t,x)\in \Xi^{\;\!\prime}.
\hfill
$
\\[1.5ex]
\indent
{\sl Proof}.
It follows from Property 5.1 that
\\[1.5ex]
\mbox{}\hfill                           
$
\displaystyle
\biggl(
\exp\biggl(\;\!
\dfrac{q_{_{\scriptstyle h_\xi^{}  f_\xi^{}} }}{\displaystyle  p^{\;\!h_\xi^{}}}+
\varphi_{_{\scriptstyle h_\xi^{}  f_\xi^{}}}\biggr),\, 
N_{h_\xi^{}  f_\xi^{}}^{}+{\sf D}\;\!\varphi_{_{\scriptstyle h_\xi^{}  f_\xi^{}}}
\biggr)\in 
\text{\rm E}_{_{\Omega_{{}_{\tiny\;\! 0}}}},
\ \ 
f_{\xi}^{}=1,\ldots, \delta_{\xi}^{},\ \xi=1,\ldots,\varepsilon.
\hfill                           
$
\\[1.5ex]
\indent
Since
\vspace{1ex}  
$\bigl(p, M\bigr)\in \text{\rm A}_{_{\Xi^{\;\!\prime}}},$ 
by Property 3.9, we see that the statement of Property 5.8 is true.$\k$

Property 5.8 is the analog of Property 5.7 in the case of 
\vspace{0.35ex} 
$\varkappa\!$-multiple partial integral 
$\bigl(p, M\bigr)\in \text{\rm A}_{_{\Xi^{\;\!\prime}}}.$
\vspace{0.5ex} 
In similarly way, we can obtain analogues of Properties 5.3, 5.4, 5.5, and 5.6 
in the case of $\varkappa\!$-multiple partial integral
$\bigl(p, M\bigr)\in \text{\rm A}_{_{\Xi^{\;\!\prime}}}\,.$
\vspace{1.25ex} 

{\bf Property 5.9.}
{\it
Let 
$k\in\N,\ \varphi\in C^1T^{\;\!\prime}.$ Then we claim that
\\[1.5ex]
\mbox{}\hfill                           
$
\bigl(p, \varphi+p^kM_{0}^{}\bigr)\in \text{\rm A}_{_{\Xi^{\;\!\prime}}}
\iff
\bigl(\bigl(p, \varphi+p^kM_{0}^{}\bigr), \bigl(k,q,-k\;\!qM_0^{}\bigr)\bigr)
\in \text{\rm B}_{_{\Xi^{\;\!\prime}}},
\hfill
$
\\[1ex]
where
\\[1ex]
\mbox{}\hfill                   % (5.4)        
$
\displaystyle
q\colon t\to\ 
\exp\biggl(\, k\int\limits_{t_0^{}}^{t}\varphi(\tau)\;\!d\tau\biggr)
$
\ for all $t\in T^{\;\!\prime},
$
\hfill {\rm(5.4)}
\\[1.5ex]
$t_0^{}$ is arbitrary fixed point from the domain $T^{\;\!\prime}.$
}
\vspace{0.5ex}

{\sl Proof}.
By Theorem 2.1, we have
\\[1.5ex]
\mbox{}\hfill                           
$
\displaystyle
{\frak d}\;\!p(t,x)=p(t,x)\;\!\bigl(\varphi(t)+p^{\;\!k}(t,x)\;\!M_0^{}(t,x)\bigr)
$
\ for all 
$(t,x)\in \Xi^{\;\!\prime},
\hfill
$
\\[1.75ex]
where the function 
\vspace{0.75ex}
$M_0^{}\in \text{\rm P}_{_{\!\Xi^{\;\!\prime}}}$
and has the degree $\deg_{\;\!x}^{} M_0^{}\leq d-k\deg_{\;\!x}^{}p-1.$

Therefore the derivative by virtue of system (0.1)
\\[1.5ex]
\mbox{}\hfill                           
$
\displaystyle
{\frak d}\,\dfrac{\displaystyle
\exp\biggl(\, k\int\limits_{t_0^{}}^{t}\varphi(\tau)\;\!d\tau\biggr)}{p^{\;\!k}(t,x)}=
\hfill
$
\\[2ex]
\mbox{}\hfill                           
$
\displaystyle
=\dfrac{\displaystyle
p^{\;\!k}(t,x)\, {\sf D}\exp\biggl(\, k\int\limits_{t_0^{}}^{t}\varphi(\tau)\;\!d\tau\biggr) -
k\;\!p^{\;\!k-1}(t,x)\;\!{\frak d}\;\!p(t,x)\, \exp\biggl(\, k\int\limits_{t_0^{}}^{t}\varphi(\tau)\;\!d\tau\biggr)}
{p^{\;\!2k}(t,x)}=
\hfill
$
\\[2ex]
\mbox{}\hfill                           
$
\displaystyle
=
\dfrac{\displaystyle
k\;\!\varphi(t)\;\!p(t,x)\;\!\exp\biggl(\, k\int\limits_{t_0^{}}^{t}\varphi(\tau)\;\!d\tau\biggr) -
k\;\!p(t,x)\;\!\bigl(\varphi(t)+p^{\;\!k}(t,x)\;\!M_0^{}(t,x)\bigr)\;\!
\exp\biggl(\, k\int\limits_{t_0^{}}^{t}\varphi(\tau)\;\!d\tau\biggr)}
{p^{\;\!k+1}(t,x)}=
\hfill
$
\\[2ex]
\mbox{}\hfill                           
$
\displaystyle
={}-k\;\!M_0^{}(t,x)\;\!\exp\biggl(\, k\int\limits_{t_0^{}}^{t}\varphi(\tau)\;\!d\tau\biggr)
$
\ for all 
$(t,x)\in\Omega_{_{0}}^{},
\hfill
$
\\[2ex]
where the set
$\Omega_{_0}\subset\Xi^{\;\!\prime}$ such that  
$p(t,x)\ne 0$ for all $(t,x)\in\Omega_{_0}$ and
$p(t,x)= 0$ for all $(t,x)\in {\sf C}_{_{\Xi^{\;\!\prime}}}\Omega_{_0}.\ \k$
\vspace{1ex}

In particular, if $\varphi(t)=0$ for all $t\in T^{\;\!\prime},$ then we get 
\vspace{1ex}

{\bf Property 5.10.}
{\it
Let 
$k\in\N,\ c\in\R\backslash\{0\}.$ Then}
\\[1.5ex]
\mbox{}\hfill                           
$
\bigl(p, p^kM_{0}^{}\bigr)\in \text{\rm A}_{_{\Xi^{\;\!\prime}}}
\iff
\bigl(\bigl(p, p^kM_{0}^{}\bigr), \bigl(k, c,{}-k\;\!c\;\!M_0^{}\bigr)\bigr)
\in \text{\rm B}_{_{\Xi^{\;\!\prime}}}.
\hfill
$
\\[2ex]
\indent
{\bf Property 5.11.}
\vspace{0.75ex}
{\it
If  
$k\in\N,\ \bigl(p, p^kM_{0}^{}\bigr)\in \text{\rm A}_{_{\Xi^{\;\!\prime}}},$ 
then the polynomial partial integral $p$ on the domain $\Xi^{\;\!\prime}$ of system {\rm (0.1)} 
is $(k+1)\!\!$-multiple and 
\\[2ex]
\mbox{}\hfill                           
$
\bigl(\bigl(p, p^kM_{0}^{}\bigr), \bigl(l, c_l^{},{}-l\;\!c_l^{}\;\!p^{\;\!k-l}M_0^{}\bigr)\bigr)
\in \text{\rm B}_{_{\Xi^{\;\!\prime}}},
\ \ 
l=1,\ldots, k,
\hfill
$
\\[2ex]
where $c_l^{},\ l=1,\ldots,k,$ 
are arbitrary fixed nonzero real numbers.}
\vspace{0.75ex}

{\sl Доказательство.}
Since  
$\bigl(p, p^kM_{0}^{}\bigr)\in \text{\rm A}_{_{\Xi^{\;\!\prime}}},$ we see that
$\bigl(p, p^{\;\!l}M_{l}^{}\bigr)\in \text{\rm A}_{_{\Xi^{\;\!\prime}}},\ l=1,\ldots, k,$
where the functions
\\[1.5ex]
\mbox{}\hfill                           
$
M_{l}^{}(t,x)=p^{\;\!k-l}(t,x)\;\!M_0^{}(t,x)
$
for all 
$(t,x)\in \Xi^{\;\!\prime},
\quad 
l=1,\ldots, k.
\hfill
$
\\[2ex]
\indent
By Property 5.10, 
\vspace{1ex}
$
\bigl(\bigl(p, p^{\;\!l}M_{l}^{}\bigr), \bigl(l, c_l^{},{}-l\;\!c_l^{}\;\!M_l^{}\bigr)\bigr)
\in \text{\rm B}_{_{\Xi^{\;\!\prime}}},
\ l=1,\ldots, k,
$
where $c_l^{}\in\R\backslash\{0\},
\linebreak 
l=1,\ldots,k.$ 
Hence, 
\vspace{1.5ex} 
$
\bigl(\bigl(p, p^{\;\!k}M_{0}^{}\bigr), \bigl(l, c_l^{},{}-l\;\!c_l^{}\;\!p^{\;\!k-l}M_0^{}\bigr)\bigr)
\in \text{\rm B}_{_{\Xi^{\;\!\prime}}},\
c_l^{}\in\R\backslash\{0\},\ l=1,\ldots,k.\ \k$

{\bf Property 5.12.}
\vspace{0.5ex} 
{\it
Suppose 
$k\in\N,\ \varphi\in C^1T^{\;\!\prime},$ 
the function $q$ is given by formula {\rm (5.4)},
the set 
$\Omega_{_0}\subset\Xi^{\;\!\prime}$ such that  
$p(t,x)\ne 0$ for all $(t,x)\in\Omega_{_0}$ and
\vspace{0.35ex} 
$p(t,x)= 0$ for all $(t,x)\in {\sf C}_{_{\Xi^{\;\!\prime}}}\Omega_{_0}.$
Then we claim that}
\\[1.5ex]
\mbox{}\hfill                           
$
\bigl(p, \varphi+p^{\;\!k}M_{0}^{}\bigr)\in \text{\rm A}_{_{\Xi^{\;\!\prime}}}
\iff
\biggl(\exp\dfrac{q}{p^{\;\!k}}\,,\, {}-k\;\!q\;\!M_0^{}\biggr)\in
\text{\rm E}_{_{\Omega_{{}_{\tiny\;\! 0}}}}.
\hfill
$
\\[2ex]
\indent
{\sl Proof.}
Taking into account Theorem 5.3, from Property 5.9 it follows that 
the statement of Property 5.12 is true.  $\k $
\vspace{1ex}

{\bf Property 5.13.}
\vspace{0.5ex} 
{\it
Suppose 
$k\in\N,\ c\in\R\backslash\{0\},$ 
the set 
$\Omega_{_0}\subset\Xi^{\;\!\prime}$ such that  
$p(t,x)\ne 0$ for all $(t,x)\in\Omega_{_0}$ and
$p(t,x)= 0$ for all $(t,x)\in {\sf C}_{_{\Xi^{\;\!\prime}}}\Omega_{_0}.$
Then we have}
\\[1.25ex]
\mbox{}\hfill                           
$
\bigl(p, p^{\;\!k}M_{0}^{}\bigr)\in \text{\rm A}_{_{\Xi^{\;\!\prime}}}
\iff
\biggl(\exp\dfrac{c}{p^{\;\!k}}\,,\, {}-k\;\!c\;\!M_0^{}\biggr)\in
\text{\rm E}_{_{\Omega_{{}_{\tiny\;\! 0}}}}.
\hfill
$
\\[1.5ex]
\indent
{\sl Proof.}
Taking into account Theorem 5.3, from Property 5.10 it follows that 
the statement of Property 5.13 is true.  $\k $
\vspace{0.75ex}

{\bf Property 5.14.}
\vspace{0.5ex} 
{\it
Suppose 
$k\in\N,\ c_l^{}\in\R\backslash\{0\},\ l=1,\ldots, k,$ 
the set
$\Omega_{_0}\subset\Xi^{\;\!\prime}$ such that   
$p(t,x)\ne 0$ for all $(t,x)\in\Omega_{_0}$ and
$p(t,x)= 0$ for all $(t,x)\in {\sf C}_{_{\Xi^{\;\!\prime}}}\Omega_{_0}.$
\vspace{0.75ex} 
Then, since
$\bigl(p, p^{\;\!k}M_{0}^{}\bigr)\in \text{\rm A}_{_{\Xi^{\;\!\prime}}},$ it follows that}
\\[1.5ex]
\mbox{}\hfill                           
$
\biggl(\exp\dfrac{c_l^{}}{p^{\;\!l}}\,,\, {}-l\;\!c_{l}^{}\;\!p^{\;\!k-l}\;\!M_0^{}\biggr)\in
\text{\rm E}_{_{\Omega_{{}_{\tiny\;\! 0}}}},
\ \ l=1,\ldots, k.
\hfill
$
\\[1.75ex]
\indent
{\sl Proof.}
Taking into account Theorem 5.3, from Property 5.11 it follows that 
the statement of Property 5.14 is true.  $\k $
\\[5.25ex]
\centerline{
{\bf  6. Complex-valued polynomial partial integrals}
}
\\[1.5ex]
\indent
The set of functions that are polynomials in the variables $x_1^{},\ldots,x_n^{}$
with continuously differentiable 
complex-valued coefficients-functions in the variable $t$ on the domain $T$ is denoted by
$\text{Z}_{_{\Xi}}.$ 
\vspace{0.75ex}

{\bf Definition 6.1.}
\vspace{0.5ex}
{\it
A function $w\in \text{\rm Z}_{_{\Xi^{\;\!\prime}}}\!$ is said to be a
\textit{\textbf{complex-valued polynomial par\-tial integral on the domain}}
$\Xi^{\;\!\prime}$ of system {\rm (0.1)} if 
the derivative by virtye of system {\rm (0.1)}
\\[2ex]
\mbox{}\hfill                           % (6.1)
$
{\frak d}\;\!w(t,x)=w(t,x)\;\!W(t,x)
$
\ for all 
$(t,x)\in \Xi^{\;\!\prime}.
$
\hfill {\rm (6.1)}
\\[2ex]
And, the function $W$ is called the \textit{\textbf{cofactor}}
\vspace{1ex}
of complex-valued polynomial partial integral $\!w.$}

{\bf Remark 6.1.}
\vspace{0.5ex}
Since the function $w\in \text{Z}_{_{\Xi^{\;\!\prime}}},$
we see that from the identity (6.1) it follows that the cofactor 
$W\in \text{Z}_{_{\Xi^{\;\!\prime}}}.$
Note also that the case $W\in \text{P}_{_{\!\Xi^{\;\!\prime}}}$ is not excluded.
\vspace{0.5ex}

\newpage

After the introduction of the notion of complex-valued polynomial partial integral 
about a polynomial partial integral in the sense of Definition 2.1 we can speak (if necessary) 
as about real polynomial partial integral.
\vspace{0.25ex}

By $\text{H}_{_{\Xi^{\;\!\prime}}}$ 
\vspace{0.25ex}
denote the set of complex-valued polynomial partial integrals on the domain 
$\Xi^{\;\!\prime}$ of system (0.1).
The phrase "the function $w$ is a complex-valued polynomial partial integral with cofactor $W$ on
the domain $\Xi^{\;\!\prime}$ of system (0.1)" is denoted by $(w, W)\in \text{H}_{_{\Xi^{\;\!\prime}}}.$
\vspace{0.75ex}

{\bf Property 6.1.} 
{\it
Let $\eta\in\C\backslash\{0\}.$ Then}
\\[1.5ex]
\mbox{}\hfill
$
(w, W)\in \text{\rm H}_{_{\Xi^{\;\!\prime}}}
\iff 
(\eta\;\! w\;\!,\;\! W)\in \text{\rm H}_{_{\Xi^{\;\!\prime}}}.
\hfill
$
\\[1.25ex]
\indent
{\sl Proof}\;\!
is based on Definition 6.1 and the identity 
\\[1.5ex]
\mbox{}\hfill                          
$
{\frak d}\;\!\bigl(\eta\;\!w(t,x)\bigr)=\eta\,{\frak d}\;\!w(t,x)
$
for all 
$(t,x)\in \Xi^{\;\!\prime},$
\ for any $\eta\in\C\backslash\{0\}.
\ \k
\hfill
$
\\[2.25ex]
\indent
{\bf Property 6.2.} 
{\it
If $\eta\in\C,\ w\in \text{\rm Z}_{_{\Xi^{\;\!\prime}}},$ and the identity holds
\\[1.5ex]
\mbox{}\hfill
$
{\frak d}\;\!w(t,x)=\bigl(w(t,x)+\eta\bigr)\;\!W(t,x)
$
\ for all 
$(t,x)\in \Xi^{\;\!\prime},
\hfill
$
\\[1.5ex]
then $(w+\eta, W)\in \text{\rm H}_{_{\Xi^{\;\!\prime}}}.$}
\vspace{0.5ex}

{\sl Indeed},
the derivative by virtue of system (0.1)
\\[1.5ex]
\mbox{}\hfill                          
$
{\frak d}\;\!\bigl(w(t,x)+\eta\bigr)=
{\frak d}\;\!w(t,x)=\bigl(w(t,x)+\eta\bigr)\;\!W(t,x)
$
for all 
$(t,x)\in \Xi^{\;\!\prime}.
\hfill
$
\\[1.5ex]
\indent
Now, by Definition 6.1, we get
$(w+\eta, W)\in \text{\rm H}_{_{\Xi^{\;\!\prime}}}.\ \k$
\vspace{1.25ex}

{\bf Property 6.3.} 
{\it
Let $k\in\N.$ Then}
\\[1.5ex]
\mbox{}\hfill
$
(w, W)\in \text{\rm H}_{_{\Xi^{\;\!\prime}}}
\iff 
(w^k,\;\! k\;\!W)\in \text{\rm H}_{_{\Xi^{\;\!\prime}}}.
\hfill
$
\\[1.5ex]
\indent
{\sl Proof}\;\!
follows from Definition 6.1 and the identity
\\[1.5ex]
\mbox{}\hfill                          
$
{\frak d}\;\!w^k(t,x)=k\;\!w^{k-1}(t,x)\,{\frak d}\;\!w(t,x)
$
for all 
$
(t,x)\in \Xi^{\;\!\prime}.\ \k
\hfill
$
\\[2ex]
\indent
{\bf Property 6.4.} 
{\it
If $(w_j^{}, W)\in \text{\rm H}_{_{\Xi^{\;\!\prime}}},\ 
\eta_j^{}\in\C\backslash\{0\},\ j=1,\ldots,m,$ then}
\\[1.5ex]
\mbox{}\hfill
$
\displaystyle
\biggl(\,\sum\limits_{j=1}^{m} \eta_j^{}\;\!w_j^{}\;\!,\;\! W\biggl)\;\! 
\in \text{\rm H}_{_{\Xi^{\;\!\prime}}}.
\hfill
$
\\[1.5ex]
\indent
{\sl Indeed}, the derivative by virtue of system (0.1)
\\[1.5ex]
\mbox{}\hfill                         
$
\displaystyle
{\frak d}\;\!\sum\limits_{j=1}^{m} \eta_j^{}\;\!w_j^{}(t,x)  = 
\sum\limits_{j=1}^{m}
\eta_j^{}\, {\frak d}\;\!w_j^{}(t,x)=
\sum\limits_{j=1}^{m}
\eta_j^{}\;\!w_j^{}(t,x)\;\!W(t,x)
$
for all 
$
(t,x)\in \Xi^{\;\!\prime}.
\hfill
$
\\[1.5ex]
\indent
Finally, by Definition 6.1, we obtain  
$\biggl(\,\sum\limits_{j=1}^{m}\! \eta_j^{}\;\!w_j^{}\;\!,\;\! W\biggl) 
\in \text{\rm H}_{_{\Xi^{\;\!\prime}}}.\ \k$
\vspace{1.25ex}

{\bf Property 6.5.}
{\it
Let $w_j^{}\in\text{\rm Z}_{_{\Xi^{\;\!\prime}}},\ j=1,\ldots,m.$ Then
\\[1.5ex]
\mbox{}\hfill                         
$
\displaystyle
\biggl(\,\prod\limits_{j=1}^{m} w_j^{}\;\!,\;\! W\biggl)\;\! 
\in \text{\rm H}_{_{\Xi^{\;\!\prime}}}
\iff
\bigl( w_j^{}\;\!,\;\! W_j^{}\bigl)\;\! \in \text{\rm H}_{_{\Xi^{\;\!\prime}}},
\quad
j=1,\ldots,m.
\hfill
$
\\[1.5ex]
Moreover, the cofactors $W,\, W_1^{},\ldots,W_m^{}$ such that 
\\[1.5ex]
\mbox{}\hfill
$
\displaystyle
W(t,x)=\sum\limits_{j=1}^{m} W_j^{}(t,x)
$
 for all} 
 $(t,x)\in \Xi^{\;\!\prime}.
\hfill
$
\\[1.5ex]
\indent
{\sl Proof.}
\vspace{0.35ex}
Suppose functions $w_1^{}, w_2^{}\in\text{\rm Z}_{_{\Xi^{\;\!\prime}}}.$
Then, by Definition 6.1,
$\bigl(w_1^{}\;\!w_2^{}\;\!,\;\! W\bigl)\;\! \in \text{\rm H}_{_{\Xi^{\;\!\prime}}}$
if and only if the derivative by virtue of system (0.1)
\\[1.75ex]
\mbox{}\hfill                          
$
\displaystyle
{\frak d}\;\! \bigl(w_1^{}(t,x)\;\!w_2^{}(t,x)\bigr) =
w_1^{}(t,x)\;\!w_2^{}(t,x)\;\!W(t,x)
$
for all 
$(t,x)\in \Xi^{\;\!\prime}.
\hfill
$
\\[1.25ex]
\indent
Since
\\[1.25ex]
\mbox{}\hfill                          
$
\displaystyle
{\frak d}\;\! \bigl(w_1^{}(t,x)\;\!w_2^{}(t,x)\bigr)\;\! =
w_2^{}(t,x)\, {\frak d}\;\!w_1^{}(t,x)+
w_1^{}(t,x)\, {\frak d}\;\!w_2^{}(t,x)
$
for all $(t,x)\in \Xi^{\;\!\prime},
\hfill
$
\\[1.75ex]
we see that 
$\bigl(w_1^{}\;\!w_2^{}\;\!,\;\! W\bigl)\;\! \in \text{\rm H}_{_{\Xi^{\;\!\prime}}}$
if and only if
\\[1.75ex]
\mbox{}\hfill                          
$
\displaystyle
w_2^{}(t,x)\, {\frak d}\;\!w_1^{}(t,x)+
w_1^{}(t,x)\, {\frak d}\;\!w_2^{}(t,x)=
w_1^{}(t,x)\;\!w_2^{}(t,x)\;\!W(t,x)
$
for all $(t,x)\in \Xi^{\;\!\prime}.
\hfill
$
\\[1.75ex]
\indent
Hence, 
$\bigl(w_1^{}\;\!w_2^{}\;\!,\;\! W\bigl)\;\! \in \text{\rm H}_{_{\Xi^{\;\!\prime}}}$
if and only if 
\\[1.5ex]
\mbox{}\hfill                          % (6.2)
$
\displaystyle
{\frak d}\;\! w_1^{}(t,x) =
w_1^{}(t,x)\;\! 
\biggl(
W(t,x)-\dfrac{{\frak d}\;\! w_2^{}(t,x)}{w_2^{}(t,x)}
\biggr)
$
for all 
$(t,x)\in\Omega_{_0},
$
\hfill (6.2)
\\[2ex]
where the set $\Omega_{_0}\subset \Xi^{\;\!\prime}$ such that
\\[1.5ex]
\mbox{}\hfill   
$|w_2^{}(t,x)|\ne 0$ for all $(t,x)\in \Omega_{_0}$ 
\ \ and \ \  
$|w_2^{}(t,x)|=0$ for all $(t,x)\in {\sf C}_{_{\Xi^{\;\!\prime}}}\Omega_{_0}.
\hfill
$
\\[1.5ex]
\indent
Assume that
\\[1.25ex]
\mbox{}\hfill                            
$
W(t,x)-\dfrac{{\frak d}\;\! w_2^{}(t,x)}{w_2^{}(t,x)}=
W_1^{}(t,x)
$ 
for all $(t,x)\in \Omega_{_0}.
\hfill 
$
\\[1.75ex]
\indent
Since $w_1^{},\ {\frak d}\;\! w_1^{}, \ W\in \text{\rm Z}_{_{\Xi^{\;\!\prime}}},$
from the identity (6.2), we get $W_1^{}\in \text{\rm Z}_{_{\Xi^{\;\!\prime}}}.$
\vspace{0.75ex}

Therefore we have
$\bigl(w_1^{}\;\!w_2^{}\;\!,\;\! W\bigl)\;\! \in \text{\rm H}_{_{\Xi^{\;\!\prime}}}$
if and only if the identities hold
\\[1.75ex]
\mbox{}\hfill                          % (6.3)
$
\displaystyle
{\frak d}\;\! w_1^{}(t,x)=
w_1^{}(t,x)\;\! W_1^{}(t,x)
$ 
for all 
$
(t,x)\in \Xi^{\;\!\prime}
$
\hfill (6.3)
\\[0.5ex]
and
\\[0.5ex]
\mbox{}\hfill                           % (6.4)
$
\displaystyle
{\frak d}\;\! w_2^{}(t,x)=
w_2^{}(t,x)\;\! \bigl(W(t,x)-W_1^{}(t,x)\bigr)
$ 
for all 
$
(t,x)\in \Xi^{\;\!\prime}.
$
\hfill (6.4)
\\[2ex]
\indent
By Definition 6.1, the identity (6.3) is true
\vspace{0.5ex}
if and only if
$\bigl(w_1^{}\;\!,\;\! W_1^{}\bigl)\;\! \in \text{\rm H}_{_{\Xi^{\;\!\prime}}},$
and the identity (6.4) is true 
if and only if 
$\bigl(w_2^{}\;\!,\;\! W_2^{}\bigl)\;\! \in \text{\rm H}_{_{\Xi^{\;\!\prime}}},\ 
W_2^{}=W-W_1^{}.$
\vspace{0.5ex}

Thus Property 6.5 under condition $m=2$ is proved. 
\vspace{0.35ex}

For $m> 2,$ we prove by induction. $\k$
\vspace{0.75ex}

{\bf Property 6.6.} 
{\it
Let 
$k_j^{}\in\N, \ w_j^{}\in \text{\rm Z}_{_{\Xi^{\;\!\prime}}},\ j=1,\ldots,m.$ Then  
\\[1.5ex]
\mbox{}\hfill                         
$
\displaystyle
\biggl(\,\prod\limits_{j=1}^{m} w_j^{k_j^{}}\;\!,\;\! W\biggl)\;\! 
\in \text{\rm H}_{_{\Xi^{\;\!\prime}}}
\iff
\bigl( w_j^{}\;\!,\;\! W_j^{}\bigl)\;\! \in \text{\rm H}_{_{\Xi^{\;\!\prime}}},
\quad
j=1,\ldots,m.
\hfill
$
\\[1.5ex]
Moreover, the cofactors $W,\, W_1^{},\ldots,W_m^{}$ such that
\\[1.5ex]
\mbox{}\hfill
$
\displaystyle
W(t,x)=\sum\limits_{j=1}^{m} k_j^{}\;\!W_j^{}(t,x)
$
for all} 
$(t,x)\in \Xi^{\;\!\prime}.
\hfill
$
\\[1.5ex]
\indent
{\sl Proof.} 
Taking into account Property 6.3, 
from Property 6.5, we get the following the statement of this property is true. $\k$
\vspace{0.75ex}

{\bf Property 6.7.} 
{\it
If 
$p_\tau^{}\in \text{\rm P}_{_{\!\Xi^{\;\!\prime}}},\, l_{\tau}^{}\!\in\N,\, \tau\!=\!1,\ldots,s,\
w_j^{}\in \text{\rm Z}_{_{\Xi^{\;\!\prime}}},\, k_j^{}\!\in\N, \, j\!=\!1,\ldots,m,$ then  
\\[1.5ex]
\mbox{}\hfill                         
$
\displaystyle
\biggl(\,
\prod\limits_{\tau=1}^{s}\! p_\tau^{\,l_\tau^{}}
\prod\limits_{j=1}^{m}\! w_j^{k_j^{}}, W\!\biggl) 
\in\! \text{\rm H}_{_{\Xi^{\;\!\prime}}}
\!\iff\!
\bigl(p_\tau^{}\;\!, M_\tau^{}\bigl)\;\! \in\! \text{\rm A}_{_{\Xi^{\;\!\prime}}},\,
\tau\!=\!1,\ldots,s,
\ \&\ 
\bigl(w_j^{}\;\!, W_j^{}\bigl)\;\!\in\! \text{\rm H}_{_{\Xi^{\;\!\prime}}},\,
j\!=\!1,\ldots,m.
\hfill
$
\\[1.5ex]
Moreover, the cofactors $W,\, M_1^{},\ldots,M_s^{},\, W_1^{},\ldots,W_m^{}$ such that 
\\[1.5ex]
\mbox{}\hfill
$
\displaystyle
W(t,x)=
\sum\limits_{\tau=1}^{s} l_\tau^{}\;\!M_\tau^{}(t,x)+
\sum\limits_{j=1}^{m} k_j^{}\;\!W_j^{}(t,x)
$
for all} $(t,x)\in \Xi^{\;\!\prime}.
\hfill
$
\\[1.5ex]
\indent
{\sl Proof.} 
\vspace{0.5ex}
If $p_\tau^{}\in \text{\rm P}_{_{\!\Xi^{\;\!\prime}}},$ then 
$p_\tau^{}\in \text{\rm Z}_{_{\Xi^{\;\!\prime}}},\ \tau=1,\ldots,s.$
Now, from Property 6.6, we obtain the statement of Property 6.7. $\k$ 
\vspace{1ex}

Separate the real and imaginary parts of functions 
$w, W\in \text{\rm Z}_{_{\Xi^{\;\!\prime}}}\colon$
\\[1.5ex]
\mbox{}\hfill
$
u\colon (t,x)\to\ {\rm Re}\;\!w(t,x)
$
for all 
$(t,x)\in \Xi^{\;\!\prime},
\qquad
v\colon (t,x)\to\ {\rm Im}\;\!w(t,x)
$
for all $(t,x)\in \Xi^{\;\!\prime},
\hfill
$
\\[2ex]
\mbox{}\hfill
$
U\colon (t,x)\to\ {\rm Re}\;\!W(t,x)
$
for all 
$(t,x)\in \Xi^{\;\!\prime},
\qquad
V\colon (t,x)\to\ {\rm Im}\;\!W(t,x)
$
for all 
$(t,x)\in \Xi^{\;\!\prime}.
\hfill
$
\\[2ex]
\indent
{\bf Property 6.8.}\! 
\vspace{0.15ex}
{\it
Suppose 
$\!p\in\! \text{\rm P}_{_{\!\Xi^{\;\!\prime}}},\!$ 
functions $\!u_1^{},v_1^{}\!\in\! \text{\rm P}_{_{\!\Xi^{\;\!\prime}}}\!$ 
are relatively prime, $\!w_1^{}\!=\!u_1^{}\!+i\;\!v_1^{},$ $w=p\;\!w_1^{}.$
Then we have
\\[1.25ex]
\mbox{}\hfill                         
$
\displaystyle
\bigl(w, W\bigl)\, \in \text{\rm H}_{_{\Xi^{\;\!\prime}}}
\iff
\bigl(p, M\bigl)\, \in \text{\rm A}_{_{\Xi^{\;\!\prime}}}
\ \&\ 
\bigl(w_1^{}, W_1^{}\bigl)\, \in \text{\rm H}_{_{\Xi^{\;\!\prime}}}.
\hfill
$
\\[1.75ex]
Moreover, the cofactors $W,\, M,\, W_1^{}$ such that
\\[1.5ex]
\mbox{}\hfill
$
\displaystyle
W(t,x)=M(t,x)+W_1^{}(t,x)
$
for all} 
$(t,x)\in \Xi^{\;\!\prime}.
\hfill
$
\\[1.75ex]
\indent
{\sl Proof}\;\! follows from Property 6.7 under the condition $s=m=1.\ \k$
\vspace{0.75ex}

Using Property 6.8, from the set $\text{\rm H}_{_{\Xi^{\;\!\prime}}}$ 
we can extract the complex-valued polynomial partial integrals 
such that 
their real and imaginary parts are relatively prime.
\vspace{0.5ex}

{\bf Theorem 6.1}
\vspace{0.15ex}
(existence criterion of complex-valued polynomial partial integral). 
{\it
Suppose 
$u,v\in \text{\rm P}_{_{\!\Xi^{\;\!\prime}}}.$ 
Then $(u+i\;\!v,\;\! U+i\;\!V)\in \text{\rm H}_{_{\Xi^{\;\!\prime}}}$
if and only if the system of identities holds
\\[1.5ex]
\mbox{}\hfill                            % (6.5)
$
{\frak d}\;\!u(t,x)=
u(t,x)\;\!U(t,x)-v(t,x)\;\!V(t,x)
$
\ for all $(t,x)\in \Xi^{\;\!\prime},
$
\hfill {\rm (6.5)}
\\[2.25ex]
\mbox{}\hfill                           % (6.6)
$
{\frak d}\;\!v(t,x)=
u(t,x)\;\!V(t,x)+v(t,x)\;\!U(t,x)
$
\ for all $(t,x)\in \Xi^{\;\!\prime},
$
\hfill {\rm (6.6)}
\\[1.75ex]
where the functions $U, V\in \text{\rm P}_{_{\!\Xi^{\;\!\prime}}}.$}
\vspace{0.5ex}

{\sl Proof}\;\!
is that the identity (6.1) is true if and only if 
the system of identities $(6.5)\& (6.6)$ is true, where 
$U, V\in \text{\rm P}_{_{\!\Xi^{\;\!\prime}}}.\ \k$
\vspace{0.75ex}

{\bf Property 6.9.}
{\it
We claim that}
\\[1.5ex]
\mbox{}\hfill                        
$
(u+i\;\!v,\;\! U+i\;\!V)\in \text{\rm H}_{_{\Xi^{\;\!\prime}}}
\iff
(u-i\;\!v,\;\! U-i\;\!V)\in \text{\rm H}_{_{\Xi^{\;\!\prime}}}.
\hfill 
$
\\[1.75ex]
\indent
{\sl Proof}\;\!
is based on Theorem 6.1 and the fact that the system of identities $(6.5)\& (6.6)$ 
is invariant under simultaneous replacement $v$ by ${}-v$ and $V$ by ${}-V.\ \k$
\vspace{0.75ex}

{\bf Lemma 6.1.}
\vspace{0.35ex}
{\it
Suppose  $u,v\in \text{\rm P}_{_{\!\Xi^{\;\!\prime}}}.$ 
Then, we have $(u+i\;\!v,\;\! U+i\;\!V)\in \text{\rm H}_{_{\Xi^{\;\!\prime}}}$
if and only if 
$(u^2+v^2,\;\! 2\;\!U)\in \text{\rm A}_{_{\Xi^{\;\!\prime}}}$
and the identity {\rm (6.5)} {\rm(}or {\rm (6.6))} holds.
}
\vspace{0.5ex}

{\sl Proof.}
By Theorem 6.1, $(u+i\;\!v,\;\! U+i\;\!V)\in \text{\rm H}_{_{\Xi^{\;\!\prime}}}$
if and only if the identities (6.5) and (6.6) are true.

Multiplying both sides of the identity (6.5) by $2u,$ we get 
\\[1.5ex]
\mbox{}\hfill                           
$
2\;\!u(t,x)\,{\frak d}\;\!u(t,x)=
2\;\!u^2(t,x)\;\!U(t,x)-2\;\!u(t,x)\;\!v(t,x)\;\!V(t,x)
$
\ for all $(t,x)\in \Xi^{\;\!\prime}.
\hfill 
$
\\[1.5ex]
\indent
Multiplying both sides of the identity (6.6) by $2v,$ we get
\\[2.25ex]
\mbox{}\hfill                          
$
2\;\!v(t,x)\,{\frak d}\;\!v(t,x)=
2\;\!v(t,x)\;\!u(t,x)\;\!V(t,x)+2\;\!v^2(t,x)\;\!U(t,x)
$
\ for all $(t,x)\in \Xi^{\;\!\prime}.
\hfill
$
\\[1.75ex]
\indent
Now, summing the obtained identities, we see that 
the system of identities $(6.5)\& (6.6)$ is reduced to the system, which contains the identity
\\[1.5ex]
\mbox{}\hfill                          
$
{\frak d}\;\!\bigl(u^2(t,x)+v^2(t,x)\bigr)=
2\;\!\bigl(u^2(t,x)+v^2(t,x)\bigr)\;\!U(t,x)
$
\ for all 
$(t,x)\in \Xi^{\;\!\prime}
\hfill
$
\\[1.5ex]
and at least one of the identities (6.5) or (6.6).

Since the last identity holds (by Theorem 2.1)
if and only if 
$(u^2+v^2,\;\! 2\;\!U)\in \text{\rm A}_{_{\Xi^{\;\!\prime}}},$
we see that the statement of Lemma 6.1 is true. $\k$
\vspace{0.75ex}

{\bf Lemma 6.2.}
\vspace{0.5ex}
{\it
Suppose $u,v\!\in\! \text{\rm P}_{_{\!\Xi^{\;\!\prime}}}\!\!$ are relatively prime, 
the set $\Omega_{_0}\!\subset\! \Xi^{\;\!\prime}$ such that   
$u(t,x)\!\ne\! 0\!$ for all $(t,x)\in \Omega_{_0}$ and
$u(t,x)=0$ for all $(t,x)\in {\sf C}_{_{\Xi^{\;\!\prime}}}\Omega_{_0}.$
\vspace{0.5ex}
Then $(u+i\;\!v,\;\! U+i\;\!V)\in \text{\rm H}_{_{\Xi^{\;\!\prime}}}$
if and only if
$\Bigl(\exp\arctan\dfrac{v}{u}\,,\;\! V\Bigr)\in 
\text{\rm E}_{_{\Omega_{{}_{\tiny\;\! 0}}}},$
\vspace{0.35ex}
where the function $U\in \text{\rm P}_{_{\!\Xi^{\;\!\prime}}}$
is obtained from the identity {\rm (6.5)} {\rm(}or {\rm (6.6)).}
}
\vspace{0.5ex}

{\sl Proof}\;\!
follows from Theorem 6.1 and Property 3.11. $\k$
\vspace{0.75ex}

{\bf Theorem 6.2}
\vspace{0.35ex}
(existence criterion of complex-valued polynomial partial integral). 
{\it
Suppose functions
\vspace{0.75ex}
$u,v\in \text{\rm P}_{_{\!\Xi^{\;\!\prime}}}$ are relatively prime, 
the set $\Omega_{_0}\subset \Xi^{\;\!\prime}$ such that  
$u(t,x)\ne 0$ for all $(t,x)\in \Omega_{_0}$ and
$u(t,x)=0$ for all $(t,x)\in {\sf C}_{_{\Xi^{\;\!\prime}}}\Omega_{_0}.$
Then}
\\[1.75ex]
\mbox{}\hfill                          
$
(u+i\;\!v,\;\! U+i\;\!V)\in \text{\rm H}_{_{\Xi^{\;\!\prime}}}
\iff
(u^2+v^2,\;\! 2\;\!U)\in \text{\rm A}_{_{\Xi^{\;\!\prime}}}
\ \&\ 
\Bigl(\exp\arctan\dfrac{v}{u}\,,\;\! V\Bigr)\in 
\text{\rm E}_{_{\Omega_{{}_{\tiny\;\! 0}}}}.
\hfill                          
$
\\[1.5ex]
\indent
{\sl Proof}\;\!
follows from Lemmas 6.1 and 6.2. $\k$
\vspace{1ex}

{\bf Property 6.10.}
{\it
Suppose functions
\vspace{0.75ex}
$u,v\in \text{\rm P}_{_{\!\Xi^{\;\!\prime}}}$ are relatively prime, 
the set $\Omega_{_0}\subset \Xi^{\;\!\prime}$ such that   
$u(t,x)\ne 0$ for all $(t,x)\in \Omega_{_0}$ and 
$u(t,x)=0$ for all $(t,x)\in {\sf C}_{_{\Xi^{\;\!\prime}}}\Omega_{_0}.$
Then  
\\[2ex]
\mbox{}\hfill
$
\Bigl(\exp\arctan\dfrac{v}{u}\,,\;\! V\Bigr)\in 
\text{\rm E}_{_{\Omega_{{}_{\tiny\;\! 0}}}}
\hfill
$
\\[1.75ex]
if and only if
$(u^2+v^2,\;\! 2\;\!U)\in \text{\rm A}_{_{\Xi^{\;\!\prime}}}$
and the identity {\rm (6.5)} {\rm(}or {\rm (6.6)) holds.}
}
\vspace{0.75ex}

{\sl Proof}\;\!
follows from Lemmas 6.1 and 6.2. $\k$
\vspace{1ex}

{\bf Remark 6.2.}
\vspace{0.35ex}
Taking into account Remark 2.1 and Theorem  3.1, 
from Theorem 6.2, we get the following
$\deg_x^{} U\leq d-1$ and $\deg_x^{} V\leq d-1.$
Therefore, in Definition 6.1 we have 
\\[1.5ex]
\mbox{}\hfill
$
\deg_x^{} W\leq d-1,
\hfill
$
\\[1.5ex]
and in the identities (3.5), (6.5), and (6.6) we obtain
\\[1.5ex]
\mbox{}\hfill
$
\deg_x^{} U\leq d-1,
\quad 
\deg_x^{} V\leq d-1.
\hfill
$
\\[2ex]
\indent
{\bf Property 6.11.} 
\vspace{0.15ex}
{\it
Let 
$u, v, U\in \text{\rm P}_{_{\!\Xi^{\;\!\prime}}}.$
Then we have}
\\[2ex]
\mbox{}\hfill                         
$
\displaystyle
(u+i\;\!v,\;\! U)\in \text{\rm H}_{_{\Xi^{\;\!\prime}}}
\iff
\bigl(u, U\bigl)\, \in \text{\rm A}_{_{\Xi^{\;\!\prime}}}
\ \&\ 
\bigl(v, U\bigl)\, \in \text{\rm A}_{_{\Xi^{\;\!\prime}}}.
\hfill
$
\\[2.25ex]
\indent
{\sl Proof.}
If $V(t,x)=0$ for all $(t,x)\in \Xi^{\;\!\prime},$
then identities (6.5) and (6.6) have the form 
\\[2ex]
\mbox{}\hfill                           
$
{\frak d}\;\!u(t,x)=
u(t,x)\;\!U(t,x)
$
for all 
$(t,x)\in \Xi^{\;\!\prime},
\quad
{\frak d}\;\!v(t,x)=
v(t,x)\;\!U(t,x)
$
for all 
$(t,x)\in \Xi^{\;\!\prime}.
\hfill
$
\\[2ex]
\indent
Now, using these identities and Theorem 2.1, from Theorem 6.1 it follows that 
the statement of Property 6.11 is true. $\k$
\vspace{1ex}

{\bf Property 6.12.}
{\it
Suppose functions
\vspace{0.75ex}
$u,v\in \text{\rm P}_{_{\!\Xi^{\;\!\prime}}}$ are relatively prime, 
the set $\Omega_{_0}\subset \Xi^{\;\!\prime}$ such that   
$u(t,x)\ne 0$ for all $(t,x)\in \Omega_{_0}$ and
$u(t,x)=0$ for all $(t,x)\in {\sf C}_{_{\Xi^{\;\!\prime}}}\Omega_{_0}.$
Then}
\\[2ex]
\mbox{}\hfill
$
(u+i\;\!v,\;\! U)\in \text{\rm H}_{_{\Xi^{\;\!\prime}}}
\iff
(u^2+v^2,\;\! 2\;\!U)\in \text{\rm A}_{_{\Xi^{\;\!\prime}}}
\ \&\ \,
\dfrac{v}{u}\in \text{\rm I}_{_{\Omega_{{}_{\tiny\;\! 0}}}}.
\hfill
$
\\[2.25ex]
\indent
{\sl Proof.}
\vspace{0.35ex}
Taking into account Theorem 0.2 and Property 0.1, 
from Theorem 6.2 under the condition 
$V(t,x)=0$ for all $(t,x)\in \Xi^{\;\!\prime},$
we get the statement of Property 6.12 is true. $\k$
\vspace{1.25ex}

{\bf Property 6.13.} 
\vspace{0.15ex}
{\it
Let 
$u, v, V\in \text{\rm P}_{_{\!\Xi^{\;\!\prime}}}.$
Then $(u+i\;\!v,\;\! i\;\!V)\in \text{\rm H}_{_{\Xi^{\;\!\prime}}}$
if and only if
\\[2ex]
\mbox{}\hfill                         
$
{\frak d}\;\!u(t,x)=
{}-v(t,x)\;\!V(t,x)
$
for all $(t,x)\!\in\! \Xi^{\;\!\prime},
\quad
{\frak d}\;\!v(t,x)=
u(t,x)\;\!V(t,x)
$
for all} $(t,x)\!\in\! \Xi^{\;\!\prime}.
\hfill
$
\\[2.25ex]
\indent
{\sl Proof}\;\!
follows from Theorem 6.1 under the condition 
$U(t,x)=0$ for all $(t,x)\in \Xi^{\;\!\prime}.\ \k$
\vspace{1.25ex}

{\bf Property 6.14.}
\vspace{0.5ex}
{\it
Suppose functions
$u,v\in \text{\rm P}_{_{\!\Xi^{\;\!\prime}}}$ are relatively prime, 
$V\in \text{\rm P}_{_{\!\Xi^{\;\!\prime}}},$
the set $\Omega_{_0}\subset \Xi^{\;\!\prime}$ such that   
$u(t,x)\ne 0$ for all $(t,x)\in \Omega_{_0}$ and 
\vspace{0.5ex}
$u(t,x)=0$ for all $(t,x)\in {\sf C}_{_{\Xi^{\;\!\prime}}}\Omega_{_0}.$
Then we claim that}
\\[1.5ex]
\mbox{}\hfill
$
(u+i\;\!v,\;\! i\;\!V)\in \text{\rm H}_{_{\Xi^{\;\!\prime}}}
\iff
\Bigl(\exp\arctan\dfrac{v}{u}\,,\, V\Bigr)\in \text{\rm E}_{_{\Omega_{{}_{\tiny\;\! 0}}}}
\ \&\ \,
u^2+v^2 \in \text{\rm I}_{_{\Xi^{\;\!\prime}}}.
\hfill
$
\\[2ex]
\indent
{\sl Proof.}
\vspace{0.35ex}
Taking into account Theorem 0.2, 
from Theorem 6.2 under the condition $\!U(t,x)\!=\!0$ for all $(t,x)\in \Xi^{\;\!\prime},$
we obtain the statement of Property 6.14 is true. $\k$
\vspace{1.25ex}

\newpage

{\bf Property 6.15.}
\vspace{0.75ex}
{\it
Suppose numbers $\gamma_1^{}, \gamma_2^{}\in\R\backslash\{0\},$
functions 
$u,v\in \text{\rm P}_{_{\!\Xi^{\;\!\prime}}}$ are relatively prime, 
$(u+i\;\!v,\;\! U+i\;\!V)\in \text{\rm H}_{_{\Xi^{\;\!\prime}}},$
\vspace{0.75ex}
the set $\Omega_{_0}\!\subset\! \Xi^{\;\!\prime}\!$ such that  
$u(t,x)\ne 0$ for all $(t,x)\!\in\! \Omega_{_0}$ and
$u(t,x)=0$ for all $(t,x)\in {\sf C}_{_{\Xi^{\;\!\prime}}}\Omega_{_0}.$
Then 
\\[2ex]
\mbox{}\hfill
$
\Bigl((u^2+v^2)^{{}^{\scriptsize \gamma_1^{}}}
\exp\Bigl(\gamma_2^{}\;\!\arctan\dfrac{v}{u}\Bigr)\,,\, M\Bigr)\in 
\text{\rm J}_{_{\Omega_{{}_{\tiny\;\! 0}}}}
\hfill
$
\\[1ex]
if and only if
\\[1.25ex]
\mbox{}\hfill
$
M(t,x)=2\gamma_1^{}\;\!U(t,x)+\gamma_2^{}\;\!V(t,x)
$
\ for all} $(t,x)\in \Xi^{\;\!\prime}.
\hfill
$
\\[2ex]
\indent
{\sl Proof}\;\!
follows from Theorem 6.2 and Property 1.9. $\k$
\\[5.25ex]
\centerline{
{\bf  7. Multiple complex-valued polynomial partial integrals}
}
\\[1.5ex]
\indent
{\bf Definition 7.1.}
\vspace{0.25ex}
{\it
A complex-valued polynomial partial integral $w$ with cofactor $W$ 
on the domain $\Xi^{\;\!\prime}$ of system {\rm (0.1)} 
is said to be  
\textit{\textbf{multiple}}
\vspace{0.15ex}
if there exist a natural number $h$ and a function
$z\in \text{\rm Z}_{_{\Xi^{\;\!\prime}}}$
\vspace{0.35ex}
relatively prime to the function $w$
such that the derivative by virtue of system {\rm (0.1)} 
\\[1.5ex]
\mbox{}\hfill                             % (7.1)
$
\displaystyle
{\frak d}\,\dfrac{z(t,x)}{w^{\;\!h}(t,x)}=Q(t,x)
$
\ for all 
$(t,x)\in \Omega_{_0},
$
\hfill {\rm (7.1)}
\\[2ex]
where the function $Q\in \text{\rm Z}_{_{\Xi^{\;\!\prime}}}$ 
has the degree $\deg_{\;\!x}^{} Q\leq d-1,$
\vspace{0.75ex}
the set 
$\Omega_{_0}\subset\Xi^{\;\!\prime}$ such that   
$
|w(t,x)|\ne 0$ for all $(t,x)\in\Omega_{_0}$ and
$|w(t,x)|= 0$ for all $(t,x)\in {\sf C}_{_{\Xi^{\;\!\prime}}}\Omega_{_0}\;\!.$
}
\vspace{1.25ex}

By $\text{\rm G}_{_{\Xi^{\;\!\prime}}}$ 
\vspace{0.5ex}
denote the set of multiple complex-valued polynomial partial integrals on the domain 
$\Xi^{\;\!\prime}$ of system (0.1).
From Definition 7.1 it follows that the set
$\text{\rm G}_{_{\Xi^{\;\!\prime}}}\subset
\text{\rm H}_{_{\Xi^{\;\!\prime}}}.$
\vspace{1ex}

The phrase "the complex-valued polynomial partial integral $w$ with cofactor $W$ 
\vspace{0.25ex}
on the domain $\Xi^{\;\!\prime}$ of system (0.1) is multiple and the identity (7.1) holds" 
is denoted by
\\[1.5ex]
\mbox{}\hfill
$
\bigl((w, W), (h,z,Q)\bigr)\in \text{\rm G}_{_{\Xi^{\;\!\prime}}}.
\hfill
$
\\[2ex]
\indent
{\bf Theorem 7.1}
\vspace{0.35ex}
(existence criterion of multiple complex-valued polynomial partial integral).
{\it
$\bigl((w, W), (h,z,Q)\bigr)\in \text{\rm G}_{_{\Xi^{\;\!\prime}}}$
\vspace{0.5ex}
if and only if the identities {\rm(6.1)} and {\rm(7.1)} are true,
where the number $h\in\N,$ the functions
\vspace{0.75ex}
$z, Q\in \text{\rm Z}_{_{\Xi^{\;\!\prime}}},$
the functions $w$ and $z$ are relatively prime, $\deg_{\;\!x}^{} Q\leq d-1,$
\vspace{0.5ex}
the set 
$\Omega_{_0}\subset\Xi^{\;\!\prime}$ such that  
$|w(t,x)|\ne 0$ for all $(t,x)\in\Omega_{_0}$ and
$|w(t,x)|= 0$ for all $(t,x)\in {\sf C}_{_{\Xi^{\;\!\prime}}}\Omega_{_0}\;\!.$
}
\vspace{0.75ex}

{\sl Proof.} 
Using Remarks 6.1 and 6.2, from Definitions 6.1 and 7.1, we get  
the statement of Theorem 7.1. $\k$
\vspace{1ex}

{\bf Theorem 7.2}
\vspace{0.25ex}
(existence criterion of multiple complex-valued polynomial partial integral).
{\it
Suppose $h\in\N,$ functions 
\vspace{0.5ex}
$w, z\in \text{\rm Z}_{_{\Xi^{\;\!\prime}}}\!$ are relatively prime, 
the set 
$\Omega_{_0}\subset\Xi^{\;\!\prime}$ such that
$|w(t,x)|\ne 0$ for all $(t,x)\in\Omega_{_0}$ and
$|w(t,x)|= 0$ for all $(t,x)\in {\sf C}_{_{\Xi^{\;\!\prime}}}\Omega_{_0}.$
Then}
\\[1.75ex]
\mbox{}\hfill                           
$
\displaystyle
\bigl((w, W), (h,z,Q)\bigr)\in \text{\rm G}_{_{\Xi^{\;\!\prime}}}
\iff
(w, W)\in \text{\rm H}_{_{\Xi^{\;\!\prime}}}
\ \& \
\hfill                           
$
\\[2ex]
\mbox{}\hfill                           
$
\& \
\Bigl(\exp\;\!{\rm Re}\;\!\dfrac{z}{w^{\;\!h}}\,, {\rm Re}\;\!Q\Bigr)\in 
\text{\rm E}_{_{\Omega_{{}_{\tiny\;\! 0}}}}
\ \,\& \ \,
\Bigl(\exp\;\!{\rm Im}\;\!\dfrac{z}{w^{\;\!h}}\,, {\rm Im}\;\!Q\Bigr)\in 
\text{\rm E}_{_{\Omega_{{}_{\tiny\;\! 0}}}}.
\hfill
$
\\[1.75ex]
\indent
{\sl Proof.}
First we separate the real and imaginary parts in the identity (7.1).
Then, taking into account 
the existence criterion of exponential partial integral (Theorem 3.1), we have 
the identity (7.1) is true if and only if 
\\[1.25ex]
\mbox{}\hfill
$
\Bigl(\exp\;\!{\rm Re}\;\!\dfrac{z}{w^{\;\!h}}\,, {\rm Re}\;\!Q\Bigr)\in 
\text{\rm E}_{_{\Omega_{{}_{\tiny\;\! 0}}}}
$
\ \ and \ \
$
\Bigl(\exp\;\!{\rm Im}\;\!\dfrac{z}{w^{\;\!h}}\,, {\rm Im}\;\!Q\Bigr)\in 
\text{\rm E}_{_{\Omega_{{}_{\tiny\;\! 0}}}}.
\hfill
$
\\[1.75ex]
\indent
Now, using Definition 6.1, from Theorem 7.1, 
\vspace{0.75ex}
we get the statement of Theorem 7.2. $\k$

{\bf Theorem 7.3}
\vspace{0.25ex}
(existence criterion of multiple complex-valued polynomial partial integral).
{\it
Suppose $h\in\N,$ functions
\vspace{0.75ex}
$u, v\in \text{\rm P}_{_{\!\Xi^{\;\!\prime}}}$ are relatively prime, 
a function $z\in \text{\rm Z}_{_{\Xi^{\;\!\prime}}}$ is relatively prime to the function 
 $u+i\;\!v,$ the set
\vspace{0.75ex}
$\Omega_{_0}\subset\Xi^{\;\!\prime}$ such that   
$u(t,x)\ne 0$ for all $(t,x)\in\Omega_{_0}$ and
$u(t,x)= 0$ for all $(t,x)\in {\sf C}_{_{\Xi^{\;\!\prime}}}\Omega_{_0}.$
Then we claim that}
\\[2ex]
\mbox{}\hfill                           
$
\displaystyle
\bigl((u+i\;\!v,\;\! U+i\;\!V), (h,z,Q)\bigr)\in \text{\rm G}_{_{\Xi^{\;\!\prime}}}
\iff
\Bigl(\bigl(u^2+v^2,\;\! 2\;\!U\bigr), 
\bigl(h,\;\!{\rm Re}\;\!\bigl(z(u-i\;\!v)^h\bigr),\;\! {\rm Re}\;\!Q\bigr)\Bigr)\in 
\text{\rm B}_{_{\Xi^{\;\!\prime}}}
\ \& \
\hfill                           
$
\\[2.25ex]
\mbox{}\hfill                           
$
\& \
\Bigl(\bigl(u^2+v^2,\;\! 2\;\!U\bigr), 
\bigl(h,\;\!{\rm Im}\;\!\bigl(z(u-i\;\!v)^h\bigr),\;\! {\rm Im}\;\!Q\bigr)\Bigr)\in 
\text{\rm B}_{_{\Xi^{\;\!\prime}}}
\ \,\& \ \,
\Bigl(\exp\arctan\dfrac{v}{u}\,,\;\! V\Bigr)\in 
\text{\rm E}_{_{\Omega_{{}_{\tiny\;\! 0}}}}.
\hfill
$
\\[1.75ex]
\indent
{\sl Proof}\;\!
is  based on Theorem 7.2 and the following facts.
\vspace{0.5ex}

By Theorem 6.2,
\\[1.25ex]
\mbox{}\hfill                           
$
\displaystyle
(u+i\;\!v,\;\! U+i\;\!V)\in \text{\rm H}_{_{\Xi^{\;\!\prime}}}
\iff
\bigl(u^2+v^2,\;\! 2\;\!U\bigr)\in\text{\rm A}_{_{\Xi^{\;\!\prime}}}
\ \& \
\Bigl(\exp\arctan\dfrac{v}{u}\,,\;\! V\Bigr)\in 
\text{\rm E}_{_{\Omega_{{}_{\tiny\;\! 0}}}}.
\hfill                           
$
\\[1.5ex]
\indent
Using 
\\[1.25ex]
\mbox{}\hfill                           
$
\displaystyle
\bigl(u^2+v^2,\;\! 2\;\!U\bigr)\in\text{\rm A}_{_{\Xi^{\;\!\prime}}},
\quad
\dfrac{z}{w^h}=
\dfrac{z}{(u+i\;\!v)^h}=
\dfrac{z\;\!(u-i\;\!v)^h}{(u^2+v^2)^h}\,,
\hfill                           
$
\\[1.75ex]
from Theorem 5.3, we get
\\[1.75ex]
\mbox{}\hfill                           
$
\displaystyle
\Bigl(\exp\;\!{\rm Re}\;\!\dfrac{z}{w^h}\,,\, {\rm Re}\;\!Q\Bigr)\in 
\text{\rm E}_{_{\Omega_{{}_{\tiny\;\! 0}}}}
\iff
\biggl(\exp\;\!\dfrac{{\rm Re}\;\!\bigl(z\;\!(u-i\;\!v)^h\bigr)}{(u^2+v^2)^h}\,,\, {\rm Re}\;\!Q\biggr)\in 
\text{\rm E}_{_{\Omega_{{}_{\tiny\;\! 0}}}}
\iff
\hfill                           
$
\\[2ex]
\mbox{}\hfill                           
$
\iff
\Bigl(\bigl(u^2+v^2,\;\! 2\;\!U\bigr), 
\bigl(h,\;\!{\rm Re}\;\!\bigl(z(u-i\;\!v)^h\bigr),\;\! {\rm Re}\;\!Q\bigr)\Bigr)\in 
\text{\rm B}_{_{\Xi^{\;\!\prime}}};
\hfill                           
$
\\[2ex]
\mbox{}\hfill                           
$
\displaystyle
\Bigl(\exp\;\!{\rm Im}\;\!\dfrac{z}{w^h}\,,\, {\rm Im}\;\!Q\Bigr)\in 
\text{\rm E}_{_{\Omega_{{}_{\tiny\;\! 0}}}}
\iff
\biggl(\exp\;\!\dfrac{{\rm Im}\;\!\bigl(z\;\!(u-i\;\!v)^h\bigr)}{(u^2+v^2)^h}\,,\, {\rm Im}\;\!Q\biggr)\in 
\text{\rm E}_{_{\Omega_{{}_{\tiny\;\! 0}}}}
\iff
\hfill                           
$
\\[2ex]
\mbox{}\hfill                           
$
\iff
\Bigl(\bigl(u^2+v^2,\;\! 2\;\!U\bigr), 
\bigl(h,\;\!{\rm Im}\;\!\bigl(z(u-i\;\!v)^h\bigr),\;\! {\rm Im}\;\!Q\bigr)\Bigr)\in 
\text{\rm B}_{_{\Xi^{\;\!\prime}}}.\ \k
\hfill                           
$
\\[2.75ex]
\indent
{\bf Corollary 7.1.}
\vspace{0.5ex}
{\it
Suppose $u, v, a, b\in \text{\rm P}_{_{\!\Xi^{\;\!\prime}}},$ 
the functions $u,\ v$ are relatively prime, 
the functions $u+i\;\!v$ and $a+i\;\!b$ are relatively prime, the set
\vspace{0.75ex}
$\Omega_{_0}\subset\Xi^{\;\!\prime}$ such that  
$u(t,x)\ne 0$ for all $(t,x)\in\Omega_{_0}$ and
$u(t,x)= 0$ for all $(t,x)\in {\sf C}_{_{\Xi^{\;\!\prime}}}\Omega_{_0}.$
Then we have}
\\[2.25ex]
\mbox{}\hfill                           
$
\displaystyle
\bigl((u+i\;\!v,\;\! U+i\;\!V), (1, a+i\;\!b, K+i\;\!L)\bigr)\in \text{\rm G}_{_{\Xi^{\;\!\prime}}}
\iff
\Bigl(\bigl(u^2+v^2,\;\! 2\;\!U\bigr), 
\bigl(1,\;\!a\;\!u+b\;\!v,\;\! K\bigr)\Bigr)\in 
\text{\rm B}_{_{\Xi^{\;\!\prime}}}
\ \& \
\hfill                           
$
\\[2.5ex]
\mbox{}\hfill                           
$
\& \
\Bigl(\bigl(u^2+v^2,\;\! 2\;\!U\bigr), 
\bigl(1,\;\!b\;\!u-a\;\!v,\;\! L\bigr)\Bigr)\in 
\text{\rm B}_{_{\Xi^{\;\!\prime}}}
\ \,\& \ \,
\Bigl(\exp\arctan\dfrac{v}{u}\,,\;\! V\Bigr)\in 
\text{\rm E}_{_{\Omega_{{}_{\tiny\;\! 0}}}}.
\hfill
$
\\[2.25ex]
\indent
{\sl Proof}\;\!
\vspace{1ex}
follows from Theorem 7.3 and the identity  
$(a+i\;\!b)(u-i\;\!v)=a\;\!u+b\;\!v+(b\;\!u-a\;\!v)\;\!i.\ \k$

{\bf Corollary 7.2.}
\vspace{0.5ex}
{\it
Suppose $u, v, a, b\in \text{\rm P}_{_{\!\Xi^{\;\!\prime}}},$ 
the functions $u,\ v$ are relatively prime, 
the functions $u+i\;\!v$ and $a+i\;\!b$ are relatively prime, the set
\vspace{0.75ex}
$\Omega_{_0}\subset\Xi^{\;\!\prime}$ such that  
$u(t,x)\ne 0$ for all $(t,x)\in\Omega_{_0}$ and
$u(t,x)= 0$ for all $(t,x)\in {\sf C}_{_{\Xi^{\;\!\prime}}}\Omega_{_0}.$
Then we have}
\\[1.75ex]
\mbox{}\hfill                           
$
\displaystyle
\bigl((u+i\;\!v,\;\! U+i\;\!V), (2, a+i\;\!b, K+i\;\!L)\bigr)\in \text{\rm G}_{_{\Xi^{\;\!\prime}}}
\iff
\hfill                           
$
\\[1.75ex]
\mbox{}\hfill                           
$
\iff
\Bigl(\bigl(u^2+v^2,\;\! 2\;\!U\bigr), 
\bigl(2,\;\!a(u^2-v^2)+2\;\!b\;\!u\;\!v,\;\! K\bigr)\Bigr)\in 
\text{\rm B}_{_{\Xi^{\;\!\prime}}}
\ \& \
\hfill                           
$
\\[1.75ex]
\mbox{}\hfill                           
$
\& \
\Bigl(\bigl(u^2+v^2,\;\! 2\;\!U\bigr), 
\bigl(2,\;\!b\;\!(u^2-v^2)-2a\;\!u\;\!v,\;\! L\bigr)\Bigr)\in 
\text{\rm B}_{_{\Xi^{\;\!\prime}}}
\ \,\& \ \,
\Bigl(\exp\arctan\dfrac{v}{u}\,,\;\! V\Bigr)\in 
\text{\rm E}_{_{\Omega_{{}_{\tiny\;\! 0}}}}.
\hfill
$
\\[1.5ex]
\indent
{\sl Proof}\;\!
follows from Theorem 7.3 and the identity  
\\[1.5ex]
\mbox{}\hfill
$
(a+i\;\!b)(u-i\;\!v)^2=a\;\!(u^2-v^2)+2\;\!b\;\!u\;\!v+
\bigl(b\;\!(u^2-v^2)-2\;\!a\;\!u\;\!v\bigr)\;\!i.\ \k
\hfill
$
\\[2ex]
\indent
{\bf Definition 7.2.}
\vspace{0.35ex} 
{\it
A complex-valued polynomial partial integral $w$  on the domain $\Xi^{\;\!\prime}$ 
of system {\rm (0.1)} is said to be
\textit{\textbf{multiple with multiplicity}}
\vspace{0.75ex} 
$
\varkappa =1 + \sum\limits_{\xi=1}^{\varepsilon} \delta_{\xi}^{}
$
if there exist natural numbers
\vspace{0.75ex} 
$h_{\xi}^{},\, \xi\! =\!1,\ldots, \varepsilon,\!$ 
and functions
\vspace{0.35ex} 
$
z_{_{\scriptstyle h_\xi^{}f_\xi^{}}}\!\in \text{\rm Z}_{_{\Xi^{\;\!\prime}}},\ 
f_{\xi}^{}\!=\!1,\ldots,\delta_{\xi}^{},\ \xi\!=\!1,\ldots,\varepsilon,$
that cor\-res\-pond to these numbers and relatively prime to the function $w,$
such that the identities hold 
\\[2ex]
\mbox{}\hfill                  
$
\displaystyle
{\frak d}\, \dfrac{z_{_{\scriptstyle h_\xi^{}  f_\xi^{}} }(t,x)}{\displaystyle  w^{\;\!h_\xi^{}} (t,x)}=
Q_{h_\xi^{}  f_\xi^{}}^{} (t,x)
$ 
\ for all 
$(t,x)\in \Omega_{_0},
\quad 
f_{\xi}^{}=1,\ldots, \delta_{\xi}^{}, \ \   
\xi=1,\ldots,\varepsilon,
\hfill
$
\\[2.5ex]
where the functions 
\vspace{1ex}
$Q_{h_\xi^{}  f_\xi^{}}^{}\in\text{\rm Z}_{_{\Xi^{\;\!\prime}}}$ 
and have the degrees
$
\deg_{\;\!x}^{}  Q_{h_\xi^{}  f_\xi^{}}^{}\leq d-1, \ 
f_{\xi}^{}=1,\ldots, \delta_{\xi}^{}, 
\linebreak  
\xi=1,\ldots,\varepsilon,
$
the set 
\vspace{0.75ex}
$\Omega_{_0}\subset\Xi^{\;\!\prime}$ such that   
$|w(t,x)|\ne 0$ for all $(t,x)\in\Omega_{_0}$ and
$|w(t,x)|= 0$ for all $(t,x)\in {\sf C}_{_{\Xi^{\;\!\prime}}}\Omega_{_0}.$
}
\vspace{1.5ex}

Using Definition 7.2, we get
\vspace{1.5ex}

{\bf Proposition 7.1.}\!
\vspace{0.75ex}
{\it
If  
$\!\Bigl((w, W), 
\Bigl(h_\xi^{},z_{_{\scriptstyle h_\xi^{}  f_\xi^{}}}, Q_{h_\xi^{}  f_\xi^{}}^{}\Bigr)\Bigr)
\!\in\! \text{\rm G}_{_{\Xi^{\;\!\prime}}},\,
f_{\xi}^{}\!=\!1,\ldots, \delta_{\xi}^{},\, \xi\!=\!1,\ldots,\varepsilon,\!$
then the complex-valued polynomial partial integral $w$ with cofactor $W$
\vspace{0.5ex}
on the domain $\Xi^{\;\!\prime}$ of system {\rm (0.1)} is 
multiple with multiplicity
$
\varkappa =1 + \sum\limits_{\xi=1}^{\varepsilon} \delta_{\xi}^{}.
$
}
\vspace{1.25ex}

{\bf Property 7.1.}
{\it
Let 
$k\in\N,\ c\in\C\backslash\{0\},\ 
\bigl(w, w^kW_{0}^{}\bigr)\in \text{\rm H}_{_{\Xi^{\;\!\prime}}}.
$ 
Then we claim that}
\\[2.25ex]
\mbox{}\hfill                           
$
\bigl(\bigl(w, w^kW_{0}^{}\bigr), \bigl(k, c,{}-k\;\!c\;\!W_0^{}\bigr)\bigr)
\in \text{\rm G}_{_{\Xi^{\;\!\prime}}}.
\hfill
$
\\[2ex]
\indent
{\sl Proof}. 
From Definition 6.1 it follows that
\\[2ex]
\mbox{}\hfill                           
$
{\frak d}\;\!w(t,x)=w^{k+1}(t,x)\;\!W_0^{}(t,x)
$
\ for all $(t,x)\in \Xi^{\;\!\prime}.
\hfill
$
\\[2ex]
\indent
Then the derivative by virtue of system (0.1)
\\[2ex]
\mbox{}\hfill                           
$
{\frak d}\;\!\dfrac{c}{w^{k}(t,x)}=
{}-k\;\!c\,
\dfrac{{\frak d}\;\!w(t,x)}{w^{k+1}(t,x)}=
{}-k\;\!c\;\!W_0^{}(t,x)
$
\ for all 
$(t,x)\in \Omega_{{}_0}^{},
\hfill
$
\\[2.5ex]
where the set 
\vspace{0.5ex}
$\Omega_{_0}\subset\Xi^{\;\!\prime}$ such that   
$|w(t,x)|\ne 0$ for all $(t,x)\in\Omega_{_0}$ and 
$|w(t,x)|= 0$ for all $(t,x)\in {\sf C}_{_{\Xi^{\;\!\prime}}}\Omega_{_0}.$
\vspace{1ex}

By Definition 7.1, 
$
\bigl(\bigl(w, w^kW_{0}^{}\bigr), \bigl(k, c,{}-k\;\!c\;\!W_0^{}\bigr)\bigr)
\in \text{\rm G}_{_{\Xi^{\;\!\prime}}}.\ \k
$
\vspace{1.5ex}

{\bf Property 7.2.}
\vspace{0.75ex}
{\it
Suppose  
$k\in\N,\ 
c_l^{}\in\C\backslash\{0\},\ l=1,\ldots,k,$ and
$\bigl(w, w^kW_{0}^{}\bigr)\in \text{\rm H}_{_{\Xi^{\;\!\prime}}}.$ 
\linebreak
Then  the complex-valued polynomial partial integral $w$ 
\vspace{0.5ex}
on the domain $\Xi^{\;\!\prime}$ of system {\rm (0.1)} 
is $(k+1)\!$-multiple and}
\\[1.75ex]
\mbox{}\hfill                           
$
\bigl(\bigl(w, w^kW_{0}^{}\bigr), \bigl(l, c_l^{},{}-l\;\!c_l^{}\;\!w^{\;\!k-l}W_0^{}\bigr)\bigr)
\in \text{\rm G}_{_{\Xi^{\;\!\prime}}},
\ \ 
l=1,\ldots, k.
\hfill
$
\\[2.75ex]
\indent
{\sl Proof.}
If  
$\bigl(w, w^kW_{0}^{}\bigr)\in \text{\rm H}_{_{\Xi^{\;\!\prime}}},$ then
$\bigl(w, w^{\;\!l}W_{l}^{}\bigr)\in \text{\rm H}_{_{\Xi^{\;\!\prime}}},$
where the functions 
\\[2.5ex]
\mbox{}\hfill                           
$
W_{l}^{}(t,x)=w^{\;\!k-l}(t,x)\;\!W_0^{}(t,x)
$
\ for all 
$(t,x)\in \Xi^{\;\!\prime},
\quad 
l=1,\ldots, k.
\hfill
$
\\[2.25ex]
\indent
By Property 7.1, 
\\[1.5ex]
\mbox{}\hfill                           
$
\bigl(\bigl(w, w^{\;\!l}W_{l}^{}\bigr), \bigl(l, c_l^{},{}-l\;\!c_l^{}\;\!W_l^{}\bigr)\bigr)
\in \text{\rm G}_{_{\Xi^{\;\!\prime}}},
\ \ c_l^{}\in\C\backslash\{0\},
\ \ l=1,\ldots, k.
\hfill
$
\\[2ex]
\indent
So, we obtain 
\\[1.75ex]
\mbox{}\hfill                           
$
\bigl(\bigl(w, w^{\;\!k}W_{0}^{}\bigr), \bigl(l, c_l^{},{}-l\;\!c_l^{}\;\!w^{\;\!k-l}W_0^{}\bigr)\bigr)
\in \text{\rm G}_{_{\Xi^{\;\!\prime}}},\
c_l^{}\in\C\backslash\{0\},\ l=1,\ldots,k.\ \k
\hfill
$
\\[-2.25ex]

\newpage

\mbox{}
\\[-2.5ex]
\centerline{
{\bf\large \S\;\!2. Last multipliers}}
\\[1.75ex]
\centerline{
{\bf  8. 
Last multiplier as partial integral}
}
\\[1.5ex]
\indent
{\bf Theorem 8.1}
(criterion of last multiplier).
{\it 
A continuously differentiable function is a last multiplier on the domain $\Omega$ 
of system {\rm (0.1)} 
if and only if
this function is a partial integral with the factor $-\,{\rm div}\;\!{\frak d}$ 
on the domain $\Omega$ of system {\rm (0.1)}.
}
\vspace{0.35ex}

{\sl Proof}.
\vspace{0.25ex}
Since the functions
$X_i^{}\in \text{\rm P}_{_{\!\Xi}},\ i=1,\ldots,n,$ we see that  
\\[1ex]
\mbox{}\hfill
$
{\rm div}\;\!{\frak d}\in \text{\rm P}_{_{\!\Xi}}
$
\ \ and \ \ 
$
\deg_{\;\!x}^{}\;\! {\rm div}\;\!{\frak d}\leq d-1.
\hfill
$ 
\\[1ex]
\indent
If ${\rm g}=\mu$ and $M=-\,{\rm div}\;\!{\frak d},$ then
the identity (1.2) from the existence criterion of partial integral (Theorem 1.1) and 
the identity (0.6) from the existence criterion of last multiplier (Theorem 0.5) are equivalent.

Therefore,
\\[0.35ex]
\mbox{}\hfill
$
\mu\in \text{\rm M}_{_{\Omega}}
\iff
(\mu, -\,{\rm div}\;\!{\frak d})\in \text{\rm J}_{_{\Omega}}.
\ \k
\hfill
$
\\[1.75ex]
\indent
By Theorem 8.1, we have
$
\text{\rm M}_{_{\Omega}}\subset\text{\rm J}_{_{\Omega}}.
$
\vspace{0.5ex}

If a last multiplier of system (0.1) is
a) polynomial partial integral;
b) multiple polynomial  partial integral;
c) exponential partial integral;
d) conditional  partial integral;
e) complex-valued polynomial  partial integral;
f) multiple complex-valued polynomial  partial integral,
then we say that this last multiplier is 
a) polynomial;
b) multiple polynomial;
c) exponential;
d) conditional;
e) complex-valued polynomial;
f) multiple complex-valued polynomial.
\vspace{0.35ex}

Now we introduce the following notations:
\vspace{0.25ex}

a) $\text{\rm MA}_{_{\Xi^{\;\!\prime}}}$ is
the set of polynomial last multipliers on the domain $\Xi^{\;\!\prime}$ for system (0.1);
\vspace{0.5ex}

b) $\text{\rm MB}_{_{\Xi^{\;\!\prime}}}$ is 
the set of multiple polynomial last multipliers on $\Xi^{\;\!\prime}$ for system (0.1);
\vspace{0.35ex}

c) $\text{\rm ME}_{_{\Omega}}$ is
the set of exponential last multipliers on the domain $\Omega$ for system (0.1); 
\vspace{0.5ex}

d) $\text{\rm MF}_{_{\!\Xi^{\;\!\prime}}}$ is
the set of  conditional last multipliers on the domain $\Xi^{\;\!\prime}$ for system (0.1); 
\vspace{0.5ex}

e) $\text{\rm MH}_{_{\Xi^{\;\!\prime}}}$ is
\vspace{0.35ex}
the set of complex-valued polynomial last multipliers on the domain $\Xi^{\;\!\prime}$ for system (0.1); 
\vspace{0.5ex}

f) $\text{\rm MG}_{_{\Xi^{\;\!\prime}}}$ is
\vspace{0.35ex}
the set of multiple complex-valued polynomial last multipliers on the domain $\Xi^{\;\!\prime}$ for system (0.1).
\vspace{0.75ex}

{\bf Theorem 8.2}
(geometric meaning of last multiplier).
{\it 
If a  last multiplier $\mu$ of the dif\-fe\-ren\-ti\-al system {\rm (0.1)} determines the manifold $\mu(t,x)=0,$ 
then this manifold is an integral manifold of the differential system {\rm (0.1)}.
}
\vspace{0.25ex}

{\sl Proof}\;\!  follows from the definition of integral manifold (Definition 0.3)
and the definition of last multiplier  (Definition 0.4). $\k$
\vspace{0.5ex}

{\bf Theorem 8.3.}
\vspace{0.25ex}
{\it
The continuously differentiable function {\rm (1.4)}
is a last multiplier on the domain $\Omega$ of system {\rm (0.1)}
if and only if 
\vspace{0.25ex}
there exist functions $M_j^{}\in C^1\Omega,\ j=1,\ldots,m,$ 
such that these functions satisfies the identities {\rm (1.5)} and}
\\[1ex]
\mbox{}\hfill
$
\displaystyle
\sum\limits_{j=1}^{m} M_j^{}(t,x)= -\,{\rm div}\;\!{\frak d}(t,x)
$
\ \, \text{\it for all}\  \,
$
(t,x)\in\Omega.
\hfill
$
\\[1.5ex]
\indent
{\sl Proof}  follows from Theorem 1.3 and Theorem 8.1. $\k$
\vspace{0.75ex}

{\bf  Property 8.1.}
\vspace{0.25ex}
{\it
Suppose a function $\varphi\in C^1T^{\;\!\prime},$ 
a set $T_{_0}\subset T^{\;\!\prime}$ 
such that 
$\varphi(t)\ne 0$ for all $t\in T_{_0}$ and 
$\varphi(t)=0$ for all $t\in {\sf C}_{{}_{T^{\;\!\prime}}}T_{_0},$
the set $\Omega_{_0}=T_{_0}\times X^{\;\!\prime}.$ Then}
\\[1.75ex]
\mbox{}\hfill
$
\varphi\;\! {\rm g}\in \text{\rm M}_{_{\Omega}}
\iff 
\bigl({\rm g}\;\!,\;\! {}-{\sf D} \ln |\varphi| -\,{\rm div}\;\!{\frak d}\bigr)\in 
\text{J}_{_{\Omega_{{}_{\tiny\;\! 0}}}}.
\hfill
$
\\[1.75ex]
\indent
{\sl Proof}\,  follows from Theorem 8.1 and Property 1.1 with 
\\[1.5ex]
\mbox{}\hfill                        
$
M(t,x)+ {\sf D} \ln |\varphi(t) | =
-\,{\rm div}\;\!{\frak d}(t,x)
$
\ for all\  $(t,x)\in  \Omega_{_0}.\ \k
\hfill 
$
\\[1.5ex]
\indent
As a consequence of Property 8.1, we obtain the following
\vspace{0.5ex}

{\bf Property 8.2.} 
{\it
If $\lambda\in\R\backslash\{0\},$ then}
\\[1.5ex]
\mbox{}\hfill
$
\mu\in \text{\rm M}_{_{\Omega}}
\iff 
\lambda\;\!\mu\in \text{\rm M}_{_{\Omega}}.
\hfill
$
\\[1.75ex]
\indent
By Property 8.2, if we speaking about two or more last multipliers of
system (0.1), then we assume that these last multipliers are pairwise linearly independent.
\vspace{1.25ex}

{\bf Property 8.3.} 
\vspace{0.25ex}
{\it
Suppose a set $\Omega_{_0}\subset\Omega$ such that  
$\mu(t,x)\ne 0$ for all $(t,x)\!\in\! \Omega_{_0}$ and
$\mu(t,x)\!=\!0$ for all $(t,x)\in {\sf C}_{_\Omega}\Omega_{_0}.$
Then}
\\[1.25ex]
\mbox{}\hfill  
$
\mu\in \text{\rm M}_{_{\Omega}}
\iff 
|\mu|\in \text{\rm M}_{_{\Omega_{{}_{\tiny\;\! 0}}}}.
\hfill
$
\\[1.75ex]
\indent
{\sl Proof}  follows from Property 1.3 and Theorem 8.1. $\k$
\vspace{1.25ex}

{\bf Property 8.4.} 
\vspace{0.5ex}
{\it
Suppose
$\lambda_j^{}\in\R\backslash\{0\},\ j=1,\ldots,m.$ Then} 
\\[1.25ex]
\mbox{}\hfill  
$
\displaystyle
\mu_j^{}\in \text{\rm M}_{_{\Omega}}, \
j=1,\ldots, m,
\ \ \Longrightarrow\ \
\sum\limits_{j=1}^{m} \lambda_j^{}\;\!\mu_j^{}\in \text{\rm M}_{_{\Omega}}\;\!.
\hfill
$
\\[1.75ex]
\indent
{\sl Proof}\;\!  follows from Property 1.5 and Theorem 8.1. $\k$
\vspace{1.25ex}

{\bf Property 8.5.} 
{\it
Suppose $\gamma\in\R\backslash\{0\},\ f^{\;\!\gamma}\in C^1\Omega.$
Then
\\[1.5ex]
\mbox{}\hfill  
$
f^{\;\!\gamma}\in \text{\rm M}_{_{\Omega}}
\iff 
\Bigl( f,\;\! {}-\dfrac{1}{\gamma}\ {\rm div}\;\!{\frak d}\Bigr)\in \text{\rm J}_{_{\Omega}}
\hfill 
$
\\[1ex]
and}
\\[1ex]
\mbox{}\hfill  
$
f\in \text{\rm M}_{_{\Omega}}
\iff 
\bigl( f^{\;\!\gamma},\;\! {}-\gamma\, {\rm div}\;\!{\frak d}\bigr)\in \text{\rm J}_{_{\Omega}}\;\!.
\hfill 
$
\\[1.75ex]
\indent
{\sl Proof}  follows from Property 1.6 and Theorem 8.1. $\k$
\vspace{1.5ex}

{\bf Property 8.6.} 
\vspace{0.5ex}
{\it
Suppose $\gamma\in\R\backslash\{0\},$ 
a set $\Omega_{_0}\!\subset\Omega$ such that   
$f(t,x)\ne 0$ for all $(t,x)\!\in\! \Omega_{_0}$ and 
$f(t,x)=0$ for all $(t,x)\in {\sf C}_{_\Omega}\Omega_{_0}.$
Then 
\\[1.5ex]
\mbox{}\hfill  
$
|f|^{{}^{\scriptstyle \gamma}}\in \text{\rm M}_{_{\Omega_{{}_{\tiny\;\! 0}}}}
\iff 
\Bigl( f,\;\! {}-\dfrac{1}{\gamma}\ {\rm div}\;\!{\frak d}\Bigr)\in \text{\rm J}_{_{\Omega}}
\hfill 
$
\\[1ex]
and}
\\[1ex]
\mbox{}\hfill  
$
f\in \text{\rm M}_{_{\Omega}}
\iff 
\bigl(\;\! |f|^{{}^{\scriptstyle \gamma}},\;\! {}-\gamma\, {\rm div}\;\!{\frak d}\bigr)\in 
\text{\rm J}_{_{\Omega_{{}_{\tiny\;\! 0}}}}.
\hfill 
$
\\[1.75ex]
\indent
{\sl Proof}  follows from Property 1.7 and Theorem 8.1. $\k$
\vspace{1.25ex}

{\bf Property 8.7.} 
\vspace{0.5ex}
{\it
Suppose  $\rho_j^{},\;\! \lambda_j^{}\in\R\backslash\{0\},\ 
{\rm g}_j^{{}^{\scriptstyle 1/\rho_{\!j}^{}}}\in C^1\Omega, ,\ j=1,\ldots,m.$ Then}
\\[1.25ex]
\mbox{}\hfill  
$
\displaystyle
\bigl({\rm g}_j^{}\;\!,\;\! \rho_j^{}\, {\rm div}\;\!{\frak d}\bigr)\in 
\text{\rm J}_{_{\Omega}}, \
j=1,\ldots, m,
\ \ \Longrightarrow\ \
\sum\limits_{j=1}^{m} \lambda_j^{}\;\!
{\rm g}_j^{{}^{\scriptstyle {}-1/\rho_{\!j}^{}}}
\in \text{\rm M}_{_{\Omega}}\;\!.
\hfill
$
\\[1.75ex]
\indent
{\sl Proof}  follows from Property 1.8 and Theorem 8.1. $\k$
\vspace{1.25ex}

{\bf Property 8.8.} 
\vspace{0.5ex}
{\it
Suppose
$({\rm g}_j^{}\;\!,\;\! M_{j}^{})\in \text{\rm J}_{_{\Omega}},\ 
\gamma_j^{}\in\R\backslash\{0\},\
{\rm g}_j^{{}^{\scriptsize \gamma_{j}^{}}}\in C^1\Omega, \ j=1,\ldots,m.$ 
Then  
\\[1.75ex]
\mbox{}\hfill
$
\displaystyle
\prod\limits_{j=1}^{m} 
{\rm g}_j^{{}^{\scriptsize \gamma_{j}^{}}}
\in \text{\rm M}_{_{\Omega}}
\hfill
$
\\[1.75ex]
if and only if the linear combination of cofactors}
\\[1.75ex]
\mbox{}\hfill                      % (8.1)
$
\displaystyle
\sum\limits_{j=1}^{m}
\gamma_j^{}\;\!M_{j}^{}(t,x)=
{}-\,{\rm div}\;\!{\frak d}(t,x)
\quad
\text{for all}\ \ (t,x)\in\Xi^{\;\!\prime}.
$
\hfill (8.1)
\\[1.5ex]
\indent
{\sl Proof}\,  follows from Theorem 8.1 and Property 1.9 with  
$M=-\,{\rm div}\;\!{\frak d}.\ \k$ 
\vspace{1.25ex}

\newpage

{\bf Property 8.9.} 
\vspace{0.75ex}
{\it
Let 
$({\rm g}_j^{}\;\!,\;\!\rho_{j}^{}\,{\rm div}\;\!{\frak d})\in \text{\rm J}_{_{\Omega}},\ 
\rho_{j}^{},\gamma_j^{}\in\R\backslash\{0\},\
{\rm g}_j^{{}^{\scriptsize \gamma_{j}^{}}}\in C^1\Omega, \ j=1,\ldots,m.$ 
Then   
\\[1.75ex]
\mbox{}\hfill
$
\displaystyle
\prod\limits_{j=1}^{m} 
{\rm g}_j^{{}^{\scriptsize \gamma_{j}^{}}}
\in \text{\rm M}_{_{\Omega}}
\hfill
$
\\[1.5ex]
if and only if}
\\[1.5ex]
\mbox{}\hfill
$
\displaystyle
\sum\limits_{j=1}^{m}\,
\rho_{j}^{}\;\!\gamma_j^{}={}-1.
\hfill
$
\\[1.75ex]
\indent
{\sl Proof}\,  follows from Property 8.8 with  
$M_j^{}=\rho_{j}^{}\,{\rm div}\;\!{\frak d},\ j=1,\ldots,m.\ \k$ 
\vspace{1.25ex}

{\bf Property 8.10.} 
\vspace{0.5ex}
{\it
Suppose 
$\mu_j^{}\in \text{\rm M}_{_{\Omega}}\;\!,\
\mu_j^{{}^{\scriptsize \gamma_{j}^{}}}\in C^1\Omega, \ 
\gamma_j^{}\in\R\backslash\{0\},\ j=1,\ldots,m.$ 
Then  
\\[1.5ex]
\mbox{}\hfill
$
\displaystyle
\biggl(\,\prod\limits_{j=1}^{m} 
\mu_j^{{}^{\scriptsize \gamma_{j}^{}}}, M\biggl) 
\in \text{\rm J}_{_{\Omega}}
\hfill
$
\\[1.5ex]
if and only if the cofactor}
\\[1.5ex]
\mbox{}\hfill                      % (8.2)
$
\displaystyle
M(t,x)={}-\sum\limits_{j=1}^{m}
\gamma_j^{}\, {\rm div}\;\!{\frak d}(t,x)
\quad
\text{\it for all}\ \  (t,x)\in\Xi^{\;\!\prime}.
$
\hfill (8.2)
\\[1.5ex]
\indent
{\sl Proof}\,  follows from Theorem 8.1 and Property 1.9 with
$M_j^{}={}-\,{\rm div}\;\!{\frak d},\ j=1,\ldots,m.\ \k$ 
\vspace{1.25ex}

{\bf Property 8.11.} 
\vspace{0.5ex}
{\it
Suppose  
$\mu_j^{}\in \text{\rm M}_{_{\Omega}}\;\!,\
\mu_j^{{}^{\scriptsize \gamma_{j}^{}}}\in C^1\Omega, \ 
\gamma_j^{}\in\R\backslash\{0\},\ j=1,\ldots,m.$ 
Then  
\\[1.5ex]
\mbox{}\hfill                      
$
\displaystyle
\prod\limits_{j=1}^{m} 
\mu_j^{{}^{\scriptsize \gamma_{j}^{}}}
\in \text{\rm M}_{_{\Omega}}
\hfill
$
\\[1ex]
if and only if}
\\[1ex]
\mbox{}\hfill                      
$
\displaystyle
\sum\limits_{j=1}^{m}\gamma_j^{}=1.
\hfill
$
\\[1.75ex]
\indent
{\sl Proof}\,  follows from Theorem 8.1 and Property 8.10 with
$M={}-\,{\rm div}\;\!{\frak d}.\ \k$ 
\vspace{1.75ex}

{\bf Property 8.12.} 
\vspace{0.5ex}
{\it
Suppose  
$\mu_j^{}\in \text{\rm M}_{_{\Omega}}\;\!,\
{\rm g}^{{}^{\scriptsize \gamma}},\mu_j^{{}^{\scriptsize \gamma_{j}^{}}}\in C^1\Omega, \ 
\gamma,\gamma_j^{}\in\R\backslash\{0\},\ j=1,\ldots,m.$ 
Then}  
\\[1ex]
\mbox{}\hfill                      
$
\displaystyle
\biggl(\,{\rm g}^{{}^{\scriptsize \gamma}}\prod\limits_{j=1}^{m} 
\mu_j^{{}^{\scriptsize \gamma_{j}^{}}}, M\biggl) 
\in \text{\rm J}_{_{\Omega}}
\iff
\biggl({\rm g}, \,
\dfrac{1}{\gamma}\;\!
\biggl(M+\sum\limits_{j=1}^{m}
\gamma_{j}^{}\;\!{\rm div}\;\!{\frak d}\biggl)\biggl) 
\in \text{\rm J}_{_{\Omega}}\;\!.
\hfill
$
\\[1.5ex]
\indent
{\sl Proof}\,  follows from Property 1.10 and Theorem 8.1. $\k$
\vspace{1.5ex}

{\bf Property 8.13.} 
\vspace{0.5ex}
{\it
Suppose  
$\mu_j^{}\in \text{\rm M}_{_{\Omega}}\;\!,\
{\rm g}^{{}^{\scriptsize \gamma}},\mu_j^{{}^{\scriptsize \gamma_{j}^{}}}\in C^1\Omega, \ 
\gamma,\gamma_j^{}\in\R\backslash\{0\},\ j=1,\ldots,m.$ 
Then}  
\\[1ex]
\mbox{}\hfill                      
$
\displaystyle
{\rm g}^{{}^{\scriptsize \gamma}}\prod\limits_{j=1}^{m} 
\mu_j^{{}^{\scriptsize \gamma_{j}^{}}} 
\in \text{\rm M}_{_{\Omega}}
\iff
\biggl({\rm g}, \,
\dfrac{1}{\gamma}\;\!
\biggl(\,\sum\limits_{j=1}^{m}
\gamma_{j}^{}-1\biggr)\, {\rm div}\;\!{\frak d}\biggl) 
\in \text{\rm J}_{_{\Omega}}\;\!.
\hfill
$
\\[1.5ex]
\indent
{\sl Proof}\,  follows from Theorem 8.1 and Property 8.12 with
$M={}-\,{\rm div}\;\!{\frak d}.\ \k$ 
\vspace{1.5ex}

{\bf Property 8.14.} 
\vspace{0.75ex}
{\it
Suppose  
$\mu_j^{}\in \text{\rm M}_{_{\Omega}}\;\!,\
\mu_j^{{}^{\scriptsize \gamma_{j}^{}}}\in C^1\Omega, \ 
\gamma_j^{}\in\R\backslash\{0\},\ j=1,\ldots,m,\ 
{\rm g}_\tau^{{}^{\scriptsize \xi_{\tau}^{}}}\!\in C^1\Omega,$
$\xi_\tau^{}\in\R\backslash\{0\},\ \tau=1,\ldots,l.$ 
Then  we claim that
\\[1ex]
\mbox{}\hfill                      
$
\displaystyle
\biggl(\,\prod\limits_{\tau=1}^{l} 
{\rm g}_\tau^{{}^{\scriptsize \xi_{\tau}^{}}}
\,\prod\limits_{j=1}^{m} 
\mu_j^{{}^{\scriptsize \gamma_{j}^{}}}, M\biggl)\, 
\in \text{\rm J}_{_{\Omega}}
\iff
\biggl(\,\prod\limits_{\tau=1}^{l} 
{\rm g}_\tau^{{}^{\scriptsize \xi_{\tau}^{}}}, \,
M+\sum\limits_{j=1}^{m}
\gamma_{j}^{}\, {\rm div}\;\!{\frak d}\biggl) \,
\in \text{\rm J}_{_{\Omega}}\;\!.
\hfill
$
\\[1.5ex]
Moreover, there exist functions
$M_\tau^{}\in C^1\Omega,\ \tau=1,\ldots, l,$
such that 

\newpage

\mbox{}
\\[-1.75ex]
\mbox{}\hfill                        % (8.3)
$
{\frak d}\;\!{\rm g}_\tau^{}(t,x)=
{\rm g}_\tau^{}(t,x)\;\!M_\tau^{}(t,x)
\quad
\text{for all}\ \  (t,x)\in\Omega,
\quad
\tau=1,\ldots, l,
$
\hfill {\rm (8.3)}
\\[1ex]
and}
\\[1ex]
\mbox{}\hfill                     
$
\displaystyle
\sum\limits_{\tau=1}^{l}
\xi_{\tau}^{}\;\!M_\tau^{}(t,x)=
M(t,x)+\sum\limits_{j=1}^{m}
\gamma_{j}^{}\, {\rm div}\;\!{\frak d}(t,x)
\quad
\text{\it for all}\ \  (t,x)\in\Omega.
\hfill
$
\\[1.5ex]
\indent
{\sl Proof}\,  follows from Property 1.11 and Theorem 8.1. $\k$
\vspace{1.5ex}

{\bf Property 8.15.} 
\vspace{0.75ex}
{\it
Suppose  
$\mu_j^{}\in \text{\rm M}_{_{\Omega}}\;\!,\
\mu_j^{{}^{\scriptsize \gamma_{j}^{}}}\in C^1\Omega, \ 
\gamma_j^{}\in\R\backslash\{0\},\ j=1,\ldots,m,\ 
{\rm g}_\tau^{{}^{\scriptsize \xi_{\tau}^{}}}\!\in C^1\Omega,$
$\xi_\tau^{}\in\R\backslash\{0\},\ \tau=1,\ldots,l.$ 
Then we have   
\\[1.5ex]
\mbox{}\hfill                      
$
\displaystyle
\prod\limits_{\tau=1}^{l} 
{\rm g}_\tau^{{}^{\scriptsize \xi_{\tau}^{}}}
\prod\limits_{j=1}^{m} 
\mu_j^{{}^{\scriptsize \gamma_{j}^{}}} 
\in \text{\rm M}_{_{\Omega}}
\iff
\biggl(\,\prod\limits_{\tau=1}^{l} 
{\rm g}_\tau^{{}^{\scriptsize \xi_{\tau}^{}}}, \,
\biggl(\,\sum\limits_{j=1}^{m}\gamma_{j}^{}-1\biggr)\, {\rm div}\;\!{\frak d}\biggl) \,
\in \text{\rm J}_{_{\Omega}}\;\!.
\hfill
$
\\[1.5ex]
Moreover, there exist functions
$M_\tau^{}\in C^1\Omega,\ \tau=1,\ldots, l,$
such that these functions satisfies the identities {\rm (8.3)} and}
\\[1ex]
\mbox{}\hfill                     
$
\displaystyle
\sum\limits_{\tau=1}^{l}
\xi_{\tau}^{}\;\!M_\tau^{}(t,x)=
\biggl(\,\sum\limits_{j=1}^{m}\gamma_{j}^{}-1\biggr)\, {\rm div}\;\!{\frak d}(t,x)
\quad
\text{\it for all}\ \  (t,x)\in\Omega.
\hfill
$
\\[1.5ex]
\indent
{\sl Proof}\,  follows from Theorem 8.1 and Property 8.14 with 
\vspace{1.5ex}
$M={}-\,{\rm div}\;\!{\frak d}.\ \k$

{\bf Remark 8.1.}
\vspace{1ex}
If necessary, it is advisable to introduce the concept of last pseudomultiplier.

{\bf Definition 8.1.}
{\it 
A function $\nu\in C^1\Omega$ is called a
\textit{\textbf{last pseudomultiplier with coefficient}} $\rho$
or \textit{\textbf{last $\rho\!$-pseudomultiplier $(\rho\in\R)$ 
on the domain}} $\Omega$ of system {\rm (0.1)} if 
the derivative by virtue of system {\rm(0.1)} 
\\[2ex]
\mbox{}\hfill                        
$
{\frak d}\;\!\nu(t,x)=\rho\;\!\nu(t,x)\;{\rm div}\, {\frak d}(t,x)
$ 
\ for all} $(t,x)\in \Omega.
\hfill
$
\\[2.25ex]
\indent
From the notion of partial integral (Definition 1.1 and Theorem 1.1), 
we get the following criterion of last pseudomultiplier.
\vspace{0.75ex}

{\bf Property 8.16.}
 \vspace{0.35ex}
 {\it 
A function $\nu\in C^1\Omega$ is a  last $\rho\!$-pseudomultiplier 
on the domain $\Omega$ of system {\rm (0.1)} if and only if 
the function $\nu$ is a partial integral with cofactor 
$\rho\;\!{\rm div}\, {\frak d}$ 
on the domain $\Omega$ of system {\rm (0.1)}. 
}
\vspace{0.75ex}

Using Property 8.5 and taking into account Definition 8.1 and Theorem 0.5, 
we obtain the relationship between last multiplier and last pseudomultiplier.
 \vspace{0.75ex}

{\bf Property 8.17A.}
 \vspace{0.35ex}
 {\it 
 A function $\nu\in C^1\Omega$ is a last pseudomultiplier with coefficient 
 $\rho\ne 0$ on the domain $\Omega$ of system {\rm (0.1)}
if and only if the function
\\[1.75ex]
\mbox{}\hfill
$
\mu\colon (t,x)\to\ \nu^{\,{}-\;\!1/\rho}(t,x)
$
\ for all $(t,x)\in\Omega
\hfill
$
\\[1.5ex]
is a last multiplier on the domain $\Omega$ of system {\rm (0.1)}. 
}
\vspace{1.25ex}

{\bf Property 8.17B.}
 \vspace{0.35ex}
 {\it 
 A function $\mu\in C^1\Omega$ is a last multiplier on the domain 
$\Omega$ of system {\rm (0.1)} if and only if the function
\\[1.5ex]
\mbox{}\hfill
$
\nu\colon (t,x)\to\ \mu^{\,{}-\;\!\rho}(t,x)
$
\ for all 
$
(t,x)\in\Omega
\hfill
$
\\[1.75ex]
is a last pseudomultiplier with coefficient $\rho\ne 0$ 
on the domain $\Omega$ of system {\rm (0.1)}. 
}
\vspace{0.75ex}

Thus, 
\vspace{0.5ex}
along with last multipliers we can consider last pseudomultipliers of system (0.1).

Finally note that the last $0\!$-pseudomultiplier of system (0.1) is a first integral.

\newpage

\mbox{}
\\[-1.75ex]
\centerline{\bf  9. Building of last multipliers on base of polynomial partial integrals} 
\\[1.75ex]
\indent
{\bf Property 9.1.} 
\vspace{0.5ex}
{\it
Suppose 
$p_j^{}\in \text{\rm P}_{_{\!\Xi^{\;\!\prime}}},\ 
\gamma_j^{}\in\R\backslash\{0\},$ 
the set $\Omega_{_0}\subset\Xi^{\;\!\prime}$ such that  
$p_j^{{}^{\scriptsize \gamma_j^{}}}\in C^1\Omega_{_0},$ $j=1,\ldots, m.$
Then we claim that
\\[1.5ex]
\mbox{}\hfill  
$
\displaystyle
\prod\limits_{j=1}^{m} 
p_j^{{}^{\scriptsize \gamma_{j}^{}}}
\in \text{\rm M}_{_{\Omega_{{}_{\tiny\;\! 0}}}}
\iff
\bigl(p_j^{}, M_j^{}\bigr)\in \text{\rm A}_{_{\Xi^{\;\!\prime}}},
\ \ 
j=1,\ldots, m,
\hfill
$
\\[1.5ex]
where the cofactors $M_j^{},\ j=1,\ldots,m,$
\vspace{0.5ex}
such that the identity {\rm (8.1)} holds.}

{\sl Proof}.\;\!
Taking into account Theorem 8.1, from Theorem 2.2 
under the condition 
\linebreak
$M={}-\,{\rm div}\;\!{\frak d},$ 
we get the statement of this property. $\k$ 
\vspace{1ex}

{\bf Property 9.2}
\vspace{0.5ex}
(existence criterion of rational last multiplier).
{\it
Suppose 
%functions
$p_1^{},p_2^{}\in \text{\rm P}_{_{\!\Xi^{\;\!\prime}}}$
are relatively prime, 
\vspace{0.5ex}
the set $\Omega_{_0}\subset\Xi^{\;\!\prime}$ such that   
$p_2^{}(t,x)\ne 0$ for all $(t,x)\in \Omega_{_0}$ and
$p_2^{}(t,x)= 0$ for all $(t,x)\in {\sf C}_{_{\Xi^{\;\!\prime}}}\Omega_{_0}.$
Then we have
\\[1.5ex]
\mbox{}\hfill  
$
\dfrac{p_1^{}}{p_2^{}} 
\in \text{\rm M}_{_{\Omega_{{}_{\tiny\;\! 0}}}}
\iff
\bigl(p_1^{}, M_1^{}\bigr)\in \text{\rm A}_{_{\Xi^{\;\!\prime}}}
\ \& \ 
\bigl(p_2^{}, M_2^{}\bigr)\in \text{\rm A}_{_{\Xi^{\;\!\prime}}},
\hfill
$
\\[1.75ex]
where the cofactors $M_1^{}$ and $M_2^{}$ such that
\\[1.75ex]
\mbox{}\hfill  
$
M_2^{}(t,x)-M_1^{}(t,x)={\rm div}\;\!{\frak d}(t,x)
$
\ for all} $(t,x)\in \Xi^{\;\!\prime}.
\hfill
$
\\[1.5ex]
\indent
{\sl Proof.}\;\! 
Taking into account Theorem 8.1, from Theorem 2.3
under the condition 
\linebreak
$M={}-\,{\rm div}\;\!{\frak d},$ 
we get the statement of this property. $\k$ 
\vspace{1ex}

{\bf  Property 9.3.}
\vspace{0.5ex}
{\it
Let $\lambda_j^{},\;\! c_j^{}\in\R, \ j=1,\ldots,m,\ m\leq n,\ 
\sum\limits_{j=1}^{m}|\lambda_j^{}|\ne 0.$
Then 
\\[1ex]
\mbox{}\hfill
$
\displaystyle
\sum\limits_{j=1}^{m}\lambda_j^{}(x_j^{}+c_j^{})
\in\text{\rm MA}_{_{\Xi}}
\hfill
$
\\[1.5ex]
if and only if the identity holds
\\[1.5ex]
\mbox{}\hfill
$
\displaystyle
\sum\limits_{j=1}^{m}
\lambda_j^{}\;\!X_j^{}(t,x)=
{}-\sum\limits_{j=1}^{m}\lambda_j^{}(x_j^{}+c_j^{})\;\!{\rm div}\;\!{\frak d}(t,x)
$
\ for all} $(t,x)\in \Xi.
\hfill
$
\\[1.5ex]
\indent
{\sl Proof.}\;\! 
Using Theorem 8.1, from Property 2.2 under the condition $M={}-\,{\rm div}\;\!{\frak d},$
we obtain the statement of Property 9.3. $\k$ 
\vspace{1ex}

{\bf  Property 9.4.}
\vspace{0.5ex}
{\it
If $c\in\R,$ then
$
x_i^{}+c \in \text{\rm MA}_{_{\Xi}},\ 
i\in \{1,\ldots, n\},
$
if and only if the function $X_i^{}$ from right side of system {\rm (0.1)} has the form
\\[1.5ex]
\mbox{}\hfill
$
X_i^{}(t,x)=
{}- (x_i^{}+c)\;\!{\rm div}\;\!{\frak d}(t,x)
$
\ for all} $(t,x)\in \Xi.
\hfill
$
\\[1.5ex]
\indent
{\sl Proof}.\;\!
Using Theorem 8.1, from Property 2.3 under the condition $M={}-\,{\rm div}\;\!{\frak d},$
we obtain the statement of Property 9.4. $\k$ 
\vspace{1ex}

%\newpage

{\bf Property 9.5}
\vspace{0.25ex}
(existence criterion of multiple polynomial last multiplier).
{\it
Suppose a fun\-c\-tion $\nu\in \text{\rm MA}_{_{\!\Xi^{\;\!\prime}}},$ an number $h\in\N,$ functions
\vspace{0.5ex}
$\nu, q\in \text{\rm P}_{_{\!\Xi^{\;\!\prime}}}$ are relatively prime, 
a function $N\in \text{\rm P}_{_{\!\Xi^{\;\!\prime}}}$ and has the degree  
\vspace{0.5ex}
$\deg_{x}^{} N\leq d-1,$
the set 
$\Omega_{_0}\subset\Xi^{\;\!\prime}$ such that  
$\nu(t,x)\ne 0$ for all $(t,x)\in\Omega_{_0}$ and 
$\nu(t,x)= 0$ for all $(t,x)\in {\sf C}_{_{\Xi^{\;\!\prime}}}\Omega_{_0}.$
\vspace{0.5ex}
Then $\bigl(\nu, (h,q,N)\bigr)\in \text{\rm MB}_{_{\Xi^{\;\!\prime}}}$
if and only if one of the following conditions is true{\rm:}
\\[1ex]
\indent
$
1)\ {\frak d}\,\dfrac{q(t,x)}{\nu^{\;\!h}(t,x)}=N(t,x)
$
\ for all $(t,x)\in \Omega_{_0};
$
\\[1.25ex]
\indent
$
2)\ {\frak d}\;\!q(t,x)=\nu^{\;\!h}(t,x)\;\!N(t,x)-h\;\!q(t,x)\;\!{\rm div}\;\!{\frak d}(t,x)
$
\ for all} $(t,x)\in \Xi^{\;\!\prime};
$
\\[1.5ex]
\indent
$
3)\ 
\Bigl(\exp\dfrac{q}{\nu^{\;\!h}}\,, N\Bigr)\in 
\text{\rm E}_{_{\Omega_{{}_{\tiny\;\! 0}}}}.
$
\\[1.5ex]
\indent
{\sl Proof.}\;\!
\vspace{0.25ex}
Taking into account Theorem 8.1, from Definition 5.1 and Theorems 5.2, 5.3 under the 
conditions $p=\nu$ and $M={}-\,{\rm div}\;\!{\frak d},$ we get the statement of this property. $\k$ 
\vspace{1ex}

{\bf Property\! 9.6.}\!
{\it
Suppose 
$\!\bigl(\nu_j^{}, (h_j^{},q_j^{},\rho_{j}^{}N_0^{})\bigr)\!\in\! 
\text{\rm MB}_{_{\Xi^{\;\!\prime}}},\, 
\rho_j^{}\!\in\!\R\backslash\{0\},\, 
\lambda_j^{},\gamma_j^{}, c_j^{}\!\in\!\R,\, 
\sum\limits_{j=1}^{m}\!|\lambda_j^{}|\!\ne\! 0,\!$ 
$\nu_j^{\gamma_j^{}}\in C^1\Omega_{_0}, \ j=1,\ldots, m,\
\varphi\in C^1T^{\;\!\prime},$ 
the set 
$\Omega_{_0}\subset\Xi^{\;\!\prime}$ such that   
$\prod\limits_{j=1}^{m}\nu_j^{}(t,x)\ne 0$
for all $(t,x)\in\Omega_{_0}$ and
$\prod\limits_{j=1}^{m}\nu_j^{}(t,x)= 0$
for all $(t,x)\in {\sf C}_{_{\Xi^{\;\!\prime}}}\Omega_{_0}.$
Then we have the statements{\rm:}
\\[1.5ex]
\indent
$
1)\
\displaystyle
\biggl(\ \prod\limits_{j=1}^{m}\nu_j^{\gamma_j^{}}
\sum\limits_{j=1}^{m}\lambda_j^{}
\exp\biggl(\;\!\dfrac{q_j^{}}{\rho_j^{}\;\!\nu_j^{\;\!h_j^{}}}+\varphi\biggr), L\biggr)\in 
\text{\rm J}_{_{\Omega_{{}_{\tiny\;\! 0}}}}
\hfill
$
\\[2ex]
if and only if the cofactor
\\[1.75ex]
\mbox{}\hfill                           
$
\displaystyle
L(t,x)={\sf D}\;\!\varphi(t)+N_0^{}(t,x)-
\sum\limits_{j=1}^{m}\gamma_j^{}\;\!{\rm div}\;\!{\frak d}(t,x)
$
\ for all 
$(t,x)\in \Xi^{\;\!\prime};
\hfill
$
\\[1.5ex]
\indent
$
2)\ 
\displaystyle
\biggl(\ \prod\limits_{j=1}^{m}\nu_j^{\gamma_j^{}}
\sum\limits_{j=1}^{m}\lambda_j^{}
\exp\biggl(\;\!\dfrac{q_j^{}}{\rho_j^{}\;\!\nu_j^{\;\!h_j^{}}}+c_j^{}\biggr), L\biggr)\in 
\text{\rm J}_{_{\Omega_{{}_{\tiny\;\! 0}}}}
\hfill
$
\\[2ex]
if and only if the cofactor
\\[1.75ex]
\mbox{}\hfill                           
$
\displaystyle
L(t,x)=N_0^{}(t,x)-
\sum\limits_{j=1}^{m}\gamma_j^{}\;\!{\rm div}\;\!{\frak d}(t,x)
$ 
\ for all} $(t,x)\in \Xi^{\;\!\prime}.
\hfill
$
\\[1.5ex]
\indent
{\sl Proof.}\;\!
Using Theorem 8.1, from Properties 5.5 and 5.6 under the conditions 
\\[1ex]
\mbox{}\hfill
$p_j^{}=\nu_j^{}\;\!,\ \ M_j^{}={}-\,{\rm div}\;\!{\frak d},\ \, j=1,\ldots,m,
\hfill
$ 
\\[1.5ex]
it follows that the statements of Property 9.6 are true. $\k$
\vspace{1ex}

{\bf Property 9.7.}\!
\vspace{0.75ex}
{\it
Suppose 
$\!\bigl(\nu_j^{}, (h_j^{},q_j^{},N_j^{})\bigr)\!\in \text{\rm MB}_{_{\Xi^{\;\!\prime}}},\, 
\gamma_j^{},\xi_j^{}\!\in\!\R,\, \varphi_j^{}\!\in\! C^1T^{\;\!\prime},\, 
j\!=\!1,\ldots, m,\!$ 
the set
\vspace{0.75ex}
$\!\Omega_{_0}\!\subset\!\Xi^{\;\!\prime}\!$ such that  
$\prod\limits_{j=1}^{m}\nu_j^{}(t,x)\!\ne\! 0$ for all $(t,x)\!\in\!\Omega_{_0}$ and
$\prod\limits_{j=1}^{m}\nu_j^{}(t,x)= 0$
for all $(t,x)\!\in\! {\sf C}_{_{\Xi^{\;\!\prime}}}\Omega_{_0},\!$ 
$\nu_j^{\gamma_j^{}}\in C^1\Omega_{_0}, \ j=1,\ldots, m.$
Then
\\[1.5ex]
\mbox{}\hfill                           
$
\displaystyle
\biggl(\ \prod\limits_{j=1}^{m}\nu_j^{\gamma_j^{}}
\exp\sum\limits_{j=1}^{m}\xi_j^{}
\biggl(\;\!\dfrac{q_j^{}}{\nu_j^{\;\! h_j^{}}}+\varphi_j^{}\biggr), L\biggr)\in 
\text{\rm J}_{_{\Omega_{{}_{\tiny\;\! 0}}}}
\hfill
$
\\[2ex]
if and only if the cofactor
\\[1.75ex]
\mbox{}\hfill                           
$
\displaystyle
L(t,x)=
\sum\limits_{j=1}^{m}
\xi_j^{}\;\!\bigl(N_j^{}(t,x)+{\sf D}\;\!\varphi_j^{}(t)\bigr)- 
\sum\limits_{j=1}^{m}
\gamma_j^{}\;\!{\rm div}\;\!{\frak d}(t,x)
$
\ for all} 
$(t,x)\in \Xi^{\;\!\prime}.
\hfill
$
\\[1.5ex]
\indent
{\sl Proof.}\;\!
Taking into account Theorem 8.1, from Property 5.7 under the conditions 
\\[1ex]
\mbox{}\hfill
$
p_j^{}=\nu_j^{},
\quad 
M_j^{}={}-{\rm div}\;\!{\frak d},
\quad 
j=1,\ldots,m,
\hfill
$
\\[1.5ex]
it follows that the statement of Property 9.7 is true. $\k$ 
\vspace{1ex}

By Definition 5.2, we introduce the notion of multiplicity for multiple polynomial  last multiplier.
\vspace{0.5ex}

{\bf Property 9.8.}
\vspace{0.75ex}
{\it
Suppose 
$\Bigl(\nu, 
\Bigl(h_\xi^{},q_{_{\scriptstyle h_\xi^{}  f_\xi^{}}}, N_{h_\xi^{}  f_\xi^{}}^{}\Bigr)\Bigr)
\in \text{\rm MB}_{_{\Xi^{\;\!\prime}}},\
\gamma,\gamma_{_{\scriptstyle h_\xi^{}  f_\xi^{}}}\in\R,\ 
\varphi_{_{\scriptstyle h_\xi^{}  f_\xi^{}}}\in C^1T^{\;\!\prime},$ 
$f_{\xi}^{}=1,\ldots, \delta_{\xi}^{},\ \xi=1,\ldots,\varepsilon,$
\vspace{1ex}
the set 
$\Omega_{_0}\subset\Xi^{\;\!\prime}$ such that  
$\nu(t,x)\ne 0$ for all $(t,x)\in\Omega_{_0}$ and
$\nu(t,x)= 0$ for all $(t,x)\in {\sf C}_{_{\Xi^{\;\!\prime}}}\Omega_{_0},\ 
\nu^{\gamma}\in C^1\Omega_{_0}.$
Then we claim that
\\[1.5ex]
\mbox{}\hfill                           
$
\displaystyle
\biggl(\ \nu^{\gamma}
\exp\sum\limits_{\xi=1}^{\varepsilon}\sum\limits_{f_{\xi}^{}=1}^{\delta_{\xi}^{}}
\biggl(\;\!\gamma_{_{\scriptstyle h_\xi^{}  f_\xi^{}}}\biggl(\,
\dfrac{q_{_{\scriptstyle h_\xi^{}  f_\xi^{}} }}{\displaystyle  \nu^{\;\!h_\xi^{}}}+
\varphi_{_{\scriptstyle h_\xi^{}  f_\xi^{}}}\biggr)\biggr),\, L\biggr)\in 
\text{\rm J}_{_{\Omega_{{}_{\tiny\;\! 0}}}}
\hfill
$
\\[2ex]
if and only if the cofactor
\\[1.75ex]
\mbox{}\hfill                           
$
\displaystyle
L(t,x)=
\sum\limits_{\xi=1}^{\varepsilon}\sum\limits_{f_{\xi}^{}=1}^{\delta_{\xi}^{}}
\biggl(\gamma_{_{\scriptstyle h_\xi^{}  f_\xi^{}}}\Bigl(
N_{h_\xi^{}  f_\xi^{}}^{} (t,x)+{\sf D}\;\!\varphi_{_{\scriptstyle h_\xi^{}  f_\xi^{}}}(t)\Bigr)\biggr)
- \gamma\;\!{\rm div}\;\!{\frak d}(t,x)
$
\ for all} 
$(t,x)\in \Xi^{\;\!\prime}.
\hfill
$
\\[2.25ex]
\indent
{\sl Proof}.\;\!
\vspace{0.35ex}
Using Theorem 8.1, from Property 5.8 under the conditions 
$M={}-\,{\rm div}\;\!{\frak d}, \ p=\nu,$
we obtain the statement of Property 9.8. $\k$ 
\vspace{1.5ex}

{\bf Property 9.9.}\!
{\it
Suppose 
$(p, M)\!\in\! \text{\rm A}_{_{\!\Xi^{\;\!\prime}}},\!$ an number $h\in\N,$ functions
\vspace{0.75ex}
$p, q\!\in\! \text{\rm P}_{_{\!\Xi^{\;\!\prime}}}\!$ are relatively prime, 
the set 
\vspace{0.75ex}
$\Omega_{_0}\subset\Xi^{\;\!\prime}$ such that   
$p(t,x)\ne 0$ for all $(t,x)\in\Omega_{_0}$ and
$p(t,x)= 0$ for all $(t,x)\in {\sf C}_{_{\Xi^{\;\!\prime}}}\Omega_{_0}.$
Then we have 
\\[1.5ex]
\mbox{}\hfill
$
\bigl((p, M), (h,q, -\;\!{\rm div}\;\!{\frak d})\bigr)\in \text{\rm B}_{_{\Xi^{\;\!\prime}}}
\hfill
$
\\[1.75ex]
if and only if one of the following conditions is true{\rm:}
\\[1.75ex]
\indent
$
1)\ {\frak d}\,\dfrac{q(t,x)}{p^{\;\!h}(t,x)}={}-\,{\rm div}\;\!{\frak d}(t,x)
$
\ for all 
$(t,x)\in \Omega_{_0};
$
\\[2.25ex]
\indent
$
2)\ {\frak d}\;\!q(t,x)=h\;\!q(t,x)\;\!M(t,x)-
p^{\;\!h}(t,x)\;\!{\rm div}\;\!{\frak d}(t,x)
$
\ for all} $(t,x)\in \Xi^{\;\!\prime};
$
\\[2.5ex]
\indent
$
3)\ 
\exp\dfrac{q}{p^{\;\!h}}\in 
\text{\rm ME}_{_{\Omega_{{}_{\tiny\;\! 0}}}}.
$
\\[2.5ex]
\indent
{\sl Proof.}\;\! 
\vspace{0.35ex}
Taking into account Theorem 8.1, from Definition 5.1 and Theorems 5.2, 5.3 
under the condition $N={}-\,{\rm div}\;\!{\frak d},$
we get the statement of Property 9.9 is true. $\k$ 
\vspace{1.5ex}

{\bf Property 9.10.}
\vspace{0.5ex}
{\it
Suppose
$\bigl((p_j^{}, M_j^{}), (h_j^{},q_j^{},\rho_{j}^{}N_0^{})\bigr)\in 
\text{\rm B}_{_{\Xi^{\;\!\prime}}},\ 
\rho_j^{}\in\R\backslash\{0\},\ 
\lambda_j^{},\gamma_j^{}, c_j^{}\in\R,
\linebreak
j=1,\ldots, m,\
\sum\limits_{j=1}^{m}\!|\lambda_j^{}|\!\ne\! 0,\
\varphi\in C^1T^{\;\!\prime},$ 
\vspace{0.5ex}
the set 
$\Omega_{_0}\subset\Xi^{\;\!\prime}$ such that   
$\prod\limits_{j=1}^{m}p_j^{}(t,x)\ne 0$
for all $(t,x)\in\Omega_{_0}$ and
$\prod\limits_{j=1}^{m}p_j^{}(t,x)= 0$
for all $(t,x)\in {\sf C}_{_{\Xi^{\;\!\prime}}}\Omega_{_0},\
p_j^{\gamma_j^{}}\in C^1\Omega_{_0}, \ j=1,\ldots, m.$
Then{\rm:}
\\[1.5ex]
\indent
$
1)\
\displaystyle
\prod\limits_{j=1}^{m}p_j^{\gamma_j^{}}
\sum\limits_{j=1}^{m}\lambda_j^{}
\exp\biggl(\;\!\dfrac{q_j^{}}{\rho_j^{}\;\!p_j^{\;\!h_j^{}}}+\varphi\biggr)\in 
\text{\rm M}_{_{\Omega_{{}_{\tiny\;\! 0}}}}
\hfill
$
\\[2ex]
if and only if
\\[1.75ex]
\mbox{}\hfill                           
$
\displaystyle
{\rm div}\;\!{\frak d}(t,x)={}-{\sf D}\;\!\varphi(t)-N_0^{}(t,x)-
\sum\limits_{j=1}^{m}\gamma_j^{}\;\!M_j^{}(t,x)
$
\ for all $(t,x)\in \Xi^{\;\!\prime};
\hfill
$
\\[1.5ex]
\indent
$
2)\ 
\displaystyle
\prod\limits_{j=1}^{m}p_j^{\gamma_j^{}}
\sum\limits_{j=1}^{m}\lambda_j^{}
\exp\biggl(\;\!\dfrac{q_j^{}}{\rho_j^{}\;\!p_j^{\;\!h_j^{}}}+c_j^{}\biggr)\in 
\text{\rm M}_{_{\Omega_{{}_{\tiny\;\! 0}}}}
\hfill
$
\\[2ex]
if and only if
\\[1.75ex]
\mbox{}\hfill                           
$
\displaystyle
{\rm div}\;\!{\frak d}(t,x)={}-N_0^{}(t,x)-
\sum\limits_{j=1}^{m}\gamma_j^{}\;\!M_j^{}(t,x)
$
\ for all} $(t,x)\in \Xi^{\;\!\prime}.
\hfill
$
\\[1.5ex]
\indent
{\sl Proof.}\;\! 
\vspace{0.25ex}
Taking into account Theorem 8.1, from Properties 5.5 and 5.6 
under the condition $L={}-\,{\rm div}\;\!{\frak d},$
we obtain the statement of Property 9.10. $\k$ 
\vspace{1ex}

\newpage

{\bf Property 9.11.}\!
\vspace{0.75ex}
{\it
Suppose 
$\!\bigl((p_j^{}, M_j^{}), (h_j^{},q_j^{},N_j^{})\bigr)\!\in\! \text{\rm B}_{_{\Xi^{\;\!\prime}}},\, 
\gamma_j^{},\xi_j^{}\!\in\!\R,\, \varphi_j^{}\!\in\! C^1T^{\;\!\prime},\, 
j\!=\!1,\ldots, m,\!$ 
the set
\vspace{0.75ex}
$\Omega_{_0}\!\subset\!\Xi^{\;\!\prime}$ such that 
$\prod\limits_{j=1}^{m}p_j^{}(t,x)\ne 0$ for all $(t,x)\in\Omega_{_0}$ and
$\prod\limits_{j=1}^{m}p_j^{}(t,x)= 0$
for all $(t,x)\!\in\! {\sf C}_{_{\Xi^{\;\!\prime}}}\Omega_{_0},\!$ 
$p_j^{\gamma_j^{}}\in C^1\Omega_{_0}, \ j=1,\ldots, m.$
Then
\\[2ex]
\mbox{}\hfill                           
$
\displaystyle
\prod\limits_{j=1}^{m}p_j^{\gamma_j^{}}
\exp\sum\limits_{j=1}^{m}\xi_j^{}
\biggl(\;\!\dfrac{q_j^{}}{p_j^{\;\! h_j^{}}}+\varphi_j^{}\biggr)\in 
\text{\rm M}_{_{\Omega_{{}_{\tiny\;\! 0}}}}
\hfill
$
\\[2ex]
if and only if
\\[2ex]
\mbox{}\hfill                           
$
\displaystyle
{\rm div}\;\!{\frak d}(t,x)=
{}-\sum\limits_{j=1}^{m}
\Bigl(\gamma_j^{}\;\!M_j^{}(t,x)+
\xi_j^{}\;\!\bigl(N_j^{}(t,x)+{\sf D}\;\!\varphi_j^{}(t)\bigr)\Bigr) 
$
\ for all} $(t,x)\in \Xi^{\;\!\prime}.
\hfill
$
\\[2ex]
\indent
{\sl Proof.}
\vspace{0.25ex}
Using Theorem 8.1, from Property 5.7 under the condition
$L={}-\, {\rm div}\;\!{\frak d},$ we get the statement of this property. $\k$
\vspace{1.5ex}

{\bf Property 9.12.}
\vspace{0.75ex}
{\it
Suppose 
$\!\Bigl((p, M), 
\Bigl(h_\xi^{},q_{_{\scriptstyle h_\xi^{}  f_\xi^{}}}, N_{h_\xi^{}  f_\xi^{}}^{}\Bigr)\Bigr)
\!\in\! \text{\rm B}_{_{\Xi^{\;\!\prime}}},\
\gamma,\gamma_{_{\scriptstyle h_\xi^{}  f_\xi^{}}}\!\in\!\R,\ 
\varphi_{_{\scriptstyle h_\xi^{}  f_\xi^{}}}\!\in\! C^1T^{\;\!\prime},$ 
$f_{\xi}^{}=1,\ldots, \delta_{\xi}^{},\ \xi=1,\ldots,\varepsilon,$
\vspace{1ex}
the set 
$\Omega_{_0}\subset\Xi^{\;\!\prime}$ such that
$p(t,x)\ne 0$ for all $(t,x)\in\Omega_{_0}$ and
$p(t,x)= 0$ for all $(t,x)\in {\sf C}_{_{\Xi^{\;\!\prime}}}\Omega_{_0},\ 
p^{\gamma}\in C^1\Omega_{_0}.$
Then
\\[2ex]
\mbox{}\hfill                           
$
\displaystyle
p^{\gamma}
\exp\sum\limits_{\xi=1}^{\varepsilon}\sum\limits_{f_{\xi}^{}=1}^{\delta_{\xi}^{}}
\biggl(\;\!\gamma_{_{\scriptstyle h_\xi^{}  f_\xi^{}}}\biggl(\,
\dfrac{q_{_{\scriptstyle h_\xi^{}  f_\xi^{}} }}{\displaystyle  p^{\;\!h_\xi^{}}}+
\varphi_{_{\scriptstyle h_\xi^{}  f_\xi^{}}}\biggr)\biggr)\in 
\text{\rm M}_{_{\Omega_{{}_{\tiny\;\! 0}}}}
\hfill
$
\\[2.25ex]
if and only if
\\[1.75ex]
\mbox{}\hfill                           
$
\displaystyle
{\rm div}\;\!{\frak d}(t,x)=
{}-\sum\limits_{\xi=1}^{\varepsilon}\sum\limits_{f_{\xi}^{}=1}^{\delta_{\xi}^{}}
\biggl(\gamma_{_{\scriptstyle h_\xi^{}  f_\xi^{}}}\Bigl(
N_{h_\xi^{}  f_\xi^{}}^{} (t,x)+{\sf D}\;\!\varphi_{_{\scriptstyle h_\xi^{}  f_\xi^{}}}(t)\Bigr)\biggr)
- \gamma\;\!M(t,x)
$
\ for all} $(t,x)\in \Xi^{\;\!\prime}.
\hfill
$
\\[2ex]
\indent
{\sl Proof.}\;\!
\vspace{0.25ex}
Using Theorem 8.1, from Property 5.8 under the condition 
$L={}-\,{\rm div}\;\!{\frak d},$ we obtain the statement of this property. $\k$ 
\vspace{1.25ex}

{\bf Property 9.13.} 
{\it
We claim that}
\\[1.75ex]
\mbox{}\hfill                         
$
\displaystyle
u+i\;\!v\in \text{\rm MH}_{_{\Xi^{\;\!\prime}}}
\iff
u\in \text{\rm MA}_{_{\Xi^{\;\!\prime}}}
\ \&\ 
v\in \text{\rm MA}_{_{\Xi^{\;\!\prime}}}.
\hfill
$
\\[2ex]
\indent
{\sl Proof.}\;\!
\vspace{0.25ex}
Taking into account Theorem 8.1, from Property 6.11 under   
$U={}-\,{\rm div}\;\!{\frak d},$ we  get the statement of Property 9.13. $\k$
\vspace{1.25ex}

{\bf Property 9.14.}
{\it
Suppose functions
\vspace{0.75ex}
$u,v\in \text{\rm P}_{_{\!\Xi^{\;\!\prime}}}$ are relatively prime, 
the set $\Omega_{_0}\subset \Xi^{\;\!\prime}$ such that
$u(t,x)\ne 0$ for all $(t,x)\in \Omega_{_0}$ and
$u(t,x)=0$ for all $(t,x)\in {\sf C}_{_{\Xi^{\;\!\prime}}}\Omega_{_0}.$
Then}
\\[2.5ex]
\mbox{}\hfill
$
u+i\;\!v\in \text{\rm MH}_{_{\Xi^{\;\!\prime}}}
\iff
(u^2+v^2, {}-2\;\!{\rm div}\;\!{\frak d})\in \text{\rm A}_{_{\Xi^{\;\!\prime}}}
\ \&\ \,
\dfrac{v}{u}\in \text{\rm I}_{_{\Omega_{{}_{\tiny\;\! 0}}}}.
\hfill
$
\\[2.25ex]
\indent
{\sl Proof.}\;\!
\vspace{0.25ex}
Taking into account Theorem 8.1, from Property 6.12 under   
$U={}-\,{\rm div}\;\!{\frak d},$ we  obtain the statement of Property 9.14. $\k$
\vspace{1.25ex}

{\bf Property 9.15.} 
{\it
Let $\eta\in\C\backslash\{0\}.$ Then we have}
\\[2ex]
\mbox{}\hfill
$
w\in \text{\rm MH}_{_{\Xi^{\;\!\prime}}}
\iff 
\eta\;\! w\in \text{\rm MH}_{_{\Xi^{\;\!\prime}}}.
\hfill
$
\\[1.75ex]
\indent
{\sl Proof.}\;\!
\vspace{0.25ex}
Using Theorem 8.1, from Property 6.1 under the condition 
$W={}-\,{\rm div}\;\!{\frak d}$ it follows that 
the statement of Property 9.15 is true. $\k$

\newpage

{\bf Property 9.16.} 
{\it
Let $k\in\N.$ Then we claim that
\\[1.75ex]
\mbox{}\hfill
$
w^k\in \text{\rm MH}_{_{\Xi^{\;\!\prime}}}
\iff 
\Bigl(w, {}-\dfrac{1}{k}\, \;\!{\rm div}\;\!{\frak d}\Bigr)\in \text{\rm H}_{_{\Xi^{\;\!\prime}}}
\hfill
$
\\[1ex]
and}
\\[1ex]
\mbox{}\hfill
$
w\in \text{\rm MH}_{_{\Xi^{\;\!\prime}}}
\iff 
(w^k,{}- k\,{\rm div}\;\!{\frak d})\in \text{\rm H}_{_{\Xi^{\;\!\prime}}}.
\hfill
$
\\[1.75ex]
\indent
{\sl Proof.}\;\! 
\vspace{1.5ex}
Using Theorem 8.1, from Property 6.3, we get the statement of Property 9.16. $\k$

{\bf Property 9.17.} 
{\it
Let $\eta_j^{}\in\R\backslash\{0\},\ j=1,\ldots,m.$ Then}
\\[1.5ex]
\mbox{}\hfill
$
\displaystyle
w_j^{}\in \text{\rm MH}_{_{\Xi^{\;\!\prime}}}, \ j=1,\ldots, m,
\ \ \Longrightarrow\ \ 
\sum\limits_{j=1}^{m} \eta_j^{}\;\!w_j^{}
\in \text{\rm MH}_{_{\Xi^{\;\!\prime}}}.
\hfill
$
\\[1.5ex]
\indent
{\sl Proof.}
\vspace{0.25ex}
Taking into account Theorem 8.1, from Property 6.4 under
$W={}-\,{\rm div}\;\!{\frak d},$ we obtain the statement of Property 9.17. $\k$ 
\vspace{1.5ex}

{\bf Property 9.18.} 
{\it
Let 
$k_j^{}\in\N, \ w_j^{}\in \text{\rm Z}_{_{\Xi^{\;\!\prime}}},\ j=1,\ldots,m.$ Then
\\[1.5ex]
\mbox{}\hfill                         
$
\displaystyle
\prod\limits_{j=1}^{m} w_j^{k_j^{}} 
\in \text{\rm MH}_{_{\Xi^{\;\!\prime}}}
\iff
\bigl( w_j^{}\;\!,\;\! W_j^{}\bigl)\;\! \in \text{\rm H}_{_{\Xi^{\;\!\prime}}},
\quad
j=1,\ldots,m.
\hfill
$
\\[1.5ex]
Moreover, the cofactors $W_1^{},\ldots,W_m^{}$ such that
\\[1.5ex]
\mbox{}\hfill
$
\displaystyle
\sum\limits_{j=1}^{m}\! k_j^{}\;\!{\rm Re}\;\!W_j^{}(t,x)=
-\;\!{\rm div}\;\!{\frak d}(t,x)
$
for all $(t,x)\!\in\! \Xi^{\;\!\prime},
\quad
\sum\limits_{j=1}^{m}\! k_j^{}\;\!{\rm Im}\;\!W_j^{}(t,x)\!=\!0
$
for all} $(t,x)\!\in\! \Xi^{\;\!\prime}.
\hfill
$
\\[1.5ex]
\indent
{\sl Proof.}
\vspace{0.25ex}
Taking into account Theorem 8.1, from Property 6.6 under
$W={}-\,{\rm div}\;\!{\frak d},$ we get the statement of Property 9.18. $\k$ 
\vspace{1.25ex}

{\bf Property 9.19.} 
{\it
If 
$p_\tau^{}\in \text{\rm P}_{_{\!\Xi^{\;\!\prime}}},\, l_{\tau}^{}\!\in\N,\, \tau\!=\!1,\ldots,s,\
w_j^{}\in \text{\rm Z}_{_{\Xi^{\;\!\prime}}},\, k_j^{}\!\in\N, \, j\!=\!1,\ldots,m,$ then  
\\[1.5ex]
\mbox{}\hfill                         
$
\displaystyle
\prod\limits_{\tau=1}^{s}\! p_\tau^{\,l_\tau^{}}
\prod\limits_{j=1}^{m}\! w_j^{k_j^{}} 
\in\! \text{\rm MH}_{_{\Xi^{\;\!\prime}}}
\!\iff\!
\bigl(p_\tau^{}\;\!, M_\tau^{}\bigl)\;\! \in\! \text{\rm A}_{_{\Xi^{\;\!\prime}}},\,
\tau\!=\!1,\ldots,s,
\ \&\ 
\bigl(w_j^{}\;\!, W_j^{}\bigl)\;\!\in\! \text{\rm H}_{_{\Xi^{\;\!\prime}}},\,
j\!=\!1,\ldots,m.
\hfill
$
\\[1.5ex]
Moreover, the cofactors $M_1^{},\ldots,M_s^{},\, W_1^{},\ldots,W_m^{}$ such that
\\[1.5ex]
\mbox{}\hfill
$
\displaystyle
\sum\limits_{\tau=1}^{s} l_\tau^{}\;\!M_\tau^{}(t,x)+
\sum\limits_{j=1}^{m} k_j^{}\;\!{\rm Re}\;\!W_j^{}(t,x)=
{}-\;\! {\rm div}\;\!{\frak d}(t,x),
\ \
\sum\limits_{j=1}^{m} k_j^{}\;\!{\rm Im}\;\!W_j^{}(t,x)=0
$
for all} $(t,x)\!\in\! \Xi^{\;\!\prime}.
\hfill
$
\\[1.5ex]
\indent
{\sl Proof.}
\vspace{0.25ex}
Taking into account Theorem 8.1, from Property 6.7 under
$W={}-\,{\rm div}\;\!{\frak d},$ we obtain the statement of Property 9.19. $\k$ 
\vspace{1ex}

{\bf Property 9.20.} 
\vspace{0.15ex}
{\it
Suppose 
$p\in \text{\rm P}_{_{\!\Xi^{\;\!\prime}}},$ 
functions $u,v\in\! \text{\rm P}_{_{\!\Xi^{\;\!\prime}}}$ 
are relatively prime, $w=u+i\;\!v.$ Then we have  
\\[0.75ex]
\mbox{}\hfill                         
$
\displaystyle
p\;\!w\in \text{\rm MH}_{_{\Xi^{\;\!\prime}}}
\iff
(p, M) \in \text{\rm A}_{_{\Xi^{\;\!\prime}}}
\ \&\ 
(w, W) \in \text{\rm H}_{_{\Xi^{\;\!\prime}}}.
\hfill
$
\\[2ex]
Moreover, the cofactors $M$ and $W$ such that
\\[1.75ex]
\mbox{}\hfill
$
\displaystyle
M(t,x)+{\rm Re}\;\!W(t,x)={}-\,{\rm div}\;\!{\frak d}(t,x)
$
for all $(t,x)\in \Xi^{\;\!\prime},
\quad
{\rm Im}\;\!W(t,x)=0
$
for all} $(t,x)\in \Xi^{\;\!\prime}.
\hfill
$
\\[1.75ex]
\indent
{\sl Proof.}\;\!
\vspace{1.25ex}
Using Theorem 8.1, from Property 6.8, we get the statement of Property 9.20. $\k$

{\bf Property 9.21.} 
{\it
We claim that}
\\[1.5ex]
\mbox{}\hfill                         
$
\displaystyle
u+i\;\!v\in \text{\rm MH}_{_{\Xi^{\;\!\prime}}}
\iff
u-i\;\!v\in \text{\rm MH}_{_{\Xi^{\;\!\prime}}}.
\hfill
$
\\[1.5ex]
\indent
{\sl Proof.}
\vspace{0.25ex}
Taking into account Theorem 8.1, from Property 6.9 under the conditions
$V=0,$ $U={}-\,{\rm div}\;\!{\frak d},$ 
we obtain the statement of Property 9.21. $\k$

\newpage

{\bf Property 9.22}
\vspace{0.35ex}
(existence criterion of multiple complex-valued polynomial last multiplier).
{\it
Suppose $h\in\N,$ functions
\vspace{1ex}
$w, z\in \text{\rm Z}_{_{\Xi^{\;\!\prime}}}\!$ are relatively prime, 
the set 
$\Omega_{_0}\subset\Xi^{\;\!\prime}$ such that  
$|w(t,x)|\ne 0$ for all $(t,x)\in\Omega_{_0}$ and
$|w(t,x)|= 0$ for all $(t,x)\in {\sf C}_{_{\Xi^{\;\!\prime}}}\Omega_{_0}.$
Then}
\\[2.5ex]
\mbox{}\hfill                           
$
\displaystyle
\bigl(w, (h,z,Q)\bigr)\in \text{\rm MG}_{_{\Xi^{\;\!\prime}}}
\iff
w\in \text{\rm MH}_{_{\Xi^{\;\!\prime}}}
\ \& \
\hfill                           
$
\\[2.5ex]
\mbox{}\hfill                           
$
\& \
\Bigl(\exp\;\!{\rm Re}\;\!\dfrac{z}{w^{\;\!h}}\,, {\rm Re}\;\!Q\Bigr)\in 
\text{\rm E}_{_{\Omega_{{}_{\tiny\;\! 0}}}}
\ \,\& \ \,
\Bigl(\exp\;\!{\rm Im}\;\!\dfrac{z}{w^{\;\!h}}\,, {\rm Im}\;\!Q\Bigr)\in 
\text{\rm E}_{_{\Omega_{{}_{\tiny\;\! 0}}}}.
\hfill
$
\\[2.5ex]
\indent
{\sl Proof.}\;\!
\vspace{0.25ex}
Using Theorem 8.1, from Theorem 7.2 under the condition
$M={}-\,{\rm div}\;\!{\frak d}$
it follows that 
the statement of Property 9.22 is true. $\k$ 
\vspace{1.5ex}

{\bf Property 9.23}
\vspace{0.35ex}
(existence criterion of multiple complex-valued polynomial last multiplier).
{\it
Suppose $h\in\N,$ functions
\vspace{0.75ex}
$u, v\in \text{\rm P}_{_{\!\Xi^{\;\!\prime}}}$ are relatively prime, 
a function $z\in \text{\rm Z}_{_{\Xi^{\;\!\prime}}}$ 
is relatively prime to the function $u+i\;\!v,$  the set
\vspace{1ex}
$\Omega_{_0}\subset\Xi^{\;\!\prime}$ such that   
$u(t,x)\ne 0$ for all $(t,x)\in\Omega_{_0}$ and
$u(t,x)= 0$ for all $(t,x)\in {\sf C}_{_{\Xi^{\;\!\prime}}}\Omega_{_0}.$
Then}
\\[2.25ex]
\mbox{}\hfill                           
$
\displaystyle
\bigl(u+i\;\!v, (h,z,Q)\bigr)\in \text{\rm MG}_{_{\Xi^{\;\!\prime}}}
\iff
\Bigl(\bigl(u^2+v^2, {}-2\,{\rm div}\;\!{\frak d}\bigr), 
\bigl(h,\;\!{\rm Re}\;\!\bigl(z(u-i\;\!v)^h\bigr),\;\! {\rm Re}\;\!Q\bigr)\Bigr)\in 
\text{\rm B}_{_{\Xi^{\;\!\prime}}}
\ \& \
\hfill                           
$
\\[2.75ex]
\mbox{}\hfill                           
$
\& \
\Bigl(\bigl(u^2+v^2, {}- 2\,{\rm div}\;\!{\frak d}\bigr), 
\bigl(h,\;\!{\rm Im}\;\!\bigl(z(u-i\;\!v)^h\bigr),\;\! {\rm Im}\;\!Q\bigr)\Bigr)\in 
\text{\rm B}_{_{\Xi^{\;\!\prime}}}
\ \,\& \ \,
\dfrac{v}{u}\in 
\text{\rm I}_{_{\Omega_{{}_{\tiny\;\! 0}}}}.
\hfill
$
\\[2.5ex]
\indent
{\sl Proof.}\;\!
\vspace{0.25ex}
Using Theorem 8.1, from Theorem 7.3 under the condition
$M={}-\,{\rm div}\;\!{\frak d}$
it follows that 
the statement of Property 9.23 is true. $\k$ 
\vspace{1ex}

By Definition 7.2, we introduce the notion of multiplicity 
for multiple complex-valued polynomial  last multiplier.
\\[4.25ex]
\centerline{
{\bf  10. Exponential last multipliers
}
}
\\[1.75ex]
\indent
{\bf Theorem\! 10.1}\!
\vspace{0.35ex}
(existence criterion of exponential last multiplier).\!
{\it
A function $\exp\omega\!\in\! \text{\rm ME}_{_{\Omega}}\!$ 
if and only if the following identity holds
\\[1.5ex]
\mbox{}\hfill               
$
{\frak d}\;\!\omega(t,x)={}-\,{\rm div}\;\!{\frak d}(t,x)
$
\ for all} 
$(t,x)\in \Omega.
\hfill
$
\\[1.5ex]
\indent
{\sl Proof}.
\vspace{0.25ex}
Taking into account Theorem 8.1, from Theorem 3.1 under 
$M={}-\,{\rm div}\;\!{\frak d},$ 
we get the statement of Theorem 10.1 is true. $\k$ 
\vspace{1ex}

{\bf Property 10.1.}
{\it
If
\\[1ex]
\mbox{}\hfill               
$
{\rm div}\;\!{\frak d}(t,x)=\varphi(t)
$
\ for all $(t,x)\in \Xi,
\hfill
$
\\[1ex]
then the function
\\[1ex]
\mbox{}\hfill               
$
\displaystyle
\mu\colon t\to\ 
\exp\biggl({}-\int\limits_{t_0^{}}^{t}\varphi(\tau)\;\!d\tau\biggr)
$
\ for all $t\in T,
\hfill
$
\\[1.5ex]
where $t_{_0}^{}$ is an arbitrary fixed point from the domain $T,$
\vspace{0.35ex}
is a conditional  last multiplier on the domain $\Xi$ of system {\rm (0.1).}
}
\vspace{0.5ex}

{\sl Proof.}
The derivative by virtue of system (0.1)
\\[1.5ex]
\mbox{}\hfill               
$
\displaystyle
{\frak d}\;\!\biggl({}-\int\limits_{t_0^{}}^{t}\varphi(\tau)\;\!d\tau\biggr)=
{\sf D}\;\!\biggl({}-\int\limits_{t_0^{}}^{t}\varphi(\tau)\;\!d\tau\biggr)=
{}-\varphi(t)=
{}-\;\!{\rm div}\;\!{\frak d}(t,x)
$
\ for all 
$(t,x)\in \Xi.
\hfill
$
\\[1.5ex]
\indent
By Theorem 10.1, $\mu\in  \text{\rm MF}_{_{\!\Xi}}.\ \k$
\vspace{1ex}

In particular, from Property 10.1, we get
\vspace{0.75ex}

{\bf Property 10.2.}
{\it
If
\\[1ex]
\mbox{}\hfill               
$
{\rm div}\;\!{\frak d}(t,x)=\lambda
$
\ for all 
$(t,x)\in \Xi,
\quad
\lambda\in\R,
\hfill
$
\\[1.5ex]
then the function $e^{{}-\lambda\;\!t}\in  \text{\rm MF}_{_{\!\Xi}}.$
}
\vspace{1.25ex}

{\bf  Property 10.3.}
{\it
Let $\gamma_j^{}\in\R\backslash\{0\},\ c_j^{}\in\R, \ j=1,\ldots,m.$
Then an exponential function
\\[1.5ex]
\mbox{}\hfill
$
\displaystyle
\exp\sum\limits_{j=1}^{m}
\gamma_j^{}\;\!(\omega_j^{}+c_j^{})\in \text{\rm ME}_{_{\Omega}}
\hfill
$
\\[0.75ex]
if and only if
\\[0.75ex]
\mbox{}\hfill
$
\displaystyle
\sum\limits_{j=1}^{m}\gamma_j^{}\;\!{\frak d}\;\!\omega_j^{}(t,x)=
{}-\;\!{\rm div}\;\!{\frak d}(t,x)
$
\ for all} $(t,x)\in \Omega.
\hfill
$
\\[1.5ex]
\indent
{\sl Proof}. 
Taking into account Theorem 8.1, from Theorem 3.3 it follows that 
the statement of Property 10.3 is true. $\k$
\vspace{1ex}

{\bf  Property 10.4.}
{\it
Let $\varphi\in C^1T^{\;\!\prime}.$ Then we claim that 
\\[1.5ex]
\mbox{}\hfill
$
\mu\in \text{\rm M}_{_{\Omega}}
\iff
\bigl(\mu\exp\varphi, {\sf D} \varphi-\;\!{\rm div}\;\!{\frak d}\bigr)\in \text{\rm J}_{_{\Omega}}
\hfill
$
\\[0.5ex]
and}
\\[0.5ex]
\mbox{}\hfill
$
{\rm g}\exp\varphi\in \text{\rm M}_{_{\Omega}}
\iff
\bigl({\rm g}, {}-{\sf D} \varphi-\;\!{\rm div}\;\!{\frak d}\bigr)\in \text{\rm J}_{_{\Omega}}.
\hfill
$
\\[1.5ex]
\indent
{\sl Proof}. 
Taking into account Theorem 8.1, from Property 3.1 it follows that 
the statement of Property 10.4 is true. $\k$
\vspace{1ex}

{\bf  Property 10.5.}
{\it
Let $\varphi\in C^1T^{\;\!\prime}.$ Then we have 
\\[1.5ex]
\mbox{}\hfill
$
\exp\omega\in \text{\rm ME}_{_{\Omega}}
\iff
\bigl(\exp(\omega+\varphi),\, {\sf D} \varphi-\;\!{\rm div}\;\!{\frak d}\bigr)\in \text{\rm E}_{_{\Omega}}
\hfill
$
\\[0.5ex]
and}
\\[0.5ex]
\mbox{}\hfill
$
\exp(\omega+\varphi)\in \text{\rm ME}_{_{\Omega}}
\iff
\bigl(\exp\omega,\, {}-{\sf D} \varphi-\;\!{\rm div}\;\!{\frak d}\bigr)\in \text{\rm E}_{_{\Omega}}.
\hfill
$
\\[1.5ex]
\indent
{\sl Proof}. 
Taking into account Theorem 8.1, from Property 3.2 it follows that 
the statement of Property 10.5 is true. $\k$
\vspace{1ex}

{\bf  Property 10.6.}
{\it
If $c\in\R,$ then} 
\\[1.5ex]
\mbox{}\hfill
$
\displaystyle
\exp(\omega+c)\in \text{\rm ME}_{_{\Omega}}
\iff
\exp\omega\in \text{\rm ME}_{_{\Omega}}.
\hfill
$
\\[1.5ex]
\indent
{\sl Proof}. 
\vspace{0.25ex}
Using Theorem 8.1, from Property 3.3 under the condition
$M={}-\;\!{\rm div}\;\!{\frak d},$
we see that the statement of Property 10.6 is true. $\k$
\vspace{1ex}

{\bf  Property 10.7.}
{\it
Suppose  $\lambda_j^{}\in\R\backslash\{0\}, \ j=1,\ldots,m.$ Then}
\\[1.5ex]
\mbox{}\hfill
$
\displaystyle
\exp\omega_j^{}\in \text{\rm ME}_{_{\Omega}},\  j=1,\ldots,m, \ \ 
\Longrightarrow\ \ 
\sum\limits_{j=1}^{m}\lambda_j^{}\exp\omega_j^{}\in \text{\rm M}_{_{\Omega}}.
\hfill
$
\\[1.5ex]
\indent
{\sl Proof}. 
\vspace{0.25ex}
Using Theorem 8.1, from Property 3.4 under the condition
$M={}-\;\!{\rm div}\;\!{\frak d},$
we get that the statement of Property 10.7 is true. $\k$
\vspace{1ex}

{\bf  Property 10.8.}
{\it
Suppose
$\rho_j^{},\lambda_j^{}\in\R\backslash\{0\}, \ j=1,\ldots,m.$ Then} 
\\[1.25ex]
\mbox{}\hfill
$
\displaystyle
\bigl(\exp\omega_j^{}\;\!, {}-\rho_j^{}\;\!{\rm div}\;\!{\frak d}\bigr)\in \text{\rm E}_{_{\Omega}},
\ j=1,\ldots,m,
\ \ \Longrightarrow\ \ 
\sum\limits_{j=1}^{m}\lambda_j^{}
\exp\dfrac{\omega_j^{}}{\rho_j^{}}\in \text{\rm M}_{_{\Omega}}.
\hfill
$
\\[1.25ex]
\indent
{\sl Proof}. 
\vspace{0.25ex}
Using Theorem 8.1, from Property 3.5 under the condition
$M={}-\;\!{\rm div}\;\!{\frak d},$
we have that the statement of Property 10.8 is true. $\k$
\vspace{1ex}

{\bf  Property 10.9.}
\vspace{0.5ex}
{\it
Suppose  
$\bigl(\exp\omega_j^{}, M_j^{}\bigr)\in \text{\rm E}_{_{\Omega}},\ 
\gamma_j^{}\in\R\backslash\{0\}, \ j=1,\ldots,m.$ Then the fun\-c\-tion
$
\exp\sum\limits_{j=1}^{m}\gamma_j^{}\;\!\omega_j^{}\in \text{\rm ME}_{_{\Omega}}
$
if and only if the identity {\rm(8.1)} is true}.
\vspace{0.75ex}

{\sl Proof}. 
\vspace{0.25ex}
Using Theorem 8.1, from Property 3.6 under the condition
$M={}-\;\!{\rm div}\;\!{\frak d},$
we obtain that the statement of Property 10.9 is true. $\k$
\vspace{1.25ex}

{\bf  Property 10.10.}
\vspace{0.5ex}
{\it
Suppose  
$\bigl(\exp\omega_j^{}, \rho_j^{}M_{_0}^{}\bigr)\in \text{\rm E}_{_{\Omega}},\ 
\rho_j^{},\gamma_j^{}\in\R\backslash\{0\}, \ j=1,\ldots,m.$ Then 
\\[1.5ex]
\mbox{}\hfill
$
\displaystyle
\exp\sum\limits_{j=1}^{m}\gamma_j^{}\;\!\omega_j^{}\in \text{\rm ME}_{_{\Omega}}
\hfill
$
\\[1.5ex]
if and only if the following identity holds
\\[1.5ex]
\mbox{}\hfill
$
\displaystyle
\sum\limits_{j=1}^{m}
\rho_j^{}\;\!\gamma_j^{}\;\!M_{_0}^{}(t,x)=
{}-\;\!{\rm div}\;\!{\frak d}(t,x)
$
\ for all} 
$(t,x)\in \Xi^{\;\!\prime}.
\hfill
$
\\[1.5ex]
\indent
{\sl Indeed}, it follows from Property 10.9 under the conditions 
$M_j^{}=\rho_j^{}\;\!M_{_0}^{},\ j=1,\ldots, m.\ \k$
\vspace{1.5ex}

{\bf  Property 10.11.}
\vspace{1ex}
{\it
Suppose  
$\bigl(\exp\omega_j^{}, {}-\rho_j^{}\;\!{\rm div}\;\!{\frak d}\bigr)\in \text{\rm E}_{_{\Omega}},\ 
\rho_j^{},\gamma_j^{}\in\R\backslash\{0\}, \ j=1,\ldots,m.$ Then 
$
\exp\sum\limits_{j=1}^{m}\gamma_j^{}\;\!\omega_j^{}\in \text{\rm ME}_{_{\Omega}},
$
if and only if} 
$\sum\limits_{j=1}^{m}\rho_j^{}\;\!\gamma_j^{}=1.$
\vspace{1ex}

{\sl Indeed}, it follows from Property 10.10 under the condition
$M_{_0}^{}={}-\;\!{\rm div}\;\!{\frak d}.\ \k$
\vspace{1.25ex}

{\bf  Property 10.12.}
\vspace{1ex}
{\it
Suppose  
$\exp\omega_j^{}\in \text{\rm ME}_{_{\Omega}},\ 
\gamma_j^{}\in\R\backslash\{0\}, \ j=1,\ldots,m.$ Then we claim that
$
\Bigl(\exp\sum\limits_{j=1}^{m}\gamma_j^{}\;\!\omega_j^{}, M\Bigr)\in \text{\rm E}_{_{\Omega}}
$
if and only if the identity {\rm (8.2)} is true.} 
\vspace{1ex}

{\sl Proof}. 
\vspace{0.35ex}
Taking into account Theorem 8.1, from Property 3.6 under the conditions
\linebreak
$M_{j}^{}={}-\;\!{\rm div}\;\!{\frak d}, \ j=1,\ldots,m,$
\vspace{1.5ex}
we see that the statement of Property 10.12 is true. $\k$

{\bf  Property 10.13.}
\vspace{1ex}
{\it
Let  
$\exp\omega_j^{}\in \text{\rm ME}_{_{\Omega}},\ 
\gamma_j^{}\in\R\backslash\{0\}, \ j=1,\ldots,m.$ Then the exponential function
$
\exp\sum\limits_{j=1}^{m}\gamma_j^{}\;\!\omega_j^{}\in \text{\rm ME}_{_{\Omega}}
$
if and only if}
$\sum\limits_{j=1}^{m}\gamma_j^{}=1.$
\vspace{1ex}

{\sl Indeed}, it follows from Property 10.11 under the conditions
$\rho_{j}^{}=1,\ j=1,\ldots,m.\ \k$
\vspace{1.25ex}

{\bf  Property 10.14.}
{\it
If $\gamma\in\R\backslash\{0\},$ then 
\\[2ex]
\mbox{}\hfill
$
\displaystyle
\exp\omega\in \text{\rm ME}_{_{\Omega}}
\iff
\bigl(\exp(\gamma\;\!\omega), {}-\gamma\;\!{\rm div}\;\!{\frak d}\bigr)\in \text{\rm E}_{_{\Omega}}
\hfill
$
\\[0.75ex]
and}
\\[0.75ex]
\mbox{}\hfill
$
\displaystyle
\exp(\gamma\;\!\omega)\in \text{\rm ME}_{_{\Omega}}
\iff
\Bigl(\exp\omega, {}-\dfrac{1}{\gamma}\ {\rm div}\;\!{\frak d}\bigr)\in \text{\rm E}_{_{\Omega}}.
\hfill
$
\\[1.5ex]
\indent
{\sl Proof}. 
\vspace{0.15ex}
Taking into account Theorem 8.1, from Corollary 3.2 
it follows that the statement of Property 10.14 is true. $\k$ 
\vspace{1ex}

{\bf  Property 10.15.}
{\it
If 
$
\bigl(\exp\omega_j^{}, M_j^{}\bigr)\in \text{\rm E}_{_{\Omega}},\ 
\gamma,\gamma_j^{}\in\R\backslash\{0\}, \ j=1,\ldots,m, \
{\rm g}^{\gamma}\in C^1\Omega,$ then}
\\[1.75ex]
\mbox{}\hfill
$
\displaystyle
{\rm g}^{\gamma}\exp\sum\limits_{j=1}^{m}\gamma_j^{}\;\!\omega_j^{}\in 
\text{\rm M}_{_{\Omega}}
\iff
\biggl(\,
{\rm g}, {}-\dfrac{1}{\gamma}\,\biggl({\rm div}\;\!{\frak d}+
\sum\limits_{j=1}^{m}\gamma_j^{}\;\!M_j^{}\biggr)\biggr)\in 
\text{\rm J}_{_{\Omega}}.
\hfill
$
\\[1.5ex]
\indent
{\sl Proof.} 
\vspace{0.25ex}
Using Theorem 8.1, from Property 3.7 under the condition 
$M={}-\;\!{\rm div}\;\!{\frak d},$ we obtain the statement of Property 10.15. $\k$
\vspace{1ex}

{\bf  Property 10.16.}
{\it
Let  
$
\exp\omega_j^{}\in \text{\rm ME}_{_{\Omega}},\ 
\gamma,\gamma_j^{}\in\R\backslash\{0\}, \ j=1,\ldots,m, \
{\rm g}^{\gamma}\in C^1\Omega.$ Then}
\\[1.5ex]
\mbox{}\hfill
$
\displaystyle
\biggl(\,
{\rm g}^{\gamma}\exp\sum\limits_{j=1}^{m}\gamma_j^{}\;\!\omega_j^{}\;\!,\;\! M\biggr)\in 
\text{\rm J}_{_{\Omega}}
\iff
\biggl(\,
{\rm g},\, \dfrac{1}{\gamma}\,\biggl(M+
\sum\limits_{j=1}^{m}\gamma_j^{}\;\!{\rm div}\;\!{\frak d}\biggr)\biggr)\in 
\text{\rm J}_{_{\Omega}}.
\hfill
$
\\[1.5ex]
\indent
{\sl Proof.} 
\vspace{0.25ex}
Taking into account Theorem 8.1, from Property 3.7 under the conditions 
\linebreak
$M_j^{}={}-\;\!{\rm div}\;\!{\frak d},\ j=1,\ldots, m,$ 
we get the statement of Property 10.16. $\k$
\vspace{1ex}

{\bf  Property 10.17.}
{\it
Let  
$
\exp\omega_j^{}\in \text{\rm ME}_{_{\Omega}},\ 
\gamma,\gamma_j^{}\in\R\backslash\{0\}, \ j=1,\ldots,m, \
{\rm g}^{\gamma}\in C^1\Omega.$ Then}
\\[1.5ex]
\mbox{}\hfill
$
\displaystyle
{\rm g}^{\gamma}\exp\sum\limits_{j=1}^{m}\gamma_j^{}\;\!\omega_j^{}\in 
\text{\rm M}_{_{\Omega}}
\iff
\biggl(\,
{\rm g},\, \dfrac{1}{\gamma}\,\biggl(\
\sum\limits_{j=1}^{m}\gamma_j^{}-1\biggr)\;\!{\rm div}\;\!{\frak d}\biggr)\in 
\text{\rm J}_{_{\Omega}}.
\hfill
$
\\[1.5ex]
\indent
{\sl Proof.} 
\vspace{0.25ex}
Using Theorem 8.1, from Property 10.16 under the condition 
$M={}-\;\!{\rm div}\;\!{\frak d},$ we obtain the statement of Property 10.17. $\k$
\vspace{1ex}

{\bf  Property 10.18.}
{\it
Suppose
$
\bigl(\exp\omega_\tau^{}, M_\tau^{}\bigr)\!\in \text{\rm E}_{_{\Omega}},\
\tau\!=\!1,\ldots, m\!-\!1,\ 
\gamma_j^{}\in\R\backslash\{0\}, \ j\!=\!1,\ldots,m.$ Then we claim that}
\\[1.5ex]
\mbox{}\hfill
$
\displaystyle
\exp\sum\limits_{j=1}^{m}\gamma_j^{}\;\!\omega_j^{}\in 
\text{\rm ME}_{_{\Omega}}
\iff
\biggl(\,
\exp\omega_{m}^{}\;\!,{}-
 \dfrac{1}{\gamma_{m}^{}}\,\biggl({\rm div}\;\!{\frak d}+
\sum\limits_{\tau=1}^{m-1}\gamma_\tau^{}\;\!M_\tau^{}\biggr)\biggr)\in 
\text{\rm E}_{_{\Omega}}.
\hfill
$
\\[1.5ex]
\indent
{\sl Proof.} 
\vspace{0.25ex}
Using Theorem 8.1, from Corollary 3.3 under the condition 
$M={}-\;\!{\rm div}\;\!{\frak d}$ 
it follows that the statement of Property 10.17 is true. $\k$
\vspace{1ex}

{\bf  Property 10.19.}
{\it
If
$
\exp\omega_\tau^{}\!\in \text{\rm ME}_{_{\Omega}},\
\tau\!=\!1,\ldots, m\!-\!1,\ 
\gamma_j^{}\in\R\backslash\{0\}, \ j\!=\!1,\ldots,m,$ then}
\\[1.5ex]
\mbox{}\hfill
$
\displaystyle
\biggl(\,
\exp\sum\limits_{j=1}^{m}\gamma_j^{}\;\!\omega_j^{}\;\!,\;\! M\biggr)\in 
\text{\rm E}_{_{\Omega}}
\iff
\biggl(\,
\exp\omega_{m}^{}\;\!,\,
 \dfrac{1}{\gamma_{m}^{}}\,\biggl(M+
\sum\limits_{\tau=1}^{m-1}\gamma_\tau^{}\;\!{\rm div}\;\!{\frak d}\biggr)\biggr)\in 
\text{\rm E}_{_{\Omega}}.
\hfill
$
\\[1.5ex]
\indent
{\sl Proof.} 
\vspace{0.25ex}
Taking into account Theorem 8.1, from Corollary 3.3 under the conditions 
\linebreak
$M_{\tau}^{}={}-\;\!{\rm div}\;\!{\frak d},\ \tau=1,\ldots,m-1,$
it follows that the statement of Property 10.19 is true. $\k$
\vspace{1ex}

{\bf  Property 10.20.}
{\it
If
$
\exp\omega_\tau^{}\!\in \text{\rm ME}_{_{\Omega}},\
\tau\!=\!1,\ldots, m\!-\!1,\ 
\gamma_j^{}\in\R\backslash\{0\}, \ j\!=\!1,\ldots,m,$ then}
\\[1.5ex]
\mbox{}\hfill
$
\displaystyle
\exp\sum\limits_{j=1}^{m}\gamma_j^{}\;\!\omega_j^{}\in 
\text{\rm ME}_{_{\Omega}}
\iff
\biggl(\,
\exp\omega_{m}^{}\;\!,\,
 \dfrac{1}{\gamma_{m}^{}}\,\biggl(\,
\sum\limits_{\tau=1}^{m-1}\gamma_\tau^{}-1\biggr)\;\!{\rm div}\;\!{\frak d}
\biggr)\in 
\text{\rm E}_{_{\Omega}}.
\hfill
$
\\[1.5ex]
\indent
{\sl Proof.} 
\vspace{0.25ex}
Using Theorem 8.1, from Property 10.19 under the condition 
$M={}-\;\!{\rm div}\;\!{\frak d},$ we get the statement of Property 10.20. $\k$
\vspace{1ex}

{\bf  Property 10.21.}
\vspace{0.5ex}
{\it
Suppose 
$
\bigl(\exp\omega_\nu^{}, M_\nu^{}\bigr)\in \text{\rm E}_{_{\Omega}},\ 
\nu=1,\ldots, s,\ s\leqslant m-2,\  
\gamma_j^{}\in\R\backslash\{0\}, 
\linebreak 
j=1,\ldots,m,\ 
{\rm g}_k^{{}^{\scriptsize \gamma_k^{}}}\!\in C^1\Omega,\ 
k=s+1,\ldots,m.$ Then we have
\\[1.75ex]
\mbox{}\hfill
$
\displaystyle
\prod\limits_{k=s+1}^{m}
{\rm g}_k^{{}^{\scriptsize \gamma_k^{}}}
\exp\sum\limits_{\nu=1}^{s}
\gamma_\nu^{}\;\!\omega_\nu^{}\in 
\text{\rm M}_{_{\Omega}}
\iff
\biggl(\,
\prod\limits_{k=s+1}^{m}
{\rm g}_k^{{}^{\scriptsize \gamma_k^{}}},\,-\;\!
{\rm div}\;\!{\frak d}-\sum\limits_{\nu=1}^{s}\gamma_\nu^{}\;\!M_\nu^{}\biggr)\in 
\text{\rm J}_{_{\Omega}}.
\hfill
$
\\[1.5ex]
Moreover, there exist functions $M_k^{}\in C^1\Omega,\ k=s+1,\ldots, m,$ 
such that the identities {\rm (1.14)} are true and
\\[1ex]
\mbox{}\hfill                     
$
\displaystyle
\sum\limits_{k=s+1}^{m}
\gamma_{k}^{}\;\!M_k^{}(t,x)=
{}-\;\!{\rm div}\;\!{\frak d}(t,x)-
\sum\limits_{\nu=1}^{s}
\gamma_{\nu}^{}\;\!M_\nu^{}(t,x)
$
\ for all} $(t,x)\in\Omega.
\hfill
$
\\[1.75ex]
\indent
{\sl Proof.} 
\vspace{0.25ex}
Using Theorem 8.1, from Property 3.8 under the condition 
$M={}-\;\!{\rm div}\;\!{\frak d},$ we obtain the statement of Property 10.21. $\k$
\vspace{1ex}

{\bf  Property 10.22.}
\vspace{0.5ex}
{\it
Suppose
$
\exp\omega_\nu^{}\in \text{\rm ME}_{_{\Omega}},\ 
\nu=1,\ldots, s,\ s\leqslant m-2,\  
\gamma_j^{}\in\R\backslash\{0\}, 
\linebreak 
j=1,\ldots,m, \
{\rm g}_k^{{}^{\scriptsize \gamma_k^{}}}\!\in C^1\Omega,\ 
k=s+1,\ldots,m.$ Then we claim that
\\[1.5ex]
\mbox{}\hfill
$
\displaystyle
\biggl(\,
\prod\limits_{k=s+1}^{m}
{\rm g}_k^{{}^{\scriptsize \gamma_k^{}}}
\exp\sum\limits_{\nu=1}^{s}
\gamma_\nu^{}\;\!\omega_\nu^{}\;\!,\;\! M\biggr)\in 
\text{\rm J}_{_{\Omega}}
\iff
\biggl(\,
\prod\limits_{k=s+1}^{m}
{\rm g}_k^{{}^{\scriptsize \gamma_k^{}}}\;\!,\;\! M+
\sum\limits_{\nu=1}^{s}\gamma_\nu^{}\;\!{\rm div}\;\!{\frak d}\biggr)\in 
\text{\rm J}_{_{\Omega}}.
\hfill
$
\\[1.25ex]
Moreover, there exist functions $M_k^{}\in C^1\Omega,\ k=s+1,\ldots, m,$
such that the identities {\rm (1.14)} are true and
\\[1ex]
\mbox{}\hfill                     
$
\displaystyle
\sum\limits_{k=s+1}^{m}
\gamma_{k}^{}\;\!M_k^{}(t,x)=
M(t,x)+
\sum\limits_{\nu=1}^{s}
\gamma_{\nu}^{}\;\!{\rm div}\;\!{\frak d}(t,x)
$
\ for all} 
$(t,x)\in\Omega.
\hfill
$
\\[-2.5ex]

\newpage

{\sl Proof.} 
\vspace{0.35ex}
Taking into account Theorem 8.1, from Property 3.8 under the conditions 
\linebreak
$M_{\nu}^{}={}-\;\!{\rm div}\;\!{\frak d},\ \nu=1,\ldots, s,$
we get the statement of Property 10.22. $\k$
\vspace{1.25ex}

{\bf  Property 10.23.}
\vspace{0.5ex}
{\it
Suppose 
$
\exp\omega_\nu^{}\in \text{\rm ME}_{_{\Omega}},\ 
\nu=1,\ldots, s,\ s\leqslant m-2,\  
\gamma_j^{}\in\R\backslash\{0\}, 
\linebreak 
j=1,\ldots,m, \
{\rm g}_k^{{}^{\scriptsize \gamma_k^{}}}\!\in C^1\Omega,\ 
k=s+1,\ldots,m.$ Then we have
\\[2.25ex]
\mbox{}\hfill
$
\displaystyle
\prod\limits_{k=s+1}^{m}
{\rm g}_k^{{}^{\scriptsize \gamma_k^{}}}
\exp\sum\limits_{\nu=1}^{s}
\gamma_\nu^{}\;\!\omega_\nu^{}\in 
\text{\rm M}_{_{\Omega}}
\iff
\biggl(\,
\prod\limits_{k=s+1}^{m}
{\rm g}_k^{{}^{\scriptsize \gamma_k^{}}}\;\!,\;\! 
\biggl(\,
\sum\limits_{\nu=1}^{s}\gamma_\nu^{}-1\biggr)\;\!{\rm div}\;\!{\frak d}\biggr)\in 
\text{\rm J}_{_{\Omega}}.
\hfill
$
\\[2.25ex]
Moreover, there exist functions $M_k^{}\in C^1\Omega,\ k=s+1,\ldots, m,$
such that the identities {\rm (1.14)} are true and
\\[1.5ex]
\mbox{}\hfill                     
$
\displaystyle
\sum\limits_{k=s+1}^{m}
\gamma_{k}^{}\;\!M_k^{}(t,x)=
\biggl(\,
\sum\limits_{\nu=1}^{s}\gamma_\nu^{}-1\biggr)\;\!{\rm div}\;\!{\frak d}(t,x)
$
\ for all} $(t,x)\in\Omega.
\hfill
$
\\[1.75ex]
\indent
{\sl Proof.} 
\vspace{0.25ex}
Using Theorem 8.1, from Property 10.22 under the condition 
$M={}-\;\!{\rm div}\;\!{\frak d},$ we obtain the statement of Property 10.23. $\k$
\vspace{1.25ex}

{\bf  Property 10.24.}
\vspace{0.5ex}
{\it
Let 
$
\bigl(\exp\omega_\nu^{}, M_\nu^{}\bigr)\in \text{\rm E}_{_{\Omega}},\ 
\nu=1,\ldots, s,\ s\leqslant m-1,\  
({\rm g}_\xi^{}, M_\xi^{})\in \text{\rm J}_{_{\Omega}},
\linebreak 
{\rm g}_\xi^{{}^{\scriptsize \gamma_\xi^{}}}\!\in C^1\Omega,
\ 
\xi=s+1,\ldots,m,\ 
\gamma_j^{}\in\R\backslash\{0\}, 
\ j=1,\ldots,m. 
$ Then
\\[1.75ex]
\mbox{}\hfill
$
\displaystyle
\prod\limits_{\xi=s+1}^{m}
{\rm g}_\xi^{{}^{\scriptsize \gamma_\xi^{}}}
\exp\sum\limits_{\nu=1}^{s}
\gamma_\nu^{}\;\!\omega_\nu^{}\in 
\text{\rm M}_{_{\Omega}}
\hfill
$
\\[1.75ex]
if and only if the identity {\rm (8.1)} holds.
}
\vspace{0.5ex}

{\sl Proof.} 
\vspace{0.25ex}
Using Theorem 8.1, from Property 3.9 under the condition 
$M={}-\;\!{\rm div}\;\!{\frak d},$ we get the statement of Property 10.24. $\k$
\vspace{1ex}

{\bf  Property 10.25.}
\vspace{0.5ex}
{\it
Let
$
\exp\omega_\nu^{}\in \text{\rm ME}_{_{\Omega}},\ 
\nu=1,\ldots, s,\ s\leqslant m-1,\  
({\rm g}_\xi^{}, M_\xi^{})\in \text{\rm J}_{_{\Omega}},
\linebreak 
{\rm g}_\xi^{{}^{\scriptsize \gamma_\xi^{}}}\!\in C^1\Omega,
\ 
\xi=s+1,\ldots,m,\ 
\gamma_j^{}\in\R\backslash\{0\}, 
\ j=1,\ldots,m. 
$ Then
\\[1.75ex]
\mbox{}\hfill
$
\displaystyle
\biggl(\,
\prod\limits_{\xi=s+1}^{m}
{\rm g}_\xi^{{}^{\scriptsize \gamma_\xi^{}}}
\exp\sum\limits_{\nu=1}^{s}
\gamma_\nu^{}\;\!\omega_\nu^{}\;\!,\;\! M\biggr)\in 
\text{\rm J}_{_{\Omega}}
\hfill
$
\\[1.75ex]
if and only if the following identity holds
\\[1.5ex]
\mbox{}\hfill                      
$
\displaystyle
M(t,x)=
\sum\limits_{\xi=s+1}^{m}
\gamma_\xi^{}\;\!M_{\xi}^{}(t,x)-
\sum\limits_{\nu=1}^{s}
\gamma_\nu^{}\;\!{\rm div}\;\!{\frak d}(t,x)
$
\ for all} 
$(t,x)\in\Xi^{\;\!\prime}.
\hfill
$
\\[1.5ex]
\indent
{\sl Proof.} 
\vspace{0.25ex}
Taking into account Theorem 8.1, from Property 3.9 under the conditions 
\linebreak
$M_{\nu}^{}={}-\;\!{\rm div}\;\!{\frak d},\ \nu=1,\ldots, s,$
\vspace{1.25ex}
it follows that the statement of Property 10.25 is true. $\k$

{\bf  Property 10.26.}
\vspace{0.5ex}
{\it
Suppose 
$
\exp\omega_\nu^{}\in \text{\rm ME}_{_{\Omega}},\ 
\nu=1,\ldots, s,\ s\leqslant m-1,\  
({\rm g}_\xi^{}, M_\xi^{})\in \text{\rm J}_{_{\Omega}},
\linebreak 
{\rm g}_\xi^{{}^{\scriptsize \gamma_\xi^{}}}\!\in C^1\Omega,
\ 
\xi=s+1,\ldots,m,\ 
\gamma_j^{}\in\R\backslash\{0\}, 
\ j=1,\ldots,m. 
$ Then
\\[1.75ex]
\mbox{}\hfill
$
\displaystyle
\prod\limits_{\xi=s+1}^{m}
{\rm g}_\xi^{{}^{\scriptsize \gamma_\xi^{}}}
\exp\sum\limits_{\nu=1}^{s}
\gamma_\nu^{}\;\!\omega_\nu^{}\in 
\text{\rm M}_{_{\Omega}}
\hfill
$
\\[1.75ex]
if and only if the following identity holds
\\[1.75ex]
\mbox{}\hfill                      
$
\displaystyle
\sum\limits_{\xi=s+1}^{m}
\gamma_\xi^{}\;\!M_{\xi}^{}(t,x)=
\biggl(\,\sum\limits_{\nu=1}^{s}\gamma_\nu^{}-1\biggr)
\;\!{\rm div}\;\!{\frak d}(t,x)
$
\ for all} 
$(t,x)\in\Xi^{\;\!\prime}.
\hfill
$
\\[1.75ex]
\indent
{\sl Proof.} 
\vspace{0.25ex}
Using Theorem 8.1, from Property 10.25 under the condition 
$M={}-\;\!{\rm div}\;\!{\frak d},$ we obtain the statement of Property 10.26. $\k$
\vspace{1ex}

{\bf Property 10.27.} 
\vspace{0.75ex}
{\it
Suppose functions
$p,q\in \text{\rm P}_{_{\!\Xi^{\;\!\prime}}}$
are relatively prime, the set
$\Omega_{_0}\subset\Xi^{\;\!\prime}$ such that   
$
p(t,x)\ne 0$ for all $(t,x)\in\Omega_{_0}$ and
$p(t,x)= 0$ for all $(t,x)\in {\sf C}_{_{\!\Xi^{\;\!\prime}}}\Omega_{_0}\;\!.
$
Then
\\[2ex]
\mbox{}\hfill
$
\exp\dfrac{q}{p}\in 
\text{\rm ME}_{_{\Omega_{{}_{\tiny\;\! 0}}}}
\hfill
$
\\[2ex]
if and only if  
$(p,M)\in \text{\rm A}_{_{\!\Xi^{\;\!\prime}}}$
and the identity holds
\\[2ex]
\mbox{}\hfill                          
$
{\frak d}\;\!q(t,x)=
q(t,x)\;\!M(t,x)-p(t,x)\;\!{\rm div}\;\!{\frak d}(t,x)
$
\ for all} $(t,x)\in \Xi^{\;\!\prime}.
\hfill
$
\\[2ex]
\indent
{\sl Proof.} 
\vspace{0.25ex}
Taking into account Theorem 8.1, from Property 3.10 under the condition 
\linebreak
$M={}-\;\!{\rm div}\;\!{\frak d},$ we get the statement of Property 10.27. $\k$
\vspace{1.25ex}

{\bf Property 10.28.} 
\vspace{0.75ex}
{\it
Suppose functions  
$p,q\in \text{\rm P}_{_{\!\Xi^{\;\!\prime}}}$
are relatively prime, the set
$\Omega_{_0}\subset\Xi^{\;\!\prime}$ such that   
$p(t,x)\ne 0$ for all $(t,x)\in\Omega_{_0}$ and
$p(t,x)= 0$ for all $(t,x)\in {\sf C}_{_{\!\Xi^{\;\!\prime}}}\Omega_{_0}\;\!.
$
Then
\\[2ex]
\mbox{}\hfill
$
\exp\arctan\dfrac{q}{p}\in 
\text{\rm ME}_{_{\Omega_{{}_{\tiny\;\! 0}}}}
\hfill
$
\\[2ex]
if and only if there exists a function
\vspace{0.5ex}
$M\in \text{\rm P}_{_{\!\Xi^{\;\!\prime}}}$
such that the degree $\deg_{x}^{}M\leq d-1$
and the following system of identities holds
\\[1.75ex]
\mbox{}\hfill                           
$
{\frak d}\;\!p(t,x)=
p(t,x)\;\!M(t,x)+q(t,x)\;\!{\rm div}\;\!{\frak d}(t,x)
$
\ for all $(t,x)\in \Xi^{\;\!\prime},
\hfill
$
\\[2.5ex]
\mbox{}\hfill                           
$
{\frak d}\;\!q(t,x)=
q(t,x)\;\!M(t,x)-p(t,x)\;\!{\rm div}\;\!{\frak d}(t,x)
$
\ for all} $(t,x)\in \Xi^{\;\!\prime}.
\hfill
$
\\[2ex]
\indent
{\sl Proof}. 
\vspace{0.15ex}
Using Remark 6.2 and Theorem 8.1, from Property 3.11 
it follows that the statement of Property 10.28 is true. $\k$
\vspace{1.25ex}

{\bf Property 10.29.} 
{\it
Let  
$p\in \text{\rm P}_{_{\!\Xi^{\;\!\prime}}}.$
Then we have
\\[2ex]
\mbox{}\hfill                       
$
\exp\arctan p\in
\text{\rm ME}_{_{\Xi^{\;\!\prime}}}
\hfill
$
\\[1.25ex]
if and only if 
\\[1.5ex]
\mbox{}\hfill                           
$
{\frak d}\;\!p(t,x)=
{}-\bigl(1+p^2(t,x)\bigr)\;\!{\rm div}\;\!{\frak d}(t,x)
$
\ for all} $(t,x)\in \Xi^{\;\!\prime}.
\hfill
$
\\[2.25ex]
\indent
{\sl Proof}. 
\vspace{0.35ex}
Using Theorem 8.1, from Property 3.12 under the condition 
$M={}-\;\!{\rm div}\;\!{\frak d},$ we obtain the statement of Property 10.29. $\k$
\vspace{1ex}

{\bf Property 10.30.} 
{\it
Let $f\in C^1\Omega.$ Then we claim that
\\[2ex]
\mbox{}\hfill                       
$
\exp\arctan f\in \text{\rm ME}_{_{\Omega}}
\iff
\bigl(\exp{\rm arccotan} f,\, {\rm div}\;\!{\frak d}\bigr)\in
\text{\rm E}_{_{\Omega}}
\hfill
$
\\[0.75ex]
and}
\\[0.75ex]
\mbox{}\hfill                       
$
\exp{\rm arccotan} f\in \text{\rm ME}_{_{\Omega}}
\iff
\bigl(\exp\arctan f,\, -\;\!{\rm div}\;\!{\frak d}\bigr)\in
\text{\rm E}_{_{\Omega}}.
\hfill
$
\\[2.75ex]
\indent
{\bf Property 10.31.} 
\vspace{0.75ex}
{\it
Suppose  
$k\in\N,\ F\in \text{\rm I}_{_{\Omega}},$
a function $f\in C^1\Omega$ is relatively prime to the function $F,$ 
the set $\Omega_{_0}\subset\Omega$ such that  
\vspace{0.75ex}
$
F(t,x)\ne 0$ for all $(t,x)\in\Omega_{_0}$ and
$F(t,x)= 0$ for all $(t,x)\in{\sf C}_{_{\Omega}}\Omega_{_0}.$
Then we have
\\[2ex]
\mbox{}\hfill                       
$
\exp\dfrac{f}{F^k}\in 
\text{\rm ME}_{_{\Omega_{{}_{\tiny\;\! 0}}}}
\hfill
$
\\[2ex]
if and only if the following identity holds
\\[2ex]
\mbox{}\hfill                      
$
{\frak d}\;\!f(t,x)={}-
F^k(t,x)\;\!{\rm div}\;\!{\frak d}(t,x)
$
\ for all} $(t,x)\in \Omega.
\hfill
$
\\[2ex]
\indent
{\sl Proof.} 
\vspace{0.35ex}
Taking into account Theorem 8.1, from Property 3.14 under the condition 
\linebreak
$M={}-\;\!{\rm div}\;\!{\frak d}$ it follows that the statement of Property 10.31 is true. $\k$

\newpage

\mbox{}
\\[-1.75ex]
\centerline{
{\bf\large \S\;\!3. First integrals}}
\\[2ex]
\centerline{\bf  11. First integrals defined by partial integrals and last multipliers}
\\[1.5ex]
\indent
{\bf Theorem 11.1}
(criterion of first integral).
{\it 
A partial integral on the domain $\Omega$ of system {\rm (0.1)} is 
a first integral on the domain $\Omega$ of system {\rm (0.1)}
if and only if 
the cofactor of the partial integral is identically equal to zero on the domain $\Omega.$
}
\vspace{0.35ex}

{\sl Proof}\;\!
is based on the definition of partial integral (Definition 1.1) and 
the existence criterion of first integral (Theorem 0.1).
Indeed, the identity (1.1) under the conditions 
${\rm g}=F$ and $M=0$ coincides with the identity (0.2). $\k$
\vspace{0.75ex}

Thus, we have 
$\text{\rm I}_{_{\Omega}}\subset \text{\rm J}_{_{\Omega}}$ 
and the statement of Theorem 11.1 is expressed by the equivalence
\\[1ex]
\mbox{}\hfill
$
F\in \text{\rm I}_{_{\Omega}}
\iff
(F,\;\! 0)\in \text{\rm J}_{_{\Omega}}.
\hfill
$
\\[1.75ex]
\indent
{\bf Property 11.1.}
{\it 
We claim that
\\[1.5ex]
\mbox{}\hfill
for all}
$C\in\R
\ \ 
\Longrightarrow
\ \ 
C\in\text{\rm I}_{_{\Xi}}.
\hfill
$
\\[1.5ex]
\indent
{\sl Indeed},
since ${\frak d}\;\!C=0$ for all $(t,x)\in\Xi,$
by Theorem 0.2, we have $C\in\text{\rm I}_{_{\Xi}}.\ \k$
\vspace{1ex}

{\bf Property 11.2.}
{\it 
If $C\in\R,$ then}
\\[1.5ex]
\mbox{}\hfill
$
F\in \text{\rm I}_{_{\Omega}}
\iff
F+C\in \text{\rm I}_{_{\Omega}}.
\hfill
$
\\[1.5ex]
\indent
{\sl Proof.}\!
\vspace{1ex}
Using Property 11.1, from Property 0.1, we get the statement of Property 11.2.$\k$

{\bf Theorem 11.2}
(geometric sense of first integral).
{\it 
A first integral $F$ on the domain $\Omega$ of system {\rm (0.1)} 
defines the family of integral manifolds $F(t,x)=C$ of system {\rm (0.1)},
where $C$ is  arbitrary constant from the range of function $F$ 
{\rm(}cases $C={}-\infty,\ C={}+\infty,$ and $C=\infty$ are not excluded{\rm)}. 
}
\vspace{0.35ex}

{\sl Proof}\;\!
follows from 
the definition of integral manifold (Definition 0.3), 
the existence criterion of first integral (Theorem 0.1), and Property 11.2.  $\k$
\vspace{1ex}

{\bf Property 11.3.}
{\it 
A linear combination over the field $\R$ 
of first integrals on the domain $\Omega$ of system {\rm (0.1)} is 
a first integral on the domain $\Omega$ of system {\rm (0.1)}, i.e.,} 
\\[1.5ex]
\mbox{}\hfill
$
\displaystyle
\lambda_j^{}\in\R,
\ \ 
F_j^{}\in \text{\rm I}_{_{\Omega}},
\ j=1,\ldots,m,
\ \ \Longrightarrow\ \ 
\sum\limits_{j=1}^{m}
\lambda_j^{}\;\!F_j^{}\in \text{\rm I}_{_{\Omega}}.
\hfill
$
\\[1.5ex]
\indent
{\sl Proof}\;\!
follows from Property 0.1. $\k$
\vspace{1ex}

Properties 11.399 and 11.1 define the mathematical structure of the set of first integrals for system (0.1).
\vspace{0.5ex}

{\bf Theorem 11.3.}
\vspace{0.35ex}
{\it 
The set of first integrals on the domain $\Omega$ of system {\rm (0.1)} 
is a linear space over the field of real numbers}:
\vspace{1ex}
$\bigl(\text{\rm I}_{_{\Omega}},\,\R,\,+\;\!,\,\cdot\;\!,\,=\bigr).$

{\bf Property 11.4.} 
\vspace{0.25ex}
{\it
Suppose the set $\Omega_{_0}\subset\Omega$ such that   
$F(t,x)\ne 0$ for all $(t,x)\in \Omega_{_0}$ and 
$F(t,x)=0$ for all $(t,x)\in {\sf C}_{_\Omega}\Omega_{_0}.$
Then}
\\[1.5ex]
\mbox{}\hfill  
$
F\in \text{I}_{_{\Omega}}
\iff 
|F|\in \text{I}_{_{\Omega_{{}_{\tiny\;\! 0}}}}.
\hfill
$
\\[1.5ex]
\indent
{\sl Proof.}\!
\vspace{1ex}
Using Theorem 11.1, from Property 1.3, we get the statement of Property 11.4.$\k$

{\bf Property 11.5.} 
{\it
If $F\in \text{\rm I}_{_{\Omega}}$ and 
\\[1.5ex]
\mbox{}\hfill                           
$
{\frak d}\;\!f(t,x)=
\bigl(f(t,x)+F(t,x)\bigr)\;\!M(t,x)
$
\ for all $(t,x)\in \Omega,
\hfill
$
\\[1.5ex]
where the function $M\in \text{\rm P}_{_{\!\Xi^{\;\!\prime}}}$
\vspace{0.5ex}
such that the degree $\deg_{x}^{}M\leq d-1,$ 
then $(f+F, M)\in \text{\rm J}_{_{\Omega}}.$
}

{\sl Proof.} 
Taking into account Theorem 0.2, we obtain
\\[1.5ex]
\mbox{}\hfill                           
$
{\frak d}\;\!\bigl(f(t,x)+F(t,x)\bigr)=
{\frak d}\;\!f(t,x)=
\bigl(f(t,x)+F(t,x)\bigr)\;\!M(t,x)
$
\ for all $(t,x)\in \Omega.
\hfill
$
\\[1.5ex]
\indent
By Theorem 1.1, $(f+F, M)\in \text{\rm J}_{_{\Omega}}.\ \k$
\vspace{1.25ex}

{\bf Property 11.6.} 
{\it
Let $F\in \text{\rm I}_{_{\Omega}}.$ Then} 
\\[1.5ex]
\mbox{}\hfill                           
$
({\rm g}\;\!F,\;\! M)\in \text{\rm J}_{_{\Omega}}
\iff
({\rm g},\;\! M)\in \text{\rm J}_{_{\Omega}}.
\hfill
$
\\[1.5ex]
\indent
{\sl Proof.} 
Using the identity (0.3), we get 
\\[1.5ex]
\mbox{}\hfill                           
$
{\frak d}\;\!\bigl({\rm g}(t,x)\;\!F(t,x)\bigr)=
F(t,x)\, {\frak d}\;\!{\rm g}(t,x)
$
\ for all $(t,x)\in \Omega.
\hfill
$
\\[1.5ex]
\indent
Now, if we apply Theorem 1.1 to the functions ${\rm g}\;\!F$ and ${\rm g},$ then 
from this identity, we obtain the statement of Property 11.6. $\k$
\vspace{1.25ex}

{\bf Property 11.7.} 
{\it
Suppose $({\rm g},\;\! M)\in \text{\rm J}_{_{\Omega}}.$ Then we have} 
\\[1.5ex]
\mbox{}\hfill                           
$
({\rm g}\;\!F,\;\! M)\in \text{\rm J}_{_{\Omega}}
\iff
F\in \text{\rm I}_{_{\Omega}}.
\hfill
$
\\[1.5ex]
\indent
{\sl Proof.} 
Taking into account the identity (1.2), from Theorem 1.1 it follows that 
the derivative by virtue of the differential system (0.1)
\\[1.5ex]
\mbox{}\hfill                           
$
{\frak d}\;\!\bigl({\rm g}(t,x)\;\!F(t,x)\bigr)=
{\rm g}(t,x)\;\!
\bigl(F(t,x)\;\!M(t,x)+ {\frak d}\;\!F(t,x)\bigr)
$
\ for all $(t,x)\in \Omega.
\hfill
$
\\[1.5ex]
\indent
If $({\rm g}\;\!F,\;\! M)\in \text{\rm J}_{_{\Omega}},$ then
\\[1.5ex]
\mbox{}\hfill                           
$
F(t,x)\;\!{\rm g}(t,x)\;\!M(t,x)=
{\rm g}(t,x)\;\!
\bigl(F(t,x)\;\!M(t,x)+ {\frak d}\;\!F(t,x)\bigr)
$
\ for all $(t,x)\in \Omega.
\hfill
$
\\[1.5ex]
This yields that ${\frak d}\;\!F(t,x)=0$ for all $(t,x)\in \Omega.$
\vspace{0.5ex}
Therefore, by Theorem 0.2, we get $F\in \text{\rm I}_{_{\Omega}}.$

If $F\in \text{\rm I}_{_{\Omega}},$ then 
\\[1.5ex]
\mbox{}\hfill                           
$
{\frak d}\;\!\bigl({\rm g}(t,x)\;\!F(t,x)\bigr)=
{\rm g}(t,x)\;\!F(t,x)\;\!M(t,x)
$
\ for all $(t,x)\in \Omega.
\hfill
$
\\[1.5ex]
\indent
By Theorem 1.1, we obtain $({\rm g}\;\!F,\;\! M)\in \text{\rm J}_{_{\Omega}}.\ \k$
\vspace{1ex}

{\bf Corollary 11.1.} 
{\it
We claim that}\,\!: 
\\[1.5ex]
\mbox{}\hfill                           
$
F\in \text{\rm I}_{_{\Omega}}
\ \ \Longrightarrow\ \ 
\bigl(\mu\;\!F\in \text{\rm M}_{_{\Omega}}
\iff 
\mu\in \text{\rm M}_{_{\Omega}}\bigr);
\hfill                           
$
\\[2ex]
\mbox{}\hfill                           
$
\mu\in \text{\rm M}_{_{\Omega}}
\ \ \Longrightarrow\ \ 
\bigl(\mu\;\!F\in \text{\rm M}_{_{\Omega}}
\iff 
F\in \text{\rm I}_{_{\Omega}}\bigr).
\hfill
$
\\[2ex]
\indent
{\bf Property 11.8.} 
{\it
Let $({\rm g},\;\!\varphi)\in \text{\rm J}_{_{\Omega}},\ 
\varphi\in C^1T^{\;\!\prime}.$ Then the function 
\\[1.5ex]
\mbox{}\hfill                           
$
\displaystyle
F\colon (t,x)\to\ 
{\rm g}(t,x)\exp\biggl({}-\int\limits_{t_0^{}}^{t}\varphi(\tau)\,d\tau\biggr)
$
\ for all $(t,x)\in \Omega,
\hfill
$
\\[1.5ex]
where $t_{_0}^{}$ is an arbitrary fixed point from the domain $T^{\;\!\prime},$
\vspace{0.25ex}
is a first integral on the domain $\Omega$ of the differential system {\rm (0.1)}.
}
\vspace{0.5ex}

{\sl Proof.} 
The derivative by virtue of system (0.1)
\\[1.5ex]
\mbox{}\hfill                           
$
\displaystyle
{\frak d}\;\!F(t,x)=
{\frak d}\;\!
\biggl({\rm g}(t,x)\exp\biggl({}-\int\limits_{t_0^{}}^{t}\varphi(\tau)\,d\tau\biggr)\biggr)=
\hfill                           
$
\\[1.5ex]
\mbox{}\hfill                           
$
\displaystyle
=\exp\biggl({}-\int\limits_{t_0^{}}^{t}\varphi(\tau)\,d\tau\biggr)\;\!
{\frak d}\;\!{\rm g}(t,x)\;\!+\;\!
{\rm g}(t,x)\, {\sf D}\exp\biggl({}-\int\limits_{t_0^{}}^{t}\varphi(\tau)\,d\tau\biggr)=
\hfill                           
$
\\[1.5ex]
\mbox{}\hfill                           
$
\displaystyle
=\varphi(t)\;\! {\rm g}(t,x)\exp\biggl({}-\int\limits_{t_0^{}}^{t}\varphi(\tau)\,d\tau\biggr)\;\!-\;\!
\varphi(t)\;\! {\rm g}(t,x)\exp\biggl({}-\int\limits_{t_0^{}}^{t}\varphi(\tau)\,d\tau\biggr)=0
$
\ for all $(t,x)\in \Omega.
\hfill
$
\\[1.75ex]
\indent
By Theorem 0.2, we have $F\in \text{\rm I}_{_{\Omega}}.\ \k$
\vspace{1.5ex}

{\bf Corollary 11.2.} 
{\it
If $({\rm g},\;\!\lambda)\in \text{\rm J}_{_{\Omega}},\ 
\lambda \in\R\backslash\{0\},$ then 
${\rm g}\;\!e^{{}-\lambda\;\!t}\in \text{\rm I}_{_{\Omega}}.$
} 
\vspace{1.5ex}

{\bf  Property 11.9.}
{\it
Let $\lambda_j^{}\in\R, \ j=1,\ldots,m,\ m\leq n.$ Then
$\sum\limits_{j=1}^{m}\lambda_j^{}\;\!x_j^{}\in \text{\rm I}_{_{\Xi}}$
if and only if
\\[1.5ex]
\mbox{}\hfill
$
\displaystyle
\sum\limits_{j=1}^{m}
\lambda_j^{}\;\!X_j^{}(t,x)=0
$
\ for all} $(t,x)\in \Xi.
\hfill
$
\\[1.5ex]
\indent
{\sl Proof.}\!
\vspace{1.5ex}
Using Theorem 11.1, from Property 2.2, we get the statement of Property 11.9.$\k$

{\bf Property 11.10.} 
{\it
Let $({\rm g}_j^{},\;\!M_j^{})\in \text{\rm J}_{_{\Omega}},\ 
\lambda_j^{}\in\R\backslash\{0\}, \ j=1,\ldots,m.$ Then
$\sum\limits_{j=1}^{m}\lambda_j^{}\;\!{\rm g}_j^{}\in \text{\rm I}_{_{\Omega}}$
if and only if the following identity holds
\\[1.5ex]
\mbox{}\hfill
$
\displaystyle
\sum\limits_{j=1}^{m}
\lambda_j^{}\;\!{\rm g}_j^{}(t,x)\;\!M_j^{}(t,x)=0
$
\ for all} $(t,x)\in \Omega.
\hfill
$
\\[1.5ex]
\indent
{\sl Proof}\;\!
is based on Theorem 11.1 and consists in the fact that
\\[1.5ex]
\mbox{}\hfill
$
\displaystyle
{\frak d}\;\!\sum\limits_{j=1}^{m}
\lambda_j^{}\;\!{\rm g}_j^{}(t,x)=
\sum\limits_{j=1}^{m}
\lambda_j^{}\;\!{\rm g}_j^{}(t,x)\;\!M_j^{}(t,x)
$
\ for all $(t,x)\in \Omega. \ \k
\hfill
$
\\[2ex]
\indent
{\bf Corollary 11.3.} 
{\it
Suppose $({\rm g}_j^{},\;\!\rho_j^{}\;\!M_{_0}^{})\in \text{\rm J}_{_{\Omega}},\ 
\rho_j^{},\lambda_j^{}\in\R\backslash\{0\}, \ j=1,\ldots,m.$ Then}
\\[1.5ex]
\mbox{}\hfill
$
\displaystyle
\sum\limits_{j=1}^{m}\lambda_j^{}\;\!{\rm g}_j^{}\in \text{\rm I}_{_{\Omega}}
\iff
\sum\limits_{j=1}^{m}
\rho_j^{}\;\!\lambda_j^{}=0.
\hfill
$
\\[2ex]
\indent
{\bf Corollary 11.4.} 
{\it
Suppose $\mu_j^{}\in \text{\rm M}_{_{\Omega}},\ 
\lambda_j^{}\in\R\backslash\{0\}, \ j=1,\ldots,m.$ Then}
\\[1.5ex]
\mbox{}\hfill
$
\displaystyle
\sum\limits_{j=1}^{m}\lambda_j^{}\;\!\mu_j^{}\in \text{\rm I}_{_{\Omega}}
\iff
\sum\limits_{j=1}^{m}
\lambda_j^{}=0.
\hfill
$
\\[2ex]
\indent
{\bf Theorem 11.4.}
\vspace{0.25ex}
{\it
The continuously differentiable function {\rm (1.4)}
is a first integral on the domain $\Omega$ of system {\rm (0.1)}
\vspace{0.25ex}
if and only if there exist functions
$M_j^{}\in C^1\Omega,\ j=1,\ldots,m,$ 
such that these functions satisfies the identities {\rm (1.5)} and
\\[1.5ex]
\mbox{}\hfill
$
\displaystyle
\sum\limits_{j=1}^{m} M_j^{}(t,x)= 0
$
\ for all} $(t,x)\in\Omega.
\hfill
$
\\[1.5ex]
\indent
{\sl Proof.}\!
\vspace{1ex}
Using Theorem 11.1, from Theorem 1.3, we get the statement of Theorem 11.4.$\k$

{\bf Property 11.11.} 
\vspace{0.5ex}
{\it
Suppose 
$({\rm g}_j^{}, M_{j}^{})\in \text{\rm J}_{_{\Omega}},\ 
\gamma_j^{}\in\R\backslash\{0\},\
{\rm g}_j^{{}^{\scriptsize \gamma_{j}^{}}}\in C^1\Omega, \ j=1,\ldots,m.$ 
Then  
$\prod\limits_{j=1}^{m} 
{\rm g}_j^{{}^{\scriptsize \gamma_{j}^{}}} 
\in \text{\rm I}_{_{\Omega}}$
if and only if 
the linear combination of cofactors
\\[1ex]
\mbox{}\hfill                      % (11.1)
$
\displaystyle
\sum\limits_{j=1}^{m}
\gamma_j^{}\;\!M_{j}^{}(t,x)=0
$
\ for all} $(t,x)\in\Xi^{\;\!\prime}.
$
\hfill (11.1)
\\[1.5ex]
\indent
{\sl Proof.}
Taking into account Theorem 11.1, 
from Property 1.9 under the condition $M=0$
it follows that the statement of Property 11.11 is true. $\k$
\vspace{0.5ex}

Property 11.11 also follows from Theorem 11.4 if we take into account Property 1.6.
\vspace{1ex}

{\bf Corollary 11.5.} 
\vspace{0.5ex}
{\it
Suppose 
$({\rm g}_j^{}, \rho_{j}^{}\;\!M_{_0}^{})\in \text{\rm J}_{_{\Omega}},\ 
\rho_{j}^{},\gamma_j^{}\in\R\backslash\{0\},\
{\rm g}_j^{{}^{\scriptsize \gamma_{j}^{}}}\in C^1\Omega, \ j=1,\ldots,m.$ 
Then  
$\prod\limits_{j=1}^{m} 
{\rm g}_j^{{}^{\scriptsize \gamma_{j}^{}}} 
\in \text{\rm I}_{_{\Omega}}$
if and only if}
$
\sum\limits_{j=1}^{m}
\rho_{j}^{}\;\!\gamma_j^{}=0.
$
\vspace{1ex}

Let us remark that 
if the conditions of Corollary 11.5 are fulfilled, then pairwise taken
partial integrals are first integrals
\\[1.5ex]
\mbox{}\hfill                   
$
{\rm g}_\xi^{{}^{\scriptsize \gamma_{\xi}^{}}}\;\!
{\rm g}_\zeta^{{}^{\scriptsize \gamma_{\zeta}^{}}} 
\in \text{\rm I}_{_{\Omega}}
\iff
\rho_{\xi}^{}\;\!\gamma_\xi^{}+
\rho_{\zeta}^{}\;\!\gamma_\zeta^{}=0,
\ \ \xi,\zeta=1,\ldots, m,\ \ \xi\ne\zeta.
\hfill
$
\\[2.5ex]
\indent
{\bf Example 11.1.}
The differential system
\\[2ex]
\mbox{}\hfill
$
\dfrac{dx}{dt}=1,
\qquad
\dfrac{dy}{dt}={}-2\;\!xy+z^2,
\qquad
\dfrac{dz}{dt}={}-2\;\!xz
\hfill
$
\\[2ex]
has the autonomous polynomial partial integral 
\\[1.5ex]
\mbox{}\hfill
$
{\rm g}_1^{}\colon (x,y,z)\to\ z
$
\ for all $(x,y,z)\in\R^3
\hfill
$
\\[2ex]
with cofactor 
$
M_1^{}\colon (x,y,z)\to {}-2\;\!x
$ for all $(x,y,z)\in\R^3,$ 
the autonomous partial integral 
\\[2ex]
\mbox{}\hfill
$
\displaystyle
{\rm g}_2^{}\colon (x,y,z)\to\ 
y-z^2e^{\;\!x^2}\int\limits_{0}^{x}e^{\;\!-\;\!\tau^2}\;\!d\tau
$
\ for all $(x,y,z)\in\R^3
\hfill
$
\\[1.5ex]
with cofactor  
$
M_2^{}\colon (x,y,z)\to {}-2\;\!x
$ for all $(x,y,z)\in\R^3,$
and the autonomous conditional partial integral 
\\[1.5ex]
\mbox{}\hfill
$
\displaystyle
{\rm g}_3^{}\colon (x,y,z)\to\ 
e^{\;\!x^2}
$
\ for all $(x,y,z)\in\R^3
\hfill
$
\\[2ex]
with cofactor  
$
M_3^{}\colon (x,y,z)\to\, 2\;\!x
$ for all $(x,y,z)\in\R^3.
$
\vspace{1ex}

By Corollary 11.5 
\vspace{0.35ex}
(under $\rho_1^{}={}-1,\ \rho_2^{}=1,\ M_{_0}^{}=2\;\!x,\ \gamma_1^{}=\gamma_2^{}=1),$
on the base of partial integrals ${\rm g}_1^{}$ and ${\rm g}_3^{},$ 
${\rm g}_2^{}$ and ${\rm g}_3^{}$ we can build the autonomous first integrals
 \\[2ex]
\mbox{}\hfill
$
\displaystyle
F_1^{}={\rm g}_1^{}\;\!{\rm g}_3^{}\colon (x,y,z)\to\ 
z\;\!e^{\;\!x^2}
$
\ for all 
$(x,y,z)\in\R^3
\hfill
$
\\[0.5ex]
and
\\[0.5ex]
\mbox{}\hfill
$
\displaystyle
F_2^{}={\rm g}_2^{}\;\!{\rm g}_3^{}\colon (x,y,z)\to\ 
\Bigl(y-z^2e^{\;\!x^2}\int\limits_{0}^{x}e^{\;\!-\;\!\tau^2}\;\!d\tau\Bigr)\;\!
e^{\;\!x^2}
$
\ for all $(x,y,z)\in\R^3.
\hfill
$
\\[1.5ex]
\indent
Since the functions $F_1^{}$ and $F_2^{}$ are functionally independent, 
we see that these functions are an autonomous integral basis on $\R^3.$
\vspace{0.5ex}

The function
\\[1ex]
\mbox{}\hfill
$
\displaystyle
F_3^{}\colon (t,x,y,z)\to\ x-t
$
\ for all $(t,x,y,z)\in\R^4
\hfill
$
\\[1.5ex]
is an nonautonomous first integral.
\vspace{0.5ex}

The set $\{F_1^{},\;\! F_2^{},\;\! F_3^{}\}$ is a basis of first integrals on $\R^4.$
\vspace{1.5ex}

{\bf Corollary 11.6.} 
\vspace{0.5ex}
{\it
Let 
$\mu_j^{}\in \text{\rm M}_{_{\Omega}},\ 
\gamma_j^{}\in\R\backslash\{0\},\
\mu_j^{{}^{\scriptsize \gamma_{j}^{}}}\in C^1\Omega, \ j=1,\ldots,m.$ 
Then the function
$\prod\limits_{j=1}^{m} 
\mu_j^{{}^{\scriptsize \gamma_{j}^{}}} 
\in \text{\rm I}_{_{\Omega}}$
if and only if}
$
\sum\limits_{j=1}^{m}
\gamma_j^{}=0.
$
\vspace{1ex}

{\bf Corollary 11.7.}
 \vspace{0.35ex}
 {\it 
 Suppose the set $\Omega_{_0}\subset \Omega$ such that 
${\rm g}_2^{}(t,x)\ne 0$ for all $(t,x)\in \Omega_{_0}$ and
${\rm g}_2^{}(t,x)= 0$ for all $(t,x)\in {\sf C}_{_\Omega}\Omega_{_0}.$
Then} 
\\[1.5ex]
\mbox{}\hfill
$
({\rm g}_1^{}, M)\in \text{\rm J}_{_{\Omega}}
\ \ \&\ \ 
({\rm g}_2^{}, M)\in \text{\rm J}_{_{\Omega}}
\ \ \Longrightarrow\ \ 
\dfrac{{\rm g}_1^{}}{{\rm g}_2^{}}\in
\text{\rm I}_{_{ \Omega_{_0}}}.
\hfill
$
\\[2ex]
\indent
Note that if $M={}-\;\!{\rm div}\;\!{\frak d},$ then from Corollary 11.7, 
we get the Jacobi property for last multipliers (Property 0.2).
\vspace{1ex}

{\bf Property 11.12.} 
\vspace{0.5ex}
{\it
Suppose 
$({\rm g}_\tau^{}, M_{\tau}^{})\in \text{\rm J}_{_{\Omega}},\ 
\tau=1,\ldots, m-1,\
\gamma_j^{}\in\R\backslash\{0\},\
{\rm g}_j^{{}^{\scriptsize \gamma_{j}^{}}}\in C^1\Omega, 
\linebreak 
j=1,\ldots,m.$ 
Then we claim that}
\\[1.25ex]
\mbox{}\hfill                      
$
\displaystyle
\prod\limits_{j=1}^{m} 
{\rm g}_j^{{}^{\scriptsize \gamma_{j}^{}}}
\in \text{\rm I}_{_{\Omega}}
\iff
\biggl({\rm g}_m^{},
{}-\dfrac{1}{\gamma_m^{}}\;\!
\sum\limits_{\tau=1}^{m-1}
\gamma_{\tau}\;\!M_{\tau}^{}\biggl) 
\in \text{\rm J}_{_{\Omega}}\;\!.
\hfill
$
\\[1.5ex]
\indent
{\sl Indeed}, 
taking into account Theorem 11.1, from Property 1.10, we obtain the statement of Property 11.12. $\k$
\vspace{1ex}

{\bf Corollary 11.8.} 
\vspace{0.5ex}
{\it
Let 
$\mu_j^{}\in \text{\rm M}_{_{\Omega}},\ 
\gamma,\gamma_j^{}\in\R\backslash\{0\},\ 
{\rm g}^{{}^{\scriptsize \gamma}},
\mu_j^{{}^{\scriptsize \gamma_{j}^{}}}\in C^1\Omega, 
\ j=1,\ldots,m.$ 
Then}
\\[1.25ex]
\mbox{}\hfill                      
$
\displaystyle
{\rm g}^{{}^{\scriptsize \gamma}}\;\!
\prod\limits_{j=1}^{m} 
\mu_j^{{}^{\scriptsize \gamma_{j}^{}}}
\in \text{\rm I}_{_{\Omega}}
\iff
\biggl({\rm g},\
\dfrac{1}{\gamma}\;\!
\sum\limits_{j=1}^{m}
\gamma_{j}^{}\,{\rm div}\;\!{\frak d}\biggl) 
\in \text{\rm J}_{_{\Omega}}\;\!.
\hfill
$
\\[1.5ex]
\indent
{\bf Property 11.13.} 
\vspace{0.5ex}
{\it
Suppose 
$({\rm g}_\nu^{}, M_{\nu}^{})\!\in \text{\rm J}_{_{\Omega}},\, 
\nu\!=\!1,\ldots, s,\ s\leq m\!-\!2, \ 
\gamma_j^{}\!\in\R\backslash\{0\},\
{\rm g}_j^{{}^{\scriptsize \gamma_{j}^{}}}\!\in C^1\Omega, 
\linebreak 
j=1,\ldots,m.$ 
Then we have
\\[1.25ex]
\mbox{}\hfill                      
$
\displaystyle
\prod\limits_{j=1}^{m} 
{\rm g}_j^{{}^{\scriptsize \gamma_{j}^{}}} 
\in \text{\rm I}_{_{\Omega}}
\iff
\biggl(
\,\prod\limits_{k=s+1}^{m} 
{\rm g}_k^{{}^{\scriptsize \gamma_{k}^{}}}, 
{}-\sum\limits_{\nu=1}^{s}
\gamma_{\nu}\;\!M_{\nu}^{}\biggl) 
\in \text{\rm J}_{_{\Omega}}\;\!.
\hfill
$
\\[1.5ex]
Moreover, there exist functions
$M_k^{}\in C^1\Omega,\ k=s+1,\ldots, m,$
such that the identities {\rm (1.14)} are true and
\\[1.25ex]
\mbox{}\hfill                      % (11.2)
$
\displaystyle
\sum\limits_{k=s+1}^{m}
\gamma_{k}^{}\;\!M_k^{}(t,x)=
{}-\sum\limits_{\nu=1}^{s}
\gamma_{\nu}^{}\;\!M_\nu^{}(t,x)
$
\ for all}
$(t,x)\in\Omega.
$
\hfill (11.2)
\\[2ex]
\indent
{\sl Proof}.\! 
\vspace{1ex}
Using Theorem 11.1, from Property 1.11, we get the statement of Property 11.13.$\k$
%Using Theorem 11.1, from Property 1.11 it follows that the statement of Pro\-perty~11.13 is true. $\k$

{\bf Corollary 11.9.} 
\vspace{0.5ex}
{\it
Suppose 
$\mu_\nu^{}\in \text{\rm M}_{_{\Omega}},\
\gamma_\nu^{}\in\R\backslash\{0\},\ 
\mu_\nu^{{}^{\scriptsize \gamma_{\nu}^{}}}\in C^1\Omega, \
\nu=1,\ldots, s,\ s\leq m-2,$
${\rm g}_k^{{}^{\scriptsize \gamma_{k}^{}}}\!\in C^1\Omega,\
\gamma_k^{}\in\R\backslash\{0\},\
k=s+1,\ldots,m.$ 
Then
\\[1.5ex]
\mbox{}\hfill                      
$
\displaystyle
\prod\limits_{\nu=1}^{s} 
\mu_\nu^{{}^{\scriptsize \gamma_{\nu}^{}}}\, 
\prod\limits_{k=s+1}^{m} 
{\rm g}_k^{{}^{\scriptsize \gamma_{k}^{}}} 
\in \text{\rm I}_{_{\Omega}}
\iff
\biggl(
\,\prod\limits_{k=s+1}^{m} 
{\rm g}_k^{{}^{\scriptsize \gamma_{k}^{}}}, \
\sum\limits_{\nu=1}^{s}
\gamma_{\nu}\,{\rm div}\;\!{\frak d}\biggl) 
\in \text{\rm J}_{_{\Omega}}\;\!.
\hfill
$
\\[1.25ex]
Moreover, there exist functions
$M_k^{}\in C^1\Omega,\ k=s+1,\ldots, m,$
such that the identities {\rm (1.14)} are true and
\\[1.25ex]
\mbox{}\hfill                      % (11.3)
$
\displaystyle
\sum\limits_{k=s+1}^{m}
\gamma_{k}^{}\;\!M_k^{}(t,x)=
\sum\limits_{\nu=1}^{s}
\gamma_{\nu}^{}\,{\rm div}\;\!{\frak d}(t,x)
$
\ for all} $(t,x)\in\Omega.
$
\hfill (11.3)
\\[2ex]
\indent
{\bf Property 11.14.} 
\vspace{0.5ex}
{\it
Suppose 
$p_j^{}\in \text{\rm P}_{_{\!\Xi^{\;\!\prime}}},\ 
\gamma_j^{}\in\R\backslash\{0\},$ 
a set $\Omega_{_0}\subset\Xi^{\;\!\prime}\!$ such that  
$p_j^{{}^{\scriptsize \gamma_j^{}}}\in C^1\Omega_{_0},$ $j=1,\ldots, m.$
Then we claim that
\\[1.5ex]
\mbox{}\hfill  
$
\displaystyle
\prod\limits_{j=1}^{m} 
p_j^{{}^{\scriptsize \gamma_{j}^{}}}
\in \text{\rm I}_{_{\Omega_{{}_{\tiny\;\! 0}}}}
\iff
\bigl(p_j^{}, M_j^{}\bigr)\in \text{\rm A}_{_{\Xi^{\;\!\prime}}},
\ \ 
j=1,\ldots, m.
\hfill
$
\\[1.5ex]
Moreover, there exist the cofactors $M_1^{},\ldots, M_m^{}$
\vspace{0.5ex}
such that the identity {\rm (11.1)} holds.}

{\sl Proof}. 
\vspace{1ex}
Using Theorem\! 11.1, from Theorem\! 2.2, we get the statement of Pro\-perty\! 11.14.$\k$

{\bf  Property 11.15.}
{\it
If 
$
\bigl(\exp\omega_j^{}, M_j^{}\bigr)\in \text{\rm E}_{_{\Omega}},\ 
\gamma,\gamma_j^{}\in\R\backslash\{0\}, \ j=1,\ldots,m, \
{\rm g}^{\gamma}\in C^1\Omega,$ then}
\\[1.5ex]
\mbox{}\hfill
$
\displaystyle
{\rm g}^{\gamma}\exp\sum\limits_{j=1}^{m}\gamma_j^{}\;\!\omega_j^{}\in 
\text{\rm I}_{_{\Omega}}
\iff
\biggl(\,
{\rm g},{}- \dfrac{1}{\gamma}\,
\sum\limits_{j=1}^{m}\gamma_j^{}\;\!M_j^{}\biggr)\in 
\text{\rm J}_{_{\Omega}}.
\hfill
$
\\[1.5ex]
\indent
{\sl Proof}. 
\vspace{1ex}
Using Theorem\! 11.1, from Property\! 3.7, we get the statement of Pro\-perty\! 11.15.$\k$

{\bf  Corollary 11.10.}
{\it
If
$
\exp\omega_j^{}\in \text{\rm ME}_{_{\Omega}},\ 
\gamma,\gamma_j^{}\in\R\backslash\{0\}, \ j=1,\ldots,m, \
{\rm g}^{\gamma}\in C^1\Omega,$ then}
\\[1.5ex]
\mbox{}\hfill
$
\displaystyle
{\rm g}^{\gamma}\exp\sum\limits_{j=1}^{m}\gamma_j^{}\;\!\omega_j^{}\in 
\text{\rm I}_{_{\Omega}}
\iff
\biggl(\,
{\rm g},\ \dfrac{1}{\gamma}\,
\sum\limits_{j=1}^{m}\gamma_j^{}\,{\rm div}\;\!{\frak d}\biggr)\in 
\text{\rm J}_{_{\Omega}}.
\hfill
$
\\[2ex]
\indent
{\bf  Property\! 11.16.}\!
{\it
If
$\!\bigl(\exp\omega_\tau^{}, M_\tau^{}\bigr)\!\in\! \text{\rm E}_{_{\Omega}},\,
\tau\!=\!1,\ldots, m\!-\!1,\ 
\gamma_j^{}\in\R\backslash\{0\}, \, j\!=\!1,\ldots,m,$ then}
\\[1.25ex]
\mbox{}\hfill
$
\displaystyle
\sum\limits_{j=1}^{m}\gamma_j^{}\;\!\omega_j^{}\in 
\text{\rm I}_{_{\Omega}}
\iff
\biggl(\,
\exp\omega_{m}^{}\;\!,{}-
 \dfrac{1}{\gamma_{m}^{}}\,
\sum\limits_{\tau=1}^{m-1}\gamma_\tau^{}\;\!M_\tau^{}\biggr)\in 
\text{\rm E}_{_{\Omega}}.
\hfill
$
\\[1.25ex]
\indent
{\sl Proof}. 
\vspace{1ex}
Using Theorem\! 11.1, from Corollary\! 3.3, we get the statement of Pro\-perty\! 11.16.$\k$

{\bf  Corollary\! 11.11.}\!
{\it
If
$\exp\omega_\tau^{}\in \text{\rm ME}_{_{\Omega}},\
\tau\!=\!1,\ldots, m\!-\!1,\ 
\gamma_j^{}\in\R\backslash\{0\}, \ j=1,\ldots,m,$ then}
\\[1.5ex]
\mbox{}\hfill
$
\displaystyle
\sum\limits_{j=1}^{m}\gamma_j^{}\;\!\omega_j^{}\in 
\text{\rm I}_{_{\Omega}}
\iff
\biggl(\,
\exp\omega_{m}^{}\;\!,\
 \dfrac{1}{\gamma_{m}^{}}\,
\sum\limits_{\tau=1}^{m-1}\gamma_\tau^{}\,{\rm div}\;\!{\frak d}\biggr)\in 
\text{\rm E}_{_{\Omega}}.
\hfill
$
\\[2.25ex]
\indent
{\bf  Property 11.17.}
\vspace{0.75ex}
{\it
Suppose 
$
\bigl(\exp\omega_\nu^{}, M_\nu^{}\bigr)\in \text{\rm E}_{_{\Omega}},\ 
\nu=1,\ldots, s,\ s\leqslant m-2,\  
\gamma_j^{}\in\R\backslash\{0\}, 
\linebreak 
j=1,\ldots,m,\ 
{\rm g}_k^{{}^{\scriptsize \gamma_k^{}}}\!\in C^1\Omega,\ 
k=s+1,\ldots,m.$ Then we have
\\[2.25ex]
\mbox{}\hfill
$
\displaystyle
\prod\limits_{k=s+1}^{m}
{\rm g}_k^{{}^{\scriptsize \gamma_k^{}}}
\exp\sum\limits_{\nu=1}^{s}
\gamma_\nu^{}\;\!\omega_\nu^{}\in 
\text{\rm I}_{_{\Omega}}
\iff
\biggl(\,
\prod\limits_{k=s+1}^{m}
{\rm g}_k^{{}^{\scriptsize \gamma_k^{}}},{}-
\sum\limits_{\nu=1}^{s}\gamma_\nu^{}\;\!M_\nu^{}\biggr)\in 
\text{\rm J}_{_{\Omega}}.
\hfill
$
\\[2.25ex]
Moreover, there exist functions $M_k^{}\in C^1\Omega,\ k=s+1,\ldots, m,$
\vspace{0.5ex}
such that the identities {\rm (1.14)} and {\rm (11.2)} hold.}
\vspace{0.5ex}

{\sl Proof}. 
\vspace{1.5ex}
Using Theorem\! 11.1, from Property\! 3.8, we get the statement of Pro\-perty\! 11.17.$\k$

{\bf  Corollary 11.12.}
\vspace{0.5ex}
{\it
Suppose 
$
\exp\omega_\nu^{}\in \text{\rm ME}_{_{\Omega}},\ 
\nu=1,\ldots, s,\ s\leqslant m-2,\  
\gamma_j^{}\in\R\backslash\{0\}, 
\linebreak 
j=1,\ldots,m,\ 
{\rm g}_k^{{}^{\scriptsize \gamma_k^{}}}\!\in C^1\Omega,\ 
k=s+1,\ldots,m.$ Then we claim that
\\[2ex]
\mbox{}\hfill
$
\displaystyle
\prod\limits_{k=s+1}^{m}
{\rm g}_k^{{}^{\scriptsize \gamma_k^{}}}
\exp\sum\limits_{\nu=1}^{s}
\gamma_\nu^{}\;\!\omega_\nu^{}\in 
\text{\rm I}_{_{\Omega}}
\iff
\biggl(\,
\prod\limits_{k=s+1}^{m}
{\rm g}_k^{{}^{\scriptsize \gamma_k^{}}},\ 
\sum\limits_{\nu=1}^{s}\gamma_\nu^{}\,{\rm div}\;\!{\frak d}\biggr)\in 
\text{\rm J}_{_{\Omega}}.
\hfill
$
\\[1.75ex]
Moreover, there exist functions $M_k^{}\in C^1\Omega,\ k=s+1,\ldots, m,$
\vspace{0.5ex}
such that the identities {\rm (1.14)} and {\rm (11.3)} hold.
}
\vspace{1.5ex}

{\bf  Property 11.18.}
\vspace{1ex}
{\it
Let 
$
\bigl(\exp\omega_\nu^{}, M_\nu^{}\bigr)\in \text{\rm E}_{_{\Omega}},\ 
\nu=1,\ldots, s,\ s\leqslant m-1,\  
({\rm g}_k^{}, M_k^{})\in \text{\rm J}_{_{\Omega}},
\linebreak 
{\rm g}_k^{{}^{\scriptsize \gamma_k^{}}}\!\in C^1\Omega,
\ 
k=s+1,\ldots,m,\ 
\gamma_j^{}\in\R\backslash\{0\}, 
\ j=1,\ldots,m. 
$ 
Then 
\vspace{0.5ex}
$
\prod\limits_{k=s+1}^{m}
{\rm g}_k^{{}^{\scriptsize \gamma_k^{}}}
\exp\sum\limits_{\nu=1}^{s}
\gamma_\nu^{}\;\!\omega_\nu^{}\in 
\text{\rm I}_{_{\Omega}}
$
if and only if the identity {\rm (11.1)} holds.
}
\vspace{0.75ex}

{\sl Proof}. 
\vspace{1.5ex}
Using Theorem\! 11.1, from Property\! 3.9, we get the statement of Pro\-perty\! 11.18.$\k$

{\bf  Corollary 11.13.}
\vspace{1ex}
{\it
Suppose 
$
\exp\omega_\nu^{}\in \text{\rm ME}_{_{\Omega}},\ 
\nu=1,\ldots, s,\ s\leqslant m-1,\  
({\rm g}_k^{}, M_k^{})\in \text{\rm J}_{_{\Omega}},
\linebreak 
{\rm g}_k^{{}^{\scriptsize \gamma_k^{}}}\!\in C^1\Omega,
\ 
k=s+1,\ldots,m,\ 
\gamma_j^{}\in\R\backslash\{0\}, 
\ j=1,\ldots,m. 
$ 
Then 
\vspace{0.75ex}
$
\prod\limits_{k=s+1}^{m}
{\rm g}_k^{{}^{\scriptsize \gamma_k^{}}}
\exp\sum\limits_{\nu=1}^{s}
\gamma_\nu^{}\;\!\omega_\nu^{}\in 
\text{\rm I}_{_{\Omega}}
$
if and only if the identity {\rm (11.3)} holds on the domain $\Xi^{\;\!\prime}.$
}
\vspace{1.5ex}

{\bf  Corollary 11.14.}
\vspace{0.75ex}
{\it
Suppose 
$
\bigl(\exp\omega_\nu^{}, M_\nu^{}\bigr)\in \text{\rm E}_{_{\Omega}},\ 
\nu=1,\ldots, s,\ s\leqslant m-1,\  
\mu_k^{}\in \text{\rm M}_{_{\Omega}},
\linebreak 
\mu_k^{{}^{\scriptsize \gamma_k^{}}}\!\in C^1\Omega,
\ 
k=s+1,\ldots,m,\ 
\gamma_j^{}\in\R\backslash\{0\}, 
\ j=1,\ldots,m. 
$ 
Then 
\\[1.75ex]
\mbox{}\hfill
$
\displaystyle
\prod\limits_{k=s+1}^{m}
\mu_k^{{}^{\scriptsize \gamma_k^{}}}
\exp\sum\limits_{\nu=1}^{s}
\gamma_\nu^{}\;\!\omega_\nu^{}\in 
\text{\rm I}_{_{\Omega}}
\hfill
$
\\[1.75ex]
if and only if the linear combination of cofactors
\\[1.75ex]
\mbox{}\hfill                     
$
\displaystyle
\sum\limits_{\nu=1}^{s}\gamma_{\nu}^{}\;\!M_\nu^{}(t,x)=
\sum\limits_{k=s+1}^{m}
\gamma_{k}^{}\,{\rm div}\;\!{\frak d}(t,x)
$
\ for all} 
$(t,x)\in\Xi^{\;\!\prime}.
\hfill
$
\\[2.5ex]
\indent
{\bf  Corollary 11.15.}
\vspace{1ex}
{\it
Suppose 
$
\exp\omega_\nu^{}\in \text{\rm ME}_{_{\Omega}},\ 
\nu=1,\ldots, s,\ s\leqslant m-1,\  
\mu_k^{}\in \text{\rm M}_{_{\Omega}},
\linebreak 
\mu_k^{{}^{\scriptsize \gamma_k^{}}}\!\in C^1\Omega,
\ 
k=s+1,\ldots,m,\ 
\gamma_j^{}\in\R\backslash\{0\}, 
\ j=1,\ldots,m. 
$ 
Then we have}
\\[1.5ex]
\mbox{}\hfill                     
$
\displaystyle
\prod\limits_{k=s+1}^{m}
\mu_k^{{}^{\scriptsize \gamma_k^{}}}
\exp\sum\limits_{\nu=1}^{s}
\gamma_\nu^{}\;\!\omega_\nu^{}\in 
\text{\rm I}_{_{\Omega}}
\iff
\sum\limits_{j=1}^{m}\gamma_{j}^{}= 0.
\hfill
$
\\[-2ex]

\newpage

{\bf Property\! 11.19.}\!
\vspace{0.75ex}
{\it
Let 
$\!\bigl((p_j^{},\! M_j^{}), (h_j^{},q_j^{},\rho_{j}^{}N_0^{})\bigr)\!\!\in\! 
\text{\rm B}_{_{\Xi^{\;\!\prime}}},\, 
\rho_j^{}\!\in\!\R\backslash\{0\},\, 
\lambda_j^{},\!\gamma_j^{}\!\in\!\R, j\!=\!1,\ldots, m,\!$ 
$\sum\limits_{j=1}^{m}|\lambda_j^{}|\ne 0,\ 
\varphi\in C^1T^{\;\!\prime},$ 
the set
\vspace{0.75ex}
$\Omega_{_0}\subset\Xi^{\;\!\prime}$ such that   
$\prod\limits_{j=1}^{m}p_j^{}(t,x)\ne 0$ for all $(t,x)\in\Omega_{_0}$ and
$\prod\limits_{j=1}^{m}p_j^{}(t,x)= 0$
for all $(t,x)\in {\sf C}_{_{\Xi^{\;\!\prime}}}\Omega_{_0},\ 
p_j^{\gamma_j^{}}\in C^1\Omega_{_0}, \ j=1,\ldots, m.$
Then we claim that
\\[1.25ex]
\mbox{}\hfill                           
$
\displaystyle
\prod\limits_{j=1}^{m}p_j^{\gamma_j^{}}
\sum\limits_{j=1}^{m}\lambda_j^{}
\exp\biggl(\;\!\dfrac{q_j^{}}{\rho_j^{}\;\!p_j^{\;\!h_j^{}}}+\varphi\biggr)\in 
\text{\rm I}_{_{\Omega_{{}_{\tiny\;\! 0}}}}
\hfill
$
\\[1.5ex]
if and only if the following identity holds
\\[1.5ex]
\mbox{}\hfill                           
$
\displaystyle
{\sf D}\;\!\varphi(t)+N_0^{}(t,x)+
\sum\limits_{j=1}^{m}\gamma_j^{}\;\!M_j^{}(t,x)=0
$
\ for all} 
$(t,x)\in \Xi^{\;\!\prime}.
\hfill
$
\\[1.5ex]
\indent
{\sl Proof}. 
\vspace{1.25ex}
Using Theorem\! 11.1, from Property\! 5.5, we get the statement of Pro\-perty\! 11.19.$\k$

{\bf Corollary 11.16.}\!
\vspace{0.75ex}
{\it
Let 
$\!\bigl(\nu_j^{}, (h_j^{},q_j^{},\rho_{j}^{}N_0^{})\!\bigr)\!\in\! 
\text{\rm MB}_{_{\Xi^{\;\!\prime}}},\ 
\rho_j^{}\!\in\!\R\backslash\{0\},\ 
\lambda_j^{},\gamma_j^{}\!\in\!\R,\ j\!=\!1,\ldots, m,\!$ 
$\sum\limits_{j=1}^{m}|\lambda_j^{}|\ne 0,\ 
\varphi\in C^1T^{\;\!\prime},$ 
the set 
\vspace{0.75ex}
$\Omega_{_0}\subset\Xi^{\;\!\prime}$ such that  
$\prod\limits_{j=1}^{m}\nu_j^{}(t,x)\ne 0$ for all $(t,x)\in\Omega_{_0}$ and
$\prod\limits_{j=1}^{m}\nu_j^{}(t,x)= 0$
for all $(t,x)\in {\sf C}_{_{\Xi^{\;\!\prime}}}\Omega_{_0},\ 
\nu_j^{\gamma_j^{}}\in C^1\Omega_{_0}, \ j=1,\ldots, m.$
Then we have
\\[1.25ex]
\mbox{}\hfill                           
$
\displaystyle
\prod\limits_{j=1}^{m}\nu_j^{\gamma_j^{}}
\sum\limits_{j=1}^{m}\lambda_j^{}
\exp\biggl(\;\!\dfrac{q_j^{}}{\rho_j^{}\;\!\nu_j^{\;\!h_j^{}}}+\varphi\biggr)\in 
\text{\rm I}_{_{\Omega_{{}_{\tiny\;\! 0}}}}
\hfill
$
\\[1.5ex]
if and only if the following identity holds
\\[1.5ex]
\mbox{}\hfill                           
$
\displaystyle
{\sf D}\;\!\varphi(t)+N_0^{}(t,x)-
\sum\limits_{j=1}^{m}\gamma_j^{}\,{\rm div}\;\!{\frak d}(t,x)=0
$
\ for all} 
$(t,x)\in \Xi^{\;\!\prime}.
\hfill
$
\\[1.75ex]
\indent
{\bf Property 11.20.}\!
\vspace{0.75ex}
{\it
Let 
$\bigl((p_j^{}, M_j^{}), (h_j^{},q_j^{},N_j^{})\bigr)\!\in\! \text{\rm B}_{_{\Xi^{\;\!\prime}}},\ 
\gamma_j^{},\xi_j^{}\!\in\!\R,\ \varphi_j^{}\!\in\! C^1T^{\;\!\prime},\ 
j\!=\!1,\ldots, m,$ 
the set 
\vspace{0.75ex}
$\!\Omega_{_0}\!\subset\!\Xi^{\;\!\prime}\!$ such that   
$\prod\limits_{j=1}^{m}p_j^{}(t,x)\!\ne\! 0$ for all $(t,x)\!\in\!\Omega_{_0}$ and
$\prod\limits_{j=1}^{m}p_j^{}(t,x)= 0$
for all $(t,x)\!\in\! {\sf C}_{_{\Xi^{\;\!\prime}}}\Omega_{_0},\!$ 
$p_j^{\gamma_j^{}}\in C^1\Omega_{_0}, \ j=1,\ldots, m.$
Then we claim that
\\[1.5ex]
\mbox{}\hfill                           
$
\displaystyle
\prod\limits_{j=1}^{m}p_j^{\gamma_j^{}}
\exp\sum\limits_{j=1}^{m}\xi_j^{}
\biggl(\;\!\dfrac{q_j^{}}{p_j^{\;\! h_j^{}}}+\varphi_j^{}\biggr)\in 
\text{\rm I}_{_{\Omega_{{}_{\tiny\;\! 0}}}}
\hfill
$
\\[1.75ex]
if and only if the following identity holds
\\[1.5ex]
\mbox{}\hfill                           
$
\displaystyle
\sum\limits_{j=1}^{m}\bigl(
\gamma_j^{}\;\!M_j^{}(t,x)+
\xi_j^{}\;\!\bigl(N_j^{}(t,x)+{\sf D}\;\!\varphi_j^{}(t)\bigr)\bigr)=0
$
\ for all} 
$(t,x)\in \Xi^{\;\!\prime}.
\hfill
$
\\[1.5ex]
\indent
{\sl Proof}. 
\vspace{1.5ex}
Using Theorem\! 11.1, from Property\! 5.7, we get the statement of Pro\-perty\! 11.20.$\k$

{\bf Corollary 11.17.}\!
\vspace{0.75ex}
{\it
Let 
$\!\bigl(\nu_j^{}, (h_j^{},q_j^{},N_j^{})\bigr)\!\in \text{\rm MB}_{_{\Xi^{\;\!\prime}}},\, 
\gamma_j^{},\xi_j^{}\in\R,\ \varphi_j^{}\in C^1T^{\;\!\prime},\, 
j=1,\ldots, m,\!$ 
the set 
\vspace{0.75ex}
$\!\Omega_{_0}\!\subset\!\Xi^{\;\!\prime}\!$ such that   
$\prod\limits_{j=1}^{m}\nu_j^{}(t,x)\!\ne\! 0$ for all $(t,x)\!\in\!\Omega_{_0}$ and
$\prod\limits_{j=1}^{m}\nu_j^{}(t,x)= 0$ 
for all $(t,x)\!\in\! {\sf C}_{_{\Xi^{\;\!\prime}}}\Omega_{_0},$ 
$\nu_j^{\gamma_j^{}}\in C^1\Omega_{_0}, \ j=1,\ldots, m.$
Then we have
\\[1.5ex]
\mbox{}\hfill                           
$
\displaystyle
\prod\limits_{j=1}^{m}\nu_j^{\gamma_j^{}}
\exp\sum\limits_{j=1}^{m}\xi_j^{}
\biggl(\;\!\dfrac{q_j^{}}{\nu_j^{\;\! h_j^{}}}+\varphi_j^{}\biggr)\in 
\text{\rm I}_{_{\Omega_{{}_{\tiny\;\! 0}}}}
\hfill
$
\\[1.5ex]
if and only if the following identity holds
\\[1.75ex]
\mbox{}\hfill                           
$
\displaystyle
\sum\limits_{j=1}^{m}
\xi_j^{}\;\!\bigl(N_j^{}(t,x)+{\sf D}\;\!\varphi_j^{}(t)\bigr)=
\sum\limits_{j=1}^{m}
\gamma_j^{}\,{\rm div}\;\!{\frak d}(t,x)
$
\ for all} 
$(t,x)\in \Xi^{\;\!\prime}.
\hfill
$
\\[2ex]
\indent
{\bf Property 11.21.}
\vspace{1ex}
{\it
Suppose 
$\Bigl((p, M), 
\Bigl(h_\xi^{},q_{_{\scriptstyle h_\xi^{}  f_\xi^{}}}, N_{h_\xi^{}  f_\xi^{}}^{}\Bigr)\Bigr)
\in \text{\rm B}_{_{\Xi^{\;\!\prime}}},\
\gamma_{_{\scriptstyle h_\xi^{}  f_\xi^{}}}\in\R,\ 
\varphi_{_{\scriptstyle h_\xi^{}  f_\xi^{}}}\in C^1T^{\;\!\prime},$ 
$f_{\xi}^{}=1,\ldots, \delta_{\xi}^{},\ \xi=1,\ldots,\varepsilon,\ \gamma\in\R\backslash\{0\},$
\vspace{1ex}
the set
$\Omega_{_0}\subset\Xi^{\;\!\prime}$ such that   
$p(t,x)\ne 0$ for all $(t,x)\in\Omega_{_0}$ and
$p(t,x)= 0$ for all $(t,x)\in {\sf C}_{_{\Xi^{\;\!\prime}}}\Omega_{_0},\ 
p^{\gamma}\in C^1\Omega_{_0}.$
Then we claim that
\\[1.5ex]
\mbox{}\hfill                           
$
\displaystyle
p^{\gamma}
\exp\sum\limits_{\xi=1}^{\varepsilon}\sum\limits_{f_{\xi}^{}=1}^{\delta_{\xi}^{}}
\biggl(\;\!\gamma_{_{\scriptstyle h_\xi^{}  f_\xi^{}}}\biggl(\,
\dfrac{q_{_{\scriptstyle h_\xi^{}  f_\xi^{}} }}{\displaystyle  p^{\;\!h_\xi^{}}}+
\varphi_{_{\scriptstyle h_\xi^{}  f_\xi^{}}}\biggr)\biggr)\in 
\text{\rm I}_{_{\Omega_{{}_{\tiny\;\! 0}}}}
\hfill
$
\\[1.5ex]
if and only if the following identity holds
\\[1.5ex]
\mbox{}\hfill                           
$
\displaystyle
\gamma\;\!M(t,x)+
\sum\limits_{\xi=1}^{\varepsilon}\sum\limits_{f_{\xi}^{}=1}^{\delta_{\xi}^{}}
\biggl(\gamma_{_{\scriptstyle h_\xi^{}  f_\xi^{}}}\Bigl(
N_{h_\xi^{}  f_\xi^{}}^{} (t,x)+{\sf D}\;\!\varphi_{_{\scriptstyle h_\xi^{}  f_\xi^{}}}(t)\Bigr)\biggr)=0
$
\ for all}
$(t,x)\in \Xi^{\;\!\prime}.
\hfill
$
\\[1.5ex]
\indent
{\sl Proof}. 
\vspace{1.25ex}
Using Theorem\! 11.1, from Property\! 5.8, we get the statement of Pro\-perty\! 11.21.$\k$

{\bf Corollary 11.18.}
\vspace{1ex}
{\it
Suppose 
$\Bigl(\nu, 
\Bigl(h_\xi^{},q_{_{\scriptstyle h_\xi^{}  f_\xi^{}}}, N_{h_\xi^{}  f_\xi^{}}^{}\Bigr)\Bigr)
\in \text{\rm MB}_{_{\Xi^{\;\!\prime}}},\
\gamma_{_{\scriptstyle h_\xi^{}  f_\xi^{}}}\in\R,\ 
\varphi_{_{\scriptstyle h_\xi^{}  f_\xi^{}}}\in C^1T^{\;\!\prime},$ 
$f_{\xi}^{}=1,\ldots, \delta_{\xi}^{},\ \xi=1,\ldots,\varepsilon,\ \gamma\in\R\backslash\{0\},$
\vspace{1ex}
the set
$\Omega_{_0}\subset\Xi^{\;\!\prime}$ such that  
$\nu(t,x)\ne 0$ for all $(t,x)\in\Omega_{_0}$ and
$\nu(t,x)= 0$ for all $(t,x)\in {\sf C}_{_{\Xi^{\;\!\prime}}}\Omega_{_0},\ 
\nu^{\gamma}\in C^1\Omega_{_0}.$
Then we have
\\[1.5ex]
\mbox{}\hfill                           
$
\displaystyle
\nu^{\gamma}
\exp\sum\limits_{\xi=1}^{\varepsilon}\sum\limits_{f_{\xi}^{}=1}^{\delta_{\xi}^{}}
\biggl(\;\!\gamma_{_{\scriptstyle h_\xi^{}  f_\xi^{}}}\biggl(\,
\dfrac{q_{_{\scriptstyle h_\xi^{}  f_\xi^{}} }}{\displaystyle  \nu^{\;\!h_\xi^{}}}+
\varphi_{_{\scriptstyle h_\xi^{}  f_\xi^{}}}\biggr)\biggr)\in 
\text{\rm I}_{_{\Omega_{{}_{\tiny\;\! 0}}}}
\hfill
$
\\[1.5ex]
if and only if the following identity holds
\\[1.25ex]
\mbox{}\hfill                           
$
\displaystyle
{\rm div}\;\!{\frak d}(t,x)=
\dfrac{1}{\gamma}\
\sum\limits_{\xi=1}^{\varepsilon}\sum\limits_{f_{\xi}^{}=1}^{\delta_{\xi}^{}}
\biggl(\gamma_{_{\scriptstyle h_\xi^{}  f_\xi^{}}}\Bigl(
N_{h_\xi^{}  f_\xi^{}}^{} (t,x)+{\sf D}\;\!\varphi_{_{\scriptstyle h_\xi^{}  f_\xi^{}}}(t)\Bigr)\biggr)
$
\ for all} 
$(t,x)\in \Xi^{\;\!\prime}.
\hfill
$
\\[2ex]
\indent
{\bf Property 11.22.}
\vspace{0.75ex}
{\it
Suppose numbers $\gamma_1^{}, \gamma_2^{}\in\R,$
functions 
$u,v\in \text{\rm P}_{_{\!\Xi^{\;\!\prime}}}$ are relatively prime, 
$(u+i\;\!v,\;\! U+i\;\!V)\in \text{\rm H}_{_{\Xi^{\;\!\prime}}},$
\vspace{0.75ex}
the set $\Omega_{_0}\!\subset\! \Xi^{\;\!\prime}\!$ such that   
$u(t,x)\ne 0$ for all $(t,x)\!\in\! \Omega_{_0}$ and
$u(t,x)=0$ for all $(t,x)\in {\sf C}_{_{\Xi^{\;\!\prime}}}\Omega_{_0}.$
Then
$(u^2+v^2)^{{}^{\scriptsize \gamma_1^{}}}
\exp\Bigl(\gamma_2^{}\;\!\arctan\dfrac{v}{u}\Bigr)\in 
\text{\rm I}_{_{\Omega_{{}_{\tiny\;\! 0}}}}$
if and only if 
\\[2ex]
\mbox{}\hfill
$
2\gamma_1^{}\;\!U(t,x)+\gamma_2^{}\;\!V(t,x)=0
$
\ for all} 
$(t,x)\in \Xi^{\;\!\prime}.
\hfill
$
\\[1.5ex]
\indent
{\sl Proof}. 
\vspace{1.25ex}
Using Theorem\! 11.1, from Property\! 6.15, we get the statement of Pro\-perty\! 11.22.$\k$

{\bf Property 11.23.}
\vspace{0.75ex}
{\it
Suppose 
$\bigl((u+i\;\!v,\;\! U+i\;\!V), (h,z,Q)\bigr)\in \text{\rm G}_{_{\Xi^{\;\!\prime}}},$
functions
$u, v\in \text{\rm P}_{_{\!\Xi}}$ are relatively prime, 
\vspace{0.75ex}
a function $z\in \text{\rm Z}_{_{\Xi}}$ is relatively prime with the function 
$u+i\;\!v,$ numbers $\gamma_1^{},\gamma_2^{},\gamma_3^{},\gamma_4^{}\in\R,$  the set
\vspace{1ex}
$\Omega_{_0}\subset\Xi^{\;\!\prime}$ such that   
$u(t,x)\ne 0$ for all $(t,x)\in\Omega_{_0}$ and
$u(t,x)= 0$ for all $(t,x)\in {\sf C}_{_{\Xi^{\;\!\prime}}}\Omega_{_0}.\!$
Then we claim that
\\[2ex]
\mbox{}\hfill                           
$
\displaystyle
\bigl(u^2+v^2\bigr)^{\!{}^{\scriptstyle\gamma_1^{}}}
\exp\biggl(
\dfrac{\gamma_2^{}\;\!{\rm Re}\;\!\bigl(z(u-i\;\!v)^h\bigr)+
\gamma_3^{}\;\!{\rm Im}\;\!\bigl(z(u-i\;\!v)^h\bigr)}{\bigl(u^2+v^2\bigr)^{h}}\,+\,
\gamma_4^{}\;\!\arctan\dfrac{v}{u}
\biggr)\in \text{\rm I}_{_{\Omega_{{}_{\tiny\;\! 0}}}}
\hfill                           
$
\\[1.75ex]
if and only if
\\[1.75ex]
\mbox{}\hfill
$
2\gamma_1^{}\;\!U(t,x)+
\gamma_2^{}\;\!{\rm Re}\;\!Q(t,x)+
\gamma_3^{}\;\!{\rm Im}\;\!Q(t,x)+
\gamma_4^{}\;\!V(t,x)=0
$
\ for all} 
$(t,x)\in \Xi^{\;\!\prime}.
\hfill
$
\\[1.5ex]
\indent
{\sl Proof}. 
\vspace{0.25ex}
Taking into account Property 1.9 and Theorem 11.1, 
from Theorem 7.3 it follows that the statement of Pro\-perty 11.23 is true. $\k$

\newpage

\mbox{}
\\[-0.5ex]
\centerline{
{\bf  12. Applications
}
}
\\[1.5ex]
\indent
{\bf 12.1.}
{\sl Darboux's problem}.
Consider a differential equation of first order
\\[1.5ex]
\mbox{}\hfill                     % (12.1)
$
Y(x,y)\;\!dx-X(x,y)\;\!dy=0,
$
\hfill (12.1)
\\[1.75ex]
where the functions $X\colon\R^2\to\R$ and $Y\colon\R^2\to\R$ 
\vspace{0.35ex}
are polynomials in the variables $x,\, y$  such that 
$
\max\bigl\{\deg X,\, \deg Y\bigr\}=d\geq 1.
$ 
\vspace{0.5ex}

The Darboux problem [8, 13] consists in constructing general integral 
of the differential equation (12.1) by known partial integrals.

The equation (12.1) is the equation of trajectories for 
the autonomous differential system of second order
\\[1.5ex]
\mbox{}\hfill                     % (12.2)
$
\dfrac{dx}{dt}=X(x,y),
\qquad
\dfrac{dy}{dt}=Y(x,y).
$
\hfill (12.2)
\\[1.75ex]
\indent
The general integral of equation (12.1) is the autonomous first integral of system (12.2) and vice versa.

Case $d=2.$
\vspace{0.35ex}
Let the function 
${\rm g}_j^{}\colon (x,y)\to {\rm g}_j^{}(x,y)$ for all $(x,y)\in D$
be a partial integral with cofactor $M_j^{}$ on the domain $D\subset\R^2$ 
\vspace{0.35ex}
of equation (12.1) (equivalently, of system (12.2)), $j=1,2,3.$
By definition 1.1,
$\deg M_j^{}\leq 1,\ j=1,2,3,$ i.e.
\\[1.5ex]
\mbox{}\hfill                    
$
M_j^{}(x,y)=
\alpha_{1j}^{}+\alpha_{2j}^{}\;\!x+\alpha_{3j}^{}\;\!y
$
\ for all 
$(x,y)\in\R^2,
\ \ 
j=1,2,3,
\hfill
$
\\[1.75ex]
where $\alpha_{1j}^{}\;\!,\;\!\alpha_{2j}^{}\;\!,\;\!\alpha_{3j}^{}\in\R,\ 
\sum\limits_{\zeta=1}^{3}|\alpha_{\zeta j}^{}|\ne 0,\ j=1,2,3.$
\vspace{0.5ex}

The set $\!D_{_{\!0}}^{}\!\subset\! D\!$ such that  
$\!{\rm g}_j^{{}^{\!\scriptsize\gamma_j^{}}}\!\in\! C^{1}D_{_0}^{}, 
\gamma_j^{}\!\in\!\R,\;\! j\!=\!1,2,3.\!$
By Property 11.11, the function 
\\[1.5ex]
\mbox{}\hfill                     % (12.3)
$
\displaystyle
F\colon (x,y)\to\ 
\prod\limits_{j=1}^3
{\rm g}_j^{{}^{\!\scriptsize\gamma_j^{}}}(x,y)
$
\ for all 
$
(x,y)\in D_{_{\!0}}^{}
$
\hfill (12.3)
\\[1.75ex]
is the general integral on the domain $D_{_{\!0}}^{}$ of equation (12.1) if and only if
\\[1.5ex]
\mbox{}\hfill                    
$
\displaystyle
\sum\limits_{j=1}^3
\gamma_j^{}\;\!
\bigl(\alpha_{1j}^{}+\alpha_{2j}^{}\;\!x+\alpha_{3j}^{}\;\!y\bigr)=0
$
\ for all 
$(x,y)\in\R^2.
\hfill
$
\\[1.75ex]
The latter is possible when the determinant
\\[2ex]
\mbox{}\hfill                    
$
\triangle=\left|
\begin{array}{ccc}
\alpha_{11}^{} & \alpha_{12}^{} & \alpha_{13}^{} 
\\[0.75ex]
\alpha_{21}^{} & \alpha_{22}^{} & \alpha_{23}^{} 
\\[0.75ex]
\alpha_{31}^{} & \alpha_{32}^{} & \alpha_{33}^{} 
\end{array}
\right|=0.
\hfill
$
\\[1.75ex]
\indent
Thus Darboux's problem is solved.

In the paper [7], we proposed to generalize the formulation of the Darboux problem.
In our case, along with the construction of general integral 
we are finding the integrating factor of equation (12.1).

By Property 8.8, the function (12.3) is the  integrating factor of equation (12.1)
if and only if  the linear combination of cofactors
\\[1.5ex]
\mbox{}\hfill              % (12.4)               
$
\displaystyle
\sum\limits_{j=1}^3
\gamma_j^{}\;\!M_j^{}(x,y)=
{}-\partial_x^{}X(x,y)-\partial_y^{}Y(x,y)
$
\ for all $(x,y)\in\R^2.
$
\hfill (12.4)
\\[1.75ex]
\indent
The divergence of vector field defined by the equation (12.1) with $d=2$ is
\\[1.5ex]
\mbox{}\hfill                    
$
\partial_x^{}X(x,y)+\partial_y^{}Y(x,y)=
\beta_{1}^{}+\beta_{2}^{}\;\!x+\beta_{3}^{}\;\!y
$
\ for all 
$(x,y)\in\R^2,
\hfill
$
\\[1.75ex]
where $\beta_{1}^{}\;\!,\;\!\beta_{2}^{}\;\!,\;\!\beta_{3}^{}$  are real numbers.
\vspace{0.5ex}

The identity (12.4) holds if and only if the system of equations
\\[1.5ex]
\mbox{}\hfill              % (12.5)               
$
\begin{array}{c}
\alpha_{11}^{}\;\!\gamma_1^{}+
\alpha_{12}^{}\;\!\gamma_2^{}+
\alpha_{13}^{}\;\!\gamma_3^{}={}-\beta_1^{},
\\[1.25ex]
\alpha_{21}^{}\;\!\gamma_1^{}+
\alpha_{22}^{}\;\!\gamma_2^{}+
\alpha_{23}^{}\;\!\gamma_3^{}={}-\beta_2^{},
\\[1.25ex]
\alpha_{31}^{}\;\!\gamma_1^{}+
\alpha_{32}^{}\;\!\gamma_2^{}+
\alpha_{33}^{}\;\!\gamma_3^{}={}-\beta_3^{}
\end{array}
$
\hfill (12.5)
\\[2ex]
is compatible.

If $\triangle\ne 0,\ |\beta_1^{}|+|\beta_2^{}|+|\beta_3^{}|\ne 0,$
\vspace{0.35ex}
then 
from the linear non-homogeneous system (12.5), we get
$\gamma_{1}^{}, \gamma_{2}^{}, \gamma_{3}^{}.$
\vspace{0.5ex}
So we can build the integrating factor (12.3) on the set
$\!D_{_{\!0}}^{}\!$ for equation (12.1).

If $\beta_1^{}=\beta_2^{}=\beta_3^{}= 0,$ then
\vspace{0.35ex}
the equation (12.1) is an exact differential equation
(the integrating factor for this equation is arbitrary real constant).

Thus, if we have three partial integrals of equation (12.1) with $d=2,$
then the generalized Darboux problem  is solved in closed form
\vspace{0.5ex}
(we can build general integral or integrating factor).

General case $d\geq 2.$
In the same way as above, we obtain:
the number of partial integrals that we need to build general integral or 
integrating factor (to solve the generalized Darboux problem) of equation (12.1) 
is equal to the number of monomials in
a degree $d-1$ polynomial of two variables $x$ and $y$ 
(the highest degree of cofactors of partial integrals).
\vspace{0.35ex}

Since two-variable polynomial of degree $d-1$ has 
$\dfrac{d\;\!(d+1)}{2}$ monomials, we have
\vspace{0.5ex}

{\bf Theorem 12.1.}
\vspace{0.35ex}
{\it 
Suppose we know $m=\dfrac{d\;\!(d+1)}{2}$
partial integrals ${\rm g}_j^{}$ on the set $D$ of equation {\rm (12.1)},
$\gamma_j^{}\in\R, $
\vspace{0.5ex}
and the set $D_{_{\!0}}^{}\subset D$ such that 
${\rm g}_j^{{}^{\!\scriptsize\gamma_j^{}}}\in C^{1}D_{_0}^{},\ 
j=1,\ldots, m.$
Then the function
\\[1ex]
\mbox{}\hfill                   
$
\displaystyle
F\colon (x,y)\to\ 
\prod\limits_{j=1}^m
{\rm g}_j^{{}^{\!\scriptsize\gamma_j^{}}}(x,y)
$
\ for all $(x,y)\in D_{_{\!0}}^{}
\hfill
$
\\[1.5ex]
is a general integral or an integrating factor on the set 
$D_{_{\!0}}^{}$ of equation} (12.1).
\vspace{1ex}

{\bf 12.2.}
The differential equation of the first order
[14, pp. 136 -- 139; 23]
\\[1.5ex]
\mbox{}\hfill                   
$
\bigl(x+a^2x^2+2\;\!axy-(1+2\;\!a^2)y^2\bigr)\;\!dx+
\bigl(y+ax^2+(3+a^2)xy-ay^2\bigr)\;\!dy=0
\quad
(a\in\R\backslash\{0\})
\hfill
$
\\[1.5ex]
has two polynomial partial integrals:
\\[1.5ex]
\mbox{}\hfill                   
$
p_1^{}\colon (x,y)\to\ 
1+(1+a^2)\bigl(3\;\!x+3\;\!ax(y+ax)+a(y+ax)^3\bigr)
$
\ for all $(x,y)\in\R^2
\hfill
$
\\[1.5ex]
with cofactor  
\\[1ex]
\mbox{}\hfill                   
$
M_1^{}\colon (x,y)\to\ 
{}-3\;\!(1+a^2)\;\!y
$
\ for all $(x,y)\in\R^2
\hfill
$
\\[1ex]
and 
\\[1ex]
\mbox{}\hfill                   
$
p_2^{}\colon (x,y)\to\ 
1+(1+a^2)\bigl(2\;\!x+(y+ax)^2\bigr)
$
\ for all $(x,y)\in\R^2
\hfill
$
\\[1.5ex]
with cofactor 
\\[1ex]
\mbox{}\hfill                   
$
M_2^{}\colon (x,y)\to\ 
{}-2\;\!(1+a^2)\;\!y
$
\ for all $(x,y)\in\R^2.
\hfill
$
\\[1.5ex]
\indent
Using Corollary 11.5 under $M_0^{}={}-(1+a^2)\;\!y,\ \rho_1^{}=3,\ \rho_2^{}=2,$
\vspace{0.5ex}
from the equality $3\gamma_1^{}+2\gamma_2^{}=0,$ for example,
we get $\gamma_1^{}=2,\ \gamma_2^{}={}-3.$
Therefore the rational function
\\[1.75ex]
\mbox{}\hfill                   
$
F\colon (x,y)\to\ 
\dfrac{p_1^{\,2}(x,y)}{p_2^{\,3}(x,y)}
$
\ for all 
$(x,y)\in D_{_{\!0}}^{}
\hfill
$
\\[1.75ex]
is a general integral on the set $D_{_{\!0}}^{}\subset\R^2,$
\vspace{0.35ex}
where the set $D_{_{\!0}}^{}$ such that
$p_2^{}(x,y)\ne 0$ for all $(x,y)\in D_{_{\!0}}^{}$ and
$p_2^{}(x,y)= 0$ for all $(x,y)\in {\sf C}_{\R^2}^{}D_{_{\!0}}^{}.$
\vspace{0.5ex}

This example is one of those cases when a general integral of equations (12.1) under $d = 2$ 
is constructed not by three partial integrals but 
\vspace{1ex}
by a smaller number (two) of partial integrals.

{\bf 12.3.}
The system [28, p. 46]
\\[1.5ex]
\mbox{}\hfill                   
$
\dfrac{dx}{dt}={}-y+\dfrac{x^2}{2}-\dfrac{y^2}{2}\,,
\qquad
\dfrac{dy}{dt}=x\;\!(1+y)
\hfill
$
\\[2ex]
has 
$(y+1, x)\in\text{A}_{\R^2}^{}$ (Property 2.3) and 
$(x^2+y^2, x)\in\text{A}_{\R^2}^{}$ (Theorem 2.1).
\vspace{0.5ex}

By Corollary 11.7, 
$\dfrac{y+1}{x^2+y^2}\in \text{I}_{\R^2\backslash\{(0,0)\}}^{}.$
\vspace{1.5ex}

{\bf 12.4.}
The system [23]
\\[1.75ex]
\mbox{}\hfill                   % (12.6)
$
\dfrac{dx}{dt}={}-2+y+x^2+xy,
\qquad
\dfrac{dy}{dt}=4+2\;\!x+xy+y^2
$
\hfill (12.6)
\\[2ex]
has 
\vspace{0.75ex}
$\bigl((2+2\;\!x+y,\;\! x+y),\;\! (1, x+y, 1)\bigr)\in\text{B}_{\R^2}^{}.$ 
By Theorem 5.3, 
$
(2+2\;\!x+y,\;\! x+y)\in\text{A}_{\R^2}^{},$
$\Bigl(\exp\dfrac{x+y}{2+2\;\!x+y}\,,\;\! 1\Bigr)\in\text{E}_{\!D}^{}\;\!,
$ 
where the set $D=\{(x,y)\colon 2+2\;\!x+y\ne 0\}.$
\vspace{0.75ex}

Taking into account Theorem 1.1, we obtain 
\vspace{0.75ex}
$\bigl(12+8\;\!x+4\;\!y+4\;\!xy+3\;\!y^2,\;\! 2\;\!(x+y)\bigr)\in\text{A}_{\R^2}^{}.$

Using the polynomial partial integrals, by  
\vspace{0.35ex}
Corollary~11.5 (under the conditions
$\rho_1^{}=1,$ $\rho_2^{}=2,\ M_0^{}=x+y,\ \gamma_1^{}={}-2,\ \gamma_2^{}=1),$
\vspace{0.25ex}
we can build the autonomous first integral
\\[1.5ex]
\mbox{}\hfill                   
$
F\colon (x,y)\to\ 
\dfrac{12+8\;\!x+4\;\!y+4\;\!xy+3\;\!y^2}{(2+2\;\!x+y)^2}
$
\ for all 
$(x,y)\in D.
\hfill
$
\\[1.75ex]
\indent
Using the exponential partial integral, 
by Corollary 11.2 (under $\lambda=1),$
we can build nonautonomous first integral
\\[1.75ex]
\mbox{}\hfill                   
$
F_1^{}\colon (t,x,y)\to\ 
e^{{}-t}\exp\dfrac{x+y}{2+2\;\!x+y}=
\exp\Bigl({}-t+\dfrac{x+y}{2+2\;\!x+y}\Bigr)
$
\ for all 
$(x,y)\in\Omega_{_0}^{},
\hfill
$
\\[2ex]
where the set $\Omega_{_0}^{}=\{(t,x,y)\colon 2+2\;\!x+y\ne 0\}.$
\vspace{0.5ex}

Now, using the functional ambiguity of first integral (Property 0.1),
we can represent the first integral $F_1^{}$ in the form
\\[1.5ex]
\mbox{}\hfill                   
$
\Psi\colon (t,x,y)\to\ 
{}-t+\dfrac{x+y}{2+2\;\!x+y}
$
\ for all $(x,y)\in\Omega_{_0}^{}.
\hfill
$
\\[1.75ex]
\indent
Since the first integrals $F$ and $\Psi$ are functionally independent, 
we see that the first integrals $F$ and $\Psi$ are 
an integral basis on the set $\Omega_{_0}^{}$ for system (12.6).
\vspace{1ex}

{\bf 12.5.}
The system [24, 23]
\\[1.75ex]
\mbox{}\hfill                   % (12.7)
$
\dfrac{dx}{dt}={}-y\;\!\bigl(2\;\!x^2+y^2+(x^2+y^2)^2\bigr),
\qquad
\dfrac{dy}{dt}=x\;\!\bigl(2\;\!x^2+y^2+2\;\!(x^2+y^2)^2\bigr)
$
\hfill (12.7)
\\[2ex]
has 
$
(2\;\!x^2+y^2,\;\! {}-2\;\!xy)\in\text{A}_{\R^2}^{}$ and (by Property 5.10)
\\[1.75ex]
\mbox{}\hfill                   
$
\bigl((x^2+y^2,\;\! 2\;\!xy(x^2+y^2)),\;\! (1, 1, {}-2\;\!xy)\bigr)\in\text{B}_{\R^2}^{}.
\hfill
$ 
\\[1.5ex]
The aliquant polynomial partial integral and 
the exponential partial integral (Theorem 5.3) have the same cofactor.
Therefore, by Corollary 11.7, the function
\\[1.5ex]
\mbox{}\hfill                   
$
F\colon (x,y)\to\ 
(2\;\!x^2+y^2)\exp\dfrac{{}-1}{x^2+y^2}
$
\ for all 
$(x,y)\in\R^2\backslash\{(0,0)\}
\hfill
$
\\[1.75ex]
is an autonomous first integral on the domain
$\R^2\backslash\{(0,0)\}$ of system (12.7).
\vspace{1ex}

{\bf 12.6.}
The system [26]
\\[1.75ex]
\mbox{}\hfill                  
$
\dfrac{dx}{dt}=x+xy+y^2,
\qquad
\dfrac{dy}{dt}=y+x^2-xy+2\;\!y^2
\hfill
$
\\[2ex]
has  
$\bigl(x+i\;\!y,\;\! 1+2\;\!y+i\;\!(x-y)\bigr)\in\text{H}_{\R^2}^{}$ (Theorem 6.1).
Hence (Theorem 6.2)
\\[1.75ex]
\mbox{}\hfill                   
$
\bigl(x^2+y^2,\;\! 2\;\!(1+2\;\!y)\bigr)\in\text{A}_{\R^2}^{},\ \
\Bigl(\exp\arctan\dfrac{y}{x}\,,\, x-y\Bigr)\in\text{E}_{\!D}^{}\;\!,\ \
D=\{(x,y)\colon x\ne 0\}\;\!.
\hfill
$
\\[1.75ex]
\indent
Also this system has 
$\bigl(x-y,\;\! {}-(x-y-1)\bigr)\in\text{A}_{\R^2}^{}$ 
\vspace{1ex}
(Theorem 2.1). Hence (Property 1.4)
$\bigl(x-y-1,\;\! {}-(x-y)\bigr)\in\text{A}_{\R^2}^{}.$ 
\vspace{1ex}

By Property 11.11, 
$(x-y-1)\exp\arctan\dfrac{y}{x}\in \text{I}_{D}^{}\;\!.$
\vspace{0.75ex}

By Property 1.9, 
\vspace{0.75ex}
$\Bigl(\dfrac{x-y-1}{x-y}\,,\;\! {}-1\Bigr)\in \text{J}_{\!D_{_{0}}}^{},\ 
D_{_{\!0}}=\{(x,y)\colon x-y\ne 0\}.$
Therefore (Corollary~11.2),
$\dfrac{x-y-1}{x-y}\,e^{\,t}\in  \text{I}_{\R\times D_{_{0}}}^{}.$
\vspace{1.75ex}

{\bf 12.7.}
The differential equation
\\[1.5ex]
\mbox{}\hfill                   % (12.8)
$
\Bigl({}-x+\dfrac{1}{2}\,y^2\Bigr)\;\!dx-
\bigl(xy+2\;\!y^3\bigr)\;\!dy=0
$
\hfill (12.8)
\\[1.5ex]
on $\R^2$ has the complex-valued polynomial partial integral
\\[1.5ex]
\mbox{}\hfill                   
$
w\colon (x,y)\to\ 
x+i\;\!y^2
$
\ for all $(x,y)\in\R^2
\hfill
$
\\[1ex]
with cofactor 
\\[1ex]
\mbox{}\hfill                   
$
W\colon y\to\ 
y-2\;\!i\;\!y
$
\ for all $y\in\R.
\hfill
$
\\[1ex]
\indent
By Theorem 6.2,
\\[1ex]
\mbox{}\hfill                   
$
\bigl(x^2+y^4,\;\! 2\;\!y\bigr)\in\text{J}_{\R^2}^{}\;\!,
\quad
\Bigl(\exp\arctan\dfrac{y^2}{x}\,,\;\! {}-2\;\!y\Bigr)\in\text{E}_{\!D}^{}\;\!,\ \
D=\{(x,y)\colon x\ne 0\}.
\hfill
$
\\[1.75ex]
By Corollary 11.5 (under 
$\rho_1^{}=1,\ \rho_2^{}={}-1,\ M_0^{}=2\;\!y,\ \gamma_1^{}=\gamma_2^{}=1),$
the function 
\\[1.75ex]
\mbox{}\hfill                   
$
F\colon (x,y)\to\ 
\bigl(x^2+y^4\bigr)\;\!
\exp\arctan\dfrac{y^2}{x}
$
\ for all $(x,y)\in D
\hfill
$
\\[1.75ex]
is an general integral on the set $D$ of equation (12.8).
\vspace{1ex}

{\bf 12.8.}
For the first time we paid attention to the application of conditional partial integrals 
in the articles [15, 16, 12], when we did the qualitative investigation 
of the behavior of trajectories on the projective phase plane for system (12.2) [17 -- 20].

Consider the differential equation from [12] 
\\[1.5ex]
\mbox{}\hfill                   % (12.9)
$
(x+x^2-y^2+a)\;\!dx-
(y+x^2-y^2+a)\;\!dy=0
\quad 
(a\in\R)
$
\hfill (12.9)
\\[1.5ex]
as an example of the differential equation with the conditional partial integral and 
the conditional integrating factor.

The polynomial
\\[0.5ex]
\mbox{}\hfill                   
$
p\colon (x,y)\to\ 
x^2-y^2+a
$
\ for all $(x,y)\in\R^2
\hfill
$
\\[1.5ex]
is a partial integral on $\R^2$ with cofactor 
\\[1.5ex]
\mbox{}\hfill                   
$
M_1^{}\colon (x,y)\to\ 
2\;\!(x-y)
$
\ for all $(x,y)\in\R^2.
\hfill
$
\\[1.5ex]
\indent
The exponential function
\\[1.25ex]
\mbox{}\hfill                   
$
{\rm g}\colon (x,y)\to\ 
\exp(x-y)
$
\ for all $(x,y)\in\R^2
\hfill
$
\\[1.5ex]
is a conditional partial integral on $\R^2$ with cofactor 
\\[1.5ex]
\mbox{}\hfill                   
$
M_2^{}\colon (x,y)\to\ 
{}-(x-y)
$
\ for all $(x,y)\in\R^2.
\hfill
$
\\[1.5ex]
\indent
By Corollary 11.5 (under 
$\rho_1^{}=2,\ \rho_2^{}={}-1,\ M_0^{}=x-y,\ \gamma_1^{}=1,\ \gamma_2^{}=2),$
the function
\\[1.75ex]
\mbox{}\hfill                   
$
F\colon (x,y)\to\ 
\bigl(x^2-y^2+a\bigr)\;\!
\exp\bigl(2\;\!(x-y)\bigr)
$
\ for all $(x,y)\in \R^2
\hfill
$
\\[1.75ex]
is a general integral on $\R^2$ of equation (12.9).
\vspace{0.5ex}

Note that
\\[0.35ex]
\mbox{}\hfill                   
$
{\rm div}\;\!{\frak d}(x,y)=2\;\!(x-y)
$
\ for all $(x,y)\in \R^2.
\hfill
$
\\[1.75ex]
\indent
Since the cofactor $M_1^{}={\rm div}\;\!{\frak d},$
we see that, by Property 8.9 (under $\rho_1^{}=1,\ \gamma_1^{}={}-1),$
the rational function
\\[1.75ex]
\mbox{}\hfill                   
$
\mu_1^{}\colon (x,y)\to\ 
\dfrac{1}{x^2-y^2+a}
$
\ for all $(x,y)\in D=\{(x,y)\colon x^2-y^2+a\ne 0\}
\hfill
$
\\[1.75ex]
is an integrating factor on the set $D$ of equation (12.9).
\vspace{0.75ex}

Since the cofactor $M_2^{}={}-\dfrac{1}{2}\,{\rm div}\;\!{\frak d},$
we see that, by Property 10.11 (under $\rho_1^{}=\dfrac{1}{2}\,,\ \gamma_1^{}=2),$
the exponential function
\\[1.25ex]
\mbox{}\hfill                   
$
\mu_2^{}\colon (x,y)\to\ 
\exp\bigl(2\;\!(x-y)\bigr)
$
\ for all $(x,y)\in \R^2
\hfill
$
\\[1.5ex]
is a conditional integrating factor on $\R^2$ of equation (12.9).
\vspace{0.5ex}

Using Jacobi's property of integrating factors (Property 0.2), we obtain 
the general integral $F=\mu_2^{}/\mu_1^{}$ of equation (12.9).
\vspace{1.25ex}

{\bf 12.9.}
{\sl The generalized Darboux problem}.
Consider an $n\!$-th order autonomous differential system
\\[1ex]
\mbox{}\hfill                     % (12.10)
$
\dfrac{dx_i^{}}{dt}=X_i^{}(x_1^{},\ldots, x_n^{}),
\quad
i=1,\ldots, n,
$
\hfill (12.10)
\\[2ex]
where $X_i^{}\colon\R^n\to\R,\ i=1,\ldots,n,$
\vspace{0.35ex}
are polynomials in the dependent variables $x_1^{},\ldots, x_n^{}$
with constant coefficients such that 
$\max\{\deg X_i^{}\colon i=1,\ldots, n\}=d\geq 1.$
\vspace{1.25ex}

{\bf Theorem 12.2.}
{\it 
Suppose we know\footnote{
$\tbinom{k}{n}=\dfrac{n!}{k!\;\!(n-k)!}$ is 
a binomial coefficient.
} 
$m=\tbinom{n+d-1}{n}$
\vspace{0.5ex}
partial integrals ${\rm g}_j^{}$ on the domain $X^{\;\!\prime}$ 
of  the differential system {\rm (12.10)},
$\gamma_j^{}\in\R, $
\vspace{0.5ex}
and the set  $X_{_{\!0}}^{}\subset X^{\;\!\prime}$ such that 
${\rm g}_j^{{}^{\!\scriptsize\gamma_j^{}}}\in C^{1}X_{_0}^{},
\linebreak 
j=1,\ldots, m.$
Then the function
\\[1.25ex]
\mbox{}\hfill                   
$
\displaystyle
F\colon x\to\ 
\prod\limits_{j=1}^m
{\rm g}_j^{{}^{\!\scriptsize\gamma_j^{}}}(x)
$
\ for all $x\in X_{_{\!0}}^{}
\hfill
$
\\[1ex]
is an autonomous first integral 
\vspace{0.15ex}
or an autonomous last multiplier 
on the set $X_{_{\!0}}^{}$ of the differential system} (12.10).
\vspace{0.35ex}

{\sl Proof}
\vspace{0.5ex}
 is analogous to the proof of Theorem 12.1 given that  an
$n\!$-variable polynomial of degree $d-1$ has $\tbinom{n+d-1}{n}$ monomials. $\k$
\vspace{1ex}

In [9, pp. 45 -- 47] and [21], 
\vspace{0.5ex}
the generalized Darboux problem (Theorem 12.2) for 
the differential system (12.10) 
was solved in the presence of  $\tbinom{n+d-1}{n}$
\vspace{0.5ex}
autonomous polynomial partial integrals that define integral manifolds. 

In [25, 31], 
the generalized Darboux problem for the differential system (12.10) 
and also for a multidimensional differential system 
and the differential system (0.1) 
\vspace{0.5ex}
was solved in the presence of  
$\tbinom{n+d-1}{n}$
\vspace{0.5ex}
polynomial partial integrals, taking into account their multiplicities, and 
conditional partial integrals.

Let us remark that the paper [25]
is the first publication in which we used the term-composite 
"conditional partial integral"
(based on the fact that this partial integral does not define an integral manifold). 
Also, note that in this article, we first introduced the concept of multiple polynomial partial integral. And instead of words "$\!\varkappa\!$-multiple polynomial partial integral"
was used word turnover 
"polynomial partial integral with weight $\varkappa\!$" 
(based on the fact that this partial integral is considered 
to be $\varkappa$ partial integrals in the construction of first integral or last multiplier).
The compositional term
"multiple polynomial partial integral"
we began to use since the paper [30].

The generalized Darboux problem and 
methods of its solution were developed in the researches
[1, 2, 12, 23, 25, 30 -- 56].
\vspace{1ex}

{\bf 12.10.}
{\sl The generalized Darboux problem} for 
the generalized Riccati-Abel equation
\\[1.5ex]
\mbox{}\hfill                     % (12.11)
$
\displaystyle
\dfrac{dx}{dt}=
\sum\limits_{i=0}^n
a_i^{}(t)\;\!x^{n-i}
\quad
(n\geq 2)
$
\hfill (12.11)
\\[1.5ex]
with coefficients $a_i^{}\in C^1T,\ i=0,\ldots,n-1,\ a_{_0}^{}\!(t)\not\equiv 0.$
\vspace{0.5ex}

Let the functions
\\[1.5ex]
\mbox{}\hfill                     % (12.12)
$
\displaystyle
{\rm g}_j^{}\colon (t,x)\to\, {\rm g}_j^{}(t,x)
\quad
\forall (t,x)\in\Omega,
\quad
j=1,\ldots, n,
$
\hfill (12.12)
\\[1.5ex]
be partial integrals with cofactors $M_j^{},\ j=1,\ldots,n,$ respectively,
on the domain $\Omega\subset\R^2$ of equation (12.11).
\vspace{0.35ex}

By Definition 1.1, $\deg_x^{} M_j^{}\leq n-1,\ j=1,\ldots, n,$ i.e.,
\\[1.5ex]
\mbox{}\hfill                     % (12.13)
$
\displaystyle
M_j^{}\colon (t,x)\to\ 
\sum\limits_{s=0}^{n-1}\alpha_{sj}^{}(t)\;\!x^{n-s-1}
$
\ for all $(t,x)\in\Xi^{\;\!\prime},
\quad
j=1,\ldots, n,
$
\hfill (12.13)
\\[1.5ex]
where the coefficients 
$\alpha_{sj}^{}\in C^1T^{\;\!\prime},\ s=0,\ldots,n-1,\ j=1,\ldots, n.$
\vspace{0.5ex}

The set $\Omega_{_0}^{}\subset\Omega$
such that  
${\rm g}_j^{{}^{\!\scriptsize\gamma_j^{}}}\in C^{1}\Omega_{_0}^{},\ 
\gamma_{j}^{}\in\R,\  j=1,\ldots, n.$
\vspace{0.75ex}

From Property 11.11 it follows that the function 
\\[1.5ex]
\mbox{}\hfill                   % (12.14)
$
\displaystyle
F\colon (t,x)\to\ 
\prod\limits_{j=1}^n
{\rm g}_j^{{}^{\!\scriptsize\gamma_j^{}}}(t,x)
$
\ for all $(t,x)\in \Omega_{_0}^{}
$
\hfill (12.14)
\\[1.5ex]
is a general integral on the set $\Omega_{_0}^{}$ of equation (12.11)
if and only if  
\\[1.5ex]
\mbox{}\hfill                  
$
\displaystyle
\sum\limits_{j=1}^{n}\gamma_{j}^{}\,
\sum\limits_{s=0}^{n-1}\alpha_{sj}^{}(t)\;\!x^{n-s-1}=0
$
\ for all $(t,x)\in\Xi^{\;\!\prime}.
\hfill
$
\\[1.5ex]
\indent
This identity holds if the functional determinant 
\\[1.75ex]
\mbox{}\hfill                  
$
\displaystyle
\triangle(t)=
\begin{vmatrix}
\alpha_{01}^{}(t) & \alpha_{02}^{}(t) & \dots & \alpha_{0n}^{}(t) 
\\[1.25ex]
\alpha_{11}^{}(t) & \alpha_{12}^{}(t) & \dots & \alpha_{1n}^{}(t) 
\\[0.75ex]
\hdotsfor{4} 
\\[0.75ex]
\alpha_{n-1{,}1}^{}(t) & \alpha_{n-1{,}2}^{}(t) & \dots & \alpha_{n-1{,}n}^{}(t) 
\end{vmatrix}
$
\ for all $t\in T^{\;\!\prime}
\hfill
$
\\[1.75ex]
is equal to the identical zero on the domain $T^{\;\!\prime}.$
\vspace{0.35ex}

Thus, we have
\vspace{0.5ex}

{\bf Theorem 12.3.}\!\!
\vspace{0.25ex}
{\it
The function {\rm (12.14)} is a general integral on the set $\Omega_{_0}^{}\!$
of equation {\rm (12.11)} if and only if the functions {\rm (12.12)} are
\vspace{0.25ex}
partial integrals on the domain $\Omega$ of equation {\rm (12.11)}
with cofactors {\rm (12.13)} such that the determinant  
$\triangle(t)=0$ for all $t\in T^{\;\!\prime}.$
}
\vspace{0.5ex}

From Property 8.8 it follows that the function (12.14) is 
an integrating factor of equation (12.11) if and only if
\\[1.25ex]
\mbox{}\hfill                  
$
\displaystyle
\sum\limits_{j=1}^{n}\gamma_{j}^{}\,
\sum\limits_{s=0}^{n-1}\alpha_{sj}^{}(t)\;\!x^{n-s-1}={}- {\rm div}\,{\frak d}(t,x)
$
\ for all $(t,x)\in\Xi^{\;\!\prime}.
\hfill
$
\\[1.25ex]
\indent
Taking into account that  
\\[1.25ex]
\mbox{}\hfill                  
$
\displaystyle
{\rm div}\,{\frak d}(t,x)=
\sum\limits_{i=0}^{n-1}(n-i)\;\!a_{i}^{}(t)\;\!x^{n-i-1} 
$
\ for all $(t,x)\in\Xi,
\hfill
$
\\[1.25ex]
we obtain that this identity is true if and only if the following system of identities holds
\\[1.5ex]
\mbox{}\hfill                  % (12.15)
$
\displaystyle
\sum\limits_{j=1}^{n}\alpha_{sj}^{}(t)\;\!\gamma_{j}^{}={}- 
(n-s)\;\!a_{s}^{}(t)
$
\ for all $t\in T^{\;\!\prime},
\quad
s=0,\ldots, n-1.
$
\hfill (12.15)
\\[1.5ex]
\indent
Using the determinant $\triangle,$ we can construct $n$ functional determinants 
$\triangle_j^{}$ by replacing $j\!$-th column on the column of  functions
\\[1.5ex]
\mbox{}\hfill
$
{}-n\;\!a_{_0}^{}(t),{}-(n-1)\;\!a_1^{}(t),\ldots,{}-a_{n-1}^{}(t).
\hfill
$
\\[1.5ex]
\indent
If the functional determinants $\triangle,\ \triangle_1^{},\ldots, \triangle_n^{}$ 
such that 
\\[1.5ex]
\mbox{}\hfill                  % (12.16)
$
\displaystyle
\dfrac{\triangle_j^{}(t)}{\triangle(t)}\in\R,
\quad
j=1,\ldots, n,
$
\hfill (12.16)
\\[1.5ex]
then the system of identities (12.15) is solvable with respect 
$\gamma_1^{},\ldots,\gamma_n^{},$
and we got the integrating factor of equation (12.11) in the form (12.14).
\vspace{0.35ex}

So, under $\triangle\not\equiv 0$ we proved
\vspace{0.5ex}

{\bf Theorem 12.4.}
\vspace{0.35ex}
{\it
The function {\rm (12.14)} is an integrating factor on the set $\Omega_{_0}^{}$
of equation {\rm (12.11)} if and only if the functions {\rm (12.12)} are
\vspace{0.25ex}
partial integrals on the domain $\Omega$ of equation {\rm (12.11)} 
with cofactors {\rm (12.13)} such that the determinant $\triangle$ 
\vspace{0.25ex}
is not the identity zero on the domain $T^{\;\!\prime}$ and
the condition {\rm (12.16)} holds.
}
\vspace{0.5ex}

In the paper [10], the generalized Darboux problem (Theorems 12.3 and 12.4) for 
equation (12.11) was solved in the presence of $n$ integral curves defined by polynomial partial integrals. And if $n=3,$ then the similar problem is solved in [11, pp. 54 -- 61].
\vspace{1ex}

{\bf 12.11.}
Since the Abel equation of first kind
\\[1.5ex]
\mbox{}\hfill                     % (12.17)
$
\displaystyle
\dfrac{dx}{dt}=
\dfrac{1}{t}\,x-\dfrac{1}{2\;\!t^2}\,x^3
$
\hfill (12.17)
\\[1.5ex]
is the Bernoulli equation, we see that this equation is integrated by classical methods [3].
\vspace{0.15ex}

Let us construct a general integral of equation (12.17) 
\vspace{0.15ex}
by partial integrals of this equation.

The polynomial 
\\[1ex]
\mbox{}\hfill
$
p\colon (t,x)\to\ {}-t+x^2
$
\ for all 
$(t,x)\in\R^2
\hfill
$
\\[1.5ex]
is a partial integral with cofactor 
\\[1.5ex]
\mbox{}\hfill
$
M\colon (t,x)\to\ 
\dfrac{1}{t}-\dfrac{1}{t^2}\,x^2
$
\ for all 
$(t,x)\in\Xi_{_0}^{}
\hfill
$
\\[1.75ex]
on the set $\Xi_{_0}^{}=T_{_{\!0}}^{}\times\R,\ T_{_{\!0}}^{}=\R\backslash\{0\}.$
\vspace{0.75ex}

Since the cofactor $M={}-\dfrac{1}{t^2}\,p,$ we see that,
\vspace{0.5ex}
by Property 5.10 (under $k=1,\ M_0^{}={}-1/t^2),$
the polynomial partial integral $p$ is double
$\bigl((p, p\;\!M_0^{}), (1, 1, {}-M_0^{})\bigr)\in \text{B}_{_{\Xi_{_0}}}.$ 
\vspace{0.75ex}

By Property 5.13, the function 
\\[1.5ex]
\mbox{}\hfill
$
{\rm g}\colon (t,x)\to\ 
\exp\dfrac{1}{x^2-t}
$
\ for all 
$(t,x)\in\Omega_{_0}^{},
\quad \Omega_{_0}^{}=\{(t,x)\colon x^2-t\ne 0\},
\hfill
$
\\[2ex]
is an exponential partial integral with cofactor 
\\[1.5ex]
\mbox{}\hfill
$
N\colon t\to\ \dfrac{1}{t^2}
$
\ for all 
$t\in T_{_0}^{}
\hfill
$
\\[1.5ex]
on the set $\Omega^{\;\!\prime}=\Omega_{_0}^{}\!\cap\Xi_{_0}^{}.$
\vspace{0.5ex}

By Property 11.8 (under $\varphi=N),$ the function 
\\[1.5ex]
\mbox{}\hfill
$
F\colon (t,x)\to\ 
\exp\dfrac{1}{t}\,
\exp\dfrac{1}{x^2-t}=
\exp\dfrac{x^2}{t\;\!(x^2-t)}
$
\ for all $(t,x)\in\Omega^{\;\!\prime},
\hfill
$
\\[1.75ex]
and taking into account Property 0.1, we have also the function
\\[1.5ex]
\mbox{}\hfill
$
\Psi\colon (t,x)\to\ 
\dfrac{x^2}{t\;\!(x^2-t)}
$
\ for all $(t,x)\in\Omega^{\;\!\prime}
\hfill
$
\\[2ex]
are a general integral on the set $\Omega^{\;\!\prime}$ of equation (12.17).
\vspace{1ex}

{\bf 12.12.}
{\sl Inverse problem}:
construction of differential systems by partial integrals.

Suppose $({\rm g}_j^{}, M_j^{})\in \text{J}_{_{\Omega}},\ j=1,\ldots, n.$
\vspace{0.5ex}
Then, using the existence criterion of partial integral (Theorem 1.1), we have
\\[1.5ex]
\mbox{}\hfill                     % (12.18)
$
\displaystyle
\sum\limits_{i=1}^n\, X_{i}^{}(t,x)\, 
\partial_{{}_{\scriptstyle x_i^{}}}{\rm g}_j^{}(t,x)=
{\rm g}_j^{}(t,x)\;\!M_j^{}(t,x)-
\partial_{{}_{\scriptstyle t}}{\rm g}_j^{}(t,x)
$ 
for all 
$(t,x)\in \Omega,
\
j\!=\!1,\ldots, n.
$
\hfill (12.18)
\\[1.75ex]
\indent
Let the Jacobian of partial integrals $\!{\rm g}_1^{},\ldots, {\rm g}_n^{}\!$
with respect to the variables $\!x_1^{},\ldots, x_n^{}\!$ be
\\[1.5ex]
\mbox{}\hfill                     % (12.19)
$
\displaystyle
\triangle(t,x)=\dfrac{{\rm D}({\rm g}_1^{},\ldots, {\rm g}_n^{})}{{\rm D}(x_1^{},\ldots, x_n^{})}
\not\equiv 0
\ \ \text{on}\ \ \Omega.
$
\hfill (12.19)
\\[1.5ex]
\indent
Using the Jacobian $\triangle,$ we can build $n$ 
functional determinants $\triangle_j^{}$ by replacing $j\!$-th column on the column of the functions
\\[1.25ex]
\mbox{}\hfill                   
$
{\rm g}_1^{}\;\!M_1^{}- 
\partial_{{}_{\scriptstyle t}}\;\!{\rm g}_1^{},\ 
{\rm g}_2^{}\;\!M_2^{}- 
\partial_{{}_{\scriptstyle t}}\;\!{\rm g}_2^{}\;\!,\ \ldots\ ,\
{\rm g}_n^{}\;\!M_n^{}- 
\partial_{{}_{\scriptstyle t}}\;\!{\rm g}_n^{}. 
\hfill
$
\\[1.5ex]
\indent
From the system of identities (12.18) it follows that 
\\[1.5ex]
\mbox{}\hfill               
$
X_i^{}(t,x)=\dfrac{\triangle_i^{}(t,x)}{\triangle(t,x)}
$
\ for all $(t,x)\in \Omega_{_0}\subset\Omega,
\quad
i=1,\ldots,n.
\hfill
$
\\[1.5ex]
\indent
If the quotients
\\[1.5ex]
\mbox{}\hfill                     % (12.20)
$
\dfrac{\triangle_i^{}(t,x)}{\triangle(t,x)}\in 
\text{P}_{_{\Xi_{_0}}},
\quad
i=1,\ldots,n.
\quad
\Xi_{_0}\subset \Xi^{\;\!\prime},
$
\hfill (12.20)
\\[1.5ex]
then we obtain
\vspace{0.5ex}

{\bf Theorem 12.5.}
\vspace{0.75ex}
{\it
If $({\rm g}_j^{}, M_j^{})\in \text{\rm J}_{_{\Omega}},\ j=1,\ldots, n,$
such that the conditions {\rm (12.19)} and {\rm (12.20)} are true,
then the system {\rm (0.1)} has the form}
\\[1.5ex]
\mbox{}\hfill               
$
\dfrac{dx_i^{}}{dt}=\dfrac{\triangle_i^{}(t,x)}{\triangle(t,x)}\,,
\quad
i=1,\ldots,n.
\hfill
$
\\[2ex]
\indent
{\bf Corollary 12.1.}\!\!
\vspace{0.35ex}
{\it
If the autonomous differential system {\rm (12.10)}
has $n$ functionally indepen\-dent partial integrals
\vspace{0.35ex}
$\!{\rm g}_j^{}\colon x\!\to\! {\rm g}_j^{}(x)\!$ for all $x\!\in\! X^{\;\!\prime}\!$ 
with cofactors 
$M_j^{}\colon\R^n\!\to\!\R,\, j\!=\!1,\ldots, n,$ 
%respectively, 
\vspace{0.35ex}
such that the conditions 
$\dfrac{\triangle_i^{}(x)}{\triangle(x)}\in \text{\rm P}_{_{\scriptstyle\!\R^n}},\
i\!=\!1,\ldots,n,$ hold, 
then this system has the form}
\\[1.5ex]
\mbox{}\hfill               
$
\dfrac{dx_i^{}}{dt}=\dfrac{\triangle_i^{}(x)}{\triangle(x)}\,,
\quad
i=1,\ldots,n.
\hfill
$
\\[2ex]
\indent
{\bf Remark 12.1.}
\vspace{0.35ex}
If the partial integral ${\rm g}_j^{}$ with cofactor $M_j^{},\ j\in\{1,\ldots,n\},$
is exponential, i.e.,
$(\exp\omega_j^{},M_j^{})\in \text{\rm E}_{_{\Omega}}\;\!,$
\vspace{0.5ex}
then (taking into account Theorem 3.1) in $j\!$-th identity from the system (12.18)
and in the determinants $\triangle,\, \triangle_1^{},\ldots, \triangle_n^{}$
\vspace{0.15ex}
from Theorem 2.5 and Corollary~12.1 
will be useful to make the formal replacement of the functions 
${\rm g}_j^{}$ by the functions $\omega_j^{}$
\vspace{0.35ex}
and of the product ${\rm g}_j^{}\;\!M_j^{}$ by the function $M_j^{}\;\!.$ 
\vspace{1ex}

{\bf 12.13.}
Suppose the system (12.2) has the polynomial partial integral
\\[1.5ex]
\mbox{}\hfill
$
p\colon (x,y)\to\ 
x^2+y^2+a
$
\ for all $(x,y)\in\R^2
\quad
(a\in\R)
\hfill
$
\\[1.5ex]
with cofactor
\\[1ex]
\mbox{}\hfill
$
M_1^{}\colon (x,y)\to\ 2\;\!(x+y)
$
\ for all $(x,y)\in\R^2
\hfill
$
\\[1.5ex]
and the conditional partial integral
\\[1.5ex]
\mbox{}\hfill
$
{\rm g}\colon (x,y)\to\ 
\exp\omega(x,y)=\exp(x-y)
$
\ for all $(x,y)\in\R^2
\hfill
$
\\[1.5ex]
with cofactor
\\[1ex]
\mbox{}\hfill
$
M_2^{}\colon (x,y)\to\ {}-(x+y)
$
\ for all $(x,y)\in\R^2.
\hfill
$
\\[1.5ex]
\indent
The Jacobian
\\[1.5ex]
\mbox{}\hfill                  
$
\displaystyle
\triangle(x,y)=
\left|\!\!
\begin{array}{cc}
\partial_{x}^{}p & \partial_{y}^{}p
\\[1.25ex]
\partial_{x}^{}\omega & \partial_{y}^{}\omega
\end{array}
\!\!\right|
=
\left|\!\!
\begin{array}{cr}
2\;\!x & 2\;\!y
\\[1ex]
1 & {}-1
\end{array}
\!\!\right|
={}-2\;\!(x+y)
$
\ for all $(x,y)\in\R^2
\hfill
$
\\[1.75ex]
is not the identical zero on $\R^2$ 
\vspace{0.25ex}
(the functions $p$ and $\omega$ are functionally independent on $\R^2).$

The determinants on $\R^2$
\\[1.75ex]
\mbox{}\hfill                  
$
\displaystyle
\triangle_{{}_{X}}(x,y)=
\left|\!\!
\begin{array}{cc}
p\;\!M_1^{} & \partial_{y}^{}p
\\[1.25ex]
M_2^{} & \partial_{y}^{}\omega
\end{array}
\!\!\right|
=
\left|\!\!
\begin{array}{cr}
2\;\!(x^2+y^2+a)(x+y) & 2\;\!y
\\[1ex]
{}-(x+y) & {}-1
\end{array}
\!\!\right|
= 2\;\!(x+y)(y-x^2-y^2-a)
%\quad
%\forall (x,y)\in\R^2
\hfill
$
\\[1.75ex]
and
\\[1.75ex]
\mbox{}\hfill                  
$
\displaystyle
\triangle_{{}_{Y}}(x,y)=
\left|\!\!
\begin{array}{cc}
\partial_{x}^{}p & p\;\!M_1^{} 
\\[1.25ex]
\partial_{x}^{}\omega & M_2^{}  
\end{array}
\!\!\right|
=
\left|\!\!
\begin{array}{cc}
2\;\!x & 2\;\!(x^2+y^2+a)
\\[1ex]
1 & {}-(x+y)
\end{array}
\!\!\right|
= {}-2\;\!(x+y)(x+x^2+y^2+a).
%\quad
%\forall (x,y)\in\R^2
\hfill
$
\\[1.75ex]
\indent
Using Remark 12.1, by Corollary 12.1, we have
\\[1.5ex]
\mbox{}\hfill               
$
\dfrac{dx}{dt}=\dfrac{\triangle_{{}_{X}}^{}(x,y)}{\triangle(x,y)}=
{}-y+x^2+y^2+a,
\hfill
$
\\[0.15ex]
\mbox{}\hfill (12.21)
\\[0.15ex]
\mbox{}\hfill               
$
\dfrac{dy}{dt}=\dfrac{\triangle_{{}_{Y}}^{}(x,y)}{\triangle(x,y)}=
x+x^2+y^2+a.\ \ \ \,
\hfill
$
\\[1.5ex]
\indent
From Corollary 11.5 (under $\rho_1^{}=2,\ \rho_2^{}={}-1,\ M_0^{}=x+y,\ 
\gamma_1^{}=1,\ \gamma_2^{}=2)$ it follows that the transcendental function
\\[1.5ex]
\mbox{}\hfill
$
F\colon (x,y)\to\ (x^2+y^2+a)\exp\bigl(\;\!2\;\!(x-y)\bigr)
$
\ for all $(x,y)\in\R^2
\hfill
$
\\[1.5ex]
is an autonomous first integral on $\R^2$ of system (12.21).

We stress that, if $a>0,$ then the polynomial $p(x,y)\ne 0$ for all $(x,y)\in\R^2.$
In other words, if $a>0,$ then the polynomial partial integral  $p$ 
does not define trajectories of the differential system (12.21).

Thus, if $a>0,$ then we construct the differential system (12.21) by two partial integrals 
$p$ and ${\rm g}$ such that  these partial integrals do not define trajectories.

In the paper [22] (see also [4, pp. 480 -- 481]),
the differential system (12.2) built by two partial integrals under the condition that 
these partial integrals define trajectories.
\vspace{1ex}

{\bf 12.14.}
\vspace{0.35ex}
Suppose the system (12.2) on the phase plane $\R^2$ has 
the multiple polynomial partial integral
$\bigl((p, M), (h, q, N)\bigr)\in \text{\rm B}_{_{\scriptstyle\R^2}}.$
By Theorem 5.2, 
\\[1.5ex]
\mbox{}\hfill
$
X(x,y)\;\!\partial_{x}^{}p(x,y)+
Y(x,y)\;\!\partial_{y}^{}p(x,y)=p(x,y)\;\!M(x,y)
$
\ for all $(x,y)\in\R^2,
\hfill
$
\\[2ex]
\mbox{}\hfill
$
X(x,y)\;\!\partial_{x}^{}q(x,y)+
Y(x,y)\;\!\partial_{y}^{}q(x,y)=h\;\!q(x,y)\;\!M(x,y)+p^{\;\!h}(x,y)\;\!N(x,y)
$
for all $(x,y)\!\in\!\R^2.
\hfill
$
\\[1.5ex]
\indent
Since the polynomials $p$ and $q$ are functionally independent on $\R^2,$
from this system of identities, we obtain
\\[1.5ex]
\mbox{}\hfill               
$
X(x,y)=\dfrac{\triangle_{{}_{X}}^{}(x,y)}{\triangle(x,y)}\,,
\quad
Y(x,y)=\dfrac{\triangle_{{}_{Y}}^{}(x,y)}{\triangle(x,y)}
$
\ for all $(x,y)\in D\subset\R^2,
\hfill
$
\\[1.5ex]
where the Jacobian of the polynomials $p$ and $q$ with respect to $x$ and $y$
\\[1.5ex]
\mbox{}\hfill                  % (12.22)
$
\displaystyle
\triangle(x,y)=
\left|\!\!
\begin{array}{cc}
\partial_{x}^{}p & \partial_{y}^{}p
\\[1.25ex]
\partial_{x}^{}q & \partial_{y}^{}q
\end{array}
\!\!\right|
\not\equiv 0
\ \ \text{on} \ \ \R^2,
$
\hfill (12.22)
\\[1.5ex]
the determinants
\\[1.5ex]
\mbox{}\hfill                 
$
\displaystyle
\triangle_{{}_{X}}^{}(x,y)=
\left|\!\!
\begin{array}{cc}
p\;\!M & \partial_{y}^{}p
\\[1.25ex]
h\;\!q\;\!M+p^{\;\!h}N & \partial_{y}^{}q
\end{array}
\!\!\right|,
\quad \ \
\triangle_{{}_{Y}}^{}(x,y)=
\left|\!\!
\begin{array}{cc}
\partial_{x}^{}p & p\;\!M 
\\[1.25ex]
\partial_{x}^{}q & h\;\!q\;\!M+p^{\;\!h}N
\end{array}
\!\!\right|
$
\ for all $(x,y)\in\R^2.
\hfill
$
\\[1.5ex]
\indent
If the quotients
\\[1.5ex]
\mbox{}\hfill                  % (12.23)
$
\displaystyle
\dfrac{\triangle_{{}_{X}}^{}(x,y)}{\triangle(x,y)}\in 
\text{\rm P}_{_{\!\scriptstyle\R^2}}
\quad \
\text{and} \ \,
\quad
\dfrac{\triangle_{{}_{Y}}^{}(x,y)}{\triangle(x,y)}\in 
\text{\rm P}_{_{\!\scriptstyle\R^2}},
$
\hfill (12.23)
\\[1.5ex]
then we obtain
\vspace{0.5ex}

{\bf Theorem 12.6.}
\vspace{0.35ex}
{\it
If the differential system {\rm (12.2)} on $\R^2$ has 
the multiple polynomial partial integral
$\bigl((p, M), (h, q, N)\bigr)\in \text{\rm B}_{_{\scriptstyle\R^2}}$
\vspace{0.35ex}
such that the conditions {\rm (12.22)} and {\rm (12.23)} are true,
then this system has the form}
\\[1.5ex]
\mbox{}\hfill               
$
\dfrac{dx}{dt}=
\dfrac{\triangle_{{}_{X}}^{}(x,y)}{\triangle(x,y)}\,,
\qquad
\dfrac{dy}{dt}=
\dfrac{\triangle_{{}_{Y}}^{}(x,y)}{\triangle(x,y)}\,.
\hfill
$
\\[2ex]
\indent
{\bf 12.15.}
\vspace{0.35ex}
Suppose the system (12.2) on the phase plane $\R^2$ 
has the multiple polynomial partial integral
$\bigl((p, M), (h, q, N)\bigr)\in \text{\rm B}_{_{\scriptstyle\R^2}}$ such that
\\[1.5ex]
\mbox{}\hfill
$
p=x^2+y^2-1,
\quad
M=xy,
\quad
h=1,
\quad
q=x^2-y^2-1,
\quad
N={}-xy.
\hfill
$
\\[1.5ex]
\indent
The Jacobian
\\[1.5ex]
\mbox{}\hfill                
$
\displaystyle
\triangle(x,y)=
\left|\!\!
\begin{array}{cc}
\partial_{x}^{}p & \partial_{y}^{}p
\\[1.25ex]
\partial_{x}^{}q & \partial_{y}^{}q
\end{array}
\!\!\right|
=
\left|\!\!
\begin{array}{cr}
2\;\!x & 2\;\!y
\\[1ex]
2\;\!x & {}-2\;\!y
\end{array}
\!\!\right|=
{}-8\;\!xy
$
\ for all $(x,y)\in \R^2
\hfill
$
\\[1.5ex]
is not the identical zero on $\R^2.$

The determinants on $\R^2$
\\[1.5ex]
\mbox{}\hfill                 
$
\displaystyle
\triangle_{{}_{X}}^{}(x,y)=
\left|\!\!
\begin{array}{cc}
p\;\!M & \partial_{y}^{}p
\\[1.25ex]
h\;\!q\;\!M+p^{\;\!h}N & \partial_{y}^{}q
\end{array}
\!\!\right|
=
\left|\!\!
\begin{array}{cr}
(x^2+y^2-1)\;\!xy & 2\;\!y
\\[1ex]
{}-2\;\!xy^3 & {}-2\;\!y
\end{array}
\!\!\right|=
{}-2\;\!xy^2(x^2-y^2-1),
%\quad
%\forall (x,y)\in \R^2
\hfill
$
\\[2.5ex]
\mbox{}\hfill
$
\triangle_{{}_{Y}}^{}(x,y)=
\left|\!\!
\begin{array}{cc}
\partial_{x}^{}p & p\;\!M 
\\[1.25ex]
\partial_{x}^{}q & h\;\!q\;\!M+p^{\;\!h}N
\end{array}
\!\!\right|
=
\left|\!\!
\begin{array}{cc}
2\;\!x & (x^2+y^2-1)\;\!xy 
\\[1ex]
2\;\!x & {}-2\;\!xy^3  
\end{array}
\!\!\right|=
{}-2\;\!x^2y\;\!(x^2+3y^2-1).
%\quad
%\forall (x,y)\in\R^2.
\hfill
$
\\[1.75ex]
\indent
By Theorem 12.6,
\\[1.5ex]
\mbox{}\hfill               % (12.24)
$
\dfrac{dx}{dt}=\dfrac{\triangle_{{}_{X}}^{}(x,y)}{\triangle(x,y)}=
\dfrac{1}{4}\,y\;\!(x^2-y^2-1),
\quad
\dfrac{dy}{dt}=\dfrac{\triangle_{{}_{Y}}^{}(x,y)}{\triangle(x,y)}=
\dfrac{1}{4}\,x\;\!(x^2+3y^2-1).
$
\hfill (12.24)
\\[1.75ex]
\indent
From Theorem 5.3 it follows that the system (12.24) has two partial integrals
\\[1.5ex]
\mbox{}\hfill
$
(x^2+y^2-1,\;\! xy)\in \text{\rm A}_{_{\scriptstyle\R^2}}
\quad \text{and}\quad
\Bigl(\exp\dfrac{x^2-y^2-1}{x^2+y^2-1}\,,\;\! {}-xy\Bigr)\in \text{\rm E}_{_{D}},
\hfill
$
\\[1.5ex]
were the set $D=\{(x,y)\colon x^2+y^2-1\ne 0\}.$
\vspace{0.5ex}

By Corollary 11.5 (under $\rho_1^{}=1,\ \rho_2^{}={}-1,\ M_0^{}=xy,\ 
\gamma_1^{}=\gamma_2^{}=1),$ the transcendental function
\\[1.5ex]
\mbox{}\hfill
$
F\colon (x,y)\to\ 
(x^2+y^2-1)\exp\dfrac{x^2-y^2-1}{x^2+y^2-1}
$
\ for all $(x,y)\in D
\hfill
$
\\[1.5ex]
is an autonomous first integral on the set $D$ of system (12.24).
\vspace{1ex}

{\bf 12.16.}
\vspace{0.35ex}
Suppose the system (12.2) on the phase plane $\R^2$ has 
the complex-valued polynomial partial integral
$(u+i\;\!v, U+i\;\!V)\in \text{\rm H}_{_{\scriptstyle\R^2}}.$
By Theorem 6.1, 
\\[1.5ex]
\mbox{}\hfill
$
X(x,y)\;\!\partial_{x}^{}u(x,y)+
Y(x,y)\;\!\partial_{y}^{}u(x,y)=u(x,y)\;\!U(x,y)-v(x,y)\;\!V(x,y)
$
\ for all $(x,y)\in\R^2,
\hfill
$
\\[2ex]
\mbox{}\hfill
$
X(x,y)\;\!\partial_{x}^{}v(x,y)+
Y(x,y)\;\!\partial_{y}^{}v(x,y)=
u(x,y)\;\!V(x,y)+v(x,y)\;\!U(x,y)
$
\ for all $(x,y)\in\R^2.
\hfill
$
\\[1.5ex]
\indent
Since the polynomials $u$ and $v$ are functionally independent on $\R^2,$
from this system of identities, we obtain
\\[1.5ex]
\mbox{}\hfill               
$
X(x,y)=\dfrac{\triangle_{{}_{X}}^{}(x,y)}{\triangle(x,y)}\,,
\quad
Y(x,y)=\dfrac{\triangle_{{}_{Y}}^{}(x,y)}{\triangle(x,y)}
$
\ for all $(x,y)\in D\subset\R^2,
\hfill
$
\\[1.5ex]
where the Jacobian of the polynomials $u$ and $v$ with respect to $x$ and $y$
\\[1.5ex]
\mbox{}\hfill                  % (12.25)
$
\displaystyle
\triangle(x,y)=
\left|\!\!
\begin{array}{cc}
\partial_{x}^{}u & \partial_{y}^{}u
\\[1.25ex]
\partial_{x}^{}v & \partial_{y}^{}v
\end{array}
\!\!\right|
\not\equiv 0
\ \ \text{on} \ \ \R^2,
$
\hfill (12.25)
\\[1.5ex]
the determinants
\\[1.5ex]
\mbox{}\hfill                 
$
\displaystyle
\triangle_{{}_{X}}^{}(x,y)=
\left|\!\!
\begin{array}{cc}
u\;\!U-v\;\!V & \partial_{y}^{}u
\\[1.25ex]
u\;\!V+v\;\!U  & \partial_{y}^{}v
\end{array}
\!\!\right|,
\quad \ \
\triangle_{{}_{Y}}^{}(x,y)=
\left|\!\!
\begin{array}{cc}
\partial_{x}^{}u & u\;\!U-v\;\!V 
\\[1.25ex]
\partial_{x}^{}v & u\;\!V+v\;\!U
\end{array}
\!\!\right|
$
\ for all $(x,y)\in\R^2.
\hfill
$
\\[1.5ex]
\indent
If the quotients 
\\[1.5ex]
\mbox{}\hfill                  % (12.26)
$
\displaystyle
\dfrac{\triangle_{{}_{X}}^{}(x,y)}{\triangle(x,y)}\in 
\text{\rm P}_{_{\!\scriptstyle\R^2}}
\quad \
\text{and} \ \,
\quad
\dfrac{\triangle_{{}_{Y}}^{}(x,y)}{\triangle(x,y)}\in 
\text{\rm P}_{_{\!\scriptstyle\R^2}},
$
\hfill (12.26)
\\[1ex]
then we have
\vspace{0.5ex}

{\bf Theorem 12.7.}
\vspace{0.35ex}
{\it
If the differential system {\rm (12.2)} on $\R^2$ has 
the  complex-valued polynomial partial integral 
$(u+i\;\!v, U+i\;\!V)\in \text{\rm H}_{_{\scriptstyle\R^2}}$
such that the conditions {\rm (12.25)} and {\rm (12.26)} are true,
then this system has the form}
\\[1.5ex]
\mbox{}\hfill               
$
\dfrac{dx}{dt}=
\dfrac{\triangle_{{}_{X}}^{}(x,y)}{\triangle(x,y)}\,,
\qquad
\dfrac{dy}{dt}=
\dfrac{\triangle_{{}_{Y}}^{}(x,y)}{\triangle(x,y)}\,.
\hfill
$
\\[2ex]
\indent
{\bf Corollary 12.2.}
\vspace{0.35ex}
{\it
The differential system {\rm (12.2)} has
the  complex-valued polynomial partial integral  
$(x+i\;\!y, U+i\;\!V)\in \text{\rm H}_{_{\scriptstyle\R^2}}$
if and only if this system has the form}
\\[1.5ex]
\mbox{}\hfill               
$
\dfrac{dx}{dt}=
x\;\!U(x,y)-y\;\!V(x,y), 
\qquad
\dfrac{dy}{dt}=
x\;\!V(x,y)+y\;\!U(x,y).
\hfill
$
\\[2ex]
\indent
{\bf 12.17.}
\vspace{0.35ex}
Suppose the system (12.2) on the phase plane $\R^2$ has 
the  complex-valued polynomial partial integral 
$(x+i\;\!y, 1+i)\in \text{\rm H}_{_{\scriptstyle\R^2}}.$ By Corollary 12.2,
\\[1.5ex]
\mbox{}\hfill     % (12.27)
$
\dfrac{dx}{dt}=x-y,
\qquad
\dfrac{dy}{dt}=x+y.
$
\hfill (12.27)
\\[1.5ex]
\indent
From Theorem 6.2 it follows that 
the system (12.27) has two partial integrals
\\[1.5ex]
\mbox{}\hfill
$
(x^2+y^2,\;\! 2)\in \text{\rm A}_{_{\scriptstyle\R^2}}
\quad \text{and}\quad
\Bigl(\exp\arctan\dfrac{y}{x}\,,\;\! 1\Bigr)\in \text{\rm E}_{_{D}},
\hfill
$
\\[1.5ex]
where the set $D=\{(x,y)\colon x\ne 0\}.$
\vspace{0.5ex}

By Corollary 11.5 
\vspace{0.15ex}
(under $\rho_1^{}=2,\ \rho_2^{}=1,\ M_0^{}=1,\ 
\gamma_1^{}=1,\ \gamma_2^{}={}-2),$ 
taking into account Property 3.13, we obtain the function 
\\[1.5ex]
\mbox{}\hfill
$
F\colon (x,y)\to\ 
(x^2+y^2)\exp\Bigl(2\arctan\dfrac{y}{x}\Bigr)
$
\ for all $(x,y)\in D
\hfill
$
\\[1.5ex]
is an autonomous first integral on the set $D$ of system (12.27).
\vspace{0.35ex}

Using the partial integrals and Corollary 11.2 (under $\lambda=2$ and $\lambda=1),$ 
we can build the non-autonomous first integrals of system (12.27)
\\[1.5ex]
\mbox{}\hfill
$
F_1^{}\colon (t,x,y)\to\ 
e^{{}-2\;\!t}(x^2+y^2)
$
\ for all $(t,x,y)\in\R^3 
\hfill
$
\\[0.75ex]
and
\\[0.75ex]
\mbox{}\hfill
$
F_2^{}\!\colon\! (t,x,y)\to
e^{{}-t}\exp\arctan\dfrac{y}{x}=
\exp\Bigl(\!\!{}-t+\arctan\dfrac{y}{x}\Bigr)\!
$
for all 
$(t,x,y)\!\in\!\Omega_{{}_0}\!=\! \{(t,x,y)\colon x\!\ne\! 0\}.
\hfill
$
\\[1.5ex]
\indent
Taking into account the functional ambiguity of first integral (Property 0.1), 
we can present the first integral $F_2^{}$ in the form 
\\[1.5ex]
\mbox{}\hfill
$
\Psi\colon (t,x,y)\to\ 
{}-t+\arctan\dfrac{y}{x}
$
\ for all $(t,x,y)\in\Omega_{{}_0}.
\hfill
$
\\[1.5ex]
\indent
Since the first integrals $F,\ F_1^{},\ \Psi$
\vspace{0.25ex}
are pairwise functional independence on $\R^3,$
we see that each of the collections of functions 
$\{F, F_1^{}\}, \ \{F, \Psi\},\ \{F_1^{}, \Psi\}$
\vspace{0.25ex}
is an integral basis on the set $\Omega_{{}_0}$ 
of the differential system (12.27).
\vspace{1ex}

{\bf 12.18.}
By Corollary 12.2, the system (12.2) has 
$\bigl(x+i\;\!y,\;\! x-y+i\;\!(x+y)\bigr)\in \text{\rm H}_{_{\scriptstyle\R^2}}$
if and only if 
\\[1.5ex]
\mbox{}\hfill     % (12.28)
$
\dfrac{dx}{dt}=x^2-2\;\!xy-y^2,
\qquad
\dfrac{dy}{dt}=x^2+2\;\!xy-y^2.
$
\hfill (12.28)
\\[1.75ex]
\indent
From Theorem 6.2 it follows that the system (12.28) has
\\[1.75ex]
\mbox{}\hfill
$
\bigl(x^2+y^2,\;\! 2\;\!(x-y)\bigr)\in \text{\rm A}_{_{\scriptstyle\R^2}},
\quad 
\Bigl(\exp\arctan\dfrac{y}{x}\,,\;\! x+y\Bigr)\in \text{\rm E}_{_{D}},
\quad 
D=\{(x,y)\colon x\ne 0\}.
\hfill
$
\\[1.75ex]
\indent
Taking into account Theorem  1.1, 
\vspace{0.25ex}
we obtain that the differential system 
of the second order (12.28) has
$
\bigl(x+y,\;\! 2\;\!(x-y)\bigr)\in \text{\rm A}_{_{\scriptstyle\R^2}}.
$
\vspace{0.75ex}

The polynomial partial integrals have the same cofactor.
By Corollary~11.7, the rational function
\\[1.5ex]
\mbox{}\hfill
$
F\colon (x,y)\to\ 
\dfrac{x+y}{x^2+y^2}
$
\ for all $(x,y)\in\R^2\backslash\{(0,0)\} 
\hfill
$
\\[1.5ex]
is an autonomous first integral on the domain $\R^2\backslash\{(0,0)\}$
of system (12.28).
\vspace{1.5ex}

{\bf 12.19.}
Using the form of right parts of system
([18, p. 88; 30] under $\lambda={}-1,\ \eta=1)$ 
\\[1.5ex]
\mbox{}\hfill     % (12.29)
$
\dfrac{dx}{dt}=x\;\!(x^2+y^2+\lambda)-y\;\!(x^2+y^2+\eta),
\quad
\dfrac{dy}{dt}=x\;\!(x^2+y^2+\eta)+y\;\!(x^2+y^2+\lambda),
$
\hfill (12.29)
\\[1.75ex]
where $\lambda,\eta\in\R, \ |\lambda|+|\eta|\ne 0,\ \lambda\ne\eta,$
by Corollary 12.2, we obtain 
\\[1.75ex]
\mbox{}\hfill
$
\bigl(x+i\;\!y,\;\! x^2+y^2+\lambda+i\;\!(x^2+y^2+\eta)\bigr)\in 
\text{\rm H}_{_{\scriptstyle\R^2}}.
\hfill
$
\\[1.5ex]
\indent
From Theorem 5.3, we get the following
\\[1.75ex]
\mbox{}\hfill
$
\bigl(x^2+y^2,\;\! 2\;\!(x^2+y^2+\lambda)\bigr)\in \text{\rm A}_{_{\scriptstyle\R^2}},
\quad 
\Bigl(\exp\arctan\dfrac{y}{x}\,,\;\! x^2+y^2+\eta\Bigr)\in \text{\rm E}_{_{D}},
\quad 
D=\{(x,y)\colon x\ne 0\}.
\hfill
$
\\[2ex]
\indent
{\sl Case} $\lambda\ne 0.$
\vspace{0.35ex}
Since 
$
\bigl(x^2+y^2,\;\! 2\;\!(x^2+y^2+\lambda)\bigr)\in \text{\rm A}_{_{\scriptstyle\R^2}},
$
we see that, by Property 1.4, 
$
\bigl(x^2+y^2+\lambda,\;\! 2\;\!(x^2+y^2)\bigr)\in \text{\rm A}_{_{\scriptstyle\R^2}}.
$
\vspace{0.75ex}

The identity
\\[1.5ex]
\mbox{}\hfill
$
2\;\!\gamma_1^{}\;\!(x^2+y^2+\lambda)+
\gamma_2^{}\;\!(x^2+y^2+\eta)+
2\;\!\gamma_3^{}\;\!(x^2+y^2)=0
$ 
\ for all $(x,y)\in\R^2
\hfill
$
\\[1.5ex]
is true, for example, under $\gamma_1^{}={}-\eta,\ \gamma_2^{}=2\;\!\lambda,\ 
\gamma_3^{}=\eta-\lambda.$
By Property 11.11, the function
\\[2ex]
\mbox{}\hfill
$
F\colon (x,y)\to\ 
\dfrac{(x^2+y^2+\lambda)^{\eta-\lambda}}{(x^2+y^2)^{\eta}}\,
\exp\Bigl(2\;\!\lambda\;\!\arctan\dfrac{y}{x}\;\!\Bigr)
$
\ for all $(x,y)\in D
\hfill
$
\\[1.5ex]
is an autonomous first integral on the set $D$ of system (12.29) under $\lambda\ne 0.$
\vspace{0.35ex}

The difference of cofactors
\\[1.5ex]
\mbox{}\hfill
$
2\;\!(x^2+y^2)-2\;\!(x^2+y^2+\lambda)={}-2\;\!\lambda
$ 
\ for all $(x,y)\in\R^2.
\hfill
$
\\[2ex]
\indent
It follows from Property 1.9 that 
$
\Bigl(\dfrac{x^2+y^2+\lambda}{x^2+y^2}\,,{}-2\;\!\lambda\Bigr)\in 
\text{\rm J}_{_{\!\scriptstyle D_{_0}}},\ 
D_{_0}=\R^2\backslash\{(0,0)\}.
$
This yields that (Corollary 11.2) the function
\\[1.75ex]
\mbox{}\hfill
$
F_1^{}\colon (t,x,y)\to\ 
\dfrac{x^2+y^2+\lambda}{x^2+y^2}\,e^{\;\!2\;\!\lambda\;\!t}
$
\ for all $(t,x,y)\in\Omega_{_{0}}
\hfill
$
\\[1.5ex]
is an non-autonomous first integral on the domain 
\vspace{0.35ex}
$\Omega_{_{0}}=\R\times D_{_{0}}$ 
of the differential system (12.29) under the condition $\lambda\ne 0.$
\vspace{0.35ex}

The linear combination of cofactors
\\[1.5ex]
\mbox{}\hfill
$
2\;\!(x^2+y^2)-2\;\!(x^2+y^2+\eta)={}-2\;\!\eta
$ 
\ for all $(x,y)\in\R^2.
\hfill
$
\\[1.75ex]
\indent
It follows from Property 1.9 that 
\vspace{0.75ex}
$
\Bigl((x^2+y^2+\lambda)\;\!\exp\Bigl({}-2\;\!\arctan\dfrac{y}{x}\;\!\Bigr),{}-2\;\!\eta\Bigr)\in 
\text{\rm J}_{_{\!\scriptstyle D}}.
$
This implies that (Corollary 11.2) the function
\\[1.75ex]
\mbox{}\hfill
$
F_2^{}\colon (t,x,y)\to\ 
(x^2+y^2+\lambda)\;\!\exp\Bigl(2\;\!\eta\;\!t-2\;\!\arctan\dfrac{y}{x}\;\!\Bigr)
$
\ for all $(t,x,y)\in\Omega^{\;\!\prime}
\hfill
$
\\[1.5ex]
is an non-autonomous first integral on the set 
$\Omega^{\;\!\prime}=\R\times D$ 
of the differential system (12.29) under the condition $\lambda\ne 0.$
\vspace{0.35ex}

The linear combination of cofactors
\\[1.5ex]
\mbox{}\hfill
$
2\;\!(x^2+y^2+\lambda)-2\;\!(x^2+y^2+\eta)=2\;\!(\lambda-\eta)
$ 
\ for all $(x,y)\in\R^2.
\hfill
$
\\[1.75ex]
\indent
It follows from Property 1.9 that 
$
\Bigl((x^2+y^2)\;\!\exp\Bigl({}-2\;\!\arctan\dfrac{y}{x}\;\!\Bigr),\, 2\;\!(\lambda-\eta)\Bigr)\in 
\text{\rm J}_{_{\!\scriptstyle D}}.
$
Therefore (Corollary 11.2) the function
\\[1.75ex]
\mbox{}\hfill
$
F_3^{}\colon (t,x,y)\to\ 
(x^2+y^2)\;\!\exp\Bigl(2\;\!(\eta-\lambda)\;\!t-2\;\!\arctan\dfrac{y}{x}\;\!\Bigr)
$
\ for all $(t,x,y)\in\Omega^{\;\!\prime}
\hfill
$
\\[1.5ex]
is an non-autonomous first integral on the set 
$\Omega^{\;\!\prime}$ 
of the differential system (12.29) under the condition $\lambda\ne 0.$
\vspace{0.35ex}

Since the first integrals $F,\ F_1^{},\ F_2^{},\ F_3^{}$
\vspace{0.35ex}
are pairwise functional independence on $\R^3,$
we see that each of the collections of functions 
\vspace{0.35ex}
$\{F, F_1^{}\},\ \{F, F_2^{}\},\ \{F, F_3^{}\},\ \{F_1^{}, F_2^{}\},\ 
\{F_1^{}, F_3^{}\},$ $\{F_2^{}, F_3^{}\}$
is an integral basis on the corresponding set of 
system (12.29) under $\lambda\ne 0.$
\vspace{1ex}

{\sl Case} $\lambda= 0.$ Then,
\\[1.5ex]
\mbox{}\hfill     % (12.30)
$
\dfrac{dx}{dt}=x\;\!(x^2+y^2)-y\;\!(x^2+y^2+\eta),
\quad
\dfrac{dy}{dt}=x\;\!(x^2+y^2+\eta)+y\;\!(x^2+y^2)
\quad
(\eta\ne 0).
$
\hfill (12.30)
\\[2ex]
\indent
We have 
\vspace{0.5ex}
$\bigl(x+i\;\!y,\;\! x^2+y^2+i\;\!(x^2+y^2+\eta)\bigr)\in \text{\rm H}_{_{\scriptstyle\R^2}},$
i.e., (Theorem 5.3)
\\[1.75ex]
\mbox{}\hfill
$
\bigl(x^2+y^2,\;\! 2\;\!(x^2+y^2)\bigr)\in \text{\rm A}_{_{\scriptstyle\R^2}},
\quad 
\Bigl(\exp\arctan\dfrac{y}{x}\,,\;\! x^2+y^2+\eta\Bigr)\in \text{\rm E}_{_{D}}.
\hfill
$
\\[1.75ex]
\indent
Since
$
\bigl(x^2+y^2,\;\! 2\;\!(x^2+y^2)\bigr)\in \text{\rm A}_{_{\scriptstyle\R^2}},
$
we see that, by Property 5.10, 
\\[1.75ex]
\mbox{}\hfill
$
\bigl((x^2+y^2,\;\! 2\;\!(x^2+y^2)),\, (1, 1, {}-2)\bigr)\in 
\text{\rm B}_{_{\scriptstyle\R^2}}.
\hfill
$
\\[1.5ex]
Therefore (Theorem 5.3),
$
\Bigl(\exp\dfrac{1}{x^2+y^2}\,,{}-2\Bigr)\in 
\text{\rm E}_{_{\!\scriptstyle D_{_0}}}.
$
\vspace{0.5ex}

The identity 
\\[1.5ex]
\mbox{}\hfill
$
2\;\!\gamma_1^{}\;\!(x^2+y^2)+
\gamma_2^{}\;\!(x^2+y^2+\eta)-2\;\!\gamma_3^{}=0
$ 
\ for all $(x,y)\in\R^2
\hfill
$
\\[1.75ex]
is true, for example, under 
$\gamma_1^{}={}-1,\ \gamma_2^{}=2,\ \gamma_3^{}=\eta.$
By Property 11.11, the function
\\[2ex]
\mbox{}\hfill
$
F_4^{}\colon (x,y)\to\ 
\dfrac{1}{x^2+y^2}\,
\exp\Bigl(\;\!\dfrac{\eta}{x^2+y^2}+2\;\!\arctan\dfrac{y}{x}\;\!\Bigr)
$
\ for all $(x,y)\in D
\hfill
$
\\[1.5ex]
is an autonomous first integral on the set $D$ of system (12.30).
\vspace{0.35ex}

The linear combination of cofactors
\\[1.5ex]
\mbox{}\hfill
$
2\;\!(x^2+y^2)-2\;\!(x^2+y^2+\eta)={}-2\;\!\eta
$ 
\ for all $(x,y)\in\R^2.
\hfill
$
\\[1.75ex]
\indent
From Property 1.9, we get
$
\Bigl((x^2+y^2)\;\!\exp\Bigl({}-2\;\!\arctan\dfrac{y}{x}\;\!\Bigr),{}-2\;\!\eta\Bigr)\in 
\text{\rm J}_{_{\!\scriptstyle D}}.
$
Then (by Corollary~11.2) the function
\\[1.75ex]
\mbox{}\hfill
$
F_5^{}\colon (t,x,y)\to\ 
(x^2+y^2)\;\!\exp\Bigl(2\;\!\eta\;\!t-2\;\!\arctan\dfrac{y}{x}\;\!\Bigr)
$
\ for all $(t,x,y)\in\Omega^{\;\!\prime}
\hfill
$
\\[1.5ex]
is an non-autonomous first integral on the set 
$\Omega^{\;\!\prime}$ of system (12.30).
\vspace{1ex}

Since
\vspace{0.75ex}
$
\Bigl(\exp\dfrac{1}{x^2+y^2}\,,{}-2\Bigr)\in 
\text{\rm E}_{_{\!\scriptstyle D_{_0}}},
$
we see that, by Corollary 11.2, 
$
\exp\Bigl(2\;\!t+\dfrac{1}{x^2+y^2}\;\!\Bigr)\in 
\text{\rm I}_{_{\scriptstyle \Omega_{_0}}}.
$
Thus (Property 0.1) the function
\\[1.75ex]
\mbox{}\hfill
$
F_6^{}\colon (t,x,y)\to\ 
2\;\!t+\dfrac{1}{x^2+y^2}
$
\ for all $(t,x,y)\in\Omega_{_{0}}
\hfill
$
\\[1.5ex]
is an non-autonomous first integral on the domain 
$\Omega_{_{0}}$ of system (12.30).
\vspace{0.5ex}

Since the first integrals $F_4^{},\ F_5^{},\ F_6^{}$
\vspace{0.5ex}
are pairwise functional independence on $\R^3,$
we see that each of the collections of functions 
\vspace{0.35ex}
$\{F_4^{}, F_5^{}\},\ \{F_4^{}, F_6^{}\},\ \{F_5^{}, F_6^{}\}$
is an integral basis on the set $\Omega^{\;\!\prime}$ of 
the differential system (12.30).
\vspace{1ex}

{\bf 12.20.}
Using the form of rights parts of the Darboux system [9]
\\[1.5ex]
\mbox{}\hfill     % (12.31)
$
\dfrac{dx}{dt}=\lambda\;\!x-\eta\;\!y+x\;\!P(x,y),
\quad
\dfrac{dy}{dt}=\eta\;\!x+\lambda\;\!y+y\;\!P(x,y),
$
\hfill (12.31)
\\[2ex]
where $\lambda\in\R,\ \eta\in\R\backslash\{0\},$ and $P$ is a 
polynomial of degree $\deg P\geq 1,$
\vspace{0.75ex}
by Corollary 12.2, we obtain 
$\bigl(x+i\;\!y,\;\! \lambda+P+\eta\;\!i\bigr)\in \text{\rm H}_{_{\scriptstyle\R^2}}.$
From Theorem 5.3 it follows that
\\[2ex]
\mbox{}\hfill
$
\bigl(x^2+y^2,\;\! 2\;\!(\lambda+P)\bigr)\in \text{\rm A}_{_{\scriptstyle\R^2}},
\quad 
\Bigl(\exp\arctan\dfrac{y}{x}\,,\;\! \eta\Bigr)\in \text{\rm E}_{_{D}},
\quad 
D=\{(x,y)\colon x\ne 0\}.
\hfill
$
\\[2ex]
\indent
By Corollary 11.2, 
$
\exp\Bigl({}-\eta\;\!t+\arctan\dfrac{y}{x}\;\!\Bigr)\in 
\text{\rm I}_{_{\scriptstyle \Omega_{_0}}}.
$
Then (Property 0.1) the function
\\[2ex]
\mbox{}\hfill
$
\Psi\colon (t,x,y)\to\ 
\eta\;\!t-\arctan\dfrac{y}{x}
\quad
\forall (t,x,y)\in\Omega_{_{0}}
\hfill
$
\\[1.75ex]
is an non-autonomous first integral on the set 
\vspace{1.25ex}
$\Omega_{_{0}}=\R\times D$ of system (12.31).

{\bf 12.20.A.}
\vspace{0.35ex}
Let $\lambda\ne 0,\ P=x^2+y^2.$
Consider the differential system (the particular case see in [28, p. 43; 29])
\\[1.75ex]
\mbox{}\hfill     % (12.32)
$
\dfrac{dx}{dt}=\lambda\;\!x-\eta\;\!y+x\;\!(x^2+y^2),
\qquad
\dfrac{dy}{dt}=\eta\;\!x+\lambda\;\!y+y\;\!(x^2+y^2).
$
\hfill (12.32)
\\[2ex]
\indent
We have
\vspace{1ex}
$\bigl(x^2+y^2,\;\! 2\;\!(x^2+y^2+\lambda)\bigr)\in \text{\rm A}_{_{\scriptstyle\R^2}},\ 
\Bigl(\exp\arctan\dfrac{y}{x}\,,\;\! \eta\Bigr)\in \text{\rm E}_{_{D}}.
$

Since
$
\bigl(x^2+y^2,\;\! 2\;\!(x^2+y^2+\lambda)\bigr)\in \text{\rm A}_{_{\scriptstyle\R^2}},
$
we see that, by Property 1.4, 
\\[1.75ex]
\mbox{}\hfill
$
\bigl(x^2+y^2+\lambda,\;\! 2\;\!(x^2+y^2)\bigr)\in \text{\rm A}_{_{\scriptstyle\R^2}}.
\hfill
$
\\[1ex]
\indent
The identity 
\\[1.5ex]
\mbox{}\hfill
$
2\;\!\gamma_1^{}\;\!(x^2+y^2+\lambda)+\gamma_2^{}\;\!\eta+
2\;\!\gamma_3^{}\;\!(x^2+y^2)=0
$
\ for all $(x,y)\in\R^2
\hfill
$
\\[1.5ex]
is true, for example, under 
$\gamma_1^{}={}-1,\ \gamma_2^{}=\dfrac{2\;\!\lambda}{\eta}\,,\ 
\gamma_3^{}=1.$
By Property 11.11, the function
\\[1.75ex]
\mbox{}\hfill
$
F\colon (x,y)\to\ 
\dfrac{x^2+y^2+\lambda}{x^2+y^2}\,
\exp\Bigl(\;\!\dfrac{2\;\!\lambda}{\eta}\;\!\arctan\dfrac{y}{x}\;\!\Bigr)
$
\ for all $(x,y)\in D
\hfill
$
\\[1.5ex]
is an autonomous first integral on the set  $D$ of system (12.32).
\vspace{0.35ex}

Since the first integrals $\Psi$ and $F$ are functional independence,
we see that they  are an integral basis on the set 
$\Omega_{_{0}}$ of system (12.32).
\vspace{0.75ex}

{\bf 12.20.B.}
Let $\lambda= 0,\ P=x^2+y^2.$ Then
\\[1.75ex]
\mbox{}\hfill     % (12.33)
$
\dfrac{dx}{dt}={}-\eta\;\!y+x\;\!(x^2+y^2),
\qquad
\dfrac{dy}{dt}=\eta\;\!x+y\;\!(x^2+y^2).
$
\hfill (12.33)
\\[2ex]
\indent
We have
\vspace{0.75ex}
$\bigl(x^2+y^2,\;\! 2\;\!(x^2+y^2)\bigr)\in \text{\rm A}_{_{\scriptstyle\R^2}}$ and 
$\Bigl(\exp\arctan\dfrac{y}{x}\,,\;\! \eta\Bigr)\in \text{\rm E}_{_{D}}.
$

Since 
$
\bigl(x^2+y^2,\;\! 2\;\!(x^2+y^2)\bigr)\in \text{\rm A}_{_{\scriptstyle\R^2}},
$
we see that, by Property 5.10, 
\\[2ex]
\mbox{}\hfill
$
\bigl((x^2+y^2,\;\! 2\;\!(x^2+y^2)),\;\! (1, 1, {}-2)\bigr)\in 
\text{\rm B}_{_{\scriptstyle\R^2}}.
\hfill
$
\\[1.5ex]
Thus (Theorem 5.3),
$
\Bigl(\exp\dfrac{1}{x^2+y^2}\,,\;\! {}-2\Bigr)\in \text{\rm E}_{_{D}}.
$
\vspace{1ex}

From 
$
\gamma_1^{}\;\!\eta -2\;\!\gamma_2^{}=0
$
it follows that, for example, $\gamma_1^{}=2,\ \gamma_2^{}=\eta.$ 
Then, by Property~11.11, we obtain 
\\[1.25ex]
\mbox{}\hfill
$
\exp\Bigl(\;\!\dfrac{\eta}{x^2+y^2}+2\;\!\arctan\dfrac{y}{x}\;\!\Bigr)\in \text{\rm I}_{_{D}},
\hfill
$
\\[1.5ex]
and using the Property 0.1, we get that the function
\\[1.75ex]
\mbox{}\hfill
$
\Psi_1^{}\colon (x,y)\to\ 
\dfrac{\eta}{x^2+y^2}+2\;\!\arctan\dfrac{y}{x}
$
\ for all $(x,y)\in D
\hfill
$
\\[1.5ex]
is an autonomous first integral on the set $D$ of system (12.33).
\vspace{0.35ex}

Since the first integrals $\Psi$ and $\Psi_1^{}$
\vspace{0.25ex}
are functional independence,
we see that they  are an integral basis on the set 
$\Omega_{_{0}}$ of system (12.33).
\vspace{1ex}

{\bf 12.20.C.} 
Let $\lambda= 0,\ P=y^2(x^2+y^2-1).$ Then ([28, p. 43] under $\eta=1)$
\\[1.75ex]
\mbox{}\hfill     % (12.34)
$
\dfrac{dx}{dt}={}-\eta\;\!y+xy^2\;\!(x^2+y^2-1),
\qquad
\dfrac{dy}{dt}=\eta\;\!x+y^3\;\!(x^2+y^2-1).
$
\hfill (12.34)
\\[2ex]
\indent
We have
$\bigl(x+i\;\!y,\;\! y^2(x^2+y^2-1)+\eta\;\!i\bigr)\in \text{\rm H}_{_{\scriptstyle\R^2}},$
i.e. (Theorem 5.3),
\\[1.75ex]
\mbox{}\hfill
$
\bigl(x^2+y^2,\;\! 2\;\!y^2\;\!(x^2+y^2-1)\bigr)\in \text{\rm A}_{_{\scriptstyle\R^2}},
\quad 
\Bigl(\exp\arctan\dfrac{y}{x}\,,\;\! \eta\Bigr)\in \text{\rm E}_{_{D}}.
\hfill
$
\\[1.75ex]
\indent
Since
$
\bigl(x^2+y^2,\;\! 2\;\!y^2\;\!(x^2+y^2-1)\bigr)\in \text{\rm A}_{_{\scriptstyle\R^2}},
$
we see that, by Property 1.4, 
\\[1.5ex]
\mbox{}\hfill
$
\bigl(x^2+y^2-1,\;\! 2\;\!y^2\;\!(x^2+y^2)\bigr)\in \text{\rm A}_{_{\scriptstyle\R^2}}.
\hfill
$
\\[1.5ex]
\indent
The divergence
\\[1.25ex]
\mbox{}\hfill
$
{\rm div}\;\!{\frak d}(t,x,y)=
6\;\!y^2\;\!(x^2+y^2)-4\;\!y^2
$
\ for all $(t,x,y)\in\R^3.
\hfill
$
\\[1.25ex]
\indent
The identity
\\[1.5ex]
\mbox{}\hfill
$
2\;\!\gamma_1^{}\;\!y^2\;\!(x^2+y^2-1)+
2\;\!\gamma_2^{}\;\!y^2\;\!(x^2+y^2)=
{}-6\;\!y^2\;\!(x^2+y^2)+4\;\!y^2
$ 
\ for all $(x,y)\in\R^2
\hfill
$
\\[1.75ex]
is true, for example, under $\gamma_1^{}={}-2,\ \gamma_2^{}={}-1.$
By Property 8.8, the function
\\[2ex]
\mbox{}\hfill
$
\mu\colon (x,y)\to\ 
\dfrac{1}{(x^2+y^2)^2(x^2+y^2-1)}
$
\ for all $(x,y)\in D_{_{0}}
\hfill
$
\\[1.75ex]
is an integrating factor on the set
\vspace{0.25ex}
$\!D_{_{0}}\!=\!\{(x,y)\colon (x^2+y^2)(x^2+y^2-1)\!\ne\! 0\}\!$ 
for the equation of trajectories of system (12.34).
\vspace{1ex}

{\bf 12.21.}
Using the form of right parts of system [27]
\\[1.5ex]
\mbox{}\hfill     % (12.35)
$
\dfrac{dx}{dt}=x-x^3-x^2y-y^3,
\qquad
\dfrac{dy}{dt}=y+x^3-x^2y+xy^2,
$
\hfill (12.35)
\\[2ex]
by Corollary 12.2, we obtain
\\[1.5ex]
\mbox{}\hfill
$
\bigl(x+i\;\!y,\;\! 1-x^2+i\;\!(x^2+y^2)\bigr)\in \text{\rm H}_{_{\scriptstyle\R^2}}.
\hfill
$
\\[1.5ex]
\indent
Therefore (Theorem~5.3),
\\[1.75ex]
\mbox{}\hfill
$
\bigl(x^2+y^2,\;\! 2\;\!(1-x^2)\bigr)\in \text{\rm A}_{_{\scriptstyle\R^2}},
\quad 
\Bigl(\exp\arctan\dfrac{y}{x}\,,\;\! x^2+y^2\Bigr)\in \text{\rm E}_{_{\!D}},
\quad 
D=\{(x,y)\colon x\ne 0\}.
\hfill
$
\\[1.75ex]
\indent
From Theorem 3.1 it follows that  
\\[1.5ex]
\mbox{}\hfill
$
\Bigl(\exp\dfrac{xy}{x^2+y^2}\,,\, x^2-y^2\Bigr)\in 
\text{\rm E}_{_{\!\scriptstyle D_{_0}}}
$
\ and \
$
\Bigl(\exp\dfrac{{}-y^2}{x^2+y^2}\,,\, {}-2\;\!xy\Bigr)\in 
\text{\rm E}_{_{\!\scriptstyle D_{_0}}},
\ D_{_0}=\R^2\backslash\{(0,0)\}.
\hfill
$
\\[1.5ex]
\indent
Hence (Theorem 5.3),
\\[1.5ex]
\mbox{}\hfill
$
\bigl((x^2\!+y^2, 2(1-x^2)), (1, xy, x^2\!-y^2)\bigr)\!\in\! 
\text{\rm B}_{_{\scriptstyle\R^2}}\!
$
and
$
\bigl((x^2\!+y^2, 2(1-x^2)), (1,\! {}-y^2,\!{}-2xy)\bigr)\!\in\! \text{\rm B}_{_{\scriptstyle\R^2}}.
\hfill
$
\\[1.5ex]
\indent
Thus,
\vspace{0.75ex}
$\!\bigl((x+iy, 1-x^2+i(x^2+y^2)), (1, y, x^2-y^2-2ixy)\bigr)\!\in\! 
\text{\rm G}_{_{\scriptstyle\R^2}}\!\!$
(Corollary 7.1).

Since the sum of cofactors
\\[1.5ex]
\mbox{}\hfill     
$
2\;\!(1-x^2)+(x^2+y^2)+(x^2-y^2)=2,
\hfill 
$
\\[1.5ex]
we see that, by Property 1.9,
\\[1.5ex]
\mbox{}\hfill
$
\Bigl((x^2+y^2)\exp\Bigl(\;\!\dfrac{xy}{x^2+y^2}+\arctan\dfrac{y}{x}\;\!\Bigr),\, 2\Bigr)\in 
\text{\rm J}_{_{\!\scriptstyle D}}.
\hfill
$
\\[1.5ex]
\indent
By Corollary 11.2, the function 
\\[1.5ex]
\mbox{}\hfill
$
F\colon (t,x,y)\to\ 
(x^2+y^2)\exp\Bigl({}-2\;\!t+\dfrac{xy}{x^2+y^2}+\arctan\dfrac{y}{x}\;\!\Bigr)
\hfill
$
\\[1.5ex]
is an non-autonomous first integral on the set $\R\times D$ of system (12.35).
\vspace{0.35ex}

The divergence
\\[1.25ex]
\mbox{}\hfill
$
{\rm div}\;\!{\frak d}(t,x,y)=
2\;\!(1-2\;\!x^2)
$
\ for all $(t,x,y)\in\R^3.
\hfill
$
\\[1.25ex]
\indent
The identity
\\[1.5ex]
\mbox{}\hfill
$
2\;\!\gamma_1^{}\;\!(1-x^2)+
\gamma_2^{}\;\!(x^2+y^2)+
\gamma_3^{}\;\!(x^2-y^2)=
{}-2\;\!(1-2\;\!x^2)
$ 
\ for all $(x,y)\in\R^2
\hfill
$
\\[1.75ex]
is true, for example, under $\gamma_1^{}={}-1,\ \gamma_2^{}=\gamma_3^{}=1.$
By Property 8.8, the function
\\[2ex]
\mbox{}\hfill
$
\mu\colon (x,y)\to\ 
\dfrac{1}{x^2+y^2}\,
\exp\Bigl(\;\!\dfrac{xy}{x^2+y^2}+\arctan\dfrac{y}{x}\;\!\Bigr)
$
\ for all $(x,y)\in D
\hfill
$
\\[1.75ex]
is an integrating factor on the set $\!D\!$ 
\vspace{1ex}
for the equation of trajectories of system (12.35).

{\bf 12.22.}
{\it The Jacobi equation} has the form 
\\[2ex]
\mbox{}\hfill        % (12.36)
$
\dfrac{dy}{dx}=
\dfrac{l_2^{}(x,y)-y\;\!l_3^{}(x,y)}{l_1^{}(x,y)-x\;\!l_3^{}(x,y)}\equiv 
\dfrac{Y(x,y)}{X(x,y)}\,,
\hfill
$
\\
\mbox{}\hfill (12.36)
\\
\mbox{}\hfill
$
l_i^{}(x,y)=a_{i}^{}\;\!x+b_{i}^{}\;\!y+c_i^{}\;\!,
\quad
a_{i}^{},b_{i}^{},c_i^{}\in\R, 
\quad
i=1,2,3.
\hfill
$
\\[2ex]
\indent
The linear differential operator of first-order
\\[2ex]
\mbox{}\hfill
$
{\frak I}(x,y)=X(x,y)\;\!\partial_{{}_{\scriptstyle x}}+Y(x,y)\;\!\partial_{{}_{\scriptstyle y}}
$
\ for all $(x,y)\in\R^2
\hfill
$
\\[2ex]
is the operator of differentiation by virtue of the Jacobi equation (12.36).

The linear function 
\\[1.25ex]
\mbox{}\hfill
$
p\colon (x,y)\to\ 
\alpha\;\! x+\beta y+\gamma
$
\ for all 
$
(x,y)\in\R^2
\quad
(|\alpha|+|\beta|\ne 0)
\hfill
$
\\[1.5ex]
is a polynomial partial integral of equation (12.36)
if and only if 
the following identity holds
\\[2ex]
\mbox{}\hfill
$
{\frak I}\;\!(\alpha\;\! x+\beta y+\gamma)=
(\alpha\;\! x+\beta y+\gamma)(\mu\;\! x+\nu\;\! y+\theta)
$
\ for all 
$(x,y)\in\R^2.
\hfill
$
\\[2ex]
\indent
This identity holds if and only if the system of equations
\\[2ex]
\mbox{}\hfill        % (12.37)
$
(a_1^{}-\lambda)\;\!\alpha+a_2^{}\beta +a_3^{}\gamma=0,
\ \ \
b_1^{}\alpha+(b_2^{}-\lambda)\;\!\beta +b_3^{}\gamma=0,
\ \ \
c_1^{}\alpha+c_2^{}\beta +(c_3^{}-\lambda)\;\!\gamma=0
$
\hfill (12.37)
\\[2ex]
is consistent, where $\mu={}-a_3^{},\ \nu={}-b_3^{},\ \theta={}-c_3^{}+\lambda.$
\vspace{0.5ex}

The linear homogeneous  system (12.37) has an nontrivial solution
$(\alpha,\beta,\gamma)$ if and only if 
the determinant of this system equals zero, i.e.,
\\[1.75ex]
\mbox{}\hfill    % (12.38)
$
\det(A-\lambda E)=0,
$
\hfill (12.38)
\\[1.5ex]
where $E$ is the identity matrix, and the matrix
\\[1.5ex]
\mbox{}\hfill 
$
A=\left(\!\!
\begin{array}{ccc}
a_1^{} & a_2^{} & a_3^{}
\\[0.75ex]
b_1^{} & b_2^{} & b_3^{}
\\[0.5ex]
c_1^{} & c_2^{} & c_3^{}
\end{array}
\!\!\right).
\hfill
$
\\[1.5ex]
\indent
The roots of equation (12.38) are  eigenvalues of matrix $A.$
Further, a solution $(\alpha,\beta,\gamma)$ of system (12.37) is
an eigenvector of matrix $A$ corresponding to the eigenvalue $\lambda.$ 
\vspace{0.35ex}

Thus we proved
\vspace{0.5ex}

{\bf Lemma 12.1.}
\vspace{0.35ex}
{\it
If $(\alpha,\beta,\gamma)$ is an eigenvector corresponding to the eigenvalue $\lambda$ 
of matrix $A,$ then the linear function
\vspace{0.5ex}
$
p\colon (x,y)\to 
\alpha\;\! x+\beta y+\gamma\;\;
$
for all $(x,y)\in\R^2$
$
(|\alpha|+|\beta|\ne 0)
$
is a polynomial partial integral with cofactor
\vspace{0.15ex}
$M\colon (x,y)\to \lambda-l_3^{}(x,y)$ for all $(x,y)\in\R^2$
for the Jacobi equation} (12.36).
\vspace{0.5ex}

Note that Lemma 12.1 does not exclude the possibility, when the eigenvalue $\lambda$ is complex.

We will build a general integral of the Jacobi equation (12.36) by eigenvectors and 
eigenvalues of matrix $A$ with using orders of elementary divisors.
\vspace{0.75ex}

{\small\bf 
12.22.1. Case of three simple elementary divisors.}
\vspace{0.35ex}

{\bf Theorem 12.8.}
\vspace{0.25ex}
{\it
Suppose $\lambda-\lambda_1^{},\ \lambda-\lambda_2^{},\ \lambda-\lambda_3^{}$ are
simple elementary divisors and 
$
(\alpha_1^{},\beta_1^{},\gamma_1^{}),\ 
(\alpha_2^{},\beta_2^{},\gamma_2^{}),\ 
(\alpha_3^{},\beta_3^{},\gamma_3^{}) 
$
\vspace{0.35ex}
are eigenvectors of matrix $A,$ which corresponding to divisors.
Then a general integral on $D$ of the Jacobi equation {\rm(12.36)} is the function
\\[2ex]
\mbox{}\hfill  % (12.39)
$
F\colon (x,y)\to\ 
p_1^{{}^{\scriptstyle h_1^{}}}\!(x,y)\;\!
p_2^{{}^{\scriptstyle h_2^{}}}(x,y)\;\!
p_3^{{}^{\scriptstyle h_3^{}}}(x,y)
$
\ for all $(x,y)\in D,
$
\hfill {\rm (12.39)}
\\[2.25ex]
where $p_i^{}(x,y)=\alpha_i^{}\;\!x+\beta_{i}^{}\;\!y+\gamma_i^{}$ for all $(x,y)\in\R^2,\ i=1,2,3,$
\vspace{0.25ex}
the power exponents $h_1^{},\, h_2^{},\, h_3^{}$ such that 
\\[1ex]
\mbox{}\hfill  % (12.40)
$
h_1^{}+h_2^{}+h_3^{}=0,
\quad
\lambda_1^{}h_1^{}+\lambda_2^{}h_2^{}+\lambda_3^{}h_3^{}=0,
$
\hfill {\rm (12.40)}
\\[2ex]
and the set $D$ such that the functions $p_i^{{}^{\scriptstyle h_i^{}}}\in C^1D,\ i=1,2,3.$}
\vspace{0.75ex}

{\sl Proof.}
\vspace{0.25ex}
By Lemma 12.1, the linear functions $p_i^{}$ are 
polynomial (real or complex-valued) partial integrals with cofactors
$\lambda_i^{}-l_3^{},\ i=1,2,3,$ of the Jacobi equation (12.36). 
\vspace{0.25ex}
Suppose the numbers $h_1^{},\, h_2^{},\, h_3^{}$ 
\vspace{0.25ex}
such that the equalities (12.40) are true.
Then the cofactors of these partial integrals such that the linear combination
\\[2ex]
\mbox{}\hfill
$
h_1^{}\bigl(\lambda_1^{}-l_3^{}(x,y)\bigr)+
h_2^{}\bigl(\lambda_2^{}-l_3^{}(x,y)\bigr)+
h_3^{}\bigl(\lambda_3^{}-l_3^{}(x,y)\bigr)=0
$
\ for all $(x,y)\in\R^2.
\hfill
$
\\[2ex]
\indent
Taking into account Properties 11.11 and 6.6, we see that 
the function (12.39) is a general integral on the set $D$ of the Jacobi equation (12.36). $\k$
\vspace{0.5ex}

Let us remark that in Theorem 12.8 among eigenvalues 
$\lambda_1^{},\, \lambda_2^{},\, \lambda_3^{}$
can be eigenvalues either equal or complex.
\vspace{0.75ex}

{\small\bf 12.22.1.1.\! Case of three distinct real eigenvalues.}\!
\vspace{0.25ex}
If the eigenvalues $\lambda_1^{},\, \lambda_2^{},\, \lambda_3^{}$ of mat\-rix $A$
are real and distinct, then they have simple elementary divisors 
\vspace{0.25ex}
$\lambda-\lambda_1^{},\ \lambda-\lambda_2^{},\ \lambda-\lambda_3^{}.$
In this case, from Theorem 12.8 under 
$h_1^{}=\lambda_2^{}-\lambda_3^{},\ 
h_2^{}=\lambda_3^{}-\lambda_1^{},\ 
h_3^{}=\lambda_1^{}-\lambda_2^{},$ we get 
\vspace{1ex}

{\bf Theorem 12.9.}
\vspace{0.35ex}
{\it
Suppose the eigenvalues $\lambda_1^{},\, \lambda_2^{},\, \lambda_3^{}$ 
of mat\-rix $A$ are real and distinct, and 
$
(\alpha_1^{},\beta_1^{},\gamma_1^{}),\ 
(\alpha_2^{},\beta_2^{},\gamma_2^{}),\ 
(\alpha_3^{},\beta_3^{},\gamma_3^{}) 
$
\vspace{0.5ex}
are eigenvectors corresponding to these eigenvalues.
Then a general integral on $D$ of the Jacobi equation {\rm(12.36)} 
is the function
\\[2ex]
\mbox{}\hfill  % (12.41)
$
F\colon (x,y)\to\ 
p_1^{{}^{\scriptstyle \lambda_2^{}-\lambda_3^{}}}\!(x,y)\;\!
p_2^{{}^{\scriptstyle \lambda_3^{}-\lambda_1^{}}}\!(x,y)\;\!
p_3^{{}^{\scriptstyle \lambda_1^{}-\lambda_2^{}}}\!(x,y)
$
\ for all $(x,y)\in D,
$
\hfill {\rm (12.41)}
\\[2.25ex]
where $p_i^{}(x,y)=\alpha_i^{}\;\!x+\beta_{i}^{}\;\!y+\gamma_i^{}$ for all $(x,y)\in\R^2,\ i=1,2,3,$
\vspace{0.75ex}
the set $D\subset\R^2$ such that the functions  
$p_1^{{}^{\scriptstyle \lambda_2^{}-\lambda_3^{}}},\,
p_2^{{}^{\scriptstyle \lambda_3^{}-\lambda_1^{}}},\,
p_3^{{}^{\scriptstyle \lambda_1^{}-\lambda_2^{}}}$
are continuously differentiable on $D.$}
\vspace{1.25ex}

{\bf Example 12.1.}
The Jacobi equation 
\\[2ex]
\mbox{}\hfill  % (12.42)
$
\dfrac{dy}{dx}=\dfrac{{}-x+5y-1-y\;\!(x-y+3)}{3x-y+1-x\;\!(x-y+3)}
$
\hfill (12.42)
\\[2ex]
such that the matrix
\vspace{0.75ex}
$
A=\left(\!\!\!
\begin{array}{rrr}
3 & {}-1 & 1
\\
{}-1 & 5 & {}-1
\\
1 & {}-1 & 3
\end{array}
\!\!\right)
$
has the eigenvalues $\lambda_1^{}=2,\ \lambda_2^{}=3,\ \lambda_3^{}=6,$ 
and corresponding eigenvectors $(1, 0,{}-1),\ (1, 1, 1),\ (1,{}-2, 1).$ 
\vspace{0.35ex}

By Theorem 12.9, the function 
\\[2ex]
\mbox{}\hfill  % (12.43)
$
F\colon (x,y)\to\ 
\dfrac{(x+y+1)^4}{(x-1)^3(x-2y+1)}
$
\ for all $(x,y)\in D
$
\hfill {\rm (12.43)}
\\[2.25ex]
is a general integral on  
$D=\{(x,y)\colon (x-1)(x-2y+1)\ne 0\}$
\vspace{1.25ex}
of the Jacobi equation (12.42).

{\small\bf 12.22.1.2. 
Case of mul\-ti\-p\-le eigenvalue with simple elementary divisors.}
Suppose a multiple eigenvalue $\!\lambda_1^{}\!$ of matrix $\!A\!$ has 
two simple elementary divisors
$\!\lambda-\lambda_1^{}\!$ and $\lambda-\lambda_1^{}.$
\vspace{0.25ex}
Then this eigenvalue $\lambda_1^{}$ has two 
linearly independent eigenvectors
$(\alpha_1^{},\beta_1^{},\gamma_1^{})$ and
$(\alpha_2^{},\beta_2^{},\gamma_2^{}).$ 
\vspace{0.75ex}

{\bf Theorem 12.10.}
\vspace{0.25ex}
{\it
If a multiple eigenvalue of matrix $A$ corresponding to 
simple elementary divisors and to the linearly independent eigenvectors
$(\alpha_1^{},\beta_1^{},\gamma_1^{})$ and
$(\alpha_2^{},\beta_2^{},\gamma_2^{}),$ 
then the homographic function
\\[2ex]
\mbox{}\hfill  % (12.44)
$
F\colon (x,y)\to\ 
\dfrac{\alpha_1^{}\;\!x+\beta_{1}^{}\;\!y+\gamma_1^{}}{\alpha_2^{}\;\!x+\beta_{2}^{}\;\!y+\gamma_2^{}}
$
\ for all 
$
(x,y)\in D
$
\hfill {\rm (12.44)}
\\[2.25ex]
is a general integral on  
$D=\{(x,y)\colon \alpha_2^{}\;\!x+\beta_{2}^{}\;\!y+\gamma_2^{}\ne 0\}$ 
of the Jacobi equation {\rm(12.36)}.
}
\vspace{0.75ex}

{\sl Proof.}
\vspace{0.25ex}
Let $\lambda$ be a multiple eigenvalue of matrix $A$ with 
corresponding linearly independent eigenvectors
$(\alpha_1^{},\beta_1^{},\gamma_1^{})$ and
$(\alpha_2^{},\beta_2^{},\gamma_2^{}).$
By Lemma 12.1, the linear functions  
\\[1.75ex]
\mbox{}\hfill
$
p_1^{}\colon (x,y)\to\, \alpha_1^{}\;\!x+\beta_{1}^{}\;\!y+\gamma_1^{}
$
\ for all 
$(x,y)\in\R^2
\hfill
$ 
\\[0.5ex] 
and
\\[0.5ex]
\mbox{}\hfill 
$
p_2^{}\colon (x,y)\to\, \alpha_2^{}\;\!x+\beta_{2}^{}\;\!y+\gamma_2^{}
$
\ for all 
$(x,y)\in\R^2
\hfill
$
\\[1.5ex]
are polynomial partial integrals of the Jacobi equation (12.36)
with the same cofactor 
\\[1.5ex]
\mbox{}\hfill
$
M\colon (x,y)\to\ \lambda-l_3^{}(x,y)$ 
\ for all $(x,y)\in\R^2.
\hfill
$
\\[1.5ex]
\indent
From Сorollary 11.7 it follows that the homographic function (12.44) is 
a general integral on the set $D$ of the Jacobi equation (12.36). $\k$
\vspace{0.35ex}

Note also that the statement of Theorem 12.10 we can obtain from Theorem 12.8 under 
the conditions $\lambda_2^{}=\lambda_1^{},$ $h_1^{}=1,\ h_2^{}={}-1,\ h_3^{}=0.$
\vspace{0.35ex}

If the triple eigenvalue $\lambda_1^{}$ of matrix $A$ 
\vspace{0.25ex}
has three simple elementary divisors
$\lambda-\lambda_1^{},\ \lambda-\lambda_1^{},$ $\lambda-\lambda_1^{},$
then the matrix $A=\lambda_1^{}E.$
\vspace{0.25ex}
In this case, the Jacobi equation (12.36) is degenerate, i.e.,
$X(x,y)=Y(x,y)=0$ for all $(x,y)\in\R^2.$
\vspace{1ex}

{\bf Example 12.2.}
The Jacobi equation
\\[1.75ex]
\mbox{}\hfill  % (12.45)
$
\dfrac{dy}{dx}=\dfrac{x-y+1-y\;\!(x+y-1)}{{}-x+y+1-x\;\!(x+y-1)}
$
\hfill (12.45)
\\[2ex]
such that the matrix
\vspace{0.75ex}
$
A=\left(\!\!\!
\begin{array}{rrr}
{}-1 & 1 & 1
\\
1 & {}-1 & 1
\\
1 & 1 & {}-1
\end{array}
\!\!\right)
$
has the eigenvalues
$\lambda_1^{}=\lambda_2^{}={}-2,\ \lambda_3^{}=1$
with the simple elementary divisors
$\lambda+2,\ \lambda+2,\ \lambda-1.$
\vspace{0.25ex}
The double eigenvalue $\lambda_1^{}={}-2$ has 
the linearly independent eigenvectors
\vspace{0.35ex}
$(0, 1,{}-1)$ and $(1, 0, {}-1).$ 

By Theorem 12.10, the homographic function
\\[2ex]
\mbox{}\hfill  % (12.46)
$
F\colon (x,y)\to\ 
\dfrac{y-1}{x-1}
$
\ for all 
$(x,y)\in D
$
\hfill {\rm (12.46)}
\\[2.25ex]
is a general integral on the set  
$D=\{(x,y)\colon x\ne 1\}$
\vspace{0.5ex}
of the Jacobi equation (12.45).

The general integral (12.46) defines the pencil of the linear trajectories
\\[1.5ex]
\mbox{}\hfill
$
\dfrac{y-1}{x-1}=C,
\quad 
C\in [{}-\infty;{}+\infty],
\hfill
$
\\[1.5ex]
of the Jacobi equation (12.45).
\vspace{0.25ex}

The vector $(1, 1, 1)$ is  an eigenvector corresponding to the 
the eigenvalue $\lambda_3^{}=1.$  
\vspace{0.35ex}
By Lemma 12.1, the linear function
$p_3^{}\colon (x,y)\to x+y+1$ for all $(x,y)\in\R^2$
\vspace{0.5ex}
is a polynomial partial integral with cofactor 
$M_3^{}\colon (x,y)\to 2-x-y$ for all $(x,y)\in\R^2$
\vspace{0.35ex}
of the Jacobi equation (12.45). 
\vspace{0.5ex}
So the straight line $x+y+1=0$ is a trajectory of the Jacobi equation (12.45).

However the linear trajectory $x+y+1=0$ 
\vspace{0.5ex}
is not contained in the the pencil of the linear trajectories
$\dfrac{y-1}{x-1}=C,\ \ C\in [{}-\infty;{}+\infty].$
\vspace{0.75ex}
It follows from the fact that the Jacobi equation (12.45) we can reduce to the form
\\[1.5ex]
\mbox{}\hfill
$
\dfrac{dy}{dx}=\dfrac{(1-y)(x+y+1)}{(1-x)(x+y+1)}\,.
\hfill
$
\\[1.5ex]
And from this representation, we see that every point of the linear trajectory 
$x+y+1=0$ is singular.
\vspace{0.25ex}

Thus, for the Jacobi equation (12.45), 
\vspace{0.25ex}
along with the general integral (12.46) we must specified and its partial integral 
$p_3^{}\colon (x,y)\to x+y+1$ for all $(x,y)\in\R^2.$
\vspace{1ex}

{\small\bf 12.22.1.3. Case of complex eigenvalue.}
\vspace{0.5ex}
If  the complex number  
$\lambda_1^{}=\xi_1^{}+\zeta_1^{}i\ (\zeta_1^{}\ne 0)$ 
is a root of equation (12.38), 
then the complex conjugate number
\vspace{0.35ex}
$\overline{\lambda}_1^{}=\xi_1^{}-\zeta_1^{}i$ 
is also a root of this equation. 
Moreover, the equation (12.38) has the real root $\lambda_3^{}.$
\vspace{0.5ex}

If the linear system (12.37) under $\lambda=\lambda_1^{}$ has the solution  
\vspace{0.5ex}
$\alpha=\alpha_1^{},\ \beta=\beta_1^{},\ \gamma=\gamma_1^{},$
then  the linear system (12.37) under $\lambda=\overline{\lambda}_1^{}$ 
has the solution  
\vspace{0.55ex}
$\alpha=\overline{\alpha}_1^{},\ \beta=\overline{\beta}_1^{},\ 
\gamma=\overline{\gamma}_1^{}.$

The solutions $\lambda_1^{}$ and $\overline{\lambda}_1^{}$ 
of equation (12.38)
\vspace{0.5ex}
are eigenvalues of matrix $A.$
Furthermore, the solutions 
$(\alpha_1^{},\beta_1^{},\gamma_1^{})$ and
$(\overline{\alpha}_1^{},\overline{\beta}_1^{},\overline{\gamma}_1^{})$ 
\vspace{0.75ex}
of system (12.37) are eigenvectors of matrix $A$ corresponding to 
the eigenvalues $\lambda_1^{}$ and $\overline{\lambda}_1^{},$ respectively.
\vspace{0.55ex}

Therefore if the matrix $A$ has the complex eigenvalue
$\lambda_1^{}=\xi_1^{}+\zeta_1^{}i\ (\zeta_1^{}\ne 0)$
\vspace{0.75ex}
with the corresponding eigenvector $(\alpha_1^{},\beta_1^{},\gamma_1^{}),$
then the vector 
$(\overline{\alpha}_1^{},\overline{\beta}_1^{},\overline{\gamma}_1^{})$ 
\vspace{0.75ex}
is an eigenvector of matrix $A$ corresponding to the eigenvalue
$\overline{\lambda}_1^{}=\xi_1^{}-\zeta_1^{}i.$
\vspace{0.75ex}

Furthermore, the matrix $A$ has the real eigenvalue $\lambda_3^{}.$
The eigenvalues $\lambda_1^{},\ \overline{\lambda}_1^{},\ \lambda_3^{}$
\vspace{0.5ex}
of matrix $A$ have the simple elementary divisors
$\lambda-\lambda_1^{},\ \lambda-\overline{\lambda}_1^{},\ \lambda-\lambda_3^{}.$
\vspace{1.25ex}

{\bf Theorem 12.11.}
\vspace{0.75ex}
{\it 
Suppose $(\alpha,\beta,\gamma)$ is an eigenvector of matrix $A$ 
corresponding to the complex eigenvalue $\lambda=\xi+\zeta\;\!i\ (\zeta\ne 0)$
and $(\alpha_{3}^{},\beta_{3}^{},\gamma_{3}^{})$
\vspace{0.75ex}
is an eigenvector of matrix $A$ 
cor\-res\-pon\-ding to the real eigenvalue $\lambda_3^{}.$
Then the Jacobi equation {\rm (12.36)} has the general integral
\\[1.5ex]
\mbox{}\hfill    % (12.47)
$
F\colon (x,y)\to\ 
H^{{}^{\scriptstyle \zeta}}\!(x,y)\exp S(x,y)
$
\ for all $(x,y)\in D,
$
\hfill {\rm (12.47)}
\\[2.25ex]
where the functions 
\\[1.5ex]
\mbox{}\hfill
$
H(x,y)=\dfrac{{\rm Re}^2p(x,y)+{\rm Im}^2p(x,y)}{p_3^2(x,y)}
$
\ for all 
$(x,y)\in D,
\hfill
$
\\[2ex]
\mbox{}\hfill
$
S(x,y)=2\;\!(\lambda_3^{}-\xi)\arctan\dfrac{{\rm Im}\;\!p(x,y)}{{\rm Re}\;\!p(x,y)}
$
\ for all 
$
(x,y)\in D,\
\hfill
$
\\[1.75ex]
\mbox{}\hfill
$
p(x,y)=\alpha x+\beta y+\gamma,
\quad
p_3^{}(x,y)=\alpha_3^{} x+\beta_3^{} y+\gamma_3^{}
$
\ for all 
$
(x,y)\in\R^2,
\hfill
$
\\[2ex]
the set $D\subset\R^2$ such that   
$H^{{}^{\scriptstyle \zeta}}$ and 
$\arctan\dfrac{{\rm Im}\;\!p}{{\rm Re}\;\!p}$
are continuously differentiable on $D.$}
\vspace{1ex}

{\sl Proof}.
From Lemma  12.1 it follows that the linear function $p_3^{}$
\vspace{0.5ex}
is a real polynomial partial integral with cofactor
$M_3^{}\colon (x,y)\to \lambda_3^{}-l_3^{}(x,y)$ for all $(x,y)\in\R^2$
\vspace{0.5ex}
of  the Jacobi equation (12.36) and
the complex linear function $p$ is a complex-valued 
\vspace{0.5ex}
polynomial partial integral with cofactor
$
M\colon (x,y)\to \lambda-l_3^{}(x,y)
$ 
for all 
$(x,y)\in\R^2
$
\vspace{1.25ex}
of  the Jacobi equation (12.36).

Suppose $p(x,y)=u(x,y)+i\;\!v(x,y)$ for all $(x,y)\in\R^2,$
where the functions
\\[1.75ex]
\mbox{}\hfill
$
u(x,y)={\rm Re}\;\!p(x,y)$ 
\ for all $(x,y)\in\R^2,
\quad 
v(x,y)={\rm Im}\;\!p(x,y)$ \ 
for all $(x,y)\in\R^2. 
\hfill
$
\\[1.75ex]
\indent
By Theorem 6.2, the function $p$ is 
\vspace{0.75ex}
a complex-valued polynomial partial integral with cofactor
$M\colon (x,y)\to\xi-l_3^{}(x,y)+\zeta i$ for all $(x,y)\in\R^2$
\vspace{0.75ex}
if and only if the function $u^2+v^2$ 
is a polynomial partial integral with cofactor
\vspace{0.75ex}
$M_1^{}\colon (x,y)\to 2(\xi-l_3^{}(x,y))$ for all $(x,y)\in\R^2$
and the function
$\exp\arctan\dfrac{v}{u}$ 
\vspace{0.75ex}
is an exponential partial integral with cofactor
$M_2^{}\colon (x,y)\to \zeta$ for all $(x,y)\in\R^2$
on the set $D_0^{}=\{(x,y)\colon u(x,y)\ne 0\}.$
\vspace{1ex}

The cofactors $M_1^{},\ M_2^{},\ M_3^{}$ of the partial integrals
$u^2+v^2,\ \exp\arctan\dfrac{v}{u}\,,\ p_3^{}$
such that the linear combination
\\[1.5ex]
\mbox{}\hfill
$
\zeta\;\!M_1^{}(x,y)+2\;\!(\lambda_3^{}-\xi)\;\!M_2^{}(x,y)-2\zeta\;\!M_3^{}(x,y)=0
$
\ for all $(x,y)\in\R^2.
\hfill
$
\\[1.75ex]
\indent
Thus (Property 11.11)
a general integral on the set $D$ of the Jacobi equation (12.36) is the function
\\[2ex]
\mbox{}\hfill
$
F\colon (x,y)\to 
\bigl(u^2(x,y)+v^2(x,y)\bigr)^{\zeta}\;\!
\biggl(\exp\arctan\dfrac{v(x,y)}{u(x,y)}\biggr)^{\!2\;\!(\lambda_3^{}-\xi)}
p_3^{{}-2\zeta}(x,y)
$
for all $(x,y)\in D.
\k
\hfill
$
\\[2.5ex]
\indent
Note also that the statement of Theorem 12.11 we can obtain from Theorem 12.8 under
\\[1.5ex]
\mbox{}\hfill
$
p_2^{}=\overline{p}_1^{},
\ \ 
\lambda_1^{}=\xi+\zeta\;\! i, 
\ \
\lambda_2^{}=\xi-\zeta\;\! i=\overline{\lambda}_1^{}, 
\ \
h_1^{}=(\overline{\lambda}_1^{}-\lambda_3^{})\;\! i=\zeta+(\xi-\lambda_3^{})\;\! i, 
\hfill
$
\\[2ex]
\mbox{}\hfill
$
h_2^{}=(\lambda_3^{}-\lambda_1^{})\;\! i=\zeta-(\xi-\lambda_3^{})\;\! i=
\overline{h}_1^{}, 
\quad
h_3^{}=(\lambda_1^{}-\overline{\lambda}_1^{})\;\! i={}-2\;\!\zeta.
\hfill
$
\\[2.25ex]
\indent
{\bf Example 12.3.}
Consider the Jacobi equation [9, pp. 24 -- 25]
\\[2ex]
\mbox{}\hfill  % (12.48)
$
\dfrac{dy}{dx}=\dfrac{6x-2y+1-y\;\!({}-2x+y)}{4x-3y-1-x\;\!({}-2x+y)}\,.
$
\hfill (12.48)
\\[2.25ex]
\indent
The matrix
\vspace{0.75ex}
$
A=\left(\!\!\!
\begin{array}{rrr}
4 & 6 & {}-2
\\
{}-3 & {}-2 & 1
\\
{}-1 & 1 & 0\end{array}
\!\!\right)
$
of the Jacobi equation (12.48) has the complex eigenvalue
$\lambda=1+\sqrt{6}\;\!i$ and the real eigenvalue $\lambda_3^{}=0.$
\vspace{1ex}

The vector $\bigl(-\;\!2\;\!(2+\sqrt{6}\;\!i),\;\! 5,\;\! 3-\sqrt{6}\;\!i\bigr)$
\vspace{0.5ex}
is an eigenvector corresponding to the complex eigenvalue
$\lambda=1+\sqrt{6}\;\!i$
and the vector $(1, 1, 5)$ is a real 
\vspace{0.35ex}
eigenvector corresponding to the eigenvalue $\lambda_3^{}=0.$
\vspace{0.5ex}

By Theorem 12.11, 
\vspace{0.5ex}
a general integral of the Jacobi equation (12.48) on the set
$D=
\linebreak
=\{(x,y)\colon (x+y+5)(4x-5y-3)\ne 0\}$ is the function
\\[2ex]
\mbox{}\hfill  % (12.49)
$
F\colon (x,y)\to\ 
\biggl(\dfrac{(4x-5y-3)^2+6(2x+1)^2}{(x+y+5)^2}\biggr)^{\!\sqrt{6}}
\exp\biggl({}-2\;\!\arctan\dfrac{\sqrt{6}\;\!(2x+1)}{4x-5y-3}\biggr).
%\quad
%\forall (x,y)\in D.
$
\hfill {\rm (12.49)}
\\[2.5ex]
\indent
{\bf Example 12.4.}
Consider the Jacobi equation 
\\[2ex]
\mbox{}\hfill  % (12.50)
$
\dfrac{dy}{dx}=\dfrac{x+y-3-y\;\!(x+2y+1)}{x+2y+3-x\;\!(x+2y+1)}\,.
$
\hfill (12.50)
\\[2ex]
\indent
The matrix
\vspace{0.75ex}
$
A=\left(\!\!\!
\begin{array}{rrr}
1 & 1 & 1
\\
2 & 1 & 2
\\
3 & {}-3 & 1\end{array}
\!\!\right)
$
has the eigenvalues $\lambda=1+i, \ \lambda_3^{}=1$ 
and the corresponding eigenvectors
$(3-i,\;\! 4, {}-3+3\;\!i),\ (1, 1, {}-1).$
\vspace{0.75ex}
By Theorem 12.11 (under $\xi=\lambda_3^{}),$ 
a general integral on the set 
\vspace{0.5ex}
$D=\{(x,y)\colon x+y-1\ne 0\}$ 
of the Jacobi equation (12.50) is the rational function
\\[1.5ex]
\mbox{}\hfill  % (12.51)
$
F\colon (x,y)\to\ 
\dfrac{(3x+4y-3)^2+(x-3)^2}{(x+y-1)^2}
$
\ for all $(x,y)\in D.
$
\hfill {\rm (12.51)}
\\[2ex]
\indent
{\small\bf 
12.22.2. Case of multiple elementary devisor.}
\vspace{0.35ex}

Suppose $\lambda_1^{}$ is an eigenvalue of matrix $A$ 
\vspace{0.35ex}
with multiple elementary devisor.
Then this eigenvalue has the eigenvector
$\theta_1^{}=(\alpha_1^{},\beta_1^{},\gamma_1^{})$ and also
\vspace{0.35ex}
the first generalized iegenvector 
$\theta_1^{(1)}=\bigl(\alpha_1^{(1)},\beta_1^{(1)},\gamma_1^{(1)}\bigr).$
\vspace{0.35ex}

The first generalized eigenvector $\theta_1^{(1)}$
is a solution of the matrix equation
\\[1.5ex]
\mbox{}\hfill % (12.52)
$
(A-\lambda_1^{}E)\;\! {\rm colon}\,\theta_1^{(1)} =
{\rm colon}\;\!\theta_1^{}.
$
\hfill (12.52)
\\[2.25ex]
\indent
{\bf Lemma 12.2.}
\vspace{0.35ex}
{\it 
Let $\lambda_1^{}$ be the eigenvalue of matrix $A$
 with multiple elementary devisor correspon\-ding to the eigenvector 
\vspace{0.5ex}
$\theta_1^{}=(\alpha_1^{},\beta_1^{},\gamma_1^{})$ and to 
the first generalized eigenvector 
\vspace{0.75ex}
$\theta_1^{(1)}=\bigl(\alpha_1^{(1)},\beta_1^{(1)},\gamma_1^{(1)}\bigr).\!$
Then the Jacobi equation {\rm (12.36)} has 
the multiple polynomial partial in\-teg\-ral
$p_1^{}\colon (x,y)\to\alpha_1^{}x+\beta_1^{}y+\gamma_1^{}$
for all $(x,y)\in\R^2$
\vspace{1ex}
with cofactor
\linebreak
$M_1^{}\colon (x,y)\to\! \lambda_1^{}\!-l_3^{}(x,y)$ for all $(x,y)\in\R^2$
\vspace{1ex}
on the set 
$D=\{(x,y)\colon \alpha_1^{}x+\beta_1^{}y+\gamma_1^{}\ne 0\}$
such that
$\bigl((p_1^{}, M_1^{}), (1,p_1^{(1)}, 1)\bigr)\in \text{\rm B}_{_{\!D}},$
where 
\vspace{1ex}
$p_1^{(1)}(x,y)= \alpha_1^{(1)}x+\beta_1^{(1)}y+\gamma_1^{(1)}$ 
for all $(x,y)\in\R^2.$
}

{\sl Proof.}
The matrix equation (12.52) is equivalent to the linear system
\\[1.5ex]
\mbox{}\hfill              % (12.53)               
$\!\!\!
\begin{array}{c}
(a_1^{}\!-\!\lambda_1^{})\;\!\alpha_{1}^{(1)}+
a_2^{}\;\!\beta_{1}^{(1)}+
a_3^{}\;\!\gamma_{1}^{(1)}=\alpha_1^{},
\\[1.5ex]
b_1^{}\;\!\alpha_{1}^{(1)}+
(b_2^{}\!-\!\lambda_1^{})\;\!\beta_{1}^{(1)}+
b_3^{}\;\!\gamma_{1}^{(1)}=\beta_1^{},
\\[1.5ex]
c_1^{}\;\!\alpha_{1}^{(1)}+
c_2^{}\;\!\beta_{1}^{(1)}+
(c_3^{}\!-\!\lambda_1^{})\;\!\gamma_{1}^{(1)}=\gamma_1^{}
\end{array}
%\Leftrightarrow
\!\!\!\iff\!\!\!
\begin{array}{c}
a_1^{}\;\!\alpha_{1}^{(1)}+
a_2^{}\;\!\beta_{1}^{(1)}=
\alpha_1^{}+\lambda_1^{}\;\!\alpha_{1}^{(1)}-a_3^{}\;\!\gamma_{1}^{(1)},
\\[1.5ex]
b_1^{}\;\!\alpha_{1}^{(1)}+
b_2^{}\;\!\beta_{1}^{(1)}= 
\beta_1^{}+\lambda_1^{}\;\!\beta_{1}^{(1)}-b_3^{}\;\!\gamma_{1}^{(1)},
\\[1.5ex]
c_1^{}\;\!\alpha_{1}^{(1)}+
c_2^{}\;\!\beta_{1}^{(1)}=
\gamma_1^{}+\lambda_1^{}\;\!\gamma_{1}^{(1)}-c_3^{}\;\!\gamma_{1}^{(1)}.
\end{array}\!\!\!
$
\hfill (12.53)
\\[2ex]
\indent
The derivative by virtue of the Jacobi equation (12.36)
\\[2ex]
\mbox{}\hfill
$
{\frak I}\;\!p_1^{(1)}(x,y)=
{\frak I}\;\!\bigl(\alpha_1^{(1)}x+\beta_1^{(1)}y+\gamma_1^{(1)}\bigr)=
\alpha_{1}^{(1)}\;\!\bigl(l_1^{}(x,y)-x\;\!l_3^{}(x,y)\bigr)+
\beta_{1}^{(1)}\;\!\bigl(l_2^{}(x,y)-y\;\!l_3^{}(x,y)\bigr)=
\hfill
$
\\[2ex]
\mbox{}\hfill
$
=
\alpha_{1}^{(1)}\bigl(a_1^{}x+b_1^{}y+c_1^{}\bigr)+
\beta_{1}^{(1)}\bigl(a_2^{}x+b_2^{}y+c_2^{}\bigr)-
\bigl(\alpha_{1}^{(1)}x+\beta_{1}^{(1)}y\bigr)\;\!l_3^{}(x,y)=
\bigl(a_1^{}\alpha_{1}^{(1)}+a_2^{}\beta_{1}^{(1)}\bigr) x\ +
\hfill
$
\\[2ex]
\mbox{}\hfill
$
+\
\bigl(b_1^{}\alpha_{1}^{(1)}+b_2^{}\beta_{1}^{(1)}\bigr) y+
\bigl(c_1^{}\alpha_{1}^{(1)}+c_2^{}\beta_{1}^{(1)}\bigr)-
\bigl(\alpha_{1}^{(1)}x+\beta_{1}^{(1)}y\bigr)\;\!l_3^{}(x,y)
$
\ for all $(x,y)\in\R^2.
\hfill
$
\\[2ex]
\indent
Now, taking into account the system (12.53), we obtain
\\[2ex]
\mbox{}\hfill
$
{\frak I}\;\!p_1^{(1)}(x,y)=
\bigl(\alpha_1^{}+\lambda_1^{}\;\!\alpha_{1}^{(1)}-a_3^{}\;\!\gamma_{1}^{(1)}\bigr)\;\! x+
\bigl(\beta_1^{}+\lambda_1^{}\;\!\beta_{1}^{(1)}-b_3^{}\;\!\gamma_{1}^{(1)}\bigr)\;\! y+
\bigl(\gamma_1^{}+\lambda_1^{}\;\!\gamma_{1}^{(1)}-c_3^{}\;\!\gamma_{1}^{(1)}\bigr)\ -
\hfill
$
\\[2ex]
\mbox{}\hfill
$
-\ \bigl(\alpha_{1}^{(1)}x+\beta_{1}^{(1)}y\bigr)\;\!l_3^{}(x,y)=
\alpha_1^{}x+\beta_1^{}y+\gamma_1^{}
+\lambda_1^{}\;\! \bigl(\alpha_1^{(1)}x+\beta_1^{(1)}y+\gamma_1^{(1)}\bigr)-
\gamma_{1}^{(1)}(a_3^{}x+b_3^{}y+c_3^{})\ -
\hfill
$
\\[2ex]
\mbox{}\hfill
$
-\ \bigl(\alpha_{1}^{(1)}x+\beta_{1}^{(1)}y\bigr)\;\!l_3^{}(x,y)=
\alpha_1^{}x+\beta_1^{}y+\gamma_1^{}
+\lambda_1^{}\;\! \bigl(\alpha_1^{(1)}x+\beta_1^{(1)}y+\gamma_1^{(1)}\bigr)\ -
\hfill
$
\\[2ex]
\mbox{}\hfill
$
-\ 
\bigl(\alpha_{1}^{(1)}x+\beta_{1}^{(1)}y+\gamma_{1}^{(1)}\bigr)\;\!l_3^{}(x,y)=
\alpha_1^{}x+\beta_1^{}y+\gamma_1^{}+
\bigl(\alpha_{1}^{(1)}x+\beta_{1}^{(1)}y+\gamma_{1}^{(1)}\bigr)
\bigl(\lambda_1^{}-l_3^{}(x,y)\bigr)=
\hfill
$
\\[2ex]
\mbox{}\hfill
$
=p_1^{}(x,y)+p_1^{(1)}(x,y)\;\!M_1^{}(x,y)
$ 
\ for all $(x,y)\in\R^2.
\hfill
$
\\[2ex]
\indent
Hence the derivative by virtue of the Jacobi equation (12.36)
\\[2ex]
\mbox{}\hfill  % (12.54)
$
{\frak I}\;\!p_1^{(1)}(x,y)=
p_1^{}(x,y)+p_1^{(1)}(x,y)\;\!M_1^{}(x,y)
$ 
\ for all $(x,y)\in\R^2.
$
\hfill (12.54)
\\[2ex]
\indent
Finally, from the criterion of multiple polynomial partial integral
\vspace{0.35ex}
(Theorem 5.2 under $h=1,\, N=1,$ $p=p_1^{},\ q=p_1^{(1)},\ M=M_1^{})$
\vspace{0.35ex}
it follows that 
the polynomial partial integral $p_1^{}$ 
of the Jacobi equation (12.36) such that 
$\bigl((p_1^{}, M_1^{}), (1,p_1^{(1)}, 1)\bigr)\in \text{B}_{_{\!D}}.\ \k$
\vspace{1.25ex}

{\small\bf 12.22.2.1. Case of double elementary devisor.}
\vspace{0.5ex}
If the multiple eigenvalue $\lambda_1^{}$ of matrix $A$
has the double elementary devisor $(\lambda-\lambda_1^{})^2,$
\vspace{0.5ex}
then the matrix $A$ has also the real eigenvalue  $\lambda_3^{}$
with the simple elementary devisor $\lambda-\lambda_3^{}.$
\vspace{0.75ex}

{\bf Theorem 12.12.}
\vspace{0.35ex}
{\it
Let $\lambda_1^{}$ be the eigenvalue of matrix $A$ 
with double elementary devisor corresponding to the eigenvector
\vspace{0.35ex}
$(\alpha_1^{},\beta_1^{},\gamma_1^{})$ and to the first order generalized 
eigenvector
$\bigl(\alpha_1^{(1)},\beta_1^{(1)},\gamma_1^{(1)}\bigr),$
\vspace{0.55ex}
let   $\lambda_3^{}$ be the eigenvalue of matrix $A$ corresponding to the 
eigenvector 
$(\alpha_3^{},\beta_3^{},\gamma_3^{}).$
\vspace{0.55ex}
Then a general integral on  
$D=\{(x,y)\colon (\alpha_1^{}x+\beta_1^{}y+\gamma_1^{})
(\alpha_3^{}x+\beta_3^{}y+\gamma_3^{})\ne 0\}$
of the Jacobi equation {\rm (12.36)} is the function
\\[1.5ex]
\mbox{}\hfill   % (12.55)
$
F\colon (x,y)\to \
\dfrac{p_1^{}(x,y)}{p_3^{}(x,y)}\,
\exp\biggl((\lambda_3^{}-\lambda_1^{})\,\dfrac{p_1^{(1)}(x,y)}{p_1^{}(x,y)}\biggr)
$
\ for all $(x,y)\in D,
$
\hfill {\rm (12.55)}
\\[2ex]
where 
\vspace{0.75ex}
$
p_1^{}(x,y)= \alpha_1^{}x+\beta_1^{}y+\gamma_1^{}$ for all $(x,y)\in\R^2,\ 
p_1^{(1)}(x,y)= \alpha_1^{(1)}x+\beta_1^{(1)}y+\gamma_1^{(1)}$ 
for all $(x,y)\in\R^2,$ 
$
p_3^{}(x,y)= \alpha_3^{}x+\beta_3^{}y+\gamma_3^{}$ for all $(x,y)\in\R^2.$
}
\vspace{0.75ex}

{\sl Proof.}
\vspace{0.75ex}
Taking into account Lemma 12.2, we have
$\bigl((p_1^{}, \lambda_1^{}-l_3^{}), (1,p_1^{(1)}, 1)\bigr)\in 
\text{B}_{_{\!\scriptstyle D_{0^{}}}},$
where the set
\vspace{0.35ex}
$D_0^{}=\{(x,y)\colon \alpha_1^{}x+\beta_1^{}y+\gamma_1^{}\ne 0\}.$
From the criterion of multiple polynomial partial integral 
\vspace{0.5ex}
(Theorem 5.3 under $p=p_1^{},\ M=\lambda_1^{}-l_3^{},\ h=1,\ q=p_1^{(1)},\ N=1),$
we get the Jacobi equation {\rm (12.36)} has
\vspace{0.5ex}
the polynomial partial integral $p_1^{}$
with cofactor $M_1^{}\colon (x,y)\to \lambda_1^{}-l_3^{}(x,y)$ for all $(x,y)\in\R^2$
\vspace{0.5ex}
and also the exponential partial integral
$\exp\dfrac{p_1^{(1)}}{p_1^{}}$ 
with cofactor 
\vspace{1ex}
$M_2^{}\colon (x,y)\to 1$ for all $(x,y)\in\R^2.$

By Lemma 12.1, the linear function $p_3^{}$ 
\vspace{1.25ex}
is a polynomial partial integral with cofactor
$M_3^{}\colon (x,y)\to \lambda_3^{}-l_3^{}(x,y)$ for all $(x,y)\in\R^2$
\vspace{0.75ex}
of the Jacobi equation (12.36).

The cofactors $M_1^{},\ M_2^{},\ M_3^{}$ of the partial integrals
$p_1^{},\ \exp\dfrac{p_1^{(1)}}{p_1^{}}\,,\ p_3^{}$
such that the linear combination
\\[1.25ex]
\mbox{}\hfill
$
M_1^{}(x,y)+(\lambda_3^{}-\lambda_1^{})\;\!M_2^{}(x,y)-M_3^{}(x,y)=0
$
\ for all $(x,y)\in\R^2.
\hfill
$
\\[1.5ex]
\indent
From Property 11.11 it follows that 
\vspace{0.15ex}
the function (12.55) is a general integral on the set $D$ of 
the Jacobi equation (12.36). $\k$
\vspace{0.75ex}

{\bf Example 12.5.}
Consider the Jacobi equation  [9, p. 25]
\\[2ex]
\mbox{}\hfill  % (12.56)
$
\dfrac{dy}{dx}=\dfrac{x-y-1-y\;\!({}-x+y)}{{}-x+y-x\;\!({}-x+y)}\,.
$
\hfill (12.56)
\\[2ex]
\indent
The matrix 
\vspace{0.75ex}
$
A=\left(\!\!\!
\begin{array}{rrr}
{}-1 & 1 & {}-1
\\
1 & {}-1 & 1
\\
0 & {}-1 & 0\end{array}
\!\!\right)
$
of this system has the eigenvalue $\lambda_1^{}={}-1$
with  double elementary devisor $(\lambda+1)^2,$
\vspace{0.35ex}
the eigenvector $(1,{}-1, {}-1),$ and
the first order generalized eigenvector $(1,{}-1, {}-2).$
\vspace{0.35ex}
Also, the matrix $A$ has the eigenvalue $\lambda_3^{}=0$ 
with the corresponding eigenvector $(1, 0, {}-1).$
\vspace{0.75ex}
By Theorem 12.12, a general integral on the set  
$D=\{(x,y)\colon (x-1)(x-y-1)\ne 0\}$ of the Jacobi equation (12.56) 
is the function
\\[1.75ex]
\mbox{}\hfill  % (12.57)
$
F\colon (x,y)\to\ 
\dfrac{x-y-1}{x-1}\,\exp\dfrac{x-y-2}{x-y-1}
$
\ for all $(x,y)\in D.
$
\hfill {\rm (12.57)}
\\[2ex]
\indent
{\bf Example 12.6.}
The Jacobi equation
\\[2ex]
\mbox{}\hfill  % (12.58)
$
\dfrac{dy}{dx}=\dfrac{y-y\;\!(y+1)}{x-x\;\!(y+1)}
$
\hfill (12.58)
\\[2ex]
such that the matrix
\vspace{0.35ex}
$
A=\left(\!\!\!
\begin{array}{ccc}
1 & 0 & 0
\\
0 & 1 & 1
\\
0 & 0 & 1\end{array}
\!\!\right)
$
has the eigenvalues $\lambda_1^{}=\lambda_2^{}=\lambda_3^{}=1$
with the elementary devisors $(\lambda-1)^2$ and $(\lambda-1).$
\vspace{0.35ex}
The eigenvalue $\lambda_1^{}=1$ with the elementary devisor $(\lambda-1)^2$
\vspace{0.35ex}
has the eigenvector $(0, 1, 0)$ and the first order generalized eigenvector $(1, 1, 1).$
The eigenvalue $\lambda_3^{}=1$ with the elementary devisor $\lambda-1$
\vspace{0.5ex}
has the eigenvector $(1, 0, 0).$
By Theorem 12.12 (under $\lambda_3^{}=\lambda_1^{}=1),$ 
\vspace{0.35ex}
a general integral on the set  
$D=\{(x,y)\colon x\ne 0\}$ of the Jacobi equation (12.58)
is the homographic function
\\[1.5ex]
\mbox{}\hfill  % (12.59)
$
F\colon (x,y)\to\ 
\dfrac{y}{x}
$
\ for all $(x,y)\in D.
$
\hfill {\rm (12.59)}
\\[2ex]
\indent
{\small\bf 
12.22.2.2. Case of triple elementary devisor.}
\vspace{0.35ex}
If the multiple eigenvalue $\lambda_1^{}$ of the matrix $A$ has 
the triple elementary devisor
$(\lambda-\lambda_1^{})^3,$ then the eigenvalue $\lambda_1^{}$ has
\vspace{0.5ex}
the eigenvector $\theta_1^{}=(\alpha_1^{},\beta_1^{},\gamma_1^{}),$
the first order generalized eigenvector 
\vspace{0.5ex}
$\theta_1^{(1)}=\bigl(\alpha_1^{(1)},\beta_1^{(1)},\gamma_1^{(1)}\bigr),$ 
and the second order generalized eigenvector 
$\theta_1^{(2)}=\bigl(\alpha_1^{(2)},\beta_1^{(2)},\gamma_1^{(2)}\bigr).$
\vspace{0.5ex}
The first order generalized eigenvector $\theta_1^{(1)}$ is 
\vspace{0.5ex}
a solution of the matrix equation (12.52). 
The second order generalized eigenvector $\theta_1^{(2)}$ is 
a solution of the matrix equation
\\[1.5ex]
\mbox{}\hfill % (12.60)
$
(A-\lambda_1^{}E)\;\! {\rm colon}\,\theta_1^{(2)} =
2\, {\rm colon}\;\!\theta_1^{(1)}.
$
\hfill (12.60)
\\[2ex]
\indent
{\bf Theorem 12.13.}
\vspace{0.25ex}
{\it
Suppose $\lambda_1^{}$ is the eigenvalue of the matrix $A$ corresponding to 
the triple elementary devisor, the eigenvector
\vspace{0.5ex}
$(\alpha_1^{},\beta_1^{},\gamma_1^{}),$ 
the first order generalized eigenvector 
\vspace{0.5ex}
$\bigl(\alpha_1^{(1)},\beta_1^{(1)},\gamma_1^{(1)}\bigr),$ 
and the second order generalized eigenvector  
$\bigl(\alpha_1^{(2)},\beta_1^{(2)},\gamma_1^{(2)}\bigr).$
Then a general integral on the set 
\vspace{0.35ex}
$D=\{(x,y)\colon \alpha_1^{}x+\beta_1^{}y+\gamma_1^{}\ne 0\}$
of the Jacobi equation {\rm (12.36)} is the rational function
\\[2ex]
\mbox{}\hfill   % (12.61)
$
F\colon (x,y)\to \
\dfrac{\bigl(p_1^{(1)}(x,y)\bigr)^2-p_1^{}(x,y)\;\!p_1^{(2)}(x,y)}{p_1^{2}(x,y)}
$
\ for all $(x,y)\in D,
$
\hfill {\rm (12.61)}
\\[2ex]
where 
\vspace{0.75ex}
$
p_1^{}(x,y)= \alpha_1^{}x+\beta_1^{}y+\gamma_1^{}$ for all $(x,y)\in\R^2,\ 
p_1^{(1)}(x,y)= \alpha_1^{(1)}x+\beta_1^{(1)}y+\gamma_1^{(1)}$ 
for all $(x,y)\in\R^2,
$ 
$
p_1^{(2)}(x,y)= \alpha_1^{(2)}x+\beta_1^{(2)}y+\gamma_1^{(2)}$ 
for all $(x,y)\in\R^2. 
$
}
\vspace{0.75ex}

{\sl Proof.}
\vspace{0.75ex}
Taking into account Lemma 12.2, we get 
$\bigl((p_1^{}, M_1^{}), (1,p_1^{(1)}, 1)\bigr)\in \text{B}_{_{\!\scriptstyle D}},$
where the cofactor  
$M_1^{}\colon (x,y)\to \lambda_1^{}-l_3^{}(x,y)$ for all $(x,y)\in\R^2.$
\vspace{0.35ex}
From the criterion of multiple polynomial partial integral  
\vspace{0.35ex}
(Theorem 5.1 under $p=p_1^{},\ M=M_1^{},\ h=1,\ q=p_1^{(1)},$ $N=1),$
we obtain the following identities
\\[1.5ex]
\mbox{}\hfill  % (12.62)
$
{\frak I}\;\!p_1^{}(x,y)=p_1^{}(x,y)\;\!M_1^{}(x,y)
$
\ for all $(x,y)\in\R^2
$
\hfill (12.62)
\\[1ex]
and 
\\[1ex]
\mbox{}\hfill  % (12.63)
$
{\frak I}\,\dfrac{p_1^{(1)}(x,y)}{p_1^{}(x,y)}=1
$
\ for all $(x,y)\in D.
$
\hfill (12.63)
\\[2ex]
\indent
The matrix equation (12.60) is equivalent to the linear system
\\[1.75ex]
\mbox{}\hfill              % (12.64)               
$\!\!\!\!\!
\begin{array}{c}
(a_1^{}\!-\!\lambda_1^{})\;\!\alpha_{1}^{(2)}\!+
a_2^{}\;\!\beta_{1}^{(2)}\!+
a_3^{}\;\!\gamma_{1}^{(2)}=2\alpha_1^{(1)},
\\[2ex]
b_1^{}\;\!\alpha_{1}^{(2)}\!+
(b_2^{}\!-\!\lambda_1^{})\;\!\beta_{1}^{(2)}\!+
b_3^{}\;\!\gamma_{1}^{(2)}=2\beta_1^{(1)},
\\[2ex]
c_1^{}\;\!\alpha_{1}^{(2)}\!+
c_2^{}\;\!\beta_{1}^{(2)}\!+
(c_3^{}\!-\!\lambda_1^{})\;\!\gamma_{1}^{(2)}=2\gamma_1^{(1)}
\end{array}
\!\!\Leftrightarrow\!\!
%\!\!\!\!\iff\!\!\!\!
\begin{array}{c}
a_1^{}\;\!\alpha_{1}^{(2)}\!+
a_2^{}\;\!\beta_{1}^{(2)}=
2\alpha_1^{(1)}\!+\lambda_1^{}\;\!\alpha_{1}^{(2)}\!-a_3^{}\;\!\gamma_{1}^{(2)},
\\[2ex]
b_1^{}\;\!\alpha_{1}^{(2)}\!+
b_2^{}\;\!\beta_{1}^{(2)}= 
2\beta_1^{(1)}\!+\lambda_1^{}\;\!\beta_{1}^{(2)}\!-b_3^{}\;\!\gamma_{1}^{(2)},
\\[2ex]
c_1^{}\;\!\alpha_{1}^{(2)}\!+
c_2^{}\;\!\beta_{1}^{(2)}=
2\gamma_1^{(1)}\!+\lambda_1^{}\;\!\gamma_{1}^{(2)}\!-c_3^{}\;\!\gamma_{1}^{(2)}.
\end{array}\!\!\!
$
\hfill (12.64)
\\[2ex]
\indent
The derivative by virtue of the Jacobi equation (12.36) on $\R^2$ is
\\[2ex]
\mbox{}\hfill
$
{\frak I}\;\!p_1^{(2)}(x,y)=
{\frak I}\;\!\bigl(\alpha_1^{(2)}x+\beta_1^{(2)}y+\gamma_1^{(2)}\bigr)=
\alpha_{1}^{(2)}\;\!\bigl(l_1^{}(x,y)-x\;\!l_3^{}(x,y)\bigr)+
\beta_{1}^{(2)}\;\!\bigl(l_2^{}(x,y)-y\;\!l_3^{}(x,y)\bigr)=
\hfill
$
\\[2ex]
\mbox{}\hfill
$
=
\alpha_{1}^{(2)}\bigl(a_1^{}x+b_1^{}y+c_1^{}\bigr)+
\beta_{1}^{(2)}\bigl(a_2^{}x+b_2^{}y+c_2^{}\bigr)-
\bigl(\alpha_{1}^{(2)}x+\beta_{1}^{(2)}y\bigr)\;\!l_3^{}(x,y)=
\hfill
$
\\[2ex]
\mbox{}\hfill
$
=
\bigl(a_1^{}\alpha_{1}^{(2)}+a_2^{}\beta_{1}^{(2)}\bigr) x+
\bigl(b_1^{}\alpha_{1}^{(2)}+b_2^{}\beta_{1}^{(2)}\bigr) y+
\bigl(c_1^{}\alpha_{1}^{(2)}+c_2^{}\beta_{1}^{(2)}\bigr)-
\bigl(\alpha_{1}^{(2)}x+\beta_{1}^{(2)}y\bigr)\;\!l_3^{}(x,y).
\hfill
$
\\[2ex]
\indent
Now, taking into account the system (12.64), we obtain
\\[2ex]
\mbox{}\hfill
$
{\frak I}\;\!p_1^{(2)}(x,y)=
\bigl(2\;\!\alpha_1^{(1)}\!+\lambda_1^{}\;\!\alpha_{1}^{(2)}-a_3^{}\;\!\gamma_{1}^{(2)}\bigr)\;\! x+
\bigl(2\beta_1^{(1)}\!+\lambda_1^{}\;\!\beta_{1}^{(2)}-b_3^{}\;\!\gamma_{1}^{(2)}\bigr)\;\! y+
\bigl(2\gamma_1^{(1)}\!+\lambda_1^{}\;\!\gamma_{1}^{(2)}-c_3^{}\;\!\gamma_{1}^{(2)}\bigr)\ -
\hfill
$
\\[2ex]
\mbox{}\hfill
$
- \bigl(\alpha_{1}^{(2)}x+\beta_{1}^{(2)}y\bigr)\;\!l_3^{}(x,y)\!=\!
2\bigl(\alpha_1^{(1)}x+\beta_1^{(1)}y+\gamma_1^{(1)}\bigr)
+\lambda_1^{}\bigl(\alpha_1^{(2)}x+\beta_1^{(2)}y+\gamma_1^{(2)}\bigr)-
\gamma_{1}^{(2)}(a_3^{}x+b_3^{}y+c_3^{}) -
\hfill
$
\\[2ex]
\mbox{}\hfill
$
-\ \bigl(\alpha_{1}^{(2)}x+\beta_{1}^{(2)}y\bigr)\;\!l_3^{}(x,y)=
2\bigl(\alpha_1^{(1)}x+\beta_1^{(1)}y+\gamma_1^{(1)}\bigr)
+\lambda_1^{}\;\! \bigl(\alpha_1^{(2)}x+\beta_1^{(2)}y+\gamma_1^{(2)}\bigr)\ -
\hfill
$
\\[2ex]
\mbox{}\hfill
$
-\ 
\bigl(\alpha_{1}^{(2)}x+\beta_{1}^{(2)}y+\gamma_{1}^{(2)}\bigr)\;\!l_3^{}(x,y)=
2\bigl(\alpha_1^{(1)}x+\beta_1^{(1)}y+\gamma_1^{(1)}\bigr)+
\bigl(\alpha_{1}^{(2)}x+\beta_{1}^{(2)}y+\gamma_{1}^{(2)}\bigr)
\bigl(\lambda_1^{}-l_3^{}(x,y)\bigr)=
\hfill
$
\\[2ex]
\mbox{}\hfill
$
=2\;\!p_1^{(1)}(x,y)+p_1^{(2)}(x,y)\;\!M_1^{}(x,y)
$ 
\ for all $(x,y)\in\R^2.
\hfill
$
\\[2ex]
\indent
Hence the derivative by virtue of the Jacobi equation (12.36)
\\[2ex]
\mbox{}\hfill  % (12.65)
$
{\frak I}\;\!p_1^{(2)}(x,y)=
2\;\!p_1^{(1)}(x,y)+p_1^{(2)}(x,y)\;\!M_1^{}(x,y)
$ 
\ for all $(x,y)\in\R^2.
$
\hfill (12.65)
\\[2ex]
\indent
Using the identities (12.62) and (12.65), we get on the set $D$ 
\\[2ex]
\mbox{}\hfill  
$
{\frak I}\,\dfrac{p_1^{(2)}(x,y)}{p_1^{}(x,y)}=
\dfrac{p_1^{}(x,y)\,{\frak I}\;\!p_1^{(2)}(x,y)-
p_1^{(2)}(x,y)\,{\frak I}\;\!p_1^{}(x,y)}{p_1^{2}(x,y)}=
\hfill  
$
\\[2ex]
\mbox{}\hfill  
$
=\dfrac{p_1^{}(x,y)\bigl(2\;\!p_1^{(1)}(x,y)+p_1^{(2)}(x,y)\;\!M_1^{}(x,y)\bigr)-
p_1^{}(x,y)\;\!p_1^{(2)}(x,y)\;\!M_1^{}(x,y)}{p_1^{2}(x,y)}=
2\,\dfrac{p_1^{(1)}(x,y)}{p_1^{}(x,y)}\,.
\hfill
$
\\[2ex]
\indent
From the identity (12.63) it follows that 
\\[1.75ex]
\mbox{}\hfill  
$
{\frak I}\,\biggl(\dfrac{p_1^{(1)}(x,y)}{p_1^{}(x,y)}\biggr)^{\!2}=
2\,\dfrac{p_1^{(1)}(x,y)}{p_1^{}(x,y)}
$ 
\ for all $(x,y)\in D.
\hfill
$
\\[2ex]
\indent
Therefore the derivative by virtue of the Jacobi equation (12.36) on the set $D$ is
\\[2ex]
\mbox{}\hfill   
$\!\!
{\frak I}F(x,y)=
{\frak I}\,\dfrac{\bigl(p_1^{(1)}(x,y)\bigr)^2\!-p_1^{}(x,y)\;\!p_1^{(2)}(x,y)}{p_1^{2}(x,y)}=
{\frak I}\biggl(\dfrac{p_1^{(1)}(x,y)}{p_1^{}(x,y)}\biggr)^{\!2}\!-
{\frak I}\,\dfrac{p_1^{(2)}(x,y)}{p_1^{}(x,y)}=0.
\hfill
$
\\[2ex]
\indent
Thus the function (12.61) is a general integral on $D$ 
\vspace{0.5ex}
of the Jacobi equation (12.36).$\k$

\newpage

{\bf Remark 12.2.}
Since the derivative by virtue of the Jacobi equation (12.36)
\\[1.5ex]
\mbox{}\hfill   
$
{\frak I}\,\dfrac{\bigl(p_1^{(1)}(x,y)\bigr)^2\!-p_1^{}(x,y)\;\!p_1^{(2)}(x,y)}{p_1^{2}(x,y)}=0
$
\ for all $(x,y)\in D,
\hfill
$
\\[1.75ex]
we see that, by Definition 5.1, 
the polynomial partial integral $p_1^{}$
\vspace{0.5ex}
of the Jacobi equation (12.36) is multiple and
\vspace{0.75ex}
$\bigl((p_1^{}, M_1^{}), \bigl(2,\bigl( p_1^{(1)}\bigr)^2-p_1^{}p_1^{(2)}, 0\bigr)\bigr)\in 
\text{B}_{_{\!\scriptstyle D}}.$
Moreover, from Lemma 12.2 it follows that
\vspace{0.75ex}
$\bigl((p_1^{}, M_1^{}), (1,p_1^{(1)}, 1)\bigr)\in \text{B}_{_{\!\scriptstyle D}}.$
Then, by Definition 5.2, 
the polynomial partial integral  $p_1^{}$ of the Jacobi equation (12.36) is double.
\vspace{0.75ex}

Thus  
\vspace{0.35ex}
{\it 
if the matrix $A$ has the eigenvalue $\lambda_1^{}$ with corresponding 
triple elementary devisor and the eigenvector
$(\alpha_1^{},\beta_1^{},\gamma_1^{}),$ 
then the polynomial partial integral
\\[1.5ex]
\mbox{}\hfill
$
p_1^{}\colon (x,y)\to \alpha_1^{}x+\beta_1^{}y+\gamma_1^{}$ 
\ for all $(x,y)\in\R^2
\hfill
$ 
\\[1.5ex]
of the Jacobi equation {\rm (12.36)} is double.
}
\vspace{1.25ex}

{\bf Example 12.7.}
Consider the Jacobi equation
\\[2ex]
\mbox{}\hfill  % (12.66)
$
\dfrac{dy}{dx}=\dfrac{x+3y+\delta-y\;\!(x+y+2)}{x-y-\delta-x\;\!(x+y+2)}\,,
$
\hfill (12.66)
\\[2ex]
where $\delta\ne 0$ is a real parameter.
The matrix
\vspace{0.5ex}
$
A=\left(\!\!\!
\begin{array}{rcc}
1 & 1 & 1
\\
{}-1 & 3 & 1
\\
{}-\delta & \delta & 2\end{array}
\!\!\right)
$
has the eigenvalue $\lambda=2$ 
\vspace{0.35ex}
with the triple elementary devisor $(\lambda-2)^3,$ 
the eigenvector $(1, 1, 0),$  
the first order generalized eigenvector $(0, 0, 1),$ 
\vspace{0.75ex}
and the second order generalized eigenvector 
$\Bigl({}-\dfrac{2}{\delta}\,,\;\! 0, {}-\dfrac{2}{\delta}\Bigr).$ 

By Theorem 12.13, a general integral on the set  
\vspace{0.25ex}
$D=\{(x,y)\colon x+y\ne 0\}$ 
of the Jacobi equation (12.66) is the rational function
\\[1.5ex]
\mbox{}\hfill  
$
F^{\ast}\colon (x,y)\to\ 
\dfrac{1-(x+y)\Bigl({}-\dfrac{2}{\delta}\,x-\dfrac{2}{\delta}\Bigr)}{(x+y)^2}
$
\ for all $(x,y)\in D
\hfill 
$
\\[1.75ex]
or (using Property 0.1 about 
functional ambiguity of general integral) the function
\\[1.75ex]
\mbox{}\hfill  % (12.67)
$
F\colon (x,y)\to\ 
\dfrac{2\;\!(x+y)(x+1)+\delta}{(x+y)^2}
$
\ for all $(x,y)\in D.
$
\hfill  (12.67)
\\[1.75ex]
\indent
Thus, as the result of solving the Darboux problem (Application 12.1)
for the Jacobi equation (12.36) we obtained:  
a general integral is constructed without quadratures.
\vspace{0.25ex}

In [9, pp. 14 -- 25], the generalized Darboux problem (Applications 12.1 and 12.9) 
for the Jacobi equation (12.36) solved. In this case, 
if the matrix $A$ has simple elementary devisors, then we build a general integral 
and 
if the matrix $A$ has multiple elementary devisor, then we 
can construct an integrating factor. 
\vspace{1ex}

{\bf 12.23.} Since 
{\it the Jacobi system}
\\[2ex]
\mbox{}\hfill        % (12.68)
$
\dfrac{dx}{dt}=l_1^{}(x,y)-x\;\!l_3^{}(x,y),
\quad
\dfrac{dy}{dt}=l_2^{}(x,y)-y\;\!l_3^{}(x,y),
\hfill
$
\\
\mbox{}\hfill (12.68)
\\
\mbox{}\hfill
$
l_i^{}(x,y)=a_{i}^{}\;\!x+b_{i}^{}\;\!y+c_i^{}\;\!,
\quad
a_{i}^{},b_{i}^{},c_i^{}\in\R, 
\quad
i=1,2,3,
\hfill
$
\\[2ex]
is an autonomous differential system of the second-order, 
we see that this system has 
two functionally independent first integrals 
(one of these first integrals is non-autonomous).

The Jacobi equation (12.36) is the equation of trajectories for 
the Jacobi system (12.68). Then a general integral of the Jacobi equation (12.36) 
is an autonomous first integral of the Jacobi system  (12.68).

Thus, to construction an integral basis of the Jacobi system (12.68), it is sufficient 
to find non-autonomous first integral of this system.

The differentiation operator by virtue of the Jacobi system (12.68) is 
the linear differential operator of first-order
\\[1.5ex]
\mbox{}\hfill
$
{\frak d}(t,x,y)=
\partial_{{}_{\scriptstyle t}}+ {\frak I}(x,y)
$
\ for all $(t,x,y)\in\R^3,
\hfill
$
\\[1.5ex]
where ${\frak I}$ is the differentiation operator 
by virtue of the Jacobi equation (12.36).

Taking into account the 
connection between the differentiation operators ${\frak d}$ and ${\frak I},$ 
we obtain that partial integrals of the Jacobi equation (12.36) are 
autonomous partial integrals of the Jacobi system  (12.68).
\vspace{0.75ex}

{\small\bf 
12.23.1. Case of distinct real eigenvalues.}
\vspace{0.5ex}

{\bf Theorem 12.14.}\!
\vspace{0.35ex}
{\it
Suppose $\lambda_1^{}\!$ and $\lambda_2^{}\!$ are
distinct real eigenvalues of the matrix $A$ with 
corresponding eigenvectors 
\vspace{0.75ex}
$
(\alpha_1^{},\beta_1^{},\gamma_1^{})
$ 
and
$
(\alpha_2^{},\beta_2^{},\gamma_2^{}). 
$
Then an non-autonomous first integ\-ral
\vspace{0.15ex}
on the set  
$\!\Omega_0^{}\!=\!\{(t,x,y)\colon\! \alpha_2^{}x+\beta_{2}^{}y+\gamma_2^{}\!\ne\! 0\}\!$ 
of the Jacobi system {\rm(12.68)} is the function 
\\[2ex]
\mbox{}\hfill  % (12.69)
$
\Psi\colon (t,x,y)\to\ 
\dfrac{p_1^{}(x,y)}{p_2^{}(x,y)}\,\exp\bigl((\lambda_2^{}-\lambda_1^{})\;\!t\bigr)
$
\ for all $(t,x,y)\in \Omega_0^{},
$
\hfill {\rm (12.69)}
\\[2.25ex]
where $p_1^{}(x,y)=\alpha_1^{}\;\!x+\beta_{1}^{}\;\!y+\gamma_1^{},$ 
$p_2^{}(x,y)=\alpha_2^{}\;\!x+\beta_{2}^{}\;\!y+\gamma_2^{}$ 
for all $(x,y)\in\R^2.$
}
\vspace{1ex}

{\sl Proof.}
\vspace{0.5ex}
From Lemma 12.1 it follows that the linear functions $p_1^{}$ and $p_2^{}$ are 
autonomous polynomial partial integrals with cofactors
\vspace{0.75ex}
$M_1^{}\colon (x,y)\to \lambda_1^{}-l_3^{}(x,y)$ for all $(x,y)\in\R^2$ and
$M_2^{}\colon (x,y)\to \lambda_2^{}-l_3^{}(x,y)$ for all $(x,y)\in\R^2$ 
of the Jacobi system (12.68).
\vspace{0.75ex}

By Property 4.1, 
\vspace{0.5ex}
the function $e^t$ is an non-autonomous conditional partial integral
with cofactor  $M_3^{}\colon t\to 1$ for all $t\in\R$ 
of the Jacobi system (12.68).
\vspace{0.5ex}

Since the cofactors $M_1^{},\ M_2^{},\ M_3^{}$ of the 
partial integrals $p_1^{},\ p_2^{},\ e^t$ such that the following identity holds
\\[1.5ex]
\mbox{}\hfill
$
M_1^{}(x,y)-M_2^{}(x,y)+(\lambda_2^{}-\lambda_1^{})M_1^{}(t)=0
$
\ for all $(t,x,y)\in \R^3,
\hfill
$ 
\\[1.5ex]
we see that (Property 11.11) 
\vspace{0.15ex}
the function (12.69) is an non-autonomous first integral on 
the set $\Omega_0^{}$ of the Jacobi system (12.68). $\k$
\vspace{1.25ex}

{\bf Example 12.8.}
The Jacobi system
\\[2ex]
\mbox{}\hfill    % (12.70)
$
\dfrac{dx}{dt}=3x-y+1-x\;\!(x-y+3),
\quad
\dfrac{dy}{dt}={}-x+5y-1-y\;\!(x-y+3)
$
\hfill (12.70)
\\[2ex]
such that the matrix $A$ has (Example 12.1) 
\vspace{0.25ex}
the distinct real eigenvalues
$\lambda_1^{}=2,\ \lambda_2^{}=3,$ $\lambda_3^{}=6$
and the corresponding eigenvectors 
$(1, 0, {}-1),\ (1, 1, 1),\ (1, {}-2, 1).$
\vspace{0.5ex}

It follows from Theorem 12.14 that 
non-autonomous first integrals of the Jacobi system (12.70) are the functions 
\\[2ex]
\mbox{}\hfill  % (12.71)
$
\Psi_{1{,}2}^{}\colon (t,x,y)\to\ 
\dfrac{x-1}{x+y+1}\, e^t
$
\ for all $(t,x,y)\in\Omega_1^{},
$
\hfill (12.71)
\\[2.25ex]
\mbox{}\hfill  % (12.72)
$
\Psi_{1{,}3}^{}\colon (t,x,y)\to\ 
\dfrac{x-1}{x-2y+1}\, e^{\;\!4t}
$
\ for all $(t,x,y)\in\Omega_2^{},
$
\hfill (12.72)
\\[2.25ex]
\mbox{}\hfill  % (12.73)
$
\Psi_{2{,}3}^{}\colon (t,x,y)\to\ 
\dfrac{x+y+1}{x-2y+1}\, e^{\;\!3t}
$
\ for all $(t,x,y)\in\Omega_2^{},
$
\hfill (12.73)
\\[2ex]
where the sets $\Omega_1^{}=\{(t,x,y)\colon x+y+1\ne 0\},\ 
\Omega_2^{}=\{(t,x,y)\colon x-2y+1\ne 0\}.$
\vspace{0.75ex}

Each of the collections of functions 
\vspace{0.35ex}
$(12.43)\& (12.71),\ (12.43)\& (12.72),\ (12.43)\& (12.73),
\linebreak 
(12.71)\&  (12.72),\ (12.71)\& (12.73),\ (12.72)\& (12.73)$
\vspace{0.35ex}
is an integral basis (on the corresponding set) 
of the Jacobi system (12.70).
\vspace{0.75ex}

{\bf Example 12.9.}
Consider the Jacobi system
\\[2ex]
\mbox{}\hfill    % (12.74)
$
\dfrac{dx}{dt}={}-x+y+1-x\;\!(x+y-1),
\quad
\dfrac{dy}{dt}=x-y+1-y\;\!(x+y-1).
$
\hfill (12.74)
\\[2ex]
\indent
The matrix $A$ (Example 12.2) has 
\vspace{0.35ex}
the distinct real eigenvalues
$\lambda_1^{}={}-2$ and $\lambda_3^{}=1.$
\vspace{0.35ex}
The eigenvalue $\lambda_1^{}={}-2$ has two
linearly independent eigenvectors
$(0, 1, -1)$ and $(1, 0, {}-1).$
\vspace{0.25ex}
The eigenvalue  $\lambda_3^{}=1$ has 
the corresponding eigenvector $(1, 1, 1).$
\vspace{0.35ex}

By Theorem 12.14, non-autonomous first integrals on the set
$\Omega_0^{}=\{(t,x,y)\colon x+y+1\ne 0\}$ 
of the Jacobi system (12.74) are the functions 
\\[2ex]
\mbox{}\hfill  % (12.75)
$
\Psi_{1}^{}\colon (t,x,y)\to\ 
\dfrac{y-1}{x+y+1}\, e^{\;\!3t}
$
\ for all $(t,x,y)\in\Omega_0^{},
$
\hfill (12.75)
\\[2ex]
\mbox{}\hfill  % (12.76)
$
\Psi_{2}^{}\colon (t,x,y)\to\ 
\dfrac{x-1}{x+y+1}\, e^{\;\!3t}
$
\ for all $(t,x,y)\in\Omega_0^{}.
$
\hfill (12.76)
\\[2ex]
\indent
Each of the collections of functions 
\vspace{0.25ex}
$(12.46)\& (12.75),\ (12.46)\& (12.76),\ (12.75)\& (12.76)$
is an integral basis (on the corresponding set) 
of the Jacobi system (12.74).
\vspace{1.25ex}

{\small\bf 
12.23.2. Case of complex eigenvalue.}
\vspace{0.5ex}

{\bf Theorem 12.15.}
\vspace{0.5ex}
{\it
Suppose $\lambda=\xi+\zeta\;\! i\ (\zeta\ne 0)$  is 
the complex eigenvalue of the matrix $A$ 
corresponding to the eigenvector
\vspace{0.5ex}
$
(\alpha,\beta,\gamma).
$ 
Then an non-autonomous first integral on the set
\vspace{0.15ex}
$\Omega_0^{}=\{(t,x,y)\colon {\rm Re}(\alpha x+\beta y+\gamma)\ne 0\}$ 
of the Jacobi system {\rm(12.68)} is the function 
\\[2ex]
\mbox{}\hfill  % (12.77)
$
\Psi\colon (t,x,y)\to\ 
\arctan\dfrac{{\rm Im}\;\!p(x,y)}{{\rm Re}\;\!p(x,y)}\,-\,\zeta\;\!t
$
\ for all $(t,x,y)\in \Omega_0^{},
$
\hfill {\rm (12.77)}
\\[2.25ex]
where $p(x,y)=\alpha\;\!x+\beta\;\!y+\gamma$ for all $(x,y)\in\R^2.$
}
\vspace{0.5ex}

{\sl Proof.}
\vspace{0.5ex}
As in the proof of Theorem 12.11, we establish that the function
$\exp\,\arctan\dfrac{{\rm Im}\;\!p}{{\rm Re}\;\!p}$
is an autonomous exponential partial integral with cofactor
\vspace{0.35ex}
$M\colon (x,y)\to \zeta$ for all $(x,y)\in\R^2$
on the set $\Omega_0^{}$ of the Jacobi system  (12.68).
\vspace{0.5ex}

By Property 4.1, 
\vspace{0.35ex}
the function $e^t$ is 
an non-autonomous conditional partial integral with cofactor
$M_1^{}\colon t\to 1$ for all $t\in\R$ of the Jacobi system (12.68).
\vspace{0.25ex}

Since the cofactors $M$ and $M_1^{}$ of the partial integrals
\vspace{0.25ex}
$\exp\,\arctan\dfrac{{\rm Im}\;\!p}{{\rm Re}\;\!p}$ and $e^t$ 
such that the following identity holds
\\[1.25ex]
\mbox{}\hfill
$
M(x,y)-\zeta\;\!M_1^{}(t) = 0
$
\ for all $(t,x,y)\in \R^3,
\hfill
$ 
\\[1.25ex]
we see that (using Property 11.11 and taking into account Property 0.1 
\vspace{0.25ex}
about functional ambiguity of first integral) 
the function (12.77) is an non-autonomous first integral on 
the set $\Omega_0^{}$ of the Jacobi system (12.68). $\k$
\vspace{1.25ex}

{\bf Example 12.10.}
The Jacobi system
\\[2ex]
\mbox{}\hfill    % (12.78)
$
\dfrac{dx}{dt}=4x-3y-1-x\;\!({}-2x+y),
\quad
\dfrac{dy}{dt}=6x-2y+1-y\;\!({}-2x+y)
$
\hfill (12.78)
\\[2ex]
such that the matrix $A$ has (Example 12.3) 
\vspace{0.5ex}
the complex eigenvalue 
$\lambda=1+\sqrt{6}\;\!i$ 
with the corresponding eigenvector 
$\bigl({}-2\bigl(2+\sqrt{6}\;\!i\bigr),\;\! 5,\;\! 3-\sqrt{6}\;\!i\bigr).$
\vspace{0.5ex}

It follows from Theorem 12.15 that 
\vspace{0.35ex}
an non-autonomous first integral of the Jacobi system (12.78) on the set 
$\Omega_0^{}=\{(t,x,y)\colon 4x-5y-3\ne 0\}$ 
 is the function 
\\[2ex]
\mbox{}\hfill  % (12.79)
$
\Psi\colon (t,x,y)\to\ 
\arctan\dfrac{\sqrt{6}\,(2x+1)}{4x-5y-3}\,-\,\sqrt{6}\,t
$
\ for all $(t,x,y)\in\Omega_0^{}.
$
\hfill (12.79)
\\[2ex]
\indent
The collection of functions $(12.49)\& (12.79)$
\vspace{0.35ex}
is an integral basis of the Jacobi system  (12.78) on the set  
$\Omega^{\ast}=\{(t,x,y)\colon (x+y+5)(4x-5y-3)\ne 0\}.$
\vspace{1ex}

{\bf Example 12.11.}
Consider the Jacobi system  
\\[2ex]
\mbox{}\hfill    % (12.80)
$
\dfrac{dx}{dt}=x+2y+3-x\;\!(x+2y+1),
\quad
\dfrac{dy}{dt}=x+y-3-y\;\!(x+2y+1).
$
\hfill (12.80)
\\[2ex]
\indent
The matrix $A$ has (Example 12.4) 
\vspace{0.35ex}
the complex eigenvalue $\lambda=1+i$ 
with the cor\-res\-pon\-ding eigenvector $(3-i,\;\! 4, {}-3+3i).$
\vspace{0.5ex}
By Theorem 12.15, 
an non-autonomous first in\-teg\-ral on the set 
$\Omega_1^{}=\{(t,x,y)\colon 3x+4y-3\ne 0\}$
of the Jacobi system (12.80) is the function
\\[2ex]
\mbox{}\hfill  % (12.81)
$
\Psi\colon (t,x,y)\to\ 
\arctan\dfrac{{}-x+3}{3x+4y-3}\,-\,t
$
\ for all $(t,x,y)\in\Omega_1^{}.
$
\hfill (12.81)
\\[2ex]
\indent
The collection of functions $(12.51)\& (12.81)$
\vspace{0.5ex}
is an integral basis of the Jacobi system  (12.80) on the set  
$\Omega_0^{}=\{(t,x,y)\colon (x+y-1)(3x+4y-3)\ne 0\}.$
\vspace{1.5ex}

{\small\bf 
12.23.3. Case of multiple elementary devisor.}
\vspace{0.5ex}

{\bf Теорема 12.16.}
\vspace{0.35ex}
{\it
Suppose the matrix $A$ has the real eigenvalue with corresponding 
multiple elementary devisor, the eigenvector 
\vspace{0.5ex}
$(\alpha_1^{},\beta_1^{},\gamma_1^{}),$ and 
the first order generalized eigenvector
$\bigl(\alpha_1^{(1)},\beta_1^{(1)},\gamma_1^{(1)}\bigr).$ 
\vspace{1ex}
Then an non-autonomous first integral 
of the Jacobi system {\rm (12.68)} on the set  
$\Omega_0^{}=\{(t,x,y)\colon \alpha_1^{}x+\beta_1^{}y+\gamma_1^{}\ne 0\}$ 
is the function 
\\[2ex]
\mbox{}\hfill   % (12.82)
$
\Psi\colon (x,y)\to \
\dfrac{p_1^{(1)}(x,y)}{p_1^{}(x,y)}\,-\,t
$
\ for all $(t,x,y)\in \Omega_0^{},
$
\hfill {\rm (12.82)}
\\[2ex]
where 
\vspace{1.25ex}
$
p_1^{}(x,y)= \alpha_1^{}x+\beta_1^{}y+\gamma_1^{}, \ 
p_1^{(1)}(x,y)= \alpha_1^{(1)}x+\beta_1^{(1)}y+\gamma_1^{(1)}$ 
for all $(x,y)\in\R^2.$
}

{\sl Proof.}
\vspace{0.5ex}
As in the proof of Theorem 12.12, we establish that the function
$\exp\dfrac{p_1^{(1)}}{p_1^{}}$
is an autonomous exponential partial integral with cofactor
\\[1.24ex]
\mbox{}\hfill
$
M_1^{}\colon (x,y)\to 1$ 
\ for all $(x,y)\in\R^2
\hfill
$
\\[1.5ex]
on the set $\Omega_0^{}$ of the Jacobi system  (12.68).
\vspace{0.75ex}

By Property 4.1, 
\vspace{0.5ex}
the function $e^t$ is 
an non-autonomous conditional partial integral with cofactor
$M_1^{(1)}\colon t\to 1$ for all $t\in\R$ of the Jacobi system (12.68).
\vspace{0.75ex}

From Corollary 11.7  
\vspace{0.35ex}
(taking into account Property 0.1 
about functional ambiguity of first integral) 
it follows that the function (12.82) is
\vspace{0.35ex}
an non-autonomous first integral 
 on the set $\Omega_0^{}$ of the Jacobi system (12.68). $\k$
\vspace{1ex}

{\bf Example 12.12.}
The Jacobi system 
\\[2ex]
\mbox{}\hfill    % (12.83)
$
\dfrac{dx}{dt}={}-x+y-x\;\!({}-x+y),
\quad
\dfrac{dy}{dt}=x-y-1-y\;\!({}-x+y)
$
\hfill (12.83)
\\[2ex]
such that the matrix $A$ has (Example 12.5) 
\vspace{0.5ex}
the distinct real eigenvalues
$\lambda_1^{}={}-1$ and $\lambda_3^{}=0$
\vspace{0.25ex}
with the corresponding eigenvectors 
$(1, {}-1, {}-1)$ and $(1, 0, {}-1).$
\vspace{0.25ex}

By Theorem 12.14, an non-autonomous first integral on the set
\vspace{0.25ex}
$\Omega_1^{}=\{(t,x,y)\colon x\ne 1\}$ 
of the Jacobi system (12.83) is the function 
\\[2ex]
\mbox{}\hfill  % (12.84)
$
\Psi_{1}^{}\colon (t,x,y)\to\ 
\dfrac{x-y-1}{x-1}\, e^{t}
$
\ for all $(t,x,y)\in\Omega_1^{}.
$
\hfill (12.84)
\\[2ex]
\indent
The eigenvalue $\lambda_1^{}={}-1$ of the matrix $A$
\vspace{0.35ex}
has double elementary devisor, the eigenvector  $(1, {}-1,{}-1),$ 
and the first order generalized eigenvector $(1, {}-1, {}-2).$
\vspace{0.5ex}

It follows from Theorem 12.16 that 
\vspace{0.5ex}
an non-autonomous first integral of the Jacobi system (12.83) 
on the set $\Omega_2^{}=\{(t,x,y)\colon x-y-1\ne 0\}$ is the function
\\[2ex]
\mbox{}\hfill  % (12.85)
$
\Psi_{2}^{}\colon (t,x,y)\to\ 
\dfrac{x-y-2}{x-y-1}\,-\,t
$
\ for all $(t,x,y)\in\Omega_2^{}.
$
\hfill (12.85)
\\[2ex]
\indent
Each of the collections of functions 
\vspace{0.5ex}
$(12.57)\& (12.84),\, (12.57)\& (12.85),\, (12.84)\& (12.85)$ is 
an integral basis on the set
$\!\Omega_0^{}\!=\!\{(t,x,y)\colon\! (x\!-\!1)(x\!-\!y\!-\!1)\!\ne 0\}\!\!$ 
\vspace{1.25ex}
of the Jacobi system (12.83).

{\bf Example 12.13.}
Consider the Jacobi system 
\\[2ex]
\mbox{}\hfill    % (12.86)
$
\dfrac{dx}{dt}=x-x\;\!(y+1),
\quad
\dfrac{dy}{dt}=y-y\;\!(y+1).
$
\hfill (12.86)
\\[2ex]
\indent
The matrix $A$ has (Example 12.6) 
\vspace{0.25ex}
the eigenvalue with the corresponding 
double ele\-men\-ta\-ry devisor, the eigenvector $(0, 1, 0),$ 
\vspace{0.5ex}
and the first order generalized eigenvector  $(1, 1, 1).$

By Theorem 12.16, 
\vspace{0.25ex}
an non-autonomous first integral on the set 
$\Omega_1^{}=\{(t,x,y)\colon y\ne 0\}$ 
of the Jacobi system (12.86) is the function 
\\[2ex]
\mbox{}\hfill  % (12.87)
$
\Psi\colon (t,x,y)\to\ 
\dfrac{x+y+1}{y}\,-\,t
$
\ for all $(t,x,y)\in\Omega_1^{}.
$
\hfill (12.87)
\\[2ex]
\indent
The collection of functions $(12.59)\& (12.87)$ is 
\vspace{0.35ex}
an integral basis of the Jacobi system (12.86)
on the set  $\Omega_0^{}=\{(t,x,y)\colon xy\ne 0\}.$ 
\vspace{1.5ex}

{\bf Example 12.14.}
Consider the Jacobi system 
\\[2ex]
\mbox{}\hfill    % (12.88)
$
\dfrac{dx}{dt}=x-y-\delta-x\;\!(x+y+2),
\quad
\dfrac{dy}{dt}=x+3y+\delta-y\;\!(x+y+2),
$
\hfill (12.88)
\\[2ex]
where $\delta\ne 0$ is a real parameter. 
\vspace{0.35ex}
The matrix $A$ has (Example 12.7) 
the eigenvalue with the corresponding 
triple elementary devisor, the eigenvector $(1, 1, 0),$ 
\vspace{0.25ex}
the first order generalized eigenvector $(0, 0, 1).$
\vspace{0.5ex}
By Theorem 12.16, 
an non-autonomous first integral of the Jacobi system (12.88) on the set 
$\Omega_0^{}=\{(t,x,y)\colon x+y\ne 0\}$ is the function 
\\[1.75ex]
\mbox{}\hfill  % (12.89)
$
\Psi\colon (t,x,y)\to\ 
\dfrac{1}{x+y}\,-\,t
$
\ for all $(t,x,y)\in\Omega_0^{}.
$
\hfill (12.89)
\\[1.75ex]
\indent
The collection of functions $(12.67)\& (12.89)$ is 
\vspace{0.15ex}
an integral basis on the set $\Omega_0^{}$ of the Jacobi system (12.88).
\vspace{1ex}

The method of constructing the integral basis for the ordinary differential Jacobi system of the second order (described in Application 12.23) 
in the monograph [1, pp. 273 -- 311] extended to the Jacobi system of 
equations in total differentials, the special case of which is the ordinary differential Jacobi system of the higher order.
\vspace{0.5ex}

Let us remark that the research [56] contains the extensive bibliography of articles that considered partial  integrals and their applications.

Note also that many other results concerning partial integrals 
can be found in [57 -- 72].

\newpage

\mbox{}
\\[-2.25ex]

{\Large\bf  References}
\vspace{1.75ex}

1.
V.N. Gorbuzov, 
{\it Integrals of  differential systems} (Russian), Grodno State University, Grodno, 2006.
\vspace{0.75ex}

2.
V.N. Gorbuzov, 
Integral equivalence of multidimensional differential systems,  
{\it Mathematics.Dynami\-cal Systems} 
(arXiv: 0909.3220v1 [math.DS]. Cornell Univ., Ithaca, New York), 2009, 1-45.
\vspace{0.75ex}

3.
N.M. Matveev, 
{\it Methods of integration of ordinary differential equations} (Rus\-si\-an),
Lan', Saint-Petersburg, 2003.
\vspace{0.75ex}

4. 
N.P. Erugin, 
{\it A book for reading the general course of differential equations} (Russian), 
Nauka and Technika, Minsk, 1979.
\vspace{0.75ex}

5.
E. Goursat, 
{\it A Course of Mathematical Analysis} (Russian), P. 2,
ONTI, Moscow-Le\-nin\-g\-rad, 1936.
\vspace{0.75ex}

6. 
E.J. Cartan,
{\it Integral invariants} (Russian), 
GITTL, Moscow-Leningrad, 1940.
\vspace{0.75ex}

7. 
V.N. Gorbuzov,
On the question of integrability in quadratures (Russian), 
{\it Doklady Akademii Nauk BSSR}, Vol. 25, 1981, No. 7,  584-585.
\vspace{0.75ex}

8.
 G. Darboux,
M\'{e}moire sur les \'{e}quations differentielles algebriques du premier ordre et du premier degr\'{e},
{\it Bulletin des Sciences Math\'{e}matiques}, 1878,  Vol. 2,  60-96.
\vspace{0.75ex}

9.
V.N. Gorbuzov  and A.A. Samodurov, 
{\it The Darboux equation and its analogs} (Russian), 
Grodno State University, Grodno, 1985.
\vspace{0.75ex}

10.
N.N. Babariko and V.N. Gorbuzov,
On the question of the integrability of first-order nonlinear differential equations (Russian), 
{\it Doklady Akademii Nauk BSSR}, Vol. 28, 1984, 
\linebreak
No. 7,  581-584.
\vspace{0.75ex}

11.
V.N. Gorbuzov  and A.A. Samodurov, 
{\it The Riccati and Abel equations} (Russian), 
Grodno State University, Grodno, 1986.
\vspace{0.75ex}

12.
V.N. Gorbuzov,
About some classes of autonomous systems with partial integral (Rus\-sian),
{\it Differential Equations}, Vol. 17, 1981, No. 9, 1685-1687.
\vspace{0.75ex}

13.
M.N. Lagutinskii, 
{\it Partial algebraic integrals}  (Russian), 
Adol'f Darre, Kharkov, 1908.
\vspace{0.75ex}

14.
K.S. Sibirsky, 
{\it  Algebraic invariants of differential equations and matrixes} (Russian), 
Shtiinca, Chisinau, 1976.
\vspace{0.75ex}

15.
V.N. Gorbuzov,
Building and the whole qualitative investigation for one class of auto\-no\-mo\-us systems (Russian),
{\it Investigation on the mathematician and the physicist} (Grodno), 1978, 26-33.
\vspace{0.75ex}

16.
V.N. Gorbuzov,
The whole qualitative investigation for one class of auto\-no\-mo\-us systems\! (Russian),\!
{\it Investigation on the mathematician and the physicist}\! (Grodno),\! 
\vspace{0.75ex}
1978, 33-37.

17.
%[3340]
V.N. Gorbuzov,
Projective atlas of trajectories of differential systems,  
{\it Mathematics.Dy\-na\-mi\-cal Sys\-tems} 
(arXiv: 1401.1000v1 [math.DS].  Cornell Univ., 
\vspace{0.75ex}
Ithaca, New York), 
2014, 1-61.

18. 
H. Poincar\'{e},
{\it On curves defined by differential equations} (Russian),
GITTL, Moscow-Leningrad, 1947.
\vspace{0.75ex}

19.
%[3327]
V.N. Gorbuzov, 
Trajectories of projective reduced differential systems (Russian),
{\it Vestnik of the Yanka Kupala Grodno State University}, 2012, Ser. 2,  No. 1(126), 39-52. 
\vspace{0.75ex}

20.
%[3339]
V.N. Gorbuzov and P.B. Pavlyuchik, 
 Linear and open limit cycles of differential systems (Russian),
{\it Vestnik of the Yanka Kupala Grodno State University}, 2013, Ser. 2,  
No. 3(159), 
\linebreak
23-32. 
\vspace{0.75ex}

21.
N.N. Babariko and V.N. Gorbuzov,
On the question of constructing the first integral or the last multiplier of a 
nonlinear system of differential equations (Russian), 
{\it Doklady Akademii Nauk BSSR}, Vol. 30, 1986, No. 9,  791-792.
\vspace{0.75ex}

22.
N.P. Erugin,  
Construction of the set of all systems of differential equations having a given integral curve (Russian), 
{\it Applied mathematics and mechanics},
 Vol. 16, 1952, 
No. 6,  659-670.
\vspace{0.75ex}

23.
V.N. Gorbuzov and V.Yu. Tyshchenko,
Particular integrals of systems of ordinary differential equations,
{\it Matematicheskii Sbornik}, 1993, Vol. 75,  No. 2, 353-369.
\vspace{0.75ex}

24.
V.V Nemytskii and  V.V Stepanov, 
{\it Qualitative theory of differential equations} (Russian), 
GITTL, Moscow-Leningrad, 1949.
\vspace{0.75ex}

25.
V.N. Gorbuzov and V.Yu. Tyshchenko,
Partial integrals of total differential systems (Rus\-sian),
{\it Differential Equations}, 1991, Vol. 27,  No. 10, 1819-1822.
\vspace{0.75ex}

26.
M.V. Dolov and V.V. Kosarev,
Darboux integrals and analytic structure of solutions for differential equations (Rus\-sian),
{\it Differential Equations}, 1983, Vol. 19,  No. 4, 697-700.
\vspace{0.75ex}

27.
M.V. Dolov and B.V. Lisin,  
An integrating factor and limit cycles  (Rus\-sian), 
{\it Differential and Integral Equations} (Gorky State University), 1984, 36-41.
\vspace{0.75ex}

28.
V.V. Amel'kin and B.S. Kalitin,  
{\it Isochronous and pulse oscillations of two-dimensional dynamical systems} (Rus\-sian), 
KomKniga, Moscow, 2006.
\vspace{0.75ex}

29.
M.V. Dolov,
A canonical integral in the neighborhood of a focus (Rus\-sian),
{\it Differential Equations}, 1976, Vol. 12,  No. 11, 1946-1953.
\vspace{0.75ex}

30.
V.N.  Gorbuzov, 
On a second-order system of differential equations and its periodic solutions  (Rus\-sian),
{\it Differential Equations}, 1994, Vol. 30,  No. 9, 1487-1497.
\vspace{0.75ex}

31.
V.N.  Gorbuzov, 
Construction of first integrals and last multipliers for polynomial autonomous 
many-dimensional differential systems (Rus\-sian),
{\it Differential Equations}, 1998, Vol. 34,  No. 4, 562-564.
\vspace{0.75ex}

32.
V.N.  Gorbuzov, 
To the question on integrals and last multipliers of
multidimensional differential systems with symmetries (Rus\-sian),
{\it Proceedings of the National Academy of Sciences of Belarus}.
Series of pfysical-mathematical sciences, 1996,  No. 1, 20-24. 
\vspace{0.75ex}

33.
V.N. Gorbuzov, 
Autonomy of a system of equations in total differentials (Rus\-sian),
{\it Differential Equations}, 1998, Vol. 34,  No. 2, 149-156.
\vspace{0.75ex}

34.
D.V. Buslyuk and V.N. Gorbuzov,
Integrals for Jacobian system of partial differential equations  (Rus\-sian), 
{\it Vestnik of the Yanka Kupala Grodno State University}, 2000, Ser. 2, 
\linebreak
 No. 2, 4-11. 
\vspace{0.75ex}

35.
V.N. Gorbuzov, 
Particular integrals of real autonomous polynomial systems of exact differential equations (Russian), \!
{\it Differential equations and control processes},\! 2000, 
\vspace{0.75ex}
No. 2, 1-36.

36.
V.N. Gorbuzov and A.F. Pranevich, 
Spectral method of construction integral basis for jacobian systems in partial equations (Russian),  
{\it Differential equations and control processes}, 2001, No. 3, 17-45.
\vspace{0.75ex}

37.
V.N. Gorbuzov and P.B. Pauliuchyk, 
Solutions, integrals and limits cycles for the Darboux system of $n\!$-th order  (Russian),  
{\it Differential equations and control processes}, 2002, No. 2, 26-46.
\vspace{0.75ex}

38.
V.N. Gorbuzov and A.F. Pranevich, 
Building of  integrals of  linear differential systems (Russian),\! 
{\it Vestnik of the\! Yanka Kupala Grodno State University}, 2003, 
\vspace{0.75ex}
Ser.\! 2, No.\! 2(22), 50-60.   

39.
V.N. Gorbuzov and A.F. Pranevich, 
Integrals of R-linear systems of exact differentials (Russian),\! 
{\it Doklady of the National Academy of Sciences of Belarus},  2004, 
\vspace{0.75ex}
Vol.\! 48, No.\! 1, 49-52.

40.
V.N. Gorbuzov and S.N. Daranchuk,
A basis of the autonomous first integrals of the Jacobi --- Fourier's system (Russian), 
{\it Vestnik of the Belarusian State  University}, 2005, Ser. 1, No. 3, 70-74.
\vspace{0.75ex}

41.
V.N. Gorbuzov, 
First integrals of Pfaff system of equations (Russian), 
{\it Vestnik of the Yanka Kupala Grodno State  University}, 2005, Ser. 2(34), No. 2, 10-29.   
\vspace{0.75ex}

42.
V.N. Gorbuzov and S.N. Daranchuk,
The integral basis of Jacobi --- Hesse partial differential system (Russian), 
{\it Izvestia of the Russian State Pedagogical University}, 2005, 
\linebreak
No. 5, 65-76.  
\vspace{0.75ex}

43.
V.N. Gorbuzov and S.N. Daranchuk,
The spectral method of building of the first integrals of Jacobi's system (Russian), 
{\it Vestnik of the Yanka Kupala Grodno State  University}, 2007, 
Ser. 2,  No. 3(57), 63-67. 
\vspace{0.75ex}

44.
V.N. Gorbuzov and S.N. Daranchuk,
The integrals and last multipliers of one class of partial differential systems (Russian),  
{\it J. Differential equations and control processes}, 2007,  No. 4, 1-16.
\vspace{0.75ex}

45.
V.N. Gorbuzov and A.F. Pranevich, 
Autonomy and cylindricality of R-differentiable integrals for systems in total differentials (Russian),  
{\it Differential equations and control processes}, 2008,  No. 1, 35-49.
\vspace{0.75ex}

46.
V.N. Gorbuzov and A.F. Pranevich, 
First integrals of linear differential systems,  
{\it Mathematics.Classical Analysis and ODEs} 
(arXiv: 0806.4155v1 [math.CA]. Cornell Univ., Ithaca, New York),  2008, 1-37.
\vspace{0.75ex}

47.
V.N. Gorbuzov and A.F. Pranevich, 
Building of first integrals for linear non-autonomous multidimensional differential systems 
with non-derogatory matrix structure (Russian), 
{\it Vestnik of the Belarusian State  University}, 2008, Ser. 1, No. 2, 75-79.
\vspace{0.75ex}

48.
V.N. Gorbuzov and A.F. Pranevich, 
First integrals for one system of Lappo-Danilevskii (Russian),\!
{\it Vestnik of the Yanka Kupala Grodno State  University}, 
\vspace{0.75ex}
2008, Ser.\! 2,  No.\! 3(73), 79-83. 

49.
V.N. Gorbuzov and S.N. Daranchuk,
The integrals and last multipliers of one class of total differential systems in complex domain (Russian),  
{\it Vestnik of the Belarusian State  University}, 2008, Ser. 1, No. 3, 59-62.
\vspace{0.75ex}

50.
V.N. Gorbuzov,
Cylindricality  and autonomy of integrals and last multipliers 
of multi\-di\-men\-si\-o\-nal differentional systems,
{\it Mathematics.Dynamical Systems} (arXiv: 0909.3234v1 [math.DS]. 
 Cornell Univ., Ithaca, New York), 2009, 1-37.
\vspace{0.75ex}

51.
V.N. Gorbuzov and A.F. Pranevich, 
R-holomorphic solutions and R-differentiable integrals 
of multidi\-mensional differential systems,
{\it Mathematics.Dynamical Systems} (arXiv: 0909.3245v1 [math.DS]. 
 Cornell Univ., Ithaca, New York), 2009, 1-29.
\vspace{0.75ex}

52.
V.N. Gorbuzov and S.N. Daranchuk,
Conditional partial integrals for system of total differential equations (Russian), 
{\it Vestnik of the Yanka Kupala Grodno State  University}, 2010, 
Ser. 2,  No. 2(96), 41-48. 
\vspace{0.75ex}

53.
V.N. Gorbuzov and A.F. Pranevich, 
Integrals of the Lappo-Danilevskii multidimensional differential system (Russian),  
{\it Differential equations and control processes}, 
\linebreak
2015,  No. 1, 1-24.
\vspace{0.75ex}

54.
V.N. Gorbuzov and A.F. Pranevich, 
First integrals of ordinary linear differential systems, 
{\it Mathematics.Dynamical Systems}
(arXiv: 1201.4141v1 [math.DS]. Cornell Univ., Ithaca, New York), 2012, 1-75.
\vspace{0.75ex}

55.
V.N. Gorbuzov,
Compact integral manifolds of  differentional systems, 
{\it Mathe\-ma\-tics.Dynamical Systems} (arXiv:1009.2998v1 [math.DS]. 
Cornell Univ., Ithaca, New York), 2010, 1-27.
\vspace{0.75ex}

56.
A. Goriely, 
{\it Integrability and nonintegrability of dynamical systems}, 
Advanced series on nonlinear dynamics, Vol. 19,  
World Scientific, 2001. 
\vspace{0.75ex}

57.
C. Christopher, 
Invariant algebraic curves and conditions for a center, 
{\it Proceedings of the Royal Society of Edinburgh}, 1994, Vol. 124A, 1209-1229.
\vspace{0.75ex}

58.
C. Christopher and  J. Llibre, 
Algebraic aspects of integrability for polynomial systems, 
{\it Qualitative theory of dynamical systems}, 1999, Vol. 1, 71-95. 
\vspace{0.75ex}

59.
C. Christopher and  J. Llibre, 
Integrability via invariant algebraic curves for planar polynomial
differential systems, {\it Annals of Differential Equations}, 2000, Vol 16, 5-19.
\vspace{0.75ex}

60.
J. Chavarriga, H. Giacomini, J. Gin\'{e}, and J. Llibre,
Darboux integrability and the inverse integrating factor, 
{\it Differential Equations}, 2003, Vol. 194, 116-139.
\vspace{0.75ex}

61.
J. Llibre, 
{\it Integrability of polynomial differential systems} in
Handbook of Differential Equations: Ordinary Differential Equations, 
Elsevier, Amsterdam, 2004.
\vspace{0.75ex}

62.
A.J. Maciejewski and M. Przybylska,
Darboux polynomials and first integrals of natural polynomial Hamiltonian systems, 
{\it Physics Letters A.}, 2004, Vol. 326, 219-226.
\vspace{0.75ex}
 
63.
A.J. Maciejewski and M. Przybylska,
Darboux points and integrability of Hamiltonian systems with homogeneous polynomial potential, 
{\it Journal of Mathematical Physics}, 2005, 
\linebreak
Vol. 46, No. 062901, 1-33.
\vspace{0.75ex}

64.
F. Dumortier,  J. Llibre, and J.C.  Art\'{e}s, 
{\it Qualitative theory of planar differential systems},
Springer,  New York, 2006.
\vspace{0.75ex}

65.
C. Christopher, J. Llibre, and J.V. Pereira, 
Multiplicity of invariant algebraic curves in polynomial vector fields, 
{\it Pacific Journal of Mathematics}, 2007, Vol. 229, 63-117.
\vspace{0.75ex}

66.
J. Llibre and  Ch. Pantazi, 
Darboux theory of integrability for a class of nonautonomous vector fields, 
{\it Journal of Mathematical Physics}, 2009, No. 50-102705, 1-19.
\vspace{0.75ex}

67.
J. Llibre and X. Zhang, 
Darboux theory of integrability for polynomial vector fields in $R^n$ taking into account the multiplicity at infinity, 
{\it Bulletin des Sciences Mathematiques}, 2009, No. 7(133), 765-778.
\vspace{0.75ex}

68.
J. Llibre and X. Zhang,
Darboux theory of integrability in $C^n$ taking into account the multiplicity, 
{\it Journal of Differential Equations}, 2009, No. 2 (246), 541-551.
\vspace{0.75ex}

69.
A.J.  Maciejewski and M. Przybylska, 
Partial integrability of Hamiltonian systems with homogeneous potential, 
{\it Regular and Chaotic Dynamics}, 2010, Vol. 15, 551-563.
\vspace{0.75ex}

70.
J. Llibre and X. Zhang,
On the Darboux integrability of the polynomial differential systems, 
{\it Qualitative Theory of Dynamical Systems}, 2012, Vol. 11, 129-144.
\vspace{0.75ex}

71.
X. Zhang,
{\it Integrability of dynamical systems: algebra and analysis},   
Springer, Singapore,~2017.
\vspace{0.75ex}

72.
V.V. Kozlov, 
Linear Hamiltonian systems: quadratic integrals, singular subspaces and stability, 
{\it Regular and Chaotic Dynamics}, 2018, Vol. 23, 26-46.  

}

\newpage

{\normalsize

%\Russian
\pagenumbering{arabic}

\sloppy

\lhead
    [\scriptsize В.Н. Горбузов]
    {\scriptsize В.Н. Горбузов}
\rhead
    [\it \scriptsize Частные интегралы обыкновенных дифференциальных систем]
    {\it \scriptsize Частные интегралы обыкновенных дифференциальных систем}

\thispagestyle{empty}

\mbox{}
\\[-0.15ex]
\centerline{
{\large
\bf
ЧАСТНЫЕ\;\! ИНТЕГРАЛЫ\;\! ОБЫКНОВЕННЫХ
}}
\\[0.5ex]
\centerline{
{\large
\bf
ДИФФЕРЕНЦИАЛЬНЫХ\;\! СИСТЕМ
}
}
\\[2.25ex]
\centerline{
\bf 
В.Н. Горбузов
}
\\[2ex]
\centerline{
\it 
Факультет математики и информатики,
}
\\[0.5ex]
\centerline{
\it 
Гродненский государственный университет имени Янки Купалы,
}
\\[0.5ex]
\centerline{
\it 
Ожешко {\rm22}, Гродно, Беларусь, {\rm 230023}
}
\\[1.5ex]
\centerline{
E-mail: gorbuzov@grsu.by
}
\\[4.25ex]
\centerline{{\large\bf Резюме}}
\\[1ex]
\indent
Для обыкновенных дифференциальных систем изучены свойства частных интегралов:
вещественных и комплекснозначных полиномиальных, кратных полиномиальных, 
экспоненциальных и условных.
Рассмотрены возможности построения первых интегралов и последних множителей 
по известным частным интегралам. Даны приложения частных интегралов для решения
задачи Дарбу и расширенной задачи Дарбу, построен интегральный базис системы Якоби,
а также решена обратная задача о построении 
дифференциальных систем на основании их частных интегралов.
\\[1.5ex]
\indent
{\it Ключевые слова}:
дифференциальная система, первый интеграл, частный интеграл, последний множитель.
\\[1.25ex]
\indent
{\it 2000 Mathematics Subject Classification}: 34A34.
%34A34 Nonlinear equations and systems, general
\\[4.25ex]
\centerline{{\large\bf Содержание}}
\\[1.25ex]
{\bf  Введение}                   \dotfill\ 2
\\[1ex]
{\bf \S 1. 
Частные интегралы}
                                                 \dotfill \ 6
\\[0.75ex]
\mbox{}\hspace{1.35em}
1. Частный интеграл. Определение. Свойства
                                                 \dotfill \ 6
\\[0.5ex]
\mbox{}\hspace{1.35em}
2. Полиномиальные частные интегралы
                                                 \dotfill \ 10
\\[0.5ex]
\mbox{}\hspace{1.35em}
3. Экспоненциальные частные интегралы
                                                 \dotfill \ 12
\\[0.5ex]
\mbox{}\hspace{1.35em}
4. Условные частные интегралы
                                                 \dotfill \ 16
\\[0.5ex]
\mbox{}\hspace{1.35em}
5. Кратные полиномиальные частные интегралы
                                                 \dotfill \ 17
\\[0.5ex]
\mbox{}\hspace{1.35em}
6. Комплекснозначные полиномиальные частные интегралы
                                                 \dotfill \ 23
\\[0.5ex]
\mbox{}\hspace{1.35em}
7. Кратные комплекснозначные полиномиальные частные интегралы
                                                 \dotfill \ 27
\\[1ex]
\noindent
{\bf \S 2. 
Последние множители
}
                                                 \dotfill \ 30
\\[0.75ex]
\mbox{}\hspace{1.35em}
8. Последний множитель как частный интеграл
                                                 \dotfill \ 30
\\[0.5ex]
\mbox{}\hspace{1.35em}
9. Построение последних множителей на основании 
\\
\mbox{}\hspace{2.6em}
полиномиальных частных интегралов
                                                 \dotfill \ 33
\\[0.5ex]
\mbox{}\hspace{1.35em}
10. Экспоненциальные последние множители
                                                 \dotfill \ 38
\\[1ex]
\noindent
{\bf \S 3. 
Первые интегралы
}
                                                 \dotfill \ 43
\\[0.75ex]
\mbox{}\hspace{1.35em}
11. Первые интегралы, определяемые частными интегралами 
\\
\mbox{}\hspace{2.9em}
и последними множителями
                                                 \dotfill \ 43
\\[0.5ex]
\mbox{}\hspace{1.35em}
12. Приложения
                                                 \dotfill \ 51
\\[1ex]
{\bf Список литературы}
                                              \dotfill \ 81

\newpage

\mbox{}
\\[-0.25ex]
\centerline{\large\bf  Введение}
\\[1.75ex]
\indent
{\bf Объект исследования и постановка задачи}.
Пусть дана нормальная обыкновенная дифференциальная система $n\!$-го порядка
\\[2ex]
\mbox{}\hfill                             %(0.1)
$
\dfrac{dx_i^{}}{dt}=X_i^{}(t,x),
\quad 
i = 1,\ldots,n,
$
\hfill(0.1)
\\[2.5ex]
у которой правые части 
$X_{i}^{}\colon \Xi\to\R,\ i=1,\ldots,n,\ \Xi = T\times\R^n,\ T\subset\R,$ 
суть полиномы по зависимым переменным $x_1^{},\ldots, x_n^{}$ 
\vspace{0.25ex}
с коэффициентами-функциями одной независимой переменной 
$t$ непрерывно дифференцируемыми на области $T,$ 
\vspace{0.35ex}
имеющие такие степени, что 
$
d=\max\bigl\{\deg_{\;\!x}^{} X_{i}^{}\colon i=1,\ldots, n\bigr\}\geq 1.
$ 
\vspace{0.5ex}

Для системы (0.1) введем понятие {\it частного интеграла}, 
относительно которого решим ряд задач, как-то:
1) существование;
2) аналитическая структура;
3) геометрическая интерпретация;
4) свойства;
5) взаимосвязь с последними множителями;
6) построение первых интегралов по известным частным интегралам.
\vspace{1ex}

{\bf Общие положения}.
С целью однозначного толкования определим используемые понятия, 
сформулируем базовые положения, оговорим принятую терминологию
и введем условные обозначения. При этом в основном будем следовать подходам,
изложенным в монографии [1] и в статье [2]. 

Под {\it областью} будем понимать открытое линейно связное множество. 
Без специальных оговорок будем считать, что области
\vspace{0.25ex}
$\Xi^{\;\!\prime}$ и $\Omega$ являются такими подмножествами области
$\Xi = T\times\R^n,\ T\subset\R,$ что 
\\[1.5ex]
\mbox{}\hfill 
$
\Xi^{\;\!\prime} = T^{\;\!\prime}\times\R^n,\ \ \,
T^{\;\!\prime}\subset T,
$ 
\ а \
$
\Omega = T^{\;\!\prime}\times X^{\;\!\prime},\ \ \,
T^{\;\!\prime}\subset T,\ \ \,
X^{\;\!\prime}\subset\R^n.
\hfill
$ 
\\[1.5ex]
\indent
Множество функций, являющихся  полиномами по переменным $x_1^{},\ldots,x_n^{}$ 
с коэффициентами-функциями одной переменной 
$t$ непрерывно дифференцируемыми на области $T,$ 
обозначим 
$\text{P}_{_{\!\Xi}}.$ 
\vspace{0.35ex}

Функции $p$ и $q$ из множества $\text{P}_{_{\!\Xi}}$ назовем {\it взаимно простыми}, 
\vspace{0.35ex}
если не существует такой 
функции $u\in \text{P}_{_{\!\Xi}},\ u\not\equiv \text{const},$ что
\\[2ex]
\mbox{}\hfill                            
$
\displaystyle
p(t,x)=u(t,x)\;\! p_1^{}(t,x),
\quad 
q(t,x)=u(t,x)\;\! q_1^{}(t,x)
\quad 
\forall (t,x)\in \Xi,
\quad
p_1^{},\, q_1^{}\in \text{P}_{_{\!\Xi}}.
\hfill
$
\\[2ex]
\indent
Система (0.1) индуцирует линейный дифференциальный оператор первого порядка
\\[2ex]
\mbox{}\hfill                            
$
\displaystyle
{\frak d}(t,x) =\partial_{{}_{\scriptstyle t}} +
\sum\limits_{i=1}^n\, X_{i}^{}(t,x)\, \partial_{{}_{\scriptstyle x_i^{}}}
\quad 
\forall (t,x)\in \Xi,
\hfill
$
\\[2ex]
который назовем 
[1, c. 20] {\it оператором дифференцирования в силу системы} (0.1).
\vspace{0.35ex}

Множество функций 
 \vspace{0.35ex}
 непрерывно дифференцируемых на области $\Omega$ обозначим $C^{1}\Omega.$
 
Пусть $f\in C^{1}\Omega.$ Тогда функцию 
\\[2ex]
\mbox{}\hfill                            
$
\displaystyle
{\frak d}f\colon (t,x)\to\ 
\partial_{{}_{\scriptstyle t}}f(t,x) +
\sum\limits_{i=1}^n\, X_{i}^{}(t,x)\, \partial_{{}_{\scriptstyle x_i^{}}}f(t,x)
\quad 
\forall (t,x)\in \Omega
\hfill
$
\\[2ex]
назовем {\it производной в силу системы} (0.1) функции $f;$ 
а если $f(t,x)> 0 \ \forall (t,x)\in \Omega,$ то функцию 
\\[2ex]
\mbox{}\hfill                            
$
\displaystyle
{\frak d}\ln f\colon (t,x)\to\ 
{\frak d}\ln f(t,x)
\quad 
\forall (t,x)\in \Omega
\hfill
$
\\[2ex]
назовем {\it логарифмической производной в силу системы} (0.1) функции $f.$
 \vspace{1ex}
 
{\bf Определение 0.1} [3, с. 256; 4, с. 129 -- 132].
\vspace{0.35ex}
{\it
Функцию $F\in C^{1}\Omega$ назовем 
\textit{\textbf{первым интегралом на области}} $\Omega$ 
системы {\rm (0.1)}, если  
 \vspace{0.35ex}
у нее сохраняется постоянное значение вдоль любого решения 
$x\colon t\to x(t)\;\; \forall t\in T_{_{0}}^{}\subset T^{\;\!\prime}$ 
 \vspace{0.75ex}
 системы {\rm (0.1)} такого, что точки $(t,x(t))\in \Omega\;\; \forall t\in T_{_{0}}^{},$  
т.е.
$F(t,x(t))=C \;\; \forall t\in T_{_{0}}^{},\ C\in\R.$
}
 \vspace{0.5ex}

Множество функций, являющихся первыми интегралами на области $\Omega$ 
системы (0.1), обозначим $\text{I}_{_{\Omega}}\;\!.$
 \vspace{0.15ex}

Укажем критерии существования первого интеграла, которые иногда 
принимаются за определение первого интеграла.
 \vspace{0.5ex}
 
{\bf Теорема 0.1} [1, с. 26; 5, с. 337].
{\it
Функция $F\in C^{1}\Omega$ является первым интегралом на области $\Omega$
системы {\rm (0.1)} тогда и только тогда, когда ее дифференциал в силу системы {\rm (0.1)}
тождественно равен нулю на области $\Omega\colon$}
\\[1.75ex]
\mbox{}\hfill
$
\displaystyle
dF(t,x)_{\displaystyle |_{(0.1)}} = \Bigl(\partial_{{}_{\scriptstyle t}}F(t,x)\;\! dt +
\sum\limits_{i=1}^n\, \partial_{{}_{\scriptstyle x_i^{}}}F(t,x)\;\! dx_i^{}\Bigr)_{\displaystyle |_{(0.1)}} =
\hfill
$
\\
\mbox{}\hfill (0.2)
\\
\mbox{}\hfill
$
\displaystyle
=\Bigl(\partial_{{}_{\scriptstyle t}}F(t,x) +
\sum\limits_{i=1}^n\, X_{i}^{}(t,x)\, \partial_{{}_{\scriptstyle x_i^{}}}F(t,x)\Bigr)dt =0
\quad
\forall (t,x)\in\Omega.
\hfill
$
\\[2ex]
\indent
Учитывая связь дифференциала в силу системы (0.1) с производной в силу системы (0.1),
дифференциальное тождество (0.2) заменим на операторное тождество.
 \vspace{0.5ex}

{\bf Теорема 0.2} [1, с. 26; 2].
{\it
Функция $F\in C^{1}\Omega$ является первым интегралом на области $\Omega$
системы {\rm (0.1)} тогда и только тогда, когда ее производная в силу системы {\rm (0.1)}
тождественно равна нулю на области $\Omega\colon$}
\\[2ex]
\mbox{}\hfill                             %(0.3)
$
{\frak d}F(t,x)=0
\quad 
\forall (t,x)\in \Omega.
$
\hfill(0.3)
\\[2.5ex]
\indent 
{\bf Свойство 0.1} [1, с. 28 -- 29; 3, с. 262].
{\it
Если $F_\nu^{}\in \text{\rm I}_{_{\Omega}}\;\!,\ \nu=1,\ldots,k,$ то функция 
\\[2ex]
\mbox{}\hfill
$
\Psi\colon (t,x)\to\Phi(F_1^{}(t,x),\ldots,F_k^{}(t,x))
\quad
\forall (t,x)\in\Omega,
\hfill
$ 
\\[2ex]
где $\Phi$ --- произвольная непрерывно дифференцируемая функция, также 
будет первым интегралом на области $\Omega$ системы}~(0.1).
 \vspace{0.5ex}
 
Этим свойством при $k=1$ выражается функциональная неоднозначность первого интеграла.
 \vspace{0.5ex}
 
{\bf Определение 0.2} [1, с. 29].
{\it
Совокупность функционально независимых на области $\Omega$ первых интегралов 
$F_\nu^{}\colon \Omega\to\R,\ \nu=1,\ldots,k,$ системы {\rm (0.1)} назовем 
\textit{\textbf{базисом первых интегралов}} {\rm(}или 
\textit{\textbf{интегральным базисом}}{\rm)} 
\textit{\textbf{на области}} $\Omega$ системы {\rm (0.1)}, если у системы {\rm (0.1)}
любой первый интеграл $F\colon \Omega\to\R$ на области $\Omega$ 
можно представить в виде 
\\[0.5ex]
\mbox{}\hfill
$
F(t,x)=\Phi(F_1^{}(t,x),\ldots,F_k^{}(t,x))
\quad
\forall (t,x)\in\Omega,
\hfill
$ 
\\[1.75ex]
где $\Phi$ --- некоторая непрерывно дифференцируемая функция. 
Число $k$ при этом назовем \textit{\textbf{размерностью}} базиса первых интегралов
на области $\Omega$ системы} (0.1).
\vspace{0.5ex} 
 
 {\bf Теорема 0.3} [1, с. 54; 3, с. 367].
{\it
Система {\rm (0.1)} на окрестности любой точки из области $\Omega$ 
имеет базис первых интегралов размерности $n$
}
\vspace{0.5ex} 
 
{\bf Определение 0.3} [2].
{\it
Гладкое многообразие ${\rm g}(t,x)=0$ назовем 
\textit{\textbf{интегральным многообразием}} системы {\rm (0.1)}, 
если дифференциал функции ${\rm g}\in C^{1}\Omega$ 
 в силу системы {\rm(0.1)} тождественно равен нулю на многообразии ${\rm g}(t,x)=0\colon$
\\[2ex]
\mbox{}\hfill                           
$
d\;\!{\rm g}(t,x)_{\displaystyle |_{(0.1)}} = \Phi(t,x)\;\! dt
\quad 
\forall (t,x)\in \Omega,
$
\hfill {\rm (0.4)}
\\[2ex]
где функция $\Phi\colon \Omega \to \R$ такая, что}
\\[1.5ex]
\mbox{}\hfill                           
$
\Phi(t,x)_{\displaystyle |_{{\rm g}(t,x)=0}} = 0
\quad 
\forall (t,x)\in \Omega.
\hfill
$
\\[-2.25ex]
\indent

\newpage

Наряду с определением 0.3 будем использовать критерий существования 
интегрального многообразия у системы (0.1).
\vspace{0.5ex}

{\bf Теорема 0.4.}
{\it
Гладкое  многообразие   ${\rm g}(t,x) = 0$ 
является  интегpальным  многообразием системы {\rm(0.1)}
тогда и только тогда, когда
производная в силу системы {\rm(0.1)} функции ${\rm g}\colon \Omega\to \R$ 
тождественно равна нулю на этом многообразии}:
\\[2ex]
\mbox{}\hfill                                    
$
{\frak d}\;\! {\rm g}(t,x)  =  \Phi(t,x), 
\quad 
\Phi(t,x)_{\displaystyle |_{{\rm g}(t,x)=0}} = 0
\quad 
\forall (t,x)\in \Omega.
$
\hfill  (0.5)
\\[2.25ex]
\indent
{\bf Определение 0.4.}
{\it
Функцию $\mu\in C^{1}\Omega$ назовем 
\textit{\textbf{последним множителем на области}} $\Omega$  системы {\rm (0.1)}, 
если дифференциал в силу системы {\rm(0.1)} 
\\[2ex]
\mbox{}\hfill                     
$
d\;\!\mu(t,x)_{\displaystyle |_{(0.1)}}={}-\mu(t,x)\;{\rm div}\, {\frak d}(t,x)\, dt
\quad 
\forall (t,x)\in \Omega,
\hfill
$
\\[1ex]
где}
\\[0ex]
\mbox{}\hfill
$
\displaystyle
{\rm div}\, {\frak d}(t,x) = 
\sum\limits_{i=1}^n\, \partial_{{}_{\scriptstyle x_i^{}}}X_i^{}(t,x)
\quad 
\forall (t,x)\in \Xi.
\hfill
$
\\[2ex]
\indent
Приведем критерий существования последнего множителя, который часто используют
в качестве определения последнего множителя.
\vspace{0.5ex}

{\bf Теорема 0.5} [5, с. 341 -- 346; 6, с. 117].
{\it
Функция $\mu\in C^{1}\Omega$ является 
последним множителем на области $\Omega$  системы {\rm (0.1)}
тогда и только тогда, когда ее производная в силу системы} {\rm(0.1)} 
\\[2ex]
\mbox{}\hfill                         %(0.6)
$
{\frak d}\;\!\mu(t,x)={}-\mu(t,x)\;{\rm div}\, {\frak d}(t,x)
\quad 
\forall (t,x)\in \Omega.
$
\hfill(0.6)
\\[2.25ex]
\indent
Множество функций, являющихся последними множителями на области $\Omega$ 
системы (0.1), обозначим $\text{M}_{_{\Omega}}\;\!.$
 \vspace{0.5ex}

{\bf Cвойство 0.2} (свойство Якоби последних множителей).
 \vspace{0.35ex}
 {\it 
Пусть $\mu_1^{},\mu_2^{}\in \text{\rm M}_{_{\Omega}}\;\!,$ 
множество $\Omega_{_0}\subset \Omega$ такое, что 
$\mu_2^{}(t,x)\ne 0\;\; \forall (t,x)\in \Omega_{_0},\
\mu_2^{}(t,x)= 0\;\; \forall (t,x)\in {\sf C}_{_\Omega}\Omega_{_0}.$
Тогда функция
\\[2ex]
\mbox{}\hfill
$
F \colon (t,x) \to\ 
\dfrac{\mu_1^{}(t,x)}{\mu_2^{}(t,x)}
\quad  
\forall(t,x)\in \Omega_{_0}
\hfill
$
\\[2.25ex]
будет первым интегралом на любой области из множества $\Omega_{_0}$ системы} (0.1).
\vspace{0.5ex}
 
{\bf Замечание 0.1.}
Определения 0.1 -- 0.4, теоремы 0.1 -- 0.5, а также свойства 0.1 и 0.2,
справедливы и в более общем случае, когда у системы (0.1) правые части 
$X_i^{}\in C^{1}\Omega,$ $i=1,\ldots, n,$
но не обязательно должны быть полиномами по зависимым переменным 
$x_1^{},\ldots, x_n^{}.$
\vspace{0.75ex}

Чтобы упростить чтение, перечислим принятые условные обозначения:
\vspace{0.5ex}

\noindent
$T$ --- область из $\R;$ 
$T^{\;\!\prime}$ --- подобласть области $T;$
\vspace{0.75ex}

\noindent
область $\Xi = T\times\R^n;$
область $\Xi^{\;\!\prime} = T^{\;\!\prime}\times\R^n;$
\vspace{0.75ex}

\noindent
$X^{\;\!\prime}$ --- область из $\R^n;$
область $\Omega = T^{\;\!\prime}\times X^{\;\!\prime};$
\vspace{0.75ex}

\noindent
${\sf C}_{_{\Omega^{\;\!\prime}}}\Omega_{_0}$ ---
дополнение множества $\Omega_{_0}\subset \Omega^{\;\!\prime}$ 
до множества $\Omega^{\;\!\prime};$
\vspace{0.75ex}

\noindent
${\sf D}$ --- символ Коши производной функции одной переменной;
\vspace{0.75ex}

\noindent
$\partial_{t}^{}$ --- частная производная по переменной $t$ 
(оператор дифференцирования 
\\
\mbox{}\hspace{2.7em}
по переменной $t);$
\vspace{0.75ex}

\noindent
${\frak d}$ --- оператор дифференцирования в силу системы (0.1);
\vspace{0.75ex}

\noindent
$d\;\!f_{\displaystyle |_{(0.1)}}$ ---
дифференциал функции $f$ в силу системы (0.1);
\vspace{1ex}

\noindent
$
d=\max\bigl\{\deg_{\;\!x}^{} X_{i}^{}\colon i=1,\ldots, n\bigr\}\;\!;
$ 
\vspace{1ex}

\noindent
$C^1\Omega$ --- множество непрерывно дифференцируемых функций на области $\Omega\;\!;$
\vspace{0.5ex}

\newpage

\noindent
$\text{P}_{_{\!\Xi}}$ --- 
множество функций, являющихся  полиномами по переменным $x_1^{},\ldots,x_n^{}$ 
\\
\mbox{}\hspace{2.7em}
с коэффициентами-функциями (вещественными) одной переменной $t$ непрерывно 
\\
\mbox{}\hspace{2.7em}
дифференцируемыми на области $T\;\!;$ 
\vspace{0.5ex}

\noindent
$\text{Z}_{_{\Xi}}$ --- 
множество функций, являющихся  полиномами по переменным $x_1^{},\ldots,x_n^{}$ 
\\
\mbox{}\hspace{2.7em}
с комплекснозначными коэффициентами-функциями одной переменной $t$ 
\\
\mbox{}\hspace{2.7em}
непрерывно дифференцируемыми на области $T\;\!;$ 
\vspace{0.5ex}

\noindent
$\text{I}_{_{\Omega}}$ ---
множество первых интегралов на области $\Omega$ системы (0.1); 
\vspace{0.5ex}

\noindent
$\text{J}_{_{\Omega}}$ ---
множество частных интегралов на области $\Omega$ системы (0.1); 
\vspace{0.75ex}

\noindent
${\rm(g}, M)\in \text{J}_{_{\Omega}}$ ---
функция ${\rm g}$ является частным интегралом с сомножителем $M$ на области 
\\[0.5ex]
\mbox{}\hspace{6.3em}
$\Omega$ системы (0.1); 
\vspace{0.5ex}

\noindent
$\text{A}_{_{\Xi}}$ ---
множество полиномиальных (вещественных) частных интегралов на области $\Xi$ 
\\
\mbox{}\hspace{2.8em}
системы (0.1); 
\vspace{0.75ex}

\noindent
$(p, M)\in \text{A}_{_{\Xi}}$ ---
функция $p$ является полиномиальным (вещественным) частным    
\\[0.35ex]
\mbox{}\hspace{6.8em}
интегралом с сомножителем $M$ на области $\Xi$ системы (0.1); 
\vspace{0.5ex}

\noindent
$\text{B}_{_{\Xi}}$ ---
множество кратных полиномиальных (вещественных) частных интегралов  
\\
\mbox{}\hspace{2.75em}
на области $\Xi$ системы (0.1); 
\vspace{0.5ex}

\noindent
$\bigl((p, M), (h,q, N)\bigr)\in \text{B}_{_{\Xi}}$ ---
полиномиальный (вещественный) частный интеграл $p$ 
\\[0.35ex]
\mbox{}\hspace{11.4em}
с сомножителем $M$ на области $\Xi$ системы (0.1) является 
\\[0.35ex]
\mbox{}\hspace{11.3em}
кратным таким,  что выполняется тождество (5.1); 
\vspace{0.5ex}

\noindent
$\text{E}_{_{\Omega}}$ ---
множество экспоненциальных частных интегралов на области $\Omega$ 
\vspace{0.75ex}
системы (0.1); 

\noindent
$(\exp\omega , M)\in \text{E}_{_{\Omega}}$ ---
функция $\exp\omega$ является экспоненциальным частным интегралом с   
\\[0.35ex]
\mbox{}\hspace{8.7em}
сомножителем $M$ на области $\Omega$ системы (0.1); 
\vspace{0.5ex}

\noindent
$\text{F}_{_{\Xi}}$ ---
множество условных частных интегралов на области $\Xi$ 
системы (0.1); 
\vspace{0.75ex}

\noindent
$(\exp p, M)\in \text{F}_{_{\Xi}}$ ---
функция $\exp p$ является условным частным интегралом с   
\\[0.35ex]
\mbox{}\hspace{8.6em}
сомножителем $M$ на области $\Xi$ системы (0.1); 
\vspace{0.5ex}

\noindent
$\text{H}_{_{\Xi}}$ ---
множество комплекснозначных полиномиальных частных интегралов  
\\
\mbox{}\hspace{2.8em}
на области $\Xi$ системы (0.1); 
\vspace{0.75ex}

\noindent
$(w, W)\in \text{H}_{_{\Xi}}$ ---
функция $w$ является комплекснозначным полиномиальным частным    
\\[0.35ex]
\mbox{}\hspace{7.25em}
интегралом с сомножителем $W$ на области $\Xi$ системы (0.1); 
\vspace{0.5ex}

\noindent
$\text{G}_{_{\Xi}}$ ---
множество кратных комплекснозначных полиномиальных частных интегралов  
\\
\mbox{}\hspace{2.8em}
на области $\Xi$ системы (0.1); 
\vspace{0.75ex}

\noindent
$\bigl((w, W), (h,z, Q)\bigr)\in \text{G}_{_{\Xi}}$ ---
комплекснозначный полиномиальный частный интеграл $w$ 
\\[0.35ex]
\mbox{}\hspace{11.7em}
с сомножителем $W$ на области $\Xi$ системы (0.1) является   
\\[0.15ex]
\mbox{}\hspace{11.7em}
кратным таким, что выполняется тождество (7.1); 
\vspace{0.5ex}

\noindent
$\text{M}_{_{\Omega}}$ ---
множество последних множителей на области $\Omega$ системы (0.1); 
\vspace{0.75ex}

\noindent
$\text{MA}_{_{\Xi}}$ ---
множество полиномиальных (вещественных) последних множителей 
\\
\mbox{}\hspace{3.7em}
на области $\Xi$ системы (0.1); 
\vspace{0.75ex}

\noindent
$\text{MB}_{_{\Xi}}$ ---
множество кратных полиномиальных (вещественных) последних множителей  
\\
\mbox{}\hspace{3.7em}
на области $\Xi$ системы (0.1); 
\vspace{0.5ex}

\noindent
$\text{ME}_{_{\Omega}}$ ---
множество экспоненциальных последних множителей на области $\Omega$ 
\\
\mbox{}\hspace{3.7em}
системы (0.1); 
\vspace{0.5ex}

\noindent
$\text{MF}_{_{\Xi}}$ ---
множество условных последних множителей на области $\Xi$ 
системы (0.1); 
\vspace{0.75ex}

\noindent
$\text{MH}_{_{\Xi}}$ ---
множество комплекснозначных полиномиальных последних множителей 
\\
\mbox{}\hspace{3.7em}
на области $\Xi$ системы (0.1); 
\vspace{0.5ex}

\noindent
$\text{MG}_{_{\Xi}}$ ---
множество кратных комплекснозначных полиномиальных последних 
\\
\mbox{}\hspace{3.7em}
множителей на области $\Xi$ системы (0.1).
\vspace{0.5ex}

\newpage

\mbox{}
\\[-0.5ex]
\centerline{
{\bf\large \S\;\!1. Частные интегралы}}
\\[2ex]
\centerline{
{\bf  1. 
Частный интеграл. Определение. Свойства
}
}
\\[1.5ex]
\indent
{\bf Определение 1.1.}
\vspace{0.5ex}
{\it
Непрерывно дифференцируемую функцию
$
{\rm g}\colon \Omega\to\R
$
назовем 
\textit{\textbf{частным интегралом на области}} 
$\Omega$ сис\-темы {\rm (0.1)}, если ее дифференциал в силу системы {\rm (0.1)} 
\\[2ex]
\mbox{}\hfill                                           % (1.1)
$
\displaystyle
d\!\;{\rm g}(t,x)_{\displaystyle |_{(0.1)}} =
{\rm g}(t,x)\;\!M(t,x)\;\!dt
\quad 
\forall (t,x)\in \Omega,
$
\hfill {\rm(1.1)}
\\[2.25ex]
где функция
\vspace{0.25ex}
$M\in \text{\rm P}_{_{\!\Xi^{\;\!\prime}}}$ и имеет степень 
$\deg_{\;\!x}^{} M\leq d-1.$ 
При этом функцию $M$ будем называть 
\textit{\textbf{сомножителем частного интеграла}} ${\rm g}.$
}
\vspace{0.5ex}

Множество функций, являющихся частными интегралами на 
области $\Omega$ системы (0.1), обозначим 
$\text{J}_{_{\Omega}}.$ 
\vspace{0.25ex}
Оборот слов <<функция ${\rm g}$ является частным интегралом
с сомножителем $M$ на области $\Omega$ сис\-темы {\rm (0.1)}>> выразим 
условной записью $({\rm g}, M)\in \text{J}_{_{\Omega}}.$
\vspace{0.35ex}
Если множество $\Omega_{_0}\subset \R^{n+1}$ не является областью, то 
запись $({\rm g}, M)\in \text{J}_{_{\Omega_{{}_{\tiny 0}}}}$ 
\vspace{0.25ex}
будем трактовать как то, что функция ${\rm g}$ 
\vspace{0.25ex}
является частным интегралом с сомножителем $M$на любой области 
из множества $\Omega_{_0}$ сис\-темы {\rm (0.1)}. 

С помощью оператора дифференцирования в силу системы (0.1) 
дифференциальное тождество (1.1) запишем в виде одного из операторных тождеств
\\[2ex]
\mbox{}\hfill                             %(1.2)
$
\displaystyle
{\frak d}\;\! {\rm g}(t,x) =
{\rm g}(t,x)\;\!M(t,x)
\quad 
\forall (t,x)\in \Omega
$
\hfill(1.2)
\\[1ex]
или
\\[1ex]
\mbox{}\hfill                             %(1.3)
$
\displaystyle
{\frak d} \ln \bigl| {\rm g}(t,x) \bigr| = M(t,x)
\quad 
\forall (t,x)\in \Omega_{_0}\;\!,
$
\hfill(1.3)
\\[2.5ex]
где множество $\Omega_{_0}\subset \Omega$ такое, что  
${\rm g}(t,x)\ne 0\;\;\forall (t,x)\in \Omega_{_0},\ 
{\rm g}(t,x)=0\;\;\forall (t,x)\in {\sf C}_{_\Omega}\Omega_{_0}.$
\vspace{0.5ex}

Тем самым, получаем два критерия существования  частного интеграла.
\vspace{0.5ex}

{\bf Теорема 1.1.}
\vspace{0.35ex}
${\rm (g}, M)\in \text{J}_{_{\Omega}}$
{\it 
тогда и только тогда, когда выполняется тождество {\rm (1.2)},
а также тогда и только тогда, когда выполняется тождество {\rm (1.3)}.
\vspace{0.25ex}
При этом как в тождестве {\rm (1.2)}, так и в тождестве {\rm (1.3)}, функция
\vspace{0.25ex}
$M\in \text{\rm P}_{_{\!\Xi^{\;\!\prime}}}$ и имеет степень 
$\deg_{\;\!x}^{} M\leq d-1.$ 
}
\vspace{0.5ex}

Основываясь на определении интегрального многообразия (определение 0.3)
и определении частного интеграла (определение 1.1), устанавливаем 
геометрический смысл частного интеграла.
\vspace{0.5ex}

{\bf Теорема 1.2.}
{\it 
Если частный интеграл ${\rm g}$ системы {\rm (0.1)} 
определяет многообразие ${\rm g}(t,x)=0,$ 
то оно будет интегральным многообразием системы {\rm (0.1)}.
}
\vspace{0.5ex}

{\bf Теорема 1.3.}
{\it
Непрерывно дифференцируемая функция
\\[1.5ex]
\mbox{}\hfill
$
\displaystyle
{\rm g}\colon (t,x)\to\ 
\prod\limits_{j=1}^{m} {\rm g}_j^{}(t,x)
\quad
\forall (t,x)\in\Omega
$
\hfill {\rm (1.4)}
\\[1.5ex]
являет\-ся частным интегралом  с 
\vspace{0.25ex}
сомножителем $M$ на области $\Omega$ системы {\rm (0.1)}
тогда и только тогда, когда существуют функции 
$M_j^{}\in C^1\Omega_{_0}\;\!,\ j=1,\ldots,m,$ 
такие, что выполняются тождества
\\[1.5ex]
\mbox{}\hfill
$
\displaystyle
{\frak d}\;\!{\rm g}_j^{}(t,x)  = 
{\rm g}_j^{}(t,x)\;\! M_j^{}(t,x)
\quad
\forall (t,x)\in\Omega_{_0}\;\!,
\quad
j=1,\ldots,m,
$
\hfill {\rm (1.5)}
\\[1ex]
и 
\\[1ex]
\mbox{}\hfill
$
\displaystyle
\sum\limits_{j=1}^{m} M_j^{}(t,x)= M(t,x)
\quad
\forall (t,x)\in\Omega_{_0}\;\!,
$
\hfill {\rm (1.6)}
\\[1.5ex]
а функция
$M\in \text{\rm P}_{_{\!\Xi^{\;\!\prime}}}$ и имеет степень 
$\deg_{\;\!x}^{} M\leq d-1.$ 
}
\vspace{0.5ex}

{\sl Доказательство. Необходимость}.
При $m=2$ для $({\rm g}_1^{}{\rm g}_2^{}, M)\in 
 \text{\rm J}_{_{\Omega}}$
тождество (1.2) из теоремы 1.1 имеет вид
\\[1.5ex]
\mbox{}\hfill                          
$
\displaystyle
{\rm g}_2^{}(t,x)\;\! {\frak d}\;\! {\rm g}_1^{}(t,x) +
{\rm g}_1^{}(t,x)\;\! {\frak d}\;\! {\rm g}_2^{}(t,x) =
{\rm g}_1^{}(t,x)\;\! {\rm g}_2^{}(t,x)\;\!M(t,x)
\quad 
\forall (t,x)\in\Omega,
\hfill
$
\\[2ex]
где функция
$M\in \text{\rm P}_{_{\!\Xi^{\;\!\prime}}}$ и  
$\deg_{\;\!x}^{} M\leq d-1,$
а значит, выполняется тождество 
\\[1.5ex]
\mbox{}\hfill                          
$
\displaystyle
{\frak d}\;\! {\rm g}_1^{}(t,x) =
{\rm g}_1^{}(t,x)\;\! 
\biggl(
M(t,x)-\dfrac{{\frak d}\;\! {\rm g}_2^{}(t,x)}{{\rm g}_2^{}(t,x)}
\biggr)
\quad 
\forall (t,x)\in\Omega_{_0},
\hfill
$
\\[2ex]
где множество $\Omega_{_0}$ такое, что  
${\rm g}_2^{}(t,x)\ne 0\;\;\forall (t,x)\in \Omega_{_0},\ 
{\rm g}_2^{}(t,x)=0\;\;\forall (t,x)\in {\sf C}_{_\Omega}\Omega_{_0}.$
\vspace{0.5ex}

Пусть
\\[1.5ex]
\mbox{}\hfill                            
$
M(t,x)-\dfrac{{\frak d}\;\! {\rm g}_2^{}(t,x)}{{\rm g}_2^{}(t,x)}=
M_1^{}(t,x)
\quad 
\forall (t,x)\in \Omega_{_0}.
\hfill 
$
\\[1.5ex]
\indent
Тогда производная в силу системы (0.1)
\\[1.5ex]
\mbox{}\hfill                          
$
\displaystyle
{\frak d}\;\! {\rm g}_1^{}(t,x)=
{\rm g}_1^{}(t,x)\;\! M_1^{}(t,x)
\quad 
\forall (t,x)\in\Omega_{_0}\;\!,
\hfill
$
\\[2.5ex]
\mbox{}\hfill                          
$
\displaystyle
{\frak d}\;\! {\rm g}_2^{}(t,x)=
{\rm g}_2^{}(t,x)\;\! \bigl(M(t,x)-M_1^{}(t,x)\bigr)
\quad 
\forall (t,x)\in\Omega_{_0}\;\!.
\hfill
$
\\[2ex]
\indent
Итак, тождества (1.5) и (1.6) при $m=2$ доказаны.

При $m>2$ тождества (1.5) и (1.6) доказываются по индукции.
\vspace{0.25ex}

{\sl Достаточность.}
Пусть выполняются тождества (1.5) и (1.6) при условии, что 
$M\in \text{\rm P}_{_{\!\Xi^{\;\!\prime}}}$ и имеет степень  
$\deg_{\;\!x}^{} M\leq d-1.$ Тогда
\\[2ex]
\mbox{}\hfill                        
$
\displaystyle
{\frak d}\;\!{\rm g}(t,x)  =
{\frak d}\;\!\prod\limits_{j=1}^{m} {\rm g}_j^{}(t,x) \, = \,
\sum\limits_{\nu=1}^{m}\, 
\prod\limits_{{\,}_{\scriptstyle j\ne \nu}^{\scriptstyle j=1,}}^{m}
{\rm g}_j^{}(t,x)\, {\frak d}\;\! {\rm g}_{\nu}^{}(t,x)  
=
\hfill
$
\\[2ex]
\mbox{}\hfill                        
$
\displaystyle
=\,
\prod\limits_{j=1}^{m}\, 
{\rm g}_j^{}(t,x)\, 
\sum\limits_{\nu=1}^{m} M_{\nu}^{}(t,x)=
{\rm g}(t,x)\;\! M(t,x)
\quad
\forall (t,x)\in \Omega.
\hfill
$
\\[1.5ex]
\indent
По теореме 1.1, функция (1.4) является частным интегралом 
с сомножителем $M$ на области $\Omega$ системы (0.1). $\k$
\vspace{0.75ex}

{\bf  Свойство 1.1.}
\vspace{0.25ex}
{\it
Пусть функция $\varphi\in C^1T^{\;\!\prime},$ 
множество $T_{_0}\subset T^{\;\!\prime}$ 
такое, что 
$\varphi(t)\ne 0$ $\forall t\in T_{_0},\ 
\varphi(t)=0\;\;\forall t\in {\sf C}_{{}_{T^{\;\!\prime}}}T_{_0},$
множество $\Omega_{_0}=T_{_0}\times X^{\;\!\prime}.$ Тогда}
\\[2ex]
\mbox{}\hfill
$
({\rm g}, M)\in \text{J}_{_{\Omega}}
\iff 
\bigl(\varphi\;\! {\rm g}\;\!,\;\! M+{\sf D} \ln |\varphi|\bigr)\in 
\text{J}_{_{\Omega_{{}_{\tiny\;\! 0}}}}.
\hfill
$
\\[2ex]
\indent
{\sl Доказательство} основано на теореме 1.1 и том, что  
\\[1.5ex]
\mbox{}\hfill                        
$
{\frak d}\bigl(\varphi(t)\;\!{\rm g}(t,x)\bigr)=
{\rm g}(t,x)\, {\sf D} \varphi(t) +\varphi(t)\, {\frak d}\;\!{\rm g}(t,x) 
=
\hfill
$
\\[2ex]
\mbox{}\hfill
$
=
\varphi(t)\;\! \bigl({\rm g}(t,x)\, {\sf D} \ln |\varphi(t) | + 
 {\frak d}\;\!{\rm g}(t,x)\bigr)
\quad
\forall (t,x)\in  \Omega_{_0}.\ \k
\hfill 
$
\\[2ex]
\indent
Частным случаем свойства 1.1 является 
\vspace{0.5ex}

{\bf Свойство 1.2.} 
{\it
Если $\lambda\in\R\backslash\{0\},$ то}
\\[1.5ex]
\mbox{}\hfill
$
({\rm g}, M)\in \text{J}_{_{\Omega}}
\iff 
(\lambda\;\! {\rm g}\;\!,\;\! M)\in \text{J}_{_{\Omega}}.
\hfill
$
\\[2ex]
\indent
В соответствии со свойством 1.2, говоря о двух и более частных интегралах 
системы (0.1), будем считать их попарно линейно независимыми.
\vspace{0.5ex}

{\bf Свойство 1.3.} 
\vspace{0.25ex}
{\it
Пусть множество $\Omega_{_0}\subset\Omega$ такое, что  
${\rm g}(t,x)\ne 0\;\;\forall (t,x)\in \Omega_{_0},$ 
${\rm g}(t,x)=0\;\;\forall (t,x)\in {\sf C}_{_\Omega}\Omega_{_0}.$
Тогда}
\\[1.5ex]
\mbox{}\hfill  % (1.7)
$
({\rm g}, M)\in \text{J}_{_{\Omega}}
\iff 
(|{\rm g}|\;\!,\;\! M)\in \text{J}_{_{\Omega_{{}_{\tiny\;\! 0}}}}.
$
\hfill (1.7)
\\[2ex]
\indent
{\sl Доказательство.}
Поскольку
\\[1.5ex]
\mbox{}\hfill 
$
({\rm g}, M)\in \text{J}_{_{\Omega}}
\iff 
({\rm g}, M)\in \text{J}_{_{\Omega_{{}_{\tiny\;\! 0}}}},
\hfill
$
\\[2ex]
\mbox{}\hfill 
$
|{\rm g}(t,x)|=\text{sgn}\;\!{\rm g}(t,x)\, {\rm g}(t,x)
\quad
\forall (t,x)\in\Omega,
\hfill
$
\\[1.5ex]
а, по свойству 1.2,
\\[1.5ex]
\mbox{}\hfill 
$
({\rm g}, M)\in \text{J}_{_{\Omega_{{}_{\tiny\;\! 0}}}}
\iff 
({\rm sgn\;\!g\,g}, M)\in \text{J}_{_{\Omega_{{}_{\tiny\;\! 0}}}},
\hfill
$
\\[1.5ex]
то в соответствии с транзитивностью эквиваленции справедливо 
утверждение (1.7). $\k$
\vspace{0.5ex}

Посредством свойства 1.3 устанавливается связь между тождествами (1.2) и (1.3).
\vspace{0.5ex}

{\bf Свойство 1.4.} 
{\it
Если производная в силу системы {\rm (0.1)}
\\[1.5ex]
\mbox{}\hfill
$
{\frak d}\;\!{\rm g}(t,x)=
\bigl({\rm g}(t,x)+c\bigr)\;\!M(t,x)
\quad
\forall (t,x)\in\Omega,
\hfill
$
\\[2ex]
где $c\in\R,\ M\in \text{\rm P}_{_{\!\Xi^{\;\!\prime}}},\  
\deg_{\;\!x}^{} M\leq d-1,$ то 
$({\rm g}+c,M)\in \text{\rm J}_{_{\Omega}}.$}
\vspace{0.5ex}

{\sl Действительно}, производная в силу системы (0.1)
\\[1.5ex]
\mbox{}\hfill
$
{\frak d}\;\!\bigl({\rm g}(t,x)+c\bigr)=
{\frak d}\;\!{\rm g}(t,x)=
\bigl({\rm g}(t,x)+c\bigr)\;\!M(t,x)
\quad
\forall (t,x)\in\Omega.
\hfill
$
\\[1.5ex]
\indent
Следовательно, по теореме 1.1, $({\rm g}+c,M)\in \text{\rm J}_{_{\Omega}}.\ \k$
\vspace{0.75ex}

{\bf Свойство 1.5.} 
\vspace{0.5ex}
{\it
Если $({\rm g}_j^{}, M)\in \text{\rm J}_{_{\Omega}},\ 
\lambda_j^{}\in\R\backslash\{0\},\ j=1,\ldots,m,$ то 
$\biggl(\,\sum\limits_{j=1}^{m}\! \lambda_j^{}\;\!{\rm g}_j^{}, M\biggl) 
\in \text{\rm J}_{_{\Omega}}.$}

{\sl Действительно}, если 
$({\rm g}_j^{}, M)\in \text{\rm J}_{_{\Omega}}\;\!,\ j=1,\ldots,m,$
то, по теореме 1.1,
\\[1.5ex]
\mbox{}\hfill
$
{\frak d}\;\!{\rm g}_j^{}(t,x)=
{\rm g}_j^{}(t,x)\;\!M(t,x)
\quad
\forall (t,x)\in\Omega,
\quad
j=1,\ldots,m,
\hfill
$
\\[1.5ex]
сомножитель $M\in \text{\rm P}_{_{\!\Xi^{\;\!\prime}}}$ и имеет степень 
$\deg_{\;\!x}^{} M\leq d-1.$
\vspace{0.35ex}

Тогда производная в силу системы (0.1)
\\[1.5ex]
\mbox{}\hfill                         
$
\displaystyle
{\frak d}\;\!\sum\limits_{j=1}^{m} \lambda_j^{}\;\!{\rm g}_j^{}(t,x)  = 
\sum\limits_{j=1}^{m}
\lambda_j^{}\;\!{\rm g}_j^{}(t,x)\;\!M(t,x)
\quad
\forall (t,x)\in\Omega.
\hfill
$
\\[1.5ex]
\indent
Стало быть, по теореме 1.1, 
$\biggl(\,\sum\limits_{j=1}^{m}\! \lambda_j^{}\;\!{\rm g}_j^{}, M\biggl) 
\in \text{\rm J}_{_{\Omega}}.\ \k$
\vspace{1.25ex}

{\bf Свойство 1.6.} 
{\it
Пусть $\gamma\in\R\backslash\{0\},\ {\rm g}^{\gamma}\in C^1\Omega.$
Тогда}
\\[1.5ex]
\mbox{}\hfill    %(1.8)
$
({\rm g}, M)\in \text{J}_{_{\Omega}}
\iff 
\bigl({\rm g}^{\gamma},\;\! \gamma\;\!M\bigr)\in \text{J}_{_{\Omega}}.
$
\hfill (1.8)
\\[2ex]
\indent
{\sl Доказательство}
основано на теореме 1.1, примененной к функциям ${\rm g}$ 
и ${\rm g}^{\gamma},$ и том, что имеет место тождество
\\[1ex]
\mbox{}\hfill
$
{\frak d}\;\!{\rm g}^{\gamma}(t,x)=
\gamma\;\!{\rm g}^{\gamma-1}(t,x)\,
{\frak d}\;\!{\rm g}(t,x)
\quad
\forall (t,x)\in\Omega. \ \k
\hfill
$
\\[2.25ex]
\indent
{\bf Свойство 1.7.} 
\vspace{0.5ex}
{\it
Пусть $\gamma\in\R\backslash\{0\},$ 
множество $\Omega_{_0}\subset\Omega$ такое, что  
${\rm g}(t,x)\ne 0$ $\forall (t,x)\in \Omega_{_0},\ 
{\rm g}(t,x)=0\;\;\forall (t,x)\in {\sf C}_{_\Omega}\Omega_{_0}.$
Тогда}
\\[2ex]
\mbox{}\hfill  % (1.9)
$
({\rm g}, M)\in \text{J}_{_{\Omega}}
\iff 
\bigl(|{\rm g}|^{\gamma},\;\! \gamma\;\!M\bigr)\in 
\text{J}_{_{\Omega_{{}_{\tiny\;\! 0}}}}.
$
\hfill (1.9)
\\[2ex]
\indent
{\sl Доказательство.}
По свойству 1.3, справедлива эквиваленция (1.7).
По свойству 1.6,
\\[1.5ex]
\mbox{}\hfill 
$
(|{\rm g}|, M)\in \text{J}_{_{\Omega_{{}_{\tiny\;\! 0}}}}
\iff 
\bigl(|{\rm g}|^{\gamma},\;\! \gamma\;\!M\bigr)\in 
\text{J}_{_{\Omega_{{}_{\tiny\;\! 0}}}}.
\hfill
$
\\[1.5ex]
\indent
Использовав транзитивность эквиваленции, получаем утверждение (1.9). $\k$ 
\vspace{0.75ex}

{\bf Свойство 1.8.} 
\vspace{0.5ex}
{\it
Если $\rho_j^{}, \lambda_j^{}\in\R\backslash\{0\},
\ ({\rm g}_j^{}, \rho_j^{}M_{0}^{})\in \text{\rm J}_{_{\Omega}},\ 
{\rm g}_j^{{}^{\scriptsize 1/\rho_{\!j}^{}}}\in C^1\Omega, ,\ j=1,\ldots,m,$ то 
$\biggl(\,\sum\limits_{j=1}^{m}\! 
\lambda_j^{}\;\!{\rm g}_j^{{}^{\scriptsize 1/\rho_{\!j}^{}}}, M_0^{}\biggl) 
\in \text{\rm J}_{_{\Omega}}.$}
%\vspace{0.5ex}

{\sl Доказательство.} 
Так как  
$({\rm g}_j^{}, \rho_j^{}M_0^{})\in \text{\rm J}_{_{\Omega}},$
то, по свойству 1.6, 
$\Bigl({\rm g}_j^{{}^{\scriptsize 1/\rho_{\!j}^{}}}, M_0^{}\Bigr)\in 
\text{\rm J}_{_{\Omega}}, 
\linebreak 
j=1,\ldots,m.$
Тогда на основании свойства 1.5 получаем, что
$\biggl(\,\sum\limits_{j=1}^{m}\! 
\lambda_j^{}\;\!{\rm g}_j^{{}^{\scriptsize 1/\rho_{\!j}^{}}}, M_0^{}\biggl) 
\in \text{\rm J}_{_{\Omega}}.\ \k$
\vspace{0.75ex}

{\bf Свойство 1.9.} 
\vspace{0.5ex}
{\it
Пусть 
$({\rm g}_j^{}, M_{j}^{})\in \text{\rm J}_{_{\Omega}},\ 
\gamma_j^{}\in\R\backslash\{0\},\
{\rm g}_j^{{}^{\scriptsize \gamma_{j}^{}}}\in C^1\Omega, \ j=1,\ldots,m.$ 
Тогда  
$\biggl(\,\prod\limits_{j=1}^{m} 
{\rm g}_j^{{}^{\scriptsize \gamma_{j}^{}}}, M\biggl) 
\in \text{\rm J}_{_{\Omega}}\;\!,$
если и только если сомножители $M,\;\!M_1^{},\ldots,M_m^{}$
такие, что выполняется тождество}
\\[1ex]
\mbox{}\hfill                      % (1.10)
$
\displaystyle
M(t,x)=\sum\limits_{j=1}^{m}
\gamma_j^{}\;\!M_{j}^{}(t,x)
\quad
\forall (t,x)\in\Xi^{\;\!\prime}.
$
\hfill (1.10)
\\[1.5ex]
\indent
{\sl Доказательство.} 
Так как
$({\rm g}_j^{}, M_{j}^{})\in \text{\rm J}_{_{\Omega}},\ j=1,\ldots,m,$
то согласно свойству 1.6 и тео\-ре\-ме 1.1 выполняются тождества 
\\[1.75ex]
\mbox{}\hfill   % (1.11)
$
{\frak d}\;\!{\rm g}_j^{{}^{\scriptsize \gamma_{j}^{}}}\!(t,x)=
\gamma_j^{}\;\!{\rm g}_j^{{}^{\scriptsize \gamma_{j}^{}}}\!(t,x)\;\!M_j^{}(t,x)
\quad
\forall (t,x)\in\Omega,
\quad
j=1,\ldots,m,
$
\hfill (1.11)
\\[2ex]
в которых функции $M_j^{}\in \text{\rm P}_{_{\!\Xi^{\;\!\prime}}}$ и имеют степени 
$\deg_{\;\!x}^{} M_j^{}\leq d-1,\ j=1,\ldots, m.$
\vspace{0.5ex}

По теореме 1.3, 
\vspace{0.35ex}
$\biggl(\,\prod\limits_{j=1}^{m} 
{\rm g}_j^{{}^{\scriptsize \gamma_{j}^{}}}, M\biggl) 
\in \text{\rm J}_{_{\Omega}}$
тогда и только тогда, когда существуют такие функции
$M_j^{\ast}\in C^1\Omega,\ j=1,\ldots, m,$
что выполняются тождества
\\[1.75ex]
\mbox{}\hfill   % (1.12)
$
{\frak d}\;\!{\rm g}_j^{{}^{\scriptsize \gamma_{j}^{}}}\!(t,x)=
{\rm g}_j^{{}^{\scriptsize \gamma_{j}^{}}}\!(t,x)\;\!M_j^{\ast}(t,x)
\quad
\forall (t,x)\in\Omega,
\quad
j=1,\ldots,m,
$
\hfill (1.12)
\\[1ex]
и
\\[1ex]
\mbox{}\hfill                      
$
\displaystyle
M(t,x)=\sum\limits_{j=1}^{m}
M_{j}^{\ast}(t,x)
\quad
\forall (t,x)\in\Omega,
\hfill
$
\\[1.5ex]
а функция $M\in \text{\rm P}_{_{\!\Xi^{\;\!\prime}}}$ и имеет степень 
$\deg_{\;\!x}^{} M\leq d-1.$
\vspace{0.5ex}

Из тождеств (1.12) с учетом тождеств (1.11) следует, что
\\[1.5ex]
\mbox{}\hfill                      
$
\displaystyle
M_{j}^{\ast}(t,x)=
\gamma_j^{}\;\!M_j^{}(t,x)
\quad
\forall (t,x)\in\Xi^{\;\!\prime},
\quad
j=1,\ldots,m,
\hfill
$
\\[1.75ex]
причем функции $M_j^{}\in \text{\rm P}_{_{\!\Xi^{\;\!\prime}}}$ и имеют степени 
$\deg_{\;\!x}^{} M_j^{}\leq d-1,\ j=1,\ldots, m.$
\vspace{0.5ex}

А значит,
\vspace{0.5ex}
$\biggl(\,\prod\limits_{j=1}^{m} 
{\rm g}_j^{{}^{\scriptsize \gamma_{j}^{}}}, M\biggl) 
\in \text{\rm J}_{_{\Omega}},$
если и только если
%тогда и только тогда, когда 
\vspace{0.5ex}
имеет место тождество (1.10). $\k$

{\bf Следствие 1.1.} 
\vspace{0.5ex}
{\it
Пусть 
$\rho_j^{},\gamma_j^{}\in\R\backslash\{0\},\
({\rm g}_j^{}, \rho_j^{}M_{0}^{})\in \text{\rm J}_{_{\Omega}},\ 
{\rm g}_j^{{}^{\scriptsize \gamma_{j}^{}}}\in C^1\Omega, \ j=1,\ldots,m.$ 
Тог\-да  
$\biggl(\,\prod\limits_{j=1}^{m} 
{\rm g}_j^{{}^{\scriptsize \gamma_{j}^{}}}, M\biggl) 
\in \text{\rm J}_{_{\Omega}}\;\!,$
если и только если сомножитель}
\\[1ex]
\mbox{}\hfill                      % (1.13)
$
\displaystyle
M(t,x)=\sum\limits_{j=1}^{m}
\rho_j^{}\;\!\gamma_j^{}\;\!M_{0}^{}(t,x)
\quad
\forall (t,x)\in\Xi^{\;\!\prime}.
$
\hfill (1.13)
\\[1.5ex]
\indent
На основании теоремы 1.3, учитывая свойства 1.6 и 1.9, получаем
\vspace{0.75ex}

{\bf Свойство 1.10.} 
\vspace{0.5ex}
{\it
Пусть 
$({\rm g}_\tau^{}, M_{\tau}^{})\in \text{\rm J}_{_{\Omega}},\ 
\tau=1,\ldots, m-1,\
\gamma_j^{}\in\R\backslash\{0\},\
{\rm g}_j^{{}^{\scriptsize \gamma_{j}^{}}}\in C^1\Omega, 
\linebreak 
j=1,\ldots,m.$ 
Тогда}
\\[1ex]
\mbox{}\hfill                      
$
\displaystyle
\biggl(\,\prod\limits_{j=1}^{m} 
{\rm g}_j^{{}^{\scriptsize \gamma_{j}^{}}}, M\biggl) 
\in \text{\rm J}_{_{\Omega}}
\iff
\biggl({\rm g}_m^{}, \,
\dfrac{1}{\gamma_m^{}}\;\!
\biggl(M-\sum\limits_{\tau=1}^{m-1}
\gamma_{\tau}\;\!M_{\tau}^{}\biggl)\biggl) 
\in \text{\rm J}_{_{\Omega}}\;\!.
\hfill
$
\\[1.5ex]
\indent
На основании теоремы 1.3 и свойства 1.10 получаем
\vspace{0.75ex}

{\bf Свойство 1.11.} 
\vspace{0.5ex}
{\it
Пусть 
$({\rm g}_\nu^{}, M_{\nu}^{})\in \text{\rm J}_{_{\Omega}},\ 
\nu=1,\ldots, s,\ s\leq m-2, \ \
\gamma_j^{}\in\R\backslash\{0\},\
{\rm g}_j^{{}^{\scriptsize \gamma_{j}^{}}}\!\in C^1\Omega, 
\linebreak 
j=1,\ldots,m.$ 
Тогда
\\[1ex]
\mbox{}\hfill                      
$
\displaystyle
\biggl(\,\prod\limits_{j=1}^{m} 
{\rm g}_j^{{}^{\scriptsize \gamma_{j}^{}}}, M\biggl) 
\in \text{\rm J}_{_{\Omega}}
\iff
\biggl(
\,\prod\limits_{k=s+1}^{m} 
{\rm g}_k^{{}^{\scriptsize \gamma_{k}^{}}}, \,
M-\sum\limits_{\nu=1}^{s}
\gamma_{\nu}\;\!M_{\nu}^{}\biggl) 
\in \text{\rm J}_{_{\Omega}}\;\!.
\hfill
$
\\[1.5ex]
\indent
Кроме этого существуют такие функции 
$M_k^{}\in C^1\Omega,\ k=s+1,\ldots, m,$
что выполняются тождества
\\[1.5ex]
\mbox{}\hfill                        % (1.14)
$
{\frak d}\;\!{\rm g}_k^{}(t,x)=
{\rm g}_k^{}(t,x)\;\!M_k^{}(t,x)
\quad
\forall (t,x)\in\Omega,
\quad
k=s+1,\ldots,m,
$
\hfill {\rm (1.14)}
\\[1ex]
и}
\\[1ex]
\mbox{}\hfill                      % (1.15)
$
\displaystyle
\sum\limits_{k=s+1}^{m}
\gamma_{k}^{}\;\!M_k^{}(t,x)=
M(t,x)-\sum\limits_{\nu=1}^{s}
\gamma_{\nu}^{}\;\!M_\nu^{}(t,x)
\quad
\forall (t,x)\in\Omega.
$
\hfill (1.15)
\\[2ex]
\indent
{\bf Замечание 1.1.}
Свойства 1.8 -- 1.11 и следствие 1.1 доказаны на основании 
эквиваленции (1.8) из свойства 1.6. 
При необходимости в свойствах 1.8 -- 1.11 и в следствии 1.1
степень ${\rm g}_j^{{}^{\scriptsize \gamma_{j}^{}}},\ j\in\{1,\ldots,m\},$
\vspace{0.35ex}
может быть заменена на степень
$|{\rm g}_j^{}|^{{}^{\scriptsize \gamma_{j}^{}}}$
с корректировкой на свойство 1.7.
\\[2.75ex]
\centerline{
{\bf  2. Полиномиальные частные интегралы}
}
\\[1.5ex]
\indent
{\bf Определение 2.1.}
{\it
Частный интеграл на области $\Xi^{\;\!\prime}$ сис\-темы {\rm (0.1)}, 
являющийся полиномом по переменным $x_1^{},\ldots, x_n^{}$
с коэффициентами-функциями одной переменной $t$ непрерывно
дифференцируемыми на области $T^{\;\!\prime},$ назовем
\textit{\textbf{полиномиальным частным интегралом на области}}
$\Xi^{\;\!\prime}$ сис\-темы {\rm (0.1)}.
}
\vspace{0.5ex}

Множество функций, являющихся полиномиальными частными интегралами 
на области $\Xi^{\;\!\prime}$ сис\-темы (0.1), обозначим 
$\text{A}_{_{\Xi^{\;\!\prime}}}.$
\vspace{0.5ex}

Согласно определению 2.1 множество
$\text{A}_{_{\Xi^{\;\!\prime}}}\subset \text{J}_{_{\Xi^{\;\!\prime}}}.$
\vspace{0.5ex}

Условной записью 
\vspace{0.35ex}
$(p, M)\in \text{A}_{_{\Xi^{\;\!\prime}}}$ 
будем выражать, что функция $p\in \text{P}_{_{\!\Xi^{\;\!\prime}}}$
является полиномиальным частным интегралом с сомножителем $M$
\vspace{0.25ex}
на области $\Xi^{\;\!\prime}$ сис\-темы (0.1).

В соответствии с определением 2.1
\\[1.5ex]
\mbox{}\hfill               % (2.1)
$
(p,M)\in \text{A}_{_{\Xi^{\;\!\prime}}}
\iff
(p,M)\in \text{J}_{_{\Xi^{\;\!\prime}}}
\ \&\ 
p\in \text{P}_{_{\!\Xi^{\;\!\prime}}}.
$
\hfill (2.1)
\\[2ex]
\indent
{\bf Теорема 2.1}
\vspace{0.35ex}
(критерий существования полиномиального частного интеграла).
{\it
$(p, M)\in \text{\rm A}_{_{\Xi^{\;\!\prime}}}$ 
тогда и только тогда, когда  производная в силу системы} (0.1)
\\[1.5ex]
\mbox{}\hfill               % (2.2)
$
{\frak d}\;\!p(t,x)=
p(t,x)\;\!M(t,x)
\quad
\forall (t,x)\in \Xi^{\;\!\prime},
\quad
p\in \text{P}_{_{\!\Xi^{\;\!\prime}}}.
$
\hfill (2.2)
\\[1.5ex]
\indent
{\sl Следует} с учетом эквиваленции (2.1) из критерия существования 
частного интеграла (теорема 1.1). $\k$
\vspace{0.5ex}

{\bf Замечание 2.1.}
\vspace{0.25ex}
Если функция $p\in \text{P}_{_{\!\Xi^{\;\!\prime}}},$
то из тождества (2.2) следует, что функция
$M\in \text{P}_{_{\!\Xi^{\;\!\prime}}}$ и имеет степень
$\deg_{\;\!x}^{} M\leq d-1.$
\vspace{0.5ex}

{\bf Теорема 2.2.}\! 
\vspace{0.5ex}
{\it
Пусть 
$p_j^{}\!\in\! \text{\rm P}_{_{\!\Xi^{\;\!\prime}}},\, 
\gamma_j^{}\!\in\!\R\backslash\{0\},\!$ 
множество $\Omega_{_0}\!\subset\!\Xi^{\;\!\prime}\!$ такое, что  
$p_j^{{}^{\scriptsize \gamma_j^{}}}\!\in\! C^1\Omega_{_0},$ $j=1,\ldots, m.$
Тогда
\\[1.5ex]
\mbox{}\hfill  
$
\displaystyle
\biggl(\,\prod\limits_{j=1}^{m} 
p_j^{{}^{\scriptsize \gamma_{j}^{}}}, M\biggl) 
\in \text{\rm J}_{_{\Omega_{{}_{\tiny\;\! 0}}}}
\iff
\bigl(p_j^{}, M_j^{}\bigr)\in \text{\rm A}_{_{\Xi^{\;\!\prime}}},
\ \ 
j=1,\ldots, m,
\hfill
$
\\[1.5ex]
где сомножители $M$ и $M_j^{},\ j=1,\ldots,m,$
\vspace{0.5ex}
такие, что выполняется тождество {\rm (1.10)}.}

{\sl Доказательство.}
\vspace{0.5ex}
По теореме 1.3,
$
\biggl(\,\prod\limits_{j=1}^{m} 
p_j^{{}^{\scriptsize \gamma_{j}^{}}}, M\biggl) 
\in \text{J}_{_{\Omega_{{}_{\tiny\;\! 0}}}}
$
тогда и только тогда, когда существуют такие функции 
 $M_j^{\ast}\in C^1\Omega_{_0},\ j=1,\ldots,m,$
что выполняются тождества
\\[1.75ex]
\mbox{}\hfill   % (2.3)
$
{\frak d}\;\!p_j^{{}^{\scriptsize \gamma_{j}^{}}}\!(t,x)=
p_j^{{}^{\scriptsize \gamma_{j}^{}}}\!(t,x)\;\!M_j^{\ast}(t,x)
\quad
\forall (t,x)\in\Omega_{_{0}},
\quad
j=1,\ldots,m,
$
\hfill (2.3)
\\[1ex]
и
\\[1ex]
\mbox{}\hfill     % (2.4)                 
$
\displaystyle
\sum\limits_{j=1}^{m}
M_{j}^{\ast}(t,x)= M(t,x)
\quad
\forall (t,x)\in \Omega_{_{0}},
$
\hfill (2.4)
\\[1.5ex]
а функция $M\in \text{\rm P}_{_{\!\Xi^{\;\!\prime}}}$ и имеет степень 
$\deg_{\;\!x}^{} M\leq d-1.$
\vspace{0.5ex}

Поскольку  производная в силу системы (0.1)
\\[1.5ex]
\mbox{}\hfill 
$
{\frak d}\;\!p_j^{{}^{\scriptsize \gamma_{j}^{}}}\!(t,x)=
\gamma_j^{}\;\!p_j^{{}^{\scriptsize \gamma_{j}^{}-1}}\!(t,x)\;\!
{\frak d}\;\!p_j^{{}}\!(t,x)
\quad
\forall (t,x)\in\Omega_{_{0}},
\quad
j=1,\ldots,m,
\hfill
$
\\[2ex]
то тождества (2.3) имеют место тогда и только тогда, когда
\\[1.5ex]
\mbox{}\hfill 
$
{\frak d}\;\!p_j^{}(t,x)=
p_j^{}(t,x)\;\!M_j^{}(t,x)
\quad
\forall (t,x)\in \Xi^{\;\!\prime},
\quad
j=1,\ldots,m,
\hfill
$
\\[1ex]
где
\\[1ex]
\mbox{}\hfill 
$
M_j^{}(t,x)=
\dfrac{1}{\gamma_j^{}}\, M_j^{\ast}(t,x)
\quad
\forall (t,x)\in \Xi^{\;\!\prime},
\quad
j=1,\ldots,m.
\hfill
$
\\[1.5ex]
\indent
Отсюда, учитывая теорему 2.1, получаем, что тождества (2.3) имеют 
место тогда и только тогда, когда 
$
\bigl(p_j^{}, M_j^{}\bigr)\in \text{\rm A}_{_{\Xi^{\;\!\prime}}},
\ 
j=1,\ldots, m.
$
\vspace{0.5ex}

Тождество (1.10) следует из тождества (2.4). $\k$
\vspace{0.5ex}

Частным случаем теоремы 2.2 является
\vspace{0.5ex}

{\bf Теорема 2.3}
(критерий существования рационального частного интеграла). 
{\it
Пусть функции
$p_1^{},p_2^{}\in \text{\rm P}_{_{\!\Xi^{\;\!\prime}}}$
взаимно простые, 
\vspace{0.75ex}
множество $\Omega_{_0}\subset\Xi^{\;\!\prime}$ такое, что  
$p_2^{}(t,x)\ne 0$ $\forall (t,x)\in \Omega_{_0},\
p_2^{}(t,x)= 0\;\; \forall (t,x)\in {\sf C}_{_{\Xi^{\;\!\prime}}}\Omega_{_0}.$
Тогда
\\[1.5ex]
\mbox{}\hfill  
$
\Bigl(\,\dfrac{p_1^{}}{p_2^{}}\,,\, M\Bigl) 
\in \text{\rm J}_{_{\Omega_{{}_{\tiny\;\! 0}}}}
\iff
\bigl(p_1^{}, M_1^{}\bigr)\in \text{\rm A}_{_{\Xi^{\;\!\prime}}}
\ \& \ 
\bigl(p_2^{}, M_2^{}\bigr)\in \text{\rm A}_{_{\Xi^{\;\!\prime}}}.
\hfill
$
\\[1.5ex]
При этом сомножители $M,\ M_1^{},\ M_2^{}$ такие, что}
\\[1.5ex]
\mbox{}\hfill  
$
M(t,x)=M_1^{}(t,x)-M_2^{}(t,x)
\quad
\forall (t,x)\in \Xi^{\;\!\prime}.
\hfill
$
\\[2ex]
\indent
{\bf  Свойство 2.1.}
\vspace{0.5ex}
{\it
Если $\varphi\in C^1T^{\;\!\prime},$ 
множество $T_{_0}\subset T^{\;\!\prime}$ 
такое, что 
$\varphi(t)\ne 0\;\; \forall t\in T_{_0},$ 
$\varphi(t)=0\;\;\forall t\in {\sf C}_{{}_{T^{\;\!\prime}}}T_{_0},$
множество $\Xi_{_0}=T_{_0}\times \R^n,$ то}
$
\bigl(\varphi, {\sf D} \ln |\varphi|\bigr)\in 
\text{\rm A}_{_{\Xi_{{}_{\tiny\;\! 0}}}}.
$
\vspace{0.25ex}

{\sl В самом деле}, производная в силу системы (0.1)  
\\[1.5ex]
\mbox{}\hfill                        
$
{\frak d}\;\!\varphi(t) =
{\sf D} \varphi(t) =
\varphi(t)\,\dfrac{{\sf D} \varphi(t)}{\varphi(t)}=
\varphi(t)\;\! {\sf D} \ln |\varphi(t) | 
\quad
\forall (t,x)\in  \Xi_{_0},
\hfill 
$
\\[1.5ex]
причем $\deg_{\;\!x}^{} {\sf D} \ln |\varphi |=0\leq d-1.$
\vspace{0.5ex}

По теореме 2.1, 
$
\bigl(\varphi, {\sf D} \ln |\varphi|\bigr)\in 
\text{\rm A}_{_{\Xi_{{}_{\tiny\;\! 0}}}}\!.\,  \k
$
\vspace{0.75ex}

{\bf  Свойство 2.2.}
\vspace{0.5ex}
{\it
Пусть $\lambda_j^{},\;\! c_j^{}\in\R, \ j=1,\ldots,m,\ m\leq n,\ 
\sum\limits_{j=1}^{m}|\lambda_j^{}|\ne 0.$
Тогда 
$
\biggl(\,\sum\limits_{j=1}^{m} 
\lambda_j^{}(x_j^{}+c_j^{}), M\biggl) 
\in\text{\rm A}_{_{\Xi}},
$
если и только если выполняется тождество}
\\[1.5ex]
\mbox{}\hfill
$
\displaystyle
\sum\limits_{j=1}^{m}
\lambda_j^{}\;\!X_j^{}(t,x)=
M(t,x)\;\!\sum\limits_{j=1}^{m}\lambda_j^{}(x_j^{}+c_j^{})
\quad
\forall (t,x)\in \Xi.
\hfill
$
\\[1.5ex]
\indent
{\sl Доказательство}
основано на теореме 2.1 и состоит в том, что производная в силу системы (0.1)
\\[1.5ex]
\mbox{}\hfill                        
$
\displaystyle
{\frak d}\;\!\sum\limits_{j=1}^{m}\lambda_j^{}(x_j^{}+c_j^{}) =
\sum\limits_{j=1}^{m}\lambda_j^{}\;\!X_j^{}(t,x)
\quad
\forall (t,x)\in  \Xi.
\ \k
\hfill 
$
\\[2ex]
\indent
Из свойства 2.2 при $m=1$ получаем
\vspace{0.5ex}

\newpage

{\bf  Свойство 2.3.}
\vspace{0.5ex}
{\it
Если $c\in\R,$ то
$
\bigl(x_i^{}+c, M\bigr) \in \text{\rm A}_{_{\Xi}},\ 
i\in \{1,\ldots, n\},
$
тогда и только тогда, когда в правой части системы {\rm (0.1)} функция}
\\[1ex]
\mbox{}\hfill
$
X_i^{}(t,x)=
(x_i^{}+c)\;\! M(t,x)
\quad
\forall (t,x)\in \Xi.
\hfill
$
\\[2.75ex]
\centerline{
{\bf  3. Экспоненциальные частные интегралы}
}
\\[1.5ex]
\indent
{\bf Определение 3.1.}
\vspace{0.25ex}
{\it
Если $\bigl(\exp\omega,M\bigr)\in \text{\rm J}_{_{\Omega}},$
то функцию $\exp\omega$ назовем
\textit{\textbf{экспоненциальным частным интегралом}}
\vspace{0.5ex}
с сомножителем $M$ на области $\Omega$ сис\-темы {\rm (0.1)}.}

Множество функций, являющихся экспоненциальными частными интегралами 
на области $\Omega$ сис\-темы (0.1), обозначим 
$\text{E}_{_{\Omega}}.$
\vspace{0.5ex}

В соответствии с определением 3.1 
$\text{E}_{_{\Omega}}\subset \text{J}_{_{\Omega}}.$
\vspace{0.5ex}

Условной записью 
$\bigl(\exp\omega,M\bigr)\in \text{\rm E}_{_{\Omega}}$
будем выражать, что функция $\exp\omega$ является 
экспоненциальным частным интегралом  
\vspace{0.75ex}
с сомножителем $\!M\!$ на области $\!\Omega\!$ сис\-темы (0.1).

{\bf Теорема 3.1}
\vspace{0.35ex}
(критерий существования экспоненциального частного интеграла).
{\it
$\bigl(\exp\omega,M\bigr)\in \text{\rm E}_{_{\Omega}}$ 
тогда и только тогда, когда выполняется тождество
\\[1.5ex]
\mbox{}\hfill               % (3.1)
$
{\frak d}\;\!\omega(t,x)=M(t,x)
\quad
\forall (t,x)\in \Omega,
$
\hfill {\rm(3.1)}
\\[1.5ex]
в котором функция $M\in \text{\rm P}_{_{\!\Xi^{\;\!\prime}}}$ 
и имеет степень $\deg_{\;\!x}^{} M\leq d-1.$
}
\vspace{0.5ex}

{\sl Доказательство} 
основано на критерии существования частного интеграла (теорема~1.1),
определении экспоненциального частного интеграла (определение 3.1)
и состоит в том, что для функции 
\vspace{0.35ex}
${\rm g}\colon (t,x)\to\exp\omega\;\;\forall (t,x)\in\Omega$
тождество (1.3), равно как и тождество (1.2), имеет вид (3.1). $\k$
\vspace{0.75ex}

{\bf Теорема 3.2}
\vspace{0.15ex}
(критерий существования экспоненциального частного интеграла).
{\it
Функция $\exp\omega\in \text{\rm E}_{_{\Omega}}$ 
\vspace{0.5ex}
тогда и только тогда, когда функция 
${\frak d}\;\!\omega\in \text{\rm P}_{_{\!\Xi^{\;\!\prime}}}$ 
и имеет степень $\deg_{\;\!x}^{} {\frak d}\;\!\omega\leq d-1.$
\vspace{0.35ex}
При этом функция ${\frak d}\;\!\omega$ является сомножителем 
экспоненциального частного интеграла $\exp\omega.$
}
\vspace{0.75ex}

{\bf  Теорема 3.3.}
\vspace{0.5ex}
{\it
Пусть $\gamma_j^{}\in\R\backslash\{0\},\ c_j^{}\in\R, \ j=1,\ldots,m.$
Тогда
\\[1.5ex]
\mbox{}\hfill
$
\displaystyle
\exp\sum\limits_{j=1}^{m}
\gamma_j^{}\;\!(\omega_j^{}+c_j^{})\in \text{\rm E}_{_{\Omega}}
\ \iff \
\sum\limits_{j=1}^{m}\gamma_j^{}\;\!{\frak d}\;\!\omega_j^{}\in 
\text{\rm P}_{_{\!\Xi^{\;\!\prime}}}
\, \ \&\ \,
\deg_{\;\!x}^{} 
\sum\limits_{j=1}^{m}
\gamma_j^{}\;\!{\frak d}\;\!\omega_j^{}\leq d-1.
\hfill
$
\\[1.5ex]
При этом функция 
$\sum\limits_{j=1}^{m}\gamma_j^{}\;\!{\frak d}\;\!\omega_j^{}$
является сомножителем экспоненциального частного интеграла 
$\exp\sum\limits_{j=1}^{m}\gamma_j^{}\;\!(\omega_j^{}+c_j^{}).$
}

{\sl Следует} из теоремы 3.2 при 
$\omega=\exp\sum\limits_{j=1}^{m}\gamma_j^{}\;\!(\omega_j^{}+c_j^{}).\ \k$
\vspace{0.5ex}

Теорема 3.3 является аналогом теоремы 1.3 в случае, когда 
${\rm g}_j^{}=\exp\bigl(\gamma_j^{}\;\!(\omega_j^{}+c_j^{})\bigr),
\linebreak 
j=1,\ldots,m.$
\vspace{0.75ex}

{\bf  Свойство 3.1.}
{\it
Пусть $\varphi\in C^1T^{\;\!\prime}.$ Тогда} 
\\[1.5ex]
\mbox{}\hfill
$
\displaystyle
({\rm g}, M)\in \text{\rm J}_{_{\Omega}}
\iff
\bigl({\rm g}\exp\varphi, M+{\sf D} \varphi\bigr)\in \text{\rm J}_{_{\Omega}}.
\hfill
$
\\[1.75ex]
\indent
{\sl Доказательство} основано на теореме 1.1 и том, что  
\\[1.5ex]
\mbox{}\hfill                        
$
{\frak d}\;\!\bigl({\rm g}(t,x)\;\!e^{\varphi(t)}\bigr) =
e^{\varphi(t)}\;\!{\frak d}\;\!{\rm g}(t,x) +
{\rm g}(t,x)\;\!{\sf D}e^{\varphi(t)} =
e^{\varphi(t)}\;\!{\rm g}(t,x)\;\!
\bigl(M(t,x)+{\sf D} \varphi(t)\bigr)
\quad
\forall (t,x)\in  \Omega. \ \k
\hfill 
$
\\[2ex]
\indent
{\bf  Свойство 3.2.}
{\it
Пусть $\varphi\in C^1T^{\;\!\prime}.$ Тогда} 
\\[1.5ex]
\mbox{}\hfill
$
\displaystyle
(\exp\omega, M)\in \text{\rm E}_{_{\Omega}}
\iff
\bigl(\exp(\omega+\varphi), M+{\sf D} \varphi\bigr)\in \text{\rm E}_{_{\Omega}}.
\hfill
$
\\[1.75ex]
\indent
{\sl Следует} из свойства 3.1 при ${\rm g}=\exp\omega.\ \k$  
\vspace{0.5ex}

Частным случаем свойства 3.2 является
\vspace{0.5ex}

\newpage

{\bf  Свойство 3.3.}
{\it
Если $c\in\R,$ то} 
\\[1.5ex]
\mbox{}\hfill
$
\displaystyle
(\exp\omega, M)\in \text{\rm E}_{_{\Omega}}
\iff
\bigl(\exp(\omega+c), M\bigr)\in \text{\rm E}_{_{\Omega}}.
\hfill
$
\\[1.75ex]
\indent
{\bf  Свойство 3.4.}
{\it
Если 
$\bigl(\exp\omega_j^{}, M\bigr)\in \text{\rm E}_{_{\Omega}},\ 
\lambda_j^{}\in\R\backslash\{0\}, \ j=1,\ldots,m,$ то}
\\[1.5ex]
\mbox{}\hfill
$
\displaystyle
\biggl(\,\sum\limits_{j=1}^{m}\lambda_j^{}\exp\omega_j^{}\;\!,\;\! M\biggr)\in \text{\rm J}_{_{\Omega}}.
\hfill
$
\\[1.5ex]
\indent
{\sl Следует} из свойства 1.5. $\k$  
\vspace{1ex}

{\bf  Свойство 3.5.}
{\it
Если 
$\rho_j^{},\lambda_j^{}\in\R\backslash\{0\}, \ 
\bigl(\exp\omega_j^{}, \rho_j^{}\;\!M_0^{}\bigr)\in \text{\rm E}_{_{\Omega}},\ 
j=1,\ldots,m,$ то}
\\[1.5ex]
\mbox{}\hfill
$
\displaystyle
\biggl(\,\sum\limits_{j=1}^{m}\lambda_j^{}
\exp\dfrac{\omega_j^{}}{\rho_j^{}}\;\!,\;\! M_0^{}\biggr)\in \text{\rm J}_{_{\Omega}}.
\hfill
$
\\[1.5ex]
\indent
{\sl Следует} из свойства 1.8. $\k$  
\vspace{1ex}

{\bf  Свойство 3.6.}
\vspace{0.5ex}
{\it
Пусть  
$\bigl(\exp\omega_j^{}, M_j^{}\bigr)\in \text{\rm E}_{_{\Omega}},\ 
\gamma_j^{}\in\R\backslash\{0\}, \ j=1,\ldots,m.$ Тогда
$
\biggl(\,\exp\sum\limits_{j=1}^{m}\gamma_j^{}\;\!\omega_j^{}\;\!,\;\! M\biggr)\in \text{\rm E}_{_{\Omega}},
$
если и только если выполняется тождество} (1.10).
\vspace{0.75ex}

{\sl Следует} из свойства 1.9. $\k$  
\vspace{1ex}

{\bf  Следствие 3.1.}
{\it
Пусть 
$\gamma_j^{}, \rho_j^{}\in\R\backslash\{0\}, \ 
\bigl(\exp\omega_j^{}, \rho_j^{}\;\!M_0^{}\bigr)\in \text{\rm E}_{_{\Omega}},\ 
j=1,\ldots,m.$ Тогда
$
\biggl(\,\exp\sum\limits_{j=1}^{m}
\gamma_j^{}\;\!\omega_j^{}\;\!,\;\! M\biggr)\in \text{\rm E}_{_{\Omega}},
$
если и только если выполняется тождество} (1.13).
\vspace{1ex}

{\bf  Следствие 3.2.}
{\it
Если $\gamma\in\R\backslash\{0\},$ то} 
\\[1.5ex]
\mbox{}\hfill
$
\displaystyle
(\exp\omega, M)\in \text{\rm E}_{_{\Omega}}
\iff
\bigl(\exp(\gamma\;\!\omega), \gamma\;\!M\bigr)\in \text{\rm E}_{_{\Omega}}.
\hfill
$
\\[1.75ex]
\indent
Следствие 3.2 является аналогом свойств 1.6 и 1.7 на случай экспоненциального 
частного интеграла.
\vspace{1ex}

{\bf  Свойство 3.7.}
{\it
Если 
$
\bigl(\exp\omega_j^{}, M_j^{}\bigr)\in \text{\rm E}_{_{\Omega}},\ 
\gamma,\gamma_j^{}\in\R\backslash\{0\}, \ j=1,\ldots,m, \
{\rm g}^{\gamma}\in C^1\Omega,$ то}
\\[1.5ex]
\mbox{}\hfill
$
\displaystyle
\biggl(\,
{\rm g}^{\gamma}\exp\sum\limits_{j=1}^{m}\gamma_j^{}\;\!\omega_j^{}\;\!,\;\! M\biggr)\in 
\text{\rm J}_{_{\Omega}}
\iff
\biggl(\,
{\rm g},\, \dfrac{1}{\gamma}\,\biggl(M-
\sum\limits_{j=1}^{m}\gamma_j^{}\;\!M_j^{}\biggr)\biggr)\in 
\text{\rm J}_{_{\Omega}}.
\hfill
$
\\[1.5ex]
\indent
{\sl Следует} из свойства 1.10. $\k$  
\vspace{1ex}

{\bf  Следствие 3.3.}
{\it
Если
$
\bigl(\exp\omega_\tau^{}, M_\tau^{}\bigr)\!\in\! \text{\rm E}_{_{\Omega}},\,
\tau\!=\!1,\ldots, m\!-\!1,\ 
\gamma_j^{}\in\R\backslash\{0\}, \, j\!=\!1,\ldots,m,$ то}
\\[1.5ex]
\mbox{}\hfill
$
\displaystyle
\biggl(\,
\exp\sum\limits_{j=1}^{m}\gamma_j^{}\;\!\omega_j^{}\;\!,\;\! M\biggr)\in 
\text{\rm E}_{_{\Omega}}
\iff
\biggl(\,
\exp\omega_{m}^{}\;\!,\,
 \dfrac{1}{\gamma_{m}^{}}\,\biggl(M-
\sum\limits_{\tau=1}^{m-1}\gamma_\tau^{}\;\!M_\tau^{}\biggr)\biggr)\in 
\text{\rm E}_{_{\Omega}}.
\hfill
$
\\[2ex]
\indent
{\bf  Свойство 3.8.}
\vspace{0.5ex}
{\it
Пусть 
$
\bigl(\exp\omega_\nu^{}, M_\nu^{}\bigr)\in \text{\rm E}_{_{\Omega}},\ 
\nu=1,\ldots, s,\ s\leqslant m-2,\  
\gamma_j^{}\in\R\backslash\{0\}, 
\linebreak 
j=1,\ldots,m,\ 
{\rm g}_k^{{}^{\scriptsize \gamma_k^{}}}\!\in C^1\Omega,\ 
k=s+1,\ldots,m.$ Тогда
\\[1.75ex]
\mbox{}\hfill
$
\displaystyle
\biggl(\,
\prod\limits_{k=s+1}^{m}
{\rm g}_k^{{}^{\scriptsize \gamma_k^{}}}
\exp\sum\limits_{\nu=1}^{s}
\gamma_\nu^{}\;\!\omega_\nu^{}\;\!,\;\! M\biggr)\in 
\text{\rm J}_{_{\Omega}}
\iff
\biggl(\,
\prod\limits_{k=s+1}^{m}
{\rm g}_k^{{}^{\scriptsize \gamma_k^{}}},\,
M-\sum\limits_{\nu=1}^{s}\gamma_\nu^{}\;\!M_\nu^{}\biggr)\in 
\text{\rm J}_{_{\Omega}}.
\hfill
$
\\[1.75ex]
Кроме этого существуют такие функции $M_k^{}\in C^1\Omega,\ k=s+1,\ldots, m,$
что выполняются тождества {\rm (1.14)} и {\rm (1.15)}.
}
\vspace{0.5ex}

{\sl Следует} из свойства 1.11. $\k$
\vspace{1ex}

\newpage

{\bf  Свойство 3.9.}
\vspace{0.5ex}
{\it
Пусть 
$
\bigl(\exp\omega_\nu^{}, M_\nu^{}\bigr)\in \text{\rm E}_{_{\Omega}},\ 
\nu=1,\ldots, s,\ s\leqslant m-1,\  
({\rm g}_\xi^{}, M_\xi^{})\in \text{\rm J}_{_{\Omega}},
\linebreak 
{\rm g}_\xi^{{}^{\scriptsize \gamma_\xi^{}}}\!\in C^1\Omega,
\ 
\xi=s+1,\ldots,m,\ 
\gamma_j^{}\in\R\backslash\{0\}, 
\ j=1,\ldots,m. 
$ Тогда
\\[1.75ex]
\mbox{}\hfill
$
\displaystyle
\biggl(\,
\prod\limits_{\xi=s+1}^{m}
{\rm g}_\xi^{{}^{\scriptsize \gamma_\xi^{}}}
\exp\sum\limits_{\nu=1}^{s}
\gamma_\nu^{}\;\!\omega_\nu^{}\;\!,\;\! M\biggr)\in 
\text{\rm J}_{_{\Omega}},
\hfill
$
\\[1.75ex]
если и только если выполняется тождество {\rm (1.10)}.
}
\vspace{0.5ex}

{\sl Следует} из свойства 1.9. $\k$
\vspace{1ex}

{\bf Свойство 3.10.} 
\vspace{0.5ex}
{\it
Пусть функции 
$p,q\in \text{\rm P}_{_{\!\Xi^{\;\!\prime}}}$
взаимно простые, множество
$\Omega_{_0}\subset\Xi^{\;\!\prime}$ такое, что  
$
p(t,x)\ne 0\;\;\forall (t,x)\in\Omega_{_0},\
p(t,x)= 0\;\;\forall (t,x)\in {\sf C}_{_{\!\Xi^{\;\!\prime}}}\Omega_{_0}\;\!.
$
Тогда
$
\Bigl(\exp\dfrac{q}{p}\,, N\Bigr)\in 
\text{\rm E}_{_{\Omega_{{}_{\tiny\;\! 0}}}}, 
$
если и только если 
$(p,M)\in \text{\rm A}_{_{\!\Xi^{\;\!\prime}}}$
и выполняется тождество
\\[1.75ex]
\mbox{}\hfill                            % (3.2)
$
{\frak d}\;\!q(t,x)=
q(t,x)\;\!M(t,x)+p(t,x)\;\!N(t,x)
\quad
\forall (t,x)\in \Xi^{\;\!\prime},
$
\hfill {\rm (3.2)}
\\[1.75ex]
в котором функция $N\in \text{\rm P}_{_{\!\Xi^{\;\!\prime}}}$ и 
имеет степень $\deg_{x}^{}N\leq d-1.$
}
\vspace{0.75ex}

{\sl Доказательство.}
В соответствии с теоремой 3.1
$
\Bigl(\exp\dfrac{q}{p}\,, N\Bigr)\in 
\text{\rm E}_{_{\Omega_{{}_{\tiny\;\! 0}}}} 
$
тогда и только тогда, когда  производная в силу системы (0.1)
\\[1.25ex]
\mbox{}\hfill  
$
{\frak d}\,\dfrac{q(t,x)}{p(t,x)}=N(t,x)
\quad
\forall (t,x)\in\Omega_{_{0}}
\hfill 
$
\\[2ex]
при условии, что функция $N\in \text{\rm P}_{_{\!\Xi^{\;\!\prime}}}$  
и имеет степень $\deg_{x}^{}N\leq d-1.$
\vspace{0.75ex}

Отсюда следует, что 
$
\Bigl(\exp\dfrac{q}{p}\,, N\Bigr)\in 
\text{\rm E}_{_{\Omega}} 
$
тогда и только тогда, когда выполняется тождество
\\[1.5ex]
\mbox{}\hfill   % (3.3)
$
p(t,x)\;\!{\frak d}\;\!q(t,x)-
q(t,x)\;\!{\frak d}\;\!p(t,x)=
p^2(t,x)\;\!N(t,x)
\quad
\forall (t,x)\in \Xi^{\;\!\prime},
$
\hfill (3.3)
\\[1.75ex]
в котором функция $N\in \text{\rm P}_{_{\!\Xi^{\;\!\prime}}}$ и 
имеет степень $\deg_{x}^{}N\leq d-1.$
\vspace{0.75ex}

{\sl Необходимость.}
Из тождества (3.3) следует, что
\\[1.5ex]
\mbox{}\hfill                       % (3.4)
$
{\frak d}\;\!q(t,x)=
p(t,x)\;\!N(t,x)+
\dfrac{{\frak d}\;\!p(t,x)}{p(t,x)}\ q(t,x)
\quad
\forall (t,x)\in\Omega_{_{0}},
$
\hfill  (3.4)
\\[1.5ex]
причем фукнция $N\in \text{\rm P}_{_{\!\Xi^{\;\!\prime}}}$
и имеет степень $\deg_{x}^{}N\leq d-1.$
\vspace{0.75ex}

Так как 
\vspace{0.5ex}
$p,q, {\frak d}\;\!p, N\!\in\! \text{\rm P}_{_{\!\Xi^{\;\!\prime}}},\!$
а функции $p$ и $q$ взаимно простые, то из тождества (3.4) следует, что 
\vspace{0.5ex}
$\dfrac{{\frak d}\;\!p}{p}\in \text{\rm P}_{_{\!\Omega_{{}_{\tiny\;\! 0}}}}.$
Стало быть, существует такая функция 
$M\in \text{\rm P}_{_{\!\Xi^{\;\!\prime}}},$
что выполняется тождество (2.2), а значит, согласно теореме 2.1 
$(p,M)\in \text{\rm A}_{_{\!\Xi^{\;\!\prime}}}.$
\vspace{0.5ex}
Из тождества (3.4) при условии (2.2) получаем тождество (3.2), 
в котором $N\in \text{\rm P}_{_{\!\Xi^{\;\!\prime}}},$  $\deg_{x}^{}N\leq d-1.$
\vspace{0.5ex}

{\sl Достаточность.}
Из тождества (3.2), учитывая тождество (2.2), получаем тождество (3.4),
которое умножением обеих его частей на $p$ приводим к тождеству (3.3). 
\vspace{0.5ex}
Тем самым доказываем, что 
$
\Bigl(\exp\dfrac{q}{p}\,, N\Bigr)\in 
\text{\rm E}_{_{\Omega_{{}_{\tiny\;\! 0}}}}. \ \k
$
\vspace{1.25ex}

{\bf Свойство 3.11.}\! 
\vspace{0.5ex}
{\it
Пусть функции 
$u,v\!\in\! \text{\rm P}_{_{\!\Xi^{\;\!\prime}}}\!\!$
взаимно простые, множество
$\Omega_{_0}\!\subset\!\Xi^{\;\!\prime}\!$ та\-кое, что  
$
u(t,x)\!\ne\! 0\;\;\forall (t,x)\!\in\!\Omega_{_0},\
u(t,x)= 0\;\;\forall (t,x)\!\in\! {\sf C}_{_{\!\Xi^{\;\!\prime}}}\Omega_{_0}\;\!.\!$
\vspace{0.5ex}
Тогда
$\!
\Bigl(\exp\arctg\dfrac{v}{u}\,, V\!\Bigr)\!\in \!
\text{\rm E}_{_{\Omega_{{}_{\tiny\;\! 0}}}}\!\!$
в том и только в том случае, когда существует такая функция 
$U\in \text{\rm P}_{_{\!\Xi^{\;\!\prime}}},$
что вы\-пол\-ня\-ет\-ся система тождеств
\\[1.25ex]
\mbox{}\hfill                            % (3.5)
$
{\frak d}\;\!u(t,x)=
u(t,x)\;\!U(t,x)-v(t,x)\;\!V(t,x)
\quad
\forall (t,x)\in \Xi^{\;\!\prime},
\hfill
$
\\
\mbox{}\hfill {\rm (3.5)}
\\
\mbox{}\hfill                           
$
{\frak d}\;\!v(t,x)=
u(t,x)\;\!V(t,x)+v(t,x)\;\!U(t,x)
\quad
\forall (t,x)\in \Xi^{\;\!\prime},
\hfill
$
\\[1.5ex]
в которой	функция $V\in \text{\rm P}_{_{\!\Xi^{\;\!\prime}}}$ и 
имеет степень $\deg_{x}^{}V\leq d-1.$
}
\vspace{0.75ex}

{\sl Доказательство.}
По теореме 3.1, 
$
\Bigl(\exp\arctg\dfrac{v}{u}\,, V\Bigr)\in 
\text{\rm E}_{_{\Omega_{{}_{\tiny\;\! 0}}}}
$
тогда и только тогда, когда
\\[1.75ex]
\mbox{}\hfill                               % (3.6)
$
{\frak d}\,\arctg\dfrac{v(t,x)}{u(t,x)}=V(t,x)
\quad
\forall (t,x)\in\Omega_{_{0}}
$
\hfill  (3.6)
\\[2ex]
при условии, что функция $V\in \text{\rm P}_{_{\!\Xi^{\;\!\prime}}}$ 
имеет степень $\deg_{x}^{}V\leq d-1.$
\vspace{0.75ex}

Учитывая, что  производная в силу системы (0.1)
\\[1.5ex]
\mbox{}\hfill                           
$
{\frak d}\,\arctg\dfrac{v(t,x)}{u(t,x)}=
\dfrac{u(t,x)\,{\frak d}\;\!v(t,x)-v(t,x)\,{\frak d}\;\!u(t,x)}{u^2(t,x)+v^2(t,x)}
\quad
\forall (t,x)\in\Omega_{_{0}},
\hfill 
$
\\[1.5ex]
тождество (3.6) приводим к виду
\\[1.5ex]
\mbox{}\hfill                           
$
{\frak d}\;\!v(t,x)=
u(t,x)\;\!V(t,x)+v(t,x)\,
\dfrac{{\frak d}\;\!u(t,x)+v(t,x)\;\!V(t,x)}{u(t,x)}
\quad
\forall (t,x)\in\Omega_{_{0}}.
\hfill 
$
\\[1.5ex]
\indent
Функции
\vspace{0.35ex}
$u,v, {\frak d}\;\!u,{\frak d}\;\!v, V\in \text{\rm P}_{_{\!\Xi^{\;\!\prime}}},$
функции $u$ и $v$ взаимно простые, поэтому существует функция 
$U\in \text{\rm P}_{_{\!\Xi^{\;\!\prime}}}$ такая, что
\\[2ex]
\mbox{}\hfill 
$
{\frak d}\;\!u(t,x)+
v(t,x)\;\!V(t,x)=
u(t,x)\;\!U(t,x)
\quad
\forall (t,x)\in \Xi^{\;\!\prime}.
\hfill
$
\\[1.75ex]
\indent
В итоге получаем, что тождество (3.6) 
\vspace{0.35ex}
выполняется тогда и только тогда, когда имеет место система тождеств (3.5),
в которой $V\in \text{\rm P}_{_{\!\Xi^{\;\!\prime}}},\ \deg_{x}^{}V\leq d-1.\ \k$
\vspace{1ex}

{\bf Свойство 3.12.} 
\vspace{0.5ex}
{\it
Пусть  
$p\in \text{\rm P}_{_{\!\Xi^{\;\!\prime}}}.$
Тогда
$\bigl(\exp\arctg p, M\bigr)\in
\text{\rm E}_{_{\Xi^{\;\!\prime}}},
$
если и только если  производная в силу системы {\rm (0.1)}
\\[1.25ex]
\mbox{}\hfill                            % (3.7)
$
{\frak d}\;\!p(t,x)=
\bigl(1+p^2(t,x)\bigr)\;\!M(t,x)
\quad
\forall (t,x)\in \Xi^{\;\!\prime}.
$
\hfill {\rm (3.7)}
\\[1.5ex]
\indent
Сомножитель $M$ экспоненциального частного интеграла 
\vspace{0.35ex}
$\exp\arctg p$ имеет степень $\deg_{x}^{}M\leq d-1-\deg_{x}^{}p.$
}
\vspace{0.75ex}

{\sl Следует}
\vspace{0.25ex}
из свойства 3.11 при $v=p,\ u=1.$
То, что $\deg_{x}^{}M\leq d-1-\deg_{x}^{}p,$
получаем из тождества (3.7). $\k$ 
\vspace{1ex}

{\bf Свойство 3.13.} 
{\it
Пусть $f\in C^1\Omega.$ Тогда}
\\[1.5ex]
\mbox{}\hfill                       
$
\bigl(\exp\arctg f, M\bigr)\in
\text{\rm E}_{_{\Omega}}
\iff
\bigl(\exp\arcctg f, {}-M\bigr)\in
\text{\rm E}_{_{\Omega}}.
\hfill
$
\\[1.5ex]
\indent
{\sl Доказательство}
основано на теореме 3.1 и том, что
\\[1.5ex]
\mbox{}\hfill                       
$
\arcctg f(t,x)=\dfrac{\pi}{2}-\arctg f(t,x)
\quad
\forall (t,x)\in\Omega. \ \k
\hfill
$
\\[2ex]
\indent
{\bf Свойство 3.14.} 
\vspace{0.5ex}
{\it
Пусть  
$k\in\N,\ F\in \text{\rm I}_{_{\Omega}},$
функция $f\in C^1\Omega$ взаимно простая с функцией $F,$ множество
$\Omega_{_0}\subset\Omega$ та\-кое, что  
\vspace{0.5ex}
$
F(t,x)\ne 0\;\;\forall (t,x)\in\Omega_{_0},\
F(t,x)= 0\;\;\forall (t,x)\in 
\linebreak
\in {\sf C}_{_{\Omega}}\Omega_{_0}\;\!.$
Тогда
$
\Bigl(\exp\dfrac{f}{F^k}\,, M\Bigr)\in 
\text{\rm E}_{_{\Omega_{{}_{\tiny\;\! 0}}}},
$
если и только если выполняется тождество
\\[2ex]
\mbox{}\hfill                      
$
{\frak d}\;\!f(t,x)=
F^k(t,x)\;\!M(t,x)
\quad
\forall (t,x)\in \Omega,
\hfill
$
\\[2ex]
в котором	функция $M\in \text{\rm P}_{_{\!\Xi^{\;\!\prime}}}$  и
имеет степень $\deg_{x}^{}M\leq d-1.$
}
\vspace{0.75ex}

{\sl Доказательство.}
\vspace{0.5ex}
Учитывая тождество (0.3) из теоремы 0.2,
на основании теоремы 3.1 получаем: 
$
\Bigl(\exp\dfrac{f}{F^k}\,, M\Bigr)\in 
\text{\rm E}_{_{\Omega_{{}_{\tiny\;\! 0}}}}
$
тогда и только тогда, когда выполняется тождество
\\[2ex]
\mbox{}\hfill                               
$
{\frak d}\;\!\dfrac{f(t,x)}{F^k(t,x)}=
\dfrac{F^k(t,x)\,{\frak d}f(t,x)-kf(t,x)\;\!F^{k-1}(t,x)\,{\frak d}F(t,x)}{F^{\;\!2k}(t,x)}=
\dfrac{{\frak d}f(t,x)}{F^k(t,x)}\!=\!M(t,x)
\ 
\forall (t,x)\!\in\!\Omega_{_{0}},
\hfill
$
\\[2.25ex]
в котором функция $M\in \text{\rm P}_{_{\!\Xi^{\;\!\prime}}}$ 
имеет степень
$\deg_{x}^{}M\leq d-1.\ \k$
\vspace{1ex}

Приложение свойства 3.14 дано в [56, с. 52 -- 55 ].

\newpage

\mbox{}
\\[-2.75ex]
\centerline{
{\bf  4. Условные частные интегралы}
}
\\[1.5ex]
\indent
{\bf Определение 4.1.}
\vspace{0.5ex}
{\it
Если $p\in \text{\rm P}_{_{\!\Xi^{\;\!\prime}}},\  
\bigl(\exp p, M\bigr)\in \text{\rm E}_{_{\Xi^{\;\!\prime}}},$
то функцию $\exp p$ назовем
\textit{\textbf{ус\-лов\-ным частным интегралом}}
\vspace{0.5ex}
с сомножителем $M$ на области $\Xi^{\;\!\prime}$ сис\-темы {\rm (0.1)}.}

Множество функций, являющихся условными частными интегралами 
на области $\Xi^{\;\!\prime}$ сис\-темы (0.1), обозначим 
$\text{\rm F}_{_{\!\Xi^{\;\!\prime}}}.$
\vspace{0.5ex}

В соответствии 
\vspace{0.15ex}
с определениями условного и экспоненциального частных интегралов
(определения 4.1 и 3.1) множество 
$\text{\rm F}_{_{\!\Xi^{\;\!\prime}}}\subset
\text{\rm E}_{_{\Xi^{\;\!\prime}}}\subset 
\text{\rm J}_{_{\!\Xi^{\;\!\prime}}}.$
\vspace{1ex}

Условной записью 
\vspace{0.5ex}
$\bigl(\exp p, M\bigr)\in \text{\rm F}_{_{\!\Xi^{\;\!\prime}}}$
будем выражать, что функция $\exp p$ является 
условным частным интегралом  
\vspace{0.75ex}
с сомножителем $M$ на области $\Xi^{\;\!\prime}$ сис\-темы (0.1).

Определение 4.1 выражается эквиваленцией
\\[1.5ex]
\mbox{}\hfill               % (4.1)
$
\bigl(\exp p, M\bigr)\in \text{\rm F}_{_{\!\Xi^{\;\!\prime}}}
\iff
\bigl(\exp p, M\bigr)\in \text{\rm E}_{_{\Xi^{\;\!\prime}}}
\ \&\ 
p\in \text{\rm P}_{_{\!\Xi^{\;\!\prime}}}.
$
\hfill {\rm(4.1)}
\\[2ex]
\indent
{\bf Теорема 4.1}\!
\vspace{0.35ex}
(критерий существования условного частного интеграла).\!
{\it
Пусть $\!p\!\in\! \text{\rm P}_{_{\!\!\Xi^{\;\!\prime}}}.$
Тогда  
$\bigl(\exp p, M\bigr)\in \text{\rm F}_{_{\!\Xi^{\;\!\prime}}},$
если и только если  производная в силу системы {\rm (0.1)}
\\[1.5ex]
\mbox{}\hfill               % (4.2)
$
{\frak d}\;\!p(t,x)=M(t,x)
\quad
\forall (t,x)\in \Xi^{\;\!\prime}
$
\hfill {\rm(4.2)}
\\[1.5ex]
при условии $\deg_{\;\!x}^{} M\leq d-1.$
}
\vspace{0.5ex}

{\sl Следует}
с учетом эквиваленции (4.1) из теоремы 3.1 при 
$\omega=p\in \text{\rm P}_{_{\!\Xi^{\;\!\prime}}}.\ \k$ 
\vspace{0.75ex}

В соответствии с теоремой 4.1 имеет место
\vspace{0.5ex}

{\bf Теорема 4.2}\!
\vspace{0.35ex}
(критерий существования условного частного интеграла).\!
{\it
Пусть $\!p\!\in\! \text{\rm P}_{_{\!\!\Xi^{\;\!\prime}}}.$
Тогда  
\\[0.5ex]
\mbox{}\hfill
$
\exp p\in \text{\rm F}_{_{\!\Xi^{\;\!\prime}}}
\iff
\deg_{\;\!x}^{}{\frak d}\;\!p\leq d-1.
\hfill
$
\\[1.5ex]
При этом функция ${\frak d}\;\!p$ является сомножителем 
\vspace{0.75ex}
условного частного интеграла $\exp p.$}

Теорема 4.2 является аналогом теоремы 3.2 на случай условного частного интеграла
(когда $\omega=p\in \text{\rm P}_{_{\!\Xi^{\;\!\prime}}}).$
\vspace{0.75ex}

{\bf Замечание 4.1.}
\vspace{0.5ex}
Если $p\in \text{\rm P}_{_{\!\Xi^{\;\!\prime}}},$ то 
${\frak d}\;\!p\in \text{\rm P}_{_{\!\Xi^{\;\!\prime}}},$ а
$\deg_{\;\!x}^{}{\frak d}\;\!p\leq d-1+\deg_{\;\!x}^{}p.$
Поэтому требование $\deg_{\;\!x}^{} M\leq d-1$ 
\vspace{0.35ex}
в тождестве (4.2) из теоремы 4.1 и требование 
$\deg_{\;\!x}^{}{\frak d}\;\!p\leq d-1$
в теореме 4.2 имеют существенное значение.
\vspace{0.75ex}

{\bf  Теорема 4.3.}
\vspace{0.5ex}
{\it
Пусть 
$p_j^{}\in \text{\rm P}_{_{\!\Xi^{\;\!\prime}}},\
\gamma_j^{}\in\R\backslash\{0\},\ c_j^{}\in\R, \ j=1,\ldots,m.$
Тогда
\\[1.5ex]
\mbox{}\hfill
$
\displaystyle
\exp\sum\limits_{j=1}^{m}
\gamma_j^{}\;\!(p_j^{}+c_j^{})\in \text{\rm F}_{_{\!\Xi^{\;\!\prime}}}
\ \iff \
\deg_{\;\!x}^{} 
\sum\limits_{j=1}^{m}
\gamma_j^{}\,{\frak d}\;\!p_j^{}\leq d-1.
\hfill
$
\\[1.5ex]
При этом функция 
$\sum\limits_{j=1}^{m}\gamma_j^{}\;\!{\frak d}\;\!p_j^{}$
является сомножителем условного частного интеграла 
$\exp\sum\limits_{j=1}^{m}\gamma_j^{}\;\!(p_j^{}+c_j^{}).$
}

{\sl Следует} из теоремы 4.2 при 
$p=\sum\limits_{j=1}^{m}\gamma_j^{}\;\!(p_j^{}+c_j^{}).\ \k$
\vspace{0.75ex}

Теорема 4.3 
\vspace{0.35ex}
является аналогом теоремы 1.3 на случай
условного частного интеграла, когда функции
${\rm g}_j^{}=\exp\bigl(\gamma_j^{}\;\!(p_j^{}+c_j^{})\bigr),
\ j=1,\ldots,m.$
\vspace{1ex}

{\bf  Свойство 4.1.}
\vspace{0.75ex}
{\it
Если $\varphi\in C^1T^{\;\!\prime},$ то} 
$
\bigl(\exp\varphi, {\sf D} \varphi\bigr)\in \text{\rm F}_{_{\!\Xi^{\;\!\prime}}}.
$

{\sl Действительно}, 
\vspace{0.75ex}
$\varphi\in \text{\rm P}_{_{\!\Xi^{\;\!\prime}}},\
{\frak d}\;\!\varphi(t)={\sf D}\varphi(t)\;\;\forall (t,x)\in \Xi^{\;\!\prime},\
\deg_{\;\!x}^{}{\sf D}\varphi=0\leq d-1.$
По тео\-ре\-ме~4.1,  
$
\bigl(\exp\varphi, {\sf D} \varphi\bigr)\in \text{\rm F}_{_{\!\Xi^{\;\!\prime}}}.\ \k
$
\vspace{0.75ex}

\newpage

{\bf  Свойство 4.2.}
{\it
Пусть $\varphi\in C^1T^{\;\!\prime}, \
\gamma_j^{}\in\R\backslash\{0\},\ c_j^{}\in\R, \ j=1,\ldots,m,\ m\leq n.$
Тогда
\\[1.5ex]
\mbox{}\hfill
$
\displaystyle
\exp\biggl(\varphi(t)+\sum\limits_{j=1}^{m}
\gamma_j^{}\;\!(x_j^{}+c_j^{})\biggr)\in \text{\rm F}_{_{\!\Xi^{\;\!\prime}}}
\ \iff \
\deg_{\;\!x}^{} 
\sum\limits_{j=1}^{m}
\gamma_j^{}\,X_j^{}\leq d-1.
\hfill
$
\\[1.5ex]
При этом функция 
${\sf D} \varphi+\sum\limits_{j=1}^{m}\gamma_j^{}\;\!X_j^{}$
является сомножителем условного частного интеграла 
$\exp\Bigl(\varphi(t)+\sum\limits_{j=1}^{m}
\gamma_j^{}\;\!(x_j^{}+c_j^{})\Bigr).$
}
\vspace{0.5ex}

{\sl Доказательство} 
основано на теореме 4.1 и том, что
\\[1.5ex]
\mbox{}\hfill
$
\displaystyle
{\frak d}\;\!\biggl(\varphi(t)+\sum\limits_{j=1}^{m}
\gamma_j^{}\;\!(x_j^{}+c_j^{})\biggr)=
{\sf D} \varphi(t)+\sum\limits_{j=1}^{m}\gamma_j^{}\;\!X_j^{}(t,x)
\quad
\forall (t,x)\in \Xi^{\;\!\prime}.\ \k
\hfill
$
\\[1.75ex]
\indent
В частности, при $m=1$ имеет место
\vspace{0.75ex}

{\bf  Свойство 4.3.}
{\it
Пусть $\varphi\in C^1T^{\;\!\prime}, \
\gamma\in\R\backslash\{0\},\ c\in\R.$
Тогда
\\[1.5ex]
\mbox{}\hfill
$
\displaystyle
\exp\bigl(\varphi(t)+\gamma\;\!(x_j^{}+c)\bigr)\in \text{\rm F}_{_{\!\Xi^{\;\!\prime}}}
\ \iff \
\deg_{\;\!x}^{} X_j^{}\leq d-1,
\quad 
j\in \{1,\ldots, n\}. 
\hfill
$
\\[1.5ex]
При этом функция 
\vspace{0.35ex}
${\sf D} \varphi+\gamma\;\!X_j^{}$
является сомножителем условного частного интеграла 
$\exp\bigl(\varphi(t)+\gamma\;\!(x_j^{}+c)\bigr).$
}
\mbox{}
\\[3.75ex]
\centerline{
{\bf  5. Кратные полиномиальные частные интегралы}
}
\\[1.5ex]
\indent
{\bf Определение 5.1.}
\vspace{0.25ex}
{\it
Полиномиальный частный интеграл $p$ с сомножителем $M$ 
на области $\Xi^{\;\!\prime}$ сис\-темы {\rm (0.1)} назовем
\textit{\textbf{кратным}},
\vspace{0.15ex}
если существуют такие натуральное число $h$ и функция
$q\in \text{\rm P}_{_{\!\Xi^{\;\!\prime}}},$
\vspace{0.35ex}
являющаяся взаимно простой с функцией $p,$
что производная в силу системы {\rm (0.1)} 
\\[1.5ex]
\mbox{}\hfill                             % (5.1)
$
\displaystyle
{\frak d}\,\dfrac{q(t,x)}{p^{\;\!h}(t,x)}=N(t,x)
\quad
\forall (t,x)\in \Omega_{_0},
$
\hfill {\rm (5.1)}
\\[1.75ex]
где функция $N\in \text{\rm P}_{_{\!\Xi^{\;\!\prime}}}$ 
и имеет степень $\deg_{\;\!x}^{} N\leq d-1,$
\vspace{0.75ex}
множество
$\Omega_{_0}\subset\Xi^{\;\!\prime}$ та\-кое, что  
$
p(t,x)\ne 0\;\;\forall (t,x)\in\Omega_{_0},\
p(t,x)= 0\;\;\forall (t,x)\in {\sf C}_{_{\Xi^{\;\!\prime}}}\Omega_{_0}\;\!.$
}
\vspace{0.75ex}

Множество функций, являющихся кратными полиномиальными частными интегралами 
на области $\Xi^{\;\!\prime}$ сис\-темы (0.1), обозначим 
$\text{\rm B}_{_{\Xi^{\;\!\prime}}}.$
\vspace{0.5ex}

В соответствии 
\vspace{0.15ex}
с определением 5.1
$\text{\rm B}_{_{\Xi^{\;\!\prime}}}\subset
\text{\rm A}_{_{\Xi^{\;\!\prime}}}\subset 
\text{\rm J}_{_{\!\Xi^{\;\!\prime}}}.$
\vspace{1ex}

Условной записью 
\vspace{0.35ex}
$\bigl((p, M), (h,q,N)\bigr)\in \text{\rm B}_{_{\Xi^{\;\!\prime}}}$
будем выражать, что полиномиальный частный интеграл $p$ с сомножителем $M$ 
\vspace{0.15ex}
на области $\Xi^{\;\!\prime}$ сис\-темы {\rm (0.1)} является кратным 
таким, что выполняется тождество (5.1).
\vspace{0.5ex}

{\bf Теорема 5.1}
\vspace{0.15ex}
(критерий существования кратного полиномиального частного интеграла).
{\it
$\bigl((p, M), (h,q,N)\bigr)\in \text{\rm B}_{_{\Xi^{\;\!\prime}}}$
\vspace{0.35ex}
тогда и только тогда, когда выполняются тождества {\rm(2.2)} и {\rm(5.1)}, 
в которых число $h\in\N,$ функции
\vspace{0.5ex}
$q, N\in \text{\rm P}_{_{\!\Xi^{\;\!\prime}}},$
функции $p$ и $q$ взаимно простые, $\deg_{\;\!x}^{} N\leq d-1,$
\vspace{0.75ex}
множество
$\Omega_{_0}\subset\Xi^{\;\!\prime}$ та\-кое, что  
$p(t,x)\ne 0\;\;\forall (t,x)\in\Omega_{_0},$
$p(t,x)= 0\;\;\forall (t,x)\in {\sf C}_{_{\Xi^{\;\!\prime}}}\Omega_{_0}\;\!.$
}
\vspace{0.5ex}

{\sl Следует} из определения 5.1 и теоремы 2.1 (с учетом замечания 2.1). $\k$
\vspace{0.75ex}

{\bf Теорема 5.2}
\vspace{0.15ex}
(критерий существования кратного полиномиального частного интеграла).
{\it
$\bigl((p, M), (h,q,N)\bigr)\in \text{\rm B}_{_{\Xi^{\;\!\prime}}},$
если и только если выполняются тождества {\rm(2.2)} и 
\\[1.5ex]
\mbox{}\hfill                             % (5.2)
$
\displaystyle
{\frak d}\;\!q(t,x)=h\;\!q(t,x)\;\!M(t,x)+p^{\;\!h}(t,x)\;\!N(t,x)
\quad
\forall (t,x)\in \Xi^{\;\!\prime},
$
\hfill {\rm (5.2)}
\\[1.75ex]
в которых число $h\in\N,$ функции
\vspace{0.5ex}
$q, N\in \text{\rm P}_{_{\!\Xi^{\;\!\prime}}},$
функции $p$ и $q$ взаимно простые, $\deg_{\;\!x}^{} N\leq d-1.$
}

\newpage

{\sl Доказательство}
основано на теореме 5.1 и состоит в том, что при условии (2.2) 
тождество (5.1) выполняется, если и только если 
выполняется тождество (5.2), так как
\\[1.5ex]
\mbox{}\hfill                           
$
\displaystyle
{\frak d}\,\dfrac{q(t,x)}{p^{\;\!h}(t,x)}=
\dfrac{p^{\;\!h}(t,x)\, {\frak d}\;\!q(t,x)-
h\;\!p^{\;\!h-1}(t,x)\;\!q(t,x)\, {\frak d}\;\!p(t,x)}{p^{\;\!2h}(t,x)}=
\hfill
$
\\[2.25ex]
\mbox{}\hfill                           
$
\displaystyle
=\dfrac{{\frak d}\;\!q(t,x)-
h\;\!q(t,x)\;\!M(t,x)}{p^{\;\!h}(t,x)}
\quad
\forall (t,x)\in \Omega_{_0},
\hfill
$
\\[1.75ex]
где множество
\vspace{1ex}
$\Omega_{_0}\subset\Xi^{\;\!\prime}$ та\-кое, что  
$p(t,x)\ne 0\;\;\forall (t,x)\in\Omega_{_0},$
$p(t,x)= 0\;\;\forall (t,x)\in {\sf C}_{_{\Xi^{\;\!\prime}}}\Omega_{_0}\;\!.\ \k$

{\bf Лемма 5.1.}
\vspace{0.15ex}
{\it
Пусть $h\in\N,$ функции
$p, q\in \text{\rm P}_{_{\!\Xi^{\;\!\prime}}}\!$ взаимно простые, 
множество
$\Omega_{_0}\!\subset\!\Xi^{\;\!\prime}\!$ та\-кое, что  
$p(t,x)\!\ne\! 0\;\;\forall (t,x)\in\Omega_{_0},\
p(t,x)= 0\;\;\forall (t,x)\in {\sf C}_{_{\Xi^{\;\!\prime}}}\Omega_{_0}.\!$
Тогда
$
\Bigl(\exp\dfrac{q}{p^{\;\!h}}\,, N\Bigr)\!\in\! 
\text{\rm E}_{_{\Omega_{{}_{\tiny\;\! 0}}}},
$
если и только если 
\vspace{0.5ex}
$(p, M)\in \text{\rm A}_{_{\Xi^{\;\!\prime}}}$
и выполняется тождество {\rm (5.2)}, 
в котором	функция $N\in \text{\rm P}_{_{\!\Xi^{\;\!\prime}}}$  и
имеет степень $\deg_{x}^{}N\leq d-1.$
}
\vspace{0.5ex}

{\sl Следует}
из свойства 3.10 с учетом того, что, по теореме 2.2,
\\[1.5ex]
\mbox{}\hfill                           
$
\displaystyle
(p, M)\in \text{\rm A}_{_{\Xi^{\;\!\prime}}}
\iff 
\bigl(p^{\;\!h}, h\;\!M\bigr)\in \text{\rm A}_{_{\Xi^{\;\!\prime}}}. \ \k
\hfill
$
\\[2ex]
\indent
{\bf Теорема 5.3}
\vspace{0.15ex}
(критерий существования кратного полиномиального частного интеграла).
{\it
Пусть $h\in\N,$ функции
\vspace{0.5ex}
$p, q\in \text{\rm P}_{_{\!\Xi^{\;\!\prime}}}\!$ взаимно простые, 
множество
$\Omega_{_0}\subset\Xi^{\;\!\prime}$ та\-кое, что  
$p(t,x)\ne 0\;\;\forall (t,x)\in\Omega_{_0},\
p(t,x)= 0\;\;\forall (t,x)\in {\sf C}_{_{\Xi^{\;\!\prime}}}\Omega_{_0}.$
Тогда}
\\[1.5ex]
\mbox{}\hfill                           
$
\displaystyle
\bigl((p, M), (h,q,N)\bigr)\in \text{\rm B}_{_{\Xi^{\;\!\prime}}}
\iff
\Bigl(\exp\dfrac{q}{p^{\;\!h}}\,, N\Bigr)\in 
\text{\rm E}_{_{\Omega_{{}_{\tiny\;\! 0}}}}.
\hfill
$
\\[2ex]
\indent
{\sl Следует}
с учетом теоремы 2.1 из теоремы 5.2 и леммы 5.1. $\k$
\vspace{0.5ex}

Согласно теореме 5.3 (с учетом леммы 5.1) у 
\vspace{0.35ex}
$\bigl((p, M), (h,q,N)\bigr)\in \text{\rm B}_{_{\Xi^{\;\!\prime}}}$
два сомножителя: 
$M$ --- сомножитель полиномиального частного интеграла $p$ и 
\vspace{0.35ex}
$N$ --- сомножитель экспоненциального частного интеграла 
$\exp\dfrac{q}{p^{\;\!h}}\,.$
\vspace{0.75ex}

{\bf Свойство 5.1.}
{\it
Пусть $h\in\N,\ \varphi\in C^1T^{\;\!\prime},$ функции
\vspace{0.5ex}
$p, q\in \text{\rm P}_{_{\!\Xi^{\;\!\prime}}}\!$ взаимно простые, 
мно\-жес\-т\-во
$\Omega_{_0}\subset\Xi^{\;\!\prime}$ та\-кое, что  
$p(t,x)\ne 0\;\;\forall (t,x)\in\Omega_{_0},\
p(t,x)= 0\;\;\forall (t,x)\in {\sf C}_{_{\Xi^{\;\!\prime}}}\Omega_{_0}.$
Тогда}
\\[1.5ex]
\mbox{}\hfill                           
$
\displaystyle
\bigl((p, M), (h,q,N)\bigr)\in \text{\rm B}_{_{\Xi^{\;\!\prime}}}
\iff
\biggl(\exp\Bigl(\;\!\dfrac{q}{p^{\;\!h}}+\varphi\Bigr), N+{\sf D}\;\!\varphi\biggr)\in 
\text{\rm E}_{_{\Omega_{{}_{\tiny\;\! 0}}}}.
\hfill
$
\\[2ex]
\indent
{\sl Следует}
с учетом свойства 3.1 из теоремы 5.3. $\k$
\vspace{0.5ex}

{\bf Следствие 5.1.}
{\it
Пусть $h\in\N,\ c\in\R,$ функции
\vspace{0.5ex}
$p, q\in \text{\rm P}_{_{\!\Xi^{\;\!\prime}}}\!$ взаимно простые, 
мно\-жес\-т\-во
$\Omega_{_0}\subset\Xi^{\;\!\prime}$ та\-кое, что  
$p(t,x)\ne 0\;\;\forall (t,x)\in\Omega_{_0},\
p(t,x)= 0\;\;\forall (t,x)\in {\sf C}_{_{\Xi^{\;\!\prime}}}\Omega_{_0}.$
Тогда}
\\[1.5ex]
\mbox{}\hfill                           
$
\displaystyle
\bigl((p, M), (h,q,N)\bigr)\in \text{\rm B}_{_{\Xi^{\;\!\prime}}}
\iff
\biggl(\exp\Bigl(\;\!\dfrac{q}{p^{\;\!h}}+c\Bigr), N\biggr)\in 
\text{\rm E}_{_{\Omega_{{}_{\tiny\;\! 0}}}}.
\hfill
$
\\[2ex]
\indent
{\bf Свойство 5.2.}
\vspace{0.75ex}
{\it
Пусть 
$\bigl((p, M), (h,q,N)\bigr)\in \text{\rm B}_{_{\Xi^{\;\!\prime}}},\ 
\gamma_1^{},\gamma_2^{}\in\R,\ \varphi\in C^1T^{\;\!\prime},$ 
мно\-жес\-т\-во
$\Omega_{_0}\subset\Xi^{\;\!\prime}$ та\-кое, что  
$p(t,x)\ne 0\;\;\forall (t,x)\in\Omega_{_0},\
p(t,x)= 0\;\;\forall (t,x)\in {\sf C}_{_{\Xi^{\;\!\prime}}}\Omega_{_0}.$
Тогда
\\[1.5ex]
\mbox{}\hfill                           
$
\displaystyle
\biggl(p^{{}^{\scriptsize \gamma_1^{}}}
\exp\Bigl(\gamma_2^{}\Bigl(\;\!\dfrac{q}{p^{\;\!h}}+\varphi\Bigr)\Bigr), L\biggr)\in 
\text{\rm J}_{_{\Omega_{{}_{\tiny\;\! 0}}}},
\hfill
$
\\[2ex]
если и только если сомножители $M,\, N,\, L$ такие, что}
\\[1.5ex]
\mbox{}\hfill                           
$
\displaystyle
L(t,x)=\gamma_1^{}\;\!M(t,x)+\gamma_2^{}\;\!
\bigl(N(t,x)+{\sf D}\;\!\varphi(t)\bigr)
\quad
\forall (t,x)\in \Xi^{\;\!\prime}.
\hfill
$
\\[1.75ex]
\indent
{\sl Следует}
с учетом свойства 5.1 из свойства 3.9. $\k$
\vspace{0.75ex}

{\bf Следствие 5.2.}
\vspace{0.75ex}
{\it
Пусть 
$\bigl((p, M), (h,q,N)\bigr)\in \text{\rm B}_{_{\Xi^{\;\!\prime}}},\ 
c, \gamma_1^{},\gamma_2^{}\in\R,$ 
мно\-жес\-т\-во
$\Omega_{_0}\subset\Xi^{\;\!\prime}$ та\-кое, что  
$p(t,x)\ne 0\;\;\forall (t,x)\in\Omega_{_0},\
p(t,x)= 0\;\;\forall (t,x)\in {\sf C}_{_{\Xi^{\;\!\prime}}}\Omega_{_0}.$
Тогда
\\[1.5ex]
\mbox{}\hfill                           
$
\displaystyle
\biggl(p^{\gamma_1^{}}
\exp\Bigl(\gamma_2^{}\Bigl(\;\!\dfrac{q}{p^{\;\!h}}+c\Bigr)\Bigr), L\biggr)\in 
\text{\rm J}_{_{\Omega_{{}_{\tiny\;\! 0}}}},
\hfill
$
\\[2ex]
если и только если сомножители $M,\, N,\, L$ такие, что}
\\[1.5ex]
\mbox{}\hfill                           
$
\displaystyle
L(t,x)=\gamma_1^{}\;\!M(t,x)+\gamma_2^{}\;\!N(t,x)
\quad
\forall (t,x)\in \Xi^{\;\!\prime}.
\hfill
$
\\[2.25ex]
\indent
{\bf Свойство 5.3.}
\vspace{0.75ex}
{\it
Пусть 
$\bigl((p_j^{}, M_j^{}), (h_j^{},q_j^{},N)\bigr)\in \text{\rm B}_{_{\Xi^{\;\!\prime}}},\ 
\lambda_j^{},\gamma_j^{}\in\R,\ j=1,\ldots, m,$ 
$\sum\limits_{j=1}^{m}|\lambda_j^{}|\ne 0,\ 
\varphi\in C^1T^{\;\!\prime},$ 
мно\-жес\-т\-во
\vspace{0.75ex}
$\Omega_{_0}\subset\Xi^{\;\!\prime}$ та\-кое, что  
$\prod\limits_{j=1}^{m}p_j^{}(t,x)\ne 0\;\;\forall (t,x)\in\Omega_{_0},$
$\prod\limits_{j=1}^{m}p_j^{}(t,x)= 0\;\;
\forall (t,x)\in {\sf C}_{_{\Xi^{\;\!\prime}}}\Omega_{_0},\ 
p_j^{\gamma_j^{}}\in C^1\Omega_{_0}, \ j=1,\ldots, m.$
Тогда
\\[1.5ex]
\mbox{}\hfill                           
$
\displaystyle
\biggl(\ \prod\limits_{j=1}^{m}p_j^{\gamma_j^{}}
\sum\limits_{j=1}^{m}\lambda_j^{}
\exp\biggl(\;\!\dfrac{q_j^{}}{p_j^{\;\!h_j^{}}}+\varphi\biggr), L\biggr)\in 
\text{\rm J}_{_{\Omega_{{}_{\tiny\;\! 0}}}},
\hfill
$
\\[2ex]
если и только если сомножители $M_1^{},\ldots, M_{m}^{},\, N,\, L$ связаны тождеством}
\\[1.75ex]
\mbox{}\hfill                           
$
\displaystyle
L(t,x)={\sf D}\;\!\varphi(t)+N(t,x)+
\sum\limits_{j=1}^{m}\gamma_j^{}\;\!M_j^{}(t,x)
\quad
\forall (t,x)\in \Xi^{\;\!\prime}.
\hfill
$
\\[1.5ex]
\indent
{\sl Доказательство}.
Последовательно используя свойства 5.1 и 3.4, получаем:
\\[1.5ex]
\mbox{}\hfill                           
$
\displaystyle
\biggl(
\exp\biggl(\;\!\dfrac{q_j^{}}{p_j^{\;\!h_j^{}}}+\varphi\biggr),\, 
N+{\sf D}\;\!\varphi\biggr)\in 
\text{\rm E}_{_{\Omega_{{}_{\tiny\;\! 0}}}},\ j=1,\ldots, m,
\ \Longrightarrow 
\hfill                           
$
\\[1.5ex]
\mbox{}\hfill                           
$
\displaystyle
\Longrightarrow \
\biggl(\ \sum\limits_{j=1}^{m}\lambda_j^{}
\exp\biggl(\;\!\dfrac{q_j^{}}{p_j^{\;\!h_j^{}}}+\varphi\biggr),\, 
N+{\sf D}\;\!\varphi\biggr)\in 
\text{\rm J}_{_{\Omega_{{}_{\tiny\;\! 0}}}}.
\hfill
$
\\[1.5ex]
\indent
Теперь, учитывая, что
\vspace{0.5ex} 
$\bigl(p_j^{}, M_j^{}\bigr)\in \text{\rm A}_{_{\Xi^{\;\!\prime}}},\ j=1,\ldots, m,$ 
на основании свойства 1.9 получаем утверждение доказываемого свойства. $\k$
\vspace{0.75ex}

{\bf Свойство 5.4.}
\vspace{0.75ex}
{\it
Пусть 
$\bigl((p_j^{}, M_j^{}), (h_j^{},q_j^{},N)\bigr)\in \text{\rm B}_{_{\Xi^{\;\!\prime}}},\ 
c_j^{}, \lambda_j^{}, \gamma_j^{}\in\R,\ j=1,\ldots, m,$ 
$\sum\limits_{j=1}^{m}|\lambda_j^{}|\ne 0,$ 
мно\-жес\-т\-во
\vspace{0.75ex}
$\Omega_{_0}\subset\Xi^{\;\!\prime}$ та\-кое, что  
$\prod\limits_{j=1}^{m}p_j^{}(t,x)\ne 0\;\;\forall (t,x)\in\Omega_{_0},\
\prod\limits_{j=1}^{m}p_j^{}(t,x)= 0$
$\forall (t,x)\in {\sf C}_{_{\Xi^{\;\!\prime}}}\Omega_{_0},\ 
p_j^{\gamma_j^{}}\in C^1\Omega_{_0}, \ j=1,\ldots, m.$
Тогда
\\[1.25ex]
\mbox{}\hfill                           
$
\displaystyle
\biggl(\ \prod\limits_{j=1}^{m}p_j^{\gamma_j^{}}
\sum\limits_{j=1}^{m}\lambda_j^{}
\exp\biggl(\;\!\dfrac{q_j^{}}{p_j^{\;\! h_j^{}}}+c_j^{}\biggr), L\biggr)\in 
\text{\rm J}_{_{\Omega_{{}_{\tiny\;\! 0}}}},
\hfill
$
\\[2ex]
если и только если сомножители $M_1^{},\ldots, M_{m}^{},\, N,\, L$ связаны тождеством}
\\[1.75ex]
\mbox{}\hfill                           
$
\displaystyle
L(t,x)=N(t,x)+
\sum\limits_{j=1}^{m}\gamma_j^{}\;\!M_j^{}(t,x)
\quad
\forall (t,x)\in \Xi^{\;\!\prime}.
\hfill
$
\\[1.5ex]
\indent
{\sl Доказательство}.
Последовательно используя следствие 5.1 и свойство 3.4, получаем:
\\[1.5ex]
\mbox{}\hfill                           
$
\displaystyle
\biggl(
\exp\biggl(\;\!\dfrac{q_j^{}}{p_j^{\;\!h_j^{}}}+c_j^{}\biggr),\, N\biggr)\in 
\text{\rm E}_{_{\Omega_{{}_{\tiny\;\! 0}}}},\ j=1,\ldots, m,
\ \Longrightarrow \
\biggl(\ \sum\limits_{j=1}^{m}\lambda_j^{}
\exp\biggl(\;\!\dfrac{q_j^{}}{p_j^{\;\!h_j^{}}}+c_j^{}\biggr),\, N\biggr)\in 
\text{\rm J}_{_{\Omega_{{}_{\tiny\;\! 0}}}}.
\hfill
$
\\[1.5ex]
\indent
Теперь, учитывая, что
\vspace{0.5ex} 
$\bigl(p_j^{}, M_j^{}\bigr)\in \text{\rm A}_{_{\Xi^{\;\!\prime}}},\ j=1,\ldots, m,$ 
на основании свойства 1.9 получаем утверждение доказываемого свойства. $\k$
\vspace{0.75ex} 

{\bf Свойство 5.5.}\!
\vspace{0.75ex}
{\it
Пусть 
$\!\bigl((p_j^{}, M_j^{}), (h_j^{},q_j^{},\rho_{j}^{}N_0^{})\bigr)\!\in\! 
\text{\rm B}_{_{\Xi^{\;\!\prime}}},\, 
\rho_j^{}\!\in\!\R\backslash\{0\},\, 
\lambda_j^{},\gamma_j^{}\!\in\!\R,\, j\!=\!1,\ldots, m,\!$ 
$\sum\limits_{j=1}^{m}|\lambda_j^{}|\ne 0,\ 
\varphi\in C^1T^{\;\!\prime},$ 
мно\-жес\-т\-во
\vspace{0.75ex}
$\Omega_{_0}\subset\Xi^{\;\!\prime}$ та\-кое, что  
$\prod\limits_{j=1}^{m}p_j^{}(t,x)\ne 0\;\;\forall (t,x)\in\Omega_{_0},$
$\prod\limits_{j=1}^{m}p_j^{}(t,x)= 0\;\;
\forall (t,x)\in {\sf C}_{_{\Xi^{\;\!\prime}}}\Omega_{_0},\ 
p_j^{\gamma_j^{}}\in C^1\Omega_{_0}, \ j=1,\ldots, m.$
Тогда
\\[1.5ex]
\mbox{}\hfill                           
$
\displaystyle
\biggl(\ \prod\limits_{j=1}^{m}p_j^{\gamma_j^{}}
\sum\limits_{j=1}^{m}\lambda_j^{}
\exp\biggl(\;\!\dfrac{q_j^{}}{\rho_j^{}\;\!p_j^{\;\!h_j^{}}}+\varphi\biggr), L\biggr)\in 
\text{\rm J}_{_{\Omega_{{}_{\tiny\;\! 0}}}},
\hfill
$
\\[2ex]
если и только если сомножитель}
\\[1.75ex]
\mbox{}\hfill                           
$
\displaystyle
L(t,x)={\sf D}\;\!\varphi(t)+N_0^{}(t,x)+
\sum\limits_{j=1}^{m}\gamma_j^{}\;\!M_j^{}(t,x)
\quad
\forall (t,x)\in \Xi^{\;\!\prime}.
\hfill
$
\\[1.5ex]
\indent
{\sl Доказательство}.
Последовательно используя теорему 5.3, свойства 3.5 (при $m=1,$ $\lambda_1^{}=1),$ 
3.2, 3.4, получаем:
\\[1.5ex]
\mbox{}\hfill                           
$
\displaystyle
\biggl(
\exp\;\!\dfrac{q_j^{}}{p_j^{\;\!h_j^{}}}\,,\, \rho_{j}^{}\;\!N_0^{}\biggr)\in 
\text{\rm E}_{_{\Omega_{{}_{\tiny\;\! 0}}}}
\ \Longrightarrow \
\biggl(
\exp\;\!\dfrac{q_j^{}}{ \rho_{j}^{}\;\! p_j^{\;\!h_j^{}}}\,,\, N_0^{}\biggr)\in 
\text{\rm E}_{_{\Omega_{{}_{\tiny\;\! 0}}}}
\ \Longrightarrow \
\hfill                           
$
\\[2ex]
\mbox{}\hfill                           
$
\displaystyle
\Longrightarrow \
\biggl(
\exp\biggl(\;\!\dfrac{q_j^{}}{ \rho_{j}^{}\;\! p_j^{\;\!h_j^{}}}+\varphi\biggr),\, 
N_0^{}+{\sf D}\;\!\varphi\biggr)\in 
\text{\rm E}_{_{\Omega_{{}_{\tiny\;\! 0}}}},\ j=1,\ldots, m,
\ \Longrightarrow 
\hfill                           
$
\\[2ex]
\mbox{}\hfill                           
$
\displaystyle
\Longrightarrow \
\biggl(\ \sum\limits_{j=1}^{m}\lambda_j^{}
\exp\biggl(\;\!\dfrac{q_j^{}}{ \rho_{j}^{}\;\! p_j^{\;\!h_j^{}}}+\varphi\biggr),\, 
N_0^{}+{\sf D}\;\!\varphi\biggr)\in 
\text{\rm J}_{_{\Omega_{{}_{\tiny\;\! 0}}}}.
\hfill
$
\\[1.5ex]
\indent
Теперь, учитывая, что
\vspace{0.5ex} 
$\bigl(p_j^{}, M_j^{}\bigr)\in \text{\rm A}_{_{\Xi^{\;\!\prime}}},\ j=1,\ldots, m,$ 
на основании свойства 1.9 получаем утверждение доказываемого свойства. $\k$
\vspace{1.25ex} 

{\bf Свойство 5.6.}
\vspace{0.75ex}
{\it
Пусть 
$\bigl((p_j^{}, M_j^{}), (h_j^{},q_j^{},\rho_{j}^{}N_0^{})\bigr)\in 
\text{\rm B}_{_{\Xi^{\;\!\prime}}},\ 
\rho_j^{}\in\R\backslash\{0\},\ 
c_j^{}, \lambda_j^{},\gamma_j^{}\in\R,$  
$j=1,\ldots, m,\
\sum\limits_{j=1}^{m}|\lambda_j^{}|\ne 0,$ 
мно\-жес\-т\-во
\vspace{0.75ex}
$\Omega_{_0}\subset\Xi^{\;\!\prime}$ та\-кое, что  
$\prod\limits_{j=1}^{m}p_j^{}(t,x)\ne 0\;\;\forall (t,x)\in\Omega_{_0},$
$\prod\limits_{j=1}^{m}p_j^{}(t,x)= 0\;\;
\forall (t,x)\in {\sf C}_{_{\Xi^{\;\!\prime}}}\Omega_{_0},\ 
p_j^{\gamma_j^{}}\in C^1\Omega_{_0}, \ j=1,\ldots, m.$
Тогда
\\[1.5ex]
\mbox{}\hfill                           
$
\displaystyle
\biggl(\ \prod\limits_{j=1}^{m}p_j^{\gamma_j^{}}
\sum\limits_{j=1}^{m}\lambda_j^{}
\exp\biggl(\;\!\dfrac{q_j^{}}{\rho_j^{}\;\!p_j^{\;\!h_j^{}}}+c_j^{}\biggr), L\biggr)\in 
\text{\rm J}_{_{\Omega_{{}_{\tiny\;\! 0}}}},
\hfill
$
\\[2ex]
если и только если сомножитель}
\\[1.75ex]
\mbox{}\hfill                           
$
\displaystyle
L(t,x)=N_0^{}(t,x)+
\sum\limits_{j=1}^{m}\gamma_j^{}\;\!M_j^{}(t,x)
\quad
\forall (t,x)\in \Xi^{\;\!\prime}.
\hfill
$
\\[1.5ex]
\indent
{\sl Доказательство}
аналогично доказательству свойства 5.5 с той лишь разницей, что вместо свойства 3.2
используется свойство 3.3. $\k$
\vspace{1.25ex} 

{\bf Свойство 5.7.}\!
\vspace{0.75ex}
{\it
Пусть 
$\!\bigl((p_j^{}, M_j^{}), (h_j^{},q_j^{},N_j^{})\bigr)\!\in\! \text{\rm B}_{_{\Xi^{\;\!\prime}}},\, 
\gamma_j^{},\xi_j^{}\!\in\!\R,\, \varphi_j^{}\!\in\! C^1T^{\;\!\prime},\, 
j\!=\!1,\ldots, m,\!$ 
мно\-жес\-т\-во
\vspace{0.75ex}
$\!\Omega_{_0}\!\subset\!\Xi^{\;\!\prime}\!$ та\-кое, что  
$\prod\limits_{j=1}^{m}p_j^{}(t,x)\!\ne\! 0\;\;\forall (t,x)\!\in\!\Omega_{_0},\
\prod\limits_{j=1}^{m}p_j^{}(t,x)= 0\;\;
\forall (t,x)\!\in\! {\sf C}_{_{\Xi^{\;\!\prime}}}\Omega_{_0},\!$ 
$p_j^{\gamma_j^{}}\in C^1\Omega_{_0}, \ j=1,\ldots, m.$
Тогда
\\[1.5ex]
\mbox{}\hfill                           
$
\displaystyle
\biggl(\ \prod\limits_{j=1}^{m}p_j^{\gamma_j^{}}
\exp\sum\limits_{j=1}^{m}\xi_j^{}
\biggl(\;\!\dfrac{q_j^{}}{p_j^{\;\! h_j^{}}}+\varphi_j^{}\biggr), L\biggr)\in 
\text{\rm J}_{_{\Omega_{{}_{\tiny\;\! 0}}}},
\hfill
$
\\[2ex]
если и только если сомножитель}
\\[1.75ex]
\mbox{}\hfill                           
$
\displaystyle
L(t,x)=
\sum\limits_{j=1}^{m}\bigl(
\gamma_j^{}\;\!M_j^{}(t,x)+
\xi_j^{}\;\!\bigl(N_j^{}(t,x)+{\sf D}\;\!\varphi_j^{}(t)\bigr)\bigr)
\quad
\forall (t,x)\in \Xi^{\;\!\prime}.
\hfill
$
\\[1.5ex]
\indent
{\sl Доказательство}.
По свойству 5.1,
\\[1.5ex]
\mbox{}\hfill                           
$
\displaystyle
\biggl(
\exp\biggl(\;\!\dfrac{q_j^{}}{p_j^{\;\!h_j^{}}}+\varphi_j^{}\biggr),\, 
N_j^{}+{\sf D}\;\!\varphi_j^{}\biggr)\in 
\text{\rm E}_{_{\Omega_{{}_{\tiny\;\! 0}}}},
\ \ 
j=1,\ldots, m.
\hfill                           
$
\\[1.5ex]
\indent
Учитывая, что
\vspace{0.5ex} 
$\bigl(p_j^{}, M_j^{}\bigr)\in \text{\rm A}_{_{\Xi^{\;\!\prime}}},\ j=1,\ldots, m,$ 
на основании свойства 3.9 получаем утверждение доказываемого свойства. $\k$
\vspace{0.75ex} 

Кратность полиномиального частного интеграла $p$ в соответствии с 
определением 5.1 зависит от количества натуральных чисел $h$ и 
соответствующих этим числам функций 
$q\in \text{\rm P}_{_{\!\Xi^{\;\!\prime}}}$ 
таких, что выполняется тождество (5.1). 
А если учесть теорему 5.3, то кратность полиномиального частного интеграла $p$
определяется количеством экспоненциальных частных интегралов
$\exp\;\!\dfrac{q}{p^{\;\!h}}\;\!.$
\vspace{0.75ex} 

{\bf Определение 5.2.}
{\it
Полиномиальный частный интеграл $p$  на области $\Xi^{\;\!\prime}$ системы {\rm (0.1)} 
назовем \textit{\textbf{кратным с кратностью}}
$
\varkappa =1 + \sum\limits_{\xi=1}^{\varepsilon} \delta_{\xi}^{},
$
если  существуют  такие натуральные числа 
\vspace{0.5ex} 
$h_{\xi}^{},\ \xi =1,\ldots, \varepsilon,$ и соответствующие этим числам функции 
$
q_{_{\scriptstyle h_\xi^{}f_\xi^{}}}\in \text{\rm P}_{_{\!\Xi^{\;\!\prime}}},\ 
f_{\xi}^{}=1,\ldots,\delta_{\xi}^{},\ \xi=1,\ldots,\varepsilon,$
\vspace{0.35ex} 
каждая из которых является взаимно простой с функцией $p,$
что выполняются тождества 
\\[1.75ex]
\mbox{}\hfill                    % (5.3)
$
\displaystyle
{\frak d}\, \dfrac{q_{_{\scriptstyle h_\xi^{}  f_\xi^{}} }(t,x)}{\displaystyle  p^{\;\!h_\xi^{}} (t,x)}=
N_{h_\xi^{}  f_\xi^{}}^{} (t,x)
\quad 
\forall (t,x)\in \Omega_{_0},
\quad 
f_{\xi}^{}=1,\ldots, \delta_{\xi}^{}, \ \   
\xi=1,\ldots,\varepsilon,
$
\hfill {\rm (5.3)}
\\[2.25ex]
где функ\-ции 
\vspace{0.75ex}
$N_{h_\xi^{}  f_\xi^{}}^{}\in\text{\rm P}_{_{\!\Xi^{\;\!\prime}}}$ 
имеют степени
$
\deg_{\;\!x}^{}  N_{h_\xi^{}  f_\xi^{}}^{}\leq d-1, \ 
f_{\xi}^{}=1,\ldots, \delta_{\xi}^{}, \,  \xi=1,\ldots,\varepsilon,
$
мно\-жес\-т\-во
$\Omega_{_0}\subset\Xi^{\;\!\prime}$ та\-кое, что  
$p(t,x)\ne 0\;\;\forall (t,x)\in\Omega_{_0},\
p(t,x)= 0\;\;\forall (t,x)\in {\sf C}_{_{\Xi^{\;\!\prime}}}\Omega_{_0}.$
}
\vspace{1.25ex}

На основании определения 5.2 получаем
\vspace{1ex}

{\bf Предложение 5.1.}\!
\vspace{0.5ex}
{\it
Если  
$\!\Bigl((p, M), 
\Bigl(h_\xi^{},q_{_{\scriptstyle h_\xi^{}  f_\xi^{}}}, N_{h_\xi^{}  f_\xi^{}}^{}\Bigr)\Bigr)
\!\in\! \text{\rm B}_{_{\Xi^{\;\!\prime}}},\,
f_{\xi}^{}\!=\!1,\ldots, \delta_{\xi}^{},\, \xi\!=\!1,\ldots,\varepsilon,\!$
то  полиномиальный частный интеграл $p$ с сомножителем $M$
на области $\Xi^{\;\!\prime}$ системы {\rm (0.1)} является кратным с 
кратностью 
$
\varkappa =1 + \sum\limits_{\xi=1}^{\varepsilon} \delta_{\xi}^{}.
$
}
\vspace{0.5ex}

Обратим внимание на то, 
\vspace{0.15ex}
что в определении 5.2 и в предложении 5.1 
предусмотрена возможность, когда одному числу $h_\xi^{}\in\N$
соответствует $\delta_{\xi}^{}$ функций
$q_{_{\scriptstyle h_{\xi}^{}  f_\xi^{}}}\in \text{\rm P}_{_{\!\Xi^{\;\!\prime}}},
\linebreak 
f_{\xi}^{}=1,\ldots, \delta_{\xi}^{}.$
Также не исключается, что
\\[1.5ex]
\mbox{}\hfill
$
q_{_{\scriptstyle h_\xi^{}  f_\xi^{}}}(t,x)=
q_{_{\scriptstyle h_\zeta^{}  f_\zeta^{}}}(t,x)
\quad 
\forall (t,x)\in \Xi^{\;\!\prime}
$
\ при \ $\zeta\ne \xi,\ \ \xi,\zeta\in\{1,\ldots,\varepsilon\}.
\hfill
$
\\[2ex]
\indent
{\bf Свойство 5.8.}
\vspace{0.75ex}
{\it
Пусть 
$\Bigl((p, M), 
\Bigl(h_\xi^{},q_{_{\scriptstyle h_\xi^{}  f_\xi^{}}}, N_{h_\xi^{}  f_\xi^{}}^{}\Bigr)\Bigr)
\in \text{\rm B}_{_{\Xi^{\;\!\prime}}},\
\gamma,\gamma_{_{\scriptstyle h_\xi^{}  f_\xi^{}}}\in\R,\ 
\varphi_{_{\scriptstyle h_\xi^{}  f_\xi^{}}}\in C^1T^{\;\!\prime},$ 
$f_{\xi}^{}=1,\ldots, \delta_{\xi}^{},\ \xi=1,\ldots,\varepsilon,$
\vspace{1ex}
мно\-жес\-т\-во
$\Omega_{_0}\subset\Xi^{\;\!\prime}$ та\-кое, что  
$p(t,x)\ne 0\;\;\forall (t,x)\in\Omega_{_0},$
$p(t,x)= 0\;\;\forall (t,x)\in {\sf C}_{_{\Xi^{\;\!\prime}}}\Omega_{_0},\ 
p^{\gamma}\in C^1\Omega_{_0}.$
Тогда
\\[1.5ex]
\mbox{}\hfill                           
$
\displaystyle
\biggl(\ p^{\gamma}
\exp\sum\limits_{\xi=1}^{\varepsilon}\sum\limits_{f_{\xi}^{}=1}^{\delta_{\xi}^{}}
\biggl(\;\!\gamma_{_{\scriptstyle h_\xi^{}  f_\xi^{}}}\biggl(\,
\dfrac{q_{_{\scriptstyle h_\xi^{}  f_\xi^{}} }}{\displaystyle  p^{\;\!h_\xi^{}}}+
\varphi_{_{\scriptstyle h_\xi^{}  f_\xi^{}}}\biggr)\biggr),\, L\biggr)\in 
\text{\rm J}_{_{\Omega_{{}_{\tiny\;\! 0}}}},
\hfill
$
\\[2ex]
если и только если сомножитель}
\\[1.75ex]
\mbox{}\hfill                           
$
\displaystyle
L(t,x)=\gamma\;\!M(t,x)+
\sum\limits_{\xi=1}^{\varepsilon}\sum\limits_{f_{\xi}^{}=1}^{\delta_{\xi}^{}}
\biggl(\gamma_{_{\scriptstyle h_\xi^{}  f_\xi^{}}}\Bigl(
N_{h_\xi^{}  f_\xi^{}}^{} (t,x)+{\sf D}\;\!\varphi_{_{\scriptstyle h_\xi^{}  f_\xi^{}}}(t)\Bigr)\biggr)
\quad
\forall (t,x)\in \Xi^{\;\!\prime}.
\hfill
$
\\[1.5ex]
\indent
{\sl Доказательство}.
По свойству 5.1,
\\[1.5ex]
\mbox{}\hfill                           
$
\displaystyle
\biggl(
\exp\biggl(\;\!
\dfrac{q_{_{\scriptstyle h_\xi^{}  f_\xi^{}} }}{\displaystyle  p^{\;\!h_\xi^{}}}+
\varphi_{_{\scriptstyle h_\xi^{}  f_\xi^{}}}\biggr),\, 
N_{h_\xi^{}  f_\xi^{}}^{}+{\sf D}\;\!\varphi_{_{\scriptstyle h_\xi^{}  f_\xi^{}}}
\biggr)\in 
\text{\rm E}_{_{\Omega_{{}_{\tiny\;\! 0}}}},
\ \ 
f_{\xi}^{}=1,\ldots, \delta_{\xi}^{},\ \xi=1,\ldots,\varepsilon.
\hfill                           
$
\\[1.5ex]
\indent
Учитывая, что
\vspace{0.35ex} 
$\bigl(p, M\bigr)\in \text{\rm A}_{_{\Xi^{\;\!\prime}}},$ 
на основании свойства 3.9 получаем утверждение доказываемого свойства. $\k$
\vspace{0.75ex} 

Свойство 5.8 является аналогом свойства 5.7 на случай $\varkappa\!$-кратного 
частного интеграла $\bigl(p, M\bigr)\in \text{\rm A}_{_{\Xi^{\;\!\prime}}}.$
\vspace{0.35ex} 
Подобный образом на случай  $\varkappa\!$-кратного 
частного интеграла $\bigl(p, M\bigr)\in \text{\rm A}_{_{\Xi^{\;\!\prime}}}$
вводятся аналоги свойств 5.3, 5.4, 5.5 и 5.6.
\vspace{0.75ex} 

{\bf Свойство 5.9.}
{\it
Пусть 
$k\in\N,\ \varphi\in C^1T^{\;\!\prime}.$ Тогда
\\[1.5ex]
\mbox{}\hfill                           
$
\bigl(p, \varphi+p^kM_{0}^{}\bigr)\in \text{\rm A}_{_{\Xi^{\;\!\prime}}}
\iff
\bigl(\bigl(p, \varphi+p^kM_{0}^{}\bigr), \bigl(k,q,-k\;\!qM_0^{}\bigr)\bigr)
\in \text{\rm B}_{_{\Xi^{\;\!\prime}}},
\hfill
$
\\[1ex]
где
\\[1ex]
\mbox{}\hfill                   % (5.4)        
$
\displaystyle
q\colon t\to\ 
\exp\biggl(\, k\int\limits_{t_0^{}}^{t}\varphi(\tau)\;\!d\tau\biggr)
\quad
\forall t\in T^{\;\!\prime},
$
\hfill {\rm(5.4)}
\\[1.5ex]
$t_0^{}$ --- произвольная фиксированная точка из области $T^{\;\!\prime}.$
}
\vspace{0.5ex}

{\sl Доказательство}.
В соответствии с теоремой 2.1
\\[1.5ex]
\mbox{}\hfill                           
$
\displaystyle
{\frak d}\;\!p(t,x)=p(t,x)\;\!\bigl(\varphi(t)+p^{\;\!k}(t,x)\;\!M_0^{}(t,x)\bigr)
\quad
\forall (t,x)\in \Xi^{\;\!\prime},
\hfill
$
\\[1.75ex]
где функция
\vspace{0.5ex}
$M_0^{}\in \text{\rm P}_{_{\!\Xi^{\;\!\prime}}}$
и имеет степень $\deg_{\;\!x}^{} M_0^{}\leq d-k\deg_{\;\!x}^{}p-1.$

Тогда производная в силу системы (0.1)
\\[1.5ex]
\mbox{}\hfill                           
$
\displaystyle
{\frak d}\,\dfrac{\displaystyle
\exp\biggl(\, k\int\limits_{t_0^{}}^{t}\varphi(\tau)\;\!d\tau\biggr)}{p^{\;\!k}(t,x)}=
\hfill
$
\\[2ex]
\mbox{}\hfill                           
$
\displaystyle
=\dfrac{\displaystyle
p^{\;\!k}(t,x)\, {\sf D}\exp\biggl(\, k\int\limits_{t_0^{}}^{t}\varphi(\tau)\;\!d\tau\biggr) -
k\;\!p^{\;\!k-1}(t,x)\;\!{\frak d}\;\!p(t,x)\, \exp\biggl(\, k\int\limits_{t_0^{}}^{t}\varphi(\tau)\;\!d\tau\biggr)}
{p^{\;\!2k}(t,x)}=
\hfill
$
\\[2ex]
\mbox{}\hfill                           
$
\displaystyle
=
\dfrac{\displaystyle
k\;\!\varphi(t)\;\!p(t,x)\;\!\exp\biggl(\, k\int\limits_{t_0^{}}^{t}\varphi(\tau)\;\!d\tau\biggr) -
k\;\!p(t,x)\;\!\bigl(\varphi(t)+p^{\;\!k}(t,x)\;\!M_0^{}(t,x)\bigr)\;\!
\exp\biggl(\, k\int\limits_{t_0^{}}^{t}\varphi(\tau)\;\!d\tau\biggr)}
{p^{\;\!k+1}(t,x)}=
\hfill
$
\\[2ex]
\mbox{}\hfill                           
$
\displaystyle
={}-k\;\!M_0^{}(t,x)\;\!\exp\biggl(\, k\int\limits_{t_0^{}}^{t}\varphi(\tau)\;\!d\tau\biggr)
\quad
\forall (t,x)\in\Omega_{_{0}}^{},
\hfill
$
\\[2ex]
где мно\-жес\-т\-во
\vspace{1ex}
$\Omega_{_0}\subset\Xi^{\;\!\prime}$ та\-кое, что  
$p(t,x)\ne 0\;\;\forall (t,x)\in\Omega_{_0},$
$p(t,x)= 0\;\;\forall (t,x)\in {\sf C}_{_{\Xi^{\;\!\prime}}}\Omega_{_0}.\ \k$

В частности, при $\varphi(t)=0\;\;\forall t\in T^{\;\!\prime}$ имеет место
\vspace{0.75ex}

{\bf Свойство 5.10.}
{\it
Пусть 
$k\in\N,\ c\in\R\backslash\{0\}.$ Тогда}
\\[1.5ex]
\mbox{}\hfill                           
$
\bigl(p, p^kM_{0}^{}\bigr)\in \text{\rm A}_{_{\Xi^{\;\!\prime}}}
\iff
\bigl(\bigl(p, p^kM_{0}^{}\bigr), \bigl(k, c,{}-k\;\!c\;\!M_0^{}\bigr)\bigr)
\in \text{\rm B}_{_{\Xi^{\;\!\prime}}}.
\hfill
$
\\[2ex]
\indent
{\bf Свойство 5.11.}
\vspace{0.5ex}
{\it
Если  
$k\in\N,\ \bigl(p, p^kM_{0}^{}\bigr)\in \text{\rm A}_{_{\Xi^{\;\!\prime}}},$ 
то полиномиальный частный интеграл $p$ на области $\Xi^{\;\!\prime}$ системы {\rm (0.1)} 
является $(k+1)\!\!$-кратным таким, что 
\\[1.5ex]
\mbox{}\hfill                           
$
\bigl(\bigl(p, p^kM_{0}^{}\bigr), \bigl(l, c_l^{},{}-l\;\!c_l^{}\;\!p^{\;\!k-l}M_0^{}\bigr)\bigr)
\in \text{\rm B}_{_{\Xi^{\;\!\prime}}},
\ \ 
l=1,\ldots, k,
\hfill
$
\\[1.5ex]
где $\!c_l^{},\, l\!=\!1,\ldots,k,\!$ суть произвольные фиксированные ненулевые действительные числа.}

{\sl Доказательство.}
Если  
$\bigl(p, p^kM_{0}^{}\bigr)\in \text{\rm A}_{_{\Xi^{\;\!\prime}}},$ то
$\bigl(p, p^{\;\!l}M_{l}^{}\bigr)\in \text{\rm A}_{_{\Xi^{\;\!\prime}}},\ l=1,\ldots, k,$
где функции
\\[2ex]
\mbox{}\hfill                           
$
M_{l}^{}(t,x)=p^{\;\!k-l}(t,x)\;\!M_0^{}(t,x)
\quad
\forall (t,x)\in \Xi^{\;\!\prime},
\quad 
l=1,\ldots, k.
\hfill
$
\\[2ex]
\indent
По свойству 5.10, 
\vspace{0.75ex}
$
\bigl(\bigl(p, p^{\;\!l}M_{l}^{}\bigr), \bigl(l, c_l^{},{}-l\;\!c_l^{}\;\!M_l^{}\bigr)\bigr)
\in \text{\rm B}_{_{\Xi^{\;\!\prime}}},
\ l=1,\ldots, k,
$
где $c_l^{}\in\R\backslash\{0\},
\linebreak 
l=1,\ldots,k.$ 
Отсюда, 
\vspace{1.25ex} 
$
\bigl(\bigl(p, p^{\;\!k}M_{0}^{}\bigr), \bigl(l, c_l^{},{}-l\;\!c_l^{}\;\!p^{\;\!k-l}M_0^{}\bigr)\bigr)
\in \text{\rm B}_{_{\Xi^{\;\!\prime}}},\
c_l^{}\in\R\backslash\{0\},\ l=1,\ldots,k.\ \k$

{\bf Свойство 5.12.}
\vspace{0.5ex} 
{\it
Пусть 
$k\in\N,\ \varphi\in C^1T^{\;\!\prime},$ 
функция $q$ задана формулой {\rm (5.4)},
мно\-жес\-т\-во
$\Omega_{_0}\subset\Xi^{\;\!\prime}$ та\-кое, что  
$p(t,x)\ne 0\;\; \forall (t,x)\in\Omega_{_0},\
p(t,x)= 0\;\;\forall (t,x)\in {\sf C}_{_{\Xi^{\;\!\prime}}}\Omega_{_0}.$
Тогда}
\\[1.5ex]
\mbox{}\hfill                           
$
\bigl(p, \varphi+p^{\;\!k}M_{0}^{}\bigr)\in \text{\rm A}_{_{\Xi^{\;\!\prime}}}
\iff
\biggl(\exp\dfrac{q}{p^{\;\!k}}\,,\, {}-k\;\!q\;\!M_0^{}\biggr)\in
\text{\rm E}_{_{\Omega_{{}_{\tiny\;\! 0}}}}.
\hfill
$
\\[2ex]
\indent
{\sl Следует}
с учетом теоремы 5.3 из свойства 5.9. $\k$
\vspace{0.75ex}

{\bf Свойство 5.13.}
\vspace{0.5ex} 
{\it
Пусть 
$k\in\N,\ c\in\R\backslash\{0\},$ 
мно\-жес\-т\-во
$\Omega_{_0}\!\subset\!\Xi^{\;\!\prime}$ та\-кое, что  
$p(t,x)\!\ne\! 0$ $\forall (t,x)\in\Omega_{_0},\
p(t,x)= 0\;\;\forall (t,x)\in {\sf C}_{_{\Xi^{\;\!\prime}}}\Omega_{_0}.$
Тогда}
\\[1.25ex]
\mbox{}\hfill                           
$
\bigl(p, p^{\;\!k}M_{0}^{}\bigr)\in \text{\rm A}_{_{\Xi^{\;\!\prime}}}
\iff
\biggl(\exp\dfrac{c}{p^{\;\!k}}\,,\, {}-k\;\!c\;\!M_0^{}\biggr)\in
\text{\rm E}_{_{\Omega_{{}_{\tiny\;\! 0}}}}.
\hfill
$
\\[1.5ex]
\indent
{\sl Следует}
с учетом теоремы 5.3 из свойства 5.10. $\k$
\vspace{0.75ex}

{\bf Свойство 5.14.}
\vspace{0.5ex} 
{\it
Пусть 
$k\in\N,\ c_l^{}\in\R\backslash\{0\},\ l=1,\ldots, k,$ 
мно\-жес\-т\-во
$\Omega_{_0}\subset\Xi^{\;\!\prime}$ та\-кое, что  
$p(t,x)\ne 0\;\;\forall (t,x)\in\Omega_{_0},\
p(t,x)= 0\;\;\forall (t,x)\in {\sf C}_{_{\Xi^{\;\!\prime}}}\Omega_{_0}.$
Тогда, если}
$\bigl(p, p^{\;\!k}M_{0}^{}\bigr)\in \text{\rm A}_{_{\Xi^{\;\!\prime}}},$ то
\\[1.5ex]
\mbox{}\hfill                           
$
\biggl(\exp\dfrac{c_l^{}}{p^{\;\!l}}\,,\, {}-l\;\!c_{l}^{}\;\!p^{\;\!k-l}\;\!M_0^{}\biggr)\in
\text{\rm E}_{_{\Omega_{{}_{\tiny\;\! 0}}}},
\ \ l=1,\ldots, k.
\hfill
$
\\[1.75ex]
\indent
{\sl Следует}
с учетом теоремы 5.3 из свойства 5.11. $\k$
\\[3.25ex]
\centerline{
{\bf  6. Комплекснозначные полиномиальные частные интегралы}
}
\\[1.25ex]
\indent
Множество функций, являющихся  полиномами по переменным $x_1^{},\ldots,x_n^{}$ 
с комплекснозначными коэффициентами-функциями одной переменной $t$ 
непрерывно дифференцируемыми на области $T,$ обозначим $\text{Z}_{_{\Xi}}.$ 
\vspace{0.5ex}

{\bf Определение 6.1.}
\vspace{0.25ex}
{\it
Функцию $w\in \text{\rm Z}_{_{\Xi^{\;\!\prime}}}$ назовем
\textit{\textbf{комплекснозначным полиномиальным частным интегралом на области}}
$\Xi^{\;\!\prime}$ сис\-темы {\rm (0.1)}, если производная в силу системы {\rm (0.1)}
\\[1ex]
\mbox{}\hfill                           % (6.1)
$
{\frak d}\;\!w(t,x)=w(t,x)\;\!W(t,x)
\quad
\forall (t,x)\in \Xi^{\;\!\prime}.
$
\hfill {\rm (6.1)}
\\[1.5ex]
При этом функцию $W$ будем называть \textit{\textbf{сомножителем}}
комплекснозначного полиномиального частного интеграла $w.$
}
\vspace{0.5ex}

{\bf Замечание 6.1.}
\vspace{0.25ex}
Так как функция $w\in \text{Z}_{_{\Xi^{\;\!\prime}}},$
то из тождества (6.1) следует, что сомножитель 
$W\in \text{Z}_{_{\Xi^{\;\!\prime}}};$
случай $W\in \text{P}_{_{\!\Xi^{\;\!\prime}}}$ не исключается.
\vspace{0.5ex}

После введения понятия комплекснозначного полиномиального частного интеграла
о полиномиальном частном интеграле в смысле определения 2.1 можно (в случае необходимости)
говорить как о вещественном полиномиальном частном интеграле.

Множество функций, являющихся комплекснозначными полиномиальными частными интегралами 
на области $\Xi^{\;\!\prime}$ сис\-темы (0.1), обозначим 
$\text{H}_{_{\Xi^{\;\!\prime}}}.$
\vspace{0.25ex}
Условной записью 
$(w, W)\in \text{H}_{_{\Xi^{\;\!\prime}}}$ 
\vspace{0.35ex}
будем выражать, что функция $w$
является комплекснозначным полиномиальным частным интегралом с сомножителем $W$
\vspace{0.5ex}
на области $\Xi^{\;\!\prime}$ сис\-темы (0.1).

{\bf Свойство 6.1.} 
{\it
Пусть $\eta\in\C\backslash\{0\}.$ Тогда}
\\[1.5ex]
\mbox{}\hfill
$
(w, W)\in \text{\rm H}_{_{\Xi^{\;\!\prime}}}
\iff 
(\eta\;\! w\;\!,\;\! W)\in \text{\rm H}_{_{\Xi^{\;\!\prime}}}.
\hfill
$
\\[1.25ex]
\indent
{\sl Доказательство}
основано на определении 6.1 и том, что при любом $\eta\in\C\backslash\{0\}$
\\[1ex]
\mbox{}\hfill                          
$
{\frak d}\;\!\bigl(\eta\;\!w(t,x)\bigr)=\eta\,{\frak d}\;\!w(t,x)
\quad
\forall (t,x)\in \Xi^{\;\!\prime}.\ \k
\hfill
$
\\[-3ex]

\newpage

{\bf Свойство 6.2.} 
{\it
Если $\eta\in\C,\ w\in \text{\rm Z}_{_{\Xi^{\;\!\prime}}}$ и выполняется тождество
\\[1.5ex]
\mbox{}\hfill
$
{\frak d}\;\!w(t,x)=\bigl(w(t,x)+\eta\bigr)\;\!W(t,x)
\quad
\forall (t,x)\in \Xi^{\;\!\prime},
\hfill
$
\\[1.5ex]
то $(w+\eta, W)\in \text{\rm H}_{_{\Xi^{\;\!\prime}}}.$}
\vspace{0.5ex}

{\sl Действительно}, 
производная в силу системы (0.1)
\\[1.5ex]
\mbox{}\hfill                          
$
{\frak d}\;\!\bigl(w(t,x)+\eta\bigr)=
{\frak d}\;\!w(t,x)=\bigl(w(t,x)+\eta\bigr)\;\!W(t,x)
\quad
\forall (t,x)\in \Xi^{\;\!\prime}.
\hfill
$
\\[1.5ex]
\indent
Согласно определению 6.1
$(w+\eta, W)\in \text{\rm H}_{_{\Xi^{\;\!\prime}}}.\ \k$
\vspace{1.25ex}

{\bf Свойство 6.3.} 
{\it
Пусть $k\in\N.$ Тогда}
\\[1.5ex]
\mbox{}\hfill
$
(w, W)\in \text{\rm H}_{_{\Xi^{\;\!\prime}}}
\iff 
(w^k,\;\! k\;\!W)\in \text{\rm H}_{_{\Xi^{\;\!\prime}}}.
\hfill
$
\\[1.5ex]
\indent
{\sl Доказательство}
основано на определении 6.1 и том, что
\\[1.5ex]
\mbox{}\hfill                          
$
{\frak d}\;\!w^k(t,x)=k\;\!w^{k-1}(t,x)\,{\frak d}\;\!w(t,x)
\quad
\forall (t,x)\in \Xi^{\;\!\prime}.\ \k
\hfill
$
\\[2ex]
\indent
{\bf Свойство 6.4.} 
{\it
Если $(w_j^{}, W)\in \text{\rm H}_{_{\Xi^{\;\!\prime}}},\ 
\eta_j^{}\in\C\backslash\{0\},\ j=1,\ldots,m,$ то}
\\[1.5ex]
\mbox{}\hfill
$
\displaystyle
\biggl(\,\sum\limits_{j=1}^{m} \eta_j^{}\;\!w_j^{}\;\!,\;\! W\biggl)\;\! 
\in \text{\rm H}_{_{\Xi^{\;\!\prime}}}.
\hfill
$
\\[1.5ex]
\indent
{\sl Действительно}, производная в силу системы (0.1)
\\[1.5ex]
\mbox{}\hfill                         
$
\displaystyle
{\frak d}\;\!\sum\limits_{j=1}^{m} \eta_j^{}\;\!w_j^{}(t,x)  = 
\sum\limits_{j=1}^{m}
\eta_j^{}\, {\frak d}\;\!w_j^{}(t,x)=
\sum\limits_{j=1}^{m}
\eta_j^{}\;\!w_j^{}(t,x)\;\!W(t,x)
\quad
\forall (t,x)\in \Xi^{\;\!\prime}.
\hfill
$
\\[1.5ex]
\indent
Согласно определению 6.1  
$\biggl(\,\sum\limits_{j=1}^{m}\! \eta_j^{}\;\!w_j^{}\;\!,\;\! W\biggl) 
\in \text{\rm H}_{_{\Xi^{\;\!\prime}}}.\ \k$
\vspace{1.25ex}

{\bf Свойство 6.5.}
{\it
Пусть $w_j^{}\in\text{\rm Z}_{_{\Xi^{\;\!\prime}}},\ j=1,\ldots,m.$ Тогда
\\[1.5ex]
\mbox{}\hfill                         
$
\displaystyle
\biggl(\,\prod\limits_{j=1}^{m} w_j^{}\;\!,\;\! W\biggl)\;\! 
\in \text{\rm H}_{_{\Xi^{\;\!\prime}}}
\iff
\bigl( w_j^{}\;\!,\;\! W_j^{}\bigl)\;\! \in \text{\rm H}_{_{\Xi^{\;\!\prime}}},
\quad
j=1,\ldots,m.
\hfill
$
\\[1.5ex]
При этом сомножители $W,\, W_1^{},\ldots,W_m^{}$ такие, что}
\\[1.5ex]
\mbox{}\hfill
$
\displaystyle
W(t,x)=\sum\limits_{j=1}^{m} W_j^{}(t,x)
\quad
\forall (t,x)\in \Xi^{\;\!\prime}.
\hfill
$
\\[1.5ex]
\indent
{\sl Доказательство.}
\vspace{0.35ex}
Пусть функции $w_1^{}, w_2^{}\in\text{\rm Z}_{_{\Xi^{\;\!\prime}}}.$
Тогда согласно определению 6.1 
$\bigl(w_1^{}\;\!w_2^{}\;\!,\;\! W\bigl)\;\! \in \text{\rm H}_{_{\Xi^{\;\!\prime}}},$
если и только если производная в силу системы (0.1)
\\[1.75ex]
\mbox{}\hfill                          
$
\displaystyle
{\frak d}\;\! \bigl(w_1^{}(t,x)\;\!w_2^{}(t,x)\bigr) =
w_1^{}(t,x)\;\!w_2^{}(t,x)\;\!W(t,x)
\quad
\forall (t,x)\in \Xi^{\;\!\prime}.
\hfill
$
\\[1.25ex]
\indent
Так как
\\[1.25ex]
\mbox{}\hfill                          
$
\displaystyle
{\frak d}\;\! \bigl(w_1^{}(t,x)\;\!w_2^{}(t,x)\bigr)\;\! =
w_2^{}(t,x)\, {\frak d}\;\!w_1^{}(t,x)+
w_1^{}(t,x)\, {\frak d}\;\!w_2^{}(t,x)
\quad
\forall (t,x)\in \Xi^{\;\!\prime},
\hfill
$
\\[1.75ex]
то $\bigl(w_1^{}\;\!w_2^{}\;\!,\;\! W\bigl)\;\! \in \text{\rm H}_{_{\Xi^{\;\!\prime}}},$
если и только если 
\\[1.75ex]
\mbox{}\hfill                          
$
\displaystyle
w_2^{}(t,x)\, {\frak d}\;\!w_1^{}(t,x)+
w_1^{}(t,x)\, {\frak d}\;\!w_2^{}(t,x)=
w_1^{}(t,x)\;\!w_2^{}(t,x)\;\!W(t,x)
\quad
\forall (t,x)\in \Xi^{\;\!\prime}.
\hfill
$
\\[1.75ex]
\indent
Отсюда, 
$\bigl(w_1^{}\;\!w_2^{}\;\!,\;\! W\bigl)\;\! \in \text{\rm H}_{_{\Xi^{\;\!\prime}}}$
тогда и только тогда, когда
\\[1.5ex]
\mbox{}\hfill                          % (6.2)
$
\displaystyle
{\frak d}\;\! w_1^{}(t,x) =
w_1^{}(t,x)\;\! 
\biggl(
W(t,x)-\dfrac{{\frak d}\;\! w_2^{}(t,x)}{w_2^{}(t,x)}
\biggr)
\quad 
\forall (t,x)\in\Omega_{_0},
$
\hfill (6.2)
\\[2ex]
где множество $\!\Omega_{_0}\!\subset\! \Xi^{\;\!\prime}\!$ такое, что  
$|w_2^{}(t,x)|\ne 0\;\;\forall (t,x)\!\in\! \Omega_{_0},\ 
|w_2^{}(t,x)|=0\;\;\forall (t,x)\!\in\! {\sf C}_{_{\Xi^{\;\!\prime}}}\Omega_{_0}.\!$

\newpage

Пусть
\\[1.25ex]
\mbox{}\hfill                            
$
W(t,x)-\dfrac{{\frak d}\;\! w_2^{}(t,x)}{w_2^{}(t,x)}=
W_1^{}(t,x)
\quad 
\forall (t,x)\in \Omega_{_0}.
\hfill 
$
\\[1.75ex]
\indent
Поскольку $w_1^{},\ {\frak d}\;\! w_1^{}, \ W\in \text{\rm Z}_{_{\Xi^{\;\!\prime}}},$
то из тождества (6.2) следует, что $W_1^{}\in \text{\rm Z}_{_{\Xi^{\;\!\prime}}}.$
\vspace{0.5ex}

В итоге получаем, что 
$\bigl(w_1^{}\;\!w_2^{}\;\!,\;\! W\bigl)\;\! \in \text{\rm H}_{_{\Xi^{\;\!\prime}}},$
если и только если выполняются тождества
\\[1.75ex]
\mbox{}\hfill                          % (6.3)
$
\displaystyle
{\frak d}\;\! w_1^{}(t,x)=
w_1^{}(t,x)\;\! W_1^{}(t,x)
\quad 
\forall (t,x)\in \Xi^{\;\!\prime}
$
\hfill (6.3)
\\[0.5ex]
и
\\[0.5ex]
\mbox{}\hfill                           % (6.4)
$
\displaystyle
{\frak d}\;\! w_2^{}(t,x)=
w_2^{}(t,x)\;\! \bigl(W(t,x)-W_1^{}(t,x)\bigr)
\quad 
\forall (t,x)\in \Xi^{\;\!\prime}.
$
\hfill (6.4)
\\[2ex]
\indent
Согласно определению 6.1 тождество (6.3)
\vspace{0.35ex}
выполняется в том и только в том случае, когда 
$\bigl(w_1^{}\;\!,\;\! W_1^{}\bigl)\;\! \in \text{\rm H}_{_{\Xi^{\;\!\prime}}},$
\vspace{0.75ex}
а тождество (6.4) выполняется в том и только в том случае, когда 
$\bigl(w_2^{}\;\!,\;\! W_2^{}\bigl)\;\! \in \text{\rm H}_{_{\Xi^{\;\!\prime}}},\ 
W_2^{}=W-W_1^{}.$
\vspace{0.35ex}

Итак, свойство 6.5 при $m=2$ доказано.
\vspace{0.35ex}

При $m>2$ доказываем по индукции. $\k$
\vspace{0.75ex}

{\bf Свойство 6.6.} 
{\it
Пусть 
$k_j^{}\in\N, \ w_j^{}\in \text{\rm Z}_{_{\Xi^{\;\!\prime}}},\ j=1,\ldots,m.$ Тогда  
\\[1.5ex]
\mbox{}\hfill                         
$
\displaystyle
\biggl(\,\prod\limits_{j=1}^{m} w_j^{k_j^{}}\;\!,\;\! W\biggl)\;\! 
\in \text{\rm H}_{_{\Xi^{\;\!\prime}}}
\iff
\bigl( w_j^{}\;\!,\;\! W_j^{}\bigl)\;\! \in \text{\rm H}_{_{\Xi^{\;\!\prime}}},
\quad
j=1,\ldots,m.
\hfill
$
\\[1.5ex]
При этом сомножители $W,\, W_1^{},\ldots,W_m^{}$ такие, что}
\\[1.5ex]
\mbox{}\hfill
$
\displaystyle
W(t,x)=\sum\limits_{j=1}^{m} k_j^{}\;\!W_j^{}(t,x)
\quad
\forall (t,x)\in \Xi^{\;\!\prime}.
\hfill
$
\\[1.5ex]
\indent
{\sl Следует} с учетом свойства 6.3 из свойства 6.5. $\k$
\vspace{0.75ex}

{\bf Свойство 6.7.} 
{\it
Если 
$p_\tau^{}\in \text{\rm P}_{_{\!\Xi^{\;\!\prime}}},\, l_{\tau}^{}\!\in\N,\, \tau\!=\!1,\ldots,s,\
w_j^{}\in \text{\rm Z}_{_{\Xi^{\;\!\prime}}},\, k_j^{}\!\in\N, \, j\!=\!1,\ldots,m,$ то  
\\[1.5ex]
\mbox{}\hfill                         
$
\displaystyle
\biggl(\,
\prod\limits_{\tau=1}^{s}\! p_\tau^{\,l_\tau^{}}
\prod\limits_{j=1}^{m}\! w_j^{k_j^{}}, W\!\biggl) 
\in\! \text{\rm H}_{_{\Xi^{\;\!\prime}}}
\!\iff\!
\bigl(p_\tau^{}\;\!, M_\tau^{}\bigl)\;\! \in\! \text{\rm A}_{_{\Xi^{\;\!\prime}}},\,
\tau\!=\!1,\ldots,s,
\ \&\ 
\bigl(w_j^{}\;\!, W_j^{}\bigl)\;\!\in\! \text{\rm H}_{_{\Xi^{\;\!\prime}}},\,
j\!=\!1,\ldots,m.
\hfill
$
\\[1.5ex]
При этом сомножители $W,\, M_1^{},\ldots,M_s^{},\, W_1^{},\ldots,W_m^{}$ такие, что}
\\[1.5ex]
\mbox{}\hfill
$
\displaystyle
W(t,x)=
\sum\limits_{\tau=1}^{s} l_\tau^{}\;\!M_\tau^{}(t,x)+
\sum\limits_{j=1}^{m} k_j^{}\;\!W_j^{}(t,x)
\quad
\forall (t,x)\in \Xi^{\;\!\prime}.
\hfill
$
\\[1.5ex]
\indent
{\sl Следует} из свойства 6.6 с учетом того, что, если 
\vspace{1ex}
$p_\tau^{}\in \text{\rm P}_{_{\!\Xi^{\;\!\prime}}},$ то 
$p_\tau^{}\in \text{\rm Z}_{_{\Xi^{\;\!\prime}}},\ \tau=1,\ldots,s.\ \k$

У функций 
$w, W\in \text{\rm Z}_{_{\Xi^{\;\!\prime}}}$
выделим вещественную и мнимую части:
\\[1.5ex]
\mbox{}\hfill
$
u\colon (t,x)\to\ {\rm Re}\;\!w(t,x)
\quad
\forall (t,x)\in \Xi^{\;\!\prime},
\qquad
v\colon (t,x)\to\ {\rm Im}\;\!w(t,x)
\quad
\forall (t,x)\in \Xi^{\;\!\prime},
\hfill
$
\\[2ex]
\mbox{}\hfill
$
U\colon (t,x)\to\ {\rm Re}\;\!W(t,x)
\quad
\forall (t,x)\in \Xi^{\;\!\prime},
\qquad
V\colon (t,x)\to\ {\rm Im}\;\!W(t,x)
\quad
\forall (t,x)\in \Xi^{\;\!\prime}.
\hfill
$
\\[2ex]
\indent
{\bf Свойство 6.8.}\! 
\vspace{0.15ex}
{\it
Пусть 
$p\in \text{\rm P}_{_{\!\Xi^{\;\!\prime}}},\!$ 
функции $u_1^{},v_1^{}\!\in\! \text{\rm P}_{_{\!\Xi^{\;\!\prime}}}$ 
взаимно простые, $w_1^{}=u_1^{}+i\;\!v_1^{},$ $w=p\;\!w_1^{}.$
Тогда  
\\[1.5ex]
\mbox{}\hfill                         
$
\displaystyle
\bigl(w, W\bigl)\, \in \text{\rm H}_{_{\Xi^{\;\!\prime}}}
\iff
\bigl(p, M\bigl)\, \in \text{\rm A}_{_{\Xi^{\;\!\prime}}}
\ \&\ 
\bigl(w_1^{}, W_1^{}\bigl)\, \in \text{\rm H}_{_{\Xi^{\;\!\prime}}}.
\hfill
$
\\[1.75ex]
При этом сомножители $W,\, M,\, W_1^{}$ такие, что}
\\[1.5ex]
\mbox{}\hfill
$
\displaystyle
W(t,x)=M(t,x)+W_1^{}(t,x)
\quad
\forall (t,x)\in \Xi^{\;\!\prime}.
\hfill
$
\\[1.75ex]
\indent
{\sl Следует} из свойства 6.7 при $s=m=1.\ \k$
\vspace{0.75ex}

\newpage

С помощью свойства 6.8 из множества $\text{\rm H}_{_{\Xi^{\;\!\prime}}}$
выделяются комплекснозначные полиномиальные частные интегралы, 
\vspace{0.5ex}
у которых вещественная и мнимая части взаимно простые.

{\bf Теорема 6.1}
\vspace{0.15ex}
(критерий существования комплекснозначного полиномиального частного интеграла). 
{\it
Пусть 
$u,v\in \text{\rm P}_{_{\!\Xi^{\;\!\prime}}}.$ 
Тогда $(u+i\;\!v,\;\! U+i\;\!V)\in \text{\rm H}_{_{\Xi^{\;\!\prime}}},$
\vspace{0.25ex}
если и только если выполняется система тождеств 
\\[1.5ex]
\mbox{}\hfill                            % (6.5)
$
{\frak d}\;\!u(t,x)=
u(t,x)\;\!U(t,x)-v(t,x)\;\!V(t,x)
\quad
\forall (t,x)\in \Xi^{\;\!\prime},
$
\hfill {\rm (6.5)}
\\[2.25ex]
\mbox{}\hfill                           % (6.6)
$
{\frak d}\;\!v(t,x)=
u(t,x)\;\!V(t,x)+v(t,x)\;\!U(t,x)
\quad
\forall (t,x)\in \Xi^{\;\!\prime},
$
\hfill {\rm (6.6)}
\\[1.75ex]
в которой	функции $U, V\in \text{\rm P}_{_{\!\Xi^{\;\!\prime}}}.$}
\vspace{0.5ex}

{\sl Доказательство} 
состоит в том, что тождество (6.1) выполняется тогда и только тогда, когда
имеет место система тождеств $(6.5)\& (6.6),$ в которой 
$U, V\in \text{\rm P}_{_{\!\Xi^{\;\!\prime}}}.\ \k$
\vspace{0.75ex}

{\bf Свойство 6.9.}
{\it
Справедливо утверждение}
\\[1.5ex]
\mbox{}\hfill                        
$
(u+i\;\!v,\;\! U+i\;\!V)\in \text{\rm H}_{_{\Xi^{\;\!\prime}}}
\iff
(u-i\;\!v,\;\! U-i\;\!V)\in \text{\rm H}_{_{\Xi^{\;\!\prime}}}.
\hfill 
$
\\[1.75ex]
\indent
{\sl Доказательство} 
основано на теореме 6.1 и том, что система тождеств $(6.5)\& (6.6)$ 
инвариантна при одновременной замене $v$ на ${}-v$ и $V$ на ${}-V.\ \k$
\vspace{0.75ex}

{\bf Лемма 6.1.}
\vspace{0.35ex}
{\it
Пусть  $u,v\in \text{\rm P}_{_{\!\Xi^{\;\!\prime}}}.$ 
Тогда $(u+i\;\!v,\;\! U+i\;\!V)\in \text{\rm H}_{_{\Xi^{\;\!\prime}}},$
если и только если $(u^2+v^2,\;\! 2\;\!U)\in \text{\rm A}_{_{\Xi^{\;\!\prime}}}$
и выполняется тождество {\rm (6.5)} {\rm(}или {\rm (6.6)).}
}
\vspace{0.5ex}

{\sl Доказательство.}
По теореме 6.1, $(u+i\;\!v,\;\! U+i\;\!V)\in \text{\rm H}_{_{\Xi^{\;\!\prime}}},$
если и только если выполняются тождества (6.5) и (6.6).

Умножим тождество (6.5) на $2u,$ а тождество (6.6) на $2v;$
затем полученные тождества почленно сложим. 
В результате систему тождеств $(6.5)\& (6.6)$ 
приведем к системе, которая состоит из тождества 
\\[1.5ex]
\mbox{}\hfill                          
$
{\frak d}\;\!\bigl(u^2(t,x)+v^2(t,x)\bigr)=
2\;\!\bigl(u^2(t,x)+v^2(t,x)\bigr)\;\!U(t,x)
\quad
\forall (t,x)\in \Xi^{\;\!\prime}
\hfill
$
\\[1.5ex]
и хотя бы одного из тождеств (6.5) или (6.6).

Учитывая, что полученное тождество в соответствии с теоремой 2.1 
выполняется тогда и только тогда, когда 
$(u^2+v^2,\;\! 2\;\!U)\in \text{\rm A}_{_{\Xi^{\;\!\prime}}},$
завершаем доказательство. $\k$
\vspace{0.75ex}

{\bf Лемма 6.2.}
\vspace{0.35ex}
{\it
Пусть  $u,v\in \text{\rm P}_{_{\!\Xi^{\;\!\prime}}}$ взаимно простые, 
множество $\Omega_{_0}\subset \Xi^{\;\!\prime}$ такое, что  
$u(t,x)\ne 0\;\;\forall (t,x)\in \Omega_{_0},\ 
u(t,x)=0\;\;\forall (t,x)\in {\sf C}_{_{\Xi^{\;\!\prime}}}\Omega_{_0}.$
\vspace{0.35ex}
Тогда $(u+i\;\!v,\;\! U+i\;\!V)\in \text{\rm H}_{_{\Xi^{\;\!\prime}}},$
если и только если 
$\Bigl(\exp\arctg\dfrac{v}{u}\,,\;\! V\Bigr)\in 
\text{\rm E}_{_{\Omega_{{}_{\tiny\;\! 0}}}},$
\vspace{0.25ex}
а функция $U\in \text{\rm P}_{_{\!\Xi^{\;\!\prime}}}$
находится из тождества {\rm (6.5)} {\rm(}или {\rm (6.6)).}
}
\vspace{0.5ex}

{\sl Следует}
из теоремы 6.1 и свойства 3.11. $\k$
\vspace{0.75ex}

{\bf Теорема 6.2}
\vspace{0.15ex}
(критерий существования комплекснозначного полиномиального частного интеграла). 
{\it
Пусть функции
\vspace{0.5ex}
$u,v\in \text{\rm P}_{_{\!\Xi^{\;\!\prime}}}$ взаимно простые, 
множество $\Omega_{_0}\subset \Xi^{\;\!\prime}$ такое, что  
$u(t,x)\ne 0\;\;\forall (t,x)\in \Omega_{_0},\ 
u(t,x)=0\;\;\forall (t,x)\in {\sf C}_{_{\Xi^{\;\!\prime}}}\Omega_{_0}.$
Тогда }
\\[1.75ex]
\mbox{}\hfill                          
$
(u+i\;\!v,\;\! U+i\;\!V)\in \text{\rm H}_{_{\Xi^{\;\!\prime}}}
\iff
(u^2+v^2,\;\! 2\;\!U)\in \text{\rm A}_{_{\Xi^{\;\!\prime}}}
\ \&\ 
\Bigl(\exp\arctg\dfrac{v}{u}\,,\;\! V\Bigr)\in 
\text{\rm E}_{_{\Omega_{{}_{\tiny\;\! 0}}}}.
\hfill                          
$
\\[1.5ex]
\indent
{\sl Следует}
из лемм 6.1 и 6.2. $\k$
\vspace{0.75ex}

{\bf Свойство 6.10.}
\vspace{0.15ex}
{\it
Пусть функции
\vspace{0.5ex}
$u,v\in \text{\rm P}_{_{\!\Xi^{\;\!\prime}}}$ взаимно простые, 
множество $\Omega_{_0}\subset \Xi^{\;\!\prime}$ такое, что  
$u(t,x)\ne 0\;\;\forall (t,x)\in \Omega_{_0},\ 
u(t,x)=0\;\;\forall (t,x)\in {\sf C}_{_{\Xi^{\;\!\prime}}}\Omega_{_0}.$
Тогда 
\\[1.5ex]
\mbox{}\hfill
$
\Bigl(\exp\arctg\dfrac{v}{u}\,,\;\! V\Bigr)\in 
\text{\rm E}_{_{\Omega_{{}_{\tiny\;\! 0}}}},
\hfill
$
\\[1.5ex]
если и только если 
$(u^2+v^2,\;\! 2\;\!U)\in \text{\rm A}_{_{\Xi^{\;\!\prime}}}$
и выполняется тождество {\rm (6.5)} {\rm(}или {\rm (6.6)).}
}
\vspace{0.5ex}

{\sl Следует}
из лемм 6.1 и 6.2. $\k$
\vspace{0.75ex}

\newpage

{\bf Замечание 6.2.}
\vspace{0.35ex}
Согласно замечанию 2.1 и теореме 3.1 из теоремы 6.2 следует, что
$\deg_x^{} U\leq d-1$ и $\deg_x^{} V\leq d-1.$
\vspace{0.35ex}
Поэтому в определении 6.1 $\deg_x^{} W\leq d-1,$
а в тождествах (3.5), (6.5) и (6.6)
$\deg_x^{} U\leq d-1,\ \deg_x^{} V\leq d-1.$
\vspace{0.75ex}

{\bf Свойство 6.11.} 
\vspace{0.15ex}
{\it
Пусть 
$u, v, U\in \text{\rm P}_{_{\!\Xi^{\;\!\prime}}}.$
Тогда}
\\[1.5ex]
\mbox{}\hfill                         
$
\displaystyle
(u+i\;\!v,\;\! U)\in \text{\rm H}_{_{\Xi^{\;\!\prime}}}
\iff
\bigl(u, U\bigl)\, \in \text{\rm A}_{_{\Xi^{\;\!\prime}}}
\ \&\ 
\bigl(v, U\bigl)\, \in \text{\rm A}_{_{\Xi^{\;\!\prime}}}.
\hfill
$
\\[1.75ex]
\indent
{\sl Следует} 
из теоремы 6.1 с учетом теоремы 2.1 и того, что при 
$V(t,x)=0\;\;\forall (t,x)\in \Xi^{\;\!\prime}$
тождества (6.5) и (6.6) соответственно примут вид
\\[1.5ex]
\mbox{}\hfill                           
$
{\frak d}\;\!u(t,x)=
u(t,x)\;\!U(t,x)
\quad
\forall (t,x)\in \Xi^{\;\!\prime},
\qquad
{\frak d}\;\!v(t,x)=
v(t,x)\;\!U(t,x)
\quad
\forall (t,x)\in \Xi^{\;\!\prime}.\ \k
\hfill
$
\\[1.75ex]
\indent
{\bf Свойство 6.12.}
\vspace{0.15ex}
{\it
Пусть функции
\vspace{0.5ex}
$u,v\in \text{\rm P}_{_{\!\Xi^{\;\!\prime}}}$ взаимно простые, 
множество $\Omega_{_0}\subset \Xi^{\;\!\prime}$ такое, что  
$u(t,x)\ne 0\;\;\forall (t,x)\in \Omega_{_0},\ 
u(t,x)=0\;\;\forall (t,x)\in {\sf C}_{_{\Xi^{\;\!\prime}}}\Omega_{_0}.$
Тогда}
\\[1.5ex]
\mbox{}\hfill
$
(u+i\;\!v,\;\! U)\in \text{\rm H}_{_{\Xi^{\;\!\prime}}}
\iff
(u^2+v^2,\;\! 2\;\!U)\in \text{\rm A}_{_{\Xi^{\;\!\prime}}}
\ \&\ \,
\dfrac{v}{u}\in \text{\rm I}_{_{\Omega_{{}_{\tiny\;\! 0}}}}.
\hfill
$
\\[1.5ex]
\indent
{\sl Следует}
из теоремы 6.2 при $V(t,x)=0\;\;\forall (t,x)\in \Xi^{\;\!\prime}$
с учетом теоремы 0.2 и свойства~0.1. $\k$
\vspace{0.75ex}

{\bf Свойство 6.13.} 
\vspace{0.15ex}
{\it
Пусть 
$u, v, V\in \text{\rm P}_{_{\!\Xi^{\;\!\prime}}}.$
Тогда $(u+i\;\!v,\;\! i\;\!V)\in \text{\rm H}_{_{\Xi^{\;\!\prime}}},$
если и только если}
\\[1.5ex]
\mbox{}\hfill                         
$
{\frak d}\;\!u(t,x)=
{}-v(t,x)\;\!V(t,x)
\quad
\forall (t,x)\in \Xi^{\;\!\prime},
\qquad
{\frak d}\;\!v(t,x)=
u(t,x)\;\!V(t,x)
\quad
\forall (t,x)\in \Xi^{\;\!\prime}.
\hfill
$
\\[1.75ex]
\indent
{\sl Следует} 
из теоремы 6.1  при 
$U(t,x)=0\;\;\forall (t,x)\in \Xi^{\;\!\prime}.\ \k$
\vspace{1.25ex}

{\bf Свойство 6.14.}
\vspace{0.15ex}
{\it
Пусть функции
\vspace{0.5ex}
$u,v\in \text{\rm P}_{_{\!\Xi^{\;\!\prime}}}$ взаимно простые, 
$V\in \text{\rm P}_{_{\!\Xi^{\;\!\prime}}},$
множество $\Omega_{_0}\subset \Xi^{\;\!\prime}$ такое, что  
$u(t,x)\ne 0\;\;\forall (t,x)\in \Omega_{_0},\ 
u(t,x)=0\;\;\forall (t,x)\in {\sf C}_{_{\Xi^{\;\!\prime}}}\Omega_{_0}.$
Тогда}
\\[1.5ex]
\mbox{}\hfill
$
(u+i\;\!v,\;\! i\;\!V)\in \text{\rm H}_{_{\Xi^{\;\!\prime}}}
\iff
\Bigl(\exp\arctg\dfrac{v}{u}\,,\, V\Bigr)\in \text{\rm E}_{_{\Omega_{{}_{\tiny\;\! 0}}}}
\ \&\ \,
u^2+v^2 \in \text{\rm I}_{_{\Xi^{\;\!\prime}}}.
\hfill
$
\\[1.5ex]
\indent
{\sl Следует}
из теоремы 6.2 при $U(t,x)=0\;\;\forall (t,x)\in \Xi^{\;\!\prime}$
с учетом теоремы 0.2. $\k$
\vspace{1.25ex}

{\bf Свойство 6.15.}
\vspace{0.5ex}
{\it
Пусть числа $\gamma_1^{}, \gamma_2^{}\in\R\backslash\{0\},$
функции
$u,v\in \text{\rm P}_{_{\!\Xi^{\;\!\prime}}}$ взаимно простые, 
$(u+i\;\!v,\;\! U+i\;\!V)\in \text{\rm H}_{_{\Xi^{\;\!\prime}}},$
\vspace{0.5ex}
множество $\Omega_{_0}\!\subset\! \Xi^{\;\!\prime}\!$ такое, что  
$u(t,x)\ne 0\;\;\forall (t,x)\!\in\! \Omega_{_0},\ 
u(t,x)=0$ $\forall (t,x)\in {\sf C}_{_{\Xi^{\;\!\prime}}}\Omega_{_0}.$
Тогда
$\Bigl((u^2+v^2)^{{}^{\scriptsize \gamma_1^{}}}
\exp\Bigl(\gamma_2^{}\;\!\arctg\dfrac{v}{u}\Bigr)\,,\, M\Bigr)\in 
\text{\rm J}_{_{\Omega_{{}_{\tiny\;\! 0}}}},$
если и только если}
\\[1.75ex]
\mbox{}\hfill
$
M(t,x)=2\gamma_1^{}\;\!U(t,x)+\gamma_2^{}\;\!V(t,x)
\quad
\forall (t,x)\in \Xi^{\;\!\prime}.
\hfill
$
\\[1.5ex]
\indent
{\sl Следует}
из теоремы 6.2 и свойства 1.9. $\k$
\\[3.75ex]
\centerline{
{\bf  7. Кратные комплекснозначные полиномиальные частные интегралы}
}
\\[1.5ex]
\indent
{\bf Определение 7.1.}
\vspace{0.25ex}
{\it
Комплекснозначный полиномиальный частный интеграл $w$ с сомножителем $W$ 
на области $\Xi^{\;\!\prime}$ сис\-темы {\rm (0.1)} назовем
\textit{\textbf{кратным}},
\vspace{0.15ex}
если существуют такие натуральное число $h$ и функция
$z\in \text{\rm Z}_{_{\Xi^{\;\!\prime}}},$
\vspace{0.35ex}
являющаяся взаимно простой с функцией $w,$
что производная в силу системы {\rm (0.1)} 
\\[1.5ex]
\mbox{}\hfill                             % (7.1)
$
\displaystyle
{\frak d}\,\dfrac{z(t,x)}{w^{\;\!h}(t,x)}=Q(t,x)
\quad
\forall (t,x)\in \Omega_{_0},
$
\hfill {\rm (7.1)}
\\[1.75ex]
где функция $Q\in \text{\rm Z}_{_{\Xi^{\;\!\prime}}}$ 
и имеет степень $\deg_{\;\!x}^{} Q\leq d-1,$
\vspace{0.75ex}
множество
$\Omega_{_0}\subset\Xi^{\;\!\prime}$ та\-кое, что  
$
|w(t,x)|\ne 0\;\;\forall (t,x)\in\Omega_{_0},\
|w(t,x)|= 0\;\;\forall (t,x)\in {\sf C}_{_{\Xi^{\;\!\prime}}}\Omega_{_0}\;\!.$
}
\vspace{0.75ex}

Множество функций, являющихся кратными 
комплекснозначными полиномиальными частными интегралами 
на области $\Xi^{\;\!\prime}$ сис\-темы (0.1), обозначим 
$\text{\rm G}_{_{\Xi^{\;\!\prime}}}.$
\vspace{0.5ex}

В соответствии 
\vspace{0.15ex}
с определением 7.1
$\text{\rm G}_{_{\Xi^{\;\!\prime}}}\subset
\text{\rm H}_{_{\Xi^{\;\!\prime}}}.$
\vspace{1ex}

\newpage

Условной записью 
\vspace{0.35ex}
$\bigl((w, W), (h,z,Q)\bigr)\in \text{\rm G}_{_{\Xi^{\;\!\prime}}}$
будем выражать, что 
комплекснозначный полиномиальный частный интеграл $w$ с сомножителем $W$ 
\vspace{0.15ex}
на области $\Xi^{\;\!\prime}$ сис\-темы {\rm (0.1)} является кратным 
таким, что выполняется тождество (7.1).
\vspace{0.5ex}

{\bf Теорема 7.1}
\vspace{0.15ex}
(критерий существования кратного комплекснозначного
полиномиального частного интеграла).
{\it
$\bigl((w, W), (h,z,Q)\bigr)\in \text{\rm G}_{_{\Xi^{\;\!\prime}}}$
\vspace{0.35ex}
тогда и только тогда, когда выполняются тождества {\rm(6.1)} и {\rm(7.1)}, 
в которых число $h\in\N,$ функции
\vspace{0.5ex}
$z, Q\in \text{\rm Z}_{_{\Xi^{\;\!\prime}}},$
функции $w$ и $z$ взаимно простые, $\deg_{\;\!x}^{} Q\leq d-1,$
\vspace{0.75ex}
множество
$\Omega_{_0}\subset\Xi^{\;\!\prime}$ та\-кое, что  
$|w(t,x)|\ne 0\;\;\forall (t,x)\in\Omega_{_0},$
$|w(t,x)|= 0\;\;\forall (t,x)\in {\sf C}_{_{\Xi^{\;\!\prime}}}\Omega_{_0}\;\!.$
}
\vspace{0.5ex}

{\sl Следует} из определения 6.1 и 7.1 (с учетом замечаний 6.1 и 6.2).$\k$
\vspace{0.75ex}

{\bf Теорема 7.2}
\vspace{0.15ex}
(критерий существования кратного комплекснозначного
полиномиального частного интеграла).
{\it
Пусть $h\in\N,$ функции
\vspace{0.5ex}
$w, z\in \text{\rm Z}_{_{\Xi^{\;\!\prime}}}\!$ взаимно простые, 
множество
$\Omega_{_0}\subset\Xi^{\;\!\prime}$ та\-кое, что  
$|w(t,x)|\ne 0\;\;\forall (t,x)\in\Omega_{_0},\
|w(t,x)|= 0\;\;\forall (t,x)\in {\sf C}_{_{\Xi^{\;\!\prime}}}\Omega_{_0}.$
Тогда}
\\[1.75ex]
\mbox{}\hfill                           
$
\displaystyle
\bigl((w, W), (h,z,Q)\bigr)\in \text{\rm G}_{_{\Xi^{\;\!\prime}}}
\iff
(w, W)\in \text{\rm H}_{_{\Xi^{\;\!\prime}}}
\ \& \
\hfill                           
$
\\[2ex]
\mbox{}\hfill                           
$
\& \
\Bigl(\exp\;\!{\rm Re}\;\!\dfrac{z}{w^{\;\!h}}\,, {\rm Re}\;\!Q\Bigr)\in 
\text{\rm E}_{_{\Omega_{{}_{\tiny\;\! 0}}}}
\ \,\& \ \,
\Bigl(\exp\;\!{\rm Im}\;\!\dfrac{z}{w^{\;\!h}}\,, {\rm Im}\;\!Q\Bigr)\in 
\text{\rm E}_{_{\Omega_{{}_{\tiny\;\! 0}}}}.
\hfill
$
\\[1.75ex]
\indent
{\sl Доказательство.}
Разделяя в тождестве (7.1) вещественные и мнимые части и учитывая 
критерий существования экспоненциального частного интеграла (теорема 3.1), 
получаем, что тождество (7.1) выполняется тогда и только тогда, когда
\\[1.25ex]
\mbox{}\hfill
$
\Bigl(\exp\;\!{\rm Re}\;\!\dfrac{z}{w^{\;\!h}}\,, {\rm Re}\;\!Q\Bigr)\in 
\text{\rm E}_{_{\Omega_{{}_{\tiny\;\! 0}}}}
$
\ \ и \ \
$
\Bigl(\exp\;\!{\rm Im}\;\!\dfrac{z}{w^{\;\!h}}\,, {\rm Im}\;\!Q\Bigr)\in 
\text{\rm E}_{_{\Omega_{{}_{\tiny\;\! 0}}}}.
\hfill
$
\\[1.5ex]
\indent
Теперь утверждение теоремы 7.2 следует 
\vspace{0.75ex}
с учетом определения 6.1 из теоремы 7.1.$\k$

{\bf Теорема 7.3}
\vspace{0.15ex}
(критерий существования кратного комплекснозначного
полиномиального частного интеграла).
{\it
Пусть $h\in\N,$ функции
\vspace{0.5ex}
$u, v\in \text{\rm P}_{_{\!\Xi^{\;\!\prime}}}$ взаимно простые, 
функция $z\in \text{\rm Z}_{_{\Xi^{\;\!\prime}}}$ взаимно простая с 
функцией $u+i\;\!v,$  множество
\vspace{0.75ex}
$\Omega_{_0}\subset\Xi^{\;\!\prime}$ та\-кое, что  
$u(t,x)\ne 0\;\;\forall (t,x)\in\Omega_{_0},\
u(t,x)= 0\;\;\forall (t,x)\in {\sf C}_{_{\Xi^{\;\!\prime}}}\Omega_{_0}.$
Тогда}
\\[1.75ex]
\mbox{}\hfill                           
$
\displaystyle
\bigl((u+i\;\!v,\;\! U+i\;\!V), (h,z,Q)\bigr)\in \text{\rm G}_{_{\Xi^{\;\!\prime}}}
\iff
\Bigl(\bigl(u^2+v^2,\;\! 2\;\!U\bigr), 
\bigl(h,\;\!{\rm Re}\;\!\bigl(z(u-i\;\!v)^h\bigr),\;\! {\rm Re}\;\!Q\bigr)\Bigr)\in 
\text{\rm B}_{_{\Xi^{\;\!\prime}}}
\ \& \
\hfill                           
$
\\[1.75ex]
\mbox{}\hfill                           
$
\& \
\Bigl(\bigl(u^2+v^2,\;\! 2\;\!U\bigr), 
\bigl(h,\;\!{\rm Im}\;\!\bigl(z(u-i\;\!v)^h\bigr),\;\! {\rm Im}\;\!Q\bigr)\Bigr)\in 
\text{\rm B}_{_{\Xi^{\;\!\prime}}}
\ \,\& \ \,
\Bigl(\exp\arctg\dfrac{v}{u}\,,\;\! V\Bigr)\in 
\text{\rm E}_{_{\Omega_{{}_{\tiny\;\! 0}}}}.
\hfill
$
\\[1.75ex]
\indent
{\sl Доказательство}
основано на теореме 7.2 и следующих обстоятельствах.

По теореме 6.2,
\\[1.25ex]
\mbox{}\hfill                           
$
\displaystyle
(u+i\;\!v,\;\! U+i\;\!V)\in \text{\rm H}_{_{\Xi^{\;\!\prime}}}
\iff
\bigl(u^2+v^2,\;\! 2\;\!U\bigr)\in\text{\rm A}_{_{\Xi^{\;\!\prime}}}
\ \& \
\Bigl(\exp\arctg\dfrac{v}{u}\,,\;\! V\Bigr)\in 
\text{\rm E}_{_{\Omega_{{}_{\tiny\;\! 0}}}}.
\hfill                           
$
\\[1.5ex]
\indent
Учитывая, что 
\\[1.25ex]
\mbox{}\hfill                           
$
\displaystyle
\bigl(u^2+v^2,\;\! 2\;\!U\bigr)\in\text{\rm A}_{_{\Xi^{\;\!\prime}}},
\quad
\dfrac{z}{w^h}=
\dfrac{z}{(u+i\;\!v)^h}=
\dfrac{z\;\!(u-i\;\!v)^h}{(u^2+v^2)^h}\,,
\hfill                           
$
\\[1.75ex]
на основании теоремы 5.3 получаем:
\\[1.75ex]
\mbox{}\hfill                           
$
\displaystyle
\Bigl(\exp\;\!{\rm Re}\;\!\dfrac{z}{w^h}\,,\, {\rm Re}\;\!Q\Bigr)\in 
\text{\rm E}_{_{\Omega_{{}_{\tiny\;\! 0}}}}
\iff
\biggl(\exp\;\!\dfrac{{\rm Re}\;\!\bigl(z\;\!(u-i\;\!v)^h\bigr)}{(u^2+v^2)^h}\,,\, {\rm Re}\;\!Q\biggr)\in 
\text{\rm E}_{_{\Omega_{{}_{\tiny\;\! 0}}}}
\iff
\hfill                           
$
\\[2ex]
\mbox{}\hfill                           
$
\iff
\Bigl(\bigl(u^2+v^2,\;\! 2\;\!U\bigr), 
\bigl(h,\;\!{\rm Re}\;\!\bigl(z(u-i\;\!v)^h\bigr),\;\! {\rm Re}\;\!Q\bigr)\Bigr)\in 
\text{\rm B}_{_{\Xi^{\;\!\prime}}};
\hfill                           
$
\\[2ex]
\mbox{}\hfill                           
$
\displaystyle
\Bigl(\exp\;\!{\rm Im}\;\!\dfrac{z}{w^h}\,,\, {\rm Im}\;\!Q\Bigr)\in 
\text{\rm E}_{_{\Omega_{{}_{\tiny\;\! 0}}}}
\iff
\biggl(\exp\;\!\dfrac{{\rm Im}\;\!\bigl(z\;\!(u-i\;\!v)^h\bigr)}{(u^2+v^2)^h}\,,\, {\rm Im}\;\!Q\biggr)\in 
\text{\rm E}_{_{\Omega_{{}_{\tiny\;\! 0}}}}
\iff
\hfill                           
$
\\[2ex]
\mbox{}\hfill                           
$
\iff
\Bigl(\bigl(u^2+v^2,\;\! 2\;\!U\bigr), 
\bigl(h,\;\!{\rm Im}\;\!\bigl(z(u-i\;\!v)^h\bigr),\;\! {\rm Im}\;\!Q\bigr)\Bigr)\in 
\text{\rm B}_{_{\Xi^{\;\!\prime}}}.\ \k
\hfill                           
$
\\[2ex]
\indent
{\bf Следствие 7.1.}
\vspace{0.5ex}
{\it
Пусть $u, v, a, b\in \text{\rm P}_{_{\!\Xi^{\;\!\prime}}},$ 
функции $u,\ v$ взаимно простые, 
функции $u+i\;\!v$ и $a+i\;\!b$ взаимно простые, множество
\vspace{0.75ex}
$\Omega_{_0}\!\subset\!\Xi^{\;\!\prime}\!$ та\-кое, что  
$u(t,x)\!\ne\! 0\;\;\forall (t,x)\!\in\!\Omega_{_0},\!$
$u(t,x)= 0\;\;\forall (t,x)\in {\sf C}_{_{\Xi^{\;\!\prime}}}\Omega_{_0}.$
Тогда}
\\[1.75ex]
\mbox{}\hfill                           
$
\displaystyle
\bigl((u+i\;\!v,\;\! U+i\;\!V), (1, a+i\;\!b, K+i\;\!L)\bigr)\in \text{\rm G}_{_{\Xi^{\;\!\prime}}}
\iff
\Bigl(\bigl(u^2+v^2,\;\! 2\;\!U\bigr), 
\bigl(1,\;\!a\;\!u+b\;\!v,\;\! K\bigr)\Bigr)\in 
\text{\rm B}_{_{\Xi^{\;\!\prime}}}
\ \& \
\hfill                           
$
\\[2ex]
\mbox{}\hfill                           
$
\& \
\Bigl(\bigl(u^2+v^2,\;\! 2\;\!U\bigr), 
\bigl(1,\;\!b\;\!u-a\;\!v,\;\! L\bigr)\Bigr)\in 
\text{\rm B}_{_{\Xi^{\;\!\prime}}}
\ \,\& \ \,
\Bigl(\exp\arctg\dfrac{v}{u}\,,\;\! V\Bigr)\in 
\text{\rm E}_{_{\Omega_{{}_{\tiny\;\! 0}}}}.
\hfill
$
\\[1.75ex]
\indent
{\sl Следует}
\vspace{1ex}
из теоремы 7.3 с учетом того, что 
$(a+i\;\!b)(u-i\;\!v)=a\;\!u+b\;\!v+(b\;\!u-a\;\!v)\;\!i.\ \k$

{\bf Следствие 7.2.}
\vspace{0.5ex}
{\it
Пусть $u, v, a, b\in \text{\rm P}_{_{\!\Xi^{\;\!\prime}}},$ 
функции $u,\ v$ взаимно простые, 
функции $u+i\;\!v$ и $a+i\;\!b$ взаимно простые, множество
\vspace{0.75ex}
$\Omega_{_0}\!\subset\!\Xi^{\;\!\prime}\!$ та\-кое, что  
$u(t,x)\!\ne\! 0\;\;\forall (t,x)\!\in\!\Omega_{_0},\!$
$u(t,x)= 0\;\;\forall (t,x)\in {\sf C}_{_{\Xi^{\;\!\prime}}}\Omega_{_0}.$
Тогда}
\\[1.75ex]
\mbox{}\hfill                           
$
\displaystyle
\bigl((u+i\;\!v,\;\! U+i\;\!V), (2, a+i\;\!b, K+i\;\!L)\bigr)\in \text{\rm G}_{_{\Xi^{\;\!\prime}}}
\iff
\hfill                           
$
\\[1.75ex]
\mbox{}\hfill                           
$
\iff
\Bigl(\bigl(u^2+v^2,\;\! 2\;\!U\bigr), 
\bigl(2,\;\!a(u^2-v^2)+2\;\!b\;\!u\;\!v,\;\! K\bigr)\Bigr)\in 
\text{\rm B}_{_{\Xi^{\;\!\prime}}}
\ \& \
\hfill                           
$
\\[1.75ex]
\mbox{}\hfill                           
$
\& \
\Bigl(\bigl(u^2+v^2,\;\! 2\;\!U\bigr), 
\bigl(2,\;\!b\;\!(u^2-v^2)-2a\;\!u\;\!v,\;\! L\bigr)\Bigr)\in 
\text{\rm B}_{_{\Xi^{\;\!\prime}}}
\ \,\& \ \,
\Bigl(\exp\arctg\dfrac{v}{u}\,,\;\! V\Bigr)\in 
\text{\rm E}_{_{\Omega_{{}_{\tiny\;\! 0}}}}.
\hfill
$
\\[1.5ex]
\indent
{\sl Следует}
из теоремы 7.3 с учетом того, что 
\\[1.5ex]
\mbox{}\hfill
$
(a+i\;\!b)(u-i\;\!v)^2=a\;\!(u^2-v^2)+2\;\!b\;\!u\;\!v+
\bigl(b\;\!(u^2-v^2)-2\;\!a\;\!u\;\!v\bigr)\;\!i.\ \k
\hfill
$
\\[2ex]
\indent
{\bf Определение 7.2.}
{\it
Комплекснозначный 
полиномиальный частный интеграл $w$  на области $\Xi^{\;\!\prime}$ системы {\rm (0.1)} 
назовем \textit{\textbf{кратным с кратностью}}
$
\varkappa =1 + \sum\limits_{\xi=1}^{\varepsilon} \delta_{\xi}^{},
$
если  существуют  такие натуральные числа 
\vspace{0.75ex} 
$h_{\xi}^{},\ \xi =1,\ldots, \varepsilon,$ и соответствующие этим числам функции 
$
z_{_{\scriptstyle h_\xi^{}f_\xi^{}}}\in \text{\rm Z}_{_{\Xi^{\;\!\prime}}},\ 
f_{\xi}^{}=1,\ldots,\delta_{\xi}^{},\ \xi=1,\ldots,\varepsilon,$
\vspace{0.35ex} 
каждая из которых является взаимно простой с функцией $w,$
что выполняются тождества 
\\[1.75ex]
\mbox{}\hfill                  
$
\displaystyle
{\frak d}\, \dfrac{z_{_{\scriptstyle h_\xi^{}  f_\xi^{}} }(t,x)}{\displaystyle  w^{\;\!h_\xi^{}} (t,x)}=
Q_{h_\xi^{}  f_\xi^{}}^{} (t,x)
\quad 
\forall (t,x)\in \Omega_{_0},
\quad 
f_{\xi}^{}=1,\ldots, \delta_{\xi}^{}, \ \   
\xi=1,\ldots,\varepsilon,
\hfill
$
\\[2.25ex]
где функ\-ции 
\vspace{0.75ex}
$Q_{h_\xi^{}  f_\xi^{}}^{}\in\text{\rm Z}_{_{\Xi^{\;\!\prime}}}$ 
имеют степени
$
\deg_{\;\!x}^{}  Q_{h_\xi^{}  f_\xi^{}}^{}\leq d-1, \ 
f_{\xi}^{}=1,\ldots, \delta_{\xi}^{}, \,  \xi=1,\ldots,\varepsilon,
$
мно\-жес\-т\-во
\vspace{1.25ex}
$\Omega_{_0}\subset\Xi^{\;\!\prime}$ та\-кое, что  
$|w(t,x)|\ne 0\;\;\forall (t,x)\in\Omega_{_0},\
|w(t,x)|= 0\;\;\forall (t,x)\in {\sf C}_{_{\Xi^{\;\!\prime}}}\Omega_{_0}.$
}

На основании определения 7.2 получаем
\vspace{1ex}

{\bf Предложение 7.1.}\!
\vspace{0.5ex}
{\it
Если  
$\!\Bigl((w, W), 
\Bigl(h_\xi^{},z_{_{\scriptstyle h_\xi^{}  f_\xi^{}}}, Q_{h_\xi^{}  f_\xi^{}}^{}\Bigr)\Bigr)
\!\in\! \text{\rm G}_{_{\Xi^{\;\!\prime}}},\,
f_{\xi}^{}\!=\!1,\ldots, \delta_{\xi}^{},\, \xi\!=\!1,\ldots,\varepsilon,\!$
то  комплекснозначный полиномиальный частный интеграл $w$ с сомножителем $W$
на области $\Xi^{\;\!\prime}$ системы {\rm (0.1)} является кратным с 
кратностью 
$
\varkappa =1 + \sum\limits_{\xi=1}^{\varepsilon} \delta_{\xi}^{}.
$
}
\vspace{0.75ex}

{\bf Свойство 7.1.}
{\it
Пусть 
$k\in\N,\ c\in\C\backslash\{0\},\ 
\bigl(w, w^kW_{0}^{}\bigr)\in \text{\rm H}_{_{\Xi^{\;\!\prime}}}.
$ 
Тогда}
\\[1.5ex]
\mbox{}\hfill                           
$
\bigl(\bigl(w, w^kW_{0}^{}\bigr), \bigl(k, c,{}-k\;\!c\;\!W_0^{}\bigr)\bigr)
\in \text{\rm G}_{_{\Xi^{\;\!\prime}}}.
\hfill
$
\\[1.5ex]
\indent
{\sl Доказательство}. 
В соответствии с определением 6.1
\\[1.5ex]
\mbox{}\hfill                           
$
{\frak d}\;\!w(t,x)=w^{k+1}(t,x)\;\!W_0^{}(t,x)
\quad
\forall (t,x)\in \Xi^{\;\!\prime}.
\hfill
$
\\[1.5ex]
\indent
Тогда производная в силу системы (0.1)
\\[1.5ex]
\mbox{}\hfill                           
$
{\frak d}\;\!\dfrac{c}{w^{k}(t,x)}=
{}-k\;\!c\,
\dfrac{{\frak d}\;\!w(t,x)}{w^{k+1}(t,x)}=
{}-k\;\!c\;\!W_0^{}(t,x)
\quad
\forall (t,x)\in \Omega_{{}_0}^{},
\hfill
$
\\[1.75ex]
где мно\-жес\-т\-во
\vspace{0.75ex}
$\Omega_{_0}\!\subset\!\Xi^{\;\!\prime}\!$ та\-кое, что  
$|w(t,x)|\ne 0\;\;\forall (t,x)\!\in\!\Omega_{_0},\
|w(t,x)|= 0\;\;\forall (t,x)\!\in\! {\sf C}_{_{\Xi^{\;\!\prime}}}\Omega_{_0}.\!$

По определению 7.1, 
$
\bigl(\bigl(w, w^kW_{0}^{}\bigr), \bigl(k, c,{}-k\;\!c\;\!W_0^{}\bigr)\bigr)
\in \text{\rm G}_{_{\Xi^{\;\!\prime}}}.\ \k
$
\vspace{1.25ex}

{\bf Свойство 7.2.}
\vspace{0.5ex}
{\it
Пусть  
$k\in\N,\ 
c_l^{}\in\C\backslash\{0\},\ l=1,\ldots,k,\ 
\bigl(w, w^kW_{0}^{}\bigr)\in \text{\rm H}_{_{\Xi^{\;\!\prime}}}.$ 
Тогда комплекснозначный полиномиальный частный интеграл $w$ 
\vspace{0.25ex}
на области $\Xi^{\;\!\prime}$ системы {\rm (0.1)} 
является $(k+1)\!$-кратным таким, что}
\\[1.75ex]
\mbox{}\hfill                           
$
\bigl(\bigl(w, w^kW_{0}^{}\bigr), \bigl(l, c_l^{},{}-l\;\!c_l^{}\;\!w^{\;\!k-l}W_0^{}\bigr)\bigr)
\in \text{\rm G}_{_{\Xi^{\;\!\prime}}},
\ \ 
l=1,\ldots, k.
\hfill
$
\\[2ex]
\indent
{\sl Доказательство.}
Если  
$\bigl(w, w^kW_{0}^{}\bigr)\in \text{\rm H}_{_{\Xi^{\;\!\prime}}},$ то
$\bigl(w, w^{\;\!l}W_{l}^{}\bigr)\in \text{\rm H}_{_{\Xi^{\;\!\prime}}},$
где функции
\\[2ex]
\mbox{}\hfill                           
$
W_{l}^{}(t,x)=w^{\;\!k-l}(t,x)\;\!W_0^{}(t,x)
\quad
\forall (t,x)\in \Xi^{\;\!\prime},
\quad 
l=1,\ldots, k.
\hfill
$
\\[2ex]
\indent
По свойству 7.1, 
\vspace{1.25ex}
$
\bigl(\bigl(w, w^{\;\!l}W_{l}^{}\bigr), \bigl(l, c_l^{},{}-l\;\!c_l^{}\;\!W_l^{}\bigr)\bigr)
\in \text{\rm G}_{_{\Xi^{\;\!\prime}}},
\ c_l^{}\in\C\backslash\{0\},
\ l=1,\ldots, k.
$
Отсюда, 
$
\bigl(\bigl(w, w^{\;\!k}W_{0}^{}\bigr), \bigl(l, c_l^{},{}-l\;\!c_l^{}\;\!w^{\;\!k-l}W_0^{}\bigr)\bigr)
\in \text{\rm G}_{_{\Xi^{\;\!\prime}}},\
c_l^{}\in\C\backslash\{0\},\ l=1,\ldots,k.\ \k$ 
\\[5.5ex]
\centerline{
{\bf\large \S\;\!2. Последние множители}}
\\[2ex]
\centerline{
{\bf  8. 
Последний множитель как частный интеграл}
}
\\[1.5ex]
\indent
{\bf Теорема 8.1}
(критерий последнего множителя).
{\it 
Непрерывно дифференцируемая функция является 
последним множителем на области $\Omega$ системы {\rm (0.1)} 
тогда и только тогда, когда она является частным интегралом 
с сомножителем $-\,{\rm div}\;\!{\frak d}$ 
на области $\Omega$ системы {\rm (0.1)}.
}
\vspace{0.35ex}

{\sl Доказательство}.
\vspace{0.25ex}
Так как функции
$X_i^{}\in \text{\rm P}_{_{\!\Xi}},\ i=1,\ldots,n,$ то  
${\rm div}\;\!{\frak d}\in \text{\rm P}_{_{\!\Xi}}$
и имеет степень
$\deg_{\;\!x}^{}\;\! {\rm div}\;\!{\frak d}\leq d-1.$ 
\vspace{0.5ex}

При ${\rm g}=\mu, \ M=-\,{\rm div}\;\!{\frak d}$
\vspace{0.15ex}
тождество (1.2) из критерия существования частного интеграла (теорема 1.1)
совпадает с тождеством (0.6) из критерия существования 
последнего множителя (теорема 0.5).

Следовательно,
\\[0.5ex]
\mbox{}\hfill
$
\mu\in \text{\rm M}_{_{\Omega}}
\iff
(\mu, -\,{\rm div}\;\!{\frak d})\in \text{\rm J}_{_{\Omega}}.
\ \k
\hfill
$
\\[1.75ex]
\indent
В соответствии с теоремой 8.1
$
\text{\rm M}_{_{\Omega}}\subset\text{\rm J}_{_{\Omega}}.
$
\vspace{0.5ex}

Если последний множитель системы (0.1)
является частным интегралом: 
а) полино\-миальным;
б) кратным полиномиальным;
в) экспоненциальным;
г) условным;
д) ком\-п\-лек\-с\-но\-значным полиномиальным;
е) кратным комплекснозначным полиномиальным,
то его соответственно назовем:
а) полиномиальным;
б) кратным полиномиальным;
\linebreak
в) экспоненциальным;
г) условным;
д) комплекснозначным полиномиальным;
е) кратным комплекснозначным полиномиальным.
\vspace{0.15ex}

Множество последних множителей на области $\Omega$ 
\vspace{0.15ex}
или $\Xi^{\;\!\prime}$ системы (0.1) обозначим:
а) $\text{\rm MA}_{_{\Xi^{\;\!\prime}}}$ --- 
полиномиальных;
\vspace{0.5ex}
\ б) $\text{\rm MB}_{_{\Xi^{\;\!\prime}}}$ --- 
кратных полиномиальных;
\ в) $\text{\rm ME}_{_{\Omega}}$ --- 
экспоненциальных;
\ г) $\text{\rm MF}_{_{\!\Xi^{\;\!\prime}}}$ --- 
условных;
\vspace{0.35ex}
\ д) $\text{\rm MH}_{_{\Xi^{\;\!\prime}}}$ --- 
комплекснозначных полиномиальных;
\ е) $\text{\rm MG}_{_{\Xi^{\;\!\prime}}}$ --- 
кратных комплекснозначных полиномиальных.
\vspace{0.75ex}

{\bf Теорема 8.2}
(геометрический смысл последнего множителя).
{\it 
Если последний множитель $\mu$ системы {\rm (0.1)} 
определяет многообразие $\mu(t,x)=0,$ 
то оно будет интегральным многообразием системы {\rm (0.1)}.
}
\vspace{0.25ex}

{\sl Следует}
из определения интегрального многообразия (определение 0.3) 
и определения последнего множителя (определение 0.4). $\k$
\vspace{0.5ex}

{\bf Теорема 8.3.}
{\it
Непрерывно дифференцируемая функция {\rm (1.4)}
являет\-ся последним множителем на области $\Omega$ системы {\rm (0.1)},
\vspace{0.25ex}
если и только если существуют такие функции 
$M_j^{}\in C^1\Omega,\ j=1,\ldots,m,$ 
что выполняются тождества {\rm (1.5)} и}
\\[1ex]
\mbox{}\hfill
$
\displaystyle
\sum\limits_{j=1}^{m} M_j^{}(t,x)= -\,{\rm div}\;\!{\frak d}(t,x)
\quad
\forall (t,x)\in\Omega.
\hfill
$
\\[1.5ex]
\indent
{\sl Следует}
из теоремы 1.3 с учетом теоремы 8.1. $\k$
\vspace{0.75ex}

{\bf  Свойство 8.1.}
\vspace{0.25ex}
{\it
Пусть функция $\varphi\in C^1T^{\;\!\prime},$ 
множество $T_{_0}\subset T^{\;\!\prime}$ 
такое, что 
$\varphi(t)\ne 0$ $\forall t\in T_{_0},\ 
\varphi(t)=0\;\;\forall t\in {\sf C}_{{}_{T^{\;\!\prime}}}T_{_0},$
множество $\Omega_{_0}=T_{_0}\times X^{\;\!\prime}.$ Тогда}
\\[2ex]
\mbox{}\hfill
$
\varphi\;\! {\rm g}\in \text{\rm M}_{_{\Omega}}
\iff 
\bigl({\rm g}\;\!,\;\! {}-{\sf D} \ln |\varphi| -\,{\rm div}\;\!{\frak d}\bigr)\in 
\text{J}_{_{\Omega_{{}_{\tiny\;\! 0}}}}.
\hfill
$
\\[2ex]
\indent
{\sl Следует} с учетом теоремы 8.1 из свойства 1.1 при   
\\[1.5ex]
\mbox{}\hfill                        
$
M(t,x)+ {\sf D} \ln |\varphi(t) | =
-\,{\rm div}\;\!{\frak d}(t,x)
\quad
\forall (t,x)\in  \Omega_{_0}.\ \k
\hfill 
$
\\[2ex]
\indent
Частным случаем свойства 8.1 является 
\vspace{0.5ex}

{\bf Свойство 8.2.} 
{\it
Если $\lambda\in\R\backslash\{0\},$ то}
\\[1.5ex]
\mbox{}\hfill
$
\mu\in \text{\rm M}_{_{\Omega}}
\iff 
\lambda\;\!\mu\in \text{\rm M}_{_{\Omega}}.
\hfill
$
\\[1.75ex]
\indent
В соответствии со свойством 8.2, говоря о двух и более последних множителях  
системы (0.1), будем считать их попарно линейно независимыми.
\vspace{0.5ex}

{\bf Свойство 8.3.} 
\vspace{0.25ex}
{\it
Пусть множество $\Omega_{_0}\subset\Omega$ такое, что  
$\mu(t,x)\ne 0\;\;\forall (t,x)\in \Omega_{_0},$ 
$\mu(t,x)=0\;\;\forall (t,x)\in {\sf C}_{_\Omega}\Omega_{_0}.$
Тогда}
\\[1.5ex]
\mbox{}\hfill  
$
\mu\in \text{\rm M}_{_{\Omega}}
\iff 
|\mu|\in \text{\rm M}_{_{\Omega_{{}_{\tiny\;\! 0}}}}.
\hfill
$
\\[1.75ex]
\indent
{\sl Следует}
из свойства 1.3 с учетом теоремы 8.1. $\k$
\vspace{0.75ex}

{\bf Свойство 8.4.} 
\vspace{0.5ex}
{\it
Пусть 
$\lambda_j^{}\in\R\backslash\{0\},\ j=1,\ldots,m.$ Тогда} 
\\[1.25ex]
\mbox{}\hfill  
$
\displaystyle
\mu_j^{}\in \text{\rm M}_{_{\Omega}}, \
j=1,\ldots, m,
\ \ \Longrightarrow\ \
\sum\limits_{j=1}^{m} \lambda_j^{}\;\!\mu_j^{}\in \text{\rm M}_{_{\Omega}}\;\!.
\hfill
$
\\[1.75ex]
\indent
{\sl Следует}
из свойства 1.5 с учетом теоремы 8.1. $\k$
\vspace{1ex}

{\bf Свойство 8.5.} 
{\it
Пусть $\gamma\in\R\backslash\{0\},\ f^{\;\!\gamma}\in C^1\Omega.$
Тогда
\\[1.5ex]
\mbox{}\hfill  
$
f^{\;\!\gamma}\in \text{\rm M}_{_{\Omega}}
\iff 
\Bigl( f,\;\! {}-\dfrac{1}{\gamma}\ {\rm div}\;\!{\frak d}\Bigr)\in \text{\rm J}_{_{\Omega}}\;\!,
\hfill 
$
\\[1ex]
а}
\\[1ex]
\mbox{}\hfill  
$
f\in \text{\rm M}_{_{\Omega}}
\iff 
\bigl( f^{\;\!\gamma},\;\! {}-\gamma\, {\rm div}\;\!{\frak d}\bigr)\in \text{\rm J}_{_{\Omega}}\;\!.
\hfill 
$
\\[1.75ex]
\indent
{\sl Следует}
из свойства 1.6 с учетом теоремы 8.1. $\k$
\vspace{1ex}

{\bf Свойство 8.6.} 
\vspace{0.5ex}
{\it
Пусть $\gamma\in\R\backslash\{0\},$ 
множество $\Omega_{_0}\subset\Omega$ такое, что  
$f(t,x)\ne 0$ $\forall (t,x)\in \Omega_{_0},\ 
f(t,x)=0\;\;\forall (t,x)\in {\sf C}_{_\Omega}\Omega_{_0}.$
Тогда
\\[1.5ex]
\mbox{}\hfill  
$
|f|^{{}^{\scriptstyle \gamma}}\in \text{\rm M}_{_{\Omega_{{}_{\tiny\;\! 0}}}}
\iff 
\Bigl( f,\;\! {}-\dfrac{1}{\gamma}\ {\rm div}\;\!{\frak d}\Bigr)\in \text{\rm J}_{_{\Omega}}\;\!,
\hfill 
$
\\[1ex]
а}
\\[1ex]
\mbox{}\hfill  
$
f\in \text{\rm M}_{_{\Omega}}
\iff 
\bigl(\;\! |f|^{{}^{\scriptstyle \gamma}},\;\! {}-\gamma\, {\rm div}\;\!{\frak d}\bigr)\in 
\text{\rm J}_{_{\Omega_{{}_{\tiny\;\! 0}}}}.
\hfill 
$
\\[1.75ex]
\indent
{\sl Следует}
из свойства 1.7 с учетом теоремы 8.1. $\k$
\vspace{1ex}

{\bf Свойство 8.7.} 
\vspace{0.5ex}
{\it
Пусть $\rho_j^{},\;\! \lambda_j^{}\in\R\backslash\{0\},\ 
{\rm g}_j^{{}^{\scriptstyle 1/\rho_{\!j}^{}}}\in C^1\Omega, ,\ j=1,\ldots,m.$ Тогда}
\\[1.25ex]
\mbox{}\hfill  
$
\displaystyle
\bigl({\rm g}_j^{}\;\!,\;\! \rho_j^{}\, {\rm div}\;\!{\frak d}\bigr)\in 
\text{\rm J}_{_{\Omega}}, \
j=1,\ldots, m,
\ \ \Longrightarrow\ \
\sum\limits_{j=1}^{m} \lambda_j^{}\;\!
{\rm g}_j^{{}^{\scriptstyle {}-1/\rho_{\!j}^{}}}
\in \text{\rm M}_{_{\Omega}}\;\!.
\hfill
$
\\[1.75ex]
\indent
{\sl Следует}
из свойства 1.8 с учетом теоремы 8.1. $\k$
\vspace{1ex}

{\bf Свойство 8.8.} 
\vspace{0.5ex}
{\it
Пусть 
$({\rm g}_j^{}\;\!,\;\! M_{j}^{})\in \text{\rm J}_{_{\Omega}},\ 
\gamma_j^{}\in\R\backslash\{0\},\
{\rm g}_j^{{}^{\scriptsize \gamma_{j}^{}}}\in C^1\Omega, \ j=1,\ldots,m.$ 
Тогда  
$\prod\limits_{j=1}^{m} 
{\rm g}_j^{{}^{\scriptsize \gamma_{j}^{}}}
\in \text{\rm M}_{_{\Omega}}\;\!,$
если и только если линейная комбинация сомножителей}
\\[1ex]
\mbox{}\hfill                      % (8.1)
$
\displaystyle
\sum\limits_{j=1}^{m}
\gamma_j^{}\;\!M_{j}^{}(t,x)=
{}-\,{\rm div}\;\!{\frak d}(t,x)
\quad
\forall (t,x)\in\Xi^{\;\!\prime}.
$
\hfill (8.1)
\\[1.5ex]
\indent
{\sl Следует}
с учетом теоремы 8.1 из свойства 1.9 при $M=-\,{\rm div}\;\!{\frak d}.\ \k$ 
\vspace{1.25ex}

{\bf Свойство 8.9.} 
\vspace{0.75ex}
{\it
Пусть 
$({\rm g}_j^{}\;\!,\;\!\rho_{j}^{}\,{\rm div}\;\!{\frak d})\in \text{\rm J}_{_{\Omega}},\ 
\rho_{j}^{},\gamma_j^{}\in\R\backslash\{0\},\
{\rm g}_j^{{}^{\scriptsize \gamma_{j}^{}}}\in C^1\Omega, \ j=1,\ldots,m.$ 
Тогда  
$\prod\limits_{j=1}^{m} 
{\rm g}_j^{{}^{\scriptsize \gamma_{j}^{}}}
\in \text{\rm M}_{_{\Omega}}\;\!,$
если и только если}
$
\sum\limits_{j=1}^{m}
\rho_{j}^{}\;\!\gamma_j^{}={}-1.
$
\vspace{0.75ex}

{\sl Следует}
из свойства 8.8 при $M_j^{}=\rho_{j}^{}\,{\rm div}\;\!{\frak d},\ j=1,\ldots,m.\ \k$ 
\vspace{1.25ex}

{\bf Свойство 8.10.} 
\vspace{0.5ex}
{\it
Пусть 
$\mu_j^{}\in \text{\rm M}_{_{\Omega}}\;\!,\
\mu_j^{{}^{\scriptsize \gamma_{j}^{}}}\in C^1\Omega, \ 
\gamma_j^{}\in\R\backslash\{0\},\ j=1,\ldots,m.$ 
Тог\-да  
$\biggl(\,\prod\limits_{j=1}^{m} 
\mu_j^{{}^{\scriptsize \gamma_{j}^{}}}, M\biggl) 
\in \text{\rm J}_{_{\Omega}}\;\!,$
если и только если сомножитель}
\\[1ex]
\mbox{}\hfill                      % (8.2)
$
\displaystyle
M(t,x)={}-\sum\limits_{j=1}^{m}
\gamma_j^{}\, {\rm div}\;\!{\frak d}(t,x)
\quad
\forall (t,x)\in\Xi^{\;\!\prime}.
$
\hfill (8.2)
\\[1.5ex]
\indent
{\sl Следует}
с учетом теоремы 8.1 
из свойства 1.9 при $M_j^{}={}-\,{\rm div}\;\!{\frak d},\ j=1,\ldots,m.\ \k$ 
\vspace{1.25ex}

{\bf Свойство 8.11.} 
\vspace{0.5ex}
{\it
Пусть 
$\mu_j^{}\in \text{\rm M}_{_{\Omega}}\;\!,\
\mu_j^{{}^{\scriptsize \gamma_{j}^{}}}\in C^1\Omega, \ 
\gamma_j^{}\in\R\backslash\{0\},\ j=1,\ldots,m.$ 
Тог\-да  
$\prod\limits_{j=1}^{m} 
\mu_j^{{}^{\scriptsize \gamma_{j}^{}}}
\in \text{\rm M}_{_{\Omega}}\;\!,$
если и только если}
$\sum\limits_{j=1}^{m}\gamma_j^{}=1.$
\vspace{0.5ex}

{\sl Следует}
с учетом теоремы 8.1 
из свойства 8.10 при $M={}-\,{\rm div}\;\!{\frak d}.\ \k$ 
\vspace{1.25ex}

{\bf Свойство 8.12.} 
\vspace{0.5ex}
{\it
Пусть 
$\mu_j^{}\in \text{\rm M}_{_{\Omega}}\;\!,\
{\rm g}^{{}^{\scriptsize \gamma}},\mu_j^{{}^{\scriptsize \gamma_{j}^{}}}\in C^1\Omega, \ 
\gamma,\gamma_j^{}\in\R\backslash\{0\},\ j=1,\ldots,m.$ 
Тог\-да}  
\\[1ex]
\mbox{}\hfill                      
$
\displaystyle
\biggl(\,{\rm g}^{{}^{\scriptsize \gamma}}\prod\limits_{j=1}^{m} 
\mu_j^{{}^{\scriptsize \gamma_{j}^{}}}, M\biggl) 
\in \text{\rm J}_{_{\Omega}}
\iff
\biggl({\rm g}, \,
\dfrac{1}{\gamma}\;\!
\biggl(M+\sum\limits_{j=1}^{m}
\gamma_{j}^{}\;\!{\rm div}\;\!{\frak d}\biggl)\biggl) 
\in \text{\rm J}_{_{\Omega}}\;\!.
\hfill
$
\\[1.5ex]
\indent
{\sl Следует}
из свойства 1.10 с учетом теоремы 8.1. $\k$
\vspace{1.25ex}

{\bf Свойство 8.13.} 
\vspace{0.5ex}
{\it
Пусть 
$\mu_j^{}\in \text{\rm M}_{_{\Omega}}\;\!,\
{\rm g}^{{}^{\scriptsize \gamma}},\mu_j^{{}^{\scriptsize \gamma_{j}^{}}}\in C^1\Omega, \ 
\gamma,\gamma_j^{}\in\R\backslash\{0\},\ j=1,\ldots,m.$ 
Тог\-да}  
\\[1ex]
\mbox{}\hfill                      
$
\displaystyle
{\rm g}^{{}^{\scriptsize \gamma}}\prod\limits_{j=1}^{m} 
\mu_j^{{}^{\scriptsize \gamma_{j}^{}}} 
\in \text{\rm M}_{_{\Omega}}
\iff
\biggl({\rm g}, \,
\dfrac{1}{\gamma}\;\!
\biggl(\,\sum\limits_{j=1}^{m}
\gamma_{j}^{}-1\biggr)\, {\rm div}\;\!{\frak d}\biggl) 
\in \text{\rm J}_{_{\Omega}}\;\!.
\hfill
$
\\[1.5ex]
\indent
{\sl Следует}
с учетом теоремы 8.1
из свойства 8.12 при $M={}-\,{\rm div}\;\!{\frak d}.\ \k$ 
\vspace{1.25ex}

{\bf Свойство 8.14.} 
\vspace{0.5ex}
{\it
Пусть 
$\mu_j^{}\in \text{\rm M}_{_{\Omega}}\;\!,\
\mu_j^{{}^{\scriptsize \gamma_{j}^{}}}\in C^1\Omega, \ 
\gamma_j^{}\in\R\backslash\{0\},\ j=1,\ldots,m,\ 
{\rm g}_\tau^{{}^{\scriptsize \xi_{\tau}^{}}}\!\in C^1\Omega,$
$\xi_\tau^{}\in\R\backslash\{0\},\ \tau=1,\ldots,l.$ 
Тог\-да  
\\[1ex]
\mbox{}\hfill                      
$
\displaystyle
\biggl(\,\prod\limits_{\tau=1}^{l} 
{\rm g}_\tau^{{}^{\scriptsize \xi_{\tau}^{}}}
\,\prod\limits_{j=1}^{m} 
\mu_j^{{}^{\scriptsize \gamma_{j}^{}}}, M\biggl)\, 
\in \text{\rm J}_{_{\Omega}}
\iff
\biggl(\,\prod\limits_{\tau=1}^{l} 
{\rm g}_\tau^{{}^{\scriptsize \xi_{\tau}^{}}}, \,
M+\sum\limits_{j=1}^{m}
\gamma_{j}^{}\, {\rm div}\;\!{\frak d}\biggl) \,
\in \text{\rm J}_{_{\Omega}}\;\!.
\hfill
$
\\[1.5ex]
Кроме этого существуют такие функции 
$M_\tau^{}\in C^1\Omega,\ \tau=1,\ldots, l,$
что выполняются тождества
\\[1.5ex]
\mbox{}\hfill                        % (8.3)
$
{\frak d}\;\!{\rm g}_\tau^{}(t,x)=
{\rm g}_\tau^{}(t,x)\;\!M_\tau^{}(t,x)
\quad
\forall (t,x)\in\Omega,
\quad
\tau=1,\ldots, l,
$
\hfill {\rm (8.3)}
\\[1ex]
и}
\\[1ex]
\mbox{}\hfill                     
$
\displaystyle
\sum\limits_{\tau=1}^{l}
\xi_{\tau}^{}\;\!M_\tau^{}(t,x)=
M(t,x)+\sum\limits_{j=1}^{m}
\gamma_{j}^{}\, {\rm div}\;\!{\frak d}(t,x)
\quad
\forall (t,x)\in\Omega.
\hfill
$
\\[1.5ex]
\indent
{\sl Следует}
из свойства 1.11 с учетом теоремы 8.1. $\k$ 
\vspace{1.25ex}

\newpage

{\bf Свойство 8.15.} 
\vspace{0.5ex}
{\it
Пусть 
$\mu_j^{}\in \text{\rm M}_{_{\Omega}}\;\!,\
\mu_j^{{}^{\scriptsize \gamma_{j}^{}}}\in C^1\Omega, \ 
\gamma_j^{}\in\R\backslash\{0\},\ j=1,\ldots,m,\ 
{\rm g}_\tau^{{}^{\scriptsize \xi_{\tau}^{}}}\!\in C^1\Omega,$
$\xi_\tau^{}\in\R\backslash\{0\},\ \tau=1,\ldots,l.$ 
Тог\-да  
\\[1ex]
\mbox{}\hfill                      
$
\displaystyle
\prod\limits_{\tau=1}^{l} 
{\rm g}_\tau^{{}^{\scriptsize \xi_{\tau}^{}}}
\prod\limits_{j=1}^{m} 
\mu_j^{{}^{\scriptsize \gamma_{j}^{}}} 
\in \text{\rm M}_{_{\Omega}}
\iff
\biggl(\,\prod\limits_{\tau=1}^{l} 
{\rm g}_\tau^{{}^{\scriptsize \xi_{\tau}^{}}}, \,
\biggl(\,\sum\limits_{j=1}^{m}\gamma_{j}^{}-1\biggr)\, {\rm div}\;\!{\frak d}\biggl) \,
\in \text{\rm J}_{_{\Omega}}\;\!.
\hfill
$
\\[1.5ex]
Кроме этого существуют такие функции 
$M_\tau^{}\in C^1\Omega,\ \tau=1,\ldots, l,$
что выполняются тождества {\rm (8.3)} и}
\\[1ex]
\mbox{}\hfill                     
$
\displaystyle
\sum\limits_{\tau=1}^{l}
\xi_{\tau}^{}\;\!M_\tau^{}(t,x)=
\biggl(\,\sum\limits_{j=1}^{m}\gamma_{j}^{}-1\biggr)\, {\rm div}\;\!{\frak d}(t,x)
\quad
\forall (t,x)\in\Omega.
\hfill
$
\\[1.5ex]
\indent
{\sl Следует}
с учетом теоремы 8.1 из свойства 8.14 при $M={}-\,{\rm div}\;\!{\frak d}.\ \k$
\vspace{1ex}

{\bf Замечание 8.1.}
При необходимости целесообразно ввести понятие последнего псевдомножителя.

{\bf Определение 8.1.}
{\it 
Функцию $\nu\in C^1\Omega$ назовем 
\textit{\textbf{последним псевдомножителем с коэффициентом}} $\rho$
или \textit{\textbf{последним $\rho\!$-псевдомножителем $(\rho\in\R)$ 
на области}} $\Omega$ системы {\rm (0.1)},
если производная в силу системы} {\rm(0.1)} 
\\[2ex]
\mbox{}\hfill                        
$
{\frak d}\;\!\nu(t,x)=\rho\;\!\nu(t,x)\;{\rm div}\, {\frak d}(t,x)
\quad 
\forall (t,x)\in \Omega.
\hfill
$
\\[2.25ex]
\indent
Если использовать понятие частного интеграла (определение 1.1 и теорема 1.1), 
то получим критерий последнего псевдомножителя.
 \vspace{0.75ex}

{\bf Cвойство 8.16.}
 \vspace{0.35ex}
 {\it 
Функция $\nu\in C^1\Omega$ является последним $\rho\!$-псевдомножителем 
на области $\Omega$ системы {\rm (0.1)} 
тогда и только тогда, когда она является частным интегралом с 
сомножителем $\rho\;\!{\rm div}\, {\frak d}$ 
на области $\Omega$ системы {\rm (0.1)}. 
}
\vspace{0.75ex}

На основании свойства 8.5 с учетом определения 8.1 и теоремы 0.5 устанавливаем взаимосвязь между последним множителем и последним псевдомножителем.
 \vspace{0.75ex}

{\bf Cвойство 8.17А.}
 \vspace{0.35ex}
 {\it 
Функция $\nu\in C^1\Omega$ является последним псевдомножителем с коэффициентом $\rho\ne 0$ на области $\Omega$ системы {\rm (0.1)},
если и только если  функция
\\[1.5ex]
\mbox{}\hfill
$
\mu\colon (t,x)\to\ \nu^{\,{}-\;\!1/\rho}(t,x)
\quad
\forall (t,x)\in\Omega
\hfill
$
\\[1.5ex]
является последним множителем на области $\Omega$ системы {\rm (0.1)}. 
}
\vspace{0.75ex}

{\bf Cвойство 8.17Б.}
 \vspace{0.35ex}
 {\it 
Функция $\mu\in C^1\Omega$ является последним множителем 
на области $\Omega$ системы {\rm (0.1)}, если и только если
\vspace{0.35ex}
функция
$
\nu\colon (t,x)\to \mu^{\,{}-\;\!\rho}(t,x)
\;\;
\forall (t,x)\in\Omega
$
является последним псевдомножителем с коэффициентом $\rho\ne 0$ 
на области $\Omega$
системы {\rm (0.1)}. 
}
\vspace{0.75ex}

Таким образом, наряду с последними множителями можно рассматривать 
последние псевдомножители системы (0.1).

Обратим внимание на то, что последний $0\!$-псевдомножитель системы (0.1) есть ее первый интеграл.
\\[2.75ex]
\centerline{\bf  9. Построение последних множителей на основании} 
\centerline{\bf полиномиальных частных интегралов}
\\[1.5ex]
\indent
{\bf Свойство 9.1.} 
\vspace{0.5ex}
{\it
Пусть 
$p_j^{}\in \text{\rm P}_{_{\!\Xi^{\;\!\prime}}},\ 
\gamma_j^{}\in\R\backslash\{0\},$ 
множество $\Omega_{_0}\subset\Xi^{\;\!\prime}$ такое, что  
\linebreak
$p_j^{{}^{\scriptsize \gamma_j^{}}}\in C^1\Omega_{_0},\ j=1,\ldots, m.$
Тогда
\\[1.5ex]
\mbox{}\hfill  
$
\displaystyle
\prod\limits_{j=1}^{m} 
p_j^{{}^{\scriptsize \gamma_{j}^{}}}
\in \text{\rm M}_{_{\Omega_{{}_{\tiny\;\! 0}}}}
\iff
\bigl(p_j^{}, M_j^{}\bigr)\in \text{\rm A}_{_{\Xi^{\;\!\prime}}},
\ \ 
j=1,\ldots, m,
\hfill
$
\\[1.5ex]
где сомножители $M_j^{},\ j=1,\ldots,m,$
\vspace{0.5ex}
такие, что выполняется тождество {\rm (8.1)}.}

{\sl Следует} 
с учетом теоремы 8.1 из теоремы 2.2 при $M={}-\,{\rm div}\;\!{\frak d}.\ \k$ 
\vspace{1ex}

{\bf Свойство 9.2}
(критерий существования рационального последнего множителя). 
{\it
Пусть функции
$p_1^{},p_2^{}\in \text{\rm P}_{_{\!\Xi^{\;\!\prime}}}$
взаимно простые, 
\vspace{0.75ex}
множество $\Omega_{_0}\subset\Xi^{\;\!\prime}$ такое, что  
$p_2^{}(t,x)\ne 0$ $\forall (t,x)\in \Omega_{_0},\
p_2^{}(t,x)= 0\;\; \forall (t,x)\in {\sf C}_{_{\Xi^{\;\!\prime}}}\Omega_{_0}.$
Тогда
\\[1.5ex]
\mbox{}\hfill  
$
\dfrac{p_1^{}}{p_2^{}} 
\in \text{\rm M}_{_{\Omega_{{}_{\tiny\;\! 0}}}}
\iff
\bigl(p_1^{}, M_1^{}\bigr)\in \text{\rm A}_{_{\Xi^{\;\!\prime}}}
\ \& \ 
\bigl(p_2^{}, M_2^{}\bigr)\in \text{\rm A}_{_{\Xi^{\;\!\prime}}},
\hfill
$
\\[1.5ex]
где сомножители $M_1^{}$ и $M_2^{}$ такие, что}
\\[1.5ex]
\mbox{}\hfill  
$
M_2^{}(t,x)-M_1^{}(t,x)={\rm div}\;\!{\frak d}(t,x)
\quad
\forall (t,x)\in \Xi^{\;\!\prime}.
\hfill
$
\\[1.5ex]
\indent
{\sl Следует} 
с учетом теоремы 8.1 из теоремы 2.3 при $M={}-\,{\rm div}\;\!{\frak d}.\ \k$ 
\vspace{1ex}

{\bf  Свойство 9.3.}
\vspace{0.5ex}
{\it
Пусть $\lambda_j^{},\;\! c_j^{}\in\R, \ j=1,\ldots,m,\ m\leq n,\ 
\sum\limits_{j=1}^{m}|\lambda_j^{}|\ne 0.$
Тогда функция
$
%\displaystyle
\sum\limits_{j=1}^{m}\lambda_j^{}(x_j^{}+c_j^{})
\in\text{\rm MA}_{_{\Xi}},
$
если и только если выполняется тождество}
\\[1.5ex]
\mbox{}\hfill
$
\displaystyle
\sum\limits_{j=1}^{m}
\lambda_j^{}\;\!X_j^{}(t,x)=
{}-\sum\limits_{j=1}^{m}\lambda_j^{}(x_j^{}+c_j^{})\;\!{\rm div}\;\!{\frak d}(t,x)
\quad
\forall (t,x)\in \Xi.
\hfill
$
\\[1.5ex]
\indent
{\sl Следует} 
с учетом теоремы 8.1 из свойства 2.2 при $M={}-\,{\rm div}\;\!{\frak d}.\ \k$ 
\vspace{1ex}

{\bf  Свойство 9.4.}
\vspace{0.35ex}
{\it
Если $c\in\R,$ то
$
x_i^{}+c \in \text{\rm MA}_{_{\Xi}},\ 
i\in \{1,\ldots, n\},
$
тогда и только тогда, когда в правой части системы {\rm (0.1)} функция}
\\[1.5ex]
\mbox{}\hfill
$
X_i^{}(t,x)=
{}- (x_i^{}+c)\;\!{\rm div}\;\!{\frak d}(t,x)
\quad
\forall (t,x)\in \Xi.
\hfill
$
\\[1.5ex]
\indent
{\sl Следует} 
с учетом теоремы 8.1 из свойства 2.3 при $M={}-\,{\rm div}\;\!{\frak d}.\ \k$ 
\vspace{1ex}

%\newpage

{\bf Свойство 9.5}
\vspace{0.15ex}
(критерий существования кратного полиномиального последнего множителя).
{\it
Пусть $\nu\in \text{\rm MA}_{_{\!\Xi^{\;\!\prime}}},$ число $h\in\N,$ функции
\vspace{0.5ex}
$\nu, q\in \text{\rm P}_{_{\!\Xi^{\;\!\prime}}}$ взаимно простые, 
функция $N\in \text{\rm P}_{_{\!\Xi^{\;\!\prime}}}$ и имеет степень 
\vspace{0.5ex}
$\deg_{x}^{} N\leq d-1,$
множество
$\Omega_{_0}\subset\Xi^{\;\!\prime}$ та\-кое, что  
$\nu(t,x)\ne 0\;\;\forall (t,x)\in\Omega_{_0},\
\nu(t,x)= 0\;\;\forall (t,x)\in {\sf C}_{_{\Xi^{\;\!\prime}}}\Omega_{_0}.$
\vspace{0.5ex}
Тогда $\bigl(\nu, (h,q,N)\bigr)\in \text{\rm MB}_{_{\Xi^{\;\!\prime}}},$
если и только если выполняется одно из условий}:
\\[1ex]
\indent
$
1)\ {\frak d}\,\dfrac{q(t,x)}{\nu^{\;\!h}(t,x)}=N(t,x)
\quad
\forall (t,x)\in \Omega_{_0};
$
\\[1.25ex]
\indent
$
2)\ {\frak d}\;\!q(t,x)=\nu^{\;\!h}(t,x)\;\!N(t,x)-h\;\!q(t,x)\;\!{\rm div}\;\!{\frak d}(t,x)
\quad
\forall (t,x)\in \Xi^{\;\!\prime};
$
\\[1.5ex]
\indent
$
3)\ 
\Bigl(\exp\dfrac{q}{\nu^{\;\!h}}\,, N\Bigr)\in 
\text{\rm E}_{_{\Omega_{{}_{\tiny\;\! 0}}}}.
$
\\[1.5ex]
\indent
{\sl Следует} 
с учетом теоремы 8.1 из определения 5.1, теорем 5.2 и 5.3 при 
$p=\nu$ и $M={}-\,{\rm div}\;\!{\frak d}.\ \k$ 
\vspace{1ex}

{\bf Свойство\! 9.6.}\!
{\it
Пусть 
$\!\bigl(\nu_j^{}, (h_j^{},q_j^{},\rho_{j}^{}N_0^{})\bigr)\!\in\! 
\text{\rm MB}_{_{\Xi^{\;\!\prime}}},\, 
\rho_j^{}\!\in\!\R\backslash\{0\},\, 
\lambda_j^{},\gamma_j^{}, c_j^{}\!\in\!\R,\, 
\sum\limits_{j=1}^{m}\!|\lambda_j^{}|\!\ne\! 0,\!$ 
$\nu_j^{\gamma_j^{}}\in C^1\Omega_{_0}, \ j=1,\ldots, m,\
\varphi\in C^1T^{\;\!\prime},$ 
мно\-жес\-т\-во
$\Omega_{_0}\subset\Xi^{\;\!\prime}$ та\-кое, что  
$\prod\limits_{j=1}^{m}\nu_j^{}(t,x)\ne 0$
$\forall (t,x)\in\Omega_{_0},\
\prod\limits_{j=1}^{m}\nu_j^{}(t,x)= 0\;\;
\forall (t,x)\in {\sf C}_{_{\Xi^{\;\!\prime}}}\Omega_{_0}.$
Тогда{\rm:}
\\[1.5ex]
\indent
$
1)\
\displaystyle
\biggl(\ \prod\limits_{j=1}^{m}\nu_j^{\gamma_j^{}}
\sum\limits_{j=1}^{m}\lambda_j^{}
\exp\biggl(\;\!\dfrac{q_j^{}}{\rho_j^{}\;\!\nu_j^{\;\!h_j^{}}}+\varphi\biggr), L\biggr)\in 
\text{\rm J}_{_{\Omega_{{}_{\tiny\;\! 0}}}},
\hfill
$
\\[2ex]
если и только если сомножитель
\\[1.75ex]
\mbox{}\hfill                           
$
\displaystyle
L(t,x)={\sf D}\;\!\varphi(t)+N_0^{}(t,x)-
\sum\limits_{j=1}^{m}\gamma_j^{}\;\!{\rm div}\;\!{\frak d}(t,x)
\quad
\forall (t,x)\in \Xi^{\;\!\prime};
\hfill
$
\\[1.5ex]
\indent
$
2)\ 
\displaystyle
\biggl(\ \prod\limits_{j=1}^{m}\nu_j^{\gamma_j^{}}
\sum\limits_{j=1}^{m}\lambda_j^{}
\exp\biggl(\;\!\dfrac{q_j^{}}{\rho_j^{}\;\!\nu_j^{\;\!h_j^{}}}+c_j^{}\biggr), L\biggr)\in 
\text{\rm J}_{_{\Omega_{{}_{\tiny\;\! 0}}}},
\hfill
$
\\[2ex]
если и только если сомножитель}
\\[1.75ex]
\mbox{}\hfill                           
$
\displaystyle
L(t,x)=N_0^{}(t,x)-
\sum\limits_{j=1}^{m}\gamma_j^{}\;\!{\rm div}\;\!{\frak d}(t,x)
\quad
\forall (t,x)\in \Xi^{\;\!\prime}.
\hfill
$
\\[1.5ex]
\indent
{\sl Следует} 
с учетом теоремы 8.1 из свойств 5.5 и 5.6 при 
\\[1ex]
\mbox{}\hfill
$p_j^{}=\nu_j^{}\;\!,\ \ M_j^{}={}-\,{\rm div}\;\!{\frak d},\ \, j=1,\ldots,m.\ \k
\hfill
$ 
\\[1.75ex]
\indent
{\bf Свойство 9.7.}\!
\vspace{0.75ex}
{\it
Пусть 
$\!\bigl(\nu_j^{}, (h_j^{},q_j^{},N_j^{})\bigr)\!\in \text{\rm MB}_{_{\Xi^{\;\!\prime}}},\, 
\gamma_j^{},\xi_j^{}\!\in\!\R,\, \varphi_j^{}\!\in\! C^1T^{\;\!\prime},\, 
j\!=\!1,\ldots, m,\!$ 
мно\-жес\-т\-во
\vspace{0.75ex}
$\!\Omega_{_0}\!\subset\!\Xi^{\;\!\prime}\!$ та\-кое, что  
$\prod\limits_{j=1}^{m}\nu_j^{}(t,x)\!\ne\! 0\;\;\forall (t,x)\!\in\!\Omega_{_0},\
\prod\limits_{j=1}^{m}\nu_j^{}(t,x)= 0\;\;
\forall (t,x)\!\in\! {\sf C}_{_{\Xi^{\;\!\prime}}}\Omega_{_0},\!$ 
$\nu_j^{\gamma_j^{}}\in C^1\Omega_{_0}, \ j=1,\ldots, m.$
Тогда
\\[1.5ex]
\mbox{}\hfill                           
$
\displaystyle
\biggl(\ \prod\limits_{j=1}^{m}\nu_j^{\gamma_j^{}}
\exp\sum\limits_{j=1}^{m}\xi_j^{}
\biggl(\;\!\dfrac{q_j^{}}{\nu_j^{\;\! h_j^{}}}+\varphi_j^{}\biggr), L\biggr)\in 
\text{\rm J}_{_{\Omega_{{}_{\tiny\;\! 0}}}},
\hfill
$
\\[2ex]
если и только если сомножитель}
\\[1.75ex]
\mbox{}\hfill                           
$
\displaystyle
L(t,x)=
\sum\limits_{j=1}^{m}
\xi_j^{}\;\!\bigl(N_j^{}(t,x)+{\sf D}\;\!\varphi_j^{}(t)\bigr)- 
\sum\limits_{j=1}^{m}
\gamma_j^{}\;\!{\rm div}\;\!{\frak d}(t,x)
\quad
\forall (t,x)\in \Xi^{\;\!\prime}.
\hfill
$
\\[1.5ex]
\indent
{\sl Следует} 
\vspace{1ex}
с учетом теоремы 8.1 из свойства 5.7 при 
$\!p_j^{}\!=\!\nu_j^{},\, M_j^{}\!=-\;\!{\rm div}\;\!{\frak d},\, j\!=\!1,\ldots,m.\k\!$

На основании определения 5.2 вводится понятие кратности для кратного 
полиномиального последнего множителя.
\vspace{0.75ex}

{\bf Свойство 9.8.}
\vspace{0.75ex}
{\it
Пусть 
$\Bigl(\nu, 
\Bigl(h_\xi^{},q_{_{\scriptstyle h_\xi^{}  f_\xi^{}}}, N_{h_\xi^{}  f_\xi^{}}^{}\Bigr)\Bigr)
\in \text{\rm MB}_{_{\Xi^{\;\!\prime}}},\
\gamma,\gamma_{_{\scriptstyle h_\xi^{}  f_\xi^{}}}\in\R,\ 
\varphi_{_{\scriptstyle h_\xi^{}  f_\xi^{}}}\in C^1T^{\;\!\prime},$ 
$f_{\xi}^{}=1,\ldots, \delta_{\xi}^{},\ \xi=1,\ldots,\varepsilon,$
\vspace{1ex}
мно\-жес\-т\-во
$\Omega_{_0}\subset\Xi^{\;\!\prime}$ та\-кое, что  
$\nu(t,x)\ne 0\;\;\forall (t,x)\in\Omega_{_0},$
$\nu(t,x)= 0\;\;\forall (t,x)\in {\sf C}_{_{\Xi^{\;\!\prime}}}\Omega_{_0},\ 
\nu^{\gamma}\in C^1\Omega_{_0}.$
Тогда
\\[1.5ex]
\mbox{}\hfill                           
$
\displaystyle
\biggl(\ \nu^{\gamma}
\exp\sum\limits_{\xi=1}^{\varepsilon}\sum\limits_{f_{\xi}^{}=1}^{\delta_{\xi}^{}}
\biggl(\;\!\gamma_{_{\scriptstyle h_\xi^{}  f_\xi^{}}}\biggl(\,
\dfrac{q_{_{\scriptstyle h_\xi^{}  f_\xi^{}} }}{\displaystyle  \nu^{\;\!h_\xi^{}}}+
\varphi_{_{\scriptstyle h_\xi^{}  f_\xi^{}}}\biggr)\biggr),\, L\biggr)\in 
\text{\rm J}_{_{\Omega_{{}_{\tiny\;\! 0}}}},
\hfill
$
\\[2ex]
если и только если сомножитель}
\\[1.75ex]
\mbox{}\hfill                           
$
\displaystyle
L(t,x)=
\sum\limits_{\xi=1}^{\varepsilon}\sum\limits_{f_{\xi}^{}=1}^{\delta_{\xi}^{}}
\biggl(\gamma_{_{\scriptstyle h_\xi^{}  f_\xi^{}}}\Bigl(
N_{h_\xi^{}  f_\xi^{}}^{} (t,x)+{\sf D}\;\!\varphi_{_{\scriptstyle h_\xi^{}  f_\xi^{}}}(t)\Bigr)\biggr)
- \gamma\;\!{\rm div}\;\!{\frak d}(t,x)
\quad
\forall (t,x)\in \Xi^{\;\!\prime}.
\hfill
$
\\[1.5ex]
\indent
{\sl Следует} 
с учетом теоремы 8.1 из свойства 5.8 при 
$M={}-\,{\rm div}\;\!{\frak d}, \ p=\nu.\ \k$ 
\vspace{1ex}

{\bf Свойство 9.9.}\!
{\it
Пусть $(p, M)\!\in\! \text{\rm A}_{_{\!\Xi^{\;\!\prime}}},\!$ число $h\!\in\!\N,$ функции
\vspace{0.5ex}
$p, q\!\in\! \text{\rm P}_{_{\!\Xi^{\;\!\prime}}}\!$ взаимно простые, 
множество
\vspace{0.5ex}
$\Omega_{_0}\subset\Xi^{\;\!\prime}$ та\-кое, что  
$p(t,x)\ne 0\;\;\forall (t,x)\in\Omega_{_0},\
p(t,x)= 0\;\;\forall (t,x)\in {\sf C}_{_{\Xi^{\;\!\prime}}}\Omega_{_0}.$
Тогда $\bigl((p, M), (h,q, -\;\!{\rm div}\;\!{\frak d})\bigr)\in \text{\rm B}_{_{\Xi^{\;\!\prime}}},$
если и только если выполняется одно из условий}:
\\[1ex]
\indent
$
1)\ {\frak d}\,\dfrac{q(t,x)}{p^{\;\!h}(t,x)}={}-\,{\rm div}\;\!{\frak d}(t,x)
\quad
\forall (t,x)\in \Omega_{_0};
$
\\[1.25ex]
\indent
$
2)\ {\frak d}\;\!q(t,x)=h\;\!q(t,x)\;\!M(t,x)-
p^{\;\!h}(t,x)\;\!{\rm div}\;\!{\frak d}(t,x)
\quad
\forall (t,x)\in \Xi^{\;\!\prime};
$
\\[1.5ex]
\indent
$
3)\ 
\exp\dfrac{q}{p^{\;\!h}}\in 
\text{\rm ME}_{_{\Omega_{{}_{\tiny\;\! 0}}}}.
$
\\[1.5ex]
\indent
{\sl Следует} 
\vspace{1.25ex}
с учетом теоремы 8.1 из определения 5.1, теорем 5.2 и 5.3 при 
$N={}-\,{\rm div}\;\!{\frak d}.\ \k$ 

\newpage

{\bf Свойство 9.10.}
\vspace{0.5ex}
{\it
Пусть 
$\bigl((p_j^{}, M_j^{}), (h_j^{},q_j^{},\rho_{j}^{}N_0^{})\bigr)\in 
\text{\rm B}_{_{\Xi^{\;\!\prime}}},\ 
\rho_j^{}\in\R\backslash\{0\},\ 
\lambda_j^{},\gamma_j^{}, c_j^{}\in\R,
\linebreak
j=1,\ldots, m,\
\sum\limits_{j=1}^{m}\!|\lambda_j^{}|\!\ne\! 0,\
\varphi\in C^1T^{\;\!\prime},$ 
\vspace{0.5ex}
мно\-жес\-т\-во
$\Omega_{_0}\subset\Xi^{\;\!\prime}$ та\-кое, что  
$\prod\limits_{j=1}^{m}p_j^{}(t,x)\ne 0$
$\forall (t,x)\in\Omega_{_0},\
\prod\limits_{j=1}^{m}p_j^{}(t,x)= 0\;\;
\forall (t,x)\in {\sf C}_{_{\Xi^{\;\!\prime}}}\Omega_{_0},\
p_j^{\gamma_j^{}}\in C^1\Omega_{_0}, \ j=1,\ldots, m.$
Тогда{\rm:}
\\[1.5ex]
\indent
$
1)\
\displaystyle
\prod\limits_{j=1}^{m}p_j^{\gamma_j^{}}
\sum\limits_{j=1}^{m}\lambda_j^{}
\exp\biggl(\;\!\dfrac{q_j^{}}{\rho_j^{}\;\!p_j^{\;\!h_j^{}}}+\varphi\biggr)\in 
\text{\rm M}_{_{\Omega_{{}_{\tiny\;\! 0}}}},
$
если и только если
\\[2ex]
\mbox{}\hfill                           
$
\displaystyle
{\rm div}\;\!{\frak d}(t,x)={}-{\sf D}\;\!\varphi(t)-N_0^{}(t,x)-
\sum\limits_{j=1}^{m}\gamma_j^{}\;\!M_j^{}(t,x)
\quad
\forall (t,x)\in \Xi^{\;\!\prime};
\hfill
$
\\[1.75ex]
\indent
$
2)\ 
\displaystyle
\prod\limits_{j=1}^{m}p_j^{\gamma_j^{}}
\sum\limits_{j=1}^{m}\lambda_j^{}
\exp\biggl(\;\!\dfrac{q_j^{}}{\rho_j^{}\;\!p_j^{\;\!h_j^{}}}+c_j^{}\biggr)\in 
\text{\rm M}_{_{\Omega_{{}_{\tiny\;\! 0}}}},
$
если и только если}
\\[1.75ex]
\mbox{}\hfill                           
$
\displaystyle
{\rm div}\;\!{\frak d}(t,x)={}-N_0^{}(t,x)-
\sum\limits_{j=1}^{m}\gamma_j^{}\;\!M_j^{}(t,x)
\quad
\forall (t,x)\in \Xi^{\;\!\prime}.
\hfill
$
\\[1.5ex]
\indent
{\sl Следует} 
\vspace{1ex}
с учетом теоремы 8.1 из свойств 5.5 и 5.6 при 
$L={}-\,{\rm div}\;\!{\frak d}.\ \k
$ 

%\newpage

{\bf Свойство 9.11.}\!
\vspace{0.75ex}
{\it
Пусть 
$\bigl((p_j^{}, M_j^{}), (h_j^{},q_j^{},N_j^{})\bigr)\in \text{\rm B}_{_{\Xi^{\;\!\prime}}},\ 
\gamma_j^{},\xi_j^{}\!\in\!\R,\ \varphi_j^{}\!\in\! C^1T^{\;\!\prime},\, 
j\!=\!1,\ldots, m,\!$ 
мно\-жес\-т\-во
\vspace{0.75ex}
$\!\Omega_{_0}\!\subset\!\Xi^{\;\!\prime}\!$ та\-кое, что  
$\prod\limits_{j=1}^{m}p_j^{}(t,x)\!\ne\! 0\;\;\forall (t,x)\!\in\!\Omega_{_0},\
\prod\limits_{j=1}^{m}p_j^{}(t,x)= 0\;\;
\forall (t,x)\!\in\! {\sf C}_{_{\Xi^{\;\!\prime}}}\Omega_{_0},\!$ 
$p_j^{\gamma_j^{}}\in C^1\Omega_{_0}, \ j=1,\ldots, m.$
Тогда
\\[1.5ex]
\mbox{}\hfill                           
$
\displaystyle
\prod\limits_{j=1}^{m}p_j^{\gamma_j^{}}
\exp\sum\limits_{j=1}^{m}\xi_j^{}
\biggl(\;\!\dfrac{q_j^{}}{p_j^{\;\! h_j^{}}}+\varphi_j^{}\biggr)\in 
\text{\rm M}_{_{\Omega_{{}_{\tiny\;\! 0}}}},
\hfill
$
\\[2ex]
если и только если}
\\[1.75ex]
\mbox{}\hfill                           
$
\displaystyle
{\rm div}\;\!{\frak d}(t,x)=
{}-\sum\limits_{j=1}^{m}
\Bigl(\gamma_j^{}\;\!M_j^{}(t,x)+
\xi_j^{}\;\!\bigl(N_j^{}(t,x)+{\sf D}\;\!\varphi_j^{}(t)\bigr)\Bigr) 
\quad
\forall (t,x)\in \Xi^{\;\!\prime}.
\hfill
$
\\[1.5ex]
\indent
{\sl Следует} 
с учетом теоремы 8.1 из свойства 5.7 при 
$L={}-\, {\rm div}\;\!{\frak d}.\ \k$ 
\vspace{1.25ex}

{\bf Свойство 9.12.}
\vspace{0.75ex}
{\it
Пусть 
$\Bigl((p, M), 
\Bigl(h_\xi^{},q_{_{\scriptstyle h_\xi^{}  f_\xi^{}}}, N_{h_\xi^{}  f_\xi^{}}^{}\Bigr)\Bigr)
\in \text{\rm B}_{_{\Xi^{\;\!\prime}}},\
\gamma,\gamma_{_{\scriptstyle h_\xi^{}  f_\xi^{}}}\in\R,\ 
\varphi_{_{\scriptstyle h_\xi^{}  f_\xi^{}}}\in C^1T^{\;\!\prime},$ 
$f_{\xi}^{}=1,\ldots, \delta_{\xi}^{},\ \xi=1,\ldots,\varepsilon,$
\vspace{1ex}
мно\-жес\-т\-во
$\Omega_{_0}\subset\Xi^{\;\!\prime}$ та\-кое, что  
$p(t,x)\ne 0\;\;\forall (t,x)\in\Omega_{_0},$
$p(t,x)= 0\;\;\forall (t,x)\in {\sf C}_{_{\Xi^{\;\!\prime}}}\Omega_{_0},\ 
p^{\gamma}\in C^1\Omega_{_0}.$
Тогда
\\[1.5ex]
\mbox{}\hfill                           
$
\displaystyle
p^{\gamma}
\exp\sum\limits_{\xi=1}^{\varepsilon}\sum\limits_{f_{\xi}^{}=1}^{\delta_{\xi}^{}}
\biggl(\;\!\gamma_{_{\scriptstyle h_\xi^{}  f_\xi^{}}}\biggl(\,
\dfrac{q_{_{\scriptstyle h_\xi^{}  f_\xi^{}} }}{\displaystyle  p^{\;\!h_\xi^{}}}+
\varphi_{_{\scriptstyle h_\xi^{}  f_\xi^{}}}\biggr)\biggr)\in 
\text{\rm M}_{_{\Omega_{{}_{\tiny\;\! 0}}}},
\hfill
$
\\[2ex]
если и только если}
\\[1.75ex]
\mbox{}\hfill                           
$
\displaystyle
{\rm div}\;\!{\frak d}(t,x)=
{}-\sum\limits_{\xi=1}^{\varepsilon}\sum\limits_{f_{\xi}^{}=1}^{\delta_{\xi}^{}}
\biggl(\gamma_{_{\scriptstyle h_\xi^{}  f_\xi^{}}}\Bigl(
N_{h_\xi^{}  f_\xi^{}}^{} (t,x)+{\sf D}\;\!\varphi_{_{\scriptstyle h_\xi^{}  f_\xi^{}}}(t)\Bigr)\biggr)
- \gamma\;\!M(t,x)
\quad
\forall (t,x)\in \Xi^{\;\!\prime}.
\hfill
$
\\[1.5ex]
\indent
{\sl Следует} 
с учетом теоремы 8.1 из свойства 5.8 при 
$L={}-\,{\rm div}\;\!{\frak d}.\ \k$ 
\vspace{1ex}

{\bf Свойство 9.13.} 
{\it
Справедливо утверждение}
\\[1.5ex]
\mbox{}\hfill                         
$
\displaystyle
u+i\;\!v\in \text{\rm MH}_{_{\Xi^{\;\!\prime}}}
\iff
u\in \text{\rm MA}_{_{\Xi^{\;\!\prime}}}
\ \&\ 
v\in \text{\rm MA}_{_{\Xi^{\;\!\prime}}}.
\hfill
$
\\[1.75ex]
\indent
{\sl Следует} 
с учетом теоремы 8.1 из свойства 6.11 при 
$U={}-\,{\rm div}\;\!{\frak d}.\ \k$ 
\vspace{1ex}

{\bf Свойство 9.14.}
\vspace{0.15ex}
{\it
Пусть функции
\vspace{0.5ex}
$u,v\in \text{\rm P}_{_{\!\Xi^{\;\!\prime}}}$ взаимно простые, 
множество $\Omega_{_0}\subset \Xi^{\;\!\prime}$ такое, что  
$u(t,x)\ne 0\;\;\forall (t,x)\in \Omega_{_0},\ 
u(t,x)=0\;\;\forall (t,x)\in {\sf C}_{_{\Xi^{\;\!\prime}}}\Omega_{_0}.$
Тогда}
\\[1.5ex]
\mbox{}\hfill
$
u+i\;\!v\in \text{\rm MH}_{_{\Xi^{\;\!\prime}}}
\iff
(u^2+v^2, {}-2\;\!{\rm div}\;\!{\frak d})\in \text{\rm A}_{_{\Xi^{\;\!\prime}}}
\ \&\ \,
\dfrac{v}{u}\in \text{\rm I}_{_{\Omega_{{}_{\tiny\;\! 0}}}}.
\hfill
$
\\[1.5ex]
\indent
{\sl Следует} 
с учетом теоремы 8.1 из свойства 6.12 при 
$U={}-\,{\rm div}\;\!{\frak d}.\ \k$ 
\vspace{1ex}

{\bf Свойство 9.15.} 
{\it
Пусть $\eta\in\C\backslash\{0\}.$ Тогда}
\\[1.5ex]
\mbox{}\hfill
$
w\in \text{\rm MH}_{_{\Xi^{\;\!\prime}}}
\iff 
\eta\;\! w\in \text{\rm MH}_{_{\Xi^{\;\!\prime}}}.
\hfill
$
\\[1.25ex]
\indent
{\sl Следует} 
с учетом теоремы 8.1 из свойства 6.1 при 
$W={}-\,{\rm div}\;\!{\frak d}.\ \k$ 
\vspace{1ex}

{\bf Свойство 9.16.} 
{\it
Пусть $k\in\N.$ Тогда
\\[1.5ex]
\mbox{}\hfill
$
w^k\in \text{\rm MH}_{_{\Xi^{\;\!\prime}}}
\iff 
\Bigl(w, {}-\dfrac{1}{k}\, \;\!{\rm div}\;\!{\frak d}\Bigr)\in \text{\rm H}_{_{\Xi^{\;\!\prime}}},
\hfill
$
\\[1ex]
а}
\\[1ex]
\mbox{}\hfill
$
w\in \text{\rm MH}_{_{\Xi^{\;\!\prime}}}
\iff 
(w^k,{}- k\,{\rm div}\;\!{\frak d})\in \text{\rm H}_{_{\Xi^{\;\!\prime}}}.
\hfill
$
\\[1.5ex]
\indent
{\sl Следует} 
из свойства 6.3 с учетом теоремы 8.1. $\k$ 
\vspace{1ex}

%\newpage

{\bf Свойство 9.17.} 
{\it
Пусть $\eta_j^{}\in\R\backslash\{0\},\ j=1,\ldots,m.$ Тогда}
\\[1.5ex]
\mbox{}\hfill
$
\displaystyle
w_j^{}\in \text{\rm MH}_{_{\Xi^{\;\!\prime}}}, \ j=1,\ldots, m,
\ \ \Longrightarrow\ \ 
\sum\limits_{j=1}^{m} \eta_j^{}\;\!w_j^{}
\in \text{\rm MH}_{_{\Xi^{\;\!\prime}}}.
\hfill
$
\\[1.5ex]
\indent
{\sl Следует} 
с учетом теоремы 8.1 из свойства 6.4 при 
$W={}-\,{\rm div}\;\!{\frak d}.\ \k$ 
\vspace{1ex}

{\bf Свойство 9.18.} 
{\it
Пусть 
$k_j^{}\in\N, \ w_j^{}\in \text{\rm Z}_{_{\Xi^{\;\!\prime}}},\ j=1,\ldots,m.$ Тогда  
\\[1.5ex]
\mbox{}\hfill                         
$
\displaystyle
\prod\limits_{j=1}^{m} w_j^{k_j^{}} 
\in \text{\rm MH}_{_{\Xi^{\;\!\prime}}}
\iff
\bigl( w_j^{}\;\!,\;\! W_j^{}\bigl)\;\! \in \text{\rm H}_{_{\Xi^{\;\!\prime}}},
\quad
j=1,\ldots,m.
\hfill
$
\\[1.5ex]
При этом сомножители $W_1^{},\ldots,W_m^{}$ такие, что}
\\[1.5ex]
\mbox{}\hfill
$
\displaystyle
\sum\limits_{j=1}^{m} k_j^{}\;\!{\rm Re}\;\!W_j^{}(t,x)=
{}-\,{\rm div}\;\!{\frak d}(t,x)
\quad
\forall (t,x)\in \Xi^{\;\!\prime},
\qquad
\sum\limits_{j=1}^{m} k_j^{}\;\!{\rm Im}\;\!W_j^{}(t,x)=0
\quad
\forall (t,x)\in \Xi^{\;\!\prime}.
\hfill
$
\\[1.5ex]
\indent
{\sl Следует} 
с учетом теоремы 8.1 из свойства 6.6 при 
$W={}-\,{\rm div}\;\!{\frak d}.\ \k$ 
\vspace{1.25ex}

{\bf Свойство 9.19.} 
{\it
Если 
$p_\tau^{}\in \text{\rm P}_{_{\!\Xi^{\;\!\prime}}},\, l_{\tau}^{}\!\in\N,\, \tau\!=\!1,\ldots,s,\
w_j^{}\in \text{\rm Z}_{_{\Xi^{\;\!\prime}}},\, k_j^{}\!\in\N, \, j\!=\!1,\ldots,m,$ то  
\\[1.5ex]
\mbox{}\hfill                         
$
\displaystyle
\prod\limits_{\tau=1}^{s}\! p_\tau^{\,l_\tau^{}}
\prod\limits_{j=1}^{m}\! w_j^{k_j^{}} 
\in\! \text{\rm MH}_{_{\Xi^{\;\!\prime}}}
\!\iff\!
\bigl(p_\tau^{}\;\!, M_\tau^{}\bigl)\;\! \in\! \text{\rm A}_{_{\Xi^{\;\!\prime}}},\,
\tau\!=\!1,\ldots,s,
\ \&\ 
\bigl(w_j^{}\;\!, W_j^{}\bigl)\;\!\in\! \text{\rm H}_{_{\Xi^{\;\!\prime}}},\,
j\!=\!1,\ldots,m.
\hfill
$
\\[1.5ex]
При этом у частных интегралов сомножители $M_1^{},\ldots,M_s^{},\, W_1^{},\ldots,W_m^{}$ такие, что}
\\[1.5ex]
\mbox{}\hfill
$
\displaystyle
\sum\limits_{\tau=1}^{s} l_\tau^{}\;\!M_\tau^{}(t,x)+
\sum\limits_{j=1}^{m} k_j^{}\;\!{\rm Re}\;\!W_j^{}(t,x)=
{}-\;\! {\rm div}\;\!{\frak d}(t,x),
\quad
\sum\limits_{j=1}^{m} k_j^{}\;\!{\rm Im}\;\!W_j^{}(t,x)=0
\quad
\forall (t,x)\in \Xi^{\;\!\prime}.
\hfill
$
\\[1.5ex]
\indent
{\sl Следует} 
с учетом теоремы 8.1 из свойства 6.7 при 
$W={}-\,{\rm div}\;\!{\frak d}.\ \k$ 
\vspace{1ex}

{\bf Свойство 9.20.} 
\vspace{0.15ex}
{\it
Пусть 
$p\in \text{\rm P}_{_{\!\Xi^{\;\!\prime}}},$ 
функции $u,v\in\! \text{\rm P}_{_{\!\Xi^{\;\!\prime}}}$ 
взаимно простые, $w=u+i\;\!v.$ Тогда  
\\[0.5ex]
\mbox{}\hfill                         
$
\displaystyle
p\;\!w\in \text{\rm MH}_{_{\Xi^{\;\!\prime}}}
\iff
(p, M) \in \text{\rm A}_{_{\Xi^{\;\!\prime}}}
\ \&\ 
(w, W) \in \text{\rm H}_{_{\Xi^{\;\!\prime}}}.
\hfill
$
\\[1.75ex]
При этом у частных интегралом сомножители $M$ и $W$ такие, что}
\\[1.5ex]
\mbox{}\hfill
$
\displaystyle
M(t,x)+{\rm Re}\;\!W(t,x)={}-\,{\rm div}\;\!{\frak d}(t,x)
\quad
\forall (t,x)\in \Xi^{\;\!\prime},
\qquad
{\rm Im}\;\!W(t,x)=0
\quad
\forall (t,x)\in \Xi^{\;\!\prime}.
\hfill
$
\\[1.75ex]
\indent
{\sl Следует} 
 из свойства 6.8 с учетом теоремы 8.1. $\k$ 
\vspace{1ex}

{\bf Свойство 9.21.} 
{\it
Справедливо утверждение}
\\[1.5ex]
\mbox{}\hfill                         
$
\displaystyle
u+i\;\!v\in \text{\rm MH}_{_{\Xi^{\;\!\prime}}}
\iff
u-i\;\!v\in \text{\rm MH}_{_{\Xi^{\;\!\prime}}}.
\hfill
$
\\[1.5ex]
\indent
{\sl Следует} 
с учетом теоремы 8.1 из свойства 6.9 при 
$V=0,\ U={}-\,{\rm div}\;\!{\frak d}.\ \k$ 
\vspace{1ex}

{\bf Свойство 9.22}
\vspace{0.15ex}
(критерий существования кратного комплекснозначного
полиномиального последнего множителя).
{\it
Пусть $h\in\N,$ функции
\vspace{0.5ex}
$w, z\in \text{\rm Z}_{_{\Xi^{\;\!\prime}}}\!$ взаимно простые, 
множество
$\Omega_{_0}\subset\Xi^{\;\!\prime}$ та\-кое, что  
$|w(t,x)|\ne 0\;\;\forall (t,x)\in\Omega_{_0},\
|w(t,x)|= 0\;\;\forall (t,x)\in {\sf C}_{_{\Xi^{\;\!\prime}}}\Omega_{_0}.$
Тогда}
\\[0.5ex]
\mbox{}\hfill                           
$
\displaystyle
\bigl(w, (h,z,Q)\bigr)\in \text{\rm MG}_{_{\Xi^{\;\!\prime}}}
\iff
w\in \text{\rm MH}_{_{\Xi^{\;\!\prime}}}
\ \& \
\hfill                           
$
\\[2ex]
\mbox{}\hfill                           
$
\& \
\Bigl(\exp\;\!{\rm Re}\;\!\dfrac{z}{w^{\;\!h}}\,, {\rm Re}\;\!Q\Bigr)\in 
\text{\rm E}_{_{\Omega_{{}_{\tiny\;\! 0}}}}
\ \,\& \ \,
\Bigl(\exp\;\!{\rm Im}\;\!\dfrac{z}{w^{\;\!h}}\,, {\rm Im}\;\!Q\Bigr)\in 
\text{\rm E}_{_{\Omega_{{}_{\tiny\;\! 0}}}}.
\hfill
$
\\[1.75ex]
\indent
{\sl Следует} 
с учетом теоремы 8.1 из теоремы 7.2 при 
$M={}-\,{\rm div}\;\!{\frak d}.\ \k$ 
\vspace{1.5ex}

{\bf Свойство 9.23}
\vspace{0.15ex}
(критерий существования кратного комплекснозначного
полиномиального последнего множителя).
{\it
Пусть $h\in\N,$ функции
\vspace{0.5ex}
$u, v\in \text{\rm P}_{_{\!\Xi^{\;\!\prime}}}$ взаимно простые, 
функция $z\in \text{\rm Z}_{_{\Xi^{\;\!\prime}}}$ взаимно простая с 
функцией $u+i\;\!v,$  множество
\vspace{0.75ex}
$\Omega_{_0}\subset\Xi^{\;\!\prime}$ та\-кое, что  
$u(t,x)\ne 0\;\;\forall (t,x)\in\Omega_{_0},\
u(t,x)= 0\;\;\forall (t,x)\in {\sf C}_{_{\Xi^{\;\!\prime}}}\Omega_{_0}.$
Тогда}
\\[1.75ex]
\mbox{}\hfill                           
$
\displaystyle
\bigl(u+i\;\!v, (h,z,Q)\bigr)\in \text{\rm MG}_{_{\Xi^{\;\!\prime}}}
\iff
\Bigl(\bigl(u^2+v^2, {}-2\,{\rm div}\;\!{\frak d}\bigr), 
\bigl(h,\;\!{\rm Re}\;\!\bigl(z(u-i\;\!v)^h\bigr),\;\! {\rm Re}\;\!Q\bigr)\Bigr)\in 
\text{\rm B}_{_{\Xi^{\;\!\prime}}}
\ \& \
\hfill                           
$
\\[1.75ex]
\mbox{}\hfill                           
$
\& \
\Bigl(\bigl(u^2+v^2, {}- 2\,{\rm div}\;\!{\frak d}\bigr), 
\bigl(h,\;\!{\rm Im}\;\!\bigl(z(u-i\;\!v)^h\bigr),\;\! {\rm Im}\;\!Q\bigr)\Bigr)\in 
\text{\rm B}_{_{\Xi^{\;\!\prime}}}
\ \,\& \ \,
\dfrac{v}{u}\in 
\text{\rm I}_{_{\Omega_{{}_{\tiny\;\! 0}}}}.
\hfill
$
\\[1.75ex]
\indent
{\sl Следует} 
с учетом теоремы 8.1 из теоремы 7.3 при 
$M={}-\,{\rm div}\;\!{\frak d}.\ \k$ 
\vspace{1ex}

На основании определения 7.2 вводится понятие кратности для кратного 
комплекснозначного полиномиального последнего множителя.
\\[3.75ex]
\centerline{
{\bf  10. Экспоненциальные последние множители 
}
}
\\[1.75ex]
\indent
{\bf Теорема 10.1}
\vspace{0.35ex}
(критерий существования экспоненциального последнего множителя).
{\it
Функция $\exp\omega\in \text{\rm ME}_{_{\Omega}}$ 
тогда и только тогда, когда выполняется тождество}
\\[1.5ex]
\mbox{}\hfill               
$
{\frak d}\;\!\omega(t,x)={}-\,{\rm div}\;\!{\frak d}(t,x)
\quad
\forall (t,x)\in \Omega.
\hfill
$
\\[1.5ex]
\indent
{\sl Следует} 
с учетом теоремы 8.1 из теоремы 3.1 при 
$M={}-\,{\rm div}\;\!{\frak d}.\ \k$ 
\vspace{1ex}

{\bf Свойство 10.1.}
{\it
Если
\\[1ex]
\mbox{}\hfill               
$
{\rm div}\;\!{\frak d}(t,x)=\varphi(t)
\quad
\forall (t,x)\in \Xi,
\hfill
$
\\[1ex]
то функция
\\[1ex]
\mbox{}\hfill               
$
\displaystyle
\mu\colon t\to\ 
\exp\biggl({}-\int\limits_{t_0^{}}^{t}\varphi(\tau)\;\!d\tau\biggr)
\quad
\forall t\in T,
\hfill
$
\\[1.5ex]
где $t_{_0}^{}$ --- произвольная фиксированная точка из области $T,$
\vspace{0.35ex}
будет условным последним множителем на области $\Xi$ системы {\rm (0.1).}
}
\vspace{0.5ex}

{\sl Доказательство.}
Производная в силу системы (0.1)
\\[1.5ex]
\mbox{}\hfill               
$
\displaystyle
{\frak d}\;\!\biggl({}-\int\limits_{t_0^{}}^{t}\varphi(\tau)\;\!d\tau\biggr)=
{\sf D}\;\!\biggl({}-\int\limits_{t_0^{}}^{t}\varphi(\tau)\;\!d\tau\biggr)=
{}-\varphi(t)=
{}-\;\!{\rm div}\;\!{\frak d}(t,x)
\quad
\forall (t,x)\in \Xi.
\hfill
$
\\[1.5ex]
\indent
По теореме 10.1, $\mu\in  \text{\rm MF}_{_{\!\Xi}}.\ \k$
\vspace{1ex}

Частным случаем свойства 10.1 является
\vspace{0.75ex}

{\bf Свойство 10.2.}
{\it
Если
\\[1ex]
\mbox{}\hfill               
$
{\rm div}\;\!{\frak d}(t,x)=\lambda
\quad
\forall (t,x)\in \Xi,
\quad
\lambda\in\R,
\hfill
$
\\[1ex]
то функция $e^{{}-\lambda\;\!t}\in  \text{\rm MF}_{_{\!\Xi}}.$
}
\vspace{1ex}

{\bf  Свойство 10.3.}
{\it
Пусть $\gamma_j^{}\in\R\backslash\{0\},\ c_j^{}\in\R, \ j=1,\ldots,m.$
Тогда экспоненциальная функция
$
\exp\sum\limits_{j=1}^{m}
\gamma_j^{}\;\!(\omega_j^{}+c_j^{})\in \text{\rm ME}_{_{\Omega}},
$
если и только если}
\\[1ex]
\mbox{}\hfill
$
\displaystyle
\sum\limits_{j=1}^{m}\gamma_j^{}\;\!{\frak d}\;\!\omega_j^{}(t,x)=
{}-\;\!{\rm div}\;\!{\frak d}(t,x)
\quad
\forall (t,x)\in \Omega.
\hfill
$
\\[1.5ex]
\indent
{\sl Следует} из теоремы 3.3 с учетом теоремы 8.1. $\k$
\vspace{0.75ex}

{\bf  Свойство 10.4.}
{\it
Пусть $\varphi\in C^1T^{\;\!\prime}.$ Тогда 
\\[1.5ex]
\mbox{}\hfill
$
\mu\in \text{\rm M}_{_{\Omega}}
\iff
\bigl(\mu\exp\varphi, {\sf D} \varphi-\;\!{\rm div}\;\!{\frak d}\bigr)\in \text{\rm J}_{_{\Omega}},
\hfill
$
\\[0.5ex]
а}
\\[0.5ex]
\mbox{}\hfill
$
{\rm g}\exp\varphi\in \text{\rm M}_{_{\Omega}}
\iff
\bigl({\rm g}, {}-{\sf D} \varphi-\;\!{\rm div}\;\!{\frak d}\bigr)\in \text{\rm J}_{_{\Omega}}.
\hfill
$
\\[1.5ex]
\indent
{\sl Следует} из свойства 3.1 с учетом теоремы 8.1. $\k$
\vspace{1ex}

{\bf  Свойство 10.5.}
{\it
Пусть $\varphi\in C^1T^{\;\!\prime}.$ Тогда 
\\[1.5ex]
\mbox{}\hfill
$
\exp\omega\in \text{\rm ME}_{_{\Omega}}
\iff
\bigl(\exp(\omega+\varphi),\, {\sf D} \varphi-\;\!{\rm div}\;\!{\frak d}\bigr)\in \text{\rm E}_{_{\Omega}},
\hfill
$
\\[0.5ex]
а}
\\[0.5ex]
\mbox{}\hfill
$
\exp(\omega+\varphi)\in \text{\rm ME}_{_{\Omega}}
\iff
\bigl(\exp\omega,\, {}-{\sf D} \varphi-\;\!{\rm div}\;\!{\frak d}\bigr)\in \text{\rm E}_{_{\Omega}}.
\hfill
$
\\[1.5ex]
\indent
{\sl Следует} из свойства 3.2 с учетом теоремы 8.1. $\k$
\vspace{1ex}

{\bf  Свойство 10.6.}
{\it
Если $c\in\R,$ то} 
\\[1.5ex]
\mbox{}\hfill
$
\displaystyle
\exp(\omega+c)\in \text{\rm ME}_{_{\Omega}}
\iff
\exp\omega\in \text{\rm ME}_{_{\Omega}}.
\hfill
$
\\[1.5ex]
\indent
{\sl Следует}  с учетом теоремы 8.1 из свойства 3.3 при 
$M={}-\;\!{\rm div}\;\!{\frak d}.\ \k$
\vspace{1ex}

{\bf  Свойство 10.7.}
{\it
Пусть  $\lambda_j^{}\in\R\backslash\{0\}, \ j=1,\ldots,m.$ Тогда}
\\[1.5ex]
\mbox{}\hfill
$
\displaystyle
\exp\omega_j^{}\in \text{\rm ME}_{_{\Omega}},\  j=1,\ldots,m, \ \ 
\Longrightarrow\ \ 
\sum\limits_{j=1}^{m}\lambda_j^{}\exp\omega_j^{}\in \text{\rm M}_{_{\Omega}}.
\hfill
$
\\[1.5ex]
\indent
{\sl Следует}  с учетом теоремы 8.1 из свойства 3.4 при 
$M={}-\;\!{\rm div}\;\!{\frak d}.\ \k$
\vspace{1ex}

{\bf  Свойство 10.8.}
{\it
Пусть
$\rho_j^{},\lambda_j^{}\in\R\backslash\{0\}, \ j=1,\ldots,m.$ Тогда} 
\\[1.5ex]
\mbox{}\hfill
$
\displaystyle
\bigl(\exp\omega_j^{}\;\!, {}-\rho_j^{}\;\!{\rm div}\;\!{\frak d}\bigr)\in \text{\rm E}_{_{\Omega}},
\ j=1,\ldots,m,
\ \ \Longrightarrow\ \ 
\sum\limits_{j=1}^{m}\lambda_j^{}
\exp\dfrac{\omega_j^{}}{\rho_j^{}}\in \text{\rm M}_{_{\Omega}}.
\hfill
$
\\[1.5ex]
\indent
{\sl Следует}  с учетом теоремы 8.1 из свойства 3.5 при 
$M_{_0}^{}={}-\;\!{\rm div}\;\!{\frak d}.\ \k$
\vspace{1ex}

{\bf  Свойство 10.9.}
\vspace{0.5ex}
{\it
Пусть  
$\bigl(\exp\omega_j^{}, M_j^{}\bigr)\in \text{\rm E}_{_{\Omega}},\ 
\gamma_j^{}\in\R\backslash\{0\}, \ j=1,\ldots,m.$ Тогда функция
$
\exp\sum\limits_{j=1}^{m}\gamma_j^{}\;\!\omega_j^{}\in \text{\rm ME}_{_{\Omega}},
$
если и только если выполняется тождество} (8.1).
\vspace{0.75ex}

{\sl Следует}  с учетом теоремы 8.1 из свойства 3.6 при 
$M={}-\;\!{\rm div}\;\!{\frak d}.\ \k$
\vspace{1ex}

{\bf  Свойство 10.10.}
\vspace{0.5ex}
{\it
Пусть  
$\bigl(\exp\omega_j^{}, \rho_j^{}M_{_0}^{}\bigr)\in \text{\rm E}_{_{\Omega}},\ 
\rho_j^{},\gamma_j^{}\in\R\backslash\{0\}, \ j=1,\ldots,m.$ Тогда 
$
\exp\sum\limits_{j=1}^{m}\gamma_j^{}\;\!\omega_j^{}\in \text{\rm ME}_{_{\Omega}},
$
если и только если выполняется тождество} 
\\[1ex]
\mbox{}\hfill
$
\displaystyle
\sum\limits_{j=1}^{m}
\rho_j^{}\;\!\gamma_j^{}\;\!M_{_0}^{}(t,x)=
{}-\;\!{\rm div}\;\!{\frak d}(t,x)
\quad
\forall (t,x)\in \Xi^{\;\!\prime}.
\hfill
$
\\[1.5ex]
\indent
{\sl Следует}  из свойства 10.9 при 
$M_j^{}=\rho_j^{}\;\!M_{_0}^{},\ j=1,\ldots, m.\ \k$
\vspace{1ex}

{\bf  Свойство 10.11.}
\vspace{0.5ex}
{\it
Пусть  
$\bigl(\exp\omega_j^{}, {}-\rho_j^{}\;\!{\rm div}\;\!{\frak d}\bigr)\in \text{\rm E}_{_{\Omega}},\ 
\rho_j^{},\gamma_j^{}\in\R\backslash\{0\}, \ j=1,\ldots,m.$ Тогда 
$
\exp\sum\limits_{j=1}^{m}\gamma_j^{}\;\!\omega_j^{}\in \text{\rm ME}_{_{\Omega}},
$
если и только если} 
$\sum\limits_{j=1}^{m}\rho_j^{}\;\!\gamma_j^{}=1.$
\vspace{0.75ex}

{\sl Следует}  из свойства 10.10 при 
$M_{_0}^{}={}-\;\!{\rm div}\;\!{\frak d}.\ \k$
\vspace{1ex}

{\bf  Свойство 10.12.}
\vspace{0.5ex}
{\it
Пусть  
$\exp\omega_j^{}\in \text{\rm ME}_{_{\Omega}},\ 
\gamma_j^{}\in\R\backslash\{0\}, \ j=1,\ldots,m.$ Тогда 
$
\Bigl(\exp\sum\limits_{j=1}^{m}\gamma_j^{}\;\!\omega_j^{}, M\Bigr)\in \text{\rm E}_{_{\Omega}},
$
если и только если выполняется тождество} (8.2). 
\vspace{0.75ex}

{\sl Следует} с учетом теоремы 8.1 из свойства 3.6 при 
$M_{j}^{}={}-\;\!{\rm div}\;\!{\frak d}, \ j=1,\ldots,m.\ \k$
\vspace{1ex}

{\bf  Свойство 10.13.}
\vspace{0.5ex}
{\it
Пусть  
$\exp\omega_j^{}\in \text{\rm ME}_{_{\Omega}},\ 
\gamma_j^{}\in\R\backslash\{0\}, \ j=1,\ldots,m.$ Тогда 
экспоненциальная функция
$
\exp\sum\limits_{j=1}^{m}\gamma_j^{}\;\!\omega_j^{}\in \text{\rm ME}_{_{\Omega}},
$
если и только если}
$\sum\limits_{j=1}^{m}\gamma_j^{}=1.$
\vspace{0.75ex}

{\sl Следует} из свойства 10.11 при 
$\rho_{j}^{}=1,\ j=1,\ldots,m.\ \k$
\vspace{1ex}

%\newpage

{\bf  Свойство 10.14.}
{\it
Если $\gamma\in\R\backslash\{0\},$ то 
\\[1.5ex]
\mbox{}\hfill
$
\displaystyle
\exp\omega\in \text{\rm ME}_{_{\Omega}}
\iff
\bigl(\exp(\gamma\;\!\omega), {}-\gamma\;\!{\rm div}\;\!{\frak d}\bigr)\in \text{\rm E}_{_{\Omega}},
\hfill
$
\\[0.35ex]
а}
\\[0.35ex]
\mbox{}\hfill
$
\displaystyle
\exp(\gamma\;\!\omega)\in \text{\rm ME}_{_{\Omega}}
\iff
\Bigl(\exp\omega, {}-\dfrac{1}{\gamma}\ {\rm div}\;\!{\frak d}\bigr)\in \text{\rm E}_{_{\Omega}}.
\hfill
$
\\[1.35ex]
\indent
{\sl Следует} с учетом теоремы 8.1 из следствия 3.2. $\k$ 
\vspace{1ex}

{\bf  Свойство 10.15.}
{\it
Если 
$
\bigl(\exp\omega_j^{}, M_j^{}\bigr)\in \text{\rm E}_{_{\Omega}},\ 
\gamma,\gamma_j^{}\in\R\backslash\{0\}, \ j=1,\ldots,m, \
{\rm g}^{\gamma}\in C^1\Omega,$ то}
\\[1.5ex]
\mbox{}\hfill
$
\displaystyle
{\rm g}^{\gamma}\exp\sum\limits_{j=1}^{m}\gamma_j^{}\;\!\omega_j^{}\in 
\text{\rm M}_{_{\Omega}}
\iff
\biggl(\,
{\rm g}, {}-\dfrac{1}{\gamma}\,\biggl({\rm div}\;\!{\frak d}+
\sum\limits_{j=1}^{m}\gamma_j^{}\;\!M_j^{}\biggr)\biggr)\in 
\text{\rm J}_{_{\Omega}}.
\hfill
$
\\[1.5ex]
\indent
{\sl Следует} с учетом теоремы 8.1 из свойства 3.7 при 
$M={}-\;\!{\rm div}\;\!{\frak d}.\ \k$
\vspace{1ex}

{\bf  Свойство 10.16.}
{\it
Пусть  
$
\exp\omega_j^{}\in \text{\rm ME}_{_{\Omega}},\ 
\gamma,\gamma_j^{}\in\R\backslash\{0\}, \ j=1,\ldots,m, \
{\rm g}^{\gamma}\in C^1\Omega.$ Тогда}
\\[1.5ex]
\mbox{}\hfill
$
\displaystyle
\biggl(\,
{\rm g}^{\gamma}\exp\sum\limits_{j=1}^{m}\gamma_j^{}\;\!\omega_j^{}\;\!,\;\! M\biggr)\in 
\text{\rm J}_{_{\Omega}}
\iff
\biggl(\,
{\rm g},\, \dfrac{1}{\gamma}\,\biggl(M+
\sum\limits_{j=1}^{m}\gamma_j^{}\;\!{\rm div}\;\!{\frak d}\biggr)\biggr)\in 
\text{\rm J}_{_{\Omega}}.
\hfill
$
\\[1.5ex]
\indent
{\sl Следует} с учетом теоремы 8.1 из свойства 3.7 при 
$M_j^{}={}-\;\!{\rm div}\;\!{\frak d},\ j=1,\ldots, m.\ \k$
\vspace{1ex}

{\bf  Свойство 10.17.}
{\it
Пусть  
$
\exp\omega_j^{}\in \text{\rm ME}_{_{\Omega}},\ 
\gamma,\gamma_j^{}\in\R\backslash\{0\}, \ j=1,\ldots,m, \
{\rm g}^{\gamma}\in C^1\Omega.$ Тогда}
\\[1.5ex]
\mbox{}\hfill
$
\displaystyle
{\rm g}^{\gamma}\exp\sum\limits_{j=1}^{m}\gamma_j^{}\;\!\omega_j^{}\in 
\text{\rm M}_{_{\Omega}}
\iff
\biggl(\,
{\rm g},\, \dfrac{1}{\gamma}\,\biggl(\
\sum\limits_{j=1}^{m}\gamma_j^{}-1\biggr)\;\!{\rm div}\;\!{\frak d}\biggr)\in 
\text{\rm J}_{_{\Omega}}.
\hfill
$
\\[1.5ex]
\indent
{\sl Следует} с учетом теоремы 8.1 из свойства 10.16 при 
$M={}-\;\!{\rm div}\;\!{\frak d}.\ \k$
\vspace{1ex}

{\bf  Свойство 10.18.}
{\it
Пусть
$
\bigl(\exp\omega_\tau^{}, M_\tau^{}\bigr)\!\in \text{\rm E}_{_{\Omega}},\
\tau\!=\!1,\ldots, m\!-\!1,\ 
\gamma_j^{}\in\R\backslash\{0\}, \ j\!=\!1,\ldots,m.$ Тогда}
\\[1ex]
\mbox{}\hfill
$
\displaystyle
\exp\sum\limits_{j=1}^{m}\gamma_j^{}\;\!\omega_j^{}\in 
\text{\rm ME}_{_{\Omega}}
\iff
\biggl(\,
\exp\omega_{m}^{}\;\!,{}-
 \dfrac{1}{\gamma_{m}^{}}\,\biggl({\rm div}\;\!{\frak d}+
\sum\limits_{\tau=1}^{m-1}\gamma_\tau^{}\;\!M_\tau^{}\biggr)\biggr)\in 
\text{\rm E}_{_{\Omega}}.
\hfill
$
\\[1.5ex]
\indent
{\sl Следует} с учетом теоремы 8.1 из следствия 3.3 при 
$M={}-\;\!{\rm div}\;\!{\frak d}.\ \k$
\vspace{1ex}

{\bf  Свойство 10.19.}
{\it
Если
$
\exp\omega_\tau^{}\!\in \text{\rm ME}_{_{\Omega}},\
\tau\!=\!1,\ldots, m\!-\!1,\ 
\gamma_j^{}\in\R\backslash\{0\}, \ j\!=\!1,\ldots,m,$ то}
\\[1.35ex]
\mbox{}\hfill
$
\displaystyle
\biggl(\,
\exp\sum\limits_{j=1}^{m}\gamma_j^{}\;\!\omega_j^{}\;\!,\;\! M\biggr)\in 
\text{\rm E}_{_{\Omega}}
\iff
\biggl(\,
\exp\omega_{m}^{}\;\!,\,
 \dfrac{1}{\gamma_{m}^{}}\,\biggl(M+
\sum\limits_{\tau=1}^{m-1}\gamma_\tau^{}\;\!{\rm div}\;\!{\frak d}\biggr)\biggr)\in 
\text{\rm E}_{_{\Omega}}.
\hfill
$
\\[1.35ex]
\indent
{\sl Следует} с учетом теоремы 8.1 из следствия 3.3 при 
\vspace{1ex}
$M_{\tau}^{}={}-\;\!{\rm div}\;\!{\frak d},\ \tau=1,\ldots,m-1.\ \k$

{\bf  Свойство 10.20.}
{\it
Если
$
\exp\omega_\tau^{}\!\in \text{\rm ME}_{_{\Omega}},\
\tau\!=\!1,\ldots, m\!-\!1,\ 
\gamma_j^{}\in\R\backslash\{0\}, \ j\!=\!1,\ldots,m,$ то}
\\[1.5ex]
\mbox{}\hfill
$
\displaystyle
\exp\sum\limits_{j=1}^{m}\gamma_j^{}\;\!\omega_j^{}\in 
\text{\rm ME}_{_{\Omega}}
\iff
\biggl(\,
\exp\omega_{m}^{}\;\!,\,
 \dfrac{1}{\gamma_{m}^{}}\,\biggl(\,
\sum\limits_{\tau=1}^{m-1}\gamma_\tau^{}-1\biggr)\;\!{\rm div}\;\!{\frak d}
\biggr)\in 
\text{\rm E}_{_{\Omega}}.
\hfill
$
\\[1.5ex]
\indent
{\sl Следует} с учетом теоремы 8.1 из свойства 10.19 при 
$M={}-\;\!{\rm div}\;\!{\frak d}.\ \k$
\vspace{1ex}

{\bf  Свойство 10.21.}
\vspace{0.5ex}
{\it
Пусть 
$
\bigl(\exp\omega_\nu^{}, M_\nu^{}\bigr)\in \text{\rm E}_{_{\Omega}},\ 
\nu=1,\ldots, s,\ s\leqslant m-2,\  
\gamma_j^{}\in\R\backslash\{0\}, 
\linebreak 
j=1,\ldots,m,\ 
{\rm g}_k^{{}^{\scriptsize \gamma_k^{}}}\!\in C^1\Omega,\ 
k=s+1,\ldots,m.$ Тогда
\\[1.75ex]
\mbox{}\hfill
$
\displaystyle
\prod\limits_{k=s+1}^{m}
{\rm g}_k^{{}^{\scriptsize \gamma_k^{}}}
\exp\sum\limits_{\nu=1}^{s}
\gamma_\nu^{}\;\!\omega_\nu^{}\in 
\text{\rm M}_{_{\Omega}}
\iff
\biggl(\,
\prod\limits_{k=s+1}^{m}
{\rm g}_k^{{}^{\scriptsize \gamma_k^{}}},\,-\;\!
{\rm div}\;\!{\frak d}-\sum\limits_{\nu=1}^{s}\gamma_\nu^{}\;\!M_\nu^{}\biggr)\in 
\text{\rm J}_{_{\Omega}}.
\hfill
$
\\[1.75ex]
Кроме этого существуют такие функции $M_k^{}\in C^1\Omega,\ k=s+1,\ldots, m,$
что выполняются тождества {\rm (1.14)} и}
\\[1.5ex]
\mbox{}\hfill                     
$
\displaystyle
\sum\limits_{k=s+1}^{m}
\gamma_{k}^{}\;\!M_k^{}(t,x)=
{}-\;\!{\rm div}\;\!{\frak d}(t,x)-
\sum\limits_{\nu=1}^{s}
\gamma_{\nu}^{}\;\!M_\nu^{}(t,x)
\quad
\forall (t,x)\in\Omega.
\hfill
$
\\[1.5ex]
\indent
{\sl Следует} с учетом теоремы 8.1 из свойства 3.8 при 
$M={}-\;\!{\rm div}\;\!{\frak d}.\ \k$
\vspace{1ex}

{\bf  Свойство 10.22.}
\vspace{0.5ex}
{\it
Пусть 
$
\exp\omega_\nu^{}\in \text{\rm ME}_{_{\Omega}},\ 
\nu=1,\ldots, s,\ s\leqslant m-2,\  
\gamma_j^{}\in\R\backslash\{0\}, 
\linebreak 
j=1,\ldots,m, \
{\rm g}_k^{{}^{\scriptsize \gamma_k^{}}}\!\in C^1\Omega,\ 
k=s+1,\ldots,m.$ Тогда
\\[1.75ex]
\mbox{}\hfill
$
\displaystyle
\biggl(\,
\prod\limits_{k=s+1}^{m}
{\rm g}_k^{{}^{\scriptsize \gamma_k^{}}}
\exp\sum\limits_{\nu=1}^{s}
\gamma_\nu^{}\;\!\omega_\nu^{}\;\!,\;\! M\biggr)\in 
\text{\rm J}_{_{\Omega}}
\iff
\biggl(\,
\prod\limits_{k=s+1}^{m}
{\rm g}_k^{{}^{\scriptsize \gamma_k^{}}}\;\!,\;\! M+
\sum\limits_{\nu=1}^{s}\gamma_\nu^{}\;\!{\rm div}\;\!{\frak d}\biggr)\in 
\text{\rm J}_{_{\Omega}}.
\hfill
$
\\[1.75ex]
Кроме этого существуют такие функции $M_k^{}\in C^1\Omega,\ k=s+1,\ldots, m,$
что выполняются тождества {\rm (1.14)} и}
\\[1.5ex]
\mbox{}\hfill                     
$
\displaystyle
\sum\limits_{k=s+1}^{m}
\gamma_{k}^{}\;\!M_k^{}(t,x)=
M(t,x)+
\sum\limits_{\nu=1}^{s}
\gamma_{\nu}^{}\;\!{\rm div}\;\!{\frak d}(t,x)
\quad
\forall (t,x)\in\Omega.
\hfill
$
\\[1.5ex]
\indent
{\sl Следует} с учетом теоремы 8.1 из свойства 3.8 при 
\vspace{1ex}
$M_{\nu}^{}={}-\;\!{\rm div}\;\!{\frak d},\ \nu=1,\ldots, s.\ \k$

{\bf  Свойство 10.23.}
\vspace{0.5ex}
{\it
Пусть 
$
\exp\omega_\nu^{}\in \text{\rm ME}_{_{\Omega}},\ 
\nu=1,\ldots, s,\ s\leqslant m-2,\  
\gamma_j^{}\in\R\backslash\{0\}, 
\linebreak 
j=1,\ldots,m, \
{\rm g}_k^{{}^{\scriptsize \gamma_k^{}}}\!\in C^1\Omega,\ 
k=s+1,\ldots,m.$ Тогда
\\[1.75ex]
\mbox{}\hfill
$
\displaystyle
\prod\limits_{k=s+1}^{m}
{\rm g}_k^{{}^{\scriptsize \gamma_k^{}}}
\exp\sum\limits_{\nu=1}^{s}
\gamma_\nu^{}\;\!\omega_\nu^{}\in 
\text{\rm M}_{_{\Omega}}
\iff
\biggl(\,
\prod\limits_{k=s+1}^{m}
{\rm g}_k^{{}^{\scriptsize \gamma_k^{}}}\;\!,\;\! 
\biggl(\,
\sum\limits_{\nu=1}^{s}\gamma_\nu^{}-1\biggr)\;\!{\rm div}\;\!{\frak d}\biggr)\in 
\text{\rm J}_{_{\Omega}}.
\hfill
$
\\[1.75ex]
Кроме этого существуют такие функции $M_k^{}\in C^1\Omega,\ k=s+1,\ldots, m,$
что выполняются тождества {\rm (1.14)} и}
\\[1.5ex]
\mbox{}\hfill                     
$
\displaystyle
\sum\limits_{k=s+1}^{m}
\gamma_{k}^{}\;\!M_k^{}(t,x)=
\biggl(\,
\sum\limits_{\nu=1}^{s}\gamma_\nu^{}-1\biggr)\;\!{\rm div}\;\!{\frak d}(t,x)
\quad
\forall (t,x)\in\Omega.
\hfill
$
\\[1.5ex]
\indent
{\sl Следует} с учетом теоремы 8.1 из свойства 10.22 при 
$M={}-\;\!{\rm div}\;\!{\frak d}.\ \k$
\vspace{1ex}

{\bf  Свойство 10.24.}
\vspace{0.5ex}
{\it
Пусть 
$
\bigl(\exp\omega_\nu^{}, M_\nu^{}\bigr)\in \text{\rm E}_{_{\Omega}},\ 
\nu=1,\ldots, s,\ s\leqslant m-1,\  
({\rm g}_\xi^{}, M_\xi^{})\in \text{\rm J}_{_{\Omega}},
\linebreak 
{\rm g}_\xi^{{}^{\scriptsize \gamma_\xi^{}}}\!\in C^1\Omega,
\ 
\xi=s+1,\ldots,m,\ 
\gamma_j^{}\in\R\backslash\{0\}, 
\ j=1,\ldots,m. 
$ Тогда
$
\prod\limits_{\xi=s+1}^{m}
{\rm g}_\xi^{{}^{\scriptsize \gamma_\xi^{}}}
\exp\sum\limits_{\nu=1}^{s}
\gamma_\nu^{}\;\!\omega_\nu^{}\in 
\text{\rm M}_{_{\Omega}},
$
если и только если выполняется тождество {\rm (8.1)}.
}
\vspace{0.5ex}

{\sl Следует} с учетом теоремы 8.1 из свойства 3.9 при 
$M={}-\;\!{\rm div}\;\!{\frak d}.\ \k$
\vspace{1ex}

{\bf  Свойство 10.25.}
\vspace{0.5ex}
{\it
Пусть 
$
\exp\omega_\nu^{}\in \text{\rm ME}_{_{\Omega}},\ 
\nu=1,\ldots, s,\ s\leqslant m-1,\  
({\rm g}_\xi^{}, M_\xi^{})\in \text{\rm J}_{_{\Omega}},
\linebreak 
{\rm g}_\xi^{{}^{\scriptsize \gamma_\xi^{}}}\!\in C^1\Omega,
\ 
\xi=s+1,\ldots,m,\ 
\gamma_j^{}\in\R\backslash\{0\}, 
\ j=1,\ldots,m. 
$ Тогда
\\[1.75ex]
\mbox{}\hfill
$
\displaystyle
\biggl(\,
\prod\limits_{\xi=s+1}^{m}
{\rm g}_\xi^{{}^{\scriptsize \gamma_\xi^{}}}
\exp\sum\limits_{\nu=1}^{s}
\gamma_\nu^{}\;\!\omega_\nu^{}\;\!,\;\! M\biggr)\in 
\text{\rm J}_{_{\Omega}},
\hfill
$
\\[1.75ex]
если и только если выполняется тождество}
\\[1.5ex]
\mbox{}\hfill                      
$
\displaystyle
M(t,x)=
\sum\limits_{\xi=s+1}^{m}
\gamma_\xi^{}\;\!M_{\xi}^{}(t,x)-
\sum\limits_{\nu=1}^{s}
\gamma_\nu^{}\;\!{\rm div}\;\!{\frak d}(t,x)
\quad
\forall (t,x)\in\Xi^{\;\!\prime}.
\hfill
$
\\[1.5ex]
\indent
{\sl Следует} с учетом теоремы 8.1 из свойства 3.9 при 
\vspace{1.5ex}
$M_{\nu}^{}={}-\;\!{\rm div}\;\!{\frak d},\ \nu=1,\ldots, s.\ \k$

{\bf  Свойство 10.26.}
\vspace{0.5ex}
{\it
Пусть 
$
\exp\omega_\nu^{}\in \text{\rm ME}_{_{\Omega}},\ 
\nu=1,\ldots, s,\ s\leqslant m-1,\  
({\rm g}_\xi^{}, M_\xi^{})\in \text{\rm J}_{_{\Omega}},
\linebreak 
{\rm g}_\xi^{{}^{\scriptsize \gamma_\xi^{}}}\!\in C^1\Omega,
\ 
\xi=s+1,\ldots,m,\ 
\gamma_j^{}\in\R\backslash\{0\}, 
\ j=1,\ldots,m. 
$ Тогда
$
\prod\limits_{\xi=s+1}^{m}
{\rm g}_\xi^{{}^{\scriptsize \gamma_\xi^{}}}
\exp\sum\limits_{\nu=1}^{s}
\gamma_\nu^{}\;\!\omega_\nu^{}\in 
\text{\rm M}_{_{\Omega}},
$
если и только если выполняется тождество}
\\[1.5ex]
\mbox{}\hfill                      
$
\displaystyle
\sum\limits_{\xi=s+1}^{m}
\gamma_\xi^{}\;\!M_{\xi}^{}(t,x)=
\biggl(\,\sum\limits_{\nu=1}^{s}\gamma_\nu^{}-1\biggr)
\;\!{\rm div}\;\!{\frak d}(t,x)
\quad
\forall (t,x)\in\Xi^{\;\!\prime}.
\hfill
$
\\[1.5ex]
\indent
{\sl Следует} с учетом теоремы 8.1 из свойства 10.25 при 
$M={}-\;\!{\rm div}\;\!{\frak d}.\ \k$
\vspace{1ex}

{\bf Свойство 10.27.} 
\vspace{0.5ex}
{\it
Пусть функции 
$p,q\in \text{\rm P}_{_{\!\Xi^{\;\!\prime}}}$
взаимно простые, множество
$\Omega_{_0}\subset\Xi^{\;\!\prime}$ такое, что  
$
p(t,x)\ne 0\;\;\forall (t,x)\in\Omega_{_0},\
p(t,x)= 0\;\;\forall (t,x)\in {\sf C}_{_{\!\Xi^{\;\!\prime}}}\Omega_{_0}\;\!.
$
Тогда
$
\exp\dfrac{q}{p}\in 
\text{\rm ME}_{_{\Omega_{{}_{\tiny\;\! 0}}}}, 
$
если и только если 
$(p,M)\in \text{\rm A}_{_{\!\Xi^{\;\!\prime}}}$
и выполняется тождество}
\\[1.75ex]
\mbox{}\hfill                          
$
{\frak d}\;\!q(t,x)=
q(t,x)\;\!M(t,x)-p(t,x)\;\!{\rm div}\;\!{\frak d}(t,x)
\quad
\forall (t,x)\in \Xi^{\;\!\prime}.
\hfill
$
\\[1.75ex]
\indent
{\sl Следует} с учетом теоремы 8.1 из свойства 3.10 при 
$N={}-\;\!{\rm div}\;\!{\frak d}.\ \k$
\vspace{1ex}

{\bf Свойство 10.28.} 
\vspace{0.5ex}
{\it
Пусть функции 
$p,q\in \text{\rm P}_{_{\!\Xi^{\;\!\prime}}}$
взаимно простые, множество
$\Omega_{_0}\subset\Xi^{\;\!\prime}$ такое, что  
$
p(t,x)\!\ne\! 0\;\;\forall (t,x)\!\in\!\Omega_{_0},\
p(t,x)\!=\! 0\;\;\forall (t,x)\!\in\! {\sf C}_{_{\!\Xi^{\;\!\prime}}}\Omega_{_0}\;\!.
$
\vspace{0.5ex}
Тогда
$\exp\arctg\dfrac{q}{p}\in\! 
\text{\rm ME}_{_{\Omega_{{}_{\tiny\;\! 0}}}}$
в том и только в том случае, когда
существует такая функция 
\vspace{0.5ex}
$M\in \text{\rm P}_{_{\!\Xi^{\;\!\prime}}},$
имеющая степень $\deg_{x}^{}M\leq d-1,$
что вы\-пол\-ня\-ет\-ся система тождеств}
\\[1.25ex]
\mbox{}\hfill                           
$
{\frak d}\;\!p(t,x)=
p(t,x)\;\!M(t,x)+q(t,x)\;\!{\rm div}\;\!{\frak d}(t,x)
\quad
\forall (t,x)\in \Xi^{\;\!\prime},
\hfill
$
\\[2ex]
\mbox{}\hfill                           
$
{\frak d}\;\!q(t,x)=
q(t,x)\;\!M(t,x)-p(t,x)\;\!{\rm div}\;\!{\frak d}(t,x)
\quad
\forall (t,x)\in \Xi^{\;\!\prime}.
\hfill
$
\\[1.5ex]
\indent
{\sl Следует}  из свойства 3.11 с учетом замечания 6.2 и теоремы 8.1. $\k$
\vspace{1ex}

{\bf Свойство 10.29.} 
\vspace{0.5ex}
{\it
Пусть  
$p\in \text{\rm P}_{_{\!\Xi^{\;\!\prime}}}.$
Тогда
$\exp\arctg p\in
\text{\rm ME}_{_{\Xi^{\;\!\prime}}},
$
если и только если}
\\[1.25ex]
\mbox{}\hfill                           
$
{\frak d}\;\!p(t,x)=
{}-\bigl(1+p^2(t,x)\bigr)\;\!{\rm div}\;\!{\frak d}(t,x)
\quad
\forall (t,x)\in \Xi^{\;\!\prime}.
\hfill
$
\\[1.5ex]
\indent
{\sl Следует} с учетом теоремы 8.1 из свойства 3.12 при 
$M={}-\;\!{\rm div}\;\!{\frak d}.\ \k$
\vspace{1ex}

{\bf Свойство 10.30.} 
{\it
Пусть $f\in C^1\Omega.$ Тогда
\\[1.5ex]
\mbox{}\hfill                       
$
\exp\arctg f\in \text{\rm ME}_{_{\Omega}}
\iff
\bigl(\exp\arcctg f,\, {\rm div}\;\!{\frak d}\bigr)\in
\text{\rm E}_{_{\Omega}},
\hfill
$
\\[0.5ex]
а}
\\[0.5ex]
\mbox{}\hfill                       
$
\exp\arcctg f\in \text{\rm ME}_{_{\Omega}}
\iff
\bigl(\exp\arctg f,\, -\;\!{\rm div}\;\!{\frak d}\bigr)\in
\text{\rm E}_{_{\Omega}}.
\hfill
$
\\[2.25ex]
\indent
{\bf Свойство 10.31.}\! 
\vspace{0.5ex}
{\it
Пусть  
$k\in\N,\, F\!\in \text{\rm I}_{_{\Omega}},$
функция $f\in C^1\Omega$ взаимно простая с функ\-цией $\!F,\!$ множество
$\!\Omega_{_0}\!\subset\!\Omega\!$ та\-кое, что  
\vspace{0.5ex}
$
F(t,x)\!\ne 0\;\forall (t,x)\!\in\!\Omega_{_0},\,
F(t,x)\!= 0\;\forall (t,x)\!\in\!{\sf C}_{_{\Omega}}\Omega_{_0}.\!$
Тогда
$
\exp\dfrac{f}{F^k}\in 
\text{\rm ME}_{_{\Omega_{{}_{\tiny\;\! 0}}}},
$
если и только если выполняется тождество}
\\[2ex]
\mbox{}\hfill                      
$
{\frak d}\;\!f(t,x)={}-
F^k(t,x)\;\!{\rm div}\;\!{\frak d}(t,x)
\quad
\forall (t,x)\in \Omega.
\hfill
$
\\[2ex]
\indent
{\sl Следует} с учетом теоремы 8.1 из свойства 3.14 при 
$M={}-\;\!{\rm div}\;\!{\frak d}.\ \k$

\newpage

\mbox{}
\\[-3.25ex]
\centerline{
{\bf\large \S\;\!3. Первые интегралы}}
\\[1.5ex]
\centerline{\bf  11. Первые интегралы, определяемые}
\centerline{\bf частными интегралами и последними множителями}
\\[1ex]
\indent
{\bf Теорема 11.1}
(критерий первого интеграла).
{\it 
Частный интеграл на области $\Omega$ системы {\rm (0.1)} 
является первым интегралом на области $\Omega$ системы {\rm (0.1)},
если и только если сомножитель частного интеграла 
\vspace{0.35ex}
тождественно равен нулю на области $\Omega.$
}

{\sl Доказательство}
основано на определении частного интеграла (определение 1.1)
и критерии существования первого интеграла (теорема 0.1).
При этом исходим из того, что тождество (1.1) при 
${\rm g}=F$ и $M=0$ совпадает с тождеством (0.2). $\k$
\vspace{0.5ex}

Таким образом, имеет место включение
$\text{\rm I}_{_{\Omega}}\subset \text{\rm J}_{_{\Omega}},$ 
а утверждение теоремы 11.1 выражается эквиваленцией 
\\[0.5ex]
\mbox{}\hfill
$
F\in \text{\rm I}_{_{\Omega}}
\iff
(F,\;\! 0)\in \text{\rm J}_{_{\Omega}}.
\hfill
$
\\[1.75ex]
\indent
{\bf Свойство 11.1.}
{\it 
Справедливо утверждение}
\\[1.35ex]
\mbox{}\hfill
$
\forall\;\! C\in\R
\ \ 
\Longrightarrow
\ \ 
C\in\text{\rm I}_{_{\Xi}}.
\hfill
$
\\[1.5ex]
\indent
{\sl Действительно},
так как ${\frak d}\;\!C=0\;\;\forall (t,x)\in\Xi,$
то, по теореме 0.2, $C\in\text{\rm I}_{_{\Xi}}.\ \k$
\vspace{1ex}

{\bf Свойство 11.2.}
{\it 
Если $C\in\R,$ то}
\\[1.5ex]
\mbox{}\hfill
$
F\in \text{\rm I}_{_{\Omega}}
\iff
F+C\in \text{\rm I}_{_{\Omega}}.
\hfill
$
\\[1.5ex]
\indent
{\sl Следует}
из свойства 0.1 с учетом свойства 11.1. $\k$
\vspace{1ex}

{\bf Теорема 11.2}
(геометрический смысл первого интеграла).
{\it 
Первый интеграл $F$ на области $\Omega$ системы {\rm (0.1)} 
определяет семейство ее интегральных многообразий 
$F(t,x)=C,$ где $C$ --- произвольная постоянная из множества 
значений функции $F$ {\rm(}случаи
$C={}-\infty,\ C={}+\infty$ и $C=\infty$ не исключаются{\rm)}. 
}
\vspace{0.35ex}

{\sl Следует}
из определения интегрального многообразия (определение 0.3),
критерия существования первого интеграла (теорема 0.1)
и свойства 11.2.  $\k$
\vspace{1ex}

{\bf Свойство 11.3.}
{\it 
Линейная комбинация над полем $\R$ первых интегралов на области $\Omega$
системы {\rm (0.1)} является первым интегралом на области $\Omega$
системы {\rm (0.1)}}:
\\[1ex]
\mbox{}\hfill
$
\displaystyle
\lambda_j^{}\in\R,
\ \ 
F_j^{}\in \text{\rm I}_{_{\Omega}},
\ j=1,\ldots,m,
\ \ \Longrightarrow\ \ 
\sum\limits_{j=1}^{m}
\lambda_j^{}\;\!F_j^{}\in \text{\rm I}_{_{\Omega}}.
\hfill
$
\\[1ex]
\indent
{\sl Следует}
из свойства 0.1. $\k$
\vspace{1ex}

Свойства 11.3 и 11.1 определяют математическую структуру множества 
первых интегралов системы (0.1).
\vspace{0.5ex}

{\bf Теорема 11.3.}
{\it 
Множество первых интегралов на области $\Omega$
системы {\rm (0.1)} является линейным пространством над 
полем действительных чисел}:
\vspace{1ex}
$\bigl(\text{\rm I}_{_{\Omega}},\,\R,\,+\;\!,\,\cdot\;\!,\,=\bigr).$

{\bf Свойство 11.4.} 
\vspace{0.25ex}
{\it
Пусть множество $\Omega_{_0}\subset\Omega$ такое, что  
$F(t,x)\ne 0\;\;\forall (t,x)\in \Omega_{_0},$ 
$F(t,x)=0\;\;\forall (t,x)\in {\sf C}_{_\Omega}\Omega_{_0}.$
Тогда}
\\[1ex]
\mbox{}\hfill  
$
F\in \text{I}_{_{\Omega}}
\iff 
|F|\in \text{I}_{_{\Omega_{{}_{\tiny\;\! 0}}}}.
\hfill
$
\\[1.5ex]
\indent
{\sl Следует}
из свойства 1.3 с учетом теоремы 11.1. $\k$
\vspace{1ex}

{\bf Свойство 11.5.} 
{\it
Если $F\in \text{\rm I}_{_{\Omega}},$ а 
\\[1.5ex]
\mbox{}\hfill                           
$
{\frak d}\;\!f(t,x)=
\bigl(f(t,x)+F(t,x)\bigr)\;\!M(t,x)
\quad
\forall (t,x)\in \Omega,
\hfill
$
\\[1.5ex]
функция $M\in \text{\rm P}_{_{\!\Xi^{\;\!\prime}}}$
и имеет степень $\deg_{x}^{}M\leq d-1,$ то
$(f+F, M)\in \text{\rm J}_{_{\Omega}}.$
}
\vspace{0.5ex}

{\sl Доказательство.} 
С учетом теоремы 0.2 
\\[1.25ex]
\mbox{}\hfill                           
$
{\frak d}\;\!\bigl(f(t,x)+F(t,x)\bigr)=
{\frak d}\;\!f(t,x)=
\bigl(f(t,x)+F(t,x)\bigr)\;\!M(t,x)
\quad
\forall (t,x)\in \Omega.
\hfill
$
\\[1ex]
\indent
По теореме 1.1, $(f+F, M)\in \text{\rm J}_{_{\Omega}}.\ \k$
\vspace{1.25ex}

{\bf Свойство 11.6.} 
{\it
Пусть $F\in \text{\rm I}_{_{\Omega}}.$ Тогда} 
\\[1.5ex]
\mbox{}\hfill                           
$
({\rm g}\;\!F,\;\! M)\in \text{\rm J}_{_{\Omega}}
\iff
({\rm g},\;\! M)\in \text{\rm J}_{_{\Omega}}.
\hfill
$
\\[1.5ex]
\indent
{\sl Доказательство.} 
С учетом тождества (0.3) 
\\[1.5ex]
\mbox{}\hfill                           
$
{\frak d}\;\!\bigl({\rm g}(t,x)\;\!F(t,x)\bigr)=
F(t,x)\, {\frak d}\;\!{\rm g}(t,x)
\quad
\forall (t,x)\in \Omega.
\hfill
$
\\[1.5ex]
\indent
Из этого тождества, применив теорему 1.1 к функциям 
${\rm g}\;\!F$ и ${\rm g},$
получаем утверждение доказываемого свойства. $\k$
\vspace{0.75ex}

{\bf Свойство 11.7.} 
{\it
Пусть $({\rm g},\;\! M)\in \text{\rm J}_{_{\Omega}}.$ Тогда} 
\\[1.5ex]
\mbox{}\hfill                           
$
({\rm g}\;\!F,\;\! M)\in \text{\rm J}_{_{\Omega}}
\iff
F\in \text{\rm I}_{_{\Omega}}.
\hfill
$
\\[1.5ex]
\indent
{\sl Доказательство.} 
С учетом тождества (1.2) из теоремы 1.1 производная в силу дифференциальной системы (0.1)
равна
\\[1.5ex]
\mbox{}\hfill                           
$
{\frak d}\;\!\bigl({\rm g}(t,x)\;\!F(t,x)\bigr)=
{\rm g}(t,x)\;\!
\bigl(F(t,x)\;\!M(t,x)+ {\frak d}\;\!F(t,x)\bigr)
\quad
\forall (t,x)\in \Omega.
\hfill
$
\\[1.5ex]
\indent
Если $({\rm g}\;\!F,\;\! M)\in \text{\rm J}_{_{\Omega}},$ то
\\[1.5ex]
\mbox{}\hfill                           
$
F(t,x)\;\!{\rm g}(t,x)\;\!M(t,x)=
{\rm g}(t,x)\;\!
\bigl(F(t,x)\;\!M(t,x)+ {\frak d}\;\!F(t,x)\bigr)
\quad
\forall (t,x)\in \Omega.
\hfill
$
\\[1.5ex]
Отсюда ${\frak d}\;\!F(t,x)=0\;\;\forall (t,x)\in \Omega,$
а значит, по теореме 0.2, $F\in \text{\rm I}_{_{\Omega}}.$
\vspace{0.5ex}

Если $F\in \text{\rm I}_{_{\Omega}},$ то 
\\[1.5ex]
\mbox{}\hfill                           
$
{\frak d}\;\!\bigl({\rm g}(t,x)\;\!F(t,x)\bigr)=
{\rm g}(t,x)\;\!F(t,x)\;\!M(t,x)
\quad
\forall (t,x)\in \Omega.
\hfill
$
\\[1.5ex]
\indent
По теореме 1.1, $({\rm g}\;\!F,\;\! M)\in \text{\rm J}_{_{\Omega}}.\ \k$
\vspace{1ex}

{\bf Следствие 11.1.} 
{\it
Справедливы утверждения}: 
\\[1.5ex]
\mbox{}\hfill                           
$
F\in \text{\rm I}_{_{\Omega}}
\ \ \Longrightarrow\ \ 
\bigl(\mu\;\!F\in \text{\rm M}_{_{\Omega}}
\iff 
\mu\in \text{\rm M}_{_{\Omega}}\bigr);
\hfill                           
$
\\[2ex]
\mbox{}\hfill                           
$
\mu\in \text{\rm M}_{_{\Omega}}
\ \ \Longrightarrow\ \ 
\bigl(\mu\;\!F\in \text{\rm M}_{_{\Omega}}
\iff 
F\in \text{\rm I}_{_{\Omega}}\bigr).
\hfill
$
\\[2ex]
\indent
{\bf Свойство 11.8.} 
{\it
Пусть $({\rm g},\;\!\varphi)\in \text{\rm J}_{_{\Omega}},\ 
\varphi\in C^1T^{\;\!\prime}.$ Тогда функция 
\\[1.5ex]
\mbox{}\hfill                           
$
\displaystyle
F\colon (t,x)\to\ 
{\rm g}(t,x)\exp\biggl({}-\int\limits_{t_0^{}}^{t}\varphi(\tau)\,d\tau\biggr)
\quad
\forall (t,x)\in \Omega,
\hfill
$
\\[1.5ex]
где $t_{_0}^{}$ --- произвольная фиксированная точка из области $T^{\;\!\prime},$
\vspace{0.25ex}
будет первым интегралом на области $\Omega$ системы {\rm (0.1)}.
}
\vspace{0.5ex}

{\sl Доказательство.} 
Производная в силу системы (0.1)
\\[1.5ex]
\mbox{}\hfill                           
$
\displaystyle
{\frak d}\;\!F(t,x)=
{\frak d}\;\!
\biggl({\rm g}(t,x)\exp\biggl({}-\int\limits_{t_0^{}}^{t}\varphi(\tau)\,d\tau\biggr)\biggr)=
\hfill                           
$
\\[1.5ex]
\mbox{}\hfill                           
$
\displaystyle
=\exp\biggl({}-\int\limits_{t_0^{}}^{t}\varphi(\tau)\,d\tau\biggr)\;\!
{\frak d}\;\!{\rm g}(t,x)\;\!+\;\!
{\rm g}(t,x)\, {\sf D}\exp\biggl({}-\int\limits_{t_0^{}}^{t}\varphi(\tau)\,d\tau\biggr)=
\hfill                           
$
\\[1.5ex]
\mbox{}\hfill                           
$
\displaystyle
=\varphi(t)\;\! {\rm g}(t,x)\exp\biggl({}-\int\limits_{t_0^{}}^{t}\varphi(\tau)\,d\tau\biggr)\;\!-\;\!
\varphi(t)\;\! {\rm g}(t,x)\exp\biggl({}-\int\limits_{t_0^{}}^{t}\varphi(\tau)\,d\tau\biggr)=0
\quad
\forall (t,x)\in \Omega.
\hfill
$
\\[1.75ex]
\indent
По теореме 0.2, $F\in \text{\rm I}_{_{\Omega}}.\ \k$
\vspace{1.5ex}

{\bf Следствие 11.2.} 
{\it
Если $({\rm g},\;\!\lambda)\in \text{\rm J}_{_{\Omega}},\ 
\lambda \in\R\backslash\{0\},$ то 
${\rm g}\;\!e^{{}-\lambda\;\!t}\in \text{\rm I}_{_{\Omega}}.$
} 
\vspace{1.5ex}

{\bf  Свойство 11.9.}
\vspace{0.5ex}
{\it
Пусть $\lambda_j^{}\in\R, \ j=1,\ldots,m,\ m\leq n.$ Тогда
$\sum\limits_{j=1}^{m}\lambda_j^{}\;\!x_j^{}\in \text{\rm I}_{_{\Xi}}$
в том и только в том случае, когда}
\\[1.5ex]
\mbox{}\hfill
$
\displaystyle
\sum\limits_{j=1}^{m}
\lambda_j^{}\;\!X_j^{}(t,x)=0
\quad
\forall (t,x)\in \Xi.
\hfill
$
\\[1.5ex]
\indent
{\sl Следует}
из свойства 2.2 с учетом теоремы 11.1. $\k$
\vspace{0.75ex}

{\bf Свойство 11.10.} 
{\it
Пусть $({\rm g}_j^{},\;\!M_j^{})\in \text{\rm J}_{_{\Omega}},\ 
\lambda_j^{}\in\R\backslash\{0\}, \ j=1,\ldots,m.$ Тогда
$\sum\limits_{j=1}^{m}\lambda_j^{}\;\!{\rm g}_j^{}\in \text{\rm I}_{_{\Omega}}$
в том и только в том случае, когда}
\\[1.5ex]
\mbox{}\hfill
$
\displaystyle
\sum\limits_{j=1}^{m}
\lambda_j^{}\;\!{\rm g}_j^{}(t,x)\;\!M_j^{}(t,x)=0
\quad
\forall (t,x)\in \Omega.
\hfill
$
\\[1.5ex]
\indent
{\sl Доказательство}
основано на теореме 11.1 и состоит в том, что
\\[1.5ex]
\mbox{}\hfill
$
\displaystyle
{\frak d}\;\!\sum\limits_{j=1}^{m}
\lambda_j^{}\;\!{\rm g}_j^{}(t,x)=
\sum\limits_{j=1}^{m}
\lambda_j^{}\;\!{\rm g}_j^{}(t,x)\;\!M_j^{}(t,x)
\quad
\forall (t,x)\in \Omega. \ \k
\hfill
$
\\[2ex]
\indent
{\bf Следствие 11.3.} 
{\it
Пусть $({\rm g}_j^{},\;\!\rho_j^{}\;\!M_{_0}^{})\in \text{\rm J}_{_{\Omega}},\ 
\rho_j^{},\lambda_j^{}\in\R\backslash\{0\}, \ j=1,\ldots,m.$ Тогда}
\\[1.5ex]
\mbox{}\hfill
$
\displaystyle
\sum\limits_{j=1}^{m}\lambda_j^{}\;\!{\rm g}_j^{}\in \text{\rm I}_{_{\Omega}}
\iff
\sum\limits_{j=1}^{m}
\rho_j^{}\;\!\lambda_j^{}=0.
\hfill
$
\\[2ex]
\indent
{\bf Следствие 11.4.} 
{\it
Пусть $\mu_j^{}\in \text{\rm M}_{_{\Omega}},\ 
\lambda_j^{}\in\R\backslash\{0\}, \ j=1,\ldots,m.$ Тогда}
\\[1.5ex]
\mbox{}\hfill
$
\displaystyle
\sum\limits_{j=1}^{m}\lambda_j^{}\;\!\mu_j^{}\in \text{\rm I}_{_{\Omega}}
\iff
\sum\limits_{j=1}^{m}
\lambda_j^{}=0.
\hfill
$
\\[2ex]
\indent
{\bf Теорема 11.4.}
{\it
Непрерывно дифференцируемая функция {\rm (1.4)}
являет\-ся первым интегралом  на области $\Omega$ системы {\rm (0.1)}
\vspace{0.25ex}
тогда и только тогда, когда существуют такие функции 
$M_j^{}\in C^1\Omega,\ j=1,\ldots,m,$ 
что выполняются тождества {\rm (1.5)} и} 
\\[1.5ex]
\mbox{}\hfill
$
\displaystyle
\sum\limits_{j=1}^{m} M_j^{}(t,x)= 0
\quad
\forall (t,x)\in\Omega.
\hfill
$
\\[1.5ex]
\indent
{\sl Следует}
из теоремы 1.3 с учетом теоремы 11.1. $\k$
\vspace{1ex}

{\bf Свойство 11.11.} 
\vspace{0.5ex}
{\it
Пусть 
$({\rm g}_j^{}, M_{j}^{})\in \text{\rm J}_{_{\Omega}},\ 
\gamma_j^{}\in\R\backslash\{0\},\
{\rm g}_j^{{}^{\scriptsize \gamma_{j}^{}}}\in C^1\Omega, \ j=1,\ldots,m.$ 
Тогда  
$\prod\limits_{j=1}^{m} 
{\rm g}_j^{{}^{\scriptsize \gamma_{j}^{}}} 
\in \text{\rm I}_{_{\Omega}}\;\!,$
если и только если линейная комбинация сомножителей}
\\[1ex]
\mbox{}\hfill                      % (11.1)
$
\displaystyle
\sum\limits_{j=1}^{m}
\gamma_j^{}\;\!M_{j}^{}(t,x)=0
\quad
\forall (t,x)\in\Xi^{\;\!\prime}.
$
\hfill (11.1)
\\[1.5ex]
\indent
{\sl Следует}
с учетом теоремы 11.1 из свойства 1.9 при $M=0.\ \k$
\vspace{0.5ex}

Свойство 11.11 также следует из теоремы 11.4, если учесть свойство 1.6. $\k$
\vspace{1ex}

{\bf Следствие 11.5.} 
\vspace{0.5ex}
{\it
Пусть 
$({\rm g}_j^{}, \rho_{j}^{}\;\!M_{_0}^{})\in \text{\rm J}_{_{\Omega}},\ 
\rho_{j}^{},\gamma_j^{}\in\R\backslash\{0\},\
{\rm g}_j^{{}^{\scriptsize \gamma_{j}^{}}}\in C^1\Omega, \ j=1,\ldots,m.$ 
Тогда  
$\prod\limits_{j=1}^{m} 
{\rm g}_j^{{}^{\scriptsize \gamma_{j}^{}}} 
\in \text{\rm I}_{_{\Omega}}\;\!,$
если и только если}
$
\sum\limits_{j=1}^{m}
\rho_{j}^{}\;\!\gamma_j^{}=0.
$
\vspace{1ex}

Обратим внимание, что при выполнении условий следствия 11.5 попарно взятые 
частные интегралы образуют первые интегралы
\\[1.5ex]
\mbox{}\hfill                   
$
{\rm g}_\xi^{{}^{\scriptsize \gamma_{\xi}^{}}}\;\!
{\rm g}_\zeta^{{}^{\scriptsize \gamma_{\zeta}^{}}} 
\in \text{\rm I}_{_{\Omega}}
\iff
\rho_{\xi}^{}\;\!\gamma_\xi^{}+
\rho_{\zeta}^{}\;\!\gamma_\zeta^{}=0,
\ \ \xi,\zeta=1,\ldots, m,\ \ \xi\ne\zeta.
\hfill
$
\\[2.5ex]
\indent
{\bf Пример 11.1.}
Дифференциальная система
\\[2ex]
\mbox{}\hfill
$
\dfrac{dx}{dt}=1,
\qquad
\dfrac{dy}{dt}={}-2\;\!xy+z^2,
\qquad
\dfrac{dz}{dt}={}-2\;\!xz
\hfill
$
\\[2ex]
имеет автономные частные интегралы: полиномиальный
\\[1.5ex]
\mbox{}\hfill
$
{\rm g}_1^{}\colon (x,y,z)\to\ z
\quad
\forall (x,y,z)\in\R^3
\hfill
$
\\[2ex]
с сомножителем 
$
M_1^{}\colon (x,y,z)\to {}-2\;\!x
\;\;\forall (x,y,z)\in\R^3;
$
\\[2ex]
\mbox{}\hfill
$
\displaystyle
{\rm g}_2^{}\colon (x,y,z)\to\ 
y-z^2e^{\;\!x^2}\int\limits_{0}^{x}e^{\;\!-\;\!\tau^2}\;\!d\tau
\quad
\forall (x,y,z)\in\R^3
\hfill
$
\\[1.5ex]
с сомножителем 
$
M_2^{}\colon (x,y,z)\to {}-2\;\!x
\;\;\forall (x,y,z)\in\R^3;
$
условный
\\[2ex]
\mbox{}\hfill
$
\displaystyle
{\rm g}_3^{}\colon (x,y,z)\to\ 
e^{\;\!x^2}
\quad
\forall (x,y,z)\in\R^3
\hfill
$
\\[2ex]
с сомножителем 
$
M_3^{}\colon (x,y,z)\to\, 2\;\!x
\;\;\forall (x,y,z)\in\R^3.
$
\vspace{1ex}

По следствию 11.5 
\vspace{0.35ex}
(при $\rho_1^{}={}-1,\ \rho_2^{}=1,\ M_{_0}^{}=2\;\!x,\ \gamma_1^{}=\gamma_2^{}=1)$
на основании частных интегралов ${\rm g}_1^{}$ и ${\rm g}_3^{},$ 
 ${\rm g}_2^{}$ и ${\rm g}_3^{}$ строим автономные первые интегралы:
 \\[2ex]
\mbox{}\hfill
$
\displaystyle
F_1^{}={\rm g}_1^{}\;\!{\rm g}_3^{}\colon (x,y,z)\to\ 
z\;\!e^{\;\!x^2}
\quad
\forall (x,y,z)\in\R^3;
\hfill
$
\\[2.5ex]
\mbox{}\hfill
$
\displaystyle
F_2^{}={\rm g}_2^{}\;\!{\rm g}_3^{}\colon (x,y,z)\to\ 
\Bigl(y-z^2e^{\;\!x^2}\int\limits_{0}^{x}e^{\;\!-\;\!\tau^2}\;\!d\tau\Bigr)\;\!
e^{\;\!x^2}
\quad
\forall (x,y,z)\in\R^3.
\hfill
$
\\[1.5ex]
\indent
Функции $F_1^{}$ и $F_2^{},$
\vspace{0.5ex}
будучи функционально независимыми, образуют автономный интегральный базис на $\R^3.$
\vspace{0.5ex}

Функция 
\\[1ex]
\mbox{}\hfill
$
\displaystyle
F_3^{}\colon (t,x,y,z)\to\ x-t
\quad
\forall (t,x,y,z)\in\R^4
\hfill
$
\\[1.5ex]
является неавтономным первым интегралом.
\vspace{0.5ex}

Совокупность $\{F_1^{},\;\! F_2^{},\;\! F_3^{}\}$ --- 
базис первых интегралов на $\R^4.$
\vspace{1.5ex}

{\bf Следствие 11.6.} 
\vspace{0.5ex}
{\it
Пусть 
$\mu_j^{}\in \text{\rm M}_{_{\Omega}},\ 
\gamma_j^{}\in\R\backslash\{0\},\
\mu_j^{{}^{\scriptsize \gamma_{j}^{}}}\in C^1\Omega, \ j=1,\ldots,m.$ 
Тогда  фун\-к\-ция
$\prod\limits_{j=1}^{m} 
\mu_j^{{}^{\scriptsize \gamma_{j}^{}}} 
\in \text{\rm I}_{_{\Omega}}\;\!,$
если и только если}
$
\sum\limits_{j=1}^{m}
\gamma_j^{}=0.
$
\vspace{1ex}

{\bf Следствие 11.7.}
 \vspace{0.35ex}
 {\it 
Пусть множество $\Omega_{_0}\subset \Omega$ такое, что 
${\rm g}_2^{}(t,x)\ne 0\;\; \forall (t,x)\in \Omega_{_0},$
${\rm g}_2^{}(t,x)= 0\;\; \forall (t,x)\in {\sf C}_{_\Omega}\Omega_{_0}.$
Тогда} 
\\[1.5ex]
\mbox{}\hfill
$
({\rm g}_1^{}, M)\in \text{\rm J}_{_{\Omega}}
\ \ \&\ \ 
({\rm g}_2^{}, M)\in \text{\rm J}_{_{\Omega}}
\ \ \Longrightarrow\ \ 
\dfrac{{\rm g}_1^{}}{{\rm g}_2^{}}\in
\text{\rm I}_{_{ \Omega_{_0}}}.
\hfill
$
\\[2ex]
\indent
Заметим, что при $M={}-\;\!{\rm div}\;\!{\frak d}$
следствие 11.7 есть свойство Якоби последних мно\-жи\-те\-лей (свойство 0.2).
\vspace{1ex}

{\bf Свойство 11.12.} 
\vspace{0.5ex}
{\it
Пусть 
$({\rm g}_\tau^{}, M_{\tau}^{})\in \text{\rm J}_{_{\Omega}},\ 
\tau=1,\ldots, m-1,\
\gamma_j^{}\in\R\backslash\{0\},\
{\rm g}_j^{{}^{\scriptsize \gamma_{j}^{}}}\in C^1\Omega, 
\linebreak 
j=1,\ldots,m.$ 
Тогда}
\\[1.25ex]
\mbox{}\hfill                      
$
\displaystyle
\prod\limits_{j=1}^{m} 
{\rm g}_j^{{}^{\scriptsize \gamma_{j}^{}}}
\in \text{\rm I}_{_{\Omega}}
\iff
\biggl({\rm g}_m^{},
{}-\dfrac{1}{\gamma_m^{}}\;\!
\sum\limits_{\tau=1}^{m-1}
\gamma_{\tau}\;\!M_{\tau}^{}\biggl) 
\in \text{\rm J}_{_{\Omega}}\;\!.
\hfill
$
\\[1.5ex]
\indent
{\sl Следует} 
из свойства 1.10 с учетом теоремы 11.1. $\k$
\vspace{1ex}

\newpage

{\bf Следствие 11.8.} 
\vspace{0.5ex}
{\it
Пусть 
$\mu_j^{}\in \text{\rm M}_{_{\Omega}},\ 
\gamma,\gamma_j^{}\in\R\backslash\{0\},\ 
{\rm g}^{{}^{\scriptsize \gamma}},
\mu_j^{{}^{\scriptsize \gamma_{j}^{}}}\in C^1\Omega, 
\ j=1,\ldots,m.$ 
Тогда}
\\[1.25ex]
\mbox{}\hfill                      
$
\displaystyle
{\rm g}^{{}^{\scriptsize \gamma}}\;\!
\prod\limits_{j=1}^{m} 
\mu_j^{{}^{\scriptsize \gamma_{j}^{}}}
\in \text{\rm I}_{_{\Omega}}
\iff
\biggl({\rm g},\
\dfrac{1}{\gamma}\;\!
\sum\limits_{j=1}^{m}
\gamma_{j}^{}\,{\rm div}\;\!{\frak d}\biggl) 
\in \text{\rm J}_{_{\Omega}}\;\!.
\hfill
$
\\[1.5ex]
\indent
{\bf Свойство 11.13.} 
\vspace{0.5ex}
{\it
Пусть 
$({\rm g}_\nu^{}, M_{\nu}^{})\!\in \text{\rm J}_{_{\Omega}},\, 
\nu\!=\!1,\ldots, s,\ s\leq m\!-\!2, \ 
\gamma_j^{}\!\in\R\backslash\{0\},\
{\rm g}_j^{{}^{\scriptsize \gamma_{j}^{}}}\!\in C^1\Omega, 
\linebreak 
j=1,\ldots,m.$ 
Тогда
\\[1.5ex]
\mbox{}\hfill                      
$
\displaystyle
\prod\limits_{j=1}^{m} 
{\rm g}_j^{{}^{\scriptsize \gamma_{j}^{}}} 
\in \text{\rm I}_{_{\Omega}}
\iff
\biggl(
\,\prod\limits_{k=s+1}^{m} 
{\rm g}_k^{{}^{\scriptsize \gamma_{k}^{}}}, 
{}-\sum\limits_{\nu=1}^{s}
\gamma_{\nu}\;\!M_{\nu}^{}\biggl) 
\in \text{\rm J}_{_{\Omega}}\;\!.
\hfill
$
\\[1.5ex]
Кроме этого существуют такие функции 
$M_k^{}\in C^1\Omega,\ k=s+1,\ldots, m,$
что выполняются тождества {\rm (1.14)} и}
\\[1.5ex]
\mbox{}\hfill                      % (11.2)
$
\displaystyle
\sum\limits_{k=s+1}^{m}
\gamma_{k}^{}\;\!M_k^{}(t,x)=
{}-\sum\limits_{\nu=1}^{s}
\gamma_{\nu}^{}\;\!M_\nu^{}(t,x)
\quad
\forall (t,x)\in\Omega.
$
\hfill (11.2)
\\[2ex]
\indent
{\sl Следует} 
из свойства 1.11 с учетом теоремы 11.1. $\k$
\vspace{1ex}

{\bf Следствие 11.9.} 
\vspace{0.5ex}
{\it
Пусть 
$\mu_\nu^{}\in \text{\rm M}_{_{\Omega}},\
\gamma_\nu^{}\in\R\backslash\{0\},\ 
\mu_\nu^{{}^{\scriptsize \gamma_{\nu}^{}}}\in C^1\Omega, \
\nu=1,\ldots, s,\ s\leq m-2,$
${\rm g}_k^{{}^{\scriptsize \gamma_{k}^{}}}\!\in C^1\Omega,\
\gamma_k^{}\in\R\backslash\{0\},\
k=s+1,\ldots,m.$ 
Тогда
\\[1.75ex]
\mbox{}\hfill                      
$
\displaystyle
\prod\limits_{\nu=1}^{s} 
\mu_\nu^{{}^{\scriptsize \gamma_{\nu}^{}}}\, 
\prod\limits_{k=s+1}^{m} 
{\rm g}_k^{{}^{\scriptsize \gamma_{k}^{}}} 
\in \text{\rm I}_{_{\Omega}}
\iff
\biggl(
\,\prod\limits_{k=s+1}^{m} 
{\rm g}_k^{{}^{\scriptsize \gamma_{k}^{}}}, \
\sum\limits_{\nu=1}^{s}
\gamma_{\nu}\,{\rm div}\;\!{\frak d}\biggl) 
\in \text{\rm J}_{_{\Omega}}\;\!.
\hfill
$
\\[1.5ex]
Кроме этого существуют такие функции 
$M_k^{}\in C^1\Omega,\ k=s+1,\ldots, m,$
что выполняются тождества {\rm (1.14)} и}
\\[1.5ex]
\mbox{}\hfill                      % (11.3)
$
\displaystyle
\sum\limits_{k=s+1}^{m}
\gamma_{k}^{}\;\!M_k^{}(t,x)=
\sum\limits_{\nu=1}^{s}
\gamma_{\nu}^{}\,{\rm div}\;\!{\frak d}(t,x)
\quad
\forall (t,x)\in\Omega.
$
\hfill (11.3)
\\[2ex]
\indent
{\bf Свойство 11.14.} 
\vspace{0.5ex}
{\it
Пусть 
$p_j^{}\in \text{\rm P}_{_{\!\Xi^{\;\!\prime}}},\ 
\gamma_j^{}\in\R\backslash\{0\},$ 
множество $\Omega_{_0}\subset\Xi^{\;\!\prime}\!$ такое, что  
$p_j^{{}^{\scriptsize \gamma_j^{}}}\in C^1\Omega_{_0},\ j=1,\ldots, m.$
Тогда
\\[1.5ex]
\mbox{}\hfill  
$
\displaystyle
\prod\limits_{j=1}^{m} 
p_j^{{}^{\scriptsize \gamma_{j}^{}}}
\in \text{\rm I}_{_{\Omega_{{}_{\tiny\;\! 0}}}}
\iff
\bigl(p_j^{}, M_j^{}\bigr)\in \text{\rm A}_{_{\Xi^{\;\!\prime}}},
\ \ 
j=1,\ldots, m.
\hfill
$
\\[1.5ex]
Кроме этого сомножители $M_1^{},\ldots, M_m^{}$
\vspace{0.75ex}
такие, что выполняется тождество {\rm (11.1)}.}

{\sl Следует} 
из теоремы 2.2 с учетом теоремы 11.1. $\k$
\vspace{1ex}

{\bf  Свойство 11.15.}
{\it
Если 
$
\bigl(\exp\omega_j^{}, M_j^{}\bigr)\in \text{\rm E}_{_{\Omega}},\ 
\gamma,\gamma_j^{}\in\R\backslash\{0\}, \ j=1,\ldots,m, \
{\rm g}^{\gamma}\in C^1\Omega,$ то}
\\[1.5ex]
\mbox{}\hfill
$
\displaystyle
{\rm g}^{\gamma}\exp\sum\limits_{j=1}^{m}\gamma_j^{}\;\!\omega_j^{}\in 
\text{\rm I}_{_{\Omega}}
\iff
\biggl(\,
{\rm g},{}- \dfrac{1}{\gamma}\,
\sum\limits_{j=1}^{m}\gamma_j^{}\;\!M_j^{}\biggr)\in 
\text{\rm J}_{_{\Omega}}.
\hfill
$
\\[1.5ex]
\indent
{\sl Следует} 
из свойства 3.7 с учетом теоремы 11.1. $\k$
\vspace{1ex}

{\bf  Следствие 11.10.}
{\it
Если
$
\exp\omega_j^{}\in \text{\rm ME}_{_{\Omega}},\ 
\gamma,\gamma_j^{}\in\R\backslash\{0\}, \ j=1,\ldots,m, \
{\rm g}^{\gamma}\in C^1\Omega,$ то}
\\[1.5ex]
\mbox{}\hfill
$
\displaystyle
{\rm g}^{\gamma}\exp\sum\limits_{j=1}^{m}\gamma_j^{}\;\!\omega_j^{}\in 
\text{\rm I}_{_{\Omega}}
\iff
\biggl(\,
{\rm g},\ \dfrac{1}{\gamma}\,
\sum\limits_{j=1}^{m}\gamma_j^{}\,{\rm div}\;\!{\frak d}\biggr)\in 
\text{\rm J}_{_{\Omega}}.
\hfill
$
\\[2ex]
\indent
{\bf  Свойство\! 11.16.}\!
{\it
Если
$\!\bigl(\exp\omega_\tau^{}, M_\tau^{}\bigr)\!\in\! \text{\rm E}_{_{\Omega}},\,
\tau\!=\!1,\ldots, m\!-\!1,\ 
\gamma_j^{}\in\R\backslash\{0\}, \, j\!=\!1,\ldots,m,$ то}
\\[1.5ex]
\mbox{}\hfill
$
\displaystyle
\sum\limits_{j=1}^{m}\gamma_j^{}\;\!\omega_j^{}\in 
\text{\rm I}_{_{\Omega}}
\iff
\biggl(\,
\exp\omega_{m}^{}\;\!,{}-
 \dfrac{1}{\gamma_{m}^{}}\,
\sum\limits_{\tau=1}^{m-1}\gamma_\tau^{}\;\!M_\tau^{}\biggr)\in 
\text{\rm E}_{_{\Omega}}.
\hfill
$
\\[2ex]
\indent
{\sl Следует} 
из следствия 3.3 с учетом теоремы 11.1. $\k$
\vspace{1ex}

{\bf  Следствие\! 11.11.}\!
{\it
Если
$\exp\omega_\tau^{}\in \text{\rm ME}_{_{\Omega}},\
\tau\!=\!1,\ldots, m\!-\!1,\ 
\gamma_j^{}\in\R\backslash\{0\}, \ j=1,\ldots,m,$ то}
\\[1.5ex]
\mbox{}\hfill
$
\displaystyle
\sum\limits_{j=1}^{m}\gamma_j^{}\;\!\omega_j^{}\in 
\text{\rm I}_{_{\Omega}}
\iff
\biggl(\,
\exp\omega_{m}^{}\;\!,\
 \dfrac{1}{\gamma_{m}^{}}\,
\sum\limits_{\tau=1}^{m-1}\gamma_\tau^{}\,{\rm div}\;\!{\frak d}\biggr)\in 
\text{\rm E}_{_{\Omega}}.
\hfill
$
\\[2ex]
\indent
{\bf  Свойство 11.17.}
\vspace{0.5ex}
{\it
Пусть 
$
\bigl(\exp\omega_\nu^{}, M_\nu^{}\bigr)\in \text{\rm E}_{_{\Omega}},\ 
\nu=1,\ldots, s,\ s\leqslant m-2,\  
\gamma_j^{}\in\R\backslash\{0\}, 
\linebreak 
j=1,\ldots,m,\ 
{\rm g}_k^{{}^{\scriptsize \gamma_k^{}}}\!\in C^1\Omega,\ 
k=s+1,\ldots,m.$ Тогда
\\[1.75ex]
\mbox{}\hfill
$
\displaystyle
\prod\limits_{k=s+1}^{m}
{\rm g}_k^{{}^{\scriptsize \gamma_k^{}}}
\exp\sum\limits_{\nu=1}^{s}
\gamma_\nu^{}\;\!\omega_\nu^{}\in 
\text{\rm I}_{_{\Omega}}
\iff
\biggl(\,
\prod\limits_{k=s+1}^{m}
{\rm g}_k^{{}^{\scriptsize \gamma_k^{}}},{}-
\sum\limits_{\nu=1}^{s}\gamma_\nu^{}\;\!M_\nu^{}\biggr)\in 
\text{\rm J}_{_{\Omega}}.
\hfill
$
\\[1.75ex]
Кроме этого существуют такие функции $M_k^{}\in C^1\Omega,\ k=s+1,\ldots, m,$
что выполняются тождества {\rm (1.14)} и {\rm (11.2)}.
}
\vspace{0.5ex}

{\sl Следует} 
из свойства 3.8 с учетом теоремы 11.1. $\k$
\vspace{1ex}

{\bf  Следствие 11.12.}
\vspace{0.5ex}
{\it
Пусть 
$
\exp\omega_\nu^{}\in \text{\rm ME}_{_{\Omega}},\ 
\nu=1,\ldots, s,\ s\leqslant m-2,\  
\gamma_j^{}\in\R\backslash\{0\}, 
\linebreak 
j=1,\ldots,m,\ 
{\rm g}_k^{{}^{\scriptsize \gamma_k^{}}}\!\in C^1\Omega,\ 
k=s+1,\ldots,m.$ Тогда
\\[1.75ex]
\mbox{}\hfill
$
\displaystyle
\prod\limits_{k=s+1}^{m}
{\rm g}_k^{{}^{\scriptsize \gamma_k^{}}}
\exp\sum\limits_{\nu=1}^{s}
\gamma_\nu^{}\;\!\omega_\nu^{}\in 
\text{\rm I}_{_{\Omega}}
\iff
\biggl(\,
\prod\limits_{k=s+1}^{m}
{\rm g}_k^{{}^{\scriptsize \gamma_k^{}}},\ 
\sum\limits_{\nu=1}^{s}\gamma_\nu^{}\,{\rm div}\;\!{\frak d}\biggr)\in 
\text{\rm J}_{_{\Omega}}.
\hfill
$
\\[1.75ex]
Кроме этого существуют такие функции $M_k^{}\in C^1\Omega,\ k=s+1,\ldots, m,$
что выполняются тождества {\rm (1.14)} и {\rm (11.3)}.
}
\vspace{1ex}

{\bf  Свойство 11.18.}
\vspace{0.75ex}
{\it
Пусть 
$
\bigl(\exp\omega_\nu^{}, M_\nu^{}\bigr)\in \text{\rm E}_{_{\Omega}},\ 
\nu=1,\ldots, s,\ s\leqslant m-1,\  
({\rm g}_k^{}, M_k^{})\in \text{\rm J}_{_{\Omega}},
\linebreak 
{\rm g}_k^{{}^{\scriptsize \gamma_k^{}}}\!\in C^1\Omega,
\ 
k=s+1,\ldots,m,\ 
\gamma_j^{}\in\R\backslash\{0\}, 
\ j=1,\ldots,m. 
$ 
Тогда 
\vspace{0.5ex}
$
\prod\limits_{k=s+1}^{m}
{\rm g}_k^{{}^{\scriptsize \gamma_k^{}}}
\exp\sum\limits_{\nu=1}^{s}
\gamma_\nu^{}\;\!\omega_\nu^{}\in 
\text{\rm I}_{_{\Omega}},
$
если и только если выполняется тождество {\rm (11.1)}.
}
\vspace{0.5ex}

{\sl Следует} 
из свойства 3.9 с учетом теоремы 11.1. $\k$
\vspace{1.25ex}

{\bf  Следствие 11.13.}
\vspace{0.75ex}
{\it
Пусть 
$
\exp\omega_\nu^{}\in \text{\rm ME}_{_{\Omega}},\ 
\nu=1,\ldots, s,\ s\leqslant m-1,\  
({\rm g}_k^{}, M_k^{})\in \text{\rm J}_{_{\Omega}},
\linebreak 
{\rm g}_k^{{}^{\scriptsize \gamma_k^{}}}\!\in C^1\Omega,
\ 
k=s+1,\ldots,m,\ 
\gamma_j^{}\in\R\backslash\{0\}, 
\ j=1,\ldots,m. 
$ 
Тогда 
\vspace{0.75ex}
$
\prod\limits_{k=s+1}^{m}
{\rm g}_k^{{}^{\scriptsize \gamma_k^{}}}
\exp\sum\limits_{\nu=1}^{s}
\gamma_\nu^{}\;\!\omega_\nu^{}\in 
\text{\rm I}_{_{\Omega}},
$
если и только если выполняется тождество {\rm (11.3)} на области $\Xi^{\;\!\prime}.$
}
\vspace{1ex}

{\bf  Следствие 11.14.}
\vspace{0.75ex}
{\it
Пусть 
$
\bigl(\exp\omega_\nu^{}, M_\nu^{}\bigr)\in \text{\rm E}_{_{\Omega}},\ 
\nu=1,\ldots, s,\ s\leqslant m-1,\  
\mu_k^{}\in \text{\rm M}_{_{\Omega}},
\linebreak 
\mu_k^{{}^{\scriptsize \gamma_k^{}}}\!\in C^1\Omega,
\ 
k=s+1,\ldots,m,\ 
\gamma_j^{}\in\R\backslash\{0\}, 
\ j=1,\ldots,m. 
$ 
Тогда 
\vspace{0.5ex}
$
\prod\limits_{k=s+1}^{m}
\mu_k^{{}^{\scriptsize \gamma_k^{}}}
\exp\sum\limits_{\nu=1}^{s}
\gamma_\nu^{}\;\!\omega_\nu^{}\in 
\text{\rm I}_{_{\Omega}},
$
если и только если линейная комбинация сомножителей}
\\[1.5ex]
\mbox{}\hfill                     
$
\displaystyle
\sum\limits_{\nu=1}^{s}\gamma_{\nu}^{}\;\!M_\nu^{}(t,x)=
\sum\limits_{k=s+1}^{m}
\gamma_{k}^{}\,{\rm div}\;\!{\frak d}(t,x)
\quad
\forall (t,x)\in\Xi^{\;\!\prime}.
\hfill
$
\\[2ex]
\indent
{\bf  Следствие 11.15.}
\vspace{0.75ex}
{\it
Пусть 
$
\exp\omega_\nu^{}\in \text{\rm ME}_{_{\Omega}},\ 
\nu=1,\ldots, s,\ s\leqslant m-1,\  
\mu_k^{}\in \text{\rm M}_{_{\Omega}},
\linebreak 
\mu_k^{{}^{\scriptsize \gamma_k^{}}}\!\in C^1\Omega,
\ 
k=s+1,\ldots,m,\ 
\gamma_j^{}\in\R\backslash\{0\}, 
\ j=1,\ldots,m. 
$ 
Тогда}
\\[1.5ex]
\mbox{}\hfill                     
$
\displaystyle
\prod\limits_{k=s+1}^{m}
\mu_k^{{}^{\scriptsize \gamma_k^{}}}
\exp\sum\limits_{\nu=1}^{s}
\gamma_\nu^{}\;\!\omega_\nu^{}\in 
\text{\rm I}_{_{\Omega}}
\iff
\sum\limits_{j=1}^{m}\gamma_{j}^{}= 0.
\hfill
$
\\[2ex]
\indent
{\bf Свойство\! 11.19.}\!
\vspace{0.75ex}
{\it
Пусть 
$\!\bigl(\!(p_j^{},\! M_j^{}), (h_j^{},q_j^{},\rho_{j}^{}N_0^{})\!\bigr)\!\!\in\! 
\text{\rm B}_{_{\Xi^{\;\!\prime}}},\, 
\rho_j^{}\!\in\!\R\backslash\{0\},\, 
\lambda_j^{},\!\gamma_j^{}\!\in\!\R, j\!=\!1,\ldots, m,\!$ 
$\sum\limits_{j=1}^{m}|\lambda_j^{}|\ne 0,\ 
\varphi\in C^1T^{\;\!\prime},$ 
мно\-жес\-т\-во
\vspace{0.75ex}
$\Omega_{_0}\subset\Xi^{\;\!\prime}$ та\-кое, что  
$\prod\limits_{j=1}^{m}p_j^{}(t,x)\ne 0\;\;\forall (t,x)\in\Omega_{_0},$
$\prod\limits_{j=1}^{m}p_j^{}(t,x)= 0\;\;
\forall (t,x)\in {\sf C}_{_{\Xi^{\;\!\prime}}}\Omega_{_0},\ 
p_j^{\gamma_j^{}}\in C^1\Omega_{_0}, \ j=1,\ldots, m.$
Тогда
\\[1.5ex]
\mbox{}\hfill                           
$
\displaystyle
\prod\limits_{j=1}^{m}p_j^{\gamma_j^{}}
\sum\limits_{j=1}^{m}\lambda_j^{}
\exp\biggl(\;\!\dfrac{q_j^{}}{\rho_j^{}\;\!p_j^{\;\!h_j^{}}}+\varphi\biggr)\in 
\text{\rm I}_{_{\Omega_{{}_{\tiny\;\! 0}}}},
\hfill
$
\\[2ex]
если и только если выполняется тождество}
\\[1.75ex]
\mbox{}\hfill                           
$
\displaystyle
{\sf D}\;\!\varphi(t)+N_0^{}(t,x)+
\sum\limits_{j=1}^{m}\gamma_j^{}\;\!M_j^{}(t,x)=0
\quad
\forall (t,x)\in \Xi^{\;\!\prime}.
\hfill
$
\\[1.5ex]
\indent
{\sl Следует} 
из свойства 5.5 с учетом теоремы 11.1. $\k$
\vspace{1.25ex}

{\bf Следствие 11.16.}\!
\vspace{0.75ex}
{\it
Пусть 
$\!\bigl(\nu_j^{}, (h_j^{},q_j^{},\rho_{j}^{}N_0^{})\!\bigr)\!\in\! 
\text{\rm MB}_{_{\Xi^{\;\!\prime}}},\, 
\rho_j^{}\!\in\!\R\backslash\{0\},\, 
\lambda_j^{},\gamma_j^{}\!\in\!\R,\, j\!=\!1,\ldots, m,\!$ 
$\sum\limits_{j=1}^{m}|\lambda_j^{}|\ne 0,\ 
\varphi\in C^1T^{\;\!\prime},$ 
мно\-жес\-т\-во
\vspace{0.75ex}
$\Omega_{_0}\subset\Xi^{\;\!\prime}$ та\-кое, что  
$\prod\limits_{j=1}^{m}\nu_j^{}(t,x)\ne 0\;\;\forall (t,x)\in\Omega_{_0},$
$\prod\limits_{j=1}^{m}\nu_j^{}(t,x)= 0\;\;
\forall (t,x)\in {\sf C}_{_{\Xi^{\;\!\prime}}}\Omega_{_0},\ 
\nu_j^{\gamma_j^{}}\in C^1\Omega_{_0}, \ j=1,\ldots, m.$
Тогда
\\[1.5ex]
\mbox{}\hfill                           
$
\displaystyle
\prod\limits_{j=1}^{m}\nu_j^{\gamma_j^{}}
\sum\limits_{j=1}^{m}\lambda_j^{}
\exp\biggl(\;\!\dfrac{q_j^{}}{\rho_j^{}\;\!\nu_j^{\;\!h_j^{}}}+\varphi\biggr)\in 
\text{\rm I}_{_{\Omega_{{}_{\tiny\;\! 0}}}},
\hfill
$
\\[2ex]
если и только если выполняется тождество}
\\[1.75ex]
\mbox{}\hfill                           
$
\displaystyle
{\sf D}\;\!\varphi(t)+N_0^{}(t,x)-
\sum\limits_{j=1}^{m}\gamma_j^{}\,{\rm div}\;\!{\frak d}(t,x)=0
\quad
\forall (t,x)\in \Xi^{\;\!\prime}.
\hfill
$
\\[2ex]
\indent
{\bf Свойство\! 11.20.}\!
\vspace{0.75ex}
{\it
Пусть 
$\!\bigl((p_j^{}, M_j^{}), (h_j^{},q_j^{},N_j^{})\bigr)\!\in\! \text{\rm B}_{_{\Xi^{\;\!\prime}}},\, 
\gamma_j^{},\xi_j^{}\!\in\!\R,\, \varphi_j^{}\!\in\! C^1T^{\;\!\prime},\, 
j\!=\!1,\ldots, m,\!$ 
мно\-жес\-т\-во
\vspace{0.75ex}
$\!\Omega_{_0}\!\subset\!\Xi^{\;\!\prime}\!$ та\-кое, что  
$\prod\limits_{j=1}^{m}p_j^{}(t,x)\!\ne\! 0\;\;\forall (t,x)\!\in\!\Omega_{_0},\
\prod\limits_{j=1}^{m}p_j^{}(t,x)= 0\;\;
\forall (t,x)\!\in\! {\sf C}_{_{\Xi^{\;\!\prime}}}\Omega_{_0},\!$ 
$p_j^{\gamma_j^{}}\in C^1\Omega_{_0}, \ j=1,\ldots, m.$
Тогда
\\[1.5ex]
\mbox{}\hfill                           
$
\displaystyle
\prod\limits_{j=1}^{m}p_j^{\gamma_j^{}}
\exp\sum\limits_{j=1}^{m}\xi_j^{}
\biggl(\;\!\dfrac{q_j^{}}{p_j^{\;\! h_j^{}}}+\varphi_j^{}\biggr)\in 
\text{\rm I}_{_{\Omega_{{}_{\tiny\;\! 0}}}},
\hfill
$
\\[2ex]
если и только если выполняется тождество}
\\[1.75ex]
\mbox{}\hfill                           
$
\displaystyle
\sum\limits_{j=1}^{m}\bigl(
\gamma_j^{}\;\!M_j^{}(t,x)+
\xi_j^{}\;\!\bigl(N_j^{}(t,x)+{\sf D}\;\!\varphi_j^{}(t)\bigr)\bigr)=0
\quad
\forall (t,x)\in \Xi^{\;\!\prime}.
\hfill
$
\\[1.5ex]
\indent
{\sl Следует} 
из свойства 5.7 с учетом теоремы 11.1. $\k$
\vspace{1.25ex}

{\bf Следствие 11.17.}\!
\vspace{0.75ex}
{\it
Пусть 
$\!\!\bigl(\nu_j^{}, (h_j^{},q_j^{},N_j^{})\bigr)\!\in \text{\rm MB}_{_{\Xi^{\;\!\prime}}},\, 
\gamma_j^{},\xi_j^{}\in\R,\ \varphi_j^{}\in C^1T^{\;\!\prime},\, 
j=1,\ldots, m,\!$ 
мно\-жес\-т\-во
\vspace{0.75ex}
$\!\Omega_{_0}\!\subset\!\Xi^{\;\!\prime}\!$ та\-кое, что  
$\prod\limits_{j=1}^{m}\nu_j^{}(t,x)\!\ne\! 0\;\;\forall (t,x)\!\in\!\Omega_{_0},\
\prod\limits_{j=1}^{m}\nu_j^{}(t,x)= 0\;\;
\forall (t,x)\!\in\! {\sf C}_{_{\Xi^{\;\!\prime}}}\Omega_{_0},\!$ 
$\nu_j^{\gamma_j^{}}\in C^1\Omega_{_0}, \ j=1,\ldots, m.$
Тогда
\\[1.5ex]
\mbox{}\hfill                           
$
\displaystyle
\prod\limits_{j=1}^{m}\nu_j^{\gamma_j^{}}
\exp\sum\limits_{j=1}^{m}\xi_j^{}
\biggl(\;\!\dfrac{q_j^{}}{\nu_j^{\;\! h_j^{}}}+\varphi_j^{}\biggr)\in 
\text{\rm I}_{_{\Omega_{{}_{\tiny\;\! 0}}}},
\hfill
$
\\[2ex]
если и только если выполняется тождество}
\\[1.75ex]
\mbox{}\hfill                           
$
\displaystyle
\sum\limits_{j=1}^{m}
\xi_j^{}\;\!\bigl(N_j^{}(t,x)+{\sf D}\;\!\varphi_j^{}(t)\bigr)=
\sum\limits_{j=1}^{m}
\gamma_j^{}\,{\rm div}\;\!{\frak d}(t,x)
\quad
\forall (t,x)\in \Xi^{\;\!\prime}.
\hfill
$
\\[2ex]
\indent
{\bf Свойство 11.21.}
\vspace{1ex}
{\it
Пусть 
$\Bigl((p, M), 
\Bigl(h_\xi^{},q_{_{\scriptstyle h_\xi^{}  f_\xi^{}}}, N_{h_\xi^{}  f_\xi^{}}^{}\Bigr)\Bigr)
\in \text{\rm B}_{_{\Xi^{\;\!\prime}}},\
\gamma_{_{\scriptstyle h_\xi^{}  f_\xi^{}}}\in\R,\ 
\varphi_{_{\scriptstyle h_\xi^{}  f_\xi^{}}}\in C^1T^{\;\!\prime},$ 
$f_{\xi}^{}=1,\ldots, \delta_{\xi}^{},\ \xi=1,\ldots,\varepsilon,\ \gamma\in\R\backslash\{0\},$
\vspace{1ex}
мно\-жес\-т\-во
$\Omega_{_0}\subset\Xi^{\;\!\prime}$ та\-кое, что  
$p(t,x)\ne 0$ $\forall (t,x)\in\Omega_{_0},$
$p(t,x)= 0\;\;\forall (t,x)\in {\sf C}_{_{\Xi^{\;\!\prime}}}\Omega_{_0},\ 
p^{\gamma}\in C^1\Omega_{_0}.$
Тогда
\\[1.5ex]
\mbox{}\hfill                           
$
\displaystyle
p^{\gamma}
\exp\sum\limits_{\xi=1}^{\varepsilon}\sum\limits_{f_{\xi}^{}=1}^{\delta_{\xi}^{}}
\biggl(\;\!\gamma_{_{\scriptstyle h_\xi^{}  f_\xi^{}}}\biggl(\,
\dfrac{q_{_{\scriptstyle h_\xi^{}  f_\xi^{}} }}{\displaystyle  p^{\;\!h_\xi^{}}}+
\varphi_{_{\scriptstyle h_\xi^{}  f_\xi^{}}}\biggr)\biggr)\in 
\text{\rm I}_{_{\Omega_{{}_{\tiny\;\! 0}}}},
\hfill
$
\\[2ex]
если и только если выполняется тождество}
\\[1.75ex]
\mbox{}\hfill                           
$
\displaystyle
\gamma\;\!M(t,x)+
\sum\limits_{\xi=1}^{\varepsilon}\sum\limits_{f_{\xi}^{}=1}^{\delta_{\xi}^{}}
\biggl(\gamma_{_{\scriptstyle h_\xi^{}  f_\xi^{}}}\Bigl(
N_{h_\xi^{}  f_\xi^{}}^{} (t,x)+{\sf D}\;\!\varphi_{_{\scriptstyle h_\xi^{}  f_\xi^{}}}(t)\Bigr)\biggr)=0
\quad
\forall (t,x)\in \Xi^{\;\!\prime}.
\hfill
$
\\[1.5ex]
\indent
{\sl Следует} 
из свойства 5.8 с учетом теоремы 11.1. $\k$
\vspace{1.25ex}

{\bf Следствие 11.18.}
\vspace{1ex}
{\it
Пусть 
$\Bigl(\nu, 
\Bigl(h_\xi^{},q_{_{\scriptstyle h_\xi^{}  f_\xi^{}}}, N_{h_\xi^{}  f_\xi^{}}^{}\Bigr)\Bigr)
\in \text{\rm MB}_{_{\Xi^{\;\!\prime}}},\
\gamma_{_{\scriptstyle h_\xi^{}  f_\xi^{}}}\in\R,\ 
\varphi_{_{\scriptstyle h_\xi^{}  f_\xi^{}}}\in C^1T^{\;\!\prime},$ 
$f_{\xi}^{}=1,\ldots, \delta_{\xi}^{},\ \xi=1,\ldots,\varepsilon,\ \gamma\in\R\backslash\{0\},$
\vspace{1ex}
мно\-жес\-т\-во
$\Omega_{_0}\subset\Xi^{\;\!\prime}$ та\-кое, что  
$\nu(t,x)\ne 0$ $\forall (t,x)\in\Omega_{_0},$
$\nu(t,x)= 0\;\;\forall (t,x)\in {\sf C}_{_{\Xi^{\;\!\prime}}}\Omega_{_0},\ 
\nu^{\gamma}\in C^1\Omega_{_0}.$
Тогда
\\[1.5ex]
\mbox{}\hfill                           
$
\displaystyle
\nu^{\gamma}
\exp\sum\limits_{\xi=1}^{\varepsilon}\sum\limits_{f_{\xi}^{}=1}^{\delta_{\xi}^{}}
\biggl(\;\!\gamma_{_{\scriptstyle h_\xi^{}  f_\xi^{}}}\biggl(\,
\dfrac{q_{_{\scriptstyle h_\xi^{}  f_\xi^{}} }}{\displaystyle  \nu^{\;\!h_\xi^{}}}+
\varphi_{_{\scriptstyle h_\xi^{}  f_\xi^{}}}\biggr)\biggr)\in 
\text{\rm I}_{_{\Omega_{{}_{\tiny\;\! 0}}}},
\hfill
$
\\[2ex]
если и только если выполняется тождество}
\\[1.75ex]
\mbox{}\hfill                           
$
\displaystyle
{\rm div}\;\!{\frak d}(t,x)=
\dfrac{1}{\gamma}\
\sum\limits_{\xi=1}^{\varepsilon}\sum\limits_{f_{\xi}^{}=1}^{\delta_{\xi}^{}}
\biggl(\gamma_{_{\scriptstyle h_\xi^{}  f_\xi^{}}}\Bigl(
N_{h_\xi^{}  f_\xi^{}}^{} (t,x)+{\sf D}\;\!\varphi_{_{\scriptstyle h_\xi^{}  f_\xi^{}}}(t)\Bigr)\biggr)
\quad
\forall (t,x)\in \Xi^{\;\!\prime}.
\hfill
$
\\[2ex]
\indent
{\bf Свойство 11.22.}
\vspace{0.5ex}
{\it
Пусть числа $\gamma_1^{}, \gamma_2^{}\in\R,$
функции
$u,v\in \text{\rm P}_{_{\!\Xi^{\;\!\prime}}}$ взаимно простые, 
$(u+i\;\!v,\;\! U+i\;\!V)\in \text{\rm H}_{_{\Xi^{\;\!\prime}}},$
\vspace{0.5ex}
множество $\Omega_{_0}\!\subset\! \Xi^{\;\!\prime}\!$ такое, что  
$u(t,x)\ne 0\;\;\forall (t,x)\!\in\! \Omega_{_0},\ 
u(t,x)=0$ $\forall (t,x)\in {\sf C}_{_{\Xi^{\;\!\prime}}}\Omega_{_0}.$
Тогда
$(u^2+v^2)^{{}^{\scriptsize \gamma_1^{}}}
\exp\Bigl(\gamma_2^{}\;\!\arctg\dfrac{v}{u}\Bigr)\in 
\text{\rm I}_{_{\Omega_{{}_{\tiny\;\! 0}}}},$
если и только если}
\\[1.75ex]
\mbox{}\hfill
$
2\gamma_1^{}\;\!U(t,x)+\gamma_2^{}\;\!V(t,x)=0
\quad
\forall (t,x)\in \Xi^{\;\!\prime}.
\hfill
$
\\[1.5ex]
\indent
{\sl Следует} 
из свойства 6.15 с учетом теоремы 11.1. $\k$
\vspace{1.25ex}

{\bf Свойство 11.23.}
\vspace{0.5ex}
{\it
Пусть 
$\bigl((u+i\;\!v,\;\! U+i\;\!V), (h,z,Q)\bigr)\in \text{\rm G}_{_{\Xi^{\;\!\prime}}},$
функции
$u, v\in \text{\rm P}_{_{\!\Xi}}$ вза\-им\-но простые, 
\vspace{0.75ex}
функция $z\in \text{\rm Z}_{_{\Xi}}$ взаимно простая с 
функцией $u+i\;\!v,\ \gamma_1^{},\gamma_2^{},\gamma_3^{},\gamma_4^{}\in\R,$  множество
$\!\Omega_{_0}\!\subset\!\Xi^{\;\!\prime}\!$ та\-кое, что  
$\!u(t,x)\!\ne\! 0\;\forall (t,x)\in\Omega_{_0},\
u(t,x)\!=\! 0\;\forall (t,x)\in {\sf C}_{_{\Xi^{\;\!\prime}}}\Omega_{_0}.\!$
Тогда
\\[1.75ex]
\mbox{}\hfill                           
$
\displaystyle
\bigl(u^2+v^2\bigr)^{\!{}^{\scriptstyle\gamma_1^{}}}
\exp\biggl(
\dfrac{\gamma_2^{}\;\!{\rm Re}\;\!\bigl(z(u-i\;\!v)^h\bigr)+
\gamma_3^{}\;\!{\rm Im}\;\!\bigl(z(u-i\;\!v)^h\bigr)}{\bigl(u^2+v^2\bigr)^{h}}\,+\,
\gamma_4^{}\;\!\arctg\dfrac{v}{u}
\biggr)\in \text{\rm I}_{_{\Omega_{{}_{\tiny\;\! 0}}}},
\hfill                           
$
\\[1.75ex]
если и только если}
\\[1.75ex]
\mbox{}\hfill
$
2\gamma_1^{}\;\!U(t,x)+
\gamma_2^{}\;\!{\rm Re}\;\!Q(t,x)+
\gamma_3^{}\;\!{\rm Im}\;\!Q(t,x)+
\gamma_4^{}\;\!V(t,x)=0
\quad
\forall (t,x)\in \Xi^{\;\!\prime}.
\hfill
$
\\[1.5ex]
\indent
{\sl Следует} 
из теоремы 7.3 с учетом свойства 1.9 и теоремы 11.1. $\k$
\vspace{1.25ex}

\newpage

\mbox{}
\\[-0.15ex]
\centerline{
{\bf  12. Приложения 
}
}
\\[2.25ex]
\indent
{\bf 12.1.}
{\sl Задача Дарбу}.
Рассмотрим дифференциальное уравнение первого порядка
\\[2ex]
\mbox{}\hfill                     % (12.1)
$
Y(x,y)\;\!dx-X(x,y)\;\!dy=0,
$
\hfill (12.1)
\\[2ex]
у которого функции $X\colon\R^2\to\R$ и $Y\colon\R^2\to\R$ 
\vspace{0.35ex}
суть полиномы по переменным $x$ и $y$ с такими степенями, что 
$
\max\bigl\{\deg X,\, \deg Y\bigr\}=d\geq 1.
$ 
\vspace{0.75ex}

Задача Дарбу [8, 13]
\vspace{0.15ex}
состоит в построении общего интеграла дифференциального уравнения (12.1)
по известным частным интегралам.
\vspace{0.25ex}

Уравнение (12.1) является уравнением траекторий автономной 
дифференциальной системы второго порядка
\\[1.75ex]
\mbox{}\hfill                     % (12.2)
$
\dfrac{dx}{dt}=X(x,y),
\qquad
\dfrac{dy}{dt}=Y(x,y).
$
\hfill (12.2)
\\[1.75ex]
\indent
Общий интеграл уравнения (12.1) является автономным первым интегралом системы (12.2) и наоборот.
\vspace{0.25ex}

Случай $d=2.$
Пусть функция 
\\[1.5ex]
\mbox{}\hfill
$
{\rm g}_j^{}\colon (x,y)\to\ 
{\rm g}_j^{}(x,y)
\quad
\forall (x,y)\in D
\hfill
$
\\[1.5ex]
является частным интегралом с сомножителем $M_j^{}$
\vspace{0.25ex}
на области $D\subset\R^2$ уравнения (12.1) 
(или, что то же, системы (12.2)), $j=1,2,3.$
\vspace{0.75ex}

Согласно определению 1.1 
$\deg M_j^{}\leq 1,\ j=1,2,3,$ т.е.
\\[1.75ex]
\mbox{}\hfill                    
$
M_j^{}(x,y)=
\alpha_{1j}^{}+\alpha_{2j}^{}\;\!x+\alpha_{3j}^{}\;\!y
\quad
\forall (x,y)\in\R^2,
\ \ 
j=1,2,3,
\hfill
$
\\[1.75ex]
где $\alpha_{1j}^{}\;\!,\;\!\alpha_{2j}^{}\;\!,\;\!\alpha_{3j}^{}\in\R,\ 
\sum\limits_{\zeta=1}^{3}|\alpha_{\zeta j}^{}|\ne 0,\ j=1,2,3.$
\vspace{0.5ex}

Множество $D_{_{\!0}}^{}\subset D$ такое, что 
${\rm g}_j^{{}^{\!\scriptsize\gamma_j^{}}}\in C^{1}D_{_0}^{},\ 
\gamma_j^{}\in\R,\ j=1,2,3.$
По свойству 11.11, функция
\\[1.5ex]
\mbox{}\hfill                     % (12.3)
$
\displaystyle
F\colon (x,y)\to\ 
\prod\limits_{j=1}^3
{\rm g}_j^{{}^{\!\scriptsize\gamma_j^{}}}(x,y)
\quad
\forall (x,y)\in D_{_{\!0}}^{}
$
\hfill (12.3)
\\[1.75ex]
является общим интегралом на множестве $D_{_{\!0}}^{}$
уравнения (12.1), если и только если 
\\[1.5ex]
\mbox{}\hfill                    
$
\displaystyle
\sum\limits_{j=1}^3
\gamma_j^{}\;\!
\bigl(\alpha_{1j}^{}+\alpha_{2j}^{}\;\!x+\alpha_{3j}^{}\;\!y\bigr)=0
\quad
\forall (x,y)\in\R^2.
\hfill
$
\\[1.75ex]
Последнее возможно, когда определитель
\\[2ex]
\mbox{}\hfill                    
$
\triangle=\left|
\begin{array}{ccc}
\alpha_{11}^{} & \alpha_{12}^{} & \alpha_{13}^{} 
\\[0.75ex]
\alpha_{21}^{} & \alpha_{22}^{} & \alpha_{23}^{} 
\\[0.75ex]
\alpha_{31}^{} & \alpha_{32}^{} & \alpha_{33}^{} 
\end{array}
\right|=0.
\hfill
$
\\[1.75ex]
\indent
Итак, задача Дарбу решена.
\vspace{0.35ex}

В статье [7] нами было предложено расширить постановку задачи Дарбу. 
Поступим таким же образом и наряду с построением общего интеграла будем 
находить интегрирующий множитель уравнения (12.1). 

По свойству 8.8, функция (12.3) является интегрирующим множителем 
уравнения (12.1), если и только если линейная комбинация сомножителей
\\[1.75ex]
\mbox{}\hfill              % (12.4)               
$
\displaystyle
\sum\limits_{j=1}^3
\gamma_j^{}\;\!M_j^{}(x,y)=
{}-\partial_x^{}X(x,y)-\partial_y^{}Y(x,y)
\quad
\forall (x,y)\in\R^2.
$
\hfill (12.4)
\\[1.75ex]
\indent
Расходимость векторного поля, определяемого уравнением (12.1) при $d=2,$
\\[2ex]
\mbox{}\hfill                    
$
\partial_x^{}X(x,y)+\partial_y^{}Y(x,y)=
\beta_{1}^{}+\beta_{2}^{}\;\!x+\beta_{3}^{}\;\!y
\quad
\forall (x,y)\in\R^2,
\hfill
$
\\[1.75ex]
где $\beta_{1}^{}\;\!,\;\!\beta_{2}^{}\;\!,\;\!\beta_{3}^{}$  --- действительные числа.
\vspace{0.75ex}

Тождество (12.4) выполняется, если и только если 
совместна система уравнений
\\[2ex]
\mbox{}\hfill              % (12.5)               
$
\begin{array}{c}
\alpha_{11}^{}\;\!\gamma_1^{}+
\alpha_{12}^{}\;\!\gamma_2^{}+
\alpha_{13}^{}\;\!\gamma_3^{}={}-\beta_1^{},
\\[1.5ex]
\alpha_{21}^{}\;\!\gamma_1^{}+
\alpha_{22}^{}\;\!\gamma_2^{}+
\alpha_{23}^{}\;\!\gamma_3^{}={}-\beta_2^{},
\\[1.5ex]
\alpha_{31}^{}\;\!\gamma_1^{}+
\alpha_{32}^{}\;\!\gamma_2^{}+
\alpha_{33}^{}\;\!\gamma_3^{}={}-\beta_3^{}.
\end{array}
$
\hfill (12.5)
\\[2.25ex]
\indent
Если $\triangle\ne 0,\ |\beta_1^{}|+|\beta_2^{}|+|\beta_3^{}|\ne 0,$
\vspace{0.35ex}
то из линейной неоднородной системы (12.5) находим 
$\gamma_{1}^{}, \gamma_{2}^{}, \gamma_{3}^{},\!$
\vspace{0.75ex}
а значит, и интегрирующий множитель\! (12.3) на множестве
$\!D_{_{\!0}}^{}\!$ уравнения\! (12.1).

Если $\beta_1^{}=\beta_2^{}=\beta_3^{}= 0,$
\vspace{0.35ex}
то  уравнение (12.1) является уравнением в полных дифференциалах 
(имеет интегрирующий множитель в виде произвольной действительной постоянной).

Таким образом, если известно три частных интеграла уравнения (12.1) при $d=2,$
то задача Дарбу в расширенной постановке решается в замкнутой форме --- 
строится общий интеграл или интегрирующий множитель.
\vspace{0.5ex}

Общий случай $d\geq 2.$
Из выше приведенных рассуждений заключаем, что 
количество частных интегралов, которые надо знать, чтобы 
по ним построить общий интеграл или интегрирующий множитель
(решить расширенную задачу Дарбу) уравнения (12.1)
равно количеству одночленов у полинома двух переменных $x$ и $y$
степени $d-1$ (наибольшая степень сомножителей частных интегралов).

Известно, что полином двух переменных степени $d-1$ имеет $\dfrac{d\;\!(d+1)}{2}$ одночленов.
Поэтому справедлива
\vspace{0.75ex}

{\bf Теорема 12.1.}
\vspace{0.25ex}
{\it 
Пусть известно $m=\dfrac{d\;\!(d+1)}{2}$
частных интегралов ${\rm g}_j^{}$ на области $D$ уравнения {\rm (12.1)},
$\gamma_j^{}\in\R, $
\vspace{0.5ex}
множество $D_{_{\!0}}^{}\subset D$ такое, что 
${\rm g}_j^{{}^{\!\scriptsize\gamma_j^{}}}\in C^{1}D_{_0}^{},\ 
j=1,\ldots, m.$
Тогда функция
\\[1.5ex]
\mbox{}\hfill                   
$
\displaystyle
F\colon (x,y)\to\ 
\prod\limits_{j=1}^m
{\rm g}_j^{{}^{\!\scriptsize\gamma_j^{}}}(x,y)
\quad
\forall (x,y)\in D_{_{\!0}}^{}
\hfill
$
\\[1.5ex]
является общим интегралом или интегрирующим множителем 
на множестве $D_{_{\!0}}^{}$ уравнения} (12.1).
\vspace{1.25ex}

{\bf 12.2.}
Дифференциальное уравнение первого порядка 
[14, с. 136 -- 139; 23]
\\[2ex]
\mbox{}\hfill                   
$
\bigl(x+a^2x^2+2\;\!axy-(1+2\;\!a^2)y^2\bigr)\;\!dx+
\bigl(y+ax^2+(3+a^2)xy-ay^2\bigr)\;\!dy=0
\quad
(a\in\R\backslash\{0\})
\hfill
$
\\[2ex]
имеет два полиномиальных частных интеграла:
\\[1.5ex]
\mbox{}\hfill                   
$
p_1^{}\colon (x,y)\to\ 
1+(1+a^2)\bigl(3\;\!x+3\;\!ax(y+ax)+a(y+ax)^3\bigr)
\quad
\forall (x,y)\in\R^2
\hfill
$
\\[1.5ex]
с сомножителем 
\\[1ex]
\mbox{}\hfill                   
$
M_1^{}\colon (x,y)\to\ 
{}-3\;\!(1+a^2)\;\!y
\quad
\forall (x,y)\in\R^2
\hfill
$
\\[1ex]
и 
\\[1ex]
\mbox{}\hfill                   
$
p_2^{}\colon (x,y)\to\ 
1+(1+a^2)\bigl(2\;\!x+(y+ax)^2\bigr)
\quad
\forall (x,y)\in\R^2
\hfill
$
\\[1.5ex]
с сомножителем 
\\[1ex]
\mbox{}\hfill                   
$
M_2^{}\colon (x,y)\to\ 
{}-2\;\!(1+a^2)\;\!y
\quad
\forall (x,y)\in\R^2.
\hfill
$
\\[2.25ex]
\indent
По следствию 11.5 при $M_0^{}={}-(1+a^2)\;\!y,\ \rho_1^{}=3,\ \rho_2^{}=2,$
\vspace{0.75ex}
из равенства $3\gamma_1^{}+2\gamma_2^{}=0$
находим, например, что $\gamma_1^{}=2,\ \gamma_2^{}={}-3.$
Следовательно, рациональная функция
\\[2ex]
\mbox{}\hfill                   
$
F\colon (x,y)\to\ 
\dfrac{p_1^{\,2}(x,y)}{p_2^{\,3}(x,y)}
\quad
\forall (x,y)\in D_{_{\!0}}^{}
\hfill
$
\\[2ex]
является общим интегралом на множестве $D_{_{\!0}}^{}\subset\R^2$
\vspace{0.35ex}
таком, что $p_2^{}(x,y)\ne 0\;\; \forall (x,y)\in D_{_{\!0}}^{},$ 
$p_2^{}(x,y)= 0\;\; \forall (x,y)\in {\sf C}_{\R^2}^{}D_{_{\!0}}^{}.$
\vspace{0.75ex}

Этот пример есть один из тех случаев, когда общий интеграл 
уравнения (12.1) при $d=2$ строится не по трем, а по меньшему числу
(двум) частным интегралам.
\vspace{1ex}

{\bf 12.3.}
Система [28, с. 46]
\\[2ex]
\mbox{}\hfill                   
$
\dfrac{dx}{dt}={}-y+\dfrac{x^2}{2}-\dfrac{y^2}{2}\,,
\qquad
\dfrac{dy}{dt}=x\;\!(1+y)
\hfill
$
\\[2.25ex]
имеет: 
$(y+1, x)\in\text{A}_{\R^2}^{}$ (свойство 2.3) и 
$(x^2+y^2, x)\in\text{A}_{\R^2}^{}$ (теорема 2.1).
\vspace{1ex}

По следствию 11.7, 
$\dfrac{y+1}{x^2+y^2}\in \text{I}_{\R^2\backslash\{(0,0)\}}^{}.$
\vspace{1.75ex}

{\bf 12.4.}
Система [23]
\\[2ex]
\mbox{}\hfill                   % (12.6)
$
\dfrac{dx}{dt}={}-2+y+x^2+xy,
\qquad
\dfrac{dy}{dt}=4+2\;\!x+xy+y^2
$
\hfill (12.6)
\\[2.25ex]
имеет: 
$\bigl((2+2\;\!x+y,\;\! x+y),\;\! (1, x+y, 1)\bigr)\in\text{B}_{\R^2}^{}.$ 
По теореме 5.3, 
\\[2ex]
\mbox{}\hfill
$
(2+2\;\!x+y,\;\! x+y)\in\text{A}_{\R^2}^{},
\qquad
\Bigl(\exp\dfrac{x+y}{2+2\;\!x+y}\,,\;\! 1\Bigr)\in\text{E}_{\!D}^{}\;\!,
\hfill
$ 
\\[1.5ex]
где множество $D=\{(x,y)\colon 2+2\;\!x+y\ne 0\}.$
\vspace{1ex}

Согласно теореме 1.1 
$\bigl(12+8\;\!x+4\;\!y+4\;\!xy+3\;\!y^2,\;\! 2\;\!(x+y)\bigr)\in\text{A}_{\R^2}^{}.$
\vspace{0.75ex}

На основании полиномиальных частных интегралов в 
\vspace{0.35ex}
соответствии со следствием~11.5 (при 
$\rho_1^{}=1,\ \rho_2^{}=2,\ M_0^{}=x+y,\ \gamma_1^{}={}-2,\ \gamma_2^{}=1)$
\vspace{0.35ex}
строим автономный первый интеграл
\\[1.5ex]
\mbox{}\hfill                   
$
F\colon (x,y)\to\ 
\dfrac{12+8\;\!x+4\;\!y+4\;\!xy+3\;\!y^2}{(2+2\;\!x+y)^2}
\quad
\forall (x,y)\in D.
\hfill
$
\\[1.75ex]
\indent
На основании экспоненциального частного интеграла, 
используя следствие 11.2 (при $\lambda=1),$
строим неавтономный первый интеграл
\\[1.75ex]
\mbox{}\hfill                   
$
F_1^{}\colon (t,x,y)\to\ 
e^{{}-t}\exp\dfrac{x+y}{2+2\;\!x+y}=
\exp\Bigl({}-t+\dfrac{x+y}{2+2\;\!x+y}\Bigr)
\quad
\forall (x,y)\in\Omega_{_0}^{},
\hfill
$
\\[2ex]
где множество $\Omega_{_0}^{}=\{(t,x,y)\colon 2+2\;\!x+y\ne 0\}.$
\vspace{0.5ex}

Учитывая функциональную неоднозначность первого интеграла (свойство 0.1),
первый интеграл $F_1^{}$ представим в виде
\\[2ex]
\mbox{}\hfill                   
$
\Psi\colon (t,x,y)\to\ 
{}-t+\dfrac{x+y}{2+2\;\!x+y}
\quad
\forall (x,y)\in\Omega_{_0}^{}.
\hfill
$
\\[2ex]
\indent
Будучи функционально независимыми, первые интегралы $F$ и $\Psi$
\vspace{0.25ex}
образуют интегральный базис на множестве $\Omega_{_0}^{}$ системы (12.6).
\vspace{1ex}

\newpage

{\bf 12.5.}
Система [24, 23]
\\[1.75ex]
\mbox{}\hfill                   % (12.7)
$
\dfrac{dx}{dt}={}-y\;\!\bigl(2\;\!x^2+y^2+(x^2+y^2)^2\bigr),
\qquad
\dfrac{dy}{dt}=x\;\!\bigl(2\;\!x^2+y^2+2\;\!(x^2+y^2)^2\bigr)
$
\hfill (12.7)
\\[2.25ex]
имеет: 
$
(2\;\!x^2+y^2,\;\! {}-2\;\!xy)\in\text{A}_{\R^2}^{}$ и (согласно свойству 5.10)
\\[2.25ex]
\mbox{}\hfill                   
$
\bigl((x^2+y^2,\;\! 2\;\!xy(x^2+y^2)),\;\! (1, 1, {}-2\;\!xy)\bigr)\in\text{B}_{\R^2}^{}.
\hfill
$ 
\\[1.5ex]
У некратного полиномиального частного интеграла и экспонциального частного интеграла (теорема 5.3)
один и тот же сомножитель. По следствию 11.7, функция
\\[1.75ex]
\mbox{}\hfill                   
$
F\colon (x,y)\to\ 
(2\;\!x^2+y^2)\exp\dfrac{{}-1}{x^2+y^2}
\quad
\forall (x,y)\in\R^2\backslash\{(0,0)\}
\hfill
$
\\[1.75ex]
является автономным первым интегралом на области 
$\R^2\backslash\{(0,0)\}$ системы (12.7).
\vspace{1ex}

{\bf 12.6.}
Система [26]
\\[1.75ex]
\mbox{}\hfill                  
$
\dfrac{dx}{dt}=x+xy+y^2,
\qquad
\dfrac{dy}{dt}=y+x^2-xy+2\;\!y^2
\hfill
$
\\[2.25ex]
имеет: 
$\bigl(x+i\;\!y,\;\! 1+2\;\!y+i\;\!(x-y)\bigr)\in\text{H}_{\R^2}^{}$ (теорема 6.1), 
а значит (теорема 6.2), 
\\[2.25ex]
\mbox{}\hfill                   
$
\bigl(x^2+y^2,\;\! 2\;\!(1+2\;\!y)\bigr)\in\text{A}_{\R^2}^{},\ \
\Bigl(\exp\arctg\dfrac{y}{x}\,,\, x-y\Bigr)\in\text{E}_{\!D}^{}\;\!,\ \
D=\{(x,y)\colon x\ne 0\}\;\!;
\hfill
$
\\[2ex]
$\bigl(x-y,\;\! {}-(x-y-1)\bigr)\in\text{A}_{\R^2}^{}$ 
\vspace{1.25ex}
(теорема 2.1), а значит (свойство 1.4), 
$\bigl(x-y-1,\;\! {}-(x-y)\bigr)\in\text{A}_{\R^2}^{}.$ 

По свойству 11.11, 
\\[1.25ex]
\mbox{}\hfill                  
$
(x-y-1)\exp\arctg\dfrac{y}{x}\in \text{I}_{D}^{}\;\!.
\hfill
$
\\[1.5ex]
\indent
По свойству 1.9, 
\\[1.5ex]
\mbox{}\hfill                  
$
\Bigl(\dfrac{x-y-1}{x-y}\,,\;\! {}-1\Bigr)\in \text{J}_{\!D_{_{0}}}^{},
\quad 
D_{_{\!0}}=\{(x,y)\colon x-y\ne 0\},
\hfill
$
\\[1.75ex]
а значит (следствие~11.2),
\\[1.25ex]
\mbox{}\hfill                  
$
\dfrac{x-y-1}{x-y}\,e^{\,t}\in  \text{I}_{\R\times D_{_{0}}}^{}.
\hfill
$
\\[2ex]
\indent
{\bf 12.7.}
Дифференциальное уравнение 
\\[1.75ex]
\mbox{}\hfill                   % (12.8)
$
\Bigl({}-x+\dfrac{1}{2}\,y^2\Bigr)\;\!dx-
\bigl(xy+2\;\!y^3\bigr)\;\!dy=0
$
\hfill (12.8)
\\[1.5ex]
на $\R^2$ имеет комплекснозначный полиномиальный частный интеграл
\\[1.5ex]
\mbox{}\hfill                   
$
w\colon (x,y)\to\ 
x+i\;\!y^2
\quad
\forall (x,y)\in\R^2
\hfill
$
\\[1ex]
с сомножителем 
\\[1ex]
\mbox{}\hfill                   
$
W\colon y\to\ 
y-2\;\!i\;\!y
\quad
\forall y\in\R.
\hfill
$
\\[1ex]
\indent
По теореме 6.2,
\\[1.5ex]
\mbox{}\hfill                   
$
\bigl(x^2+y^4,\;\! 2\;\!y\bigr)\in\text{J}_{\R^2}^{}\;\!,
\quad
\Bigl(\exp\arctg\dfrac{y^2}{x}\,,\;\! {}-2\;\!y\Bigr)\in\text{E}_{\!D}^{}\;\!,\ \
D=\{(x,y)\colon x\ne 0\}.
\hfill
$
\\[2ex]
По следствию 11.5 (при 
$\rho_1^{}=1,\ \rho_2^{}={}-1,\ M_0^{}=2\;\!y,\ \gamma_1^{}=\gamma_2^{}=1),$
функция
\\[2.25ex]
\mbox{}\hfill                   
$
F\colon (x,y)\to\ 
\bigl(x^2+y^4\bigr)\;\!
\exp\arctg\dfrac{y^2}{x}
\quad
\forall (x,y)\in D
\hfill
$
\\[1.75ex]
является общим интегралом на множестве $D$ уравнения (12.8).
\vspace{1ex}

{\bf 12.8.}
На приложение условных частных интегралов нами впервые было обращено внимание 
в статьях [15, 16, 12] при качественном исследовании поведения траекторий на 
проективной фазовой плоскости системы (12.2) [17 -- 20].

В качестве примера дифференциального уравнения с условным частным интегралом 
и условным интегрирующим множителем рассмотрим дифференциальное уравнение 
из статьи [12]
\\[1.5ex]
\mbox{}\hfill                   % (12.9)
$
(x+x^2-y^2+a)\;\!dx-
(y+x^2-y^2+a)\;\!dy=0
\quad 
(a\in\R).
$
\hfill (12.9)
\\[1.5ex]
Полином
\\[1ex]
\mbox{}\hfill                   
$
p\colon (x,y)\to\ 
x^2-y^2+a
\quad
\forall (x,y)\in\R^2
\hfill
$
\\[2ex]
является частным интегралом на $\R^2$ с сомножителем 
\\[1.75ex]
\mbox{}\hfill                   
$
M_1^{}\colon (x,y)\to\ 
2\;\!(x-y)
\quad
\forall (x,y)\in\R^2,
\hfill
$
\\[1.5ex]
а экспоненциальная функция
\\[1.5ex]
\mbox{}\hfill                   
$
{\rm g}\colon (x,y)\to\ 
\exp(x-y)
\quad
\forall (x,y)\in\R^2
\hfill
$
\\[1.75ex]
является условным частным интегралом на $\R^2$ с сомножителем 
\\[2ex]
\mbox{}\hfill                   
$
M_2^{}\colon (x,y)\to\ 
{}-(x-y)
\quad
\forall (x,y)\in\R^2.
\hfill
$
\\[2ex]
\indent
По следствию 11.5 (при 
$\rho_1^{}=2,\ \rho_2^{}={}-1,\ M_0^{}=x-y,\ \gamma_1^{}=1,\ \gamma_2^{}=2),$
функция
\\[2.25ex]
\mbox{}\hfill                   
$
F\colon (x,y)\to\ 
\bigl(x^2-y^2+a\bigr)\;\!
\exp\bigl(2\;\!(x-y)\bigr)
\quad
\forall (x,y)\in \R^2
\hfill
$
\\[1.75ex]
является общим интегралом на $\R^2$ уравнения (12.9).
\vspace{0.5ex}

Заметим, что
\\[1ex]
\mbox{}\hfill                   
$
{\rm div}\;\!{\frak d}(x,y)=2\;\!(x-y)
\quad
\forall (x,y)\in \R^2.
\hfill
$
\\[2ex]
\indent
Так как сомножитель $M_1^{}={\rm div}\;\!{\frak d},$
то, по свойству 8.9 (при 
$\rho_1^{}=1,\ \gamma_1^{}={}-1),$
рациональная функция
\\[2ex]
\mbox{}\hfill                   
$
\mu_1^{}\colon (x,y)\to\ 
\dfrac{1}{x^2-y^2+a}
\quad
\forall (x,y)\in D=\{(x,y)\colon x^2-y^2+a\ne 0\}
\hfill
$
\\[1.75ex]
является интегрирующим множителем на множестве $D$ уравнения (12.9).
\vspace{0.75ex}

Так как сомножитель $M_2^{}={}-\dfrac{1}{2}\,{\rm div}\;\!{\frak d},$
\vspace{0.5ex}
то, по свойству 10.11 (при 
$\rho_1^{}=\dfrac{1}{2}\,,\ \gamma_1^{}=2),$
экспоненциальная функция
\\[1.75ex]
\mbox{}\hfill                   
$
\mu_2^{}\colon (x,y)\to\ 
\exp\bigl(2\;\!(x-y)\bigr)
\quad
\forall (x,y)\in \R^2
\hfill
$
\\[2ex]
является условным интегрирующим множителем на $\R^2$ уравнения (12.9).
\vspace{0.5ex}

На основании свойства Якоби интегрирующих множителей (свойство 0.2),
\vspace{0.35ex}
общий интеграл $F=\mu_2^{}/\mu_1^{}.$
\vspace{1.5ex}

{\bf 12.9.}
{\sl Расширенная задача Дарбу}.
Рассмотрим автономную дифференциальную систему $n\!$-го порядка
\\[2ex]
\mbox{}\hfill                     % (12.10)
$
\dfrac{dx_i^{}}{dt}=X_i^{}(x_1^{},\ldots, x_n^{}),
\quad
i=1,\ldots, n,
$
\hfill (12.10)
\\[2.5ex]
у которой правые части $X_i^{}\colon\R^n\to\R,\ i=1,\ldots,n,$
\vspace{0.35ex}
суть полиномы по зависимым переменным $x_1^{},\ldots, x_n^{}$
\vspace{0.75ex}
с постоянными коэффициентами, имеющие такие степени, что 
$\max\{\deg X_i^{}\colon i=1,\ldots, n\}=d\geq 1.$
\vspace{0.75ex}

\newpage

{\bf Теорема 12.2.}
\vspace{0.5ex}
{\it 
Пусть известно\footnote{
$\tbinom{k}{n}=\dfrac{n!}{k!\;\!(n-k)!}$ --- биномиальные коэффициенты.
} 
$m=\tbinom{n+d-1}{n}$
\vspace{0.5ex}
частных интегралов ${\rm g}_j^{}$ на области $X^{\;\!\prime}$ 
дифференциальной системы {\rm (12.10)},
$\gamma_j^{}\in\R, $
\vspace{0.5ex}
множество $X_{_{\!0}}^{}\subset X^{\;\!\prime}$ такое, что 
${\rm g}_j^{{}^{\!\scriptsize\gamma_j^{}}}\in C^{1}X_{_0}^{},\ 
j=1,\ldots, m.$
Тогда функция
\\[1.5ex]
\mbox{}\hfill                   
$
\displaystyle
F\colon x\to\ 
\prod\limits_{j=1}^m
{\rm g}_j^{{}^{\!\scriptsize\gamma_j^{}}}(x)
\quad
\forall x\in X_{_{\!0}}^{}
\hfill
$
\\[1.5ex]
является автономным первым интегралом 
\vspace{0.15ex}
или автономным последним множителем 
на множестве $X_{_{\!0}}^{}$ системы} (12.10).
\vspace{0.5ex}

{\sl Доказательство}
аналогично доказательству теоремы 12.1 с учетом того, что у полинома 
$n$ переменных степени $d-1$ количество одночленов равно числу сочетаний из $n+d-1$
по $n.\ \k$
\vspace{1ex}

В [9, с. 45 -- 47] и [21] для системы (12.10) 
\vspace{0.75ex}
расширенная задача Дарбу (теорема 12.2) 
решена при наличии $\tbinom{n+d-1}{n}$
\vspace{0.5ex}
автономных полиномиальных частных интегралов, определяющих 
интегральные многообразия. 
\vspace{0.35ex}

В [25, 31] для дифференциальной системы (12.10), 
а также для многомерной дифференциальной системы 
и дифференциальной системы (0.1), 
\vspace{0.5ex}
расширенная задача Дарбу решена при наличии 
$\tbinom{n+d-1}{n}$
\vspace{0.5ex}
полиномиальных частных интегралов с учетом их кратностей и 
условных частных интегралов.

Обратим внимание на то, что статья [25]
является первой публикацией, в которой нами был использован
термин-композит <<условный частный интеграл>> 
(исходя из того, что этот частный интеграл не определяет интегрального многообразия)
и впервые было введено понятие кратного полиномиального частного интеграла. 
При этом вместо <<$\!\varkappa\!$-кратный полиномиальный частный интеграл>>
использовался оборот слов <<полиномиальный частный интеграл с весом $\varkappa\!$>>,
отражая тем самым, что при построении первого интеграла и последнего множителя 
он считается за $\varkappa$ частных интегралов.

Композиционный термин <<кратный полиномиальный частный интеграл>>
нами стал использоваться, начиная со статьи [30].

Обобщенная задача Дарбу и методы ее решения получили развитие в работах
[1, 2, 12, 23, 25, 30 -- 56].
\vspace{1ex}

{\bf 12.10.}
{\sl Расширенная задача Дарбу}
для обобщенного уравнения Риккати-Абеля
\\[1.75ex]
\mbox{}\hfill                     % (12.11)
$
\displaystyle
\dfrac{dx}{dt}=
\sum\limits_{i=0}^n
a_i^{}(t)\;\!x^{n-i}
\quad
(n\geq 2)
$
\hfill (12.11)
\\[1.5ex]
с коэффициентами $a_i^{}\in C^1T,\ i=0,\ldots,n-1,\ a_{_0}^{}\!(t)\not\equiv 0.$
\vspace{0.5ex}

Пусть функции
\\[1.5ex]
\mbox{}\hfill                     % (12.12)
$
\displaystyle
{\rm g}_j^{}\colon (t,x)\to\, {\rm g}_j^{}(t,x)
\quad
\forall (t,x)\in\Omega,
\quad
j=1,\ldots, n,
$
\hfill (12.12)
\\[1.75ex]
являются частными интегралами с сомножителями $M_j^{},\ j=1,\ldots,n,$
соответственно на области $\Omega\subset\R^2$ уравнения (12.11).
\vspace{0.5ex}

Согласно определению 1.1 $\deg_x^{} M_j^{}\leq n-1,\ j=1,\ldots, n,$ т.е.
\\[1.5ex]
\mbox{}\hfill                     % (12.13)
$
\displaystyle
M_j^{}\colon (t,x)\to\ 
\sum\limits_{s=0}^{n-1}\alpha_{sj}^{}(t)\;\!x^{n-s-1}
\quad
\forall (t,x)\in\Xi^{\;\!\prime},
\quad
j=1,\ldots, n,
$
\hfill (12.13)
\\[2ex]
где коэффициенты $\alpha_{sj}^{}\in C^1T^{\;\!\prime},\ s=0,\ldots,n-1,\ j=1,\ldots, n.$
\vspace{0.5ex}

Множество $\Omega_{_0}^{}\subset\Omega$
такое, что 
${\rm g}_j^{{}^{\!\scriptsize\gamma_j^{}}}\in C^{1}\Omega_{_0}^{},\ 
\gamma_{j}^{}\in\R,\  j=1,\ldots, n.$
\vspace{0.5ex}

По свойству 11.11, функция
\\[1.75ex]
\mbox{}\hfill                   % (12.14)
$
\displaystyle
F\colon (t,x)\to\ 
\prod\limits_{j=1}^n
{\rm g}_j^{{}^{\!\scriptsize\gamma_j^{}}}(t,x)
\quad
\forall (t,x)\in \Omega_{_0}^{}
$
\hfill (12.14)
\\[1.5ex]
является общим интегралом 
на множестве $\Omega_{_0}^{}$ уравнения (12.11), 
если и только если 
\\[1.5ex]
\mbox{}\hfill                  
$
\displaystyle
\sum\limits_{j=1}^{n}\gamma_{j}^{}\,
\sum\limits_{s=0}^{n-1}\alpha_{sj}^{}(t)\;\!x^{n-s-1}=0
\quad
\forall (t,x)\in\Xi^{\;\!\prime}.
\hfill
$
\\[1.5ex]
\indent
Это тождество имеет место, если функциональный определитель
\\[1.75ex]
\mbox{}\hfill                  
$
\displaystyle
\triangle(t)=
\begin{vmatrix}
\alpha_{01}^{}(t) & \alpha_{02}^{}(t) & \dots & \alpha_{0n}^{}(t) 
\\[1.25ex]
\alpha_{11}^{}(t) & \alpha_{12}^{}(t) & \dots & \alpha_{1n}^{}(t) 
\\[0.75ex]
\hdotsfor{4} 
\\[0.75ex]
\alpha_{n-1{,}1}^{}(t) & \alpha_{n-1{,}2}^{}(t) & \dots & \alpha_{n-1{,}n}^{}(t) 
\end{vmatrix}
\quad
\forall t\in T^{\;\!\prime}
\hfill
$
\\[1.75ex]
тождественно равен нулю на области $T^{\;\!\prime}.$
\vspace{0.35ex}

Тем самым, доказана
\vspace{0.5ex}

{\bf Теорема 12.3.}
{\it
Функция {\rm (12.14)} является общим интегралом на множестве $\Omega_{_0}^{}$
уравнения {\rm (12.11)} тогда и только тогда, когда функции {\rm (12.12)} 
являются частными интегралами на области $\Omega$ уравнения {\rm (12.11)}
\vspace{0.25ex}
с такими сомножителями {\rm (12.13)}, что определитель 
$\triangle(t)=0\;\;\forall t\in T^{\;\!\prime}.$
}
\vspace{0.5ex}

По свойству 8.8, функция (12.14) является 
интегрирующим множителем уравнения (12.11), если и только если
\\[1.5ex]
\mbox{}\hfill                  
$
\displaystyle
\sum\limits_{j=1}^{n}\gamma_{j}^{}\,
\sum\limits_{s=0}^{n-1}\alpha_{sj}^{}(t)\;\!x^{n-s-1}={}- {\rm div}\,{\frak d}(t,x)
\quad
\forall (t,x)\in\Xi^{\;\!\prime}.
\hfill
$
\\[1.5ex]
\indent
Учитывая, что
\\[1.5ex]
\mbox{}\hfill                  
$
\displaystyle
{\rm div}\,{\frak d}(t,x)=
\sum\limits_{i=0}^{n-1}(n-i)\;\!a_{i}^{}(t)\;\!x^{n-i-1} 
\quad
\forall (t,x)\in\Xi,
\hfill
$
\\[1.5ex]
тождество выполняется тогда и только тогда, когда имеет место система тождеств
\\[1.5ex]
\mbox{}\hfill                  % (12.15)
$
\displaystyle
\sum\limits_{j=1}^{n}\alpha_{sj}^{}(t)\;\!\gamma_{j}^{}={}- 
(n-s)\;\!a_{s}^{}(t)
\quad
\forall t\in T^{\;\!\prime},
\quad
s=0,\ldots, n-1.
$
\hfill (12.15)
\\[1.5ex]
\indent
На основе определителя $\triangle$ составим $n$ функциональных 
определителей $\triangle_j^{}$ путем замены $j\!$-го столбца на столбец, 
состоящий из функций 
\\[1.5ex]
\mbox{}\hfill
$
{}-n\;\!a_{_0}^{}(t),{}-(n-1)\;\!a_1^{}(t),\ldots,{}-a_{n-1}^{}(t).
\hfill
$
\\[1.5ex]
\indent
Если функциональные определители $\triangle,\ \triangle_1^{},\ldots, \triangle_n^{}$ такие, что
\\[1.5ex]
\mbox{}\hfill                  % (12.16)
$
\displaystyle
\dfrac{\triangle_j^{}(t)}{\triangle(t)}\in\R,
\quad
j=1,\ldots, n,
$
\hfill (12.16)
\\[1.5ex]
то система тождеств (12.15) разрешима относительно $\gamma_1^{},\ldots,\gamma_n^{},$
а значит, найден интегрирующий множитель уравнения (12.11) в виде (12.14).
\vspace{0.35ex}

Итак, при $\triangle\not\equiv 0$ доказана
\vspace{0.5ex}

{\bf Теорема 12.4.}
{\it
Функция {\rm (12.14)} является интегрирующим множителем на множестве $\Omega_{_0}^{}$
уравнения {\rm (12.11)} тогда и только тогда, когда функции {\rm (12.12)} 
являются частными интегралами на области $\Omega$ уравнения {\rm (12.11)}
с такими сомножителями {\rm (12.13)}, что определитель $\triangle$
тождественно не равен нулю на области $T^{\;\!\prime}$ и 
выполняется условие {\rm (12.16)}.
}
\vspace{0.5ex}

В статье [10] для уравнения (12.11) расширенная задача Дарбу (теоремы 12.3 и 12.4)
решена при наличии $n$ интегральных кривых, определяемых 
полиномиальными частными интегралами, а при $n=3$ 
аналогичная задача решена в [11, с. 54 -- 61].
\vspace{1ex}

{\bf 12.11.}
Уравнение Абеля первого рода
\\[1.75ex]
\mbox{}\hfill                     % (12.17)
$
\displaystyle
\dfrac{dx}{dt}=
\dfrac{1}{t}\,x-\dfrac{1}{2\;\!t^2}\,x^3
$
\hfill (12.17)
\\[1.75ex]
будучи уравнением Бернулли, интегрируется классическими методами [3].
\vspace{0.15ex}

Построим общий интеграл уравнения (12.17) по его частным интегралам.
\vspace{0.15ex}

Полином 
\\[1ex]
\mbox{}\hfill
$
p\colon (t,x)\to\ {}-t+x^2
\quad
\forall (t,x)\in\R^2
\hfill
$
\\[1.75ex]
является частным интегралом с сомножителем 
\\[1.5ex]
\mbox{}\hfill
$
M\colon (t,x)\to\ 
\dfrac{1}{t}-\dfrac{1}{t^2}\,x^2
\quad
\forall (t,x)\in\Xi_{_0}^{}
\hfill
$
\\[1.75ex]
на множестве 
$\Xi_{_0}^{}=T_{_{\!0}}^{}\times\R,\ T_{_{\!0}}^{}=\R\backslash\{0\}.$
\vspace{0.75ex}

Поскольку сомножитель $M={}-\dfrac{1}{t^2}\,p,$
\vspace{0.75ex}
то по свойству 5.10 (при $k=1,\ M_0^{}={}-1/t^2)$
полиномиальный частный интеграл $p$ двукратный
\\[1.75ex]
\mbox{}\hfill
$
\bigl((p, p\;\!M_0^{}), (1, 1, {}-M_0^{})\bigr)\in \text{B}_{_{\Xi_{_0}}}.
\hfill
$ 
\\[2ex]
\indent
По свойству 5.13, функция
\\[1.75ex]
\mbox{}\hfill
$
{\rm g}\colon (t,x)\to\ 
\exp\dfrac{1}{x^2-t}
\quad
\forall (t,x)\in\Omega_{_0}^{},
\quad
\Omega_{_0}^{}=\{(t,x)\colon x^2-t\ne 0\},
\hfill
$
\\[2.25ex]
является экспоненциальным частным интегралом с сомножителем 
\\[1.5ex]
\mbox{}\hfill
$
N\colon t\to\ \dfrac{1}{t^2}
\quad
\forall t\in T_{_0}^{}
\hfill
$
\\[1.5ex]
на множестве $\Omega^{\;\!\prime}=\Omega_{_0}^{}\!\cap\Xi_{_0}^{}.$
\vspace{1ex}

По свойству 11.8 (при $\varphi=N),$ функция
\\[2ex]
\mbox{}\hfill
$
F\colon (t,x)\to\ 
\exp\dfrac{1}{t}\,
\exp\dfrac{1}{x^2-t}=
\exp\dfrac{x^2}{t\;\!(x^2-t)}
\quad
\forall (t,x)\in\Omega^{\;\!\prime},
\hfill
$
\\[2.25ex]
а с учетом свойства 0.1, и функция
\\[2ex]
\mbox{}\hfill
$
\Psi\colon (t,x)\to\ 
\dfrac{x^2}{t\;\!(x^2-t)}
\quad
\forall (t,x)\in\Omega^{\;\!\prime}
\hfill
$
\\[2.25ex]
является общим интегралом на множестве $\Omega^{\;\!\prime}$ уравнения (12.17).
\vspace{1.25ex}

{\bf 12.12.}
{\sl Обратная задача}:
построение дифференциальных систем по частным интег\-ралам.

Пусть $({\rm g}_j^{}, M_j^{})\in \text{J}_{_{\Omega}},\ j=1,\ldots, n.$
\vspace{0.75ex}
Тогда, по критерию существования частного интеграла (теорема 1.1),
\\[2ex]
\mbox{}\hfill                     % (12.18)
$
\displaystyle
\sum\limits_{i=1}^n\, X_{i}^{}(t,x)\, 
\partial_{{}_{\scriptstyle x_i^{}}}{\rm g}_j^{}(t,x)=
{\rm g}_j^{}(t,x)\;\!M_j^{}(t,x)-
\partial_{{}_{\scriptstyle t}}{\rm g}_j^{}(t,x)
\quad 
\forall (t,x)\in \Omega,
\quad
j=1,\ldots, n.
$
\hfill (12.18)
\\[2.25ex]
\indent
Пусть якобиан частных интегралов ${\rm g}_1^{},\ldots, {\rm g}_n^{}$
по переменным $x_1^{},\ldots, x_n^{}$
\\[2ex]
\mbox{}\hfill                     % (12.19)
$
\displaystyle
\triangle(t,x)=\dfrac{{\rm D}({\rm g}_1^{},\ldots, {\rm g}_n^{})}{{\rm D}(x_1^{},\ldots, x_n^{})}
\not\equiv 0
\ \ \text{на}\ \ \Omega.
$
\hfill (12.19)
\\[2ex]
\indent
На основании якобиана $\triangle$ составим $n$ 
\vspace{0.35ex}
функциональных определителей $\triangle_j^{}$ 
путем замены $j\!$-го столбца на столбец, состоящий из функций
\\[1.5ex]
\mbox{}\hfill                   
$
{\rm g}_1^{}\;\!M_1^{}- 
\partial_{{}_{\scriptstyle t}}\;\!{\rm g}_1^{},\ 
{\rm g}_2^{}\;\!M_2^{}- 
\partial_{{}_{\scriptstyle t}}\;\!{\rm g}_2^{}\;\!,\ \ldots\ ,\
{\rm g}_n^{}\;\!M_n^{}- 
\partial_{{}_{\scriptstyle t}}\;\!{\rm g}_n^{}. 
\hfill
$
\\[2ex]
\indent
Из системы тождеств (12.18) находим
\\[2ex]
\mbox{}\hfill               
$
X_i^{}(t,x)=\dfrac{\triangle_i^{}(t,x)}{\triangle(t,x)}
\quad
\forall (t,x)\in \Omega_{_0}\subset\Omega,
\quad
i=1,\ldots,n.
\hfill
$
\\[1.75ex]
\indent
Если
\\[1.75ex]
\mbox{}\hfill                     % (12.20)
$
\dfrac{\triangle_i^{}(t,x)}{\triangle(t,x)}\in 
\text{P}_{_{\Xi_{_0}}},
\quad
i=1,\ldots,n.
\quad
\Xi_{_0}\subset \Xi^{\;\!\prime},
$
\hfill (12.20)
\\[1.75ex]
то справедлива
\vspace{1ex}

{\bf Теорема 12.5.}
\vspace{0.75ex}
{\it
Если $({\rm g}_j^{}, M_j^{})\in \text{\rm J}_{_{\Omega}},\ j=1,\ldots, n,$
такие, что выполняются условия {\rm (12.19)} и {\rm (12.20)},
то система {\rm (0.1)} имеет вид}
\\[1.75ex]
\mbox{}\hfill               
$
\dfrac{dx_i^{}}{dt}=\dfrac{\triangle_i^{}(t,x)}{\triangle(t,x)}\,,
\quad
i=1,\ldots,n.
\hfill
$
\\[2.25ex]
\indent
{\bf Следствие 12.1.}
\vspace{0.5ex}
{\it
Если автономная дифференциальная система {\rm (12.10)}
имеет $n$ таких функционально независимых частных интегралов
\vspace{0.5ex}
${\rm g}_j^{}\colon x\to {\rm g}_j^{}(x)\;\;\forall x\in X^{\;\!\prime}$ 
с сом\-ножителями
$M_j^{}\colon\R^n\to\R,\ j=1,\ldots, n,$ соответственно, 
\vspace{0.5ex}
что выполняются условия
$\dfrac{\triangle_i^{}(x)}{\triangle(x)}\in \text{\rm P}_{_{\scriptstyle\!\R^n}},\
i=1,\ldots,n,$
то она имеет вид}
\\[1.5ex]
\mbox{}\hfill               
$
\dfrac{dx_i^{}}{dt}=\dfrac{\triangle_i^{}(x)}{\triangle(x)}\,,
\quad
i=1,\ldots,n.
\hfill
$
\\[2ex]
\indent
{\bf Замечание 12.1.}
\vspace{0.5ex}
Если частный интеграл ${\rm g}_j^{}$ с сомножителем $M_j^{},\ j\in\{1,\ldots,n\},$
является экспоненциальным 
$(\exp\omega_j^{},M_j^{})\in \text{\rm E}_{_{\Omega}}\;\!,$
\vspace{0.75ex}
то (с учетом теоремы 3.1) в $j\!$-ом тождестве из системы (12.18)
и в определителях $\triangle,\, \triangle_1^{},\ldots, \triangle_n^{}$
\vspace{0.35ex}
из теоремы 2.5 и из следст\-вия~12.1 целесообразно произвести формальную замену 
функций ${\rm g}_j^{}$ на функции $\omega_j^{},$
\vspace{0.5ex}
а произведения ${\rm g}_j^{}\;\!M_j^{}$ на функцию $M_j^{}\;\!.$ 
\vspace{1.5ex}

{\bf 12.13.}
Пусть система (12.2) имеет полиномиальный частный интеграл
\\[1.5ex]
\mbox{}\hfill
$
p\colon (x,y)\to\ 
x^2+y^2+a
\quad
\forall (x,y)\in\R^2
\quad
(a\in\R)
\hfill
$
\\[1.5ex]
с сомножителем
\\[1ex]
\mbox{}\hfill
$
M_1^{}\colon (x,y)\to\ 2\;\!(x+y)
\quad
\forall (x,y)\in\R^2
\hfill
$
\\[1.5ex]
и условный частный интеграл
\\[1.5ex]
\mbox{}\hfill
$
{\rm g}\colon (x,y)\to\ 
\exp\omega(x,y)=\exp(x-y)
\quad
\forall (x,y)\in\R^2
\hfill
$
\\[1.5ex]
с сомножителем
\\[1ex]
\mbox{}\hfill
$
M_2^{}\colon (x,y)\to\ {}-(x+y)
\quad
\forall (x,y)\in\R^2.
\hfill
$
\\[1.5ex]
\indent
Якобиан
\\[1.75ex]
\mbox{}\hfill                  
$
\displaystyle
\triangle(x,y)=
\left|\!\!
\begin{array}{cc}
\partial_{x}^{}p & \partial_{y}^{}p
\\[1.25ex]
\partial_{x}^{}\omega & \partial_{y}^{}\omega
\end{array}
\!\!\right|
=
\left|\!\!
\begin{array}{cr}
2\;\!x & 2\;\!y
\\[1ex]
1 & {}-1
\end{array}
\!\!\right|
={}-2\;\!(x+y)
\quad
\forall (x,y)\in\R^2
\hfill
$
\\[2.25ex]
отличен от тождественного нуля на $\R^2$ 
\vspace{0.35ex}
(функции $p$ и $\omega$ функционально независимы на $\R^2).$
Определители на $\R^2$
\\[1.75ex]
\mbox{}\hfill                  
$
\displaystyle
\triangle_{{}_{X}}(x,y)=
\left|\!\!
\begin{array}{cc}
p\;\!M_1^{} & \partial_{y}^{}p
\\[1.25ex]
M_2^{} & \partial_{y}^{}\omega
\end{array}
\!\!\right|
=
\left|\!\!
\begin{array}{cr}
2\;\!(x^2+y^2+a)(x+y) & 2\;\!y
\\[1ex]
{}-(x+y) & {}-1
\end{array}
\!\!\right|
= 2\;\!(x+y)(y-x^2-y^2-a),
%\quad
%\forall (x,y)\in\R^2
\hfill
$
\\[2ex]
и
\\[1.75ex]
\mbox{}\hfill                  
$
\displaystyle
\triangle_{{}_{Y}}(x,y)=
\left|\!\!
\begin{array}{cc}
\partial_{x}^{}p & p\;\!M_1^{} 
\\[1.25ex]
\partial_{x}^{}\omega & M_2^{}  
\end{array}
\!\!\right|
=
\left|\!\!
\begin{array}{cc}
2\;\!x & 2\;\!(x^2+y^2+a)
\\[1ex]
1 & {}-(x+y)
\end{array}
\!\!\right|
= {}-2\;\!(x+y)(x+x^2+y^2+a).
%\quad
%\forall (x,y)\in\R^2
\hfill
$
\\[2.25ex]
\indent
По следствию 12.1 с учетом замечания 12.1, 
\\[2.25ex]
\mbox{}\hfill               
$
\dfrac{dx}{dt}=\dfrac{\triangle_{{}_{X}}^{}(x,y)}{\triangle(x,y)}=
{}-y+x^2+y^2+a,
\hfill
$
\\[0.5ex]
\mbox{}\hfill (12.21)
\\[0.5ex]
\mbox{}\hfill               
$
\dfrac{dy}{dt}=\dfrac{\triangle_{{}_{Y}}^{}(x,y)}{\triangle(x,y)}=
x+x^2+y^2+a.\ \ \ \,
\hfill
$
\\[2.25ex]
\indent
Согласно следствию 11.5 
\vspace{0.35ex}
(при $\rho_1^{}=2,\ \rho_2^{}={}-1,\ M_0^{}=x+y,\ 
\gamma_1^{}=1,\ \gamma_2^{}=2),$
трансцендентная функция 
\\[2ex]
\mbox{}\hfill
$
F\colon (x,y)\to\ (x^2+y^2+a)\exp\bigl(\;\!2\;\!(x-y)\bigr)
\quad
\forall (x,y)\in\R^2
\hfill
$
\\[2ex]
является автономным первым интегралом на $\R^2$ системы (12.21).
\vspace{0.35ex}

Обратим внимание на то, что при $a>0$ полином 
\vspace{0.35ex}
$p(x,y)\ne 0\;\;\forall (x,y)\in\R^2.$
Иначе говоря, при $a>0$ полиномиальный частный интеграл $p$ 
\vspace{0.15ex}
не определяет траекторий дифференциальной системы (12.21).
\vspace{0.35ex}

Таким образом, при $a>0$ 
\vspace{0.15ex}
дифференциальная система (12.21) построена на основании
двух частных интегралов 
\vspace{0.35ex}
$p$ и ${\rm g},$ каждый из которых не определяет траекторий.

В статье [22] (см. также [4, с. 480 -- 481])
\vspace{0.15ex}
построена дифференциальная система (12.2) по частным интегралам при условии, 
что они определяют траектории.
\vspace{1ex}

{\bf 12.14.}
\vspace{0.5ex}
Пусть система (12.2) на фазовой плоскости $\R^2$ имеет кратный 
полиномиальный частный интеграл 
$\bigl((p, M), (h, q, N)\bigr)\in \text{\rm B}_{_{\scriptstyle\R^2}}.$
По теореме 5.2, 
\\[2ex]
\mbox{}\hfill
$
X(x,y)\;\!\partial_{x}^{}p(x,y)+
Y(x,y)\;\!\partial_{y}^{}p(x,y)=p(x,y)\;\!M(x,y)
\quad
\forall (x,y)\in\R^2,
\hfill
$
\\[2.75ex]
\mbox{}\hfill
$
X(x,y)\;\!\partial_{x}^{}q(x,y)+
Y(x,y)\;\!\partial_{y}^{}q(x,y)=h\;\!q(x,y)\;\!M(x,y)+p^{\;\!h}(x,y)\;\!N(x,y)
\quad
\forall (x,y)\in\R^2.
\hfill
$
\\[2ex]
\indent
Считая взаимно простые полиномы $p$ и $q$ функционально независимыми на $\R^2,$
из системы тождеств находим
\\[1.5ex]
\mbox{}\hfill               
$
X(x,y)=\dfrac{\triangle_{{}_{X}}^{}(x,y)}{\triangle(x,y)}\,,
\quad
Y(x,y)=\dfrac{\triangle_{{}_{Y}}^{}(x,y)}{\triangle(x,y)}
\quad
\forall (x,y)\in D\subset\R^2,
\hfill
$
\\[1.75ex]
где якобиан полиномов $p$ и $q$ по $x$ и $y$
\\[1.75ex]
\mbox{}\hfill                  % (12.22)
$
\displaystyle
\triangle(x,y)=
\left|\!\!
\begin{array}{cc}
\partial_{x}^{}p & \partial_{y}^{}p
\\[1.25ex]
\partial_{x}^{}q & \partial_{y}^{}q
\end{array}
\!\!\right|
\not\equiv 0
\ \ \text{на} \ \ \R^2,
$
\hfill (12.22)
\\[1.5ex]
определители
\\[1.5ex]
\mbox{}\hfill                 
$
\displaystyle
\triangle_{{}_{X}}^{}(x,y)=
\left|\!\!
\begin{array}{cc}
p\;\!M & \partial_{y}^{}p
\\[1.25ex]
h\;\!q\;\!M+p^{\;\!h}N & \partial_{y}^{}q
\end{array}
\!\!\right|,
\quad \ \
\triangle_{{}_{Y}}^{}(x,y)=
\left|\!\!
\begin{array}{cc}
\partial_{x}^{}p & p\;\!M 
\\[1.25ex]
\partial_{x}^{}q & h\;\!q\;\!M+p^{\;\!h}N
\end{array}
\!\!\right|
\quad
\forall (x,y)\in\R^2.
\hfill
$
\\[1.75ex]
\indent
Если частные
\\[1.5ex]
\mbox{}\hfill                  % (12.23)
$
\displaystyle
\dfrac{\triangle_{{}_{X}}^{}(x,y)}{\triangle(x,y)}\in 
\text{\rm P}_{_{\!\scriptstyle\R^2}}
\quad \
\text{и} \ \,
\quad
\dfrac{\triangle_{{}_{Y}}^{}(x,y)}{\triangle(x,y)}\in 
\text{\rm P}_{_{\!\scriptstyle\R^2}},
$
\hfill (12.23)
\\[1.5ex]
то справедлива
\vspace{0.5ex}

{\bf Теорема 12.6.}
\vspace{0.5ex}
{\it
Если дифференциальная  система {\rm (12.2)} на $\R^2$ имеет такой кратный  
полиномиальный частный интеграл 
$\bigl((p, M), (h, q, N)\bigr)\in \text{\rm B}_{_{\scriptstyle\R^2}},$
\vspace{0.5ex}
что выполняются условия {\rm (12.22)} и {\rm (12.23)},
то она имеет вид}
\\[1.5ex]
\mbox{}\hfill               
$
\dfrac{dx}{dt}=
\dfrac{\triangle_{{}_{X}}^{}(x,y)}{\triangle(x,y)}\,,
\qquad
\dfrac{dy}{dt}=
\dfrac{\triangle_{{}_{Y}}^{}(x,y)}{\triangle(x,y)}\,.
\hfill
$
\\[2.25ex]
\indent
{\bf 12.15.}
\vspace{0.5ex}
Пусть система (12.2) на фазовой плоскости $\R^2$ имеет кратный 
полиномиальный частный интеграл 
$\bigl((p, M), (h, q, N)\bigr)\in \text{\rm B}_{_{\scriptstyle\R^2}}$ такой, что
\\[2ex]
\mbox{}\hfill
$
p=x^2+y^2-1,
\quad
M=xy,
\quad
h=1,
\quad
q=x^2-y^2-1,
\quad
N={}-xy.
\hfill
$
\\[1.5ex]
\indent
Якобиан
\\[1.5ex]
\mbox{}\hfill                
$
\displaystyle
\triangle(x,y)=
\left|\!\!
\begin{array}{cc}
\partial_{x}^{}p & \partial_{y}^{}p
\\[1.25ex]
\partial_{x}^{}q & \partial_{y}^{}q
\end{array}
\!\!\right|
=
\left|\!\!
\begin{array}{cr}
2\;\!x & 2\;\!y
\\[1ex]
2\;\!x & {}-2\;\!y
\end{array}
\!\!\right|=
{}-8\;\!xy
\quad
\forall (x,y)\in \R^2
\hfill
$
\\[2ex]
отличен от тождественного нуля на $\R^2.$
\vspace{0.5ex}

Определители на $\R^2$
\\[2ex]
\mbox{}\hfill                 
$
\displaystyle
\triangle_{{}_{X}}^{}(x,y)=
\left|\!\!
\begin{array}{cc}
p\;\!M & \partial_{y}^{}p
\\[1.25ex]
h\;\!q\;\!M+p^{\;\!h}N & \partial_{y}^{}q
\end{array}
\!\!\right|
=
\left|\!\!
\begin{array}{cr}
(x^2+y^2-1)\;\!xy & 2\;\!y
\\[1ex]
{}-2\;\!xy^3 & {}-2\;\!y
\end{array}
\!\!\right|=
{}-2\;\!xy^2(x^2-y^2-1),
%\quad
%\forall (x,y)\in \R^2
\hfill
$
\\[3.5ex]
\mbox{}\hfill
$
\triangle_{{}_{Y}}^{}(x,y)=
\left|\!\!
\begin{array}{cc}
\partial_{x}^{}p & p\;\!M 
\\[1.25ex]
\partial_{x}^{}q & h\;\!q\;\!M+p^{\;\!h}N
\end{array}
\!\!\right|
=
\left|\!\!
\begin{array}{cc}
2\;\!x & (x^2+y^2-1)\;\!xy 
\\[1ex]
2\;\!x & {}-2\;\!xy^3  
\end{array}
\!\!\right|=
{}-2\;\!x^2y\;\!(x^2+3y^2-1).
%\quad
%\forall (x,y)\in\R^2.
\hfill
$
\\[2.25ex]
\indent
По теореме 12.6,
\\[2.25ex]
\mbox{}\hfill               % (12.24)
$
\dfrac{dx}{dt}=\dfrac{\triangle_{{}_{X}}^{}(x,y)}{\triangle(x,y)}=
\dfrac{1}{4}\,y\;\!(x^2-y^2-1),
\quad
\dfrac{dy}{dt}=\dfrac{\triangle_{{}_{Y}}^{}(x,y)}{\triangle(x,y)}=
\dfrac{1}{4}\,x\;\!(x^2+3y^2-1).
$
\hfill (12.24)
\\[2.25ex]
\indent
Согласно теореме 5.3 система (12.24) имеет два частных интеграла
\\[1.5ex]
\mbox{}\hfill
$
(x^2+y^2-1,\;\! xy)\in \text{\rm A}_{_{\scriptstyle\R^2}}
\quad \text{и}\quad
\Bigl(\exp\dfrac{x^2-y^2-1}{x^2+y^2-1}\,,\;\! {}-xy\Bigr)\in \text{\rm E}_{_{D}},
\hfill
$
\\[1.5ex]
где множество $D=\{(x,y)\colon x^2+y^2-1\ne 0\}.$
\vspace{0.75ex}

По следствию 11.5 (при $\rho_1^{}=1,\ \rho_2^{}={}-1,\ M_0^{}=xy,\ 
\gamma_1^{}=\gamma_2^{}=1),$ трансцендентная функция
\\[1.5ex]
\mbox{}\hfill
$
F\colon (x,y)\to\ 
(x^2+y^2-1)\exp\dfrac{x^2-y^2-1}{x^2+y^2-1}
\quad
\forall (x,y)\in D
\hfill
$
\\[1.5ex]
является автономным первых интегралом на множестве $D$ системы (12.24).
\vspace{1ex}

{\bf 12.16.}
\vspace{0.35ex}
Пусть система (12.2) на фазовой плоскости $\R^2$ имеет комплекснозначный  
полиномиальный частный интеграл 
$(u+i\;\!v, U+i\;\!V)\in \text{\rm H}_{_{\scriptstyle\R^2}}.$
По теореме 6.1, 
\\[1.5ex]
\mbox{}\hfill
$
X(x,y)\;\!\partial_{x}^{}u(x,y)+
Y(x,y)\;\!\partial_{y}^{}u(x,y)=u(x,y)\;\!U(x,y)-v(x,y)\;\!V(x,y)
\quad
\forall (x,y)\in\R^2,
\hfill
$
\\[2ex]
\mbox{}\hfill
$
X(x,y)\;\!\partial_{x}^{}v(x,y)+
Y(x,y)\;\!\partial_{y}^{}v(x,y)=
u(x,y)\;\!V(x,y)+v(x,y)\;\!U(x,y)
\quad
\forall (x,y)\in\R^2.
\hfill
$
\\[1.5ex]
\indent
Считая полиномы $u$ и $v$ функционально независимыми на $\R^2,$
из системы тождеств находим
\\[1.5ex]
\mbox{}\hfill               
$
X(x,y)=\dfrac{\triangle_{{}_{X}}^{}(x,y)}{\triangle(x,y)}\,,
\quad
Y(x,y)=\dfrac{\triangle_{{}_{Y}}^{}(x,y)}{\triangle(x,y)}
\quad
\forall (x,y)\in D\subset\R^2,
\hfill
$
\\[1.5ex]
где якобиан полиномов $u$ и $v$ по $x$ и $y$
\\[1.5ex]
\mbox{}\hfill                  % (12.25)
$
\displaystyle
\triangle(x,y)=
\left|\!\!
\begin{array}{cc}
\partial_{x}^{}u & \partial_{y}^{}u
\\[1.25ex]
\partial_{x}^{}v & \partial_{y}^{}v
\end{array}
\!\!\right|
\not\equiv 0
\ \ \text{на} \ \ \R^2,
$
\hfill (12.25)
\\[1.5ex]
определители
\\[1.5ex]
\mbox{}\hfill                 
$
\displaystyle
\triangle_{{}_{X}}^{}(x,y)=
\left|\!\!
\begin{array}{cc}
u\;\!U-v\;\!V & \partial_{y}^{}u
\\[1.25ex]
u\;\!V+v\;\!U  & \partial_{y}^{}v
\end{array}
\!\!\right|,
\quad \ \
\triangle_{{}_{Y}}^{}(x,y)=
\left|\!\!
\begin{array}{cc}
\partial_{x}^{}u & u\;\!U-v\;\!V 
\\[1.25ex]
\partial_{x}^{}v & u\;\!V+v\;\!U
\end{array}
\!\!\right|
\quad
\forall (x,y)\in\R^2.
\hfill
$
\\[1.5ex]
\indent
Если частные
\\[1.5ex]
\mbox{}\hfill                  % (12.26)
$
\displaystyle
\dfrac{\triangle_{{}_{X}}^{}(x,y)}{\triangle(x,y)}\in 
\text{\rm P}_{_{\!\scriptstyle\R^2}}
\quad \
\text{и} \ \,
\quad
\dfrac{\triangle_{{}_{Y}}^{}(x,y)}{\triangle(x,y)}\in 
\text{\rm P}_{_{\!\scriptstyle\R^2}},
$
\hfill (12.26)
\\[1.5ex]
то справедлива
\vspace{0.5ex}

{\bf Теорема 12.7.}
\vspace{0.35ex}
{\it
Если дифференциальная  система {\rm (12.2)} на $\R^2$ имеет такой 
комплекснозначный полиномиальный частный интеграл 
$(u+i\;\!v, U+i\;\!V)\in \text{\rm H}_{_{\scriptstyle\R^2}},$
что выполняются условия {\rm (12.25)} и {\rm (12.26)},
то она имеет вид}
\\[1.5ex]
\mbox{}\hfill               
$
\dfrac{dx}{dt}=
\dfrac{\triangle_{{}_{X}}^{}(x,y)}{\triangle(x,y)}\,,
\qquad
\dfrac{dy}{dt}=
\dfrac{\triangle_{{}_{Y}}^{}(x,y)}{\triangle(x,y)}\,.
\hfill
$
\\[2ex]
\indent
{\bf Следствие 12.2.}
\vspace{0.35ex}
{\it
Система {\rm (12.2)} имеет  
комплекснозначный полиномиальный частный интеграл 
$(x+i\;\!y, U+i\;\!V)\in \text{\rm H}_{_{\scriptstyle\R^2}}$
тогда и только тогда, когда она имеет вид}
\\[1.5ex]
\mbox{}\hfill               
$
\dfrac{dx}{dt}=
x\;\!U(x,y)-y\;\!V(x,y), 
\qquad
\dfrac{dy}{dt}=
x\;\!V(x,y)+y\;\!U(x,y).
\hfill
$
\\[2ex]
\indent
{\bf 12.17.}
\vspace{0.35ex}
Пусть система (12.2) на фазовой плоскости $\R^2$ имеет 
комплекснозначный полиномиальный частный интеграл 
$(x+i\;\!y, 1+i)\in \text{\rm H}_{_{\scriptstyle\R^2}}.$ По следствию 12.2,
\\[1.5ex]
\mbox{}\hfill     % (12.27)
$
\dfrac{dx}{dt}=x-y,
\qquad
\dfrac{dy}{dt}=x+y.
$
\hfill (12.27)
\\[1.5ex]
\indent
Согласно теореме 6.2 система (12.27) имеет два частных интеграла
\\[1.5ex]
\mbox{}\hfill
$
(x^2+y^2,\;\! 2)\in \text{\rm A}_{_{\scriptstyle\R^2}}
\quad \text{и}\quad
\Bigl(\exp\arctg\dfrac{y}{x}\,,\;\! 1\Bigr)\in \text{\rm E}_{_{D}},
\hfill
$
\\[1.5ex]
где множество $D=\{(x,y)\colon x\ne 0\}.$
\vspace{0.5ex}

По следствию 11.5 
\vspace{0.15ex}
(при $\rho_1^{}=2,\ \rho_2^{}=1,\ M_0^{}=1,\ 
\gamma_1^{}=1,\ \gamma_2^{}={}-2)$ 
с последующим учетом свойства 3.13, функция
\\[1.5ex]
\mbox{}\hfill
$
F\colon (x,y)\to\ 
(x^2+y^2)\exp\Bigl(2\arctg\dfrac{y}{x}\Bigr)
\quad
\forall (x,y)\in D
\hfill
$
\\[1.5ex]
является автономным первых интегралом на множестве $D$ системы (12.27).
\vspace{0.35ex}

На основании частных интегралов, используя следствие 11.2 (при $\lambda=2$ и $\lambda=1),$
строим неавтономные первые интегралы системы (12.27):
\\[1.5ex]
\mbox{}\hfill
$
F_1^{}\colon (t,x,y)\to\ 
e^{{}-2\;\!t}(x^2+y^2)
\quad
\forall (t,x,y)\in\R^3 
\hfill
$
\\[1ex]
и
\\[1.25ex]
\mbox{}\hfill
$
F_2^{}\colon (t,x,y)\to\ 
e^{{}-t}\exp\arctg\dfrac{y}{x}=
\exp\Bigl({}-t+\arctg\dfrac{y}{x}\Bigr)
\quad
\forall (t,x,y)\in\Omega_{{}_0}= 
\{(t,x,y)\colon x\ne 0\}.
\hfill
$
\\[1.5ex]
\indent
Учитывая функциональную неоднозначность первого интеграла (свойство 0.1), 
первый интеграл $F_2^{}$ представим в виде
\\[1.5ex]
\mbox{}\hfill
$
\Psi\colon (t,x,y)\to\ 
{}-t+\arctg\dfrac{y}{x}
\quad
\forall (t,x,y)\in\Omega_{{}_0}.
\hfill
$
\\[1.5ex]
\indent
Поскольку первые интегралы $F,\ F_1^{},\ \Psi$
\vspace{0.25ex}
попарно функционально независимы на $\R^3,$
то каждая из совокупностей $\{F, F_1^{}\}, \ \{F, \Psi\},\ \{F_1^{}, \Psi\}$
\vspace{0.25ex}
образует интегральный базис на множестве $\Omega_{{}_0}$
дифференциальной системы (12.27).
\vspace{1ex}

{\bf 12.18.}
По следствию 12.2, система (12.2) имеет 
\\[1.75ex]
\mbox{}\hfill
$
\bigl(x+i\;\!y,\;\! x-y+i\;\!(x+y)\bigr)\in \text{\rm H}_{_{\scriptstyle\R^2}},
\hfill
$
\\[1.5ex]
если и только если 
\\[2ex]
\mbox{}\hfill     % (12.28)
$
\dfrac{dx}{dt}=x^2-2\;\!xy-y^2,
\qquad
\dfrac{dy}{dt}=x^2+2\;\!xy-y^2.
$
\hfill (12.28)
\\[2ex]
\indent
Согласно теореме 6.2 система (12.28) имеет
\\[1.75ex]
\mbox{}\hfill
$
\bigl(x^2+y^2,\;\! 2\;\!(x-y)\bigr)\in \text{\rm A}_{_{\scriptstyle\R^2}},
\quad 
\Bigl(\exp\arctg\dfrac{y}{x}\,,\;\! x+y\Bigr)\in \text{\rm E}_{_{D}},
\quad 
D=\{(x,y)\colon x\ne 0\}.
\hfill
$
\\[1.75ex]
\indent
На основании теоремы 1.1 устанавливаем, 
\vspace{0.35ex}
что дифференциальная система второго порядка (12.28)
имеет 
$
\bigl(x+y,\;\! 2\;\!(x-y)\bigr)\in \text{\rm A}_{_{\scriptstyle\R^2}}.
$
\vspace{0.75ex}

У полиномиальных частных интегралов один и тот же сомножитель.
По следст\-вию~11.7, рациональная функция
\\[1.5ex]
\mbox{}\hfill
$
F\colon (x,y)\to\ 
\dfrac{x+y}{x^2+y^2}
\quad
\forall (x,y)\in\R^2\backslash\{(0,0)\} 
\hfill
$
\\[1.5ex]
является автономным первым интегралов на области $\R^2\backslash\{(0,0)\}$
системы (12.28).
\vspace{1.5ex}

{\bf 12.19.}
У системы ([18, с. 88; 30] при $\lambda={}-1,\ \eta=1)$ 
\\[2.25ex]
\mbox{}\hfill     % (12.29)
$
\dfrac{dx}{dt}=x\;\!(x^2+y^2+\lambda)-y\;\!(x^2+y^2+\eta),
\quad
\dfrac{dy}{dt}=x\;\!(x^2+y^2+\eta)+y\;\!(x^2+y^2+\lambda),
$
\hfill (12.29)
\\[2.5ex]
где $\lambda,\eta\in\R, \ |\lambda|+|\eta|\ne 0,\ \lambda\ne\eta,$
\vspace{0.5ex}
по виду правых частей на основании следствия 12.2, находим
\\[1.25ex]
\mbox{}\hfill
$
\bigl(x+i\;\!y,\;\! x^2+y^2+\lambda+i\;\!(x^2+y^2+\eta)\bigr)\in \text{\rm H}_{_{\scriptstyle\R^2}}.
\hfill
$
\\[1.75ex]
\indent
В соответствии с теоремой 5.3
\\[1.75ex]
\mbox{}\hfill
$
\bigl(x^2+y^2,\;\! 2\;\!(x^2+y^2+\lambda)\bigr)\in \text{\rm A}_{_{\scriptstyle\R^2}},
\quad 
\Bigl(\exp\arctg\dfrac{y}{x}\,,\;\! x^2+y^2+\eta\Bigr)\in \text{\rm E}_{_{D}},
\quad 
D=\{(x,y)\colon x\ne 0\}.
\hfill
$
\\[1.75ex]
\indent
{\sl Случай} $\lambda\ne 0.$
Поскольку
\\[1.5ex]
\mbox{}\hfill
$
\bigl(x^2+y^2,\;\! 2\;\!(x^2+y^2+\lambda)\bigr)\in \text{\rm A}_{_{\scriptstyle\R^2}},
\hfill
$
\\[1.5ex]
то, по свойству 1.4, 
\\[1.5ex]
\mbox{}\hfill
$
\bigl(x^2+y^2+\lambda,\;\! 2\;\!(x^2+y^2)\bigr)\in \text{\rm A}_{_{\scriptstyle\R^2}}.
\hfill
$
\\[1.5ex]
\indent
Тождество
\\[1.5ex]
\mbox{}\hfill
$
2\;\!\gamma_1^{}\;\!(x^2+y^2+\lambda)+
\gamma_2^{}\;\!(x^2+y^2+\eta)+
2\;\!\gamma_3^{}\;\!(x^2+y^2)=0
\quad 
\forall (x,y)\in\R^2
\hfill
$
\\[2ex]
выполняется, например, при $\gamma_1^{}={}-\eta,\ \gamma_2^{}=2\;\!\lambda,\ 
\gamma_3^{}=\eta-\lambda.$
По свойству 11.11, функция
\\[2.25ex]
\mbox{}\hfill
$
F\colon (x,y)\to\ 
\dfrac{(x^2+y^2+\lambda)^{\eta-\lambda}}{(x^2+y^2)^{\eta}}\,
\exp\Bigl(2\;\!\lambda\;\!\arctg\dfrac{y}{x}\;\!\Bigr)
\quad
\forall (x,y)\in D
\hfill
$
\\[2ex]
является автономным первых интегралом на множестве $D$ системы (12.29) при $\lambda\ne 0.$
\vspace{0.35ex}

Разность сомножителей
\\[1.5ex]
\mbox{}\hfill
$
2\;\!(x^2+y^2)-2\;\!(x^2+y^2+\lambda)={}-2\;\!\lambda
\quad 
\forall (x,y)\in\R^2.
\hfill
$
\\[1.5ex]
\indent
Согласно свойству 1.9
\\[1.5ex]
\mbox{}\hfill
$
\Bigl(\dfrac{x^2+y^2+\lambda}{x^2+y^2}\,,{}-2\;\!\lambda\Bigr)\in 
\text{\rm J}_{_{\!\scriptstyle D_{_0}}},
\quad 
D_{_0}=\R^2\backslash\{(0,0)\},
\hfill
$
\\[1.5ex]
а значит (следствие 11.2), функция
\\[1.5ex]
\mbox{}\hfill
$
F_1^{}\colon (t,x,y)\to\ 
\dfrac{x^2+y^2+\lambda}{x^2+y^2}\,e^{\;\!2\;\!\lambda\;\!t}
\quad
\forall (t,x,y)\in\Omega_{_{0}}
\hfill
$
\\[1.5ex]
является неавтономным первых интегралом на области 
$\Omega_{_{0}}=\R\times D_{_{0}}$ 
дифференциальной системы (12.29) при $\lambda\ne 0.$
\vspace{0.35ex}

Линейная комбинация сомножителей
\\[1.5ex]
\mbox{}\hfill
$
2\;\!(x^2+y^2)-2\;\!(x^2+y^2+\eta)={}-2\;\!\eta
\quad 
\forall (x,y)\in\R^2.
\hfill
$
\\[1.75ex]
\indent
Согласно свойству 1.9
\\[1.5ex]
\mbox{}\hfill
$
\Bigl((x^2+y^2+\lambda)\;\!\exp\Bigl({}-2\;\!\arctg\dfrac{y}{x}\;\!\Bigr),{}-2\;\!\eta\Bigr)\in 
\text{\rm J}_{_{\!\scriptstyle D}},
\hfill
$
\\[1.5ex]
а значит (следствие 11.2), функция
\\[1.5ex]
\mbox{}\hfill
$
F_2^{}\colon (t,x,y)\to\ 
(x^2+y^2+\lambda)\;\!\exp\Bigl(2\;\!\eta\;\!t-2\;\!\arctg\dfrac{y}{x}\;\!\Bigr)
\quad
\forall (t,x,y)\in\Omega^{\;\!\prime}
\hfill
$
\\[1.5ex]
является неавтономным первых интегралом на множестве 
$\Omega^{\;\!\prime}=\R\times D$ 
дифференциальной системы (12.29) при $\lambda\ne 0.$
\vspace{0.35ex}

Линейная комбинация сомножителей
\\[1.5ex]
\mbox{}\hfill
$
2\;\!(x^2+y^2+\lambda)-2\;\!(x^2+y^2+\eta)=2\;\!(\lambda-\eta)
\quad 
\forall (x,y)\in\R^2.
\hfill
$
\\[1.75ex]
\indent
Согласно свойству 1.9
\\[1.5ex]
\mbox{}\hfill
$
\Bigl((x^2+y^2)\;\!\exp\Bigl({}-2\;\!\arctg\dfrac{y}{x}\;\!\Bigr),\, 2\;\!(\lambda-\eta)\Bigr)\in 
\text{\rm J}_{_{\!\scriptstyle D}},
\hfill
$
\\[1.5ex]
а значит (следствие 11.2), функция
\\[1.25ex]
\mbox{}\hfill
$
F_3^{}\colon (t,x,y)\to\ 
(x^2+y^2)\;\!\exp\Bigl(2\;\!(\eta-\lambda)\;\!t-2\;\!\arctg\dfrac{y}{x}\;\!\Bigr)
\quad
\forall (t,x,y)\in\Omega^{\;\!\prime}
\hfill
$
\\[1.5ex]
является неавтономным первых интегралом на множестве 
$\Omega^{\;\!\prime}$ 
дифференциальной системы (12.29) при $\lambda\ne 0.$
\vspace{0.35ex}

Поскольку первые интегралы $F,\ F_1^{},\ F_2^{},\ F_3^{}$
\vspace{0.35ex}
попарно функционально независимы на $\R^3,$
то каждая из совокупностей 
\vspace{0.35ex}
$\{F, F_1^{}\},\ \{F, F_2^{}\},\ \{F, F_3^{}\},\ 
\{F_1^{}, F_2^{}\},\ \{F_1^{}, F_3^{}\},\ \{F_2^{}, F_3^{}\}$
образует интегральный базис на соответствующем множестве
дифференциальной системы (12.29) при $\lambda\ne 0.$
\vspace{0.75ex}

{\sl Случай} $\lambda= 0.$
Тогда 
\\[2ex]
\mbox{}\hfill     % (12.30)
$
\dfrac{dx}{dt}=x\;\!(x^2+y^2)-y\;\!(x^2+y^2+\eta),
\quad
\dfrac{dy}{dt}=x\;\!(x^2+y^2+\eta)+y\;\!(x^2+y^2)
\quad
(\eta\ne 0).
$
\hfill (12.30)
\\[2.25ex]
\indent
Имеем:
\vspace{0.5ex}
$\bigl(x+i\;\!y,\;\! x^2+y^2+i\;\!(x^2+y^2+\eta)\bigr)\in \text{\rm H}_{_{\scriptstyle\R^2}},$
т.е. (теорема 5.3)
\\[1.75ex]
\mbox{}\hfill
$
\bigl(x^2+y^2,\;\! 2\;\!(x^2+y^2)\bigr)\in \text{\rm A}_{_{\scriptstyle\R^2}},
\quad 
\Bigl(\exp\arctg\dfrac{y}{x}\,,\;\! x^2+y^2+\eta\Bigr)\in \text{\rm E}_{_{D}}.
\hfill
$
\\[1.75ex]
\indent
Поскольку
\\[1.5ex]
\mbox{}\hfill
$
\bigl(x^2+y^2,\;\! 2\;\!(x^2+y^2)\bigr)\in \text{\rm A}_{_{\scriptstyle\R^2}},
\hfill
$
\\[1.5ex]
то, по свойству 5.10, 
\\[1.5ex]
\mbox{}\hfill
$
\bigl((x^2+y^2,\;\! 2\;\!(x^2+y^2)),\, (1, 1, {}-2)\bigr)\in \text{\rm B}_{_{\scriptstyle\R^2}},
\hfill
$
\\[1.5ex]
а значит (теорема 5.3),
\\[1.5ex]
\mbox{}\hfill
$
\Bigl(\exp\dfrac{1}{x^2+y^2}\,,{}-2\Bigr)\in 
\text{\rm E}_{_{\!\scriptstyle D_{_0}}}.
\hfill
$
\\[1.5ex]
\indent
Тождество
\\[1.5ex]
\mbox{}\hfill
$
2\;\!\gamma_1^{}\;\!(x^2+y^2)+
\gamma_2^{}\;\!(x^2+y^2+\eta)-2\;\!\gamma_3^{}=0
\quad 
\forall (x,y)\in\R^2
\hfill
$
\\[1.75ex]
выполняется, например, при $\gamma_1^{}={}-1,\ \gamma_2^{}=2,\ \gamma_3^{}=\eta.$
По свойству 11.11, функция
\\[2.25ex]
\mbox{}\hfill
$
F_4^{}\colon (x,y)\to\ 
\dfrac{1}{x^2+y^2}\,
\exp\Bigl(\;\!\dfrac{\eta}{x^2+y^2}+2\;\!\arctg\dfrac{y}{x}\;\!\Bigr)
\quad
\forall (x,y)\in D
\hfill
$
\\[2ex]
является автономным первых интегралом на множестве $D$ системы (12.30).
\vspace{0.35ex}

Линейная комбинация сомножителей
\\[1.5ex]
\mbox{}\hfill
$
2\;\!(x^2+y^2)-2\;\!(x^2+y^2+\eta)={}-2\;\!\eta
\quad 
\forall (x,y)\in\R^2.
\hfill
$
\\[1.75ex]
\indent
Согласно свойству 1.9
\\[1.5ex]
\mbox{}\hfill
$
\Bigl((x^2+y^2)\;\!\exp\Bigl({}-2\;\!\arctg\dfrac{y}{x}\;\!\Bigr),{}-2\;\!\eta\Bigr)\in 
\text{\rm J}_{_{\!\scriptstyle D}},
\hfill
$
\\[1.5ex]
а значит (следст\-вие~11.2), функция
\\[1.75ex]
\mbox{}\hfill
$
F_5^{}\colon (t,x,y)\to\ 
(x^2+y^2)\;\!\exp\Bigl(2\;\!\eta\;\!t-2\;\!\arctg\dfrac{y}{x}\;\!\Bigr)
\quad
\forall (t,x,y)\in\Omega^{\;\!\prime}
\hfill
$
\\[1.5ex]
является неавтономным первых интегралом на множестве 
$\Omega^{\;\!\prime}$ системы (12.30).
\vspace{0.5ex}

Поскольку
\vspace{0.5ex}
$
\Bigl(\exp\dfrac{1}{x^2+y^2}\,,{}-2\Bigr)\in 
\text{\rm E}_{_{\!\scriptstyle D_{_0}}},
$
то, по следствию 11.2, 
\\[1.5ex]
\mbox{}\hfill
$
\exp\Bigl(2\;\!t+\dfrac{1}{x^2+y^2}\;\!\Bigr)\in 
\text{\rm I}_{_{\scriptstyle \Omega_{_0}}},
\hfill
$
\\[1.5ex]
а значит (свойство 0.1), функция
\\[1.5ex]
\mbox{}\hfill
$
F_6^{}\colon (t,x,y)\to\ 
2\;\!t+\dfrac{1}{x^2+y^2}
\quad
\forall (t,x,y)\in\Omega_{_{0}}
\hfill
$
\\[1.5ex]
является неавтономным первых интегралом на области 
$\Omega_{_{0}}$ системы (12.30).
\vspace{0.35ex}

Поскольку первые интегралы $F_4^{},\ F_5^{},\ F_6^{}$
\vspace{0.35ex}
попарно функционально независимы на $\R^3,$
то каждая из совокупностей 
\vspace{0.35ex}
$\{F_4^{}, F_5^{}\},\ \{F_4^{}, F_6^{}\},\ \{F_5^{}, F_6^{}\}$
образует интегральный базис на множестве $\Omega^{\;\!\prime}$
дифференциальной системы (12.30).
\vspace{1ex}

{\bf 12.20.}
У системы Дарбу [9]
\\[2.25ex]
\mbox{}\hfill     % (12.31)
$
\dfrac{dx}{dt}=\lambda\;\!x-\eta\;\!y+x\;\!P(x,y),
\quad
\dfrac{dy}{dt}=\eta\;\!x+\lambda\;\!y+y\;\!P(x,y),
$
\hfill (12.31)
\\[2.5ex]
где $\lambda\in\R,\ \eta\in\R\backslash\{0\}, \ P$ --- полином степени $\deg P\geq 1,$
\vspace{0.75ex}
по виду правых частей на осно\-ва\-нии следствия 12.2, находим
$\bigl(x+i\;\!y,\;\! \lambda+P+\eta\;\!i\bigr)\in \text{\rm H}_{_{\scriptstyle\R^2}}.$
В соответствии с теоремой~5.3
\\[2.25ex]
\mbox{}\hfill
$
\bigl(x^2+y^2,\;\! 2\;\!(\lambda+P)\bigr)\in \text{\rm A}_{_{\scriptstyle\R^2}},
\quad 
\Bigl(\exp\arctg\dfrac{y}{x}\,,\;\! \eta\Bigr)\in \text{\rm E}_{_{D}},
\quad 
D=\{(x,y)\colon x\ne 0\}.
\hfill
$
\\[2.25ex]
\indent
По следствию 11.2, 
\\[1.5ex]
\mbox{}\hfill
$
\exp\Bigl({}-\eta\;\!t+\arctg\dfrac{y}{x}\;\!\Bigr)\in 
\text{\rm I}_{_{\scriptstyle \Omega_{_0}}},
\hfill
$
\\[1.5ex]
а значит (свойство 0.1), функция
\\[1.75ex]
\mbox{}\hfill
$
\Psi\colon (t,x,y)\to\ 
\eta\;\!t-\arctg\dfrac{y}{x}
\quad
\forall (t,x,y)\in\Omega_{_{0}}
\hfill
$
\\[1.5ex]
является неавтономным первых интегралом на множестве 
\vspace{0.75ex}
$\Omega_{_{0}}=\R\times D$ системы (12.31).

{\bf 12.20.А.}
\vspace{0.35ex}
Будем считать $\lambda\ne 0,\ P=x^2+y^2$ и 
рассмотрим систему (частный случай в [28, с. 43; 29])
\\[1.5ex]
\mbox{}\hfill     % (12.32)
$
\dfrac{dx}{dt}=\lambda\;\!x-\eta\;\!y+x\;\!(x^2+y^2),
\qquad
\dfrac{dy}{dt}=\eta\;\!x+\lambda\;\!y+y\;\!(x^2+y^2).
$
\hfill (12.32)
\\[2.25ex]
\indent
Имеем:
\vspace{1ex}
$\bigl(x^2+y^2,\;\! 2\;\!(x^2+y^2+\lambda)\bigr)\in \text{\rm A}_{_{\scriptstyle\R^2}},\ 
\Bigl(\exp\arctg\dfrac{y}{x}\,,\;\! \eta\Bigr)\in \text{\rm E}_{_{D}}.
$

Поскольку
$
\bigl(x^2+y^2,\;\! 2\;\!(x^2+y^2+\lambda)\bigr)\in \text{\rm A}_{_{\scriptstyle\R^2}},
$
то, по свойству 1.4, 
\\[2ex]
\mbox{}\hfill
$
\bigl(x^2+y^2+\lambda,\;\! 2\;\!(x^2+y^2)\bigr)\in \text{\rm A}_{_{\scriptstyle\R^2}}.
\hfill
$
\\[1.5ex]
\indent
Тождество
\\[1.5ex]
\mbox{}\hfill
$
2\;\!\gamma_1^{}\;\!(x^2+y^2+\lambda)+\gamma_2^{}\;\!\eta+
2\;\!\gamma_3^{}\;\!(x^2+y^2)=0
\quad 
\forall (x,y)\in\R^2
\hfill
$
\\[2ex]
выполняется, например, при $\gamma_1^{}={}-1,\ \gamma_2^{}=\dfrac{2\;\!\lambda}{\eta}\,,\ 
\gamma_3^{}=1.$
По свойству 11.11, функция
\\[2ex]
\mbox{}\hfill
$
F\colon (x,y)\to\ 
\dfrac{x^2+y^2+\lambda}{x^2+y^2}\,
\exp\Bigl(\;\!\dfrac{2\;\!\lambda}{\eta}\;\!\arctg\dfrac{y}{x}\;\!\Bigr)
\quad
\forall (x,y)\in D
\hfill
$
\\[1.5ex]
является автономным первых интегралом на множестве $D$ системы (12.32).
\vspace{0.35ex}

Первые интегралы $\Psi$ и $F,$
будучи функционально независимыми, образуют интегральный базис 
на множестве $\Omega_{_{0}}$ системы (12.32).
\vspace{1ex}

{\bf 12.20.Б.}
Пусть $\lambda= 0,\ P=x^2+y^2.$ Тогда
\\[1.5ex]
\mbox{}\hfill     % (12.33)
$
\dfrac{dx}{dt}={}-\eta\;\!y+x\;\!(x^2+y^2),
\qquad
\dfrac{dy}{dt}=\eta\;\!x+y\;\!(x^2+y^2).
$
\hfill (12.33)
\\[2ex]
\indent
Имеем:
\vspace{1ex}
$\bigl(x^2+y^2,\;\! 2\;\!(x^2+y^2)\bigr)\in \text{\rm A}_{_{\scriptstyle\R^2}},\ 
\Bigl(\exp\arctg\dfrac{y}{x}\,,\;\! \eta\Bigr)\in \text{\rm E}_{_{D}}.
$

Поскольку
$
\bigl(x^2+y^2,\;\! 2\;\!(x^2+y^2)\bigr)\in \text{\rm A}_{_{\scriptstyle\R^2}},
$
то, по свойству 5.10, 
\\[2ex]
\mbox{}\hfill
$
\bigl((x^2+y^2,\;\! 2\;\!(x^2+y^2)),\;\! (1, 1, {}-2)\bigr)\in \text{\rm B}_{_{\scriptstyle\R^2}},
\hfill
$
\\[1.5ex]
а значит (теорема 5.3), 
$
\Bigl(\exp\dfrac{1}{x^2+y^2}\,,\;\! {}-2\Bigr)\in \text{\rm E}_{_{D}}.
$
\vspace{1ex}

Так как 
$
\gamma_1^{}\;\!\eta -2\;\!\gamma_2^{}=0,
$
например, при $\gamma_1^{}=2,\ \gamma_2^{}=\eta,$ 
то, по свойству 11.11, 
\\[2ex]
\mbox{}\hfill
$
\exp\Bigl(\;\!\dfrac{\eta}{x^2+y^2}+2\;\!\arctg\dfrac{y}{x}\;\!\Bigr)\in \text{\rm I}_{_{D}},
\hfill
$
\\[1.5ex]
а значит (свойство 0.1), функция
\\[1.75ex]
\mbox{}\hfill
$
\Psi_1^{}\colon (x,y)\to\ 
\dfrac{\eta}{x^2+y^2}+2\;\!\arctg\dfrac{y}{x}
\quad
\forall (x,y)\in D
\hfill
$
\\[1.5ex]
является автономным первых интегралом на множестве $D$ системы (12.33).
\vspace{0.35ex}

Первые интегралы $\Psi$ и $\Psi_1^{},$
\vspace{0.25ex}
будучи функционально независимыми, образуют интегральный базис 
на множестве $\Omega_{_{0}}$ системы (12.33).
\vspace{1ex}

{\bf 12.20.В.} 
Пусть $\lambda= 0,\ P=y^2(x^2+y^2-1).$ Тогда ([28, с. 43] при $\eta=1)$
\\[1.75ex]
\mbox{}\hfill     % (12.34)
$
\dfrac{dx}{dt}={}-\eta\;\!y+xy^2\;\!(x^2+y^2-1),
\qquad
\dfrac{dy}{dt}=\eta\;\!x+y^3\;\!(x^2+y^2-1).
$
\hfill (12.34)
\\[2.25ex]
\indent
Имеем:
$\bigl(x+i\;\!y,\;\! y^2(x^2+y^2-1)+\eta\;\!i\bigr)\in \text{\rm H}_{_{\scriptstyle\R^2}},$
т.е. (теорема 5.3)
\\[1.75ex]
\mbox{}\hfill
$
\bigl(x^2+y^2,\;\! 2\;\!y^2\;\!(x^2+y^2-1)\bigr)\in \text{\rm A}_{_{\scriptstyle\R^2}},
\quad 
\Bigl(\exp\arctg\dfrac{y}{x}\,,\;\! \eta\Bigr)\in \text{\rm E}_{_{D}}.
\hfill
$
\\[1.75ex]
\indent
Поскольку
$
\bigl(x^2+y^2,\;\! 2\;\!y^2\;\!(x^2+y^2-1)\bigr)\in \text{\rm A}_{_{\scriptstyle\R^2}},
$
то, по свойству 1.4, 
\\[2ex]
\mbox{}\hfill
$
\bigl(x^2+y^2-1,\;\! 2\;\!y^2\;\!(x^2+y^2)\bigr)\in \text{\rm A}_{_{\scriptstyle\R^2}}.
\hfill
$
\\[1.5ex]
\indent
Расходимость 
\\[1.15ex]
\mbox{}\hfill
$
{\rm div}\;\!{\frak d}(t,x,y)=
6\;\!y^2\;\!(x^2+y^2)-4\;\!y^2
\quad
\forall (t,x,y)\in\R^3.
\hfill
$
\\[1ex]
\indent
Тождество
\\[1.5ex]
\mbox{}\hfill
$
2\;\!\gamma_1^{}\;\!y^2\;\!(x^2+y^2-1)+
2\;\!\gamma_2^{}\;\!y^2\;\!(x^2+y^2)=
{}-6\;\!y^2\;\!(x^2+y^2)+4\;\!y^2
\quad 
\forall (x,y)\in\R^2
\hfill
$
\\[1.75ex]
выполняется, например, при $\gamma_1^{}={}-2,\ \gamma_2^{}={}-1.$
По свойству 8.8, функция
\\[2ex]
\mbox{}\hfill
$
\mu\colon (x,y)\to\ 
\dfrac{1}{(x^2+y^2)^2(x^2+y^2-1)}
\quad
\forall (x,y)\in D_{_{0}}
\hfill
$
\\[1.75ex]
является интегрирующим множителем на множестве 
\vspace{0.25ex}
$\!D_{_{0}}\!=\!\{(x,y)\colon (x^2+y^2)(x^2+y^2-1)\!\ne\! 0\}\!$ 
уравнения траекторий системы (12.34).
\vspace{1ex}

{\bf 12.21.}
У системы [27]
\\[1.5ex]
\mbox{}\hfill     % (12.35)
$
\dfrac{dx}{dt}=x-x^3-x^2y-y^3,
\qquad
\dfrac{dy}{dt}=y+x^3-x^2y+xy^2
$
\hfill (12.35)
\\[2ex]
по виду правых частей на осно\-ва\-нии следствия 12.2, находим
\\[1.5ex]
\mbox{}\hfill
$
\bigl(x+i\;\!y,\;\! 1-x^2+i\;\!(x^2+y^2)\bigr)\in \text{\rm H}_{_{\scriptstyle\R^2}},
\hfill
$
\\[1.5ex]
а значит (теорема~5.3),
\\[1.75ex]
\mbox{}\hfill
$
\bigl(x^2+y^2,\;\! 2\;\!(1-x^2)\bigr)\in \text{\rm A}_{_{\scriptstyle\R^2}},
\quad 
\Bigl(\exp\arctg\dfrac{y}{x}\,,\;\! x^2+y^2\Bigr)\in \text{\rm E}_{_{\!D}},
\quad 
D=\{(x,y)\colon x\ne 0\}.
\hfill
$
\\[1.75ex]
\indent
По теореме 3.1, 
\\[1.5ex]
\mbox{}\hfill
$
\Bigl(\exp\dfrac{xy}{x^2+y^2}\,,\, x^2-y^2\Bigr)\in 
\text{\rm E}_{_{\!\scriptstyle D_{_0}}}
$
\ и \
$
\Bigl(\exp\dfrac{{}-y^2}{x^2+y^2}\,,\, {}-2\;\!xy\Bigr)\in 
\text{\rm E}_{_{\!\scriptstyle D_{_0}}},
\ D_{_0}=\R^2\backslash\{(0,0)\},
\hfill
$
\\[1.5ex]
а значит (теорема 5.3),
\\[1.5ex]
\mbox{}\hfill
$
\bigl((x^2+y^2, 2(1-x^2)), (1, xy, x^2-y^2)\bigr)\!\in\! \text{\rm B}_{_{\scriptstyle\R^2}}
$
и
$
\bigl((x^2+y^2, 2(1-x^2)), (1, {}-y^2,{}-2xy)\bigr)\!\in\! \text{\rm B}_{_{\scriptstyle\R^2}}.
\hfill
$
\\[1.5ex]
\indent
Таким образом,
\vspace{0.75ex}
$\!\bigl((x+iy, 1-x^2+i(x^2+y^2)), (1, y, x^2-y^2-2ixy)\bigr)\!\in\! 
\text{\rm G}_{_{\scriptstyle\R^2}}\!\!$
(следствие 7.1).

Поскольку сумма сомножителей
\\[1.5ex]
\mbox{}\hfill     
$
2\;\!(1-x^2)+(x^2+y^2)+(x^2-y^2)=2,
\hfill 
$
\\[1.5ex]
то согласно свойству 1.9
\\[1.5ex]
\mbox{}\hfill
$
\Bigl((x^2+y^2)\exp\Bigl(\;\!\dfrac{xy}{x^2+y^2}+\arctg\dfrac{y}{x}\;\!\Bigr),\, 2\Bigr)\in 
\text{\rm J}_{_{\!\scriptstyle D}}.
\hfill
$
\\[1.5ex]
\indent
По следствию 11.2, функция
\\[1.5ex]
\mbox{}\hfill
$
F\colon (t,x,y)\to\ 
(x^2+y^2)\exp\Bigl({}-2\;\!t+\dfrac{xy}{x^2+y^2}+\arctg\dfrac{y}{x}\;\!\Bigr)
\hfill
$
\\[1.5ex]
является 
\vspace{0.25ex}
неавтономным первым интегралом на множестве $\R\times D$ системы (12.35).

Расходимость 
\\[1.15ex]
\mbox{}\hfill
$
{\rm div}\;\!{\frak d}(t,x,y)=
2\;\!(1-2\;\!x^2)
\quad
\forall (t,x,y)\in\R^3.
\hfill
$
\\[0.75ex]
\indent
Тождество
\\[1.5ex]
\mbox{}\hfill
$
2\;\!\gamma_1^{}\;\!(1-x^2)+
\gamma_2^{}\;\!(x^2+y^2)+
\gamma_3^{}\;\!(x^2-y^2)=
{}-2\;\!(1-2\;\!x^2)
\quad 
\forall (x,y)\in\R^2
\hfill
$
\\[1.75ex]
выполняется, например, при $\gamma_1^{}={}-1,\ \gamma_2^{}=\gamma_3^{}=1.$
По свойству 8.8, функция
\\[2ex]
\mbox{}\hfill
$
\mu\colon (x,y)\to\ 
\dfrac{1}{x^2+y^2}\,
\exp\Bigl(\;\!\dfrac{xy}{x^2+y^2}+\arctg\dfrac{y}{x}\;\!\Bigr)
\quad
\forall (x,y)\in D
\hfill
$
\\[1.5ex]
есть интегрирующий множитель на множестве $\!D\!$ 
уравнения траекторий системы (12.35).

\newpage

{\bf 12.22.}
{\it Уравнение Якоби} имеет вид
\\[2ex]
\mbox{}\hfill        % (12.36)
$
\dfrac{dy}{dx}=
\dfrac{l_2^{}(x,y)-y\;\!l_3^{}(x,y)}{l_1^{}(x,y)-x\;\!l_3^{}(x,y)}\equiv 
\dfrac{Y(x,y)}{X(x,y)}\,,
\hfill
$
\\
\mbox{}\hfill (12.36)
\\
\mbox{}\hfill
$
l_i^{}(x,y)=a_{i}^{}\;\!x+b_{i}^{}\;\!y+c_i^{}\;\!,
\quad
a_{i}^{},b_{i}^{},c_i^{}\in\R, 
\quad
i=1,2,3.
\hfill
$
\\[2ex]
\indent
Линейный дифференциальный оператор первого порядка
\\[2ex]
\mbox{}\hfill
$
{\frak I}(x,y)=X(x,y)\;\!\partial_{{}_{\scriptstyle x}}+Y(x,y)\;\!\partial_{{}_{\scriptstyle y}}
\quad
\forall (x,y)\in\R^2
\hfill
$
\\[2ex]
суть оператор дифференцирования в силу уравнения Якоби (12.36).

Линейная функция 
\\[1.25ex]
\mbox{}\hfill
$
p\colon (x,y)\to\ 
\alpha\;\! x+\beta y+\gamma
\quad
\forall (x,y)\in\R^2
\quad
(|\alpha|+|\beta|\ne 0)
\hfill
$
\\[1.5ex]
является полиномиальным частным интегралом уравнения (12.36)
тогда и только тогда, когда выполняется тождество
\\[2ex]
\mbox{}\hfill
$
{\frak I}\;\!(\alpha\;\! x+\beta y+\gamma)=
(\alpha\;\! x+\beta y+\gamma)(\mu\;\! x+\nu\;\! y+\theta)
\quad
\forall (x,y)\in\R^2.
\hfill
$
\\[2ex]
\indent
Это тождество имеет место, если и только если совместна система уравнений
\\[2ex]
\mbox{}\hfill        % (12.37)
$
(a_1^{}-\lambda)\;\!\alpha+a_2^{}\beta +a_3^{}\gamma=0,
\ \ \
b_1^{}\alpha+(b_2^{}-\lambda)\;\!\beta +b_3^{}\gamma=0,
\ \ \
c_1^{}\alpha+c_2^{}\beta +(c_3^{}-\lambda)\;\!\gamma=0,
$
\hfill (12.37)
\\[2ex]
причем $\mu={}-a_3^{},\ \nu={}-b_3^{},\ \theta={}-c_3^{}+\lambda.$
\vspace{0.5ex}

Для того чтобы линейная однородная
система (12.37) имела нетривиальное решение $(\alpha,\beta,\gamma)$
необходимо и достаточно, чтобы был равен нулю ее определитель:
\\[1.75ex]
\mbox{}\hfill    % (12.38)
$
\det(A-\lambda E)=0,
$
\hfill (12.38)
\\[1.5ex]
где $E$ --- единичная матрица, а матрица
\\[1.5ex]
\mbox{}\hfill 
$
A=\left(\!\!
\begin{array}{ccc}
a_1^{} & a_2^{} & a_3^{}
\\[0.75ex]
b_1^{} & b_2^{} & b_3^{}
\\[0.5ex]
c_1^{} & c_2^{} & c_3^{}
\end{array}
\!\!\right).
\hfill
$
\\[1.5ex]
\indent
Корни уравнения (12.38) суть собственные числа матрицы $A.$
При этом решение $(\alpha,\beta,\gamma)$ системы (12.37) есть 
собственный вектор матрицы $A,$ соответствующий ее 
собственному числу $\lambda.$

Таким образом, доказана
\vspace{0.35ex}

{\bf Лемма 12.1.}
\vspace{0.25ex}
{\it
Если $(\alpha,\beta,\gamma)$ --- собственный вектор, 
соответствующий собственному числу $\lambda$ матрицы $A,$
то линейная функция 
\vspace{0.25ex}
$
p\colon (x,y)\to 
\alpha\;\! x+\beta y+\gamma\;\;
\forall (x,y)\in\R^2
$
$
(|\alpha|+|\beta|\ne 0)
$
\vspace{0.35ex}
является полиномиальным частным интегралом с сомножителем 
$M\colon (x,y)\to \lambda-l_3^{}(x,y)\;\;\forall (x,y)\in\R^2$
уравнения Якоби} (12.36).
\vspace{0.5ex}

Обратим внимание на то, что в лемме 12.1 не исключается возможность,
когда собственное число $\lambda$ является комплексным.

Построим общий интеграл уравнения Якоби (12.36) по собственным векторам и 
собственным числам матрицы $A$ с учетом кратности элементарных делителей.
\vspace{0.5ex}

{\small\bf 
12.22.1. Случай трех простых элементарных делителей.}
\vspace{0.35ex}

{\bf Теорема 12.8.}
\vspace{0.25ex}
{\it
Пусть $\lambda-\lambda_1^{},\ \lambda-\lambda_2^{},\ \lambda-\lambda_3^{}$ --- 
простые элементарные делители, а 
$
(\alpha_1^{},\beta_1^{},\gamma_1^{}),\ 
(\alpha_2^{},\beta_2^{},\gamma_2^{}),\ 
(\alpha_3^{},\beta_3^{},\gamma_3^{}) 
$
\vspace{0.35ex}
--- соответствующие собственные векторы матрицы $A.$
Тогда общим интегралом на $D$ уравнения Якоби {\rm(12.36)}
будет функция
\\[2ex]
\mbox{}\hfill  % (12.39)
$
F\colon (x,y)\to\ 
p_1^{{}^{\scriptstyle h_1^{}}}\!(x,y)\;\!
p_2^{{}^{\scriptstyle h_2^{}}}(x,y)\;\!
p_3^{{}^{\scriptstyle h_3^{}}}(x,y)
\quad
\forall (x,y)\in D,
$
\hfill {\rm (12.39)}
\\[2.25ex]
где $p_i^{}(x,y)=\alpha_i^{}\;\!x+\beta_{i}^{}\;\!y+\gamma_i^{}\;\;\forall (x,y)\in\R^2,\ i=1,2,3,$
показатели степеней $h_1^{},\, h_2^{},\, h_3^{}$ связаны равенствами
\\[2ex]
\mbox{}\hfill  % (12.40)
$
h_1^{}+h_2^{}+h_3^{}=0,
\quad
\lambda_1^{}h_1^{}+\lambda_2^{}h_2^{}+\lambda_3^{}h_3^{}=0,
$
\hfill {\rm (12.40)}
\\[2ex]
а множество $D$ такое, что функции $p_i^{{}^{\scriptstyle h_i^{}}}\in C^1D,\ i=1,2,3.$}
\vspace{0.35ex}

{\sl Доказательство.}
По лемме 12.1, линейные функции $p_i^{}$ являются полиномиальными 
(вещественными или комплекснозначными) частными интегралами с сомножителями 
$\lambda_i^{}-l_3^{},\ i=1,2,3,$ уравнения Якоби (12.36). 
\vspace{0.25ex}
Пусть числа $h_1^{},\, h_2^{},\, h_3^{}$ связаны равенствами (12.40).
Тогда у частных интегралов сомножители такие, что линейная комбинация
\\[2ex]
\mbox{}\hfill
$
h_1^{}\bigl(\lambda_1^{}-l_3^{}(x,y)\bigr)+
h_2^{}\bigl(\lambda_2^{}-l_3^{}(x,y)\bigr)+
h_3^{}\bigl(\lambda_3^{}-l_3^{}(x,y)\bigr)=0
\quad
\forall (x,y)\in\R^2.
\hfill
$
\\[2ex]
\indent
Согласно свойствам 11.11 и 6.6 функция (12.39) будет общим интегралом на 
множестве $D$ уравнения Якоби (12.36). $\k$
\vspace{0.35ex}

Обратим внимание на то, что в теореме 12.8 среди собственных чисел 
$\lambda_1^{},\, \lambda_2^{},\, \lambda_3^{}$
могут быть как равные, так и комплексные.
\vspace{0.5ex}

{\small\bf 12.22.1.1. Случай трех различных вещественных собственных чисел.}
Если у матрицы $A$ собственные числа $\lambda_1^{},\, \lambda_2^{},\, \lambda_3^{}$
\vspace{0.25ex}
вещественные и различные, то им соответствуют простые элементарные делители
$\lambda-\lambda_1^{},\ \lambda-\lambda_2^{},\ \lambda-\lambda_3^{}.$
\vspace{0.25ex}
Тогда из теоремы 12.8 при 
$h_1^{}=\lambda_2^{}-\lambda_3^{},\ 
h_2^{}=\lambda_3^{}-\lambda_1^{},\ 
h_3^{}=\lambda_1^{}-\lambda_2^{}$ следует
\vspace{0.5ex}

{\bf Теорема 12.9.}
\vspace{0.25ex}
{\it
Пусть у матрицы $A$ собственные числа
$\lambda_1^{},\, \lambda_2^{},\, \lambda_3^{}$ вещественные и различные, а 
$
(\alpha_1^{},\beta_1^{},\gamma_1^{}),\ 
(\alpha_2^{},\beta_2^{},\gamma_2^{}),\ 
(\alpha_3^{},\beta_3^{},\gamma_3^{}) 
$
\vspace{0.35ex}
--- соответствующие собственные векторы.
Тогда общим интегралом на $D$ уравнения Якоби {\rm(12.36)}
будет функция
\\[2ex]
\mbox{}\hfill  % (12.41)
$
F\colon (x,y)\to\ 
p_1^{{}^{\scriptstyle \lambda_2^{}-\lambda_3^{}}}\!(x,y)\;\!
p_2^{{}^{\scriptstyle \lambda_3^{}-\lambda_1^{}}}\!(x,y)\;\!
p_3^{{}^{\scriptstyle \lambda_1^{}-\lambda_2^{}}}\!(x,y)
\quad
\forall (x,y)\in D,
$
\hfill {\rm (12.41)}
\\[2.25ex]
где $p_i^{}(x,y)=\alpha_i^{}\;\!x+\beta_{i}^{}\;\!y+\gamma_i^{}\;\;\forall (x,y)\in\R^2,\ i=1,2,3,$
\vspace{0.5ex}
а множество $D\subset\R^2$ такое, что функции 
$p_1^{{}^{\scriptstyle \lambda_2^{}-\lambda_3^{}}},\,
p_2^{{}^{\scriptstyle \lambda_3^{}-\lambda_1^{}}},\,
p_3^{{}^{\scriptstyle \lambda_1^{}-\lambda_2^{}}}$
являются непрерывно дифференцируемыми на $D.$}
\vspace{1ex}

{\bf Пример 12.1.}
У уравнения Якоби
\\[2ex]
\mbox{}\hfill  % (12.42)
$
\dfrac{dy}{dx}=\dfrac{{}-x+5y-1-y\;\!(x-y+3)}{3x-y+1-x\;\!(x-y+3)}
$
\hfill (12.42)
\\[2ex]
матрица
$
A=\left(\!\!\!
\begin{array}{rrr}
3 & {}-1 & 1
\\
{}-1 & 5 & {}-1
\\
1 & {}-1 & 3
\end{array}
\!\!\right)
$
имеет собственные числа
$\lambda_1^{}=2,\ \lambda_2^{}=3,\ \lambda_3^{}=6,$
которым соответствуют собственные векторы 
$(1, 0,{}-1),\ (1, 1, 1),\ (1,{}-2, 1).$ 

По теореме 12.9, функция
\\[2ex]
\mbox{}\hfill  % (12.43)
$
F\colon (x,y)\to\ 
\dfrac{(x+y+1)^4}{(x-1)^3(x-2y+1)}
\quad
\forall (x,y)\in D
$
\hfill {\rm (12.43)}
\\[2.25ex]
является общим интегралом на множестве 
$D=\{(x,y)\colon (x-1)(x-2y+1)\ne 0\}$
уравнения Якоби (12.42).
\vspace{0.5ex}

{\small\bf 12.22.1.2. Случай, когда кратному собственному числу соответствуют
простые элементарные делители.}
Если $\lambda_1^{}$ --- кратное собственное число
матрицы $A$ и ему соответствуют два простых элементарных делителя
$\lambda-\lambda_1^{}$ и $\lambda-\lambda_1^{},$
то кратному собственному числу $\lambda_1^{}$  соответствуют два 
линейно независимых собственных вектора
$(\alpha_1^{},\beta_1^{},\gamma_1^{})$ и
$(\alpha_2^{},\beta_2^{},\gamma_2^{}).$ 
\vspace{0.5ex}

{\bf Теорема 12.10.}
\vspace{0.25ex}
{\it
Если кратному собственному числу матрицы $A$ соответствуют 
простые элементарные делители и линейно независимые собственные векторы
$(\alpha_1^{},\beta_1^{},\gamma_1^{})$ и
$(\alpha_2^{},\beta_2^{},\gamma_2^{}),$ 
то дробно-линейная функция
\\[2ex]
\mbox{}\hfill  % (12.44)
$
F\colon (x,y)\to\ 
\dfrac{\alpha_1^{}\;\!x+\beta_{1}^{}\;\!y+\gamma_1^{}}{\alpha_2^{}\;\!x+\beta_{2}^{}\;\!y+\gamma_2^{}}
\quad
\forall (x,y)\in D
$
\hfill {\rm (12.44)}
\\[2.25ex]
будет общим интегралом на множестве 
$D=\{(x,y)\colon \alpha_2^{}\;\!x+\beta_{2}^{}\;\!y+\gamma_2^{}\ne 0\}$ 
уравнения Якоби {\rm(12.36)}.
}
\vspace{0.35ex}

{\sl Доказательство.}
Пусть $\lambda$ --- кратное собственное число матрицы $A,$ которому соответствуют 
линейно независимые собственные векторы
$(\alpha_1^{},\beta_1^{},\gamma_1^{})$ и
$(\alpha_2^{},\beta_2^{},\gamma_2^{}).$
По лемме 12.1, линейные функции 
\\[1.5ex]
\mbox{}\hfill
$
p_1^{}\colon (x,y)\to\, \alpha_1^{}\;\!x+\beta_{1}^{}\;\!y+\gamma_1^{}
\;\;\forall (x,y)\in\R^2
$ 
\ и\ 
$
p_2^{}\colon (x,y)\to\, \alpha_2^{}\;\!x+\beta_{2}^{}\;\!y+\gamma_2^{}
\;\;\forall (x,y)\in\R^2
\hfill
$
\\[1.5ex]
являются полиномиальными частными интегралами уравнения Якоби (12.36)
\vspace{0.25ex}
и у них один и тот же сомножитель 
$M\colon (x,y)\to \lambda-l_3^{}(x,y)\;\;\forall (x,y)\in\R^2.$
\vspace{0.25ex}
Тогда согласно следствию 11.7 дробно-линейная функция (12.44)
будет общим интегралом на множестве $D$ уравнения Якоби (12.36). $\k$
\vspace{0.35ex}

Заметим, что утверждение теоремы 12.10 также следует из теоремы 12.8
при $\lambda_2^{}=\lambda_1^{},$ $h_1^{}=1,\ h_2^{}={}-1,\ h_3^{}=0.$
\vspace{0.35ex}

Если у матрицы $A$ трехкратному собственному числу $\lambda_1^{}$ 
\vspace{0.25ex}
соответствует три простых элементарных делителя
$\lambda-\lambda_1^{},\ \lambda-\lambda_1^{},\ \lambda-\lambda_1^{},$
то матрица $A=\lambda_1^{}E.$
\vspace{0.25ex}
В этом случае уравнение Якоби (12.36) вырождается, так как
$X(x,y)=Y(x,y)=0\;\;\forall (x,y)\in\R^2.$
\vspace{0.5ex}

{\bf Пример 12.2.}
У уравнения Якоби
\\[2ex]
\mbox{}\hfill  % (12.45)
$
\dfrac{dy}{dx}=\dfrac{x-y+1-y\;\!(x+y-1)}{{}-x+y+1-x\;\!(x+y-1)}
$
\hfill (12.45)
\\[2ex]
матрица
$
A=\left(\!\!\!
\begin{array}{rrr}
{}-1 & 1 & 1
\\
1 & {}-1 & 1
\\
1 & 1 & {}-1
\end{array}
\!\!\right)
$
имеет собственные числа
$\lambda_1^{}=\lambda_2^{}={}-2,\ \lambda_3^{}=1,$
которым соответствуют простые элементарные делители
$\lambda+2,\ \lambda+2,\ \lambda-1.$
\vspace{0.25ex}
Двукратному собственному числу $\lambda_1^{}={}-2$
соответствуют собственные векторы 
\vspace{0.25ex}
$(0, 1,{}-1)$ и $(1, 0, {}-1).$ 
По теореме 12.10, дробно-линейная функция
\\[2ex]
\mbox{}\hfill  % (12.46)
$
F\colon (x,y)\to\ 
\dfrac{y-1}{x-1}
\quad
\forall (x,y)\in D
$
\hfill {\rm (12.46)}
\\[2.25ex]
является общим интегралом на множестве 
$D=\{(x,y)\colon x\ne 1\}$
\vspace{0.35ex}
уравнения Якоби (12.45).

Общий интеграл (12.46) определяет пучок прямых-траекторий
\\[1.5ex]
\mbox{}\hfill
$
\dfrac{y-1}{x-1}=C,
\quad 
C\in [{}-\infty;{}+\infty],
\hfill
$
\\[1.5ex]
уравнения Якоби (12.45).
\vspace{0.25ex}

Собственному числу $\lambda_3^{}=1$
\vspace{0.35ex}
соответствует собственный вектор $(1, 1, 1).$ 
По лемме~12.1, линейная функция
$p_3^{}\colon (x,y)\to x+y+1\;\; \forall (x,y)\in\R^2$
\vspace{0.5ex}
является полиномиальным част\-ным интегралом с сомножителем 
$M_3^{}\colon (x,y)\to 2-x-y\;\;\forall (x,y)\in\R^2$
\vspace{0.35ex}
уравнения Яко\-би (12.45). 
\vspace{0.35ex}
Поэтому прямая $x+y+1=0$ суть траектория уравнения Якоби (12.45).

Однако, прямая-траектория $x+y+1=0$ 
\vspace{0.5ex}
не содержится в пучке прямых-траекторий
$\dfrac{y-1}{x-1}=C,\ \ C\in [{}-\infty;{}+\infty].$
\vspace{1.5ex}
Это обосновано тем, что уравнение Якоби (12.45) приводится к виду
\\[1ex]
\mbox{}\hfill
$
\dfrac{dy}{dx}=\dfrac{(1-y)(x+y+1)}{(1-x)(x+y+1)}\,,
\hfill
$
\\[1.5ex]
а значит, каждая точка прямой-траектории $x+y+1=0$ является особой.
\vspace{0.25ex}

Таким образом, для уравнения Якоби (12.45) наряду с общим интегралом (12.46)
следует указывать и его частный интеграл $p_3^{}\colon (x,y)\to x+y+1\;\; \forall (x,y)\in\R^2.$

{\small\bf 12.22.1.3. Случай комплексного собственного числа.}
\vspace{0.25ex}
Если существенно комплексное число 
$\lambda_1^{}=\xi_1^{}+\zeta_1^{}i\ (\zeta_1^{}\ne 0)$ 
\vspace{0.5ex}
является корнем уравнения (12.38), то и комплексно сопряженное число
\vspace{0.35ex}
$\overline{\lambda}_1^{}=\xi_1^{}-\zeta_1^{}i$ будет его корнем. 
Кроме этого, уравнение (12.38) имеет вещественный корень $\lambda_3^{}.$
\vspace{0.35ex}

Если при $\lambda=\lambda_1^{}$ линейная система (12.37)
\vspace{0.35ex}
имеет решение 
$\alpha=\alpha_1^{},\ \beta=\beta_1^{},\ \gamma=\gamma_1^{},$
то при $\lambda=\overline{\lambda}_1^{}$ решением линейной системы (12.37)
будет 
\vspace{0.35ex}
$\alpha=\overline{\alpha}_1^{},\ \beta=\overline{\beta}_1^{},\ 
\gamma=\overline{\gamma}_1^{}.$

Решения $\lambda_1^{}$ и $\overline{\lambda}_1^{}$ уравнения (12.38)
\vspace{0.35ex}
суть собственные числа матрицы $A,$ а решения 
$(\alpha_1^{},\beta_1^{},\gamma_1^{})$ и
$(\overline{\alpha}_1^{},\overline{\beta}_1^{},\overline{\gamma}_1^{})$ 
\vspace{0.35ex}
системы (12.37) есть собственные векторы матрицы $A,$
соответствующие собственным числам $\lambda_1^{}$ и $\overline{\lambda}_1^{}.$ 
\vspace{0.25ex}

Следовательно, если матрица $A$ имеет существенно комплексное собственное число
$\lambda_1^{}=\xi_1^{}+\zeta_1^{}i\ (\zeta_1^{}\ne 0),$
\vspace{0.5ex}
которому соответствует собственный вектор $(\alpha_1^{},\beta_1^{},\gamma_1^{}),$
то вектор $(\overline{\alpha}_1^{},\overline{\beta}_1^{},\overline{\gamma}_1^{})$ 
\vspace{0.5ex}
есть собственный вектор матрицы $A,$ соответствующий ее собственному числу
$\overline{\lambda}_1^{}=\xi_1^{}-\zeta_1^{}i.$
\vspace{0.35ex}
Кроме этого, матрица $A$ имеет вещественное собственное число $\lambda_3^{}.$
Собственным числам $\lambda_1^{},\ \overline{\lambda}_1^{},\ \lambda_3^{}$
\vspace{0.35ex}
матрицы $A$ соответствуют простые  элементарные делители
$\lambda-\lambda_1^{},\ \lambda-\overline{\lambda}_1^{},\ \lambda-\lambda_3^{}.$
\vspace{0.75ex}

{\bf Теорема 12.11.}
\vspace{0.25ex}
{\it 
Пусть собственное число $\lambda=\xi+\zeta\;\!i$
матрицы $A$ существенно комплексное $(\zeta\ne 0),$
\vspace{0.25ex}
которому соответствует собственный вектор $(\alpha,\beta,\gamma),$ а
$\lambda_3^{}$ --- вещественное собственное число матрицы $A,$
которому соответствует собственный вектор $(\alpha_{3}^{},\beta_{3}^{},\gamma_{3}^{}).$
\vspace{0.25ex}
Тогда общим интегралом на множестве $D$ уравнения Якоби {\rm (12.36)}
будет функция
\\[1.5ex]
\mbox{}\hfill    % (12.47)
$
F\colon (x,y)\to\ 
H^{{}^{\scriptstyle \zeta}}\!(x,y)\exp S(x,y)
\quad
\forall (x,y)\in D,
$
\hfill {\rm (12.47)}
\\[2.25ex]
где 
\vspace{1ex}
$
H(x,y)=\dfrac{{\rm Re}^2p(x,y)+{\rm Im}^2p(x,y)}{p_3^2(x,y)}
\;\; \forall (x,y)\in D,
\ \,
S(x,y)=2\;\!(\lambda_3^{}-\xi)\arctg\dfrac{{\rm Im}\;\!p(x,y)}{{\rm Re}\;\!p(x,y)}
$
$
\forall (x,y)\in D,\
\vspace{0.5ex}
p(x,y)=\alpha x+\beta y+\gamma
\;\;\forall (x,y)\in\R^2, \
p_3^{}(x,y)=\alpha_3^{} x+\beta_3^{} y+\gamma_3^{}
\;\;\forall (x,y)\in\R^2, \
$
множество $D\subset\R^2$ такое, что функции 
$H^{{}^{\scriptstyle \zeta}}$ и 
$\arctg\dfrac{{\rm Im}\;\!p}{{\rm Re}\;\!p}$
непрерывно дифференцируемые на $D.$
}
\vspace{0.5ex}

{\sl Доказательство}.
По лемме 12.1, для уравнения Якоби (12.36) линейная функция $p_3^{}$
\vspace{0.35ex}
является вещественным полиномиальным частным интегралом с сомножителем
$M_3^{}\colon (x,y)\to \lambda_3^{}-l_3^{}(x,y)\;\;\forall (x,y)\in\R^2,$
а комплекснозначная линейная функция $p$ является 
комплекснозначным полиномиальным частным интегралом с сомножителем 
$M\colon (x,y)\to \lambda-l_3^{}(x,y)\;\;\forall (x,y)\in\R^2.$
\vspace{0.5ex}

Пусть $p(x,y)=u(x,y)+i\;\!v(x,y)\;\;\forall (x,y)\in\R^2,$
\vspace{0.5ex}
где 
$
u(x,y)={\rm Re}\;\!p(x,y)\;\;\forall (x,y)\in\R^2,
$
$ 
v(x,y)={\rm Im}\;\!p(x,y)\;\;\forall (x,y)\in\R^2. 
$
\vspace{0.35ex}
По теореме 6.2, функция $p$ является комплекснозначным полиномиальным
частным интегралом с сомножителем 
\vspace{0.35ex}
$M\colon (x,y)\to\xi-l_3^{}(x,y)+\zeta i$ $\forall (x,y)\in\R^2$
тогда и только тогда, когда функция $u^2+v^2$ 
\vspace{0.35ex}
является полиномиальным частным
интегралом с сомножителем
\vspace{0.25ex}
$M_1^{}\colon (x,y)\to 2(\xi-l_3^{}(x,y))\;\;\forall (x,y)\in\R^2$
и функция
$\exp\arctg\dfrac{v}{u}$ 
\vspace{0.5ex}
является экспоненциальным частным интегралом 
с сомножителем
$M_2^{}\colon (x,y)\to \zeta\;\;\forall (x,y)\in\R^2$
на множестве $D_0^{}=\{(x,y)\colon u(x,y)\ne 0\}.$
\vspace{0.75ex}

Сомножители $M_1^{},\ M_2^{},\ M_3^{}$ частных интегралов
$u^2+v^2,\ \exp\arctg\dfrac{v}{u}\,,\ p_3^{}$
такие, что линейная комбинация
\\[1.5ex]
\mbox{}\hfill
$
\zeta\;\!M_1^{}(x,y)+2\;\!(\lambda_3^{}-\xi)\;\!M_2^{}(x,y)-2\zeta\;\!M_3^{}(x,y)=0
\quad
\forall (x,y)\in\R^2,
\hfill
$
\\[1.5ex]
а значит (свойство 11.11), общим интегралом на множестве $D$ уравнения Якоби (12.36)
будет функция
\\[2ex]
\mbox{}\hfill
$
F\colon (x,y)\to\ 
\bigl(u^2(x,y)+v^2(x,y)\bigr)^{\zeta}\;\!
\biggl(\exp\arctg\dfrac{v(x,y)}{u(x,y)}\biggr)^{\!2\;\!(\lambda_3^{}-\xi)}
p_3^{{}-2\zeta}(x,y)
\quad
\forall (x,y)\in D.
\ \k
\hfill
$
\\[2.5ex]
\indent
Заметим, что утверждение теоремы 12.11 
\vspace{0.35ex}
также следует из теоремы 12.8 при 
$
p_2^{}=\overline{p}_1^{},
$
$
\lambda_1^{}=\xi+\zeta\;\! i, \
\vspace{0.5ex}
\lambda_2^{}=\xi-\zeta\;\! i=\overline{\lambda}_1^{}, \
h_1^{}=(\overline{\lambda}_1^{}-\lambda_3^{})\;\! i=\zeta+(\xi-\lambda_3^{})\;\! i, \
h_2^{}=(\lambda_3^{}-\lambda_1^{})\;\! i=\zeta-(\xi-\lambda_3^{})\;\! i=
\linebreak
=\overline{h}_1^{}, \
h_3^{}=(\lambda_1^{}-\overline{\lambda}_1^{})\;\! i={}-2\;\!\zeta.
$
\vspace{0.75ex}

{\bf Пример 12.3.}
У уравнения Якоби [9, с. 24 -- 25]
\\[2ex]
\mbox{}\hfill  % (12.48)
$
\dfrac{dy}{dx}=\dfrac{6x-2y+1-y\;\!({}-2x+y)}{4x-3y-1-x\;\!({}-2x+y)}
$
\hfill (12.48)
\\[2ex]
матрица
\vspace{0.35ex}
$
A=\left(\!\!\!
\begin{array}{rrr}
4 & 6 & {}-2
\\
{}-3 & {}-2 & 1
\\
{}-1 & 1 & 0\end{array}
\!\!\right)
$
имеет существенно комплексное собственное число
$\lambda=1+\sqrt{6}\;\!i$ и вещественное собственное число $\lambda_3^{}=0.$
\vspace{0.35ex}

Комплексному собственному числу $\lambda=1+\sqrt{6}\;\!i$
\vspace{0.35ex}
соответствует собственный вектор 
$\bigl(-\;\!2\;\!(2+\sqrt{6}\;\!i),\;\! 5,\;\! 3-\sqrt{6}\;\!i\bigr),$
\vspace{0.35ex}
а вещественному собственному числу $\lambda_3^{}=0$
соответствует собственный вектор $(1, 1, 5).$
\vspace{0.35ex}

По теореме 12.11, 
\vspace{0.35ex}
общим интегралом уравнения Якоби (12.48) на множестве 
$D=
\linebreak
=\{(x,y)\colon (x+y+5)(4x-5y-3)\ne 0\}$ будет функция
\\[2ex]
\mbox{}\hfill  % (12.49)
$
F\colon (x,y)\to\ 
\biggl(\dfrac{(4x-5y-3)^2+6(2x+1)^2}{(x+y+5)^2}\biggr)^{\!\sqrt{6}}
\exp\biggl({}-2\;\!\arctg\dfrac{\sqrt{6}\;\!(2x+1)}{4x-5y-3}\biggr).
%\quad
%\forall (x,y)\in D.
$
\hfill {\rm (12.49)}
\\[2.5ex]
\indent
{\bf Пример 12.4.}
У уравнения Якоби
\\[2ex]
\mbox{}\hfill  % (12.50)
$
\dfrac{dy}{dx}=\dfrac{x+y-3-y\;\!(x+2y+1)}{x+2y+3-x\;\!(x+2y+1)}
$
\hfill (12.50)
\\[2ex]
матрица
\vspace{0.35ex}
$
A=\left(\!\!\!
\begin{array}{rrr}
1 & 1 & 1
\\
2 & 1 & 2
\\
3 & {}-3 & 1\end{array}
\!\!\right)
$
имеет собственные числа $\lambda=1+i$ и $\lambda_3^{}=1,$
которым соответствуют собственные векторы
$(3-i,\;\! 4, {}-3+3\;\!i)$ и $(1, 1, {}-1).$
\vspace{0.5ex}
По теореме 12.11 (при $\xi=\lambda_3^{}),$ 
общим интегралом  на множестве 
\vspace{0.5ex}
$D=\{(x,y)\colon x+y-1\ne 0\}$ уравнения Якоби (12.50) 
будет рациональная функция
\\[1.5ex]
\mbox{}\hfill  % (12.51)
$
F\colon (x,y)\to\ 
\dfrac{(3x+4y-3)^2+(x-3)^2}{(x+y-1)^2}
\quad
\forall (x,y)\in D.
$
\hfill {\rm (12.51)}
\\[2ex]
\indent
{\small\bf 
12.22.2. Случай кратного элементарного делителя.}
\vspace{0.35ex}

Если $\lambda_1^{}$ --- собственное число матрицы $A,$ которому соответствует 
кратный элементарный делитель, то ему кроме собственного вектора
$\theta_1^{}=(\alpha_1^{},\beta_1^{},\gamma_1^{})$
соответствует первый присоединенный вектор 
$\theta_1^{(1)}=\bigl(\alpha_1^{(1)},\beta_1^{(1)},\gamma_1^{(1)}\bigr).$
\vspace{0.35ex}

Первый присоединенный вектор $\theta_1^{(1)}$
является решением матричного уравнения
\\[1.5ex]
\mbox{}\hfill % (12.52)
$
(A-\lambda_1^{}E)\;\! {\rm colon}\,\theta_1^{(1)} =
{\rm colon}\;\!\theta_1^{}.
$
\hfill (12.52)
\\[2ex]
\indent
{\bf Лемма 12.2.}
{\it 
Пусть собственному числу $\lambda_1^{}$ матрицы $A$ 
соответствуют кратный элементарный делитель, собственный вектор 
$\theta_1^{}=(\alpha_1^{},\beta_1^{},\gamma_1^{})$ и 
присоединенный вектор 
\vspace{0.35ex}
$\theta_1^{(1)}=\bigl(\alpha_1^{(1)},\beta_1^{(1)},\gamma_1^{(1)}\bigr).$
Тогда уравнение Якоби {\rm (12.36)}
имеет кратный полиномиальный частный интеграл
$p_1^{}\colon (x,y)\to \alpha_1^{}x+\beta_1^{}y+\gamma_1^{}\;\; \forall (x,y)\in\R^2$
\vspace{0.35ex}
с сомножителем
$M_1^{}\colon (x,y)\to \lambda_1^{}-l_3^{}(x,y)\;\; \forall (x,y)\in\R^2$
\vspace{0.35ex}
на множестве 
$D=\{(x,y)\colon \alpha_1^{}x+\beta_1^{}y+\gamma_1^{}\ne 0\}$
такой, что
$\bigl((p_1^{}, M_1^{}), (1,p_1^{(1)}, 1)\bigr)\in \text{\rm B}_{_{\!D}},$
где 
$p_1^{(1)}(x,y)= \alpha_1^{(1)}x+\beta_1^{(1)}y+\gamma_1^{(1)}\;\; \forall (x,y)\in\R^2.$
}

{\sl Доказательство.}
Матричное уравнение (12.52) равносильно линейной системе
\\[1.5ex]
\mbox{}\hfill              % (12.53)               
$\!\!\!
\begin{array}{c}
(a_1^{}\!-\!\lambda_1^{})\;\!\alpha_{1}^{(1)}+
a_2^{}\;\!\beta_{1}^{(1)}+
a_3^{}\;\!\gamma_{1}^{(1)}=\alpha_1^{},
\\[1.5ex]
b_1^{}\;\!\alpha_{1}^{(1)}+
(b_2^{}\!-\!\lambda_1^{})\;\!\beta_{1}^{(1)}+
b_3^{}\;\!\gamma_{1}^{(1)}=\beta_1^{},
\\[1.5ex]
c_1^{}\;\!\alpha_{1}^{(1)}+
c_2^{}\;\!\beta_{1}^{(1)}+
(c_3^{}\!-\!\lambda_1^{})\;\!\gamma_{1}^{(1)}=\gamma_1^{}
\end{array}
%\Leftrightarrow
\!\!\!\iff\!\!\!
\begin{array}{c}
a_1^{}\;\!\alpha_{1}^{(1)}+
a_2^{}\;\!\beta_{1}^{(1)}=
\alpha_1^{}+\lambda_1^{}\;\!\alpha_{1}^{(1)}-a_3^{}\;\!\gamma_{1}^{(1)},
\\[1.5ex]
b_1^{}\;\!\alpha_{1}^{(1)}+
b_2^{}\;\!\beta_{1}^{(1)}= 
\beta_1^{}+\lambda_1^{}\;\!\beta_{1}^{(1)}-b_3^{}\;\!\gamma_{1}^{(1)},
\\[1.5ex]
c_1^{}\;\!\alpha_{1}^{(1)}+
c_2^{}\;\!\beta_{1}^{(1)}=
\gamma_1^{}+\lambda_1^{}\;\!\gamma_{1}^{(1)}-c_3^{}\;\!\gamma_{1}^{(1)}.
\end{array}\!\!\!
$
\hfill (12.53)
\\[2ex]
\indent
Производная в силу уравнения Якоби (12.36)
\\[2ex]
\mbox{}\hfill
$
{\frak I}\;\!p_1^{(1)}(x,y)=
{\frak I}\;\!\bigl(\alpha_1^{(1)}x+\beta_1^{(1)}y+\gamma_1^{(1)}\bigr)=
\alpha_{1}^{(1)}\;\!\bigl(l_1^{}(x,y)-x\;\!l_3^{}(x,y)\bigr)+
\beta_{1}^{(1)}\;\!\bigl(l_2^{}(x,y)-y\;\!l_3^{}(x,y)\bigr)=
\hfill
$
\\[2ex]
\mbox{}\hfill
$
=
\alpha_{1}^{(1)}\bigl(a_1^{}x+b_1^{}y+c_1^{}\bigr)+
\beta_{1}^{(1)}\bigl(a_2^{}x+b_2^{}y+c_2^{}\bigr)-
\bigl(\alpha_{1}^{(1)}x+\beta_{1}^{(1)}y\bigr)\;\!l_3^{}(x,y)=
\hfill
$
\\[2ex]
\mbox{}\hfill
$
=
\bigl(a_1^{}\alpha_{1}^{(1)}+a_2^{}\beta_{1}^{(1)}\bigr) x+
\bigl(b_1^{}\alpha_{1}^{(1)}+b_2^{}\beta_{1}^{(1)}\bigr) y+
\bigl(c_1^{}\alpha_{1}^{(1)}+c_2^{}\beta_{1}^{(1)}\bigr)-
\bigl(\alpha_{1}^{(1)}x+\beta_{1}^{(1)}y\bigr)\;\!l_3^{}(x,y)
\ \, \forall (x,y)\!\in\!\R^2.
\hfill
$
\\[2ex]
\indent
Отсюда, учитывая систему (12.53), получаем:
\\[2ex]
\mbox{}\hfill
$
{\frak I}\;\!p_1^{(1)}(x,y)=
\bigl(\alpha_1^{}+\lambda_1^{}\;\!\alpha_{1}^{(1)}-a_3^{}\;\!\gamma_{1}^{(1)}\bigr)\;\! x+
\bigl(\beta_1^{}+\lambda_1^{}\;\!\beta_{1}^{(1)}-b_3^{}\;\!\gamma_{1}^{(1)}\bigr)\;\! y+
\bigl(\gamma_1^{}+\lambda_1^{}\;\!\gamma_{1}^{(1)}-c_3^{}\;\!\gamma_{1}^{(1)}\bigr)\ -
\hfill
$
\\[2ex]
\mbox{}\hfill
$
-\ \bigl(\alpha_{1}^{(1)}x+\beta_{1}^{(1)}y\bigr)\;\!l_3^{}(x,y)=
\alpha_1^{}x+\beta_1^{}y+\gamma_1^{}
+\lambda_1^{}\;\! \bigl(\alpha_1^{(1)}x+\beta_1^{(1)}y+\gamma_1^{(1)}\bigr)-
\gamma_{1}^{(1)}(a_3^{}x+b_3^{}y+c_3^{})\ -
\hfill
$
\\[2ex]
\mbox{}\hfill
$
-\ \bigl(\alpha_{1}^{(1)}x+\beta_{1}^{(1)}y\bigr)\;\!l_3^{}(x,y)=
\alpha_1^{}x+\beta_1^{}y+\gamma_1^{}
+\lambda_1^{}\;\! \bigl(\alpha_1^{(1)}x+\beta_1^{(1)}y+\gamma_1^{(1)}\bigr)\ -
\hfill
$
\\[2ex]
\mbox{}\hfill
$
-\ 
\bigl(\alpha_{1}^{(1)}x+\beta_{1}^{(1)}y+\gamma_{1}^{(1)}\bigr)\;\!l_3^{}(x,y)=
\alpha_1^{}x+\beta_1^{}y+\gamma_1^{}+
\bigl(\alpha_{1}^{(1)}x+\beta_{1}^{(1)}y+\gamma_{1}^{(1)}\bigr)
\bigl(\lambda_1^{}-l_3^{}(x,y)\bigr)=
\hfill
$
\\[2ex]
\mbox{}\hfill
$
=p_1^{}(x,y)+p_1^{(1)}(x,y)\;\!M_1^{}(x,y)
\quad 
\forall (x,y)\in\R^2.
\hfill
$
\\[2ex]
\indent
Итак, производная в силу уравнения Якоби (12.36)
\\[2ex]
\mbox{}\hfill  % (12.54)
$
{\frak I}\;\!p_1^{(1)}(x,y)=
p_1^{}(x,y)+p_1^{(1)}(x,y)\;\!M_1^{}(x,y)
\quad 
\forall (x,y)\in\R^2.
$
\hfill (12.54)
\\[2ex]
\indent
Согласно критерию кратного полиномиального частного интеграла 
(теорема 5.2 при $h=1, \ N=1,\ p=p_1^{},\ q=p_1^{(1)},\ M=M_1^{})$
полиномиальный частный интеграл $p_1^{}$ 
уравнения Якоби (12.36) такой, что 
$\bigl((p_1^{}, M_1^{}), (1,p_1^{(1)}, 1)\bigr)\in \text{B}_{_{\!D}}.\ \k$
\vspace{0.5ex}

{\small\bf 12.22.2.1. Случай двукратного элементарного делителя.}
\vspace{0.25ex}
Если кратному собственному числу $\lambda_1^{}$ матрицы $A$
\vspace{0.25ex}
соответствует двукратный элементарный делитель $(\lambda-\lambda_1^{})^2,$
то матрица $A$ имеет и вещественное собственное число  $\lambda_3^{},$
которому соответствует простой элементарный делитель $\lambda-\lambda_3^{}.$
\vspace{0.5ex}

{\bf Теорема 12.12.}
\vspace{0.25ex}
{\it
Пусть $\lambda_1^{}$ --- собственное число матрицы $A,$ которому
соответствуют двукратный элементарный делитель, 
собственный вектор 
\vspace{0.25ex}
$(\alpha_1^{},\beta_1^{},\gamma_1^{})$ и первый 
присоединенный вектор 
\vspace{0.35ex}
$\bigl(\alpha_1^{(1)},\beta_1^{(1)},\gamma_1^{(1)}\bigr),$
а  $\lambda_3^{}$ --- собственное число матрицы $A,$ которому
соответствует собственный вектор 
$(\alpha_3^{},\beta_3^{},\gamma_3^{}).$
\vspace{0.5ex}
Тогда общим интегралом на множестве 
$D=\{(x,y)\colon (\alpha_1^{}x+\beta_1^{}y+\gamma_1^{})
(\alpha_3^{}x+\beta_3^{}y+\gamma_3^{})\ne 0\}$
уравнения Якоби {\rm (12.36)} будет функция
\\[2ex]
\mbox{}\hfill   % (12.55)
$
F\colon (x,y)\to \
\dfrac{p_1^{}(x,y)}{p_3^{}(x,y)}\,
\exp\biggl((\lambda_3^{}-\lambda_1^{})\,\dfrac{p_1^{(1)}(x,y)}{p_1^{}(x,y)}\biggr)
\quad
\forall (x,y)\in D,
$
\hfill {\rm (12.55)}
\\[2ex]
где 
\vspace{0.5ex}
$
p_1^{}(x,y)= \alpha_1^{}x+\beta_1^{}y+\gamma_1^{}\;\; \forall (x,y)\in\R^2,\ 
p_1^{(1)}(x,y)= \alpha_1^{(1)}x+\beta_1^{(1)}y+\gamma_1^{(1)}\;\; \forall (x,y)\in\R^2,
$ 
$
p_3^{}(x,y)= \alpha_3^{}x+\beta_3^{}y+\gamma_3^{}\;\; \forall (x,y)\in\R^2. 
$
}
\vspace{0.35ex}

{\sl Доказательство.}
\vspace{0.25ex}
По лемме 12.2, 
$\bigl((p_1^{}, \lambda_1^{}-l_3^{}), (1,p_1^{(1)}, 1)\bigr)\in \text{B}_{_{\!\scriptstyle D_{0^{}}}},$
где множество
$D_0^{}=
\linebreak
=\{(x,y)\colon \alpha_1^{}x+\beta_1^{}y+\gamma_1^{}\ne 0\}.$
В соответствии с критерием кратного полиномиального частного интеграла 
(теорема 5.3 при $p=p_1^{},\ M=\lambda_1^{}-l_3^{},\ h=1,\ q=p_1^{(1)},\ N=1)$
уравнение Якоби (12.36) 
\vspace{0.35ex}
наряду с полиномиальным частным интегралом $p_1^{}$
с сомножителем $M_1^{}\colon (x,y)\to \lambda_1^{}-l_3^{}(x,y)\;\; \forall (x,y)\in\R^2$
\vspace{0.25ex}
имеет экспоненциальный частный интеграл
$\exp\dfrac{p_1^{(1)}}{p_1^{}}$ 
с сомножителем 
\vspace{0.5ex}
$M_2^{}\colon (x,y)\to 1\;\; \forall (x,y)\in\R^2.$

По лемме 12.1, линейная функция $p_3^{}$ 
\vspace{0.35ex}
является полиномиальным частным интегралом
с сомножителем $M_3^{}\colon (x,y)\to \lambda_3^{}-l_3^{}(x,y)\;\; \forall (x,y)\in\R^2$
\vspace{0.35ex}
уравнения Якоби (12.36).

Сомножители $M_1^{},\ M_2^{},\ M_3^{}$ частных интегралов
$p_1^{},\ \exp\dfrac{p_1^{(1)}}{p_1^{}}\,,\ p_3^{}$
такие, что линейная комбинация
\\[1.5ex]
\mbox{}\hfill
$
M_1^{}(x,y)+(\lambda_3^{}-\lambda_1^{})\;\!M_2^{}(x,y)-M_3^{}(x,y)=0
\quad
\forall (x,y)\in\R^2.
\hfill
$
\\[1.25ex]
\indent
По свойству 11.11, 
\vspace{0.15ex}
функция (12.55) будет общим интегралом на 
множестве $D$ уравнения Якоби (12.36). $\k$
\vspace{0.5ex}

{\bf Пример 12.5.}
У уравнения Якоби [9, с. 25]
\\[2ex]
\mbox{}\hfill  % (12.56)
$
\dfrac{dy}{dx}=\dfrac{x-y-1-y\;\!({}-x+y)}{{}-x+y-x\;\!({}-x+y)}
$
\hfill (12.56)
\\[2ex]
матрица
\vspace{0.35ex}
$
A=\left(\!\!\!
\begin{array}{rrr}
{}-1 & 1 & {}-1
\\
1 & {}-1 & 1
\\
0 & {}-1 & 0\end{array}
\!\!\right)
$
имеет собственное число $\lambda_1^{}={}-1,$ 
которому соответствуют двукратный элементарный делитель $(\lambda+1)^2,$
\vspace{0.25ex}
собственный вектор $(1,{}-1, {}-1),$
первый присоединенный вектор $(1,{}-1, {}-2),$
\vspace{0.25ex}
а также собственное число $\lambda_3^{}=0,$
которому соответствует собственный вектор  $(1, 0, {}-1).$
\vspace{0.35ex}
По теореме 12.12, общим интегралом  на множестве 
$D=\{(x,y)\colon (x-1)(x-y-1)\ne 0\}$ уравнения Якоби (12.56) 
будет функция
\\[1.5ex]
\mbox{}\hfill  % (12.57)
$
F\colon (x,y)\to\ 
\dfrac{x-y-1}{x-1}\,\exp\dfrac{x-y-2}{x-y-1}
\quad
\forall (x,y)\in D.
$
\hfill {\rm (12.57)}
\\[2ex]
\indent
{\bf Пример 12.6.}
У уравнения Якоби 
\\[2ex]
\mbox{}\hfill  % (12.58)
$
\dfrac{dy}{dx}=\dfrac{y-y\;\!(y+1)}{x-x\;\!(y+1)}
$
\hfill (12.58)
\\[2ex]
матрица
\vspace{0.35ex}
$
A=\left(\!\!\!
\begin{array}{ccc}
1 & 0 & 0
\\
0 & 1 & 1
\\
0 & 0 & 1\end{array}
\!\!\right)
$
имеет собственные числа $\lambda_1^{}=\lambda_2^{}=\lambda_3^{}=1,$ 
которым соот\-ветствуют элементарные делители $(\lambda-1)^2$ и $(\lambda-1).$
\vspace{0.25ex}
Собственному числу $\lambda_1^{}=1,$
ко\-то\-рому соответствует кратный элементарный делитель $(\lambda-1)^2,$
\vspace{0.25ex}
соответствуют собст\-венный вектор $(0, 1, 0)$ и 
первый присоединенный вектор $(1, 1, 1).$
Собственному числу $\lambda_3^{}=1,$
которому соответствует простой элементарный делитель $\lambda-1,$
соответствует собственный вектор $(1, 0, 0).$
\vspace{0.25ex}
По теореме 12.12 (при $\lambda_3^{}=\lambda_1^{}=1),$ 
общим интегралом  на множестве 
$\!D\!=\!\{(x,y)\colon x\!\ne\! 0\}\!$ уравнения Якоби (12.58) 
будет дробно-линейная функция
\\[1.5ex]
\mbox{}\hfill  % (12.59)
$
F\colon (x,y)\to\ 
\dfrac{y}{x}
\quad
\forall (x,y)\in D.
$
\hfill {\rm (12.59)}
\\[2ex]
\indent
{\small\bf 
12.22.2.2. Случай трехкратного элементарного делителя.}
Если $\lambda_1^{}$ --- собственное число матрицы $A,$ которому соответствует 
\vspace{0.25ex}
трехкратный элементарный делитель $(\lambda-\lambda_1^{})^3,$ то ему кроме собственного вектора
\vspace{0.5ex}
$\theta_1^{}=(\alpha_1^{},\beta_1^{},\gamma_1^{})$
соответствуют первый присоединенный вектор 
\vspace{0.35ex}
$\theta_1^{(1)}=\bigl(\alpha_1^{(1)},\beta_1^{(1)},\gamma_1^{(1)}\bigr)$
и второй присоединенный вектор 
$\theta_1^{(2)}=\bigl(\alpha_1^{(2)},\beta_1^{(2)},\gamma_1^{(2)}\bigr).$
\vspace{0.35ex}
Первый присоединенный вектор $\theta_1^{(1)}$ является решением 
матричного уравнения (12.52). Второй присоединенный вектор $\theta_1^{(2)}$
является решением матричного уравнения
\\[1.5ex]
\mbox{}\hfill % (12.60)
$
(A-\lambda_1^{}E)\;\! {\rm colon}\,\theta_1^{(2)} =
2\, {\rm colon}\;\!\theta_1^{(1)}.
$
\hfill (12.60)
\\[2ex]
\indent
{\bf Теорема 12.13.}
\vspace{0.25ex}
{\it
Пусть собственному числу $\lambda_1^{}$ матрицы $A$ соответствуют 
трехкратный элементарный делитель, 
собственный вектор 
\vspace{0.5ex}
$(\alpha_1^{},\beta_1^{},\gamma_1^{}),$ 
первый присоединенный вектор 
\vspace{0.5ex}
$\bigl(\alpha_1^{(1)},\beta_1^{(1)},\gamma_1^{(1)}\bigr)$ 
и второй присоединенный вектор 
$\bigl(\alpha_1^{(2)},\beta_1^{(2)},\gamma_1^{(2)}\bigr).$
Тогда общим интегралом на множестве 
\vspace{0.35ex}
$D=\{(x,y)\colon \alpha_1^{}x+\beta_1^{}y+\gamma_1^{}\ne 0\}$
уравнения Якоби {\rm (12.36)} будет рациональная функция
\\[2ex]
\mbox{}\hfill   % (12.61)
$
F\colon (x,y)\to \
\dfrac{\bigl(p_1^{(1)}(x,y)\bigr)^2-p_1^{}(x,y)\;\!p_1^{(2)}(x,y)}{p_1^{2}(x,y)}
\quad
\forall (x,y)\in D,
$
\hfill {\rm (12.61)}
\\[2ex]
где 
\vspace{0.75ex}
$
p_1^{}(x,y)= \alpha_1^{}x+\beta_1^{}y+\gamma_1^{}\;\; \forall (x,y)\in\R^2,\ 
p_1^{(1)}(x,y)= \alpha_1^{(1)}x+\beta_1^{(1)}y+\gamma_1^{(1)}\;\; \forall (x,y)\in\R^2,
$ 
$
p_1^{(2)}(x,y)= \alpha_1^{(2)}x+\beta_1^{(2)}y+\gamma_1^{(2)}\;\; \forall (x,y)\in\R^2. 
$
}
\vspace{0.5ex}

{\sl Доказательство.}
\vspace{0.5ex}
По лемме 12.2, 
$\bigl((p_1^{}, M_1^{}), (1,p_1^{(1)}, 1)\bigr)\in \text{B}_{_{\!\scriptstyle D}},$
где сомножитель 
$M_1^{}\colon (x,y)\to \lambda_1^{}-l_3^{}(x,y)\;\; \forall (x,y)\in\R^2.$
\vspace{0.25ex}
Согласно критерию существования кратного полиномиального частного интеграла 
(теорема 5.1 при $p=p_1^{},\ M=M_1^{},\ h=1,\ q=p_1^{(1)},$ $N=1)$
выполняются тождества
\\[1.5ex]
\mbox{}\hfill  % (12.62)
$
{\frak I}\;\!p_1^{}(x,y)=p_1^{}(x,y)\;\!M_1^{}(x,y)
\quad
\forall (x,y)\in\R^2
$
\hfill (12.62)
\\[1ex]
и 
\\[1ex]
\mbox{}\hfill  % (12.63)
$
{\frak I}\,\dfrac{p_1^{(1)}(x,y)}{p_1^{}(x,y)}=1
\quad
\forall (x,y)\in D.
$
\hfill (12.63)
\\[1.5ex]
\indent
Матричное уравнение (12.60) равносильно линейной системе
\\[1.5ex]
\mbox{}\hfill              % (12.64)               
$\!\!\!\!\!
\begin{array}{c}
(a_1^{}\!-\!\lambda_1^{})\;\!\alpha_{1}^{(2)}\!+
a_2^{}\;\!\beta_{1}^{(2)}\!+
a_3^{}\;\!\gamma_{1}^{(2)}=2\alpha_1^{(1)},
\\[1.5ex]
b_1^{}\;\!\alpha_{1}^{(2)}\!+
(b_2^{}\!-\!\lambda_1^{})\;\!\beta_{1}^{(2)}\!+
b_3^{}\;\!\gamma_{1}^{(2)}=2\beta_1^{(1)},
\\[1.5ex]
c_1^{}\;\!\alpha_{1}^{(2)}\!+
c_2^{}\;\!\beta_{1}^{(2)}\!+
(c_3^{}\!-\!\lambda_1^{})\;\!\gamma_{1}^{(2)}=2\gamma_1^{(1)}
\end{array}
\!\!\Leftrightarrow\!\!
%\!\!\!\!\iff\!\!\!\!
\begin{array}{c}
a_1^{}\;\!\alpha_{1}^{(2)}\!+
a_2^{}\;\!\beta_{1}^{(2)}=
2\alpha_1^{(1)}\!+\lambda_1^{}\;\!\alpha_{1}^{(2)}\!-a_3^{}\;\!\gamma_{1}^{(2)},
\\[1.5ex]
b_1^{}\;\!\alpha_{1}^{(2)}\!+
b_2^{}\;\!\beta_{1}^{(2)}= 
2\beta_1^{(1)}\!+\lambda_1^{}\;\!\beta_{1}^{(2)}\!-b_3^{}\;\!\gamma_{1}^{(2)},
\\[1.5ex]
c_1^{}\;\!\alpha_{1}^{(2)}\!+
c_2^{}\;\!\beta_{1}^{(2)}=
2\gamma_1^{(1)}\!+\lambda_1^{}\;\!\gamma_{1}^{(2)}\!-c_3^{}\;\!\gamma_{1}^{(2)}.
\end{array}\!\!\!
$
\hfill (12.64)
\\[2ex]
\indent
Производная в силу уравнения Якоби (12.36)
\\[2ex]
\mbox{}\hfill
$
{\frak I}\;\!p_1^{(2)}(x,y)=
{\frak I}\;\!\bigl(\alpha_1^{(2)}x+\beta_1^{(2)}y+\gamma_1^{(2)}\bigr)=
\alpha_{1}^{(2)}\;\!\bigl(l_1^{}(x,y)-x\;\!l_3^{}(x,y)\bigr)+
\beta_{1}^{(2)}\;\!\bigl(l_2^{}(x,y)-y\;\!l_3^{}(x,y)\bigr)=
\hfill
$
\\[2ex]
\mbox{}\hfill
$
=
\alpha_{1}^{(2)}\bigl(a_1^{}x+b_1^{}y+c_1^{}\bigr)+
\beta_{1}^{(2)}\bigl(a_2^{}x+b_2^{}y+c_2^{}\bigr)-
\bigl(\alpha_{1}^{(2)}x+\beta_{1}^{(2)}y\bigr)\;\!l_3^{}(x,y)=
\hfill
$
\\[2ex]
\mbox{}\hfill
$
=
\bigl(a_1^{}\alpha_{1}^{(2)}+a_2^{}\beta_{1}^{(2)}\bigr) x+
\bigl(b_1^{}\alpha_{1}^{(2)}+b_2^{}\beta_{1}^{(2)}\bigr) y+
\bigl(c_1^{}\alpha_{1}^{(2)}+c_2^{}\beta_{1}^{(2)}\bigr)-
\bigl(\alpha_{1}^{(2)}x+\beta_{1}^{(2)}y\bigr)\;\!l_3^{}(x,y)
\ \, \forall (x,y)\!\in\!\R^2.
\hfill
$
\\[2ex]
\indent
Отсюда, учитывая систему (12.64), получаем:
\\[2ex]
\mbox{}\hfill
$
{\frak I}\;\!p_1^{(2)}(x,y)=
\bigl(2\;\!\alpha_1^{(1)}\!+\lambda_1^{}\;\!\alpha_{1}^{(2)}-a_3^{}\;\!\gamma_{1}^{(2)}\bigr)\;\! x+
\bigl(2\beta_1^{(1)}\!+\lambda_1^{}\;\!\beta_{1}^{(2)}-b_3^{}\;\!\gamma_{1}^{(2)}\bigr)\;\! y+
\bigl(2\gamma_1^{(1)}\!+\lambda_1^{}\;\!\gamma_{1}^{(2)}-c_3^{}\;\!\gamma_{1}^{(2)}\bigr)\ -
\hfill
$
\\[2ex]
\mbox{}\hfill
$
- \bigl(\alpha_{1}^{(2)}x+\beta_{1}^{(2)}y\bigr)\;\!l_3^{}(x,y)\!=\!
2\bigl(\alpha_1^{(1)}x+\beta_1^{(1)}y+\gamma_1^{(1)}\bigr)
+\lambda_1^{}\bigl(\alpha_1^{(2)}x+\beta_1^{(2)}y+\gamma_1^{(2)}\bigr)-
\gamma_{1}^{(2)}(a_3^{}x+b_3^{}y+c_3^{}) -
\hfill
$
\\[2ex]
\mbox{}\hfill
$
-\ \bigl(\alpha_{1}^{(2)}x+\beta_{1}^{(2)}y\bigr)\;\!l_3^{}(x,y)=
2\bigl(\alpha_1^{(1)}x+\beta_1^{(1)}y+\gamma_1^{(1)}\bigr)
+\lambda_1^{}\;\! \bigl(\alpha_1^{(2)}x+\beta_1^{(2)}y+\gamma_1^{(2)}\bigr)\ -
\hfill
$
\\[2ex]
\mbox{}\hfill
$
-\ 
\bigl(\alpha_{1}^{(2)}x+\beta_{1}^{(2)}y+\gamma_{1}^{(2)}\bigr)\;\!l_3^{}(x,y)=
2\bigl(\alpha_1^{(1)}x+\beta_1^{(1)}y+\gamma_1^{(1)}\bigr)+
\bigl(\alpha_{1}^{(2)}x+\beta_{1}^{(2)}y+\gamma_{1}^{(2)}\bigr)
\bigl(\lambda_1^{}-l_3^{}(x,y)\bigr)=
\hfill
$
\\[2ex]
\mbox{}\hfill
$
=2\;\!p_1^{(1)}(x,y)+p_1^{(2)}(x,y)\;\!M_1^{}(x,y)
\quad 
\forall (x,y)\in\R^2.
\hfill
$
\\[2ex]
\indent
Итак, производная в силу уравнения Якоби (12.36)
\\[2ex]
\mbox{}\hfill  % (12.65)
$
{\frak I}\;\!p_1^{(2)}(x,y)=
2\;\!p_1^{(1)}(x,y)+p_1^{(2)}(x,y)\;\!M_1^{}(x,y)
\quad 
\forall (x,y)\in\R^2.
$
\hfill (12.65)
\\[2ex]
\indent
Учитывая тождества (12.62) и (12.65), на множестве $D$ находим:
\\[2ex]
\mbox{}\hfill  
$
{\frak I}\,\dfrac{p_1^{(2)}(x,y)}{p_1^{}(x,y)}=
\dfrac{p_1^{}(x,y)\,{\frak I}\;\!p_1^{(2)}(x,y)-
p_1^{(2)}(x,y)\,{\frak I}\;\!p_1^{}(x,y)}{p_1^{2}(x,y)}=
\hfill  
$
\\[2ex]
\mbox{}\hfill  
$
=\dfrac{p_1^{}(x,y)\bigl(2\;\!p_1^{(1)}(x,y)+p_1^{(2)}(x,y)\;\!M_1^{}(x,y)\bigr)-
p_1^{}(x,y)\;\!p_1^{(2)}(x,y)\;\!M_1^{}(x,y)}{p_1^{2}(x,y)}=
2\,\dfrac{p_1^{(1)}(x,y)}{p_1^{}(x,y)}\,.
\hfill
$
\\[2ex]
\indent
Из тождества (12.63) следует, что
\\[1.5ex]
\mbox{}\hfill  
$
{\frak I}\,\biggl(\dfrac{p_1^{(1)}(x,y)}{p_1^{}(x,y)}\biggr)^{\!2}=
2\,\dfrac{p_1^{(1)}(x,y)}{p_1^{}(x,y)}
\quad 
\forall (x,y)\in D.
\hfill
$
\\[1.75ex]
\indent
Стало быть, производная в силу уравнения Якоби (12.36)
\\[1.75ex]
\mbox{}\hfill   
$\!\!
{\frak I}F(x,y)=
{\frak I}\,\dfrac{\bigl(p_1^{(1)}(x,y)\bigr)^2\!-p_1^{}(x,y)\;\!p_1^{(2)}(x,y)}{p_1^{2}(x,y)}=
{\frak I}\biggl(\dfrac{p_1^{(1)}(x,y)}{p_1^{}(x,y)}\biggr)^{\!2}\!-
{\frak I}\,\dfrac{p_1^{(2)}(x,y)}{p_1^{}(x,y)}=0
\ \,
\forall (x,y)\!\in\! D.
\hfill
$
\\[1.75ex]
\indent
Следовательно, рациональная функция (12.61) является общим интегралом на 
множестве $D$ уравнения Якоби (12.36). $\k$
\vspace{0.5ex}

{\bf Замечание 12.2.}
Поскольку производная в силу уравнения Якоби (12.36)
\\[1.75ex]
\mbox{}\hfill   
$
{\frak I}\,\dfrac{\bigl(p_1^{(1)}(x,y)\bigr)^2\!-p_1^{}(x,y)\;\!p_1^{(2)}(x,y)}{p_1^{2}(x,y)}=0
\quad
\forall (x,y)\in D,
\hfill
$
\\[1.75ex]
то согласно определению 5.1 полиномиальный частный интеграл $p_1^{}$
\vspace{0.35ex}
уравнения Якоби (12.36) является кратным таким, что 
\vspace{0.25ex}
$\bigl((p_1^{}, M_1^{}), \bigl(2,\bigl( p_1^{(1)}\bigr)^2-p_1^{}p_1^{(2)}, 0\bigr)\bigr)\in 
\text{B}_{_{\!\scriptstyle D}}.$
Кроме этого, по лемме 12.2,
\vspace{0.25ex}
$\bigl((p_1^{}, M_1^{}), (1,p_1^{(1)}, 1)\bigr)\in \text{B}_{_{\!\scriptstyle D}}.$
Тогда в соответствии с определением 5.2
полиномиальный частный интеграл $p_1^{}$ уравнения Якоби (12.36)
будет двукратным.
\vspace{0.25ex}

Следовательно, 
\vspace{0.15ex}
{\it 
если матрица $A$ имеет собственное число $\lambda_1^{},$ 
которому соответствует трехкратный элементарный делитель и 
собственный вектор 
\vspace{0.5ex}
$(\alpha_1^{},\beta_1^{},\gamma_1^{}),$ 
то полиномиальный частный интеграл
\vspace{0.15ex}
$p_1^{}\colon (x,y)\to \alpha_1^{}x+\beta_1^{}y+\gamma_1^{}\;\; \forall (x,y)\in\R^2$ 
уравнения Якоби {\rm (12.36)} будет двукратным.
}
\vspace{0.5ex}

{\bf Пример 12.7.}
У уравнения Якоби 
\\[2ex]
\mbox{}\hfill  % (12.66)
$
\dfrac{dy}{dx}=\dfrac{x+3y+\delta-y\;\!(x+y+2)}{x-y-\delta-x\;\!(x+y+2)}
$
\hfill (12.66)
\\[2ex]
с вещественным параметром $\delta\ne 0$
матрица
\vspace{0.35ex}
$
A=\left(\!\!\!
\begin{array}{rcc}
1 & 1 & 1
\\
{}-1 & 3 & 1
\\
{}-\delta & \delta & 2\end{array}
\!\!\right)
$
имеет собственное число $\lambda=2,$ 
\vspace{0.35ex}
которому соот\-ветствуют трехкратный элементарный делитель $(\lambda-2)^3,$ 
собст\-венный вектор $(1, 1, 0),$  
первый присоединенный вектор $(0, 0, 1)$ 
\vspace{0.75ex}
и второй присоединенный вектор $\Bigl({}-\dfrac{2}{\delta}\,,\;\! 0, {}-\dfrac{2}{\delta}\Bigr).$ 
\vspace{0.75ex}
По теореме 12.13, общим интегралом  на множестве 
$D=\{(x,y)\colon x+y\ne 0\}$ уравнения Якоби (12.66) 
будет рациональная функция
\\[1.5ex]
\mbox{}\hfill  
$
F^{\ast}\colon (x,y)\to\ 
\dfrac{1-(x+y)\Bigl({}-\dfrac{2}{\delta}\,x-\dfrac{2}{\delta}\Bigr)}{(x+y)^2}
\quad
\forall (x,y)\in D
\hfill 
$
\\[1.5ex]
или (с учетом свойства 0.1 о функциональной неоднозначности общего интеграла) функция
\\[1.5ex]
\mbox{}\hfill  % (12.67)
$
F\colon (x,y)\to\ 
\dfrac{2\;\!(x+y)(x+1)+\delta}{(x+y)^2}
\quad
\forall (x,y)\in D.
$
\hfill  (12.67)
\\[1.5ex]
\indent
Таким образом, для уравнения Якоби (12.36) в результате решения задачи Дарбу 
(приложение 12.1) общий интеграл построен без квадратур. 
В [9, с. 14 -- 25] для уравнения Якоби (12.36) решена расширенная задача Дарбу 
(приложения 12.1 и 12.9), когда в случае простых элементарных делителей у 
матрицы $A$ строится общий интеграл, а при наличии кратного элементарного делителя
указывается интегрирующий множитель.
\vspace{1ex}

{\bf 12.23.}
{\it Система Якоби}
\\[2ex]
\mbox{}\hfill        % (12.68)
$
\dfrac{dx}{dt}=l_1^{}(x,y)-x\;\!l_3^{}(x,y),
\quad
\dfrac{dy}{dt}=l_2^{}(x,y)-y\;\!l_3^{}(x,y),
\hfill
$
\\
\mbox{}\hfill (12.68)
\\
\mbox{}\hfill
$
l_i^{}(x,y)=a_{i}^{}\;\!x+b_{i}^{}\;\!y+c_i^{}\;\!,
\quad
a_{i}^{},b_{i}^{},c_i^{}\in\R, 
\quad
i=1,2,3,
\hfill
$
\\[2ex]
будучи автономной дифференциальной системой второго порядка, 
имеет два функционально независимых первых интеграла, один из которых
является неавтономным.

Уравнение Якоби (12.36) есть уравнение траекторий системы Якоби (12.68). 
Поэтому общий интеграл уравнения Якоби (12.36) является автономным 
первым интегралом системы Якоби (12.68).

Таким образом, для построения интегрального базиса системы Якоби (12.68)
достаточно найти ее неавтономный первый интеграл.

Оператором дифференцирования в силу системы Якоби (12.68) является
линейный дифференциальный оператор первого порядка
\\[1.5ex]
\mbox{}\hfill
$
{\frak d}(t,x,y)=
\partial_{{}_{\scriptstyle t}}+ {\frak I}(x,y)
\quad
\forall (t,x,y)\in\R^3,
\hfill
$
\\[1.5ex]
где ${\frak I}$ --- оператор дифференцирования в силу уравнения Якоби (12.36).

Связь операторов дифференцирования ${\frak d}$ и ${\frak I}$ позволяет 
частные интегралы уравнения Якоби (12.36) трактовать как автономные 
\vspace{0.5ex}
частные интегралы системы Якоби (12.68).

{\small\bf 
12.23.1. Случай различных вещественных собственных чисел.}
\vspace{0.35ex}

{\bf Теорема 12.14.}
\vspace{0.25ex}
{\it
Пусть $\lambda_1^{}$ и $\lambda_2^{}$ --- 
различные вещественные собственные числа матрицы $A,$
которым соответствуют собственные векторы
\vspace{0.35ex}
$
(\alpha_1^{},\beta_1^{},\gamma_1^{})
$ 
и
$
(\alpha_2^{},\beta_2^{},\gamma_2^{}). 
$
Тогда неавтономным 
\vspace{0.15ex}
первым интегралом на множестве 
$\!\Omega_0^{}\!=\!\{(t,x,y)\colon\! \alpha_2^{}x+\beta_{2}^{}y+\gamma_2^{}\!\ne\! 0\}\!$ 
системы Якоби {\rm(12.68)} будет функция
\\[2ex]
\mbox{}\hfill  % (12.69)
$
\Psi\colon (t,x,y)\to\ 
\dfrac{p_1^{}(x,y)}{p_2^{}(x,y)}\,\exp\bigl((\lambda_2^{}-\lambda_1^{})\;\!t\bigr)
\quad
\forall (t,x,y)\in \Omega_0^{},
$
\hfill {\rm (12.69)}
\\[2.25ex]
где $p_1^{}(x,y)=\alpha_1^{}\;\!x+\beta_{1}^{}\;\!y+\gamma_1^{}\;\;\forall (x,y)\in\R^2,\ 
p_2^{}(x,y)=\alpha_2^{}\;\!x+\beta_{2}^{}\;\!y+\gamma_2^{}\;\;\forall (x,y)\in\R^2.$
}
\vspace{0.5ex}

{\sl Доказательство.}
\vspace{0.25ex}
По лемме 12.1, линейные функции $p_1^{}$ и $p_2^{}$ являются автономными полиномиальными 
частными интегралами с сомножителями
\vspace{0.35ex}
$M_1^{}\colon (x,y)\to \lambda_1^{}-l_3^{}(x,y)$ $\forall (x,y)\in\R^2$ и
$M_2^{}\colon (x,y)\to \lambda_2^{}-l_3^{}(x,y)\;\;\forall (x,y)\in\R^2$ 
системы Якоби (12.68).
\vspace{0.35ex}

Согласно свойству 4.1 
\vspace{0.15ex}
функция $e^t$ является неавтономным условным частным интегралом
с сомножителем $M_3^{}\colon t\to 1\;\;\forall t\in\R$ системы Якоби (12.68).
\vspace{0.25ex}

Так как у частных интегралов $p_1^{},\ p_2^{},\ e^t$ 
сомножители $M_1^{},\ M_2^{},\ M_3^{}$ такие, что выполняется тождество
\\[1.5ex]
\mbox{}\hfill
$
M_1^{}(x,y)-M_2^{}(x,y)+(\lambda_2^{}-\lambda_1^{})M_1^{}(t)=0
\quad
\forall (t,x,y)\in \R^3,
\hfill
$ 
\\[1.5ex]
то (свойство 11.11) 
\vspace{0.15ex}
функция (12.69) будет неавтономным первым интегралом на 
множестве $\Omega_0^{}$ системы Якоби (12.68). $\k$
\vspace{0.5ex}

{\bf Пример 12.8.}
У системы Якоби
\\[2ex]
\mbox{}\hfill    % (12.70)
$
\dfrac{dx}{dt}=3x-y+1-x\;\!(x-y+3),
\quad
\dfrac{dy}{dt}={}-x+5y-1-y\;\!(x-y+3)
$
\hfill (12.70)
\\[2ex]
матрица $A$ имеет (пример 12.1) различные вещественные собственные числа
$\lambda_1^{}=2,$ $\lambda_2^{}\!=\!3,\, \lambda_3^{}\!=\!6,\!$
которым соответствуют собственные векторы 
$\!(1, 0, -1),\, (1, 1, 1),\, (1, -2, 1).$
По теореме 12.14, неавтономными первыми интегралами на соответствующих 
множествах системы Якоби (12.70) будут функции
\\[2ex]
\mbox{}\hfill  % (12.71)
$
\Psi_{1{,}2}^{}\colon (t,x,y)\to\ 
\dfrac{x-1}{x+y+1}\, e^t
\quad
\forall (t,x,y)\in\Omega_1^{},
$
\hfill (12.71)
\\[2ex]
\mbox{}\hfill  % (12.72)
$
\Psi_{1{,}3}^{}\colon (t,x,y)\to\ 
\dfrac{x-1}{x-2y+1}\, e^{\;\!4t}
\quad
\forall (t,x,y)\in\Omega_2^{},
$
\hfill (12.72)
\\[2ex]
\mbox{}\hfill  % (12.73)
$
\Psi_{2{,}3}^{}\colon (t,x,y)\to\ 
\dfrac{x+y+1}{x-2y+1}\, e^{\;\!3t}
\quad
\forall (t,x,y)\in\Omega_2^{},
$
\hfill (12.73)
\\[2ex]
где множества $\Omega_1^{}=\{(t,x,y)\colon x+y+1\ne 0\},\ 
\Omega_2^{}=\{(t,x,y)\colon x-2y+1\ne 0\}.$
\vspace{0.35ex}

Каждая из совокупностей $(12.43)\& (12.71),\ (12.43)\& (12.72),\ (12.43)\& (12.73),\ (12.71)\& 
\linebreak 
\& (12.72),\ (12.71)\& (12.73),\ (12.72)\& (12.73)$
образует интегральный базис на соответствующем множестве системы Якоби (12.70).
\vspace{0.5ex}

{\bf Пример 12.9.}
У системы Якоби
\\[2ex]
\mbox{}\hfill    % (12.74)
$
\dfrac{dx}{dt}={}-x+y+1-x\;\!(x+y-1),
\quad
\dfrac{dy}{dt}=x-y+1-y\;\!(x+y-1)
$
\hfill (12.74)
\\[2ex]
матрица $A$ имеет (пример 12.2) 
\vspace{0.15ex}
различные вещественные собственные числа
$\lambda_1^{}={}-2$ и $\lambda_3^{}=1.$
\vspace{0.25ex}
Собственному числу $\lambda_1^{}={}-2$
соответствуют два линейно независимых собственных вектора
$(0, 1, -1)$ и $(1, 0, {}-1).$
\vspace{0.15ex}
Собственному числу $\lambda_3^{}=1$
соответствует собственный вектор $(1, 1, 1).$
\vspace{0.25ex}
По теореме 12.14, неавтономными первыми интегралами на  
множестве $\Omega_0^{}=\{(t,x,y)\colon x+y+1\ne 0\}$ 
системы Якоби (12.74) будут функции
\\[2ex]
\mbox{}\hfill  % (12.75)
$
\Psi_{1}^{}\colon (t,x,y)\to\ 
\dfrac{y-1}{x+y+1}\, e^{\;\!3t}
\quad
\forall (t,x,y)\in\Omega_0^{},
$
\hfill (12.75)
\\[2ex]
\mbox{}\hfill  % (12.76)
$
\Psi_{2}^{}\colon (t,x,y)\to\ 
\dfrac{x-1}{x+y+1}\, e^{\;\!3t}
\quad
\forall (t,x,y)\in\Omega_0^{}.
$
\hfill (12.76)
\\[2ex]
\indent
Каждая из совокупностей $(12.46)\& (12.75),\ (12.46)\& (12.76),\ (12.75)\& (12.76)$
образует интегральный базис на соответствующем множестве системы Якоби (12.74).
\vspace{0.5ex}

{\small\bf 
12.23.2. Случай комплексного собственного числа.}
\vspace{0.35ex}

{\bf Теорема 12.15.}
\vspace{0.25ex}
{\it
Пусть собственное число $\lambda=\xi+\zeta\;\! i$  матрицы $A$
существенно комплексное $(\zeta\ne 0),$
которому соответствует собственный вектор
\vspace{0.35ex}
$
(\alpha,\beta,\gamma).
$ 
Тогда неавтономным 
\vspace{0.15ex}
первым интегралом на множестве 
$\Omega_0^{}=\{(t,x,y)\colon {\rm Re}(\alpha x+\beta y+\gamma)\ne 0\}$ 
системы Якоби {\rm(12.68)} будет функция
\\[2ex]
\mbox{}\hfill  % (12.77)
$
\Psi\colon (t,x,y)\to\ 
\arctg\dfrac{{\rm Im}\;\!p(x,y)}{{\rm Re}\;\!p(x,y)}\,-\,\zeta\;\!t
\quad
\forall (t,x,y)\in \Omega_0^{},
$
\hfill {\rm (12.77)}
\\[2.25ex]
где $p(x,y)=\alpha\;\!x+\beta\;\!y+\gamma\;\;\forall (x,y)\in\R^2.$
}
\vspace{0.5ex}

{\sl Доказательство.}
\vspace{0.25ex}
Как и при доказательстве теоремы 12.11 устанавливаем, что функция
$\exp\,\arctg\dfrac{{\rm Im}\;\!p}{{\rm Re}\;\!p}$
\vspace{0.5ex}
является автономным экспоненциальным частным интегралом с 
сомножителем $M\colon (x,y)\to \zeta\;\;\forall (x,y)\in\R^2$
на множестве $\Omega_0^{}$ системы Якоби (12.68).
\vspace{0.5ex}

По свойству 4.1, 
\vspace{0.15ex}
функция $e^t$ является неавтономным условным частным интегралом
с сомножителем $M_1^{}\colon t\to 1\;\;\forall t\in\R$ системы Якоби (12.68).
\vspace{0.25ex}

Так как у частных интегралов 
\vspace{0.25ex}
$\exp\,\arctg\dfrac{{\rm Im}\;\!p}{{\rm Re}\;\!p}$ и $e^t$ 
сомножители $M$ и $M_1^{}$ такие, что выполняется тождество
\\[1.25ex]
\mbox{}\hfill
$
M(x,y)-\zeta\;\!M_1^{}(t) = 0
\quad
\forall (t,x,y)\in \R^3,
\hfill
$ 
\\[1.25ex]
то (свойство 11.11 с учетом свойства 0.1 о функциональной неоднозначности 
первого интеграла) функция (12.77) будет неавтономным первым интегралом на 
множестве $\Omega_0^{}$ системы Якоби (12.68). $\k$
\vspace{0.5ex}

{\bf Пример 12.10.}
У системы Якоби
\\[2ex]
\mbox{}\hfill    % (12.78)
$
\dfrac{dx}{dt}=4x-3y-1-x\;\!({}-2x+y),
\quad
\dfrac{dy}{dt}=6x-2y+1-y\;\!({}-2x+y)
$
\hfill (12.78)
\\[2ex]
матрица $A$ имеет (пример 12.3) 
\vspace{0.15ex}
существенно комплексное собственное число
$\lambda=1+\sqrt{6}\;\!i,$ которому
соответствует  собственный вектор
$\bigl({}-2\bigl(2+\sqrt{6}\;\!i\bigr),\;\! 5,\;\! 3-\sqrt{6}\;\!i\bigr).$
\vspace{0.25ex}
По теореме 12.15, неавтономным первым интегралом на  
множестве $\Omega_0^{}=\{(t,x,y)\colon 4x-5y-3\ne 0\}$ 
системы Якоби (12.78) будет функция
\\[2ex]
\mbox{}\hfill  % (12.79)
$
\Psi\colon (t,x,y)\to\ 
\arctg\dfrac{\sqrt{6}\,(2x+1)}{4x-5y-3}\,-\,\sqrt{6}\,t
\quad
\forall (t,x,y)\in\Omega_0^{}.
$
\hfill (12.79)
\\[2ex]
\indent
Совокупность $(12.49)\& (12.79)$
\vspace{0.35ex}
образует интегральный базис системы Якоби (12.78) на множестве 
$\Omega^{\ast}=\{(t,x,y)\colon (x+y+5)(4x-5y-3)\ne 0\}.$
\vspace{1ex}

{\bf Пример 12.11.}
У системы Якоби
\\[2ex]
\mbox{}\hfill    % (12.80)
$
\dfrac{dx}{dt}=x+2y+3-x\;\!(x+2y+1),
\quad
\dfrac{dy}{dt}=x+y-3-y\;\!(x+2y+1)
$
\hfill (12.80)
\\[2ex]
матрица $A$ имеет (пример 12.4) 
\vspace{0.15ex}
существенно комплексное собственное число
$\lambda=1+i,$ которому
соответствует  собственный вектор
$(3-i,\;\! 4, {}-3+3i).$
\vspace{0.35ex}
По теореме 12.15, не\-ав\-то\-номным первым интегралом на  
множестве $\Omega_1^{}=\{(t,x,y)\colon 3x+4y-3\ne 0\}$ 
\vspace{0.25ex}
системы Якоби (12.80) будет функция
\\[2ex]
\mbox{}\hfill  % (12.81)
$
\Psi\colon (t,x,y)\to\ 
\arctg\dfrac{{}-x+3}{3x+4y-3}\,-\,t
\quad
\forall (t,x,y)\in\Omega_1^{}.
$
\hfill (12.81)
\\[2ex]
\indent
Совокупность $(12.51)\& (12.81)$
\vspace{0.35ex}
образует интегральный базис системы Якоби (12.80) на множестве 
$\Omega_0^{}=\{(t,x,y)\colon (x+y-1)(3x+4y-3)\ne 0\}.$
\vspace{1.5ex}

{\small\bf 
12.23.3. Случай кратного элементарного делителя.}
\vspace{0.35ex}

{\bf Теорема 12.16.}
\vspace{0.25ex}
{\it
Пусть вещественному собственному числу матрицы $A$ соответствуют 
кратный элементарный делитель, 
собственный вектор 
\vspace{0.5ex}
$(\alpha_1^{},\beta_1^{},\gamma_1^{})$ и 
первый присоединенный вектор 
\vspace{0.75ex}
$\bigl(\alpha_1^{(1)},\beta_1^{(1)},\gamma_1^{(1)}\bigr).$ 
Тогда неавтономным первым интегралом на множестве 
$\Omega_0^{}=\{(t,x,y)\colon \alpha_1^{}x+\beta_1^{}y+\gamma_1^{}\ne 0\}$
системы Якоби {\rm (12.68)} будет функция
\\[2ex]
\mbox{}\hfill   % (12.82)
$
\Psi\colon (x,y)\to \
\dfrac{p_1^{(1)}(x,y)}{p_1^{}(x,y)}\,-\,t
\quad
\forall (t,x,y)\in \Omega_0^{},
$
\hfill {\rm (12.82)}
\\[2ex]
где 
\vspace{1.25ex}
$
p_1^{}(x,y)= \alpha_1^{}x+\beta_1^{}y+\gamma_1^{}\;\; \forall (x,y)\in\R^2,\ 
p_1^{(1)}(x,y)= \alpha_1^{(1)}x+\beta_1^{(1)}y+\gamma_1^{(1)}\;\; \forall (x,y)\in\R^2.
$}

{\sl Доказательство.}
\vspace{0.25ex}
Как и при доказательстве теоремы 12.12 устанавливаем, что функция
$\exp\dfrac{p_1^{(1)}}{p_1^{}}$
\vspace{0.75ex}
является автономным экспоненциальным частным интегралом с 
сомножителем $M_1^{}\colon (x,y)\to 1\;\;\forall (x,y)\in\R^2$
на множестве $\Omega_0^{}$ системы Якоби (12.68).
\vspace{0.75ex}

По свойству 4.1, 
\vspace{0.35ex}
функция $e^t$ является неавтономным условным частным интегралом
с сомножителем $M_1^{(1)}\colon t\to 1\;\;\forall t\in\R$ системы Якоби (12.68).
\vspace{0.5ex}

Согласно следствию 11.7  
\vspace{0.15ex}
(с учетом свойства 0.1 о функциональной неоднозначности 
первого интеграла) функция (12.82) 
\vspace{0.25ex}
будет неавтономным первым интегралом на 
множестве $\Omega_0^{}$ системы Якоби (12.68). $\k$
\vspace{0.5ex}

\newpage

{\bf Пример 12.12.}
У системы Якоби
\\[1.75ex]
\mbox{}\hfill    % (12.83)
$
\dfrac{dx}{dt}={}-x+y-x\;\!({}-x+y),
\quad
\dfrac{dy}{dt}=x-y-1-y\;\!({}-x+y)
$
\hfill (12.83)
\\[1.75ex]
матрица $A$ имеет (пример 12.5) 
\vspace{0.15ex}
различные вещественные собственные числа
$\lambda_1^{}={}-1$ и $\lambda_3^{}=0,$
\vspace{0.25ex}
которым соответствуют  собственные векторы
$(1, {}-1, {}-1)$ и $(1, 0, {}-1).$
\vspace{0.25ex}
По теореме 12.14, неавтономным первым интегралом на  
множестве $\Omega_1^{}=\{(t,x,y)\colon x\ne 1\}$ 
системы Якоби (12.83) будет функция
\\[1.75ex]
\mbox{}\hfill  % (12.84)
$
\Psi_{1}^{}\colon (t,x,y)\to\ 
\dfrac{x-y-1}{x-1}\, e^{t}
\quad
\forall (t,x,y)\in\Omega_1^{}.
$
\hfill (12.84)
\\[1.75ex]
\indent
Собственному числу $\lambda_1^{}={}-1$ матрицы $A$
соответствуют двукратный элементарный делитель, 
собственный вектор $(1, {}-1,{}-1)$ и 
первый присоединенный вектор $(1, {}-1, {}-2).$
По теореме 12.16, 
\vspace{0.15ex}
неавтономным первым интегралом системы Якоби (12.83) на  
мно\-жест\-ве $\Omega_2^{}=\{(t,x,y)\colon x-y-1\ne 0\}$ будет функция
\\[1.75ex]
\mbox{}\hfill  % (12.85)
$
\Psi_{2}^{}\colon (t,x,y)\to\ 
\dfrac{x-y-2}{x-y-1}\,-\,t
\quad
\forall (t,x,y)\in\Omega_2^{}.
$
\hfill (12.85)
\\[1.75ex]
\indent
\!\!
Каждая из совокупностей 
\vspace{0.25ex}
$\!\!(12.57)\!\& (12.84), (12.57)\!\& (12.85), (12.84)\!\& (12.85)\!\!$
образует ин- тегральный базис на множестве
$\!\Omega_0^{}\!=\!\{(t,x,y)\colon\! (x\!-\!1)(x\!-\!y\!-\!1)\!\ne 0\}\!\!$ 
\vspace{0.5ex}
системы Якоби\! (12.83).

{\bf Пример 12.13.}
У системы Якоби
\\[1.75ex]
\mbox{}\hfill    % (12.86)
$
\dfrac{dx}{dt}=x-x\;\!(y+1),
\quad
\dfrac{dy}{dt}=y-y\;\!(y+1)
$
\hfill (12.86)
\\[1.75ex]
матрица $A$ имеет (пример 12.6) 
\vspace{0.15ex}
собственное число, которому соответствуют двукратный элементарный делитель,
собственный вектор $(0, 1, 0)$ и 
первый присоединенный вектор $(1, 1, 1).$
\vspace{0.25ex}
По теореме 12.16, неавтономным первым интегралом на  
множестве $\Omega_1^{}=\{(t,x,y)\colon y\ne 0\}$ 
системы Якоби (12.86) будет функция
\\[2ex]
\mbox{}\hfill  % (12.87)
$
\Psi\colon (t,x,y)\to\ 
\dfrac{x+y+1}{y}\,-\,t
\quad
\forall (t,x,y)\in\Omega_1^{}.
$
\hfill (12.87)
\\[2ex]
\indent
Совокупность $(12.59)\& (12.87)$
образует интегральный базис системы Якоби (12.86) 
на множестве $\Omega_0^{}=\{(t,x,y)\colon xy\ne 0\}.$ 
\vspace{0.5ex}

{\bf Пример 12.14.}
У системы Якоби
\\[1.75ex]
\mbox{}\hfill    % (12.88)
$
\dfrac{dx}{dt}=x-y-\delta-x\;\!(x+y+2),
\quad
\dfrac{dy}{dt}=x+3y+\delta-y\;\!(x+y+2)
$
\hfill (12.88)
\\[1.75ex]
с вещественным параметром $\!\delta\!\ne\! 0\!$ 
\vspace{0.15ex}
матрица $\!A\!$ имеет (пример 12.7) 
собственное число, ко\-торому соответствуют трехкратный элементарный делитель,
собственный вектор $\!\!(1, 1, 0),$ 
первый присоединенный вектор $(0, 0, 1).$
\vspace{0.25ex}
По теореме 12.16, неавтономным первым интегралом на  
множестве $\Omega_0^{}=\{(t,x,y)\colon x+y\ne 0\}$ 
системы Якоби (12.88) будет функция
\\[1.5ex]
\mbox{}\hfill  % (12.89)
$
\Psi\colon (t,x,y)\to\ 
\dfrac{1}{x+y}\,-\,t
\quad
\forall (t,x,y)\in\Omega_0^{}.
$
\hfill (12.89)
\\[1.5ex]
\indent
Совокупность $(12.67)\& (12.89)$
образует интегральный базис на множестве $\Omega_0^{}$ 
системы Якоби (12.88).
\vspace{0.25ex}

Изложенный метод построения интегрального базиса обыкновенной дифференциальной 
системы Якоби второго порядка в монографии [1, с. 273 -- 311] распространен на 
систему Якоби в полных дифференциалах, частным случаем которой является 
обыкновенная дифференциальная система Якоби высшего порядка.
\vspace{0.25ex}

Заметим, что в работе [56] приведена обширная библиография статей, 
в которых рассматривались частные интегралы и их приложения.
Более поздние результаты по теории частных интегралов и их приложений
можно найти в работах [57 -- 72].

\newpage

\mbox{}
\\[-5.75ex]

}

\end{document}